%% file: book.tex
\documentclass[12pt,reqno,twoside,openright]{amsbook}
\usepackage[all]{xy}
\UseComputerModernTips
\CompileMatrices

\usepackage{amssymb}
\usepackage{eufrak}
\usepackage[dvips]{graphics}
\usepackage{setspace}
\usepackage{fontenc}
\usepackage{mathrsfs}
\usepackage{MnSymbol}

\newcommand{\head}[1]{\medbreak\noindent\textbf{#1}}

\makeindex

\swapnumbers
\theoremstyle{plain}
\newtheorem{theorem}{Theorem}[chapter] 
\newtheorem{proposition}[theorem]{Proposition}
\newtheorem{corollary}[theorem]{Corollary}
\newtheorem{claim}[theorem]{Claim}
\newtheorem{lemma}[theorem]{Lemma}
\newtheorem{conjecture}[theorem]{Conjecture}

\theoremstyle{definition}

\DeclareMathOperator{\Sym}{Sym}
\DeclareMathOperator{\FSym}{FSym}
\DeclareMathOperator{\Aut}{Aut}

\DeclareMathOperator{\FAut}{FAut}
\DeclareMathOperator{\BSym}{BSym}
\DeclareMathOperator{\Alt}{Alt}
\DeclareMathOperator{\Out}{Out}
\DeclareMathOperator{\Spin}{Spin}
\DeclareMathOperator{\Inn}{Inn}
\DeclareMathOperator{\SAut}{SAut}
\DeclareMathOperator{\DAut}{DAut}
\DeclareMathOperator{\BAut}{BAut}
\DeclareMathOperator{\BdAut}{BdAut}

\DeclareMathOperator{\Gal}{Gal}
\DeclareMathOperator{\ind}{ind}

\DeclareMathOperator{\col(}{col(}
\DeclareMathOperator{\cl}{cl}
\DeclareMathOperator{\dom(}{dom(}

\DeclareMathOperator{\im}{im}
\DeclareMathOperator{\id}{id}

\DeclareMathOperator{\max(}{max(}

\DeclareMathOperator{\Th(}{Th(}
\DeclareMathOperator{\card}{card}
\DeclareMathOperator{\Cay}{Cay}

\DeclareMathOperator{\diag}{diag}
\DeclareMathOperator{\supp}{supp}
\DeclareMathOperator{\ker(}{ker(}

\DeclareMathOperator{\Suz}{Suz}
\DeclareMathOperator{\Cr}{Cr}
\DeclareMathOperator{\Dr}{Dr}
\DeclareMathOperator{\Cox}{Cox}
\DeclareMathOperator{\Ts}{Ts}
\DeclareMathOperator{\W(}{W(}
\DeclareMathOperator{\PSL(}{PSL(}
\DeclareMathOperator{\PGL(}{PGL(}
\DeclareMathOperator{\GL(}{GL(}
\DeclareMathOperator{\Hyp}{Hyp}
\DeclareMathOperator{\Co}{Co}
\DeclareMathOperator{\Comm}{Comm}
\DeclareMathOperator{\SL}{SL}
\DeclareMathOperator{\Iso}{Iso}
\DeclareMathOperator{\End}{End}

\DeclareMathOperator{\soc}{soc}

\DeclareMathOperator{\Mlt}{Mlt}
\DeclareMathOperator{\RMlt}{RMlt}
\DeclareMathOperator{\LMlt}{LMlt}
\DeclareMathOperator{\SA}{SAut}
\DeclareMathOperator{\DA}{DAut}
\DeclareMathOperator{\BA}{BAut}
\DeclareMathOperator{\AAut}{AAut}

\DeclareMathOperator{\TAut}{TAut}
\DeclareMathOperator{\NA}{NA}
\DeclareMathOperator{\NB}{NB}
\DeclareMathOperator{\NE}{NE}
\DeclareMathOperator{\NS}{NS}
\newcommand{\symd}{\bigtriangleup}
\newcommand{\sd}{\rtimes}
\DeclareMathOperator{\fix}{fix}
\DeclareMathOperator{\tp}{tp}
\DeclareMathOperator{\Hom}{Hom}
\DeclareMathOperator{\homo}{homo}

\DeclareMathOperator{\diam}{diam}

\DeclareMathOperator{\Cos}{Cos}
\DeclareMathOperator{\Core}{Core}
\DeclareMathOperator{\Gr(}{Gr(}
\DeclareMathOperator{\Cl(}{Cl(}

\DeclareMathOperator{\acl}{acl}
\DeclareMathOperator{\lm}{lm}

\DeclareMathOperator{\dist}{dist}

\DeclareMathOperator{\Oct(}{\Oct(}

\DeclareMathOperator{\SO(}{SO(}
\DeclareMathOperator{\PSO(}{PSO(}
\DeclareMathOperator{\Bn}{Bn}

\newcommand{\Wr}{\mathop{\mathrm{Wr}}}

\newcommand{\onto}{\twoheadrightarrow}

\def \ds {\displaystyle} 
\def \ns {{}^{\ds .}}

\begin{document}

\title{Multicoloured Random Graphs: Constructions and Symmetry}
\vspace{7cm}
\author{Sam Tarzi\\[3pt]
London, England}
\date{Draft, 15 June 2014}
\maketitle

\pagenumbering{roman}

\cleardoublepage

\chapter*{Copyright}


\bigskip

\centerline{\copyright\  2014 Sam Tarzi}
\centerline{All Rights Reserved}

\vspace{1cm}

\begin{flushleft}
Sam Tarzi is identified as the author of this work in accordance with the Copyright, Designs and Patents Act, 1988.
\end{flushleft}

\vspace{1cm}

\begin{flushleft}
No part of this publication may be reproduced or stored in a retrieval system, nor transmitted in any form or by any means - digital scanning, mechanical, photo-copying, recording, or by any other media - without prior permission of the copyright owner.
\end{flushleft}

\vspace{1cm}

\begin{flushleft}
Every effort has been made to trace copyright ownership and permissions are confirmed in the text.  Any omission brought to the notice of the publishers will be acknowledged in future impressions.
\end{flushleft}

\clearpage
\vspace{3cm}

\bigskip

\chapter*{Abstract}

\centerline{Multicoloured Random Graphs: Constructions and Symmetry}

\bigskip

This is a research monograph on constructions of, and group actions on countable homogeneous graphs, concentrating particularly on the simple random graph and its edge-coloured variants.  Taking this family of graphs as our starting point, we study various aspects of the graphs, but the emphasis has been on illuminating those groups that are supported by these graphs, together with links with other structures such as lattices, topologies $\&$ filters, rings $\&$ algebras, metric spaces, sets $\&$ models, Moufang loops and monoids.  All graphs that we consider have a countable number of vertices, but some of the associated groups have uncountable order.  The large amount of background material included serves as an introduction to the theories that are used to produce the new results.  The large number of references should help in making this a resource for anyone interested in beginning research in this or allied fields.  Our exposition is intended to be thorough rather than encyclopedic. 

The countably infinite random graph which we denote by $\mathfrak{R}$ is known to be unique up to isomorphism; in other words a countable random graph is isomorphic to $\mathfrak{R}$ with probability $1$.
Its most common variant is as a simple graph, thus having only two
possible adjacencies, edges and non-edges, no loops or multiple edges
being permitted.  As we shall see there are three equivalent ways of characterizing $\mathfrak{R}$ up to isomorphism as a countable graph, these being the so-called one-point extension property, the injection property or a combination of universality and homogeneity.

There is an obvious generalization to graphs with any finite or infinite number of different adjacency types, which we realize as colours.  Our study of these generalized graphs demonstrates that whilst the basic properties such as homogeneity, universality and various stability conditions remain the same as those of $\mathfrak{R}$, the infinite-vertex $m$-coloured random graphs $\mathfrak{R}_{m, \omega}$ for $m \geq 3$ support groups with a different structure and are associated with rings of a different structure, to those of $\mathfrak{R}$.  Some of these differences are highlighted.  We also devote some attention to uncovering properties that are unique to the graph with three adjacency types which we call the \emph{triality graph}
\index{graph ! triality}%
 denoting it $\mathfrak{R^{t}}$.  Much of the work on the random graph and its relations over the past half-century has been on combinatorial aspects of the theory, whilst our exposition focuses on complementary questions.  

\clearpage

\chapter*{Dedication}


\vspace{4cm}

\centerline{\textbf{\Large In loving memory of my late father Anwar}} 

\bigskip

\centerline{\textbf{\Large my mother Adiba}}

\bigskip



\centerline{\textbf{\Large Jacob Avner and Daniel}}

\clearpage

\tableofcontents





\chapter*{Preface}

Graph theory is regarded as having begun with Euler's
\index{Euler, L.}%
solution to the K\"onigsberg bridge problem
\index{K\"onigsberg bridge problem}%
 in $1736$; Euler was possibly the most prolific
mathematician of all time~\cite{dun}, (though both Paul Erd\H{o}s and Serge Lang may stake a claim to this description).  The historical development of infinite discrete mathematics was such that~\cite{jecha}
it took another two centuries until Erd\H{o}s and A. R\'enyi
\index{R\'enyi, A.}%
\index{Erd\H{o}s, P.}%
 published the unique defining relation for the infinite graph.  Their papers
beginning in $1959$ built upon the notion discovered by Erd\H{o}s in the
$1940$s but that had lay dormant for a while, by giving a non-constructive existence proof of what
became known as the random graph
\index{graph ! random}%
 which we denote $\mathfrak{R}$, and by proving its
uniqueness.  The binomial model of the random graph due to E. N. Gilbert
\index{Gilbert, E. N.}%
 also dates back to $1959$.  Anatol Rapoport
\index{Rapoport, A.}%
 began a study of random-biased nets
in~\cite{solomra}, in which the emergence of a giant cluster in a
certain limit is derived, though the derivation is heuristic.  Rapoport
\index{Rapoport, A.}%
 introduced \emph{biases} to provide
more realistic models than can be achieved by using pure randomness.  The theory is summarized in~\cite{rapoport}.  More on the origins of random graph theory can be found in~\cite{kar}.

In this work we focus on studying versions of the random graph with more than
two adjacency types, which for clarity we denote by colours.  We
have devoted much attention to the countably infinite graph with three
adjacency types, calling it
the \emph{triality graph},
\index{graph ! triality}%
denoting it $\mathfrak{R^{t}}$ and using
the symbols $\mathfrak{r}, \mathfrak{b},\mathfrak{g}$ for its red,
  blue and green edges.  Many of
our results apply more generally to multicoloured generalizations of
$\mathfrak{R}$ having any finite or infinite number of edge-colours.

A simple graph is a relational structure with one binary
symmetric irreflexive relation.  One feature which persists in the transition from simple graphs with edges and non-edges to graphs on any number of colours is the existence
of a structure that is countable, universal
\index{structure ! universal}%
 (that is contains all finite substructures) and homogeneous
\index{structure ! homogeneous}%
(that is for which any isomorphism between finite substructures extends to an automorphism of the entire structure)~\cite{f}~\cite{jonsson1}~\cite{morley}.  Let us say a little about these two properties.

Firstly universality.  Two examples of \emph{categories}, that is structures and accompanying natural maps, leading to familiar universal objects are the collection of countable linear orders
\index{linear order}%
 and order-preserving maps having the linear order of rationals as universal object, and countable boolean algebras
\index{boolean algebra}%
 together with boolean algebra isomorphisms that have the countable boolean algebra on countably many generators as universal.  The universal objects may be of interest in their own right.  One aim of studying a universal structure is the possibility of finding a unifying theme for the individual questions of interest as well as the chance that the universal object may have properties that influence those of the embedded substructures as well as being influenced by them.  So in addition to the properties of the universal object itself, we may gather information that is useful in classifying the finite sub-objects. For the two examples just mentioned proof of universality gives us the theory of dense linear orders
\index{linear order}%
 without endpoints and the result that a finitely generated subalgebra of a boolean algebra is finite.  
\index{boolean algebra}%

Universal objects are seldom unique, because adding one point results in a different universal object, but universal homogeneous objects are unique.

There are certainly instances where universal objects fail to exist, for example for the category of countable fields and mappings preserving 0, 1, and the field operations.  In graph theory, there is no countable planar graph which is universal with respect to one-to-one graph homomorphisms, nor is there a universal graph of fixed finite degree under \emph{weak embeddings},
\index{graph ! weak embedding}%
 this being a $1$--$1$ mapping $i$ from the vertices of graph $\Gamma_1$ to the vertices of graph $\Gamma_2$ such that if $x, y$ are neighbours in $\Gamma_1$ then they are also neighbours in $\Gamma_2$.  The paper by Moss~\cite{moss1} contains more on the existence and nonexistence of universal graphs. 
\index{Moss, L. S.}%

 A variation on this theme is the idea of \emph{homomorphism-universal} structures
\index{structure ! homomorphism-universal}%
 whose relationship with embedding-universal structures is given in~\cite{hubickanes1}.  Briefly, if a homomorphism is $1$--$1$ then it is an embedding.  Another example of a phenomenon where a global structure provides a description of the possible local structures of graphs is that of \emph{limit graph}~\cite{giudici1}; see Appendix~\ref{PermutationGroups}.
\index{graph ! limit}%

\smallskip

Next homogeneity.  To say that a structure is homogeneous
\index{structure ! homogeneous}%
 is equivalent to asserting that it has the maximum amount of symmetry; or that many parts of the structure look alike; or that the structure looks the same when viewed from many positions within the structure.

We have mentioned the uniqueness of the countable random graph;
\index{graph ! random}%
 it also has a huge amount of symmetry.  To see how remarkable and unexpected this is, contrast the situation with the finite case.  Construct a finite random graph on $n$ vertices by ordering the pairs of vertices in a countable sequence, and tossing a fair coin a total of $\frac{n(n-1)}{2}$ times where heads means join the two vertices by an edge, and tails means do nothing.  Every $n$-vertex graph (up to isomorphism) occurs with non-zero probability, and the probability of a particular graph arising is inversely proportional to its symmetry as measured by the order of its automorphism group.  The asymmetric graphs are overwhelmingly more numerous than the symmetric ones.  The probability measure associated with a countable sequence of coin tosses is discussed in Appendix~\ref{CategoryandMeasure}.  To get a feel for the meaning of homogeneity beyond its definition consider that the pentagon or 5-cycle graph is homogeneous, but that the hexagon or 6-cycle graph is not because a pair of points two steps apart and a pair three steps apart are isomorphic as induced subgraphs but are not equivalent under automorphisms of the whole graph.
  
\smallskip

The existence and uniqueness of homogeneous universal structures follow immediately from Fra\"{\i}ss\'e's Theorem~\cite{f} on relational structures.
\index{Fra\"{\i}ss\'e's Theorem}%
\index{relational structure}%
As we stated, homogenous structures are those with the maximum amount of symmetry and one of the most useful methods for constructing objects with a large amount of symmetry is based on this theorem.  In the two-colour case, thinking of a
two-coloured complete graph as a simple graph, the appropriate structure is
\emph{Rado's graph},
\index{graph ! Rado's}%
or
the Erd\H{o}s-R\'enyi \emph{random graph}
\index{graph ! random}%
~\cite{er} $\mathfrak{R}$; it
was Richard Rado
\index{Rado, R.}%
~\cite{rado} who gave the first construction of
$\mathfrak{R}$. 

Fra\"{\i}ss\'e's Theory has been extended in several directions.  For example, a study of relational structures which satisfy the analogue of homogeneity but for homomorphisms rather than isomorphisms was begun in~\cite{camnes}.  We include some material on some of the developments in an appendix on Fra\"{\i}ss\'e's Theory of Relational Structures.

\smallskip

Multi-adjacency random graphs have been analysed before, as countable universal homogeneous
$C$-coloured graphs
\index{graph ! homogeneous}%
\index{graph ! universal}%
 $\Gamma_{C}$ for $|C| \geq 2$,
using permutation group theory~\cite{truss1}~\cite{truss2}.  In another appendix we pr\'ecis some of the results from previous research
work.

The random graph is unique up to isomorphism and random graphs with
different numbers of colours will be pairwise non-isomorphic.
Allowing colours on vertices yields different graphs and in Chapter 2 we
comment briefly on random graphs with two-coloured vertices.

Much of the impetus for the renaissance of discrete mathematics
beginning in the
second half of last century has been finite or finitary in origin and
often guided by real-world problems.  This does not diminish the
contribution of questions that are infinite in their nature.
An apposite quote is the first sentence of the relevant chapter
in the book~\cite{bollobas}: ``Although the theory of random graphs is one of the
youngest branches of graph theory, in importance it is second to
none.''    The two-coloured random graph
\index{graph ! random}%
  has been studied using combinatorial~\cite{biggs}, probabilistic~\cite{bollobas1}, logical~\cite{spencer} and spectral techniques.  The study of random graphs on any number of colours greater than two adds further scope to the field.

Most of the work in random graph theory so far has blended
probability theory, combinatorics
\index{combinatorics}%
 and logic, but our focus has been on graph constructions and group actions on the graphs.  A group action can uncover some of the detailed structure and properties of an operand by seeing what actions it supports, and at the same time it can illuminate the characteristics of the acting group by giving some of its representations.  This is particularly the case if a mathematical object can be uniquely determined by the algebra of functions that it can support.

In studying the infinite counterpart of a finite structure, it is
natural to let one of the parameters go to infinity.  This 

Taking the infinite limit of a parameter of in a finite structure cannot be done arbitrarily; in discrete mathematics it may cause axiom of choice
\index{axiom of choice}%
problems and in continuous mathematics, convergence issues may arise.
Another reason is that certain properties may be lost in taking the
limit.  For example, the linear group $GL(V)$
\index{group ! linear}%
acting on the countably infinite vector space
$V=(\aleph_{0}, K)$ over a finite field $K$ is \emph{oligomorphic}
\index{group ! oligomorphic}%
 (meaning that the number of
$GL(V)$ orbits
\index{group ! orbit}%
on the set of $n$-subsets of $V=(\aleph_{0}, K)$ is
finite for all $n$); the field must
be finite, because fixing a vector fixes all scalar multiples of it, but an
infinite field gives an infinite number of scalar multiples,
contradicting oligomorphicity.  However oligomorphicity is \emph{not} lost
when the dimension of the space goes from finite to infinite.  In
fact if $K$ is finite then $V=(\aleph_{0}, K)$ is totally categorical
(these are the best behaved of the four classes in the categoricity
\index{categoricity}%
spectrum implied by Morley's theorem in \emph{model theory}),
\index{model theory}%
\index{Morley's Theorem}%
whilst it is uncountably
categorical if $K$ is countably infinite.  (A set of sentences in a model is \emph{categorical in power}
\index{categorical in power}%
 $\lambda$ (an infinite cardinal) if any two models of the set of cardinality $\lambda$ are isomorphic.)  Morley showed that a set of sentences over a countable language which is categorical in some uncountable power is categorical in all; Shelah extended the notion from first-order
\index{first-order logic}%
 to infinitary logic, where the Compactness Theorem
\index{Compactness Theorem}%
fails~\cite{baldwin2}.  These are all notions from logic and in particular model theory, and we defer discussions these concepts to the appendices.  


The word \emph{oligomorphic}
\index{group ! oligomorphic}%
is intended to capture the notion of `few shapes', where the group orbits contain a finite number of structures of any given finite size (that is, few) and each orbit contains an isomorphism class of structures (that is, shapes).

One systematic way to go from the finite to the infinite is using the model
theory of relational structures
\index{relational structure}%
 developed by Roland Fra\"{\i}ss\'e~\cite{frai}.
\index{Fra\"{\i}ss\'e, R.}%
  This theory is widely-applicable to many different types of structures, but Fra\"{\i}ss\'e's main interest was in relational ones.  If $N$ is a finite substructure of any homogeneous structure $M$, then all automorphisms of $N$ extend to automorphisms of $M$ fixing $N$.  If
the converse holds then $N$ is said to be a \emph{finite homogeneous substructure} of $M$.
\index{finite homogeneous substructure}%
 Rose and Woodrow
\index{Rose, B. I.}%
\index{Woodrow, R. E.}%
  use the term \emph{ultrahomogeneous structures} if every isomorphism between substructures of a smaller cardinality extends to an automorphism.  In~\cite{rose} they prove a sufficient condition for their existence as well as others on connections between such structures and \emph{quantifier elimination (q.e.)} which defined in the appendix on model theory.
\index{quantifier elimination (q.e.)}%


Model theory
\index{model theory}%
is a part of pure mathematics which is very general in its notions and formalism, and simultaneously very good at solving problems when specialised.  
Its methods are increasingly being used to great eefect in other areas of mathematics.   The proofs of theorems of G\"odel-Deligne, Chevalley-Tarski, Ax-Grothendieck, Tarski-Seidenberg, and Weil-Hrushovski
\index{Hrushovski, E.}%
 which were an intersection of model-theoretic methods and techniques from other areas; see reference~\cite{baldwin1}.  Other crossovers can be found in the work of Ax-Kochen-Ershov, Macintyre, Denef and du Sautoy.
 The important point is that such techniques were used
before or instead of the usual methods of say field theory or algebraic
geometry.  There is a claim that algebraic geometry is in some sense a canonical special case of model theory,
\index{model theory}%
 for while not everything can be reduced to algebraic geometry,
the classical functions and the structures they give rise to are nice
extensions of it.  So there may be a generalization or a re-writing of deep parts of mathematics, for instance the geometric Langlands programme using model-theoretic language.  The tools we use are largely group-theoretic, many of the ideas such as reducts and the small index property,
\index{small index property (SIP)}%
are at heart model-theoretic.  A comprehensive treatment of the field can be found in~\cite{hod1}.

In Zermelo-Fraenkel set theory,
\index{Zermelo-Fraenkel set theory}%
$ZF$, \emph{$\aleph_0$-categoricity} of $\mathfrak{R}$ means that there is essentially
only one countable model of set theory with symmetrised set-inclusion $\in$.  $\aleph_0$-categoricity is a property of the automorphism group
\index{group ! automorphism}%
 of a structure and implies a very high degree of symmetry.  Homogeneity is an even stronger condition.  

Relational structures
\index{relational structure}%
 with large automorphism groups have a high degree of
symmetry; this is also implied by an orbit space with a small number of orbits of the
acting automorphism group.  A countable homogeneous relational
structure of a given type is determined up to isomorphism by the
isomorphism classes of its finite substructures.  An important outcome of Fra\"{\i}ss\'e's Theorem
\index{Fra\"{\i}ss\'e's Theorem}%
 is that determining the numbers of orbits on $n$-sets or $n$-tuples of \emph{oligomorphic} permutation groups
\index{group ! permutation ! oligomorphic}%
(those having finitely many orbits in their induced action on $n$-tuples for all $n$) is equivalent to enumerating the unlabelled or labelled objects in certain classes of finite structures characterized largely by the \emph{amalgamation property}, which is a way of building larger structures from smaller ones (see Appendix~\ref{TheoryofRelationalStructures}).  Not all graph classes amalgamate; for example, $n$-colourable graphs (those having chromatic number at most $n$) do not have the amalgamation property.  

The model-theoretic equivalent of the concept of oligomorphic permutation group
\index{group ! permutation ! oligomorphic}%
is given by the Engeler--Ryll-Nardzewski--Svenonius Theorem~\cite{engeler}~\cite{ryll}~\cite{svenonius}
\index{Engeler--Ryll-Nardzewski--Svenonius ! Theorem}%
\index{Engeler, E.}%
\index{Ryll-Nardzewski, C.}%
\index{Svenonius, L.}%
which states that a countable (first-order) structure
\index{first-order structure}%
 is axiomatisable (that is characterised up to isomorphism as a countable structure by first-order sentences)
\index{first-order sentence}%
if and only if its automorphism group
\index{group ! automorphism}%
 is oligomorphic.  The surprising conclusion of this theorem is that \emph{axiomatisability is equivalent to symmetry}.  This result appears to be an \emph{exception}.  We know of no other such equivalences in model theory;
\index{model theory}%
 it does not apply to higher cardinalities.  Remarkable though the theorem is, it adds nothing to our knowledge of  finite structures, for every finite permutation group is oligomorphic, and every finite first-order structure is categorical.  Indeed, if the automorphism group of the $n$-element structure $M$ has a base of size $k$, then every point of $M$ is uniquely identified by a formula having the elements of the base as parameters.  (If $G$ is a permutation group on a set $\Omega$, a $base$
\index{group ! permutation ! base}%
 for $G$ is a list $B$ of elements of $\Omega$ whose pointwise stabilizer in $G$ is the identity.  One philosophy is that a list of points of the structure $M$ which is a base for the 
automorphism group
\index{group ! automorphism}%
 of $M$ should be in some sense a base for $M$; this holds for first-order structures.  That is, if we cannot move $x$ to $y$ by an automorphism fixing $B$ pointwise, this is because the structures of $(M, [B, x])$ and $(M, [B, y])$ are different.  For more on bases, see~\cite{bailcam} and~\cite{cameron2}.)

The occurrence of the random graph
\index{graph ! random}%
 paradigm in various areas of
mathematics is well-known~\cite{cameron}, and perhaps less famously in
models in physics, see for example~\cite{ant}~\cite{childs}~\cite{filk}~\cite{req}.  What
is typically being asked for is the form of `most' graphs in families with
certain properties.  This is the essence of the so-called
probabilistic method~\cite{alonspencer}.
\index{Spencer, J.}%
\index{Alon, N.}%
The underlying principle of random graph theory can in a heuristic way be considered to
be diametrically opposite to the principle underlying chaos theory, where neighbouring trajectories diverge exponentially for long-term behaviour, that is on an attractor.  

The random graph
\index{graph ! random}%
 construction has been extended to more general first-order relational structures~\cite{droste}.
\index{first-order structure}%
\index{relational structure}%
  It can also be viewed as a special case of the
Borel-Cantelli lemmas
\index{Borel-Cantelli lemmas}%
concerning infinite sequences of trials, for
example Bernoulli trials~\cite{feller}; see Appendix~\ref{CategoryandMeasure}.
\index{Bernoulli trials}%
That such constructions can arise both via Borel-Cantelli or from a Kolmogorov
\index{Kolmogorov, A. N.}%
zero-one law,
\index{zero-one law}%
 points to the significance of probability in proofs of existence, uniqueness as well as in building mathematical structures.  In fact it can be shown~\cite{droste} that if a
class of countable relational structures
\index{relational structure}%
 contains an infinite $\aleph_0$-categorical universal homogeneous structure $M$, then $M$
can be constructed probabilistically.  A wide range of systems that arise as complex networks have been realized as random graphs~\cite{albertbar}.
\index{graph ! random}%
  Many random graph properties have no unique vertex-independent threshold, but rather one that depends on the system size, with critical probability obeying $p_0(n \to \infty) \to 0$.  We can think of a random graph \emph{evolution} or \emph{emergence} as the graph develops from a starting set of isolated vertices by the addition of random edges, that is as a changing dynamical network.  The alternative passive approach, to which we have tacitly adhered, is to view random graphs as a family of graphs to which is added a probability distribution, thus turning it into a probability space.  

As the graphs grow they can undergo a phase transition which for example leads from small components to a giant 
\index{giant component}%
 component~\cite{spencer1}.  Threshold phenomena are neither peculiar to a particular definition of randomness nor to graphs, occurring in a variety of random combinatorial structures in the limit of large size.  The existence of a limit as the number of parts or interacting units in the system becomes very large is however crucial to the possibility of a sudden change in a certain property as a function of smooth changes of a parameter~\cite{cohenj}.  Whilst our theory is not an application of threshold effects, we do also utilize the idea of the large system limit, as we will see.  

Random graph theory offers the opportunity to work along both the finite-infinite and the discrete-continuous boundaries of mathematics.  Examples of this symbiosis abound, for example the exciting work on `continuous
 graphs'~\cite{aldous}~\cite{borgs}~\cite{diaconisjan}~\cite{lovasz}~\cite{petrovv}. 

After concepts such as function and space, the concept of \emph{group}
is one of the most ubiquitous ideas in mathematics
and permeates almost all branches, often in a fundamental way.  We have concentrated our efforts on studying the types of groups whose action is supported by random graphs, including the random graph
\emph{reducts}
\index{group ! reduct}%
 (or closed supergroups of the automorphism groups)
\index{group ! automorphism}%
 and their multicoloured variants, which we felt was a good starting point from which to uncover symmetry properties of the graphs.  Closure as the automorphism group of some first-order structure
\index{first-order structure}%
  is equivalent to closure in the topology of pointwise convergence~\cite{cameron}.

Switching of undirected graphs was introduced by Van Lint
\index{Van Lint, J. H.}%
 and Seidel~\cite{vanlint}
\index{Seidel, J. J.}%
  in their study of a problem concerning equilateral $n$-tuples of points in elliptic geometry. 

According to a theorem of Simon Thomas~\cite{thomas} the group of switching automorphisms of $\mathfrak{R}$ is one of its three non-trivial reducts,
\index{Thomas, S.}%
  a fact that can be explained by use of the \emph{two-graph} $\mathcal{T}(\mathfrak{R})$,
\index{graph ! two-graph}%
\index{two-graph}%
 a relational structure
\index{relational structure}%
  derived from $\mathfrak{R}$ whose automorphism group
\index{group ! automorphism}%
 is $\SA(\mathfrak{R})$, the group of switching
automorphisms.  The extra scope gained by consideration of more than two adjacencies offers new possibilities.  In the three-coloured case the equivalent group of
switching automorphisms is not a reduct, for it is highly transitive
but not the symmetric group;
\index{group ! symmetric}%
a group is \emph{highly transitive}
\index{group ! permutation ! highly transitive}%
if it
acts on an infinite set and is $k$-transitive, meaning that it maps any $k$-tuple to any other $k$-tuple for all integers $k \ge 1$.  However there are switching automorphism
reducts for multicoloured random graphs,
\index{graph ! random ! $m$-coloured}%
 which together with duality-type reducts are classified in the Ph.D. thesis of James
H. Bennett~\cite{bennett}.
\index{Bennett, J.H.}%

Some of the new results provide evidence that the random two-coloured graph has different switching properties to the random graph on three or more colours; broadly, the former has structure which appears to be diluted when more than two colours are considered.  So the random graph $\mathfrak{R}$
\index{graph ! random}%
 appears to have contrasting symmetry properties to those of the higher-adjacency graphs, as
witnessed in the results of several of the chapters.  This phenomenon where the low dimensional
cases of a theory are seen to exhibit special behaviour which is
lost in higher dimensions arises all over mathematics.
We tabulate some of the differences between the two-colour and $m$-colour $(m> 3)$ properties in an appendix.  It would be worth investigating whether or not this difference extends to a loss of the `strange logic'~\cite{spencer}
\index{Spencer, J.}%
that arises for $\mathfrak{R}$, in particular the phase transition in $\mathfrak{R}$ whereby
over a short range, events quickly turn from being almost certainly false to almost certainly true.  

\bigskip

J. K. Truss
\index{Truss, J. K.}%
 has studied countable universal edge-coloured graphs, and we list some of his results, leaving a more complete summary for Appendix~\ref{PreviousResults}.   In~\cite{truss1} he proved that these graphs have simple automorphism groups.
\index{group ! automorphism}%
  More concretely a certain subset $\Sigma$ of elements of these groups is defined such that: (a) If $\sigma\sb 1,\sigma\sb 2\ne 1$ then there is a conjugate $\tau$ of $\sigma\sb 1$ such that $\sigma\sb 2\tau \in \Sigma$, (b) if $\sigma\sb 1,\sigma\sb 2,\sigma\sb 3\in \Sigma$ there are conjugates $\tau\sb 1,\tau\sb 2,\tau\sb 3$ of $\sigma\sb 1,\sigma\sb 2,\sigma\sb 3$ respectively such that $\tau\sb 3\tau\sb 2\tau\sb 1=1$. From (a) and (b) it follows that if $\sigma$,$\tau$ $\ne 1$, then $\tau$ is the product of five conjugates of $\sigma$.  The result is proved by approximating the desired conjugacies by means of carefully chosen finite approximations. \par  The cycle types of elements of the automorphism groups  $\Aut(\mathfrak{R}_{m,\omega})$ are also classified.  Elements with cycle types $\infty\sp 1,\infty\sp 2$ and $p\sp{\infty}$ for prime $p$ are classified up to conjugacy. Finally a study is made of the relations between the groups $\Aut(\mathfrak{R}_{m,\omega})$,
\index{graph ! random ! $m$-coloured}%
for different $m$.

Truss looked at the highly transitive overgroup of $\Aut(\mathfrak{R}_{m,\omega})$, the \emph{group of almost automorphisms}
\index{group ! almost automorphism}%
 $\AAut(\mathfrak{R}_{m,\omega})$ of $\mathfrak{R}_{m,\omega}$, in~\cite{truss2}. A permutation of $\mathfrak{R}_{m,\omega}$ is said to be an almost automorphism if the set of edges whose colour it alters is finite.  He proved that $\AAut(\mathfrak{R}_{m,\omega})$ is not simple, described its normal subgroups and characterized the possible cycle types of members of $\AAut(\mathfrak{R}_{m,\omega})$.

\bigskip

In the main text we will study and in part extend the above-mentioned properties of the random graph to its edge-coloured relatives.

\medskip

This book is part research monograph, part introduction to and summary of ground already covered over decades written with the intention of drawing in audiences from elsewhere including beginning researchers.  The style of writing is intended to convey that multicoloured random graphs are a peg onto which to hang a great breadth of material.  Multicoloured random graphs do have distinguishable properties from the simple random graph, as listed in the Table in Figure 2 of Appendix E, and this book is just the beginning.  We have pointed out those issues where the two-coloured random graph illustrates the property under discussion.

The contents of this tractate are arranged as follows.  
In order to maintain the flow of the main chapters, we defer to the appendices an introduction to the objects of study and elements of the theories which we apply in the main text.

\bigskip

Chapter $1$
introduces the triality graph
\index{graph ! triality}%
 $\mathfrak{R^{t}}$ and proves using standard arguments some fundamental properties such as homogeneity, universality, and the pigeonhole property.  We end with some comments regarding triality graphs on vertices with two colours.

\bigskip

In $\S1$ of Chapter $2$ we define \emph{switching
groups}
\index{group ! switching}%
$S_{m,n}$ that act on the set of complete graphs with $n$ vertices and $m$ different adjacencies, as well as groups of switching automorphisms and identify the interrelationship of these two types of groups.  In the next two sections we define two-coloured switching, firstly for finite graphs and then for the infinite random graph $\mathfrak{R}$ stating Thomas' Theorem 
\index{Thomas' Theorem}%
that classifies its reducts.
\index{group ! reduct}%
After a section which lists the reducts of the triality graph, we turn our attention to random graphs with coloured vertices, commenting on (but not giving) a classification of the reducts of the \emph{random bipartite graph}
\index{graph ! random ! bipartite}%
$\mathfrak{B}^{v}$.  After a section outlining some of the known results on random graphs with forbidden substructures, we give a detailed examination of the switching group $S_{3,3}$, as motivation for the determination in Theorem~\ref{swgpfm} of the structure of the finite switching groups $S_{m,n}$.
\index{group ! switching}%
One section is devoted to the switching reduct for three colours.  We go on to study some of the properties of $S_{m, n}$, thereby extending the work of D. G. Higman
\index{Higman, D. G.}%
 and J. J. Seidel
\index{Seidel, J. J.}%
 for $m=2$.  We
prove in Theorem~\ref{extswgp} that for $m \geq 3$ the extended groups~\label{extsw}
defined by $G^W_{m,n}:= S_{m,n} \sd \Sym(n)$ generated by the
switchings and vertex permutations are primitive in their action
\index{group ! action}%
\index{group ! permutation ! primitive}%
on the set of complete graphs with $n$ vertices and
$m$ edge colours, whilst for $m=2$ these extended groups are not even
transitive in their action.  For any coloured graph $\Gamma$ the group
$\SAut(\Gamma)$ of switching automorphisms is the stabilizer
\index{group ! permutation ! stabilizer}%
 of $\Gamma$ in $G^W_{m,n}$, or
equivalently the group of permutations $g$ such that $\Gamma g =
\Gamma \sigma$ for some switching $\sigma$.  We also show, in Theorem~\ref{wte}, for $3 \leq n \leq \aleph_0$, $3 \leq m <
\aleph_0$ and for a given $n$, that $S_{m,n}$ is a weak transitive
extension of $S_{m-1,n}$.  In the case of a graph with a countably
infinite vertex set and a finite number of colours on edges, we show
that its group of switching automorphisms is the symmetric group on
the vertex set if and only if the graph is a so-called \emph{switched c-clique}
\index{graph ! switched $c$-clique}%
  for some colour $c$, where this latter concept is defined.  Different
definitions of switched $c$-clique are given and the equivalence of
these definitions is proved.  

\bigskip

In Chapter $3$ we investigate
variations on the switching group theme by allowing the random graphs
to have an infinite number of different colours,
\index{graph ! random ! $m$-coloured}%
 but where an element of the switching group can only change a finite number of
colours and such elements therefore generate what can be called a
\emph{finitary switching group}.
\index{group ! switching ! finitary}%
We prove some results characterizing profinite properties of these groups and local finiteness of related groups.  We investigate switching groups that are closed in a certain topology, identifying the parity structure which the group elements preserve.  This turns out to be different to the global parity-equivalence of graphs arising in Chapter $2$.  We end by showing the circumstances in which switching groups can form a direct limit.

Chapters $2$ and $3$ offer evidence that the theory of switching groups has developed somewhat from Seidel's
\index{Seidel, J. J.}%
original definition.  

\bigskip

In the first section of Chapter $4$ we outline the theory of \emph{local switching}~\cite{camseits} together with elements of the cohomological theory of switching classes and two-graphs.  In the next section we begin a formulation of the theory for multicoloured graphs, which involves working with a set of complete multicoloured graphs
\index{graph ! random ! $m$-coloured}%
 with a list of edge-colours, which we call an \emph{SML graph}.  We generalize the relevant closed switching group from the previous section together with the parity that it leaves invariant.  After a subsection outlining some known material on trees, we derive in two appendices two presentations for the switching group $S_{m,n}$.
 
\emph{We should mention that there are 3 types of colour parities appearing in our work.  Firstly  that in Chapter~\ref{chap2} (in the prelude to Theorem~\ref{psequiv}), which is defined on $3$-coloured complete graphs, where one of the colours is stabilized and the parity refers to that of the number of second-colour edges in triangles (vertex triples) on the second and third colours.  Secondly, in Chapter~\ref{chap3} (in the prelude to Theorem~\ref{parthm}), refers to a parity-presercing permutation of a list of colours on triangles in a complete $m$-coloured graph.  Thirdly, in Chapter~\ref{chap4} (in the prelude to Lemma~\ref{treeprop}) there is also a permutation parity of list colourings on edges of multicoloured graphs, but one which generalizes the second parity in being applicable in a slightly more general graph setting.}

\bigskip

Chapter $5$ has two parts.  The first deals with polynomial algebras
and polynomial invariants related to random graphs.  We reveal how the algebras
related to the two-coloured random graph reducts differ in their structure
to those related to random graphs with more colours.  

The second half of the chapter is devoted to the elucidation of an isomorphism between two algebras that have quite different structures, explaining why the isomorphism works for the very simplest infinite version of the algebras for pure point sets, whilst it fails for finite algebras.  We comment on the inexhaustibility of relational structures and on links between growth rates and reconstruction.  In the appendix, we begin showing the isomorphism in the case of 1, 2 and $3$-vertex graphs.

\bigskip

In Chapter $6$ various aspects of random graphs as Cayley graphs are examined.  

Firstly, we use the construction of $\mathfrak{R}$ as a Cayley graph
\index{graph ! Cayley}%
\index{Cayley, A.}%
to derive results concerning cyclic automorphisms
\index{automorphism ! cyclic}%
 and cyclic almost automorphisms of $\mathfrak{R}$.  Square root sets of elements of a structure such as a group are defined and a necessary and sufficient condition for a group $G$ to act regularly on $\mathfrak{R}$ is given, this being a translation of the one-point extension property,
\index{one-point extension property}%
  one of the equivalent ways of defining random graphs. 
 
We then construct $\mathfrak{R^{t}}$ as a homogeneous Cayley graph
\index{graph ! Cayley}%
for both an index $3$ subgroup of the Modular group
\index{group ! modular}%
 in Theorem~\ref{mdgptm} and the complex Leech Lattice
\index{lattice ! Leech ! complex}%
in Theorem~\ref{cplltm}.  More generally we construct random Cayley
$m$-coloured graphs from lattices $\mathbb{L}$
\index{lattice}%
 in $\mathbb{R}^d$, for $d \ge 2$, identifying lattice vectors with vertices
and vector pairs with edges such that the graphs are invariant under the lattice automorphism group $\Aut(\mathbb{L})$.  This is quite different to known links between groups, graphs and lattices, such as the Cayley graph for the free abelian group
\index{group ! free abelian}%
 on two generators being a lattice in the plane with vertex set $\mathbb{Z}^2$, because here the graphs are random.  
 
After a section defining so-called \emph{groups with triality}
\index{group ! with triality}%
we prove the existence of such a subgroup of $\Aut(\mathfrak{R^{t}})$, and further prove that $\Aut(\mathfrak{R^{t}})$ has a split extension.
 \index{group ! split extension}%
 Groups with triality are a generalization of what we have called Cartan triality, the algebraic triality outer automorphisms
\index{automorphism ! outer}%
 of the $8$-dimensional projective reduced 
\index{group ! orthogonal}%
 orthogonal group.
 
An extended discussion lead us to conclude that any direct approach using the Cayley object
\index{Cayley object}%
 method, to demonstrating a link between the two occurrences of triality, the normalizing action of $\Sym(\mathfrak{r} , \mathfrak{b} , \mathfrak{g})$ on $\Aut(\mathfrak{R^{t}})$ and a Cartan-type one is likely to fail.  This is because of the absence of a fixed-point-free
\index{group ! action ! fixed-point-free}%
 group action which we would require for an unambiguous graph construction.  We do make some progress in the steps required to demonstrate such a connection.  For example, in Theorem~\ref{moulthm} we give a construction of $\mathfrak{R^{t}}$ via a Moufang loop,
\index{Moufang ! loop}%
 whose \emph{multiplication group} is a
\index{group ! with triality}%
  group with triality.  The Cartan-like triality $\Sym(3)$ acts as outer automorphisms of $\Mlt(Q)$.

\bigskip

Chapter $7$ begins by defining outer automorphism groups
\index{group ! automorphism}%
 for random graphs and then proving a theorem showing that if $m$ is odd, an extension of $\Aut(\mathfrak{R}_{m,\omega})$ splits
\index{group ! split extension}%
  over this group, that is there is a \emph{complement}.  In the case that $m$ is even but not divisible by 8, then there is a subgroup $H$ of the extension which is a \emph{supplement} such that $|\Aut(\mathfrak{R}_{m,\omega}) \cap H| = 2$.  The formalism introduced by Alperin, Covington and Macpherson for sets~\cite{alcoma}
\index{Alperin, J. L.}%
\index{Covington, J.}%
\index{Macpherson, H. D.}%
is then applied to graphs, and this leads to the definition of two new types of groups, the group of near symmetries and the group of equivalence classes of near-automorphisms for the random graph.  Two further new groups are identified, the zero vertex index group and the finite vertex index group.  Maps of various types between these groups and known variations of random graph automorphism groups are indicated.

\bigskip

In Chapter $8$ we initiate a study of neighbourhood filters on sets of random graph vertices identifying some new topologies and groups whose basic features we outline.  After a motivating discussion of Sierpi\'nski's Theorem characterizing spaces homeomorphic to $\mathbb{Q}$ we study various topologies on random graphs.  A Hasse diagram of the novel groups that arise in this and the previous chapter is given.

An account of the Stone--\v{C}ech compactification
\index{Stone--\v{C}ech compactification}%
 of a topological space is given for neighbourhood filters on graphs.  

There is a section which introduces the $\mathfrak{R}$-uniform Hypergraph $\mathfrak{RHyp}$ of Claude Laflamme, Norbert Sauer and Maurice Pouzet and Robert Woodrow
\index{Laflamme, C.}%
\index{Sauer, N.}%
\index{Pouzet, M.}%
\index{Woodrow, R. E.}%
 defined on the vertex set of $\mathfrak{R}$ whose edges are those sets of vertices which induce a copy of $\mathfrak{R}$.  Groups related to this hypergraph are related to the groups that have arisen in this and the previous chapter.

\bigskip

The unifying theme of Chapter $9$ is various aspects of the interaction of homomorphisms and graphs, in both a finite and an infinite context.  Theories of algebraic structures, such as monoids, are very wide-ranging, so we focus on special classes of such structures.  

We give results on endomorphisms on graphs in general and the random graph in particular, as well as the monoids they generate.  Monoids have a more amorphous structure than groups, for example lacking the equivalent of Lagrange's Theorem
\index{Lagrange's Theorem}%
 which even exists~\cite{grish}~\cite{hallga} for Moufang loops
\index{Moufang ! loop}%
  (a finite Moufang loop $(Q,\cdot)$ has the Lagrange property if for each subloop $L$ of $Q$, $|L|$ divides $|Q|$).   As another example of the group-like properties of Moufang loops, which loops in general surely cannot be expected to obey, is the analog of the first Sylow theorem giving a criterion for the existence of a $p$-Sylow subloop of a finite Moufang loop, due to Grishkov and Zavarnitsine~\cite{grishkovza}.
\index{Grishkov, A. N.}%
\index{Zavarnitsine, A. V.}%
Finally we mention Glauberman's
\index{Glauberman, G.}%
 $Z^{*}$-theorem~\cite{glaubermana}, through which he was able to prove version of the Feit-Thompson Theorem, that is that Moufang loops of odd order are solvable~\cite{glauberman}.  For more on this aspect of Moufang loops, including the connection between Moufang loops and groups with triality, see section 4 of~\cite{hallji2}.  

Various relationships are studied between the \emph{core} of a finite graph $\Gamma$, which is the smallest graph which is homomorphically equivalent to $\Gamma$, and a concept which is in some sense dual to the core, that is the \emph{hull} of a graph $\Gamma$, which is a graph containing $\Gamma$ as a spanning subgraph, admitting all the endomorphisms of $\Gamma$, and having as core a complete graph of the same order as the core of $\Gamma$.  There is a section which comments on the equivalent theory for cores of  infinite graphs and which contains open questions.

The proof of a monoidal version of Fra\"{\i}ss\'e's Theorem
\index{Fra\"{\i}ss\'e's Theorem}%
for so-called HH structures is given.  As we shall see, a countable structure is HH if and only if it satisfies the so-called \emph{homo-extension property}, (and that a Fra\"{\i}ss\'e limit
\index{Fra\"{\i}ss\'e limit}%
 exists for its age which has the \emph{homo-amalgamation property}, though the limit may not be unique); all these terms are explained in the chapter.  

We further give an extension of the standard formalism allowing us to prove an equivalence theorem for (topological) endomorphism monoids and discuss how a weakening of the Engeler--Ryll-Nardzewski--Svenonius Theorem,
\index{Engeler--Ryll-Nardzewski--Svenonius ! Theorem}%
may lead to a characterization of oligomorphic monoids.  This discussion acts as motivation for just such a characterization given by Bodirsky and Pinsker
\index{Bodirsky, M.}%
\index{Pinsker, M.}%
  in~\cite{bodirsky}, which we briefly review. 
 

\bigskip

Chapter $10$ is on random graph constructions.  We illustrate some of the promising connections between number theory and random graph theory by giving number-theoretic constructions of both $\mathfrak{R}$
and $\mathfrak{R^{t}}$.  The next section is an account of known similarities between universal metric spaces and $\mathfrak{R}$.  This is followed by a section on homogeneous integral metric spaces, including references to previous work and the beginnings of a classification of reducts of such spaces.  After a section on miscellaneous observations focussing on aspects of set theory, we end with an appendix on results obtained by working in a model of set theory in which the axiom of choice
\index{axiom of choice}%
 is false.  

\medskip

The first ten chapters have open questions interspersed amongst them, but Chapter $11$ serves as a repository for a selection of further research directions which is intended to stimulate further work in this field.  Some of the questions are natural outgrowths of the text whilst others are more speculative.  

\medskip

The three major themes of our account afforded by the generalization to multi-coloured graphs are an extension of switching groups to a much greater depth than has thus far been studied,  constructions of random graphs in an array of other areas in mathematics, and the actions of a variety of new groups uncovered.

We have concentrated on (i) developing the theory of switching groups, their relatives and other groups supported by random graphs, (ii) constructing Cayley graphs of multicoloured random graphs
\index{graph ! random ! $m$-coloured}%
 and (iii) finding links between random graphs and other algebraic objects, in particular lattices, these being at the opposite extremes as ordered structures.

For the readers' convenience and for clarity of exposition we have taken a pedestrian approach to
some of the sections, recalling readily available theory  before derivation of a result.  This style is intended to make the material more easily digestible than if it were written down simply as
theorem, proof, theorem, proof, $\ldots$; it also helps to place a result in some context, given the episodic nature of the work.  Interesting comments on proof in mathematics can be found in~\cite{thurston1}.

The problems we have attacked represent our independent research program that is set in a
wider context of pushing the random graph concept into new
and uncharted territory, and in some chapters with a concentration on ideas rather than
results.  

There is a sense in which this work uses as its template the article
``Oligomorphic Groups and Homogeneous Graphs''~\cite{cameron} by
P.J. Cameron,
\index{Cameron, P. J.}%
in which there is a cross-fertilisation of random graph theory
with other parts of mathematics.  It has
guided the philosophy of the approach of this research program.  Furthermore, it has been over 20 years since the publication of the exposition ``Oligomorphic Permutation Groups'', which represented another key thesis in the subject, and our work is in a sense an extension of this book.  It has not been our aim to extend every topic that arose in the aforementioned paper and book, but only certain aspects.  Our pedagogic style is in keeping with its parent~\cite{cameron} having been a summer study school lecture course that included new results; this approach has been taken with a view to making the work more accessible to beginning researchers and those new to the field.  Hence approximately one third of the main text comprises background material that serves as an introduction to the theories that are used to produce the new results, giving the treatise the maximum breadth of readership, from beginning research students to experts.  Few of the new results have been published elsewhere.

We anticipate that the project will help towards an increased understanding of lively and rich areas of mathematics.  It is part of the effort to find relationships between graph theory and other fields; its multidisciplinarity should appeal to a wide range of tastes.  We
address specific problems whose solution were feasible, but also provide more
speculative longer term goals both to indicate the wide range of possible directions and because we it is not possible to know a priori the progress that will be made along a particular
path.  There is current interest in using graph-theoretic
tools in physics and we would in particular hope that symmetry
properties of random graphs find some uses there. 

Another way to view this work is as an application of the probabilistic
method.  Many of our results exemplify different aspects of (infinite) permutation group theory.  It was not possible to cover every conceivable topic; we have made little or no mention of cycle index theory, linear, circular and treelike-objects that arise as operands of oligomorphic groups, nor dynamics on random graphs.  We refer to~\cite{kechris} and~\cite{durrett} for work on the last of these.  Furthermore there are interpretations of random graph growth which may broadly be viewed as dynamic~\cite[p.~155]{alonspencer}
\index{Alon, N.}%
\index{Spencer, J.}%
 and which would be interesting to develop in future studies.

If there is one theme underlying the work as a whole, it maybe the concept of \emph{exceptionality}.  From the uniqueness of $\mathfrak{R}$ viewed as a (Bernoulli probabilistic, or Borel) measure~\cite{petrovv}
\index{Bernoulli measure}%
\index{Borel measure}%
 on a space of countable graphs, and its properties that are not shared by coloured random graphs on more than two edge-colours, to the isomorphism of Cameron and Glynn algebras in the infinite case alone, to the connection of $\mathfrak{R^{t}}$ with Cartan triality which we hope will be realised in the future, to the extension of the apparently unique Engeler--Ryll-Nardzewski--Svenonius Theorem~\cite{engeler}~\cite{ryll}~\cite{svenonius}.
\index{Engeler--Ryll-Nardzewski--Svenonius ! Theorem}%

Of course, not all possible generalizations are fruitful, witness for example the theory of higher-adjacency tournaments.  Nevertheless, the opportunities afforded the researcher for further work in all manner
of new directions are limitless.  The overarching vision is that the `random
graph' signify a mathematical method for solving problems rather than an object in discrete mathematics. 

Finally, in trying to cover so many aspects of a field, we have inevitably omitted a great many references that would be worthy of mention.  Regrettably, it is not possible to be totally comprehensive in surveying such a broad panorama as we have attempted to do.

\bigskip
\bigskip
\bigskip

This book is an account of joint work with Peter Cameron. 

\bigskip
\bigskip
\bigskip

``we require different colours, the [edges] belonging to any one substitution [i.e. generator] being of the same colour.''

\begin{flushright}
Arthur Cayley describing~\cite{cayley} Cayley graphs on two or more generators, as quoted in~\cite{meier}
\end{flushright}

\clearpage

\pagenumbering{arabic}

\chapter{Introduction: Basic Properties of $\mathfrak{R^{t}}$}
\label{fstchap}
\bigskip

Mathematics studies material things not as they are but as abstracted
from, though always existing in, matter $\ldots$.  As for logic, they
do not count it among the sciences, but rather as an instrument to
science.  Indeed it has been said: One cannot properly study or teach
except by means of the art of logic; for it is an instrument, and an
instrument of something is not a part thereof.

\begin{flushright}
Moses Maimonedes, \textit{Treatise on Logic c.1152}
\end{flushright}

\medskip

There are also many subjects of speculation, which, though not
preparing the way for metaphysics, help to train the reasoning power,
enabling it to understand the nature of a proof, and to test truth by
characteristics essential to it $\ldots$.  Consequently he who wishes
to attain human perfection, must therefore first study Logic, next the
various branches of Mathematics in their proper order, then Physics,
and lastly Metaphysics.
\begin{flushright}
Moses Maimonedes, \textit{The Guide for the Perplexed c.1190}
\end{flushright}

\medskip

Randomness is the true foundation of mathematics.
\begin{flushright}
Gregory Chaitin.
\end{flushright}

\bigskip

In the first section we restrict ourselves to the three-coloured random graph $\mathfrak{R^{t}}$,
\index{graph ! triality}%
however all the properties derived or mentioned are equally valid for the $m$-coloured random graphs denoted $\mathfrak{R}_{m,\omega}$,
\index{graph ! random ! $m$-coloured}%
where $m > 3$.  The core background is relegated to Appendix~\ref{TheRandomGraph}, and a more detailed account can found in~\cite{cameron}.

\bigskip

The \emph{random graph} $\mathfrak{R}$~\label{mathfrak{R}}
\index{graph ! random}%
is the unique countably infinite graph whose defining relation satisfies the ($*$)-condition:

($*$)  If $U$ and $V$ are finite
disjoint sets of vertices of $\mathfrak{R}$, then there exists in
$\mathfrak{R}$ a vertex $z$ joined to every vertex in $U$ and to no vertex in $V$.

We will return to the proof of this assertion later.  

The two adjacency types here are edge and non-edge, but this can readily be generalized to any number.  For example, if we choose three symmetric binary relations with any pair of vertices satisfying exactly one, they can more conveniently be thought of as colours, red, blue and green ($\mathfrak{r}, \mathfrak{b}, \mathfrak{g}$) thereby giving the \emph{triality graph} $\mathfrak{R^{t}}$.
\index{graph ! triality}%
The modified \emph{I-property} (so-called
because it is a form of injectivity)
\index{I-property}%
defining $\mathfrak{R^{t}}$~\label{mathfrak{Rt}} is:

($*_{t}$)  If $U$, $V$ and $W$ are finite
disjoint sets of vertices of $\mathfrak{R^{t}}$, then there exist in
$\mathfrak{R^{t}}$, a vertex $z$,
joined to every vertex in $U$ with a red edge, to every vertex in
$V$ with a blue edge, and joined to every vertex in
$W$ with a green edge.

The subscript $t$ and superscript $\mathfrak{t}$ stand for three or triality.  The property ($*_{t}$) is equivalent to the satisfaction of an infinite number of first-order sentences,
\index{first-order sentence}%
 so is a property of the first-order theory of $\mathfrak{R^{t}}$.
\index{first-order property}%

To demonstrate that such a graph  exists, we simply modify Rado's construction of $\mathfrak{R}$.
\index{Rado, R.}%
For $x, y \in \mathbb{N} \cup \{ 0 \}$, assuming $x < y$ (so that we can get an asymmetric joining rule for an undirected graph), express $y$ as a base $3$ expansion.  Join $x$ to $y$ respectively with a red, blue or green edge according to whether $y$ has $2, 1$ or $0$ in position $x$.  Take $3$ finite disjoint sets $U, V, W \subseteq \mathbb{N}$, and assume by adding new elements to $U$ and $V$ if necessary that $\mbox{max}(U) > \mbox{max}(V) > \mbox{max}(W)$.  Then the uniquely defined natural number $\sum_{u \in U} 2 . 3^u + \sum_{v \in V} 3^v$ can be chosen to represent $z$.

A  countable first-order structure
\index{first-order structure}%
 $M$ is \emph{$\aleph_0$-categorical}
\index{aleph@$\aleph_0$-categorical}%
or equivalently \emph{countably categorical}
\index{structure ! countably categorical}%
if any countable structure $N$ over the same language which satisfies the same first-order sentences is isomorphic to $M$.  
  
We can prove that $\mathfrak{R^{t}}$ is \emph{$\aleph_0$-categorical}
by using the \emph{back-and-forth method}.
\index{back-and-forth method}%

\begin{theorem}
Any two triality graphs with a countable infinity of vertices having property ($*_{t}$) are isomorphic.
\end{theorem}

\begin{proof} 
Let 
\[ \Gamma_1 = (a_1, a_2, \ldots)  \text{ and } \Gamma_2 = (b_1, b_2, \ldots) \]
be two such graphs with vertices enumerated.  The ($*_{t}$) condition
applies to both $\Gamma_1$ and $\Gamma_2$.  If $T$ is any graph
satisfying condition ($*_{t}$) we can show that the
following one-point extension property
\index{one-point extension property}%
 holds:

$(\dagger)$  If $A \subset B$ are finite $3$-coloured graphs
$S$, with $|B|=|A|+1$, then every embedding of $A$ into graph $T$ can be
extended to an embedding of $B$ into $T$.

Let $x$ be the vertex in $B$ and not in $A$, and $f$ the embedding of
$A$ into $T$.  Let $U, V, W$ be the sets of
vertices in $A$ joined to $x$ by red, blue and green edges
respectively.  Applying the definition with  $U_1 = f(U), V_1 = f(V),
W_1 = f(W)$, we find a vertex $z \in T$ satisfying the
conditions of the definition of $\mathfrak{R^{t}}$.  Then extend $f$
to $B$ by setting $f(x) = z$.  This shows that $f$ can be extended
so that its domain is all of $B$, proving $(\dagger)$.

Let $A_{i}$ be a finite set of vertices of $\Gamma_i$ $(i=1,2)$, and let
$$
\xymatrix{  \Theta:A_1 \ar[r]^{\cong} & A_2}$$.

Begin the back-and-forth argument with the empty isomorphism,
$\Theta_0$.  At the odd-numbered stage $n$, let $z_1$ be the first vertex after the enumeration of $A_{n-1} \subset \Gamma_1$ such that
$z_1 \notin \dom(\Theta_{n-1})$.  Then by $(\dagger)$, there exists
$\Theta_n$, which is an extension of $\Theta_{n-1}$, with domain
$A_{n-1} \cup {z_1}$.  

At even-numbered stage $n$, let $z_2$ be the first
vertex in the enumeration of $\Gamma_2$ such that
$z_2 \notin \dom(\Theta_{n-1}^{-1})$.  Then by applying $(\dagger)$ in
$\Gamma_1$, there exists $\Theta_n^{-1}$, with $z_{2} \in \dom
(\Theta_n^{-1})$.

Finally, the required isomorphism that has every vertex in $\Gamma_1$
in its domain and every vertex in $\Gamma_2$ in its range is the union
of the partial isomorphisms:

\[ \Theta = \bigcup_{n \ge 0} \Theta_{n} \]
where 
$$
\xymatrix{ \Gamma_1 \ar[r]^{\cong} & \Gamma_2 }.$$
\end{proof} 

An induced subgraph
\index{graph ! induced subgraph}%
 is obtained by omitting some of the vertices of a
graph and naturally the adjacencies on those vertices.  Embeddings of
graphs that arise for us are always as induced substructures.

The back-and-forth method
\index{back-and-forth method}%
 shows that any two countable triality graphs
 \index{graph ! triality}%
  are isomorphic.  Apart from uniqueness, it can be used to prove
homogeneity and universality, as per the next two results.
\index{graph ! homogeneous}%

\begin{theorem}
$\mathfrak{R^{t}}$ is homogeneous, that is any isomorphism between
finite subgraphs of $\mathfrak{R^{t}}$ can be extended to an
automorphism of $\mathfrak{R^{t}}$.
\end{theorem}

\begin{proof} 
Starting with any isomorphism between finite substructures of $\Gamma_1$, instead of the
empty isomorphism, and setting $\Gamma_1=\Gamma_2$ gives, via the
back-and-forth method, an automorphism of $\mathfrak{R^{t}}$.
\end{proof}

For structures such as $\mathfrak{R^{t}}$, with only finitely many $n$-element substructures up to isomorphism, homogeneity is a stronger concept than countable categoricity, (though this is no longer necessarily true if there are infinitely many such substructures).
\index{categoricity}%
  The homogeneity of $\mathfrak{R^{t}}$ immediately implies that the group $\Aut(\mathfrak{R^{t}})$~\label{Aut} is transitive on both edges of a particular colour and on finite subgraphs of
any given isomorphism type.

\begin{theorem}
$\mathfrak{R^{t}}$ is universal, that is any finite or countably
infinite $3$-edge-coloured complete graph can be embedded as an induced subgraph of
$\mathfrak{R^{t}}$.
\end{theorem}
\index{graph ! universal}%
\begin{proof} 
Let $\Gamma$ be a finite or countable graph, with vertices denoted by
${v_1, v_2, \ldots}$  We can inductively embed the
$n$-vertex subgraph $\Gamma_n = \{v_1, \ldots, v_n\}$
into $\mathfrak{R^{t}}$, using $(\dagger)$ to extend the
embedding of  $\Gamma_n$ to $\Gamma_{n+1}$.
Repeating this process we can embed the whole of $\Gamma$.
\end{proof}

\smallskip

Denote the random graph on $m$ colours and a countable infinity of
vertices by $\mathfrak{R}_{m,\omega}$;~\label{mathfrak{R_mn}} when $m = 3$ then
$\mathfrak{R}_{m,\omega}$ is simply $\mathfrak{R^{t}}$.  
\index{graph ! random ! $m$-coloured}%

Note that  $\mathfrak{R}_{2,\omega}$ is simply that graph $\mathfrak{R}$ that we began with, and the proofs of uniqueness, homogeneity and universality and the explicit construction are almost exactly the same as those we gave for $\mathfrak{R^{t}}$.  Going colourblind in two of the three colours of $\mathfrak{R^{t}}$ gives a graph that is isomorphic to $\mathfrak{R}$ and whose complement in $\mathfrak{R^{t}}$ is also isomorphic to $\mathfrak{R}$.

The
\emph{pigeonhole property}
\index{pigeonhole property}%
for any relational structure states that for
every partition of the point set of the structure into two
nonempty parts, the substructure induced on one of the parts is
isomorphic to the original structure.  For graphs, the point set
would be its vertex set.

\begin{theorem}
\label{ppforthree}
A countable $3$ edge-coloured graph in which all three colours occur
has the pigeonhole property if and only if it satisfies ($*_{t}$).
\end{theorem}

\begin{proof} 
We use the same method as that used in the corresponding proof for the
random graph $\mathfrak{R}$ to be found in \cite{cameron}.

Let $\Gamma_1$, $\Gamma_2$ be a partition
of the vertex set $\Gamma$ and suppose neither satisfies ($*_{t}$).  So for
$i=1,2$, there are disjoint finite subsets $U_i, V_i, W_i$ of $\Gamma_i$ for which
there does not exist a vertex $z_i \in \Gamma_i$ satisfying ($*_{t}$).
Then $U=U_1 \cup U_2,  V=V_1 \cup V_2, W=W_1 \cup W_2$ also fail ($*_{t}$)
in the whole graph $\Gamma$.

Conversely suppose $\Gamma$ has the pigeonhole property and partition
$\Gamma$ as $\Gamma_1 \cup \Gamma_2$ where $\Gamma_1$ consists of all
vertices lying on no green edge.  Then $\Gamma_1$ has red and blue
edges only and every vertex of $\Gamma_2$ lies on a green edge within $\Gamma_2$.  By assumption $\Gamma \ncong \Gamma_1$.  So $\Gamma \cong
\Gamma_2$, and thus every vertex lies on a green edge.  Similarly for
the other two colours.

Suppose, for a contradiction, that $\Gamma$ does not satisfy
($*_{t}$), and let $U, V, W$ be finite disjoint sets of vertices for
which no $z$ witnesses ($*_{t}$), and choose $|S|:=|U| + |V| + |W|$
minimal subject to this.  If $|S| = 1$ and $S = \{ x \}$ with $|U| =
1, V = \emptyset = W$, then no vertex can be joined to $x$ by a red
edge.  Similarly for the other two colours.  So $|S| \geq 2$.  Without
loss of generality $U \neq \emptyset$.

Choose $u \in U$.  Let $A$ consist of $u$ and those
vertices in $\Gamma \backslash (V \cup W)$ that are joined to it by blue and
green edges.  Let $B$ consist of the
remaining vertices.  Then $u$ is
red-isolated in $A$, so this subgraph cannot be isomorphic to $\Gamma$.  By the pigeonhole
property the induced subgraph on $B$ is isomorphic to $\Gamma$.
Therefore by the inductive hypothesis, that is by the minimality of $|S|$, $B$ has a vertex $z$ joined by a
red edge to every vertex in $U \backslash \{ u \}$, by a blue edge
to every vertex in $V$ and by a green edge to every
vertex in $W$.  But by construction of $B$, $z$ is also joined by a
red edge to $u$.  So $(U, V, W)$ was not
after all a counterexample to ($*_{t}$).
\end{proof}

We mention two relevant studies of partition properties of random graphs.
\index{graph ! random}%
  The first is an investigation of partition properties of edge-coloured random graphs by Pouzet and Sauer~\cite{pouzetsau}.
\index{Pouzet, M.}%
\index{Sauer, N.}%
 The second is concerned with colouring subgraphs of the Rado graph; here Sauer shows~\cite{sauer} that given a universal binary countable homogeneous structure $\mathcal{M}$ and $n \in \omega$, there is a partition of the induced $n$-element substructures of $\mathcal{M}$ into finitely many classes so that for any partition $C_0, \ldots, C_{m-1}$ of such a class $Q$ into finitely many parts there is a number $k \in m$ and a copy $\mathcal{M}^{*}$ of $\mathcal{M}$ in $\mathcal{M}$ such that all of the induced $n$-element substructures of $\mathcal{M}^{*}$ which are in $Q$ are also in $C_k$.  We refer to Sauer's paper for a precise definition of universal in this context.  Furthermore, the partition of the induced $n$-element substructures of $\mathcal{M}$ is given explicitly.

\medskip

A modification of the non-constructive existence proof of Erd\H{o}s
\index{Erd\H{o}s, P.}%
and R\'enyi
\index{R\'enyi, A.}%
for $\mathfrak{R}$ gives $\mathfrak{R^{t}}$ as follows.  Decide independently with
probability $\frac{1}{3}$, whether or not to join a pair of vertices with an edge of colour $\mathfrak{r}, \mathfrak{b}$ or $\mathfrak{g}$.  There are countably many different
choices for the triple $(U, V, W)$ in a graph $\Gamma$.  The probability
$p(z_1, \ldots, z_i)$ that
$z_1, \ldots, z_i, \in \Gamma \setminus (U, V, W)$ are not correctly
joined as $i \to \infty$ is 
\begin{displaymath}
 \lim_{i \to \infty}  \left(1 - \frac{1}{3^{s}}\right)^{i}  = 0,
\end{displaymath}
where $s=|S|=|U| \cup |V| \cup |W|$.  This leads us to the fact that ($*_{t}$) holds with probability one, because the countable union of null sets is null.  Since a set with probability 1 is certainly non-empty, this demonstrates non-constructively that $\mathfrak{R^{t}}$ exists.

We emphasize that $\aleph_0$-categoricity
\index{aleph@$\aleph_0$-categorical}%
is a weaker concept than
homogeneity over a finite relational language for the latter implies
the former, but not the other way round, as evidenced by the following three examples.  Firstly, an
infinite-dimensional vector space over a finite field cannot be made homogeneous over a finite relational language by adding relation symbols~\cite{macth}.  Secondly, two countable structures with the same \emph{age} (which is the class of all finite substructures  -- see Appendix~\ref{TheoryofRelationalStructures}) are isomorphic if they are homogeneous but not if they are $\aleph_0$-categorical.  In fact, Droste and Macpherson~\cite{drostemac}
\index{Droste, M.}%
\index{Macpherson, H. D.}%
constructed $2^{\aleph_0}$ non-isomorphic countable $\aleph_0$-categorical graphs, each of whose age being the class of all finite graphs, and further show that for every $k\in \mathbb{N}$ there are continuously many countable $\aleph_0$-categorical universal graphs which are $k$-homogeneous but not $(k+1)$-homogeneous (see Appendix~\ref{PermutationGroups}).  Thirdly, Philip Hall's countable universal homogeneous locally finite group
\index{Hall, P.}%
\index{group ! locally finite}%
\index{group ! Philip Hall's}%
(see Chapter~\ref{chapFD}), which however is not $\aleph_0$-categorical because it is not of finite exponent, and so has infinitely many $1$-types~\cite[p.56]{kayem}; this group is the Fra\"{\i}ss\'e limit
\index{Fra\"{\i}ss\'e limit}%
(see Appendix~\ref{TheoryofRelationalStructures}) of a suitable class of finite groups.

\head{Triality Graphs with Coloured Vertices}
\index{graph ! triality ! with coloured vertices}%

For the most part in this work we will assume that we have only one type of vertex, but to end this chapter let us assume that there are two types of vertices, cyan and yellow $(\mathfrak{c}, \mathfrak{y})$.  Our comments can be readily generalized to graphs with more than two vertex colours and more than three edge colours.  

So we can define a slightly more general graph denoted $\mathfrak{R^{t(v)}}$, which allows for $3$ possible edge-colours as well as two types of vertices, whose injection property states: 

($*_{t,v}$)  If $U$, $V$ and $W$ are finite disjoint sets of vertices of $\mathfrak{R^{t(v)}}$, then there exist in
$\mathfrak{R^{t(v)}}$, a cyan vertex $z^{\mathfrak{c}}$ and a yellow vertex $z^{\mathfrak{y}}$
both of which are joined to every vertex in $U$ with a red edge, to every vertex in
$V$ with a blue edge, and joined to every vertex in $W$ with a green edge. 

To demonstrate that such a graph exists, we simply modify Rado's construction of $\mathfrak{R}$, using the odd natural numbers for the cyan vertices, and the even natural numbers including zero for the yellow vertices.  The edge-formation rule is as above for $\mathfrak{R^{t}}$, except that the uniquely defined natural number $\sum_{u \in U} 2 . 3^u + \sum_{v \in V} 3^v$ can be chosen to represent $z^{\mathfrak{c}}$ if it is an odd number and $z^{\mathfrak{y}}$ if it is even.  If $m > |U \cup V \cup W|$, then we add $3^m$ to obtain a vertex of the opposite colour, also correctly joined.

\begin{theorem}
Any two triality graphs with a countable infinity of vertices of both cyan and yellow variety having property ($*_{t,v}$) are isomorphic.
\end{theorem}

The proof that $\mathfrak{R^{t(v)}}$ is \emph{$\aleph_0$-categorical} follows that for $\mathfrak{R^{t}}$, except that (a)  both a cyan vertex $z^{\mathfrak{c}} \in \mathfrak{R^{t(v)}}$ and a yellow vertex $z^{\mathfrak{y}} \in \mathfrak{R^{t(v)}}$ must be found to satisfy ($*_{t,v}$); (b)  $f$ is extended to $B$ by setting $f(x) = z^{\mathfrak{c}}$ if $x$ is cyan or $f(x) = z^{\mathfrak{y}}$ if $x$ is yellow.

Homogeneity and universality follow as they did for $\mathfrak{R^{t}}$.  An immediate corollary of the homogeneity of $\mathfrak{R^{t(v)}}$ is that the group $\Aut(\mathfrak{R^{t(v)}})$ has two orbits, namely cyan vertices and yellow vertices, and is transitive on edges with a particular colour and with endpoints having a given number of cyan vertices.

Theorem~\ref{ppforthree} is false when coloured vertices are allowed because we can have a
partition where one part consists entirely of cyan vertices and the other entirely of yellow vertices.

By working with a graph $\Gamma$ containing both cyan vertices ($z^{\mathfrak{c}}$) and yellow vertices ($z^{\mathfrak{y}}$) in the the Erd\H{o}s-R\'enyi construction, the probability $p(z^{\mathfrak{c}}_1, \ldots, z^{\mathfrak{c}}_i, z^{\mathfrak{y}}_1, \ldots, z^{\mathfrak{y}}_j)$ that $z^{\mathfrak{c}}_1, \ldots, z^{\mathfrak{c}}_i \in \Gamma \setminus (U, V, W)$ (respectively  $z^{\mathfrak{y}}_1, \ldots, z^{\mathfrak{y}}_j \in \Gamma \setminus (U, V, W)$) are not correctly joined according to ($*_{t,v}$) is 
\begin{displaymath}
 \lim_{i \to \infty} \left(1 - \frac{1}{3^{s}}\right)^{i} = 0,
\end{displaymath}
(respectively \begin{displaymath}
\lim_{j \to \infty} \left(1 - \frac{1}{3^{s}}\right)^{j} = 0),
\end{displaymath}
where $s=|S|=|U| \cup |V| \cup |W|$; this proves the existence of $\mathfrak{R^{t(v)}}$.

\chapter{Reducts of the Random Graphs}
\label{chap2}
\bigskip

The three Rs -- ``Reading, Writing and Arithmetic, the three subjects
which are taught in all elementary schools and are considered an
essential part of anyone's education.'' 
\begin{flushright}
D. M. Gulland and D. G. Hinde-Howell, \textit{The Penguin Dictionary of English Idioms -- Penguin Books (1988)}
\end{flushright}

\medskip

The Platonists as we have seen in Kepler's case, favoured in general a
trinitarian attitude in which the soul occupies an intermediary
position between mind and body.
\begin{flushright}
Wolfgang Pauli, \textit{Letter to Fierz, 9 March 1948,~\cite{enz}}
\end{flushright}

\medskip

Plus \c{c}a change, plus c'est la m$\hat{e}$me chose
\begin{flushright}
Alfonse Karr (1808-1890), \textit{Les Guepes,  Jan.1849}, (Novelist
and journalist; Oxford Dictionary of Quotations),
\end{flushright}

\medskip

We had to classify things in order to cope with the complexity of
nature.
\begin{flushright}
Robert Winston, \textit{BBC TV Program, 1999}
\end{flushright}

\bigskip

In this chapter, we introduce the reducts, concentrating our efforts
on the switching groups and deriving some of their properties.

\section{Switchings and Switching Permutations}

Let $\mathcal{G}_{m,n}$~\label{mathcal{G}_{m,n}} be the set of edge-coloured simple complete graphs
\index{graph ! simple}%
on a fixed set of $n$ vertices and $m$ edge-colours and $\Gamma \in \mathcal{G}_{m,n}$ be an arbitrary element of this set.

Let $c$ and $d$ be distinct colours, and $Y$ a subset of the vertex
set~$X$ of the graphs in $\mathcal{G}_{m,n}$. The operation
$\sigma_{c,d,Y}$~\label{sigmacdy} consists in interchanging the colours $c$ and $d$ whenever
they occur on an edge with just one vertex in $Y$, leaving all other colours
unchanged. Note that
\[\sigma_{c,d,Y} = \sigma_{d,c,Y} = \sigma_{c,d,X\setminus Y}\]
and
\[\sigma_{c,d,Y}^2=1.\]
The \emph{switching group} $S_{m,n}$~\label{switchinggroup}
\index{group ! switching}%
is
the group of permutations of $\mathcal{G}_{m,n}$ generated by the
\emph{switching operations},
\index{switching ! operation}%
which are defined to be the switchings on all vertex subsets of
$X$. A \emph{switching class}
\index{switching ! class}%
 is an orbit of
$S_{m,n}$ on $\mathcal{G}_{m,n}$.

If $g$ is any element of $\Sym(X)$, then $\sigma_{c,d,Y}^g=\sigma_{c,d,Y^g}$.
So $\Sym(X)$ normalizes $S_{m,n}$. Note
that $\Sym(X)$ also acts on the set $\mathcal{G}_{m,n}$, and hence so
does the semidirect product of these groups.

Let $\Gamma$ be an element of $\mathcal{G}_{m,n}$, that is, an $m$-coloured
graph on $X$. The \emph{group of switching automorphisms}~\label{groupofswaut}
\index{group ! switching ! automorphism}%
 $\SA(\Gamma)$ is defined by
\[\SA(\Gamma) = \{g\in\Sym(X) : \exists \sigma\in S_{m,n},\Gamma
g=\Gamma\sigma\}. \]
That is, it consists of all permutations of $X$ whose effect on $\Gamma$ can
be undone by a switching.

Another way to think about this group is as follows: let $\hat{G}$ be the
stabilizer of $\Gamma$ (then $\Gamma = \Gamma \sigma g^{-1}$, and $g$ undoes
the effect of $\sigma$ on $\Gamma$) in the semidirect product
$S_{m,n} \sd \Sym(X)$, and $\SA(\Gamma)$ the image
of $\hat{G}$ under the canonical projection of the semidirect product
onto $\Sym(X)$.

A switching class is an equivalence class of graphs under switching.  Note that multicoloured graphs in the same switching class have the same switching automorphism group.
\index{group ! automorphism}%
 For suppose that $g\in \SA(\Gamma)$, so that
$\Gamma g = \Gamma\sigma$,and let $\xi$~\label{switop} be any switching operation. Then
$(\Gamma\xi)g = (\Gamma\xi) \xi\sigma\xi^g$; since $\xi^g$ is again a
switching operation, $g$ is a switching automorphism of $\Gamma\xi$.

We end the section by emphasising the difference between switching
groups and groups of switching automorphisms.  A switching
group acts on the set of all graphs with specified sets
of vertices and colours.  The group of switching automorphisms acts
on the vertex set of \emph{one\ fixed} graph on a specific number of vertices and edges.  

\section{Two-Coloured Switching: The Finite Case}

A two-coloured complete graph can be identified with a graph on the same
vertex set; we identify the two colours $c$ and $d$ with adjacency and
nonadjacency. We write the switching operation $\sigma_{c,d,Y}$ simply as
$\sigma_Y$; it involves deleting edges with one end in $Y$ and inserting
edges between pairs of non-adjacent vertices of which one lies in $Y$.

The theory of such switching was worked out by J.~J.~Seidel
\index{Seidel, J. J.}%
in the 1960s in
connection with strongly regular graphs;
\index{graph ! strongly regular}%
it is sometimes referred to as
``Seidel switching''. 
\index{Seidel switching}%
\index{switching! Seidel}%
See, for example,~\cite{sei1}.  (A graph in which every vertex has the same degree is called a \emph{regular graph}.
\index{graph ! regular}%
A graph is \emph{strongly regular} if it is a regular graph such that the number of vertices mutually adjacent to a pair of vertices $v_1,v_2\in V(\Gamma)$ depends only on whether or not $\{v_1,v_2\}$ is an edge in the graph.)
We survey briefly parts of the theory relevant to our discussion.

Since $\sigma_Y\sigma_Z=\sigma_{Y\symd Z}$, and $\sigma_Y=1$ if and only if
$Y=\emptyset$ or $Y=X$, we see that the switching group
$S_{2,n}$ is elementary abelian
\index{group ! elementary abelian}%
of order $2^{n-1}$, and is isomorphic to $\mathcal{P}X/
\{\emptyset, X\}$, where the group operation on the power set $\mathcal{P}X$
is symmetric difference $\symd$. Moreover, this group acts semiregularly
\index{group ! permutation ! semiregular}%
 (by definition) on the set
$\mathcal{G}_{2,n}$ of graphs on the vertex set $X$. So the number of
switching classes is $2^{n(n-1)/2}/2^{n-1}=2^{(n-1)(n-2)/2}$.  This is
equal to the number of orbits by simple use of the orbit stabilizer theorem.
\index{Orbit-Stabilizer Theorem}%

Note that these are ``labelled'' switching classes;
\index{switching ! class}%
the enumeration of ``unlabelled'' switching classes (that is, the orbits of the semidirect
product $S_{2,n} \sd \Sym(X)$) was obtained by Mallows
\index{Mallows, C. L.}%
and Sloane~\cite{ms},
\index{Sloane, N. J. A.}%
a result to which we shall return shortly.

That switching is an equivalence relation is easily seen: (a)
Reflexivity: switch with respect to the whole of set $X$ and
everything outside $X$, that is $\emptyset$. (b) Symmetry: switch
with respect to $Y \subset X$; and then $Y$ again to recover the
original configuration. (c) Transitivity: switching with respect to $Y_1
\subset X$ then with respect to $Y_2 \subset X$ is the same as
switching with respect to $Y_1 \symd Y_2$, the symmetric difference of $Y_1$ and $Y_2$.

We need one further result about this case. Given a graph $\Gamma$ on the
vertex set $X$, we let $\mathcal{T}(\Gamma)$ denote the set of all
$3$-subsets of~$X$ which contain an odd number of edges of the graph
$\Gamma$. The following holds~\cite{camerona}:

(a) Graphs $\Gamma_1$ and $\Gamma_2$ on the vertex set $X$ satisfy
$\mathcal{T}(\Gamma_1) = \mathcal{T}(\Gamma_2)$ if and only if they belong to
the same switching class.

(b) A set $\mathcal{T}$ of $3$-subsets of $X$ is equal to
$\mathcal{T}(\Gamma)$ for some graph $\Gamma$ if and only if every $4$-subset
of $X$ contains an even number of members of $\mathcal{T}$.

A set $\mathcal{T}$ of $3$-sets satisfying the condition of (b) is called a
\emph{two-graph}.~\label{two-graph}
\index{graph ! two-graph}%
\index{two-graph}%
  The two-graph is trivial if either the set of $3$-subsets is empty or is the entire set of all $3$-subsets.  The terminology is due to G.~Higman
  \index{Higman, D. G.}%
(unpublished); see
Seidel~\cite{sei}.
\index{Seidel, J. J.}%
Thus there is a natural bijection between two-graphs and
switching classes, and the group of switching automorphisms of a graph is
equal to the group of automorphisms (in the obvious sense) of the two-graph
corresponding to its switching class. In particular, we see that the group of
switching automorphisms of $\Gamma$ cannot be more than $2$-transitive
\index{group ! permutation ! $2$-transitive}%
 unless $\Gamma$ is in the switching class of the complete or null graph.  The automorphism group
\index{group ! automorphism}%
  of the two-graph can act $2$-transitively on the graph vertex set, but not always.

Seidel has shown~\cite{sei} that for $n$-vertex graphs, the number
$t_n$~\label{tnswcl} of switching classes equals the number of two-graphs.  An
\emph{Euler graph}
\index{Euler graph}%
 is one in which every vertex has even degree.

\begin{theorem}[Mallows and Sloane]
\index{Mallows-Sloane Theorem}%
If $e_n$~\label{enswcl} is the number of Euler graphs on $n$-verices, then for all $n$, $e_n = t_n$.
\end{theorem}

\begin{proof}
Identify the edge set of a graph with its characteristic function, a binary vector of length $n \choose 2$.  Now let $V$ be the vector space of all such functions (with dimension $n \choose 2$) and consider the subspace $U$ spanned by the characteristic functions of the \emph{$n$ star graphs $S_i$}, where $S_i$ has all edges from $i$ to the vertices different from $i$.  Then $U$ consists of the characteristic functions of complete bipartite graphs, and has dimension $n - 1$.

Now observe:

\begin{itemize}

\item  Adding the star graph $S_i$ to a graph (modulo 2) is equivalent to switching at vertex $i$.  So the cosets of $U$ in $V$ are the switching classes of graphs.

\item A graph is Euler if and only if its characteristic vector is orthogonal to all the star graphs, and hence to all graphs in $U$.  So $U^{\perp}$ is the set of even graphs.

\end{itemize}

From this, the labelled version of the theorem follows: the numbers of switching classes and of even graphs are both $2^{{n \choose 2} - (n-1)} = 2^{(n-1) (n-2) / 2}$.

Now also there is a natural isomorphism between $(V / U)^*$ (the dual space of $V / U$) and $U^{\perp}$; this isomorphism respects the action of the symmetric group $\Sym(n)$.  By Brauer's Lemma,
\index{Brauer's Lemma}%
the numbers of orbits of a finite group on a finite vector space and its dual are equal.  Since the $\Sym(n)$-orbits on $V / U$ and on $U^{\perp}$ are isomorphism classes of switching classes and of Euler graphs respectively, the numbers of these isomorphism classes are equal; so the unlabelled version of the theorem is proved.
\end{proof}

\section{Two-Coloured Switching: The Infinite Case}

In this section we consider the case of two-coloured switching on infinite graphs.  The definition remains the same and we will consider just the random graph $\mathfrak{R}$.
\index{graph ! random}%

In mathematical logic new structures are often created from old ones
by taking the same underlying set and replacing one relation with
another.  For example we could replace a binary relation by an $n$-ary
relation for any $n > 2$.  In this way we can construct uncountably
many distinct definitions for essentially the same structure.  Often
two structures are considered to be the same if they are defined on
the same underlying set and have the same automorphism group.
\index{group ! automorphism}%
  A model-theoretic definition of reduct
\index{group ! permutation ! reduct}%
is

$\mathbf{Definition}$  Let $M$ be a countable structure for a finite
relational language.  A \emph{reduct}
\index{reduct}%
of $M$ is a permutation group $G$ for
which there is a structure $N$ and a relational language $L$ such
that:

1. $N$ has the same universe as $M$.

2.  For each relation $R \in L$, $R^N$ is definable without parameters
    in $M$.

3.  $G = \Aut(N)$.\\

Another definition which has a more group-theoretic focus is\\

$\mathbf{Definition}$  Let $M$ be a countable structure for a finite
relational language.  A \emph{reduct} of $M$ is a permutation group $G \leq
\Sym(M)$ such that:

1. $\Aut(M) \leq G$.

2.  $G$ is a closed subgroup of $\Sym(M)$, that is, $G = \overline{G}$, where $\overline{G}$ consists of all permutations fixing the $G$-orbits on $n$-tuples for all $n$.\\

So if $\bar{g} \in \Sym(M)$ and for all finite $A \subset M$ there is
a $g \in G$ such that $\bar{g}$ restricted to $A$ is the same map as
$g$ restricted to $A$, then $\bar{g} \in G$ also.

The two definitions are equivalent for $\aleph_0$-categorical
structures
\index{aleph@$\aleph_0$-categorical}%
as is proved for example in~\cite{bennett}.

The \emph{reducts} of the $m$-coloured random graph
\index{graph ! random ! $m$-coloured}%
 are those closed subgroups $G$ such that $\Aut(\mathfrak{R}_{m, \omega}) \leq G \leq \Sym(\mathfrak{R}_{m, \omega})$.  That there are
at least three types of proper reducts for $\mathfrak{R}$ was shown
in~\cite{camb}.  That these are the only ones for $\mathfrak{R}$ is
Thomas' Theorem~\cite{thomas},
\index{Thomas' Theorem}%
\index{Thomas, S.}%
whose proof gives the anatomy of the reducts of
$\mathfrak{R}$ as those in Figure~\ref{thomast}.
\begin{figure}[!ht]
$$
\xymatrix{
& \Sym(\mathfrak{R}) \ar@{-}[d]^{\infty} \\
& \BAut(\mathfrak{R}) \ar@{-}[dl]_{2} \ar@{-}[dr]^{\infty} \\
\SAut(\mathfrak{R}) \ar@{-}[dr]_{\infty} &&  \DAut(\mathfrak{R})\ar@{-}[dl]^{2} \\
& \Aut(\mathfrak{R})
}    
$$
\caption{Reducts of $\mathfrak{R}$}
\label{thomast}
\end{figure}
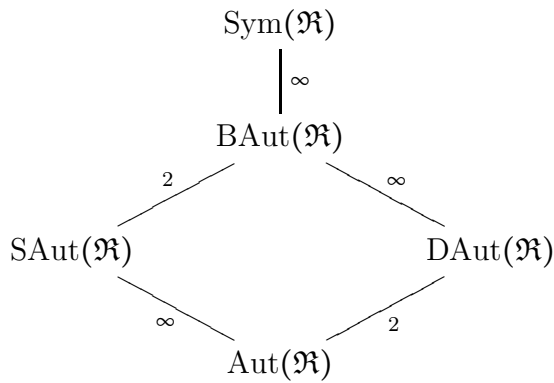  The proper reducts in this diagram are called \emph{duality} ($\DAut(\mathfrak{R})$),
\index{reduct ! duality}%
\emph{switching} ($\SAut(\mathfrak{R})$)
\index{reduct ! switching}%
and \emph{biggest} ($\BAut(\mathfrak{R})$) groups.~\label{dsbgps}
\index{reduct ! biggest}%

An \emph{anti-automorphism}
\index{graph ! anti-automorphism}%
of a graph $\Gamma$ is an isomorphism from $\Gamma$
to the complementary graph $\overline{\Gamma}$, while a \emph{switching
anti-automorphism}
\index{graph ! switching anti-automorphism}%
is a permutation $g$ such that $\Gamma g=\overline{\Gamma} \sigma$ for some switching $\sigma$.

We have already defined the group of switching automorphisms of a graph; that of $\mathfrak{R}$ is denoted $\SAut(\mathfrak{R})$.  Clearly $\mathfrak{R}$ is isomorphic to its complement $\overline{\mathfrak{R}}$.  Let $\DAut(\mathfrak{R})$ be the \emph{duality group}
\index{group ! duality}%
 of automorphisms and anti-automorphisms of $\mathfrak{R}$.  The \emph{biggest group} 
\index{group ! biggest}%
 whose elements are switching-automorphisms and switching anti-automorphisms, is defined by $\BAut(\mathfrak{R}) : = \langle \DAut(\mathfrak{R}), \SAut(\mathfrak{R}) \rangle$.
  
Thomas' Theorem (after Simon Thomas) can be stated as follows~\cite{thomas}:
\index{Thomas' Theorem}%
\index{Thomas, S.}%

\begin{theorem}[Thomas' Theorem]
There are just five reducts of the random graph $\mathfrak{R}$.
\index{graph ! random}%
 These are
\begin{itemize}
\item[(a)]  $\Aut(\mathfrak{R})$; 
\item[(b)]  $\DAut(\mathfrak{R})$; 
\item[(c)]  $\SAut(\mathfrak{R})$;
\item[(d)]  $\BAut(\mathfrak{R})$;
\item[(e)]  the symmetric group $\Sym(\mathfrak{R})$ on the vertex set of $\mathfrak{R}$.
\end{itemize}
\end{theorem}

 The duality group $\DAut(\mathfrak{R})$ acts $2$-transitively
\index{group ! permutation ! $2$-transitive}%
but not $3$-transitively on $\mathfrak{R}$ which follows from its definition, as does $\SAut(\mathfrak{R})$ because the parity of the number of edges in a 3-vertex set is invariant under switching.  The group $\BAut(\mathfrak{R})$ acts $3$-transitively but not $4$-transitively.  

We reproduce the outline of the proof that appears in~\cite[p.~344]{cam7a}.

\begin{proof}(Sketch)
Recall the Erd\H{o}s-R\'enyi
\index{Erd\H{o}s, P.}%
\index{R\'enyi, A.}%
$\mathfrak{R}$-construction.  Choose whether or not to join each pair of vertices in a countable set $S$ of them with an edge, each choice being made independently with probability $\frac{1}{2}$.  To see that property $(*)$ (given in Appendix~\ref{TheRandomGraph}) holds with probability $1$ in the random graph, firstly note that there are only countably many choices of disjoint finite sets $U, V \in S$.  Enumerate the vertices outside $U \cup V$ as $z_1, z_2, \ldots \in Z$.  If $s=|S|=|U \cup V|$ then the probability that no $z_i$ is correctly joined is
\begin{displaymath}
 \lim_{i \to \infty}  \left(1 - \frac{1}{2^{s}}\right)^{i}  = 0.
\end{displaymath}
The union of countably many null sets is a null set, so the construction follows.

So the set $Z$ of vertices satisfying $(*)$ produces a graph isomorphic to $\mathfrak{R}$.  From the homogeneity of $\mathfrak{R}$ it follows that the pointwise stabilizer in $\Aut(\mathfrak{R})$ of the set $U \cup V$ in $(*)$ acts homogeneously on $Z$, that is the induced permutation group is a dense subgroup of $\Aut(Z)$.

We call two finite graphs \emph{equivalent} if a copy of the first in $\mathfrak{R}$ can be mapped to a copy of the second by an element of a closed supergroup $G$ of $\Aut(\mathfrak{R})$.  The pointwise stabilizer $\Aut(\mathfrak{R})_{(X)}$ of a finite subset $X \subset \mathfrak{R}$ has $2^{|X|}$ orbits in $\mathfrak{R} \backslash X$, each orbit $Z$ corresponding to the subset of $X$ consisting of vertices adjacent to that orbit.  The orbits of $G_{(X)}$ are unions of these $Z$ orbits.  Of the following two cases, will discuss only (a) which applies to $\Aut(\mathfrak{R})$ and $\DAut(\mathfrak{R})$, and note that (b) which applies to $\SAut(\mathfrak{R})$ and $\BAut(\mathfrak{R})$ follows by similar arguments:\\
(a) $\exists\ \emptyset \neq X$ and an $\Aut(\mathfrak{R})_{(X)}$-orbit $Z$ which is fixed by $G_{(X)}$;\\
(b) no such $X$ and $Z$ exist.

Now assume that $G$ is a counterexample to the theorem.  Let $G$ be $n$-transitive but not $(n+1)$-transitive, and assume that $n$ is as small as possible.  It must be that $n > 1$, for if $G$ is not 2-transitive, it must preserve the $\Aut(\mathfrak{R})$-orbits on pairs, that is edges and non-edges.  Choose $X$ and $Z$ as in case (a).  By our earlier remark, the closure $A^{cl}$ of the permutation group induced by $\Aut(\mathfrak{R})_{(X)}$ on $Z$ is isomorphic to $\Aut(\mathfrak{R})$.  Now the closure $G^{cl}$~\label{G^{cl}} of $G_{(X)}$ is a closed supergroup of $A^{cl}$.  A combinatorial argument shows that $G^{cl}$ is not $n$-transitive.  By choice of $n$, $G^{cl}$ must be isomorphic to $\Aut(\mathfrak{R})$, $\DAut(\mathfrak{R})$, $\SAut(\mathfrak{R})$ or $\BAut(\mathfrak{R})$, and we can discount the latter two as we are in case (a).  By enlarging $X$ by two points of $Z$ if necessary, we may assume that $G^{cl} \cong \Aut(\mathfrak{R})$, as no anti-automorphism
\index{graph ! anti-automorphism}%
 can fix two points.  Equivalently, if $g \in G$ induces an isomorphism on $X$ then it induces an isomorphism on $Z$.  A combinatorial argument then shows that $G_{(X)} \le \Aut(\mathfrak{R})$, that is if $g \in G$ induces an isomorphism on $X$ then $g \in \Aut(\mathfrak{R})$.

We must now show that for every $g \in G$ there is a subgraph of $\mathfrak{R}$ isomorphic to $X$ on which $g$ induces an isomorphism or an anti-isomorphism.  This uses a Ramsey-type
\index{Ramsey theory}%
theorem of Deuber~\cite{deuber}
\index{Deuber, W.}%
and Ne\v{s}et\v{r}il
\index{Ne\v{s}et\v{r}il, J.}%
and R\"{o}dl~\cite{nesetril}:
\index{R\"{o}dl, V.}%
for every finite graph $U$ there is a finite graph $V$ with the property that for any 2-colouring of the edges of $V$ there is a monochromatic induced subgraph isomorphic to $U$.  By iterating we get a graph $W$ such that given 2-colourings of both its edges and non-edges, there is an induced subgraph isomorphic to $U$ in which both its edges and non-edges are monochromatic.  Apply the theorem with $U$ as the union of the $X$ constructed above, a graph not equivalent to a complete graph, and a graph not equivalent to a null graph.  Let $W$ be chosen thus, and assume that $W$ is embedded in $\mathfrak{R}$.  Give edge $e$ of $W$ colour $c_1$ if $eg$ is an edge and colour $c_2$ if $eg$ is a non-edge; and 2-colour the non-edges of $W$ similarly.  We find a copy of $U$ in $\mathfrak{R}$ for which both edges and non-edges are monochromatic, that is $Ug$ is complete or null or isomorphic to $U$ or its complement.  If we chose $U$ so that $Ug$ cannot be complete or null then $g$ induces an isomorphism or anti-isomorphism on $U$ (and hence on $X$), and the theorem is proved.
\end{proof}

\section{Three-Coloured Reducts}

That the random graph
\index{graph ! random}%
 is isomorphic to its complement is due to the
duality between edges and non-edges.  Similarly there is a triality between the types of edges in
$\mathfrak{R^{t}}$.  

An automorphism of $\mathfrak{R^{t}}$ is by definition an
adjacency-preserving permutation of its vertex set -- that is, $g \in
\Aut(\mathfrak{R^{t}})$ is a permutation if for each type (or colour) of edge
$\{x, y\}$, the edge $\{gx, gy\}$ is of the same type; so
$\Aut(\mathfrak{R^{t}})$ stabilizes the colour classes.  Let $\Aut^{*}(\mathfrak{R^{t}})$ be the group of
permutations of the vertex set which preserves the partition of the
edge set into the three colour classes.  The group $\Aut^{*}(\mathfrak{R^{t}})$ has
$\Aut(\mathfrak{R^{t}})$ as a normal subgroup, such that if
$T(\mathfrak{R^{t}})$ is the triality group
\index{group ! triality}%
(which is just $\Sym(\mathfrak{r} , \mathfrak{b} , \mathfrak{g})$) acting on the colours then $\Aut(\mathfrak{R^{t}}) \lhd  \Aut^{*}(\mathfrak{R^{t}})$ and 
\[  \Aut^{*}(\mathfrak{R^{t}}) / \Aut(\mathfrak{R^{t}}) \cong T(\mathfrak{R^{t}}) \cong \Sym(\mathfrak{r} , \mathfrak{b} , \mathfrak{g} ) \cong \Sym(3). \]
This is because $\Aut(\mathfrak{R^{t}})$ is the kernel of the action on colours; the quotient is $\Sym(\mathfrak{r} , \mathfrak{b} , \mathfrak{g})$ because permuting the colours gives a graph isomorphic to the original one.  The quotient is
isomorphic to $\Sym(3)$ because the group of permutations of the
colours has precisely the structure of $\Sym(3)$.  This triality
property of $\mathfrak{R^{t}}$ is
equivalent to the duality property of the random graph whereby
$\mathfrak{R}$ is isomorphic to its complement.

An automorphism of a graph gives orbits of different 
colours of edges.  In Chapter~\ref{chap7} we will more accurately identify this as an \emph{inner} automorphism, for clearly the cosets of $\Aut^{*}(\mathfrak{R^{t}})$ relative to $\Aut(\mathfrak{R^{t}})$ form outer automorphisms,
\index{automorphism ! outer}%
 a precise definition of which is given in this chapter.

Notice that $\Aut^{*}(\mathfrak{R^{t}})$ is $2$-transitive
\index{group ! permutation ! $2$-transitive}%
on vertices; for let $\alpha$ and $\beta$ be red-adjacent vertices of $\mathfrak{R^{t}}$.  By homogeneity these can be mapped to any other vertex pair joined by a red edge.  If $\gamma$ and $\delta$
are any vertex pair joined by a blue edge, then a triality mapping
$\mathfrak{r}$ to some $\mathfrak{b}$ in $\Aut^{*}(\mathfrak{R^{t}})$,
(either by stabilizing $\mathfrak{g}$ and transposing $\mathfrak{r}$
and $\mathfrak{b}$ or by a ($\mathfrak{r} \mathfrak{b}
\mathfrak{g}$)-cycle), takes $(\alpha, \beta)$ to some blue edge and an automorphism maps this blue edge to $(\gamma, \delta)$.

With Thomas' Theorem in mind we can define equivalent types of groups for $m$-coloured random graphs $\mathfrak{R}_{m,\omega}$ on a countable infinity of vertices, as follows:

(1) The \emph{Duality}
\index{group ! duality}%
 group reducts ($\DAut(\mathfrak{R}_{m,\omega})$).  For any subgroup $H
    \leq \Sym(m)$ there is a reduct of the $m$-coloured random graph
    corresponding to $H$, this being the group
    $\DAut^{H}(\mathfrak{R}_{m,\omega})$ of duality automorphisms which induce an element of $H$ on the set of colours.  For the usual random graph, the duality group preserves the parity of edges in every 4-subset of $\mathfrak{R}$.  We can construct such a duality-type
    reduct as follows.  For some positive integer $k$ take a $k$-ary relation $R$ (this exists at least when $k = m-1$) on the colour-set
    $\{c_1, \ldots, c_m \}$ such that $\Aut(R) = H$.  Define a
    $2k$-ary relation $\tilde{R}$ on vertices by 

$(x_1, y_1, x_2, y_2,
    \ldots, x_k, y_k) \in \tilde{R} \Leftrightarrow (c(x_1, y_1),
    \ldots, c(x_k, y_k)) \in R$,\\ where $c(x, y)$ denotes the colour of the edge between vertices $x$
and $y$.~\label{edgecol}

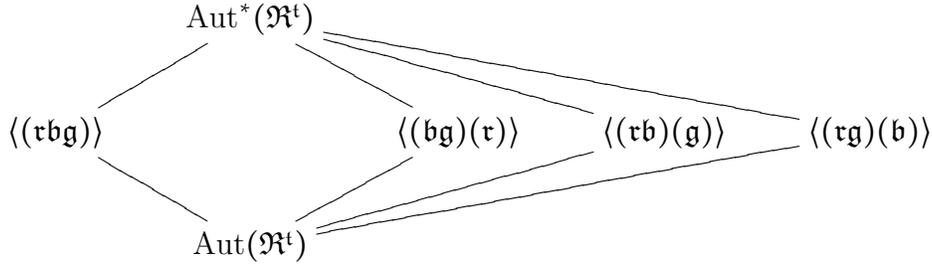
\begin{figure}[!ht]
$$
\xymatrix{
& \Aut^{*}(\mathfrak{R^{t}}) \ar@{-}[dl] \ar@{-}[dr] \ar@{-}[drr] \ar@{-}[drrr]\\
{\langle (\mathfrak{r} \mathfrak{b} \mathfrak{g}) \rangle} \ar@{-}[dr] && \langle (\mathfrak{b} \mathfrak{g})(\mathfrak{r}) \rangle \ar@{-}[dl] & \langle
(\mathfrak{r} \mathfrak{b})(\mathfrak{g}) \rangle \ar@{-}[dll] & \langle(\mathfrak{r} \mathfrak{g})(\mathfrak{b}) \rangle \ar@{-}[dlll] \\
& \Aut(\mathfrak{R^{t}})
}$$\caption{Duality-type reducts of $\mathfrak{R^{t}}$}  
\label{dualityRt}\end{figure}

Also let $P$ be a
    quaternary relation such that $(x_1, y_1, x_2, y_2) \in P
    \Leftrightarrow c(x_1, y_1) = c(x_2, y_2)$.  If a permutation $g$
    preserves $P$ then it induces a permutation $\hat{g}$ of the
    colours and if $g$ preserves $\tilde{R}$ then $\hat{g}$ preserves
    $R$ so $\hat{g} \in H$.  Therefore $g \in
    \Aut^{H}(\mathfrak{R}_{m,\omega})$,~\label{auth} which is the group of
    vertex-permutations that induces permutations $H$ of colours.
    Conversely $\Aut^{H}(\mathfrak{R}_{m,\omega})$ preserves
    $\tilde{R}$ and $P$.  Therefore $\Aut^{H}(\mathfrak{R}_{m,\omega})
    = \Aut(\tilde{R}, P)$ is a reduct.  So all such duality-type
    groups are genuine reducts.  Furthermore different groups $H$ give
    different reducts.
Specialising this to $\mathfrak{R^{t}}$, we can immediately give in Figure~\ref{dualityRt} five
of the duality-type reducts that play the role of the duality group of
$\mathfrak{R^{t}}$
where for example $\langle (\mathfrak{b} \mathfrak{g})(\mathfrak{r}) \rangle$
denotes the subgroup of $\Aut^{*}(\mathfrak{R^{t}})$ that is generated by permutations of the
vertices which induce the transposition $(\mathfrak{b} \mathfrak{g})$ on colours.  This lattice is isomorphic to the subgroup lattice of
$\Sym(3)$; so for example defining $K = \langle (\mathfrak{b} \mathfrak{g})(\mathfrak{r}) \rangle$ to be the group of isomorphisms and
    complementations (that is $\mathfrak{b}\mathfrak{g}$ colour transpositions) of $\mathfrak{R^{t}}$, gives a group acting on the vertices
    of $\mathfrak{R^{t}}$, and such that $ | K : \Aut(\mathfrak{R^{t}}) | = 2$.
 
(2) The \emph{Switching} group reducts
\index{group ! switching}%
($\SAut(\mathfrak{R}_{m,\omega})$).
    Switching a two-coloured graph with respect to a set of vertices
    transposes the colours of the edges between the vertex set and
    its complement, whilst leaving unaltered the coloured adjacencies within
    and outside the chosen vertex set.  For the usual random graph, the switching group preserves the parity of edges in every 3-subset of $\mathfrak{R}$.   If more than two colours are
    considered any subgroup $H \leq \Sym(m)$ yields a switching
    group as follows.  For any element $h \in H$, the $h$-switching
    with respect to a subset $X$ of the vertex set of  $\mathfrak{R^{t}}$ is the map which permutes the colours of edges  between $X$ and its complement according to $h$, whilst fixing the
    colours of edges within or outside $X$.  To such subgroups $H$, there corresponds a group of switching automorphisms which is
    the group of all permutations of the graph vertex set which map it
    to its image under some $h$-switching for some $h \in H$.  For example in the case of $\mathfrak{R^{t}}$, if the orbit partition
    of $H$ is of the form $\{\mathfrak{r}\} \{\mathfrak{b},
    \mathfrak{g}\}$, this gives a reduct for the closure of this group
    because, denoting a graph edge by $(x_i, y_i)$, the edges satisfy the relation $R$ for which:
$\{R = (x_1, y_1, x_2, y_2,
    \ldots, x_k, y_k):\mbox{parity\ of\ the\ $\mathfrak{b}$\ and\ $\mathfrak{g}$\ coloured\ edges\
    is\ invariant\ under\ switching}\}$; for more on this see Theorem~\ref{trreduct}.  As we will see below the
    partition $\{\mathfrak{r}, \mathfrak{b}, \mathfrak{g}\}$ yields a
    duality reduct but not a switching reduct.  
    
    The $\SAut(\mathfrak{R}_{3,\omega})$ reducts for intransitive
    subgroups $H$ of $\Sym(\mathfrak{r}, \mathfrak{b},
    \mathfrak{g})$ are not
    $2$-transitive,
   \index{group ! permutation ! $2$-transitive}%
    for if we stabilize the red colour say by switching
    only blue and green adjacencies, we would not be able to take a
    green or a blue edge to a red one by switching.  However if $H =
    \Sym(\mathfrak{r}, \mathfrak{b},
    \mathfrak{g})$ then the corresponding group of switching automorphisms
    turns out to be highly vertex transitive, as noted after the proof of Theorem~\ref{hightran} below.
 
(3) The \emph{Biggest} group reducts
\index{group ! biggest}%
($\BAut(\mathfrak{R}_{m,\omega})$) are those groups generated by both the switching reducts and by those dualities that preserve the partition of the colour set.  For the usual random graph, the biggest group preserves the parity of edges in every 5-subset of $\mathfrak{R}$.

The reducts of the random $m$-coloured graph have been classified in
the Ph.D. thesis of Bennett~\cite{bennett}.
\index{Bennett, J.H.}%
  For the rest of this section we report on his work, using our own notation for the reducts.  He finds two types of reducts:

(i) \emph{reducible} reducts
\index{reduct ! reducible}%
 ($G \leq \Sym(\mathfrak{R}_{m,\omega})$), for
which two colours are indistinguishable, (that is we are colourblind
in two colours), so that as far as the group of automorphisms is
concerned there are $m-1$ distinct colours.  The classification of
reducible reducts thus simply becomes the classification of reducts of $\mathfrak{R}_{m-1,\omega}$, which is done by induction.

(ii) \emph{irreducible}
\index{reduct ! irreducible}%
 reducts ($\Aut(\mathfrak{R}_{m,\omega}) \leq G < \Sym(\mathfrak{R}_{m,\omega})$) of $\mathfrak{R}_{m,\omega}$.
These are generated (as topological groups)
\index{topological ! group}%
 by groups of switching
automorphisms and duality automorphisms, where each group of switching
automorphims $\SAut(\mathfrak{R}_{m,\omega})$ corresponds to an
abelian
\index{group ! abelian}%
subgroup $A$ of $\Sym(m)$ and each
$\DAut(\mathfrak{R}_{m,\omega})$ corresponds to a subgroup $B$ of the normalizer of $A$
in $\Sym(m)$.

Note that the \emph{trivial} reducts,
\index{reduct ! trivial}%
$\Aut(\mathfrak{R}_{m,\omega})$ and
$\Sym(\mathfrak{R}_{m,\omega})$ have as respective structures the
original structure itself (which corresponds to taking the trivial
duality and switching groups) and a vertex set with no structure at
all (which is reducible since all colours are indistinguishable).

The motivation for Bennett's abelian characterization -- that if $\SAut(\mathfrak{R}_{m,\omega})$ is irreducible, then $\exists H \leq \Sym(m)$ corresponding to the switching group that must be abelian, is the following.  The commutator of two non-commuting permutations in $\Sym(m)$ must fail to preserve some colour $c$.  A switching of $\mathfrak{R}_{m,\omega}$ can be found that stabilizes colours on all edges except a specific one which originally had colour $c$.  Repeating this process gives a graph with no edges of colour $c$.  Another of Bennett's theorems then implies that $\SAut(\mathfrak{R}_{m,\omega})$ is reducible.  

Now let $H, D < \Sym(\mathfrak{R}_{m,\omega})$.  By the previous result, if a reduct $\Aut(\mathfrak{R}_{m,\omega}) \leq G < \Sym(\mathfrak{R}_{m,\omega})$ is irreducible then $H$ is abelian.  We can assume that $\SAut(\mathfrak{R}_{m,\omega})$ contains all 1-vertex switchings in $G$, so $H$ contains all corresponding colour permutations.  So if $\alpha \in D, \beta \in H$, but $\alpha \beta \alpha^{-1} \notin H$ then $\exists \delta \in \DAut(\mathfrak{R}_{m,\omega})$ and a single vertex switch $\sigma \in \SAut(\mathfrak{R}_{m,\omega})$ with colour permutations $\alpha$ and $\beta$ respectively, such that $\delta \sigma \delta^{-1} \notin \SAut(\mathfrak{R}_{m,\omega})$.  So $D$ normalizes $H$.

\medskip

Whereas Bennett's thesis~\cite{bennett}
\index{Bennett, J.H.}%
concentrates on classifying reducts of $\mathfrak{R}_{m,\omega}$, we have
explored slightly more widely, looking for example at groups that are
not reducts by virtue of their being highly transitive
\index{group ! permutation ! highly transitive}%
in their action
on the graph vertex set; this is detailed in later sections of
this chapter and in later chapters.  We also look at wider classes of graphs in the next section.

To illustrate the classification given
in~\cite{bennett} we find the reducts in the $3$-colour case, that is
the reducts of $\mathfrak{R^{t}}$.  There is an irreducible reduct for every combination of abelian group
$1 \leq A \leq \Sym(\mathfrak{r} , \mathfrak{b} , \mathfrak{g})$ and
their normalizers $1 \leq N(A) \leq \Sym(\mathfrak{r} , \mathfrak{b} ,
\mathfrak{g})$.  

The following table gives the $18$ irreducible reducts.  The interrelationships of all $31$ of the triality graph
\index{graph ! triality}%
 reducts that we give are illustrated in the Hasse diagram
\index{Hasse diagram}%
 given in Figure~\ref{redofr}.  An explicit proof showing that our list of reducts of $\mathfrak{R^{t}}$ is exhaustive remains an open question.

\vspace*{2cm}

\begin{figure}[!ht]
\[
\begin{array}{|c|c|c|}
\hline
A & B & Reduct \\
\hline
1  &  1  &  \Aut(\mathfrak{R^{t}}) \\
1  &  (\mathfrak{b} \mathfrak{g})  &  {\mathfrak{r}}D(\mathfrak{R^{t}})  \\
1  &  (\mathfrak{r} \mathfrak{g})  &  {\mathfrak{b}}D(\mathfrak{R^{t}})  \\
1  &  (\mathfrak{r} \mathfrak{b})  &  {\mathfrak{g}}D(\mathfrak{R^{t}})  \\
1  &  (\mathfrak{r} \mathfrak{b} \mathfrak{g})  &  \Alt(\mathfrak{r}, \mathfrak{b}, \mathfrak{g}) \\
1  &  \Sym(\mathfrak{r} , \mathfrak{b} , \mathfrak{g})  & T(\mathfrak{R^{t}})   \\
(\mathfrak{b} \mathfrak{g})  &  1  &  {\mathfrak{r}}S(\mathfrak{R^{t}})  \\
(\mathfrak{b} \mathfrak{g})  &  (\mathfrak{b} \mathfrak{g})  &  {\mathfrak{r}}B(\mathfrak{R^{t}})  \\
(\mathfrak{r} \mathfrak{g})  &  1  &  {\mathfrak{b}}S(\mathfrak{R^{t}})  \\
(\mathfrak{r} \mathfrak{g})  &  (\mathfrak{r} \mathfrak{g})  &  {\mathfrak{b}}B(\mathfrak{R^{t}})  \\
(\mathfrak{r} \mathfrak{b})  &  1  &  {\mathfrak{g}}S(\mathfrak{R^{t}})  \\
(\mathfrak{r} \mathfrak{b})  &  (\mathfrak{r} \mathfrak{b})  &  {\mathfrak{g}}B(\mathfrak{R^{t}})  \\
(\mathfrak{r} \mathfrak{b} \mathfrak{g})  &  1  &  {3}{\Aut(\mathfrak{R^{t}})}  \\
(\mathfrak{r} \mathfrak{b} \mathfrak{g})  &  (\mathfrak{b} \mathfrak{g})  &  {3\mathfrak{r}}D(\mathfrak{R^{t}})  \\
(\mathfrak{r} \mathfrak{b} \mathfrak{g})  &  (\mathfrak{r} \mathfrak{g})  &  {3\mathfrak{b}}D(\mathfrak{R^{t}})  \\
(\mathfrak{r} \mathfrak{b} \mathfrak{g})  &  (\mathfrak{r} \mathfrak{b})  &  {3\mathfrak{g}}D(\mathfrak{R^{t}}) \\
(\mathfrak{r} \mathfrak{b} \mathfrak{g})  &  (\mathfrak{r} \mathfrak{b} \mathfrak{g})  &  {3}{\Alt(\mathfrak{r}, \mathfrak{b}, \mathfrak{g})}  \\
(\mathfrak{r} \mathfrak{b} \mathfrak{g})  &  \Sym(\mathfrak{r} , \mathfrak{b} , \mathfrak{g})  &  {3}T(\mathfrak{R^{t}})  \\
\hline
\end{array}
\]
\caption{Irreducible reducts of $\mathfrak{R^{t}}$}
\end{figure} 

\vspace*{1cm}

\clearpage

\begin{figure}[ht]
\vbox to \vsize{
 \vss
 \hbox to \hsize{
  \hss
  \rotatebox{90}{
   \vbox{
    $$
    \hss
\xymatrix@R=30pt@C=5pt{
&&&&&&&&& \Sym(\mathfrak{R^{t}}) \ar@{-}[dlllll]^{\infty} \ar@{-}[dllllllll]_{\infty} \ar@{-}[dll]^{\infty} \\
& B(\mathfrak{R^{t}_{\mathfrak{r}}}) \ar@{-}[dl]_{2} \ar@{-}[dr]^{\infty} 
  &&& B(\mathfrak{R^{t}_{\mathfrak{b}}}) \ar@{-}[dl]_{2} \ar@{-}[dr]^{\infty}  
  &&& B(\mathfrak{R^{t}_{\mathfrak{g}}}) \ar@{-}[dl]_{2} \ar@{-}[dr]^{\infty} \\
S(\mathfrak{R^{t}_{\mathfrak{r}}}) \ar@{-}[dr]_{\infty} &&  D(\mathfrak{R^{t}_{\mathfrak{r}}}) \ar@{-}[dl]^{2} 
  & S(\mathfrak{R^{t}_{\mathfrak{b}}}) \ar@{-}[dr]_{\infty} && D(\mathfrak{R^{t}_{\mathfrak{b}}}) \ar@{-}[dl]^{2} \ 
  & S(\mathfrak{R^{t}_{\mathfrak{g}}}) \ar@{-}[dr]_{\infty} &&  D(\mathfrak{R^{t}_{\mathfrak{g}}}) \ar@{-}[dl]^{2} \\
& \Aut(\mathfrak{R^{t}_{\mathfrak{r}}})  \ar@{-}[d]^{\infty} &&& \Aut(\mathfrak{R^{t}_{\mathfrak{b}}})  \ar@{-}[d]^{\infty}
  &&& \Aut(\mathfrak{R^{t}_{\mathfrak{g}}}) \ar@{-}[d]_{\infty}
&& {3}T(\mathfrak{R^{t}}) \ar@{.}[d]^{2} \ar@{-}[uuu]_{\infty} \ar@{.}[dl]
\ar@{.}[dr] \ar@{.}[dllll] \ar@{.}[dlllllll]\\
& {\mathfrak{r}}B(\mathfrak{R^{t}}) \ar@{-}[dl]_{2} \ar@{-}[dr]^{\infty} 
& {3\mathfrak{r}}D(\mathfrak{R^{t}}) \ar@{.}[d]^{2}
  && {\mathfrak{b}}B(\mathfrak{R^{t}}) \ar@{-}[dl]_{2}
\ar@{-}[dr]^{\infty} 
& {3\mathfrak{b}}D(\mathfrak{R^{t}}) \ar@{.}[d]^{2}
  && {\mathfrak{g}}B(\mathfrak{R^{t}}) \ar@{-}[dl] \ar@{-}[dr]
& {3\mathfrak{g}}D(\mathfrak{R^{t}}) \ar@{.}[d]^{2}
& T(\mathfrak{R^{t}}) \ar@{-}[dl]^{\quad 3} \ar@{-}[dr]^{\quad 2}
\ar@{-}[dllll]^{3} \ar@{-}[dlllllll]_{3} 
& {3}{\Alt(\mathfrak{r}, \mathfrak{b}, \mathfrak{g})} \ar@{.}[d]^{2}\\
{\mathfrak{r}}S(\mathfrak{R^{t}}) \ar@{-}[drrrrrrrrr]_{\infty} &&  {\mathfrak{r}}D(\mathfrak{R^{t}})
\ar@{-}[drrrrrrr] & {\mathfrak{b}}S(\mathfrak{R^{t}}) \ar@{-}[drrrrrr] &&
{\mathfrak{b}}D(\mathfrak{R^{t}}) \ar@{-}[drrrr] & {\mathfrak{g}}S(\mathfrak{R^{t}})
\ar@{-}[drrr]^{\infty} &&  {\mathfrak{g}}D(\mathfrak{R^{t}})
\ar@{-}[dr]^{2} &
{3}{\Aut(\mathfrak{R^{t}})} \ar@{.}[d]^{2} \ar@{.}[ul] \ar@{.}[ur]
& {\Alt(\mathfrak{r}, \mathfrak{b}, \mathfrak{g})} \ar@{-}[dl]^{3} \\
&&&&&&&&& \Aut(\mathfrak{R^{t}}) }
    \hss
    $$
    \caption{Reducts of $\mathfrak{R^{t}}$}
    \label{redofr}  
   }
  }
 \hss}
\vss}
\end{figure}

\clearpage

We require an explanation of our notation for the groups in this diagram.  The above groups have the following actions:-

(i)  $\Sym(\mathfrak{R^{t}})$ is the group of all permutations of the vertices.

(ii)  $\Aut(\mathfrak{R^{t}_{\mathfrak{r}}}),~\label{autri}
S(\mathfrak{R^{t}_{\mathfrak{r}}}), D(\mathfrak{R^{t}_{\mathfrak{r}}}), B(\mathfrak{R^{t}_{\mathfrak{r}}})$.
These are respectively the automorphism group,
\index{group ! automorphism}%
 the group of switching automorphisms, the group of duality automorphisms, and the biggest group of automorphisms that is generated by the switching and duality
automorphisms.  These act on the vertex set of the triality graph
\index{graph ! triality}%
 which is colourblind in colour red ($\mathfrak{r}$), that is without
differentiating between red and either the green or blue.  So,
as only two colours are being considered, the existence of these
groups follows from Thomas' Theorem~\cite{thomas}.
\index{Thomas' Theorem}%
\index{Thomas, S.}%
They are what Bennett
\index{Bennett, J.H.}%
calls reducible reducts.

(iii)  ${\mathfrak{r}}D(\mathfrak{R^{t}}),
{\mathfrak{r}}S(\mathfrak{R^{t}}), {\mathfrak{r}}B(\mathfrak{R^{t}})$.~\label{stabi}
These groups act on the vertex set of the triality graph so as to stabilize the colour
red.  For example ${\mathfrak{r}}D(\mathfrak{R^{t}})$ permutes the set
of vertices with the effect of fixing red coloured edges
and either fixing or interchanging blue and green edges. Similarly
$\mathfrak{r}S(\mathfrak{R^{t}})$ switches blue and green, whilst preserving red.

(iv)  The existence of the triality group $T(\mathfrak{R^{t}})$
\index{group ! triality}%
and its index $2$ subgroup $\Alt(\mathfrak{r}, \mathfrak{b},
\mathfrak{g})$ was mentioned with at the beginning of this section.  More will be said in a later chapter.

(v)  ${3\mathfrak{r}}D(\mathfrak{R^{t}})$,
${3\mathfrak{b}}D(\mathfrak{R^{t}})$,
${3\mathfrak{g}}D(\mathfrak{R^{t}})$, ${3}T(\mathfrak{R^{t}})$,
${3}{\Alt(\mathfrak{r}, \mathfrak{b}, \mathfrak{g})}$, ${3}{\Aut(\mathfrak{R^{t}})}$.
These groups correspond to the colour permutations of the cyclic group $H$ of order $3$ in $\Sym(\mathfrak{r},
\mathfrak{b}, \mathfrak{g})$ preceding the respective actions, rather than simply being generated by
transpositions of pairs of colours.  We have prefixed these groups with $3$ to denote the
$3$-cycle action and have indicated them to be overgroups of the
other duality groups by dotted lines in the diagram.   The smallest one of these reducts is the $H$-automorphism group
\index{group ! automorphism}%
 ${3}{\Aut(\mathfrak{R^{t}})}$.  Three of the other five are
the $H$-duality groups for three non-trivial subgroups of $\Sym(\mathfrak{r}, \mathfrak{b},
\mathfrak{g})$.  The three
indices of the form $|{3\mathfrak{i}}D(\mathfrak{R^{t}}) :
{\mathfrak{i}}D(\mathfrak{R^{t}})|$ ($\mathfrak{i} = \mathfrak{r}, \mathfrak{b}, \mathfrak{g}$) equal $2$, because we can cycle in
one of two directions, either $\mathfrak{r} \to \mathfrak{g} \to
\mathfrak{b} \to \mathfrak{r}$ or $\mathfrak{r} \to \mathfrak{b} \to
\mathfrak{g} \to \mathfrak{r}$.  This inclusion also arises between the
other pairs of such $H$-duality reducts and their duality subgroups.

(vi)  $\Aut(\mathfrak{R^{t}})$ is the automorphism group
\index{group ! automorphism}%
 stabilizing the colours.

The green and blue equivalents of the red groups in (ii) and (iii) are
similarly defined and this gives us the remaining reducts in the diagram.  

\clearpage

\section{Random Graphs with Coloured Vertices}
\index{graph ! random}%
The aim of this section is to introduce simple random graphs on two different types of vertices, which we denote $\mathfrak{R}^{v}$ and to which we assign the colours cyan and yellow $(\mathfrak{c}, \mathfrak{y})$.  This section is taxonomic in nature
and indicates how the simple modification of adding an extra species of vertex leads to a proliferation of reducts.  This can be readily generalized to graphs with more than two vertex colours.  But first we briefly recall random bipartite graphs, which inevitably arise in any attempt to classify reducts of random graphs with coloured vertices.

The class of bipartite graphs is not homogeneous; for example, given two non-adjacent vertices, we do not know whether they are in the same bipartite block or not.  For a further example, it is not possible to amalgamate a path of length 2 and one of length 3 with the same endpoints in the class of bipartite graphs.  However the class of graphs with a prescribed bipartition does have a Fra\"{\i}ss\'e limit,
\index{Fra\"{\i}ss\'e limit}%
 this being a universal bipartite graph which is distance-transitive of diameter 3~\cite{cam8a}.  This \emph{random bipartite graph}
\index{graph ! random ! bipartite}%
 as defined in~\cite{cameron}, is the unique countable graph $\mathfrak{B}$~\label{randbip} such that, if we take two disjoint countable sets of vertices and join
pairs of vertices in different sets independently with probability
$\frac{1}{2}$, the resulting graph is almost surely isomorphic to
$\mathfrak{B}$.  Another way to obtain this graph, is to toss a fair coin to determine in which block a vertex goes, then to join it at random to vertices in the other block; with probability 1 we get infinitely many vertices in each bipartite block.  We recall some of the discussion of this object from~\cite{cameron}.

We assume that there are no edges between vertices within either
of the blocks.  A bipartite graph cannot be made
homogeneous unless it is complete bipartite or null, for two
non-adjacent vertices in different blocks cannot be mapped by an
automorphism to two vertices in the same block; so these pairs are not
considered equivalent.  A graph $\Gamma$ is said to be \emph{almost
homogeneous}
\index{graph ! almost homogeneous}%
if there are a finite number of relations $R_1, \ldots, R_k$~\label{Relations} such that

(i)  each $R_i$ is first-order definable (without parameters) in
$\Gamma$;

(ii)  the structure $(\Gamma, R_1, \ldots, R_k)$ is homogeneous.

For example the disjoint union of two copies of $\mathfrak{R}$ (having no edges between the copies) is not homogeneous, but becomes so when an extra binary relation is added to the language which is an equivalence relation whose classes are the two copies of $\mathfrak{R}$.  A second example is a bipartite graph with diameter at most $3$ which has a
bipartition defined by the set $\{(x, y) : (\exists z) (x
\sim z \sim y)\}$. 

The graph $\mathfrak{B}$ is universal for finite
bipartite graphs.  It is also almost homogeneous, and can be made homogeneous by
adding the above bipartition relation.  The \emph{I-property}
\index{I-property}%
for this graph is

($*_{B}$)  If $U$ and $V$ are finite
disjoint sets of vertices of $\mathfrak{B}$ in the same
bipartite block, then there exists a vertex $z$ in the opposite block of
$\mathfrak{B}$ joined to every vertex in $U$ and to no vertex in $V$.

Let $\Delta_1, \Delta_2$~\label{deltabip} denote the bipartite blocks of $\mathfrak{B}$, and $\Aut(\mathfrak{B})$ its automorphism group.
\index{group ! automorphism}%
 We will list some groups that are clearly reducts, though a proof that this list is exhaustive is left to future work.  A highly transitive group
\index{group ! permutation ! highly transitive}%
(whose closure is the symmetric group) can occur as the group induced on a definable subset in a nontrivial closed group.  We conjecture that this fact can be used in classifying the set of reducts of $\mathfrak{B}$, and though we will not give this classification we make the following comments:

1. $\mathfrak{B}$ is isomorphic to its bipartite complement $\mathfrak{B}^{*}$,~\label{bipcomp} obtained by interchanging edges and non-edges between $\Delta_1$ and $\Delta_2$.  So there is a duality group $\DAut(\mathfrak{B})$ containing $\Aut(\mathfrak{B})$ as a subgroup of index 2.

2.  Define \emph{bipartite switching}
\index{switching ! bipartite}%
with respect to $(X_1, X_2)$, where $X_1 \subseteq \Delta_1$ and $X_2 \subseteq \Delta_2$ as follows: exchange edges and non-edges between $X_1$ and $\Delta_2 \backslash X_2$, and between $X_2$ and $\Delta_1 \backslash X_1$, leaving the rest unaltered.  This leads to the switching group $\SAut(\mathfrak{B})$.

3.  There is a ``biggest group'' $\BAut(\mathfrak{B}) = \langle \DAut(\mathfrak{B}), \SAut(\mathfrak{B}) \rangle$, (except that in this case it is not the biggest nontrivial reduct).

4.  There is a group $\Sym(\Delta_1) \Wr \Sym(2)$, the automorphism group
\index{group ! automorphism}%
 of the complete bipartite graph on $(\Delta_1, \Delta_2)$.

5.  Finally, $\Sym(\Delta_1 \cup \Delta_2)$.

\begin{figure}[!ht]
$$
\xymatrix{
& \Sym(\Delta_1 \cup \Delta_2) \ar@{-}[d]^{\infty} \\
& \Sym(\Delta_1) \Wr \Sym(2) \ar@{-}[d]^{\infty} \\
& \BAut(\mathfrak{B}) \ar@{-}[dl]_{2} \ar@{-}[dr]^{\infty} \\
\SAut(\mathfrak{B}) \ar@{-}[dr]_{\infty} &&  \DAut(\mathfrak{B})\ar@{-}[dl]^{2} \\
& \Aut(\mathfrak{B})
}    
$$
\caption{Reducts of $\mathfrak{B}$}
\label{redofb}
\end{figure}
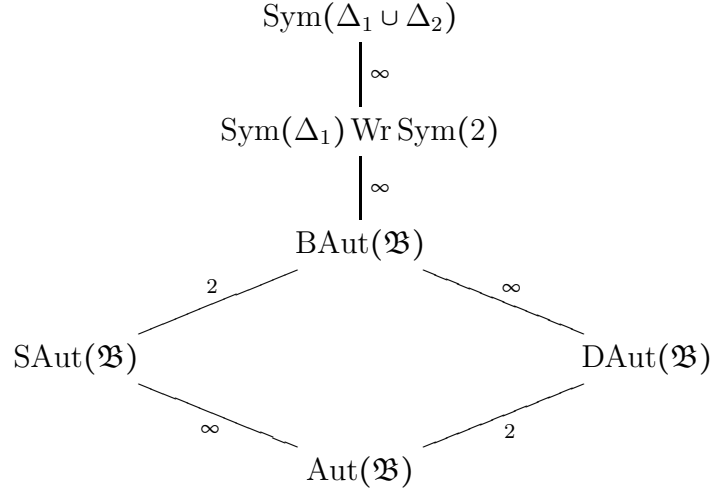  We can retain the names for the proper reducts that were used in the case of the random graph, that is \emph{duality} ($\DAut(\mathfrak{B})$),
\index{reduct ! duality}%
\emph{switching} ($\SAut(\mathfrak{B})$)
\index{reduct ! switching}%
and \emph{biggest} ($\BAut(\mathfrak{B})$) groups
\index{reduct ! biggest}%

Only the uppermost symmetric group reduct in Figure~\ref{redofb} is primitive
\index{group ! permutation ! primitive}%
 because all the others preserve the bipartite blocks.

Next we give Hasse diagrams
\index{Hasse diagram}%
of some reducts
of (a) the random bipartite graph
\index{graph ! random ! bipartite}%
 which we denote $\mathfrak{B}^{v}$ with one block comprising cyan vertices
$\mathfrak{B^c}$ and the other block yellow vertices
$\mathfrak{B^y}$, and (b)  the random graph with both cyan and
yellow vertices, $\mathfrak{R}^{v}$.  All of these groups are highly transitive.

\head{Remark}  A \emph{pseudoreduct}
\index{pseudoreduct}%
\index{reduct ! pseudoreduct}%
 is a permutation group that is not closed, but is closed on an enlarged vertex set.  Thomas
\index{Thomas, S.}%
 gave a definition for any countable structure~\cite{thomas1}.  The group induced on a bipartite block of $\mathfrak{B}$ by its stabilizer in $\Aut(\mathfrak{B})$ is highly transitive and so not a reduct, but it is a pseudoreduct.  The following groups are also pseudoreducts.

(a) Random Bipartite Graph with vertex colours, $\mathfrak{B}^{v}$.
\index{graph ! random ! bipartite}%

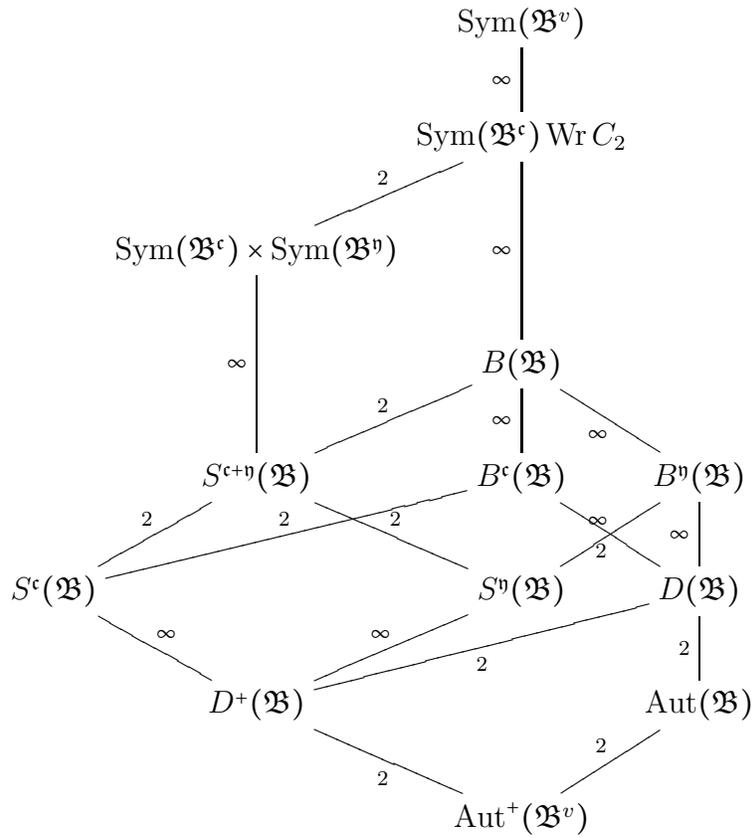
\begin{figure}[!ht]
$$
\xymatrix@C=0pt{
&& {} \ar@{}[d]\\
&& {} \ar@{}[d]\\
&& {} \ar@{}[d]\\
&& {\Sym(\mathfrak{B}^{v})} \ar@{-}[d]_{\infty}\\
&& \Sym(\mathfrak{B}^{\mathfrak{c}}) \Wr C_{2} 
\ar@{-}[dl]_{2} \ar@{-}[dd]_{\infty}\\
& \Sym(\mathfrak{B}^{\mathfrak{c}}) \times \Sym(\mathfrak{B}^{\mathfrak{y}})
\ar@{-}[dd]_{\infty}\\
&& B(\mathfrak{B}) \ar@{-}[dl]_{2} \ar@{-}[dr]_{\infty} \ar@{-}[d]_{\infty}\\
& S^{\mathfrak{c} + \mathfrak{y}}(\mathfrak{B}) \ar@{-}[dl]_{2} \ar@{-}[dr]^{2} & B^{\mathfrak{c}}(\mathfrak{B}) \ar@{-}[dll]_{2} \ar@{-}[dr]_{2} & B^{\mathfrak{y}}(\mathfrak{B}) \ar@{-}[dl]_{\infty} \ar@{-}[d]_{\infty}\\ 
S^{\mathfrak{c}}(\mathfrak{B}) \ar@{-}[dr]^{\infty} && S^{\mathfrak{y}}(\mathfrak{B})
\ar@{-}[dl]_{\infty} & D(\mathfrak{B}) \ar@{-}[d]_{2} \ar@{-}[dll]^{2}\\
& D^{+}(\mathfrak{B}) \ar@{-}[dr]_{2} && \Aut(\mathfrak{B}) \ar@{-}[dl]_{2}\\
&& \Aut^{+}(\mathfrak{B}^{v}) }
$$
\caption{Reducts of $\mathfrak{B}^{v}$ - (equivalently pseudoreducts
of $\mathfrak{R}^{v}$ (Part 1/4))}
\label{redbloops}
\end{figure}

The pseudoreduct diagram is Figure~\ref{redbloops}, with the following
meaning attached to the pseudoreducts:

$\Aut^{+}(\mathfrak{B}^{v})$ - group of automorphisms fixing the
two bipartite blocks.

$\Aut(\mathfrak{B})$ - group of automorphisms not necessarily fixing
the two bipartite blocks, but fixing the graph as a whole.

$D^{+}(\mathfrak{B})$ - group of duality permutations, interchanging
edges and non-edges between the two blocks.

$D(\mathfrak{B})$ - as $D^{+}(\mathfrak{B})$, except the permutations
either fix or interchange the two blocks.

$\Sym(\mathfrak{B^c}) \times \Sym(\mathfrak{B^y})$ -
symmetric group action independently on the two blocks.
\index{group ! action}%

$\Sym(\mathfrak{B^c}) \Wr C_{2} \cong
\Sym(\mathfrak{B^y}) \Wr C_{2}$ - wreath product
\index{wreath product}%
 group action permuting both blocks independently of each other and together with each other.

$\Sym(\mathfrak{B}^{v})$ - full symmetric group.

Now define a $\mathfrak{c}$-vertex switching as follows: for a subset $A
\subseteq B^{\mathfrak{c}}$, it interchanges edges and non-edges from $A$ to
$B^{\mathfrak{y}}$, and leaves all other adjacencies unchanged.  Similarly for a $\mathfrak{y}$-vertex switching.  Then we have that 

$S^{\mathfrak{c}}(\mathfrak{B})$ - group of $\mathfrak{c}$-vertex switching
automorphisms, which are permutations that fix blocks $B^{\mathfrak{c}}$ and
$B^{\mathfrak{y}}$ and map the graph into one that is $\mathfrak{c}$-vertex switching equivalent.

$S^{\mathfrak{y}}(\mathfrak{B})$ - as previous group with $\mathfrak{c}$ and $\mathfrak{y}$ interchanged.

$S^{\mathfrak{c} + \mathfrak{y}}(\mathfrak{B})$ - group generated by elements from both
$S^{\mathfrak{c}}(\mathfrak{B})$ and $S^{\mathfrak{y}}(\mathfrak{B})$.

$B^{\mathfrak{c}}(\mathfrak{B})$ - group generated by elements from both
$S^{\mathfrak{c}}(\mathfrak{B})$ and $D(\mathfrak{B})$.

$B^{\mathfrak{y}}(\mathfrak{B})$ - group generated by elements from both
$S^{\mathfrak{y}}(\mathfrak{B})$ and $D(\mathfrak{B})$.

$B(\mathfrak{B})$ - group generated by elements from both
$S^{\mathfrak{c} + \mathfrak{y}}(\mathfrak{B})$ and $D(\mathfrak{B})$.

Note that switching with respect to the whole of $B^{\mathfrak{c}}$ or the
whole of $B^{\mathfrak{y}}$ is complementation, which gives a pseudoreduct.  However,
we must demonstrate that more generally, the set of $\mathfrak{c}$-vertex switching
automorphisms forms a reduct.  This follows because the composition of
sequential switchings about two sets of vertices is equivalent to
switching about their symmetric difference; switching about the same
set twice gives the identity operation.  Closure in $\Sym(\mathfrak{B}^{\mathfrak{c}})
\times \Sym(\mathfrak{B}^{\mathfrak{y}})$ follows because
$S^{\mathfrak{c}}(\mathfrak{B})$ satisfies the following binary relation on
vertices, 

$rel := \{a_{1}, a_{2} \in B^{\mathfrak{c}}\ \text{and}\ b_{1}, b_{2} \in
B^{\mathfrak{y}}:$ 

$\text{\quad\quad the edges between these four vertices have the same parity} \}$.

Firstly any edge not in the switching set is stabilized.  Any pair of edges or
pair of non-edges switched together is parity-preserving, as is any switching
of an edge and a non-edge pair.  
Similarly for $S^{\mathfrak{y}}(\mathfrak{B})$ and $S^{\mathfrak{c} + \mathfrak{y}}(\mathfrak{B})$.
The wreath product action also leaves invariant a relation required
for closure of the resulting wreath product of groups, the $C_{2}$
involution corresponding to independent interchange of the two blocks.

\clearpage

\bigskip

\clearpage

Recall some theory~\cite{seideltsar} of two-graph
\index{graph ! two-graph}%
\index{two-graph}%
 representations for two-colour graphs.  Imagine an $n$-vertex graph $\Gamma = (V, E)$ with
a $(-1, 1)$-adjacency matrix $E$.
\index{adjacency matrix}%
  Switching about any vertex $v \in
V$ can be effected by addition$\pmod 2$ of the graph $K_{1,n-1}$ in
$v$.  So the switching class of $\Gamma$ is the equivalence class
$\Gamma$ modulo $\mathcal{B}_{c}$,~\label{mathcal{B}_{c}} where $\mathcal{B}_{c}$ is the set
of $n$-vertex complete bipartite graphs.  Equivalently it is the set
of graphs $(V, DED)$, for $D \in \mathcal{D}_{n}$ , where
$\mathcal{D}_{n}$ is the set of $n \times n$ diagonal matrices with
$\pm 1$ diagonal entries.  Equivalently still, it is the map $(V, E) \to (V, DED)$ where $D$ has $-1$ in the diagonal position $v$ and $1$ otherwise.

Switching classes of simple graphs can thus be explained in terms of bipartite graphs~\cite{seideltsar}.  The set of all graphs on a given vertex set is a $\mathbb{Z}_2$-vector space, where the sum of two graphs is obtained by taking the symmetric difference of their edge sets.  Complementation then corresponds to adding the complete graph, and switching to adding a complete bipartite graph. 

An example of the effect of the addition of a complete graph is the
following:

$$
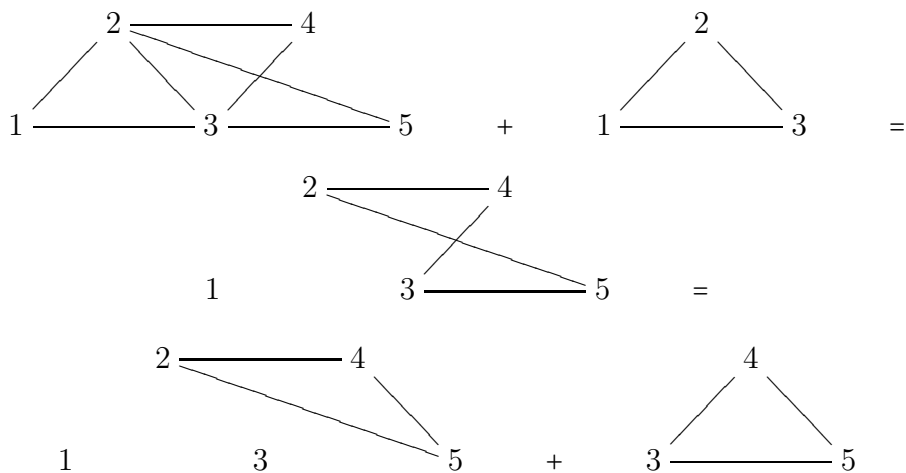
\xymatrix{
& {2} \ar@{-}[dl] \ar@{-}[dr] \ar@{-}[drrr] \ar@{-}[rr] && {4} \ar@{-}[dl] \ar@{}[rrr] &&&& {2} \ar@{-}[dl] \ar@{-}[dr]\\
 {1} \ar@{-}[rr] && {3} \ar@{-}[rr] && {5} & {+} & {1} \ar@{-}[rr] &&
{3} & {=}
}$$

$$\xymatrix{
& {2} \ar@{-}[drrr] \ar@{-}[rr] && {4} \ar@{-}[dl]\\
 {1} \ar@{}[rr] && {3} \ar@{-}[rr] && {5} & {=}
}$$

$$\xymatrix{
& {2} \ar@{-}[drrr] \ar@{-}[rr] && {4} \ar@{-}[dr] \ar@{}[rr] &&&& {4} \ar@{-}[dl] \ar@{-}[dr]\\
{1} \ar@{}[rr] && {3} && {5} & {+} & {3} \ar@{-}[rr] && {5}
}$$

\smallskip

 It follows from Thomas' Theorem
\index{Thomas' Theorem}%
\index{Thomas, S.}%
 that if $G$ is a closed supergroup of $\Aut(\mathfrak{R})$, then the set of all images of $\mathfrak{R}$ under $G$ is contained in a coset of a subspace $W(G)$ of this vector space.  (For example, $W(\BAut(\mathfrak{R}))$ consists of all complete bipartite graphs and all unions of at most two complete graphs).  These subspaces are invariant under the symmetric group.  

\clearpage 

(b) Random Graph with coloured vertices, $\mathfrak{R}^{v}$.~\label{rangphcolv}
\index{graph ! random ! with coloured vertices}%
This graph has a countable infinity of both $\mathfrak{c}$-vertices and
$\mathfrak{y}$-vertices.  The reducts can be treated as four separate cases, which also
makes the reduct diagram more transparent to follow, for these four
cases are mutually independent and lie between $\Sym(\mathfrak{R}^{v})$
and $\Aut(\mathfrak{R}^{v})$ in separate sub-diagrams, which we label
Parts $1-4$.

Firstly there are reducts that correspond to the graph
$\mathfrak{B}^{v}$, given in Figure~\ref{redbloops} above, all lying between
$\Sym(\mathfrak{R}^{v})$ and $\Aut(\mathfrak{R}^{v})$.
These permute the $\mathfrak{c}$ and $\mathfrak{y}$ vertices as though they lay in disjoint countable subgraphs.  

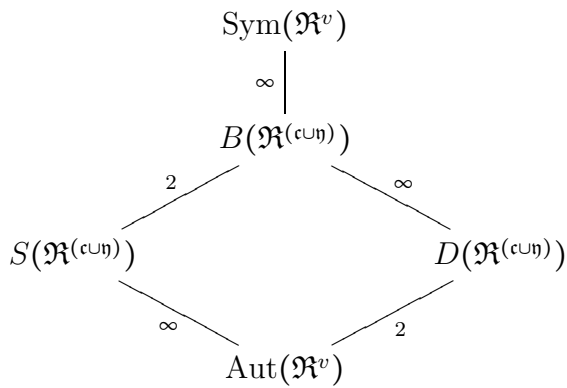
\begin{figure}[!ht]
$$
\xymatrix{
& \Sym(\mathfrak{R}^{v}) \ar@{-}[d]_{\infty} \\
& B(\mathfrak{R}^{(\mathfrak{c} \cup \mathfrak{y})}) \ar@{-}[dl]_{2} \ar@{-}[dr]^{\infty} \\
S(\mathfrak{R}^{(\mathfrak{c} \cup \mathfrak{y})}) \ar@{-}[dr]_{\infty} &&  D(\mathfrak{R}^{(\mathfrak{c} \cup
\mathfrak{y})}) \ar@{-}[dl]^{2} \\
& \Aut(\mathfrak{R}^{v})
}$$
\caption{Reducts of $\mathfrak{R}^{v}$ (Part 2/4)}
\label{part2/4}
\end{figure}

Secondly there are reducts corresponding to indistinguishability of $\mathfrak{c}$ and $\mathfrak{y}$ vertices, which we
represent by the superscript $(\mathfrak{c} \cup \mathfrak{y})$.  The resulting Hasse
diagram is given in Figure~\ref{part2/4}.

\begin{figure}[ht]
\vbox to \vsize{
 \vss
 \hbox to \hsize{
  \hss
  \rotatebox{90}{
   \vbox{
    $$
    \hss
    \xymatrix@R=20pt@C=2pt{
    &&& \Sym(\mathfrak{R}^{v}) \ar@{-}[d]^{\infty} \\
    &&& \Sym^{\mathfrak{c}} \times \Sym^{\mathfrak{y}} \ar@{-}[dll]_{\infty} \ar@{-}[d]^{\infty} \ar@{-}[drrr]^{\infty} \\
    &B^{\mathfrak{c}} \times \Sym^{\mathfrak{y}} \ar@{-}[d]^{\infty} \ar@{-}[dl]_{2}
    && B^{\mathfrak{c}} \times B^{\mathfrak{y}} \ar@{-}[d]^{\infty} \ar@{-}[dl]_{\infty} \ar@{-}[dr]^{\infty} \ar@{-}[drr]_{\infty} &&& \Sym^{\mathfrak{c}}
    \times B^{\mathfrak{y}} \ar@{-}[d]^{2} \ar@{-}[dr]_{\infty} \\
    S^{\mathfrak{c}} \times \Sym^{\mathfrak{y}} & D^{\mathfrak{c}} \times \Sym^{\mathfrak{y}} &
    B^{\mathfrak{c}} \times S^{\mathfrak{y}} & S^{\mathfrak{c}} \times B^{\mathfrak{y}} & B^{\mathfrak{c}} \times D^{\mathfrak{y}}
    & D^{\mathfrak{c}} \times B^{\mathfrak{y}}& \Sym^{\mathfrak{c}}
     \times S^{\mathfrak{y}} & \Sym^{\mathfrak{c}} \times D^{\mathfrak{y}} \\
    & \Aut^{\mathfrak{c}} \times \Sym^{\mathfrak{y}} \ar@{-}[u]_{2}
    \ar@{-}[ul]_{\infty} & S^{\mathfrak{c}} \times S^{\mathfrak{y}} \ar@{-}[u]_{2}
    \ar@{-}[ur]_{2} & D^{\mathfrak{c}} \times S^{\mathfrak{y}} \ar@{-}[ul]_{\infty}
    \ar@{-}[urr]_{2} & S^{\mathfrak{c}} \times D^{\mathfrak{y}} \ar@{-}[u]_{2}
    \ar@{-}[ul]_{\infty} &  D^{\mathfrak{c}} \times D^{\mathfrak{y}} \ar@{-}[ul]_{\infty} \ar@{-}[u]_{\infty} &
    \Sym^{\mathfrak{c}} \times \Aut^{\mathfrak{y}} \ar@{-}[u]_{\infty} \ar@{-}[ur]_{2} \\
    & \Aut^{\mathfrak{c}} \times B^{\mathfrak{y}} \ar@{-}[u]_{\infty} &&&&& B^{\mathfrak{c}} \times
    \Aut^{\mathfrak{y}} \ar@{-}[u]_{\infty} \\ 
    & \Aut^{\mathfrak{c}} \times S^{\mathfrak{y}} \ar@{-}[u]_{2} \ar@{-}[rruu]_{2}
    \ar@{-}[ruu]_{\infty} & \Aut^{\mathfrak{c}} \times
    D^{\mathfrak{y}} \ar@{-}[ul]_{\infty} \ar@{-}[uurr]_{\infty}
    \ar@{-}[uurrr]_{2} &&& S^{\mathfrak{c}}
    \times \Aut^{\mathfrak{y}} \ar@{-}[llluu]_{\infty} \ar@{-}[luu]^{2} \ar@{-}[ur]_{2} & D^{\mathfrak{c}}
    \times \Aut^{\mathfrak{y}} \ar@{-}[luu]^{2} \ar@{-}[llluu]_{\infty}
    \ar@{-}[u]_{\infty} \\ 
    &&& \Aut^{\mathfrak{c}} \times \Aut^{\mathfrak{y}} \ar@{-}[ul]_{2} \ar@{-}[ull]_{\infty}
    \ar@{-}[urr]_{\infty} \ar@{-}[urrr]_{2} \ar@{-}[d]^{\infty} \\
    &&& \Aut(\mathfrak{R}^{v}) 
    }
    \hss
    $$
    \caption{Reducts of $\mathfrak{R}^{v}$ (Part 3/4)}  
\label{part3/4}
   }
  }
 \hss}
\vss}
\end{figure}
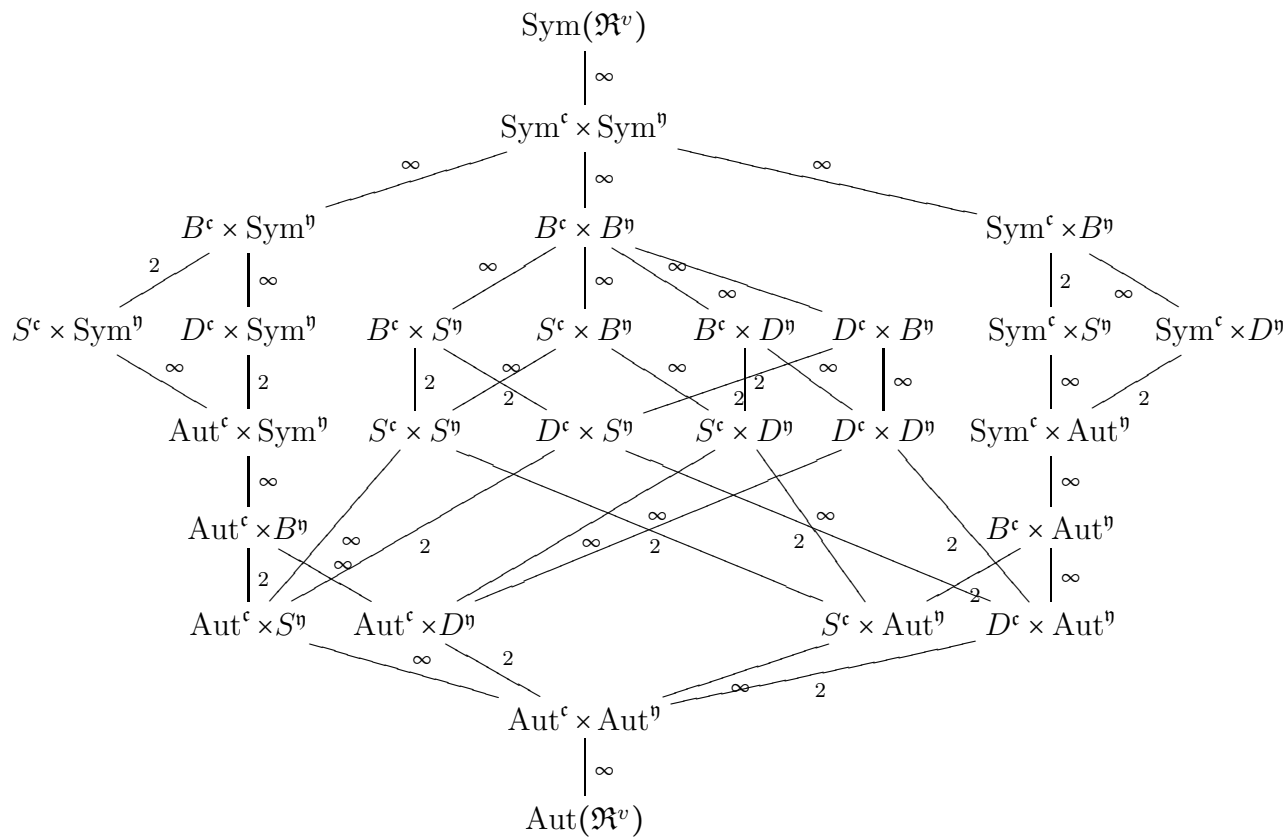

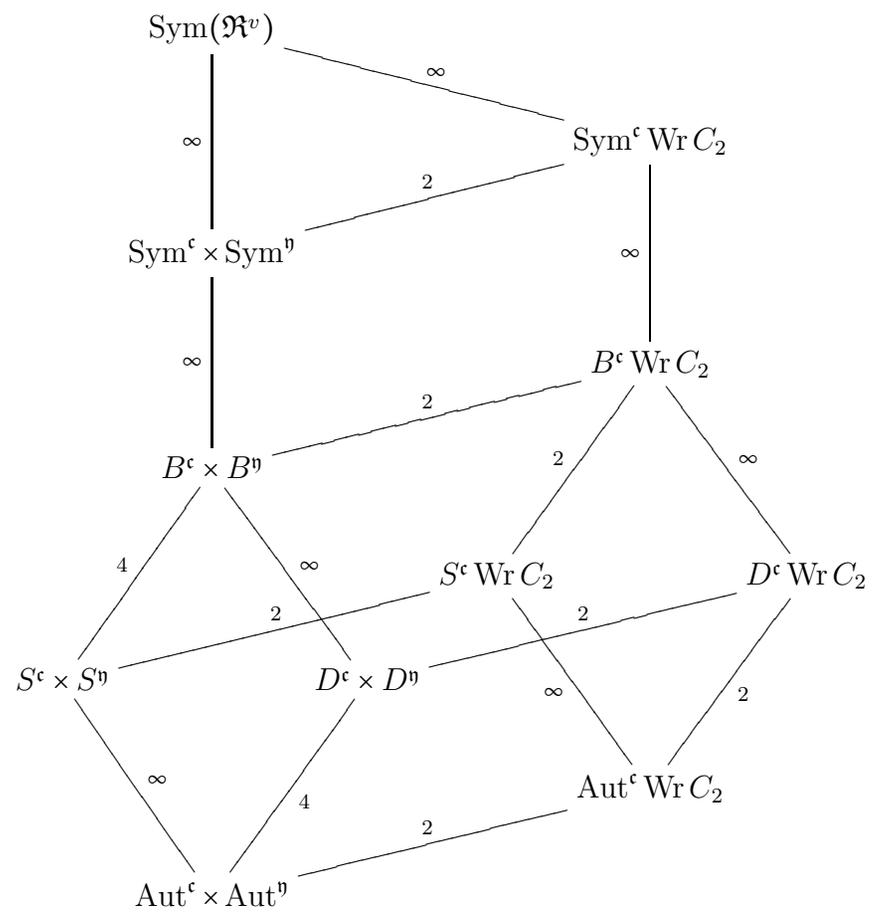
\begin{figure}[ht]
\vbox to \vsize{\vss
 \hbox to \hsize{\hss
  \rotatebox{90}{\vbox{
   $$
   \xymatrix@C=0pt{
   & {\Sym(\mathfrak{R}^{v})} \ar@{-}[dd]_{\infty}
   \ar@{-}[drrr]^{\infty}\\
   &&&& \Sym^{\mathfrak{c}} \Wr C_{2} \ar@{-}[dd]_{\infty}\\
   & \Sym^{\mathfrak{c}} \times \Sym^{\mathfrak{y}}
   \ar@{-}[dd]_{\infty} \ar@{-}[urrr]^{2}\\
   &&&& B^{\mathfrak{c}} \Wr C_{2} \ar@{-}[ddl]_{2} \ar@{-}[ddr]^{\infty}\\
   & B^{\mathfrak{c}} \times B^{\mathfrak{y}} \ar@{-}[rrru]^{2} \ar@{-}[ddl]_{4}
   \ar@{-}[ddr]^{\infty} \\
   &&& S^{\mathfrak{c}} \Wr C_{2} \ar@{-}[ddr]_{\infty} && D^{\mathfrak{c}} \Wr C_{2}
   \ar@{-}[ddl]^{2}\\
   S^{\mathfrak{c}} \times S^{\mathfrak{y}} \ar@{-}[ddr]^{\infty} \ar@{-}[rrru]^{2} &&
   D^{\mathfrak{c}} \times D^{\mathfrak{y}} \ar@{-}[rrru]^{2} \ar@{-}[ddl]^{4}\\ 
   &&&& \Aut^{\mathfrak{c}} \Wr C_{2} \ar@{-}[dlll]_{2}\\
   & \Aut^{\mathfrak{c}} \times \Aut^{\mathfrak{y}} }
   $$
   \caption{Reducts of $\mathfrak{R}^{v}$ (Part 4/4)}  
\label{part4/4}
  }
  }
 \hss}
\vss}
\end{figure}

Thirdly there are the reducts permuting the $\mathfrak{c}$ vertices amongst
themselves, and the $\mathfrak{y}$ vertices amongst themselves.
For all wreath products of any reduct $G$, the group  $G^{\mathfrak{c}} \Wr C_{2}$ is isomorphic to $G^{\mathfrak{y}} \Wr C_{2}$ as abstract
groups.  Also $\Sym(\mathfrak{B}^{\mathfrak{c}}) \times \Sym(\mathfrak{B}^{\mathfrak{y}})$
$\cong$ $\Sym(\mathfrak{R}^{\mathfrak{c}}) \times \Sym(\mathfrak{R}^{\mathfrak{y}})$.  We break
down the resulting Hasse diagram
\index{Hasse diagram}%
 into two independent parts, given on
the next two pages as Figures~\ref{part3/4} and~\ref{part4/4}, which
are joined at their common groups $\Sym(\mathfrak{R}^{v})$ and
$\Aut^{\mathfrak{c}} \times \Aut^{\mathfrak{y}}$.  

In Figure~\ref{part3/4} we have not included all possible subgroup inclusions for reasons of space, omitting for example that $\Aut^{\mathfrak{c}} \times \Aut^{\mathfrak{y}} < D^{\mathfrak{c}} \times D^{\mathfrak{y}}$ where $| D^{\mathfrak{c}} \times D^{\mathfrak{y}} : \Aut^{\mathfrak{c}} \times \Aut^{\mathfrak{y}}| = 4$.


\clearpage

\vspace{30pt}

\clearpage

\clearpage

Clearly we would not be able to comfortably fit the reduct diagram of the
triality graph
\index{graph ! triality}%
 with coloured vertices on one page.
In principle it is possible to construct a reduct diagram for graphs on any countable number of different vertex species.

\bigskip

Finally in this section, we mention a model of simple, finite random graphs in which the vertices rather than the edges are the focus.  In these \emph{random intersection graphs}~\cite{karonski}~\cite{jaworski},
\index{graph ! random ! intersection}%
 each vertex is independently assigned a random structure and the adjacency of two vertices is decided by comparing these structures.
 
More specifically, a graph $\Gamma$ is an \emph{intersection graph}
\index{graph ! intersection}%
 if each vertex $v \in V(\Gamma)$ can be assigned a set $S_v$  such that $\{v, w\} \in V(\Gamma)$ precisely when  $S_v \cap S_w \neq 0$.  Then $\Gamma$ is the intersection graph of the family of sets $\mathcal{S} = \{ S_v : v \in V(\Gamma) \}$.  If furthermore, the sets from $\mathcal{S}$ are generated in a random way then $\Gamma$ is a \emph{random intersection graph}.

\section{Random Graphs with Forbidden Substructures}
\index{graph ! random ! with coloured vertices}%

Let $K_k$~\label{K_k} be the complete graph on $k$ vertices.  Finite $K_k$-free graphs
\index{graph ! Kfree@$K_k$-free}%
amalgamate to satisfy Fra\"{\i}ss\'e's hypotheses, (see Appendix~\ref{TheoryofRelationalStructures}).  For $k \ge 3$, a \emph{Henson graph}, denoted $\mathfrak{H}_k$,
\index{graph ! Henson}%
\index{Henson, C. W.}%
is the unique countable homogeneous $K_k$-free graph which contains every finite $K_k$-free graph~\cite{hen1}.  The relevant I-property is
\index{I-property}%

($*_{k}$)  If $U$ and $V$ are finite disjoint sets of vertices, such that $U$ contains no $K_{k-1}$, then there is a vertex $z$ joined to every vertex in $U$ and to no vertex in $V$.
 
The notion of \emph{indivisibility}
\index{structure ! indivisible}%
 for relational structures originates in Fra\"{\i}ss\'e's book~\cite{frai}.
\index{Fra\"{\i}ss\'e, R.}%
 A structure $\mathcal{M}$ is \emph{Ramsey}
\index{structure ! Ramsey}%
 for $\mathcal{N}$, if for every partition $\mathcal{M}$ into two classes $C$ and $D$ there is an embedding of $\mathcal{N}$ into $C$ or an embedding of $\mathcal{N}$ into $D$; see Pouzet's~\cite{pouzeta2}.
\index{Pouzet, M.}%
 A graph $\Gamma$ is said to be \emph{indivisible}
\index{graph ! indivisible}%
if for every partition of the vertex set of $\Gamma$ into two classes $A$ and $B$ there is an isomorphic copy of $\Gamma$ either in $A$ or in $B$.  Folkman~\cite{folkman} studied the Ramsey-type
\index{Ramsey theory}%
 property that graphs have monochromatic complete subgraphs in every edge colouring.  El-Zahar and Sauer
\index{El-Zahar, M.}%
\index{Sauer, N.}%
 proved in~\cite{elz} that every homogeneous $K_k$-free graph, $\mathfrak{H}_k$ $(k \ge 3)$, is indivisible, whilst in~\cite{elzahar} they related the divisibility or indivisibility of the age of homogeneous relational structures to the way they amalgamate.  See~\cite[Chapter 6, \S 6]{frai} for more on this notion, and~\cite{covingtonmac} for further applications to countable homogeneous indivisible structures.
\index{structure ! indivisible}%

Thomas
\index{Thomas, S.}%
used the same technique as that in~\cite{thomas} to prove that the reducts of the homogeneous universal $K_k$-free graphs $\mathfrak{H}_k$ for $2 < k < \omega$ are $\Aut(\mathfrak{H}_k)$ and $\Sym(\mathfrak{H}_k)$.  Switching or complementing a subgraph of $\mathfrak{H}_k$ may not result in another $\mathfrak{H}_k$.  

The strongest known universality result about graphs with forbidden substructures is the following result of Cherlin, Shelah and Shi~\cite{cherlin}:
\index{Cherlin, G. L.}%
\index{Shelah, S.}%
\index{Shi, N.}%
\begin{theorem}
For every finite set of finite connected graphs, the class of all graphs for which there is no homomorphism into any member of the finite graph set, has a universal object.
\end{theorem}
An extension of this to relational structures
\index{relational structure}%
 is in~\cite{cherlin1}, while the main result of~\cite{hubickanes1} is that the universal object for the classes of the previous theorem can be obtained as a reduct of a so-called \emph{generic structure}, 
\index{generic structure}%
\index{structure ! generic}%
which we now define. 

\medskip

One way to construct the countable random graph
\index{graph ! random}%
 that is equivalent to the usual method is
as a countable \emph{generic structure}
\index{structure ! generic}%
where generic means construction by finite approximations.  Random graphs are defined in a space of graphs with a given property.  The potential difficulties with this process are firstly that of defining a measure
\index{measure}%
on this space and secondly that of determining whether or not a limit exists.  Even if a limiting structure exists it could be a surprising one, as we shall see presently.  That a random structure is not always equivalent to a generic structure is exemplified by triangle-free graphs.
\index{graph ! triangle-free}%
A triangle-free graph is generic if (i) it is countably infinite, (ii) a triangle of edges is not embeddable in it, (iii) it is homogeneous, and (iv) it is universal for finite triangle-free graphs.  Such graphs provide an example of a structure whose large size limit is different to that which we would expect, for a random triangle-free graph is bipartite,
\index{graph ! random ! bipartite}%
 meaning that if 
\[a_n := \frac{\sharp\ \text{$n$-vertex bipartite\
graphs}}{\sharp\ \text{triangle-free\ $n$-vertex\ graphs}}\]
then $a_n \to 1$ as $n \to \infty$~\cite{erdos}.  (A bipartite graph is necessarily triangle-free.  The generic triangle-free graph is Henson's graph $\mathfrak{H}_3$).
\index{graph ! Henson}%
  By contrast, the class of cycle-free graphs
\index{graph ! cycle-free}%
fails to be universal for cycle-lengths greater than 3.

Hage et al. look at the switching operation
\index{switching ! operation}%
 of a graph transformation, from the viewpoint that it is a global transformation of a graph, achieved by applying local transformations to the vertices in the form of group actions.
\index{group ! action}%
 Switching classes grow rapidly and in~\cite{hage}, a study is made of detecting Euler, triangle-free and bipartite graphs in these classes.

\bigskip

Another example of a structure-free graph is the Covington
graph
\index{graph ! Covington}%
\index{Covington, J.}%
 $\mathfrak{N}$~\cite{coving}~\label{cogr}.  There are no homogeneous \emph{$N$-free}
graphs, those which contain no induced path of length $3$, because in
such a graph any set of three vertices has the possibility of a
distinguished vertex.  The class $\mathbf{N}$~\label{clnfr} of $N$-free graphs
\index{graph ! N-free@$N$-free}%
does not have the amalgamation property,
\index{amalgamation property}%
so no homogeneous graph has age
$\mathbf{N}_{f}$, the class of finite $N$-free graphs.  However $\mathfrak{N}$ is almost homogeneous
\index{almost homogeneous}%
and by the judicious addition
of a ternary relation $r$ to the language of graph theory, the failure of
amalgamation can be repaired and the graph $\mathfrak{N}$ constructed
as the unique countable universal $N$-free graph.  So there is a countable homogeneous
structure with age $\mathbf{N}[r]_{f}$.  

The subject of homogenizability of structures has been studied by Kun and Ne\v{s}et\v{r}il
\index{Kun, G.}%
\index{Ne\v{s}et\v{r}il, J.}%
 in~\cite{kun}.

We can say a little more detail about $N$-free graphs.  The purpose of the ternary relation is to distinguish one of any triple of vertices.  For triples where the induced subgraph is one of the following two,
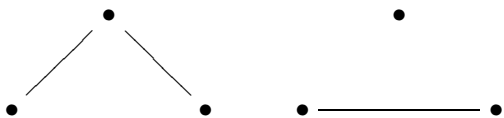
\begin{figure}[!h]$$\xymatrix{
& {\bullet} \ar@{-}[dr] \ar@{-}[dl] &&& {\bullet} \ar@{}[dr] \ar@{}[dl]\\
{\bullet} \ar@{}[rr] && {\bullet} & {\bullet} \ar@{-}[rr] && {\bullet}
}$$
\end{figure}

one vertex (the uppermost one in the diagrams) is already distinguished by the graph structure.  However consider the triangle subgraph.
\begin{figure}[!h]$$\xymatrix{
& {\bullet} \ar@{-}[dr] \ar@{-}[dl]\\
{\bullet} \ar@{-}[rr] && {\bullet}
}$$
\caption{Triangle subgraph}
\label{triangle}
\end{figure}
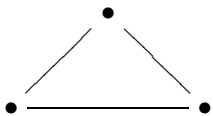 \\ 

\begin{figure}[!h]$$\xymatrix{
& {\bullet}  \ar@{-}[dl] \ar@{-}[rr] && {\bullet} &&& {\bullet}  \ar@{-}[dl] \ar@{-}[rr] && {\bullet}\\
{\bullet} \ar@{}[rr] &&&&& {\bullet}\\
& {\bullet} \ar@{}[rr] \ar@{-}[ul] \ar@{-}[uu]  \ar@{-}[rr] && {\bullet} &&&  {\bullet} \ar@{}[rr] \ar@{-}[ul] \ar@{-}[uu]  \ar@{-}[rr] && {\bullet} \ar@{-}[uu]
}$$
\caption{Non $N$-free supergraphs of the triangle subgraph}
\label{nonnfreetriangle}
\end{figure}
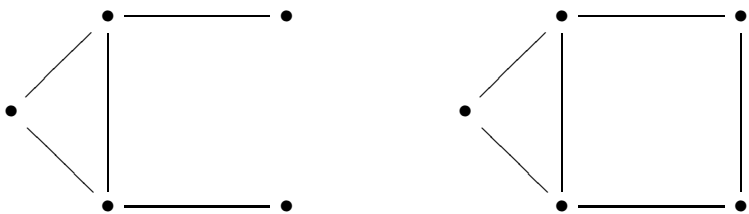

It is not possible to have two vertices, each with one neighbour in the triangle and these neighbours distinct, since neither of the the next two graphs (Figure~\ref{nonnfreetriangle}) are $N$-free.  
\index{graph ! N-free@$N$-free}%

So if there is a vertex with one neighbour in the triangle, this neighbour is distinguished.  If not, the ternary relation prescribes the only vertex which would be distinguished in this way in any $N$-free supergraph.  Similar remarks hold for the subgraph
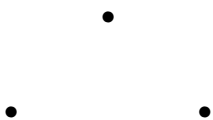
\begin{figure}[!h]$$\xymatrix{
& {\bullet} \ar@{}[dr] \ar@{}[dl]\\
{\bullet} \ar@{}[rr] && {\bullet}
}$$
\end{figure} 

Our comments show that in the universal $N$-free graph, the ternary relation is determined by the graph structure.  A list of areas where $N$-free graphs arise is given in~\cite{coving}.

In summary, although $N$-free graphs are not naturally homogeneous, they are homogenizable.  If a structure is homogenizable
\index{structure ! homogenizable}%
 then it is \emph{generic}; we briefly say a little about this and expand on a similar concept of generic automorphism.  A \emph{generic structure} 
\index{structure ! generic}%
is one that is residual in some natural complete metric space in the sense of the \emph{Baire category theorem}~\cite{baire}
\index{Baire category theorem}%
 (see Appendix~\ref{CategoryandMeasure}).  This says that we can put a natural metric structure on the class of all objects with a given \emph{age} (see Appendix~\ref{TheoryofRelationalStructures}) or smaller, so that the isomorphism class of such objects is residual, or in other words that almost all objects look like the one we are interested in.  The notion of \emph{generic structures}
\index{generic structure}%
\index{structure ! generic}%
includes homogeneous, $\aleph_{0}$-categorical as well as homogenizable structures, including $N$-free graphs.  

A \emph{generic automorphism} 
\index{automorphism ! generic}%
 is one whose conjugacy class is large in the sense of Baire category.  A homogeneous structure will have a generic automorphism, where a \emph{generic element}
\index{generic group element}%
of the automorphism group
\index{group ! automorphism}%
 is one whose conjugacy class is residual, that is it lies in a comeagre conjugacy class relative to the given topology.  Let $G$ be a permutation group acting on a countably infinite set $\Omega$.  The natural topology
\index{topology ! natural}%
  on permutation groups is the \emph{topology of pointwise convergence}
\index{topology ! of pointwise convergence}%
in which for a sequence $g_n \in G$ of permutations $G= \Sym(X)$ on $X$, $\lim_{n \to
  \infty} (g_n) = g$ if and only if $\forall x_i \in X,\ \exists n_0
\in \omega$ such that $\forall n > n_0,\ x_i g_n = x_i g$.  (We shall assume that the permutation groups have countable degree; for arbitrary degree we would have to consider limits of nets rather than sequences.)  Endowing a group $G$ with this topology turns it into a \emph{topological group},
\index{group ! topological}%
that is multiplication and inversion are continuous, so that if $g_n \to g$ and $h_n \to h$ then $g_n h_n \to gh$ and $g_n^{-1} \to g^{-1}$.  Further endowing $G$ with a metric turns it into a complete metric space, so that we can use the Baire category theoretic
\index{Baire category theorem}%
notion of `meagre' set (see Appendix~\ref{CategoryandMeasure}).  The topology is generated by basic open sets
\index{sets ! basic open}%
 of the form $[p] = \{g \in G : g\ \text{agrees with}\ p\ \text{on its domain} \}$ where $p$ is a $1$--$1$ map from a finite subset of $\Omega$ to $\Omega$.  Given any relational or first-order structure, its automorphism group is a closed subgroup of $\Sym(\Omega)$ and so is automatically a complete metric space, and therefore we can talk about generic automorphisms.  One way of capturing the notion of $g$ being a `typical' element is to require $g$ to lie in certain \emph{dense open sets}.  
\index{dense open set}%
An example of a dense open set for any $g$ is $\{G \backslash g\}$.  

Truss~\cite{truss3} 
\index{Truss, J. K.}%
 undertook a study of generic automorphisms of homogeneous structures, proving the above results and that $\Sym(\omega)$, $\Aut(\mathfrak{R}_{m,\omega})$ and $\Aut(\mathbb{Q}, <)$ the group of order automorphisms of the rationals (namely, $\{g \in \Sym(\mathbb{Q}) : p < q \Leftrightarrow p^g < q^g, \forall p, q \in \mathbb{Q}\}$), all have generic elements with specified cycle types for such elements.  More information on generic automorphisms, together with references can be found in Appendix~\ref{TopologyinPermutationGroups}.

A further sample to the literature on random graphs with forbidden substructures is~\cite{cam8a}~\cite{cherlin2}~\cite{cherlin3}~\cite{cherlin4}~\cite{furedi}~\cite{komjath}.

\section{Properties of $S_{3,3}$}

The main aim of this section is to discover the structure of the group
$S_{3,3}$.  This will help in finding the form of the more general
switching groups $S_{m,n}$, a task that we delegate to section~\ref{dersec} of
this chapter.  We also learn some of the properties of the action of
$S_{3,3}$ on complete graphs $\mathfrak{R}_{3,3} \in
\mathcal{G}_{3,3}$ with $3$ vertices (1, 2, 3) and $3$ colours ($i, j,
k)$.  

The possible switchings operations are given by

$$\xymatrix{
& {2^{2}} \ar[r] & {\sigma_{i,j,1}} & {\sigma_{i,j,2}} & {\sigma_{i,j,3}} \\
& {2^{2}} \ar[r] & {\sigma_{j,k,1}} &
{\sigma_{j,k,2}} & {\sigma_{j,k,3}} \\
& {2^{2}} \ar[r] & {\sigma_{i,k,1}} & {\sigma_{i,k,2}} & {\sigma_{i,k,3}} \\
&& {\Sym(3)} \ar[u] &{\Sym(3)} \ar[u] &{\Sym(3)} \ar[u]
}$$

Following the notation established previously in this
chapter, $\sigma_{i,j,1}$ means a switch of colours $i$
and $j$ whenever they occur on an edge about vertex $1$.  The $2^{2}$ denotes the Klein
$4$-group
\index{group ! Klein}%
 generated along each row, and $\Sym(3)$ denotes the group generated along each column.

We will derive the group generated by these $9$ involutions.  Firstly
note that five of the nine entries in this matrix are redundant as
generators.  Later, we will consider the $3$-set
embedded in an $n$-vertex graph; but for now, the three vertices of
the $3$-set are considered as comprising the whole graph.
Any two of $\sigma_{i,j,X}$, $\sigma_{j,k,X}$ and $\sigma_{i,k,X}$
give the third.  Also in both rows the product of two of the
switchings equals the third.  In other words, $\sigma_{i,j,1} \sigma_{i,k,1} \sigma_{i,j,1} =
\sigma_{j,k,1}$ and $\sigma_{i,j,1} \sigma_{i,j,2} =
\sigma_{i,j,3}$.  So we actually need only
consider four of the nine possible generators; without loss of
generality we choose $\sigma_{i,j,1}$, $\sigma_{j,k,1}$,
$\sigma_{i,j,2}$, $\sigma_{j,k,2}$.  

We use $\mathsf{GAP}$ 4
\index{gap@$\mathsf{GAP}$}%
to find the group
generated by these on the $3^{3} = 27$ graphs on ordered $3$-sets,
which are as follows.  One way to label the $27$ possible graphs is
using elements of the set $\{0, \ldots, 26\}$.  We can parametrize the
$3$ possible colours on each of the $3$ edges using ordered triples of
numbers from $\{0,1,2\}$.  So a coloured graph can be labelled using
the $27$ numbers $\{0, \ldots, 26\}$ written in base $3$ then adding
$1$, since $\mathsf{GAP}$ 4 assumes the set permuted is $\{1, \ldots, 27\}$.  By pairing-off switching-equivalent graphs, we can combine the
resulting transpositions to form the four requisite generators, where each permutation runs into a second line:-
\[ \texttt{P1} := \sigma_{i,j,1} =
(1\ 7)(2\ 4)(3\ 8)(5\ 6)(9\ 18)(10\ 25)(11\ 20)(12\ 22) \]
\[ (13\ 24)(16\ 19)(17\ 26)(21\ 27) \]
\[ \texttt{P2} := \sigma_{j,k,1} =
(1\ 14)(2\ 20)(3\ 23)(4\ 18)(5\ 26)(6\ 27)(9\ 11) \]
\[ (10\ 15)(12\ 13)(16\ 24)(17\ 21)(19\ 22) \]
\[ \texttt{P3} := \sigma_{i,j,,2} =
(1\ 6)(2\ 8)(3\ 4)(5\ 7)(9\ 27)(10\ 19)(11\ 17)(12\ 22) \]
\[ (14\ 23)(16\ 25)(18\ 21)(20\ 26) \]
\[ \texttt{P4} := \sigma_{j,k,2} =
(1\ 13)(2\ 24)(3\ 17)(4\ 19)(5\ 26)(7\ 25)(9\ 15) \]
\[ (10\ 11)(12\ 14)(16\ 20)(18\ 22)(21\ 23) \]

The $\mathsf{GAP}$ 4
\index{gap@$\mathsf{GAP}$}%
program gives information about the group
$S_{3,3}$ which we denote $\texttt{G} = \langle \texttt{P1},
\texttt{P2}, \texttt{P3}, \texttt{P4}  \rangle $.  Rather than give the $\mathsf{GAP}$ 4 code, we merely summarize the information that we obtained about $\texttt{G}$, for example that its order is 108, together with the order and generators of the subgroups of $\texttt{G}$, and also which of these are primitive,
\index{group ! permutation ! primitive}%
Sylow,
\index{group ! Sylow subgroup}%
normal, nilpotent,
\index{group ! nilpotent}%
Hall,
\index{group ! abelian}%
abelian, transitive,
\index{group ! permutation ! transitive}%
maximal and so on.  It is found that the Sylow $2$-subgroups of $G$ are not normal subgroups and so not unique, whilst $\texttt{G}$ has a normal and therefore unique Sylow $3$-subgroup.  Also $\texttt{G}$ is not nilpotent which implies that the the Sylow $2$-subgroup acts non-trivially on the Sylow $3$-subgroup because in a finite group if the Sylow
subgroups commute then the group is nilpotent.

That the orbit of $\texttt{G}$ has $27$ elements indicates a transitive
group action
\index{group ! action}%
 on these $27$ graphs, all of which lie in a single equivalence class.  In other words, starting from any of the graphs, we can sequentially switch to get any
of the other $26$.  This proves the following theorem:-

\begin{theorem}
\label{s33}
Any two $3$-vertex $3$-colour graphs are switching-equivalent.
\end{theorem}

\begin{proof}
Given by the $\mathsf{GAP}$ 4 output.
\end{proof}

By using $\mathsf{GAP}$ 4 to give the order of the Sylow subgroup elements, we finally
arrive at
\[ S_{3,3} \cong (C_{3} \times C_{3} \times C_{3}) \sd (C_{2} \times C_{2}) \]
This is the semidirect product of an elementary abelian $3$-group
\index{group ! elementary abelian}%
and
the Klein $4$-group,
\index{group ! Klein}%
the latter being the Sylow
\index{group ! Sylow subgroup}%
$2$-subgroup $C_{2}\times C_{2}$.  So
$S_{3,3}$ is metabelian,
\index{group ! metabelian}%
being an abelian extension of an abelian group.  

The theory of groups with regular normal subgroups
\index{group ! regular normal subgroup}%
is relevant here.  Now $N \cong C_{3} \times C_{3} \times C_{3}$, $N \vartriangleleft
S_{3,3}$, and from the $\mathsf{GAP}$ 4 output $N$ is transitive on the $27$ graphs.  A
transitive abelian group is regular, so $N$ is a regular normal subgroup of
$S_{3,3}$.  So
$S_{3,3} = N \sd C_{2}\times C_{2} \leq N \sd \Aut(N)$.  In our case $S_{3,3}$ is transitive on $\mathcal{G}_{3,3}$,
(a result proved above in Theorem~\ref{s33} and below for general $m, n$ in Theorem~\ref{strg}
(a)), with $C_{3} \times C_{3} \times C_{3}$ as the regular
normal subgroup being normalized by $C_{2} \times C_{2}$,
this latter assuming the role of the point stabilizer $G_{\alpha}$.

If $V$ is a vector space, $N$ its additive group, and $H = GL(V)$,
then $N \sd H \cong AGL(V)$~\label{AGL(V)} the \emph{affine general linear group}
\index{group ! affine}
(A transitive group $G$ is \emph{affine} if it has an elementary abelian regular normal subgroup.)  In our case, $N$ is a three-dimensional vector space over
$\mathbb{F}_{3}$, with $\Aut(N) \cong GL(3, {\mathbb{F}}_{3})$.
Next we derive the form of the elements of order $2$ in $GL(3,
{\mathbb{F}}_{3})$.  If a matrix $A$ is taken to represent an
involution, then $A^{2} = 1 \Rightarrow (A + I) (A - I)
= 0$.  This is the minimal polynomial, so $A$ is diagonalizable and has
eigenvalues $\pm 1$.  Hence, $A$ is similar to one of the following
three matrices, $A_{1}$, $A_{2}$ or $A_{3}$, which up to conjugacy are:-

(i) $A_{1} =$
$ \begin{pmatrix} 
1&0&0\\
0&1&0\\
0&0&-1
\end{pmatrix}$

This corresponds to $9$ fixed points and  $\sharp$\ $2$-cycles is $(27 - 9) / 2 = 9$.

(ii) $A_{2} =$
$ \begin{pmatrix} 1&0&0\\
0&-1&0\\
0&0&-1\end{pmatrix}$

This corresponds to $3$ fixed points and $\sharp$\ $2$-cycles is $(27 - 3) / 2 = 12$.

(iii) $A_{3} =$
$ \begin{pmatrix} -1&0&0\\
0&-1&0\\
0&0&-1\end{pmatrix}$

This corresponds to $1$ fixed point and $\sharp$\ $2$-cycles is $(27 - 1) / 2 = 13$.

So the possible Klein $4$-group
\index{group ! Klein}%
is $V \cong \langle I, A, B, C\rangle$
where each of the three non-identity elements $A, B, C$ is similar to
one of the elements $A_{1}, A_{2}$, or $A_{3}$.
$\mathsf{GAP}$ 4 gives $V$ as having two generators, each having $12$
orbits.  So in our problem $A = A_{2}$, with three fixed points.  By
inspection of the $\mathsf{GAP}$ 4 orbit outputs, the three fixed points of the first orbit
(that is generator) are $1, 2, 9$ and those of the second orbit are
$1, 3, 10$.  So the fixed point of the whole group is $\{1, 2, 9\}
\cap \{1, 3, 10\} = 1$.  So our particular Klein group turns out to be
that for which all of $A, B, C$ are similar to $A_{2}$.

Another representation in terms of matrices is given by

$A =$
$ \begin{pmatrix} 
1&0&0\\
0&-1&0\\
0&0&-1
\end{pmatrix}$,
$B =$
$ \begin{pmatrix} -1&0&0\\
0&1&0\\
0&0&-1\end{pmatrix}$,
$C =$
$ \begin{pmatrix} -1&0&0\\
0&-1&0\\
0&0&1\end{pmatrix}$

Each element of the Klein $4$-group has $3$ fixed points and $12$ transpositions.  This has a geometrical interpretation in the
projective plane.
\index{projective plane}%
  A collineation is an automorphism of a projective
geometry $PG(n,\mathbb{F})$, mapping point sets and subspaces to each
other.  Each of $A, B, C$ fix 3 points forming a basis, and thus one side (2 vertices) of the triangle together with the opposite vertex pointwise.

$$
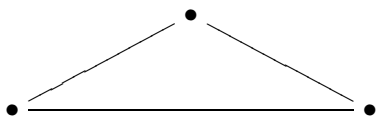
\xymatrix{
&&& {\bullet} \ar@{-}[drr] \ar@{-}[dll]\\
{} \ar@{}[r] & {\bullet} \ar@{-}[rrrr] &&&& {\bullet}
}$$

\section{The Switching Reduct for Three Colours}

To begin with we repeat some introductory material from earlier.  Given a complete graph $\Gamma$ with edges coloured red, blue and green, we
define
\emph{blue-green switching} of $\Gamma$ with respect to a set $X$ of vertices
as
follows:
\begin{itemize}
\item[(a)] for edges within $X$, or disjoint from $X$, the colours are
unchanged;
\item[(b)] for edges with one end in $X$, the colour red is unchanged,
while
the colours blue and green are interchanged.
\end{itemize}
As with the usual concept of switching, this gives an equivalence
relation on
the set of graphs on a given vertex set: switching with respect to $X$
and
then with respect to $X'$ is the same as switching with respect to the
symmetric difference of these sets.
A \emph{switching-automorphism}
\index{group ! switching ! automorphism}%
of $\Gamma$ is an isomorphism from $\Gamma$ to a graph equivalent to $\Gamma$ under switching.

The \emph{triality graph}
\index{graph ! triality}%
$\mathfrak{R^{t}}$ is the countable universal homogeneous
edge-colouring with three colours (which we take to be red, blue and
green).

\begin{lemma}
If $\mathfrak{R^{t}}$ is switched with respect to a finite set $X$ then
the resulting graph is isomorphic to $\mathfrak{R^{t}}$.
\end{lemma}

\begin{proof} 
Let $\Gamma'$ be the result of switching $\Gamma$ with respect to $X$, and let
$U,V,W$
be finite disjoint subsets of the vertex set. We seek a vertex $z$
joined to
these three sets by edges of colours red, blue, green respectively in
$\Gamma'$.
Now any vertex $z$ will do this if it has the following properties in
$\Gamma$:
\begin{itemize}
\item[(a)] $z\notin X$;
\item[(b)] $z$ is joined to $U$ by red edges;
\item[(c)] $z$ is joined to $(V\setminus X)\cup (W\cap X)$ by blue
edges.
\item[(d)] $z$ is joined to $(W\setminus X)\cup (V\cap X)$ by green
edges.
\end{itemize}
But such a vertex exists in $\Gamma$ by definition.
\end{proof}

\begin{theorem}
\label{trreduct}
Let $\mathfrak{R^{t}}$ be the triality graph.
\index{graph ! triality}%
  Then blue-green switching defines a
(proper) reduct of $\mathfrak{R^{t}}$.
\end{theorem}

\begin{proof} We have to show two things:
\begin{itemize}
\item[(a)] the group of switching-automorphisms is strictly larger than
the
group of automorphisms;
\item[(b)] this group is closed, that is, it is the automorphism group
\index{group ! automorphism}%
of some relational structure.
\end{itemize}
\index{relational structure}%

Part (a) is equivalent to showing that it is possible to switch $\Gamma$ into
a
graph which is isomorphic to it but not identical; then the isomorphism
is a
switching-automorphism of $\Gamma$ which is not an automorphism. In fact switching with respect to any non-empty finite set $X$ does this.

For part (b), we observe that switching preserves the set of
red
edges, and also the parity of the number of blue edges in a blue-green
triangle. The next result shows that a permutation which preserves both red edges and the parity of the number of blue edges in a blue-green triangle can be realized by a switching, and so that these properties define the
appropriate reduct.
\end{proof}

\bigskip

Consider the situation where there are
three colours called red, blue and green, and only blue-green switchings are
permitted.  This kind of switching has a geometrical interpretation. We are given a set
of lines in Euclidean space
\index{Euclidean space}%
making angles $\pi/2$ and $\alpha$. Choose unit
vectors along the lines; their Gram matrix
\index{Gram matrix}%
 has the form $I+(\cos\alpha)A$,
where $A$ is a matrix with entries $0$ and $\pm1$. If colours red, blue,
green correspond to entries $0$, $+1$, $-1$ respectively, then changing the
sign of a set of vectors corresponds to blue-green switching.

\smallskip

Blue-green switching clearly leaves all red edges unchanged. It also
preserves an analogue of a two-graph,
\index{graph ! two-graph}%
\index{two-graph}%
 namely, the parity of the number of green edges (say) in any blue-green triangle. Is the converse true?
Let us say that two $3$-coloured complete graphs on $V$ are
\emph{P-equivalent}~\label{P-equivalent}
\index{graph ! P-equivalent}%
 if they have the same red edges and each blue-green
triangle has the same parity of the number of green edges; and
\emph{S-equivalent}~\label{S-equivalent}
\index{graph ! S-equivalent}%
 if one can be obtained from the other
by blue-green switching.

\smallskip

P-equivalence does not imply S-equivalence in general. Suppose that $\Gamma$
consists of a blue $n$-cycle (with $n\ge4$), all other edges red. By
switching, we can make any even number of edges in the cycle green; but any
replacement of blue by green gives a P-equivalent graph. 

\bigskip

However, the following is true~\cite{camtar}:

\begin{theorem}
\label{psequiv}
\begin{itemize}
\item[(a)] Any $3$-coloured complete graph which is P-equivalent to the
countable random $3$-coloured complete graph $\mathfrak{R^{t}}$ is S-equivalent to~$\mathfrak{R^{t}}$.
\item[(b)]  Let $\Gamma$ be a random finite $3$-coloured complete graph. Then
the probability of the event that every $3$-coloured complete graph
P-equivalent to $\Gamma$ is S-equivalent to $\Gamma$ tends to~$1$ as
$n\to\infty$.
\end{itemize}
\end{theorem}

\begin{proof}
(a) Suppose that $\Gamma_1$ is the random $3$-coloured complete graph $\mathfrak{R^{t}}$,
and $\Gamma_2$ is a graph which is P-equivalent to $\Gamma_1$.
We begin with some notation. We let $c_i(xy)$ denote the colour of the edge
$\{x,y\}$ in $\Gamma_i$, and $R_i(v)$, $B_i(v)$, $G_i(v)$ the sets of vertices
joined to $v$ by red, blue, or green edges respectively in $\Gamma_i$, for
$i=1,2$. We let $BG_i(v) = B_i(v)\cup G_i(v)$. In the proof we shall modify
the graph $\Gamma_2$ so that various colours or sets become the same; once we
know that, for example, $c_1(xy)=c_2(xy)$, we drop the subscript. Note that we
can immediately write $R(v)$ and $BG(v)$, by the definition of P-equivalence.

Let $\Delta(v)$ be the symmetric difference of $B_1(v)$ and $B_2(v)$.
Switching $\Gamma_2$ with respect to $\Delta(v)$ gives a new graph $\Gamma'_2$
such that all edges containing $v$ have the same colour in $\Gamma_1$ and
$\Gamma'_2$. Now replacing $\Gamma_2$ by $\Gamma'_2$, we may assume that this
holds for $\Gamma_2$.

Now the subgraphs on $\{v\}\cup BG(v)$ are identical in $\Gamma_1$ and
$\Gamma_2$. For let $x,y\in BG(v)$. If $c_1(xy)$ is red, the result is
clear. Otherwise, $c_1(vx)=c_2(vx)$ and $c_1(vy)=c_2(vy)$, and so
$c_1(xy)=c_2(xy)$ by hypothesis.

Next we claim that, for any two vertices $x,y\in R(v)$, the edges from $x$
and $y$ to $BG(v)$ are either of the same colour in the two graphs, or differ
by an interchange of blue and green. Suppose that $c_1(xz)=c_2(xz)$ is blue
or green for some $z\in BG(v)$. Let $z'\in BG(v)$ be another point such that
$c_1(xz')$ is blue or green; we must show that $c_1(xz')=c_2(xz')$. If
$c(zz')$ is blue or green, then this assertion follows by hypothesis.
But since $\Gamma_1\cong \mathfrak{R^{t}}$, the blue-green graph on $BG(v)\cap BG(x)$ is
connected. (The induced structure on this set is isomorphic to $\mathfrak{R^{t}}$, so any
two vertices in $BG(v)\cap BG(x)$ are joined by a blue-green path of length
at most~$2$). So the claim follows.

Now $R(v)=R^+(v)\cup R^-(v)$, where, for $x\in R(v)$, the colours of edges
from $x$ to $BG(v)$ are the same in $\Gamma_1$ and $\Gamma_2$ if
$x\in R^+(v)$, and differ by a blue-green exchange if $x\in R^-(v)$. Let
$\Gamma'_2$ be obtained by switching $\Gamma_2$ with respect to $R^-(v)$.
This switching doesn't change the colours in $\{v\}\cup BG(v)$, and has the
result that $R^-(v)$ is empty in the switched graph. Replacing $\Gamma_2$ by
$\Gamma'_2$, we may assume that edges between $R(v)$ and $BG(v)$ have the
same colour in $\Gamma_1$ and $\Gamma_2$.

Finally, take $x,y\in R(v)$ with $c_1(xy)$ blue or green. Again, since
$\Gamma_1\cong \mathfrak{R^{t}}$, there exists $z\in BG(v)$ such that $c(xz)$ and $c(yz)$
are each blue or green. That is, such a configuration must exist within
the triality graph
\index{graph ! triality}%
 by universality.  Then by hypothesis we ensure that $c_1(xy)=c_2(xy)$. So
$\Gamma_1=\Gamma_2$. Since we switched the original $\Gamma_2$ twice in the
course of the proof, the proposition is proved.

\medbreak
(b) The above argument only depends on the existence of one vertex joined
to at most four given vertices by blue or green edges. It is straightforward
to show that this property holds in almost all random $3$-coloured finite
complete graphs. We omit the details.
\end{proof}

\begin{corollary}
The group of switching automorphisms corresponding to blue-green switchings
is a reduct of $\mathfrak{R^{t}}$.
\end{corollary}

\begin{proof}
This group preserves the first-order structure
\index{first-order structure}%
 with two relations, namely adjacency in the red edges in the graph, and the ternary relation picking out the
blue-green triangles with an odd number of green edges.
\end{proof}

The theorem can be expressed in another way, following~\cite{cama}.  Consider the red graph as given. Let $\mathcal{C}_2$~\label{2dimcomp} be the
$2$-dimensional complex whose simplices are the vertices, edges, and triangles
in the blue-green graph. Then P- and S-equivalence classes of the colouring
with no green edges are $1$-cocycles and $1$-coboundaries over
$\mathbb{Z}_{2}$; so the cohomology group
\index{group ! cohomology}%
 $H^1(\mathcal{C}_2,\mathbb{Z}_{2})$~\label{cohgp}
measures the extent to which P-equivalence fails to imply S-equivalence.
The theorem asserts that this cohomology vanishes in the infinite random
structure and in almost all finite structures.  The Theorem of Mallows and Sloane
\index{Mallows-Sloane Theorem}%
 whose own proof differs from that we gave above, has yet another proof~\cite{cama} using the same algebro-geometric formalism of this paragraph.

\medbreak

The most general type of restricted switching works as follows. Let $G$ be a
group of permutations on the set of colours (a subgroup of $\Sym(m)$). A
switching operation
\index{switching ! operation}%
 has the form $\sigma_{g,X}$ for $g\in G$ and $X\subseteq
V$; it applies the permutation~$g$ to the colours of edges between $X$ and
its complement, and leaves other colours unaltered. We refer to this
operation as \emph{$G$-restricted switching}.
\index{switching ! $G$-restricted}%

We record one simple observation. If the derived group $[G,G]$~\label{dergp} of $G$ is a
transitive subgroup of $G$, then all $m$-coloured graphs on the finite set
$V$ are equivalent under $G$-restricted switching. The proof is the same as
that of Theorem~\ref{strg} below.

\head{Remarks} 

1.   For ordinary switching (the two-coloured case,
corresponding to switching on blue and green edges only), it is well known that two
graphs
are switching-equivalent if and only if the parities of the numbers of
blue
edges in any triangle are the same in the two graphs. This is not true
for
three colours; whilst S-equivalence
implies P-equivalence, the example of the $n$-cycle for $n > 3$ given before Theorem~\ref{psequiv} shows
that the converse is false.

2.  For any number of colours equal to or greater than three, the
group of switchings for a given pair of colours (leaving all
others unchanged) is a reduct of the corresponding universal homogeneous
structure.  This is proved by the previous argument, by just regarding
all colours except blue or green as red.  This
gives the triality graph.
\index{graph ! triality}%

\section{Switching on More Than Two Colours: The Finite Case}

The theory of switching with more than two colours is radically different
from the two-coloured case.  We begin with a useful lemma which in
fact applies to graphs with either a finite or an infinite number of vertices.

\begin{lemma}
\label{tinded}
Consider the $m$-edge-coloured graph $\Gamma$, where $m \ge 3$, (or any of its finite subgraphs).  The colour of any edge can be switched to any other colour independently of all other edges.
\end{lemma}
\begin{proof} 
Let the three colours be $\mathfrak{r}, \mathfrak{b}, \mathfrak{g}$.
Let the edge between vertices $1$ and $2$ be coloured $\mathfrak{r}$ and say we
wanted to switch this to $\mathfrak{g}$ leaving unchanged all other colours on
all other edges.  We have the following situation:

$$
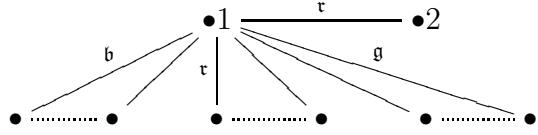
\xymatrix{ 
&& {\bullet}1 \ar@{-}[d]_{\mathfrak{r}}  \ar@{-}[dll]_{\mathfrak{b}}  \ar@{-}[dl] \ar@{-}[dr]
\ar@{-}[drr]  \ar@{-}[drrr]^{\mathfrak{g}} \ar@{-}[rr]^{\mathfrak{r}} && {\bullet}2 \\
{\bullet} \ar@{.}[r] & {\bullet} & {\bullet} \ar@{.}[r] & {\bullet} &
{\bullet} \ar@{.}[r] & {\bullet}
}$$

We want to prove that we
can make all edges joined to $1$ have colour $\mathfrak{b}$ except $\{1, 2\}$, so that we can
simply switch about vertex $1$ from $\mathfrak{r}$ to $\mathfrak{g}$, and undo all the other
unwanted colour changes.
Let $v_{{\mathfrak{r}}_{1}}, v_{{\mathfrak{r}}_{2}}, \ldots$ be those vertices
attached to vertex $1$ with red-coloured edges $e_{{\mathfrak{r}}_{1}}, e_{{\mathfrak{r}}_{2}},
\ldots$ (apart from $\mathfrak{r}$ between $\{1, 2\}$).
Let $v_{{\mathfrak{g}}_{1}}, v_{{\mathfrak{g}}_{2}}, \ldots$ be those vertices
attached to vertex $1$ with green-coloured edges  $e_{{\mathfrak{g}}_{1}}, e_{{\mathfrak{g}}_{2}},\\
\ldots$.
Perform a switch $\sigma_{\mathfrak{r}, \mathfrak{b}}$ about $v_{{\mathfrak{r}}_{1}},
v_{{\mathfrak{r}}_{2}}, \ldots$ and a switch
$\sigma_{\mathfrak{g}, \mathfrak{b}}$ about
$v_{{\mathfrak{g}}_{1}}, v_{{\mathfrak{g}}_{2}}, \ldots$.  Then we get 
$$
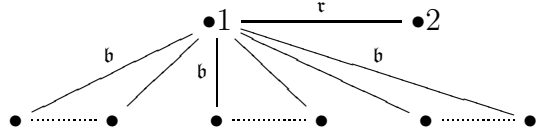
\xymatrix{ 
&& {\bullet}1 \ar@{-}[d]_{\mathfrak{b}} \ar@{-}[dll]_{\mathfrak{b}}  \ar@{-}[dl] \ar@{-}[dr]
\ar@{-}[drr]  \ar@{-}[drrr]^{\mathfrak{b}} \ar@{-}[rr]^{\mathfrak{r}} && {\bullet}2 \\
{\bullet} \ar@{.}[r] & {\bullet} & {\bullet} \ar@{.}[r] & {\bullet} &
{\bullet} \ar@{.}[r] & {\bullet}
}$$

Now perform the switch $\sigma_{\mathfrak{r}, \mathfrak{g}}^{'}$ about vertex $1$ to get the
requisite edge colour change on $\{1, 2\}$.  Next redo the switchings
$\sigma_{\mathfrak{r}, \mathfrak{b}}$ about $v_{{\mathfrak{r}}_{1}},
v_{{\mathfrak{r}}_{2}}, \ldots$ and
$\sigma_{\mathfrak{g}, \mathfrak{b}}$
about $v_{{\mathfrak{g}}_{1}}, v_{{\mathfrak{g}}_{2}}, \ldots$ in the reverse order.  Because the
switchings are involutions, this undoes the unwanted changes leaving as
the only change that edge $\{1, 2\}$ goes from $\mathfrak{r}$ to
$\mathfrak{g}$, and restoring the colours on the edges adjoining
$v_{{\mathfrak{r}}_{i}}$ and $v_{{\mathfrak{g}}_{j}}$.  
\end{proof}

Whilst the last lemma applies to graphs with either a finite or an infinite
number of vertices, the next theorem is restricted to finite graphs.  It shows that there is only one
switching class, and that the group of switching automorphisms
\index{group ! switching ! automorphism}%
of any graph is the symmetric group.

\begin{theorem}
\label{strg}
Let $m$ and $n$ be positive integers with $m>2$.
\item{(a)} The group $S_{m,n}$ acts transitively on the set
$\mathcal{G}_{m,n}$ of $m$-edge-coloured graphs on $n$ vertices.
\item{(b)} For any $m$-edge-coloured graph $\Gamma$ on $n$ vertices, the
group $\SA(\Gamma)$ of switching automorphisms of $\Gamma$ is the symmetric
group $\Sym(n)$.
\end{theorem}
\begin{proof}
(a) Observe that $[\sigma_{c,d,\{x\}}, \sigma_{d,e,\{y\}}] =
(\sigma_{c,d,\{x\}} \sigma_{d,e,\{y\}})^{2}$.  This has the effect of
permuting the triplet of colours $(c, d, e)$ to $(d, e, c)$.  So the commutator $[\sigma_{c,d,\{x\}}, \sigma_{d,e,\{y\}}]$ induces the 3-cycle $(c,d,e)$ on the colours on the edge
$\{x,y\}$ and acts as the identity on the colours on all other edges.
Now these $3$-cycles generate the direct product of alternating
groups -- a product indexed by $2$-subsets of $\{1, \ldots, n\}$.  Hence result.

(b) Obvious, for $S_{m,n}$ acts transitively so any graph
configuration can be switched to any other.  This is equivalent to proving the group is $\Sym(n)$.  
\end{proof}

\smallskip

This theorem more or less demolishes the theory of switching classes
for more than two colours. However, the problem of determining the
structure of the group of switching operations remains and is dealt
with in a separate section.

\bigskip
\bigskip
\bigskip

\begin{theorem}
\label{hightran}
The group $\SA\left(\mathfrak{R^{t}}\right)$ of switching automorphisms of $\mathfrak{R^{t}}$ is highly vertex transitive.
\index{group ! permutation ! highly transitive}%
\index{group ! switching ! automorphism}%
\end{theorem}

\begin{proof}
Take a bijection between finite sets of vertices, $\alpha : A \to B$
in $\mathfrak{R^{t}}$.  There is a switching $\sigma \in S_{3,\omega}$
which has the property that for any $x, y \in A$, $c(x \alpha, y
\alpha)^{\sigma} = c(x, y)$, where $c(x, y)$ denotes the colour of the
edge between vertices $x$ and $y$.  Then $\alpha$ is an isomorphism from $A$
(in $\mathfrak{R^{t}}$) to $B$ (in $\mathfrak{R^{t}}^{\sigma}$).  So
it extends to an isomorphism $\beta: \mathfrak{R^{t}} \to
\mathfrak{R^{t}}^{\sigma}$.  Then by definition $\beta \in
\SAut(\mathfrak{R^{t}})$ and $\beta$ extends $\alpha$.  Therefore $\SA\left(\mathfrak{R^{t}}\right)$ is highly transitive.
\end{proof}

Whilst we have proved this theorem for $\mathfrak{R^{t}}$, a similar
proof will prove that the $m$-coloured random graphs
\index{graph ! random ! $m$-coloured}%
 are highly vertex transitive for all finite $m \ge 3$.  The edge-transitivity of $\SA\left(\mathfrak{R^{t}}\right)$ follows
from its $2$-transitivity on vertices.
\index{group ! permutation ! $2$-transitive}%
 
\head{Remarks} 

1.  That $\mathfrak{R^{t}}$ is homogeneous is already a very strong symmetry
condition; it implies the transitivity of
$\Aut\left(\mathfrak{R^{t}}\right)$ on vertices and on adjacencies of
a particular colour.

2.  That the switching group is highly transitive is
equivalent to saying that its closure is $\Sym(V)$, that is the group of all permutations whose
effects can be undone by some switching.
Alternatively it is dense in $\Sym(V)$.  (A group $H$ that is a subgroup of a closed group $G \leq \Sym(X)$ is
\emph{dense}
\index{group ! dense}%
in $G$ if and only if it has the same orbits as $G$ on
finite ordered subsets of $X$).  In the next section we will see that
$\SA\left(\mathfrak{R^{t}}\right) \neq \Sym(V)$ and so it is not
a reduct.  There are switchings $\alpha$ with respect to an infinite
number of vertices so that $(\mathfrak{R^{t}})^{\alpha} =
\mathfrak{R^{t}}$.  Indeed switching with respect to a random subset
of vertices almost surely gives $\mathfrak{R^{t}}$.

\section{Switching on More Than Two Colours: The Infinite Case}

What about infinite multicoloured graphs?

 We assume that the graphs have countably many vertices. In the two-coloured case, the situation is like that
for the finite case. But for more than two colours, there are significant
differences.

We begin by observing that $S_{m,\omega}$ does not act transitively
on $\mathcal{G}_{m,\omega}$: there is no sequence of switching operations
\index{switching ! operation}%
which maps an infinite graph with all edges of one colour to one with all
edges of a different colour. This follows from the following characterisation
of the switching class of a monochromatic complete graph.

\bigskip

\begin{lemma}
\label{eqconsw}
The following conditions on a countable $m$-coloured graph $\Gamma$ are
equivalent for $m>2$:

(a) $\Gamma$ is obtained from a $c$-coloured clique by switching;

(b) there is an equivalence relation on $\Gamma$ with finitely many
equivalence classes such that, if $x\equiv y$ and $x\ne y$, then $\{x,y\}$
has colour $c$, and $\{x,z\}$ and $\{y,z\}$ have the same colour for all
$z\ne x,y$.
\end{lemma}

\begin{proof}
(a) implies (b):
Any element of the group of switching operations is a product of finitely
many switchings $\sigma_{c_i,d_i,Y_i}$ for $i=1, \ldots, n$. Now the
infinite set $X$ is partitioned into at most $2^n$ sets $Y_J$ for
$J\subseteq\{1,\ldots, n\}$, where
\[Y_J=\{x\in X:(x\in Y_i)\Leftrightarrow(i\in J)\}.\]
Now edges joining vertices in the same set $Y_J$ have their colours unaltered
by the given sequence of switchings, so have colour~$c$; and if $x$ and $y$
belong to the same one of these sets and $z$ is any other vertex, then the
colours of $\{x,z\}$ and $\{y,z\}$ experience the same sequence of changes.

(b) implies (a): suppose that (b) holds. We form a new graph $\hat{\Gamma}$
whose vertices are the equivalence classes of vertices of $\Gamma$, with
an edge coloured $d$ from $Y_1$ to $Y_2$ if the edges from vertices in $Y_1$
to vertices in $Y_2$ have colour $d$ in $\Gamma$. The graph $\hat{\Gamma}$ is
finite; so, by Theorem~\ref{strg}, we can switch it so that all edges
have colour~$c$. The corresponding sequence of switchings of $\Gamma$ changes
the colours of all edges between equivalence classes to~$c$, and does not
alter the colours of edges within classes. So a $c$-coloured complete graph
results.
\end{proof}

A graph $\Gamma$ is a \emph{switched $c$-clique}
\index{graph ! switched $c$-clique}%
 if the equivalent conditions
of Lemma~\ref{eqconsw} hold. Note that every subgraph of a switched $c$-clique
is a switched $c$-clique. Note also that an infinite graph cannot be both a
switched $c$-clique and a switched $d$-clique for $c\ne d$, since a switched
$c$-clique has no infinite clique of colour different from~$c$.

A \emph{moiety}
\index{moiety}%
 in a countably infinite set~$X$ is an infinite subset $Y$ of
$X$ such that $X\setminus Y$ is also infinite.

\begin{lemma}
\label{lmthmo}
A countably infinite multicoloured graph is a switched $c$-clique if and only if the vertex set can be partitioned into three moieties such that the induced
subgraph on the union of any two is a switched $c$-clique.
\end{lemma}

\begin{proof}
The forward implication is clear. So suppose that $X$ is the disjoint union
of $Y_1, Y_2, Y_3$, and for each $i\ne j$, there is an equivalence relation
$\equiv_{ij}$ on $Y_i\cup Y_j$ with the properties of Lemma~\ref{eqconsw}.
Extend $\equiv_{ij}$ to an equivalence relation on $X$ in which the remaining
set $Y_k$ is a single class. Let $\equiv$ be the meet of these three
equivalence relations.

We claim that $\equiv$ has the properties of Lemma~\ref{eqconsw}. Certainly it
has only finitely many classes. Take two points $x,y$ in the same class.
Then they belong to the same set $Y_i$, say $Y_1$ without loss of generality.
Since $x\equiv_{12}y$, the edge $\{x,y\}$ has colour $c$. Now let $z$ be
any point in a different equivalence class. Suppose, without loss of
generality, that $z\in Y_1\cup Y_2$. Then the properties of $\equiv_{12}$
ensure that $\{x,z\}$ and $\{y,z\}$ have the same colour.
\end{proof}

Now we can characterise multicoloured graphs for which the group of switching
automorphisms is the symmetric group.

\begin{theorem}
\label{sanesm}
Let $\Gamma$ be an $m$-coloured graph on a countable set $X$, where $m>2$.
Then the group of switching automorphisms of $\Gamma$ is the symmetric group
on $X$ if and only if $\Gamma$ is a switched $c$-clique for some colour $c$.
\end{theorem}

\begin{proof}
Since switching does not change the group of switching automorphisms, the
reverse implication is clear. So suppose that $\SA(\Gamma)$ $= \Sym(X)$. By
Ramsey's Theorem,
\index{Ramsey theorem}%
there is a moiety $Y$ of $X$ which is a $c$-clique for some
colour $c$. Since $\Sym(X)$ is transitive on moieties, every moiety is a
switched $c$-clique. Now Lemma~\ref{lmthmo} gives the result.
\end{proof}

So the group of switching automorphisms is not in general equal to the
symmetric group; in particular if the $m$-coloured graph contains
infinite monochromatic subgraphs in two different colours.
Equivalently, if $\SA(\Gamma)$ $= \Sym(X)$ there is only one colour for
which $\Gamma$ contains an infinite monochromatic clique.  In
particular $\mathfrak{R^{t}}$ contains infinite monochromatic cliques
of all $3$ colours, so $\SA\left(\mathfrak{R^{t}}\right) \neq
\Sym(\mathfrak{R^{t}})$.  Hence $\SA\left(\mathfrak{R^{t}}\right)$ is
not a reduct of $\mathfrak{R^{t}}$.

For completeness we can give some instances of countable graphs
$\Gamma$ with at least three colours for which the group of switching
automorphisms is \emph{not} highly transitive.

\head{Example}

If $V$ is the disjoint union of two infinite sets $V_1$ and $V_2$, where
edges within $V_1$ are red, edges within $V_2$ are blue, and edges between
$V_1$ and $V_2$ are blue, then $\SAut(\Gamma)$ contains $\FSym(V)$.
For although the graph is not a switched $c$-clique, so that
$\SAut(\Gamma) \ne \Sym(V)$, because $\FSym(V)$ permutes only a
finite number of vertices and only a finite number of edge colours are
present, Lemma~\ref{tinded} gives that $\FSym(V) \le \SAut(\Gamma)$.  The
following statements can easily be proved:

(a) in this case $\SAut(\Gamma)$ is the product of $\FSym(V)$ and
$\Sym(V_1)\times\Sym(V_2)$;

(b) in general, if $\Gamma$ has a finite partition $V = V_1 \cup \ldots \cup V_k$ such that
all edges within a single part of the partition have the same colour, while
the colour of an edge from $V_i$ to $V_j$ depend only on $i$ and
$j$, then $\SAut(\Gamma)$ contains $\FSym(\Gamma)$.

The converse statement of (b) would be:  if $\FSym(V) \le
\SAut(\Gamma)$ then $\Gamma$ has a finite partition of vertex sets,
such that each part is monochromatic and the colour of edges joining
any two parts is only dependent on those parts.  This is unknown.
\medskip

A permutation $g \in \Sym(X)$ of $X$ is an \emph{almost automorphism}
\index{group ! almost automorphism}%
of a multicoloured graph $\Gamma$ on $X$ if the colour of $\{xg,yg\}$ is equal to the colour of
$\{x,y\}$ for all but finitely many $2$-element subsets of $X$. The
group of such almost automorphisms is denoted $\AAut(\Gamma)$.  In the
two-colour case, this group has been investigated by Truss
\index{Truss, J. K.}%
and others~\cite{mek}~\cite{truss2}.  

\begin{theorem}
\label{contth}
Suppose that $m>2$. Then the group $\AAut(\Gamma)$ of almost automorphisms of
any $m$-edge-coloured countably infinite graph $\Gamma$ is contained in the group $\SA(\Gamma)$ of switching automorphisms.
\end{theorem}

\begin{proof}
Just as in the finite case (Theorem~\ref{strg}), we can change the colour
of any single edge (and hence of any finite number of edges) by a sequence of
switchings.  Any finite number of colour changes can be switched back.
\end{proof}

For finite graphs we have the equivalence $\SAut(\Gamma) =
\AAut(\Gamma) = \Sym(\Gamma)$.  However the containment in Theorem~\ref{contth} is proper in the case of $\Gamma = \mathfrak{R}_{m,\omega}$.
For switching
with respect to a random subset of vertices almost surely gives
$\mathfrak{R}_{m,\omega}$, that is the exceptions which do not give
this graph form a null set in the sense of measure.
\index{measure}%
  So there is a switching automorphism changing
the colours on infinitely many edges which is therefore not an almost
automorphism. 

As we shall see in the next paragraph, the group of almost automorphisms
\index{group ! almost automorphism}%
of the universal homogeneous multicoloured graph
\index{graph ! random ! $m$-coloured}%
 is highly transitive
\index{group ! permutation ! highly transitive}%
so cannot be a non-trivial
reduct, since any closed subgroup of the symmetric group is the
automorphism group
\index{group ! automorphism}%
 of a relational structure.
\index{relational structure}%
   Furthermore
$\SAut(\mathfrak{R}_{m,\omega})$ is highly transitive.  In summary, for $m \geq 3$,
\[ \AAut(\mathfrak{R}_{m,\omega}) < \SAut(\mathfrak{R}_{m,\omega}) < \Sym(\mathfrak{R}_{m,\omega}). \]

Let $f_m$~\label{fm} be the colouring function on the set of unordered vertex pairs of $\mathfrak{R}_{m,\omega}$ into the set $\mathcal{C}_m$ of $m$ edge colours, satisfying the following universality property: if $\alpha : \mathfrak{R}_{m,\omega} \to \mathcal{C}_m$ is a finite partial map then $\exists x \in \dom(\alpha)$ such that $\forall y \in \dom(\alpha)$ $\alpha(y) = f_m \{x, y\}$.  If $m = 2$, then $\mathcal{C}_2 = \{0, 1\}$ and vertices $x, y \in \mathfrak{R}_{m,\omega}$ are joined if $f_2 \{x, y\} = 1$.  The main result required to prove that $\AAut(\mathfrak{R}_{m,\omega})$ is a highly transitive permutation group on $\mathfrak{R}_{m,\omega}$ for $m \ge 2$~\cite{camb} is the following lemma, which is proved using the universality property and a back-and-forth argument,
\index{back-and-forth method}%

\begin{lemma}
Let $p : \mathfrak{R}_{m,\omega} \to \mathfrak{R}_{m,\omega}$~\label{fpb} be a finite partial bijection.  Then $\exists g \in \AAut(\mathfrak{R}_{m,\omega})$ extending $p$, which is a permutation of $\mathfrak{R}_{m,\omega}$, and is such that if $f_m \{gx, gy\} \neq f_m \{x, y \}$ then $x, y \in \dom(p)$.
\end{lemma}

An alternative form of the lemma is

\begin{lemma}
Let edges $e_1, e_2 \in E(\mathfrak{R}_{m,\omega})$ be such $f_m(e_1) \neq f_m(e_2)$.  Then $\exists g \in \AAut(\mathfrak{R}_{m,\omega})$ with $|\{e_i \in E(\mathfrak{R}_{m,\omega}) : f_m(g e_i) \neq f_m (e_i) \}| = 1$, such that $g e_1 = e_2$.
\end{lemma}

We leave the proofs of these lemmas as exercises for the interested reader.

\section{Derivation of the Form of $S_{m,n}$}
\label{dersec}
We turn our attention to determining the form of the finite switching groups.

\begin{proposition}
The switching group $S_{m,3}$ on $3$ vertices and $m$ colours is a normal subgroup of index
$2$ in the direct product of $3$ symmetric groups of degree $m$.
\index{group ! permutation ! degree}%
\end{proposition}

\begin{proof}
Denote the switching group in question by
$S_{m,3}$.  For $n=3=m$, we have seen
above that $|S_{3,3}| = 108.$ 
Also $|(\Sym(3))^{3}| = 216$, so $|(\Sym(3))^{3} :
S_{3,3}| = 2$ and
$S_{3,3} \lhd (\Sym(3))^{3}$.

Proceed inductively on $m$.  We have
$S_{m,3}$ acting on $m^{3}$ triples of
colours, and we want to show that $|(\Sym(m))^{3} :
S_{m,3}| = 2$, for which we would require
$|S_{m,3}| = (m!)^{3}/2$.  For then, by
the inductive hypothesis, $|S_{m-1,3}| =
((m-1)!)^{3}/2$, and the assertion $|(\Sym(m))^{3} :
S_{m-1,3}|
. |S_{m,3} : S_{m-1,3} |^{-1} = 2$ would follow and would give the result if
we could show that  $|S_{m,3} : S_{m-1,3} | = m^{3}$.  So it remains to
prove this last equation.  The transitivity of
$S_{m,3}$ as given by Theorem~\ref{strg} implies that this group acts
transitively on triples of colours.  Also, since the stabilizer of the
graph with all edges coloured $m$ is $S_{m-1,3}$, and each
entry in the triple can either stay the same colour $m$, or change to
one of the other $m-1$ colours, this gives $m^{3}$ possibilities, and
proves the proposition.
\end{proof}

This leads us to our desired result, which is a generalization of
Seidel's which states that $|S_{2,n}| = 2^{n-1}.$
\index{Seidel, J. J.}%

\begin{theorem}
\label{swgpfm}
The group of all switchings on $n$-vertex subgraphs having $m$ colours
is given by
\[S_{m,n} \cong (\Alt(m))^{n(n-1)/2}\sd (C_2)^{n-1}.\]
\end{theorem}

\begin{proof}
The semidirect product action here partitions the vertex set into two parts and the $C_2$ groups transpose the colours on edges crossing the partition.  

Take any edge $\{x,y\}$
and any three colours $a,b,c$. Then $\sigma_{a,b,\{x\}}\sigma_{a,c,\{y\}}$
induces the $3$-cycle $(a,c,b)$ on the colours on $\{x,y\}$, the $2$-cycle
$(a,b)$ on edges $\{x,z\}$, the $2$-cycle $(a,c)$ on edges $\{y,z\}$ (where
$z\notin\{x,y\}$), and the identity on all other edges.  (This is the
content of Lemma~\ref{tinded}).  So the square of this
element induces $(a,b,c)$ on $\{x,y\}$ and the identity elsewhere, for $\sigma^{-1} \xi^{-1} \sigma \xi = (\sigma \xi )^2$.  Since
$3$-cycles generate the alternating group,
\index{group ! alternating}%
we get $\Alt(m)$ on the colours on
$\{x,y\}$, and hence the direct product of $n \choose 2$ alternating groups
overall.

Now there is a homomorphism from $S_{m,n}$ to $S_{2,n}$,
where the parity of the permutation on the $m$ colours at an edge determines
whether the two colours are switched or not. The kernel is the product of
alternating groups just constructed, and the image is the whole group
$S_{2,n}$, whose structure we know by the result of Seidel.
\index{Seidel, J. J.}%
\end{proof}

\section{Primitivity of Extended Switching Group Actions}
\index{group ! action}%

There are two actions of $S_{m,n}$ (and $S_{m,n} \sd \Sym(n)$):-

(i)  on the set of coloured edges.  This is imprimitive
\index{group ! permutation ! imprimitive}%
 because the ``same edge'' is a congruence,
 \index{congruence}%
that is all pairs with a fixed edge form a block of imprimitivity;

(ii)  on the set of all edge-coloured graphs.  In this, more natural
action, $S_{m,n}$ is transitive.  We will show that this action of
$S_{m,n} \sd \Sym(n)$, where $S_{m,n}$ acts as $\Sym(m)$ on
the colours of each edge and $\Sym(n)$ acts transitively on the $n
\choose 2$ pairs, is primitive
\index{group ! permutation ! primitive}%
 for $m \geq 3$.  However
$S_{m,n}$ by itself is imprimitive
\index{group ! permutation ! imprimitive}%
 in this action, for take an edge
$e$ and colour $c$, then the set of all graphs in which the colour of
$e$ is $c$ forms a block of imprimitivity. 

Take the group $G^W_{m,n}$ to be a subgroup of $\Sym(m)  \Wr \Sym{n
\choose 2}$ with the product action.
\index{group ! action ! product}%
 The group $G^W_{m,n}$ acts on the set $\Omega =
\Gamma^{\Delta}$ of functions from $\Delta$ to $\Gamma$, where colour
set $|\Gamma| = m$ and the set of vertex pairs is $|\Delta| = {n
\choose 2}$, and is \emph{non-basic},
\index{group ! permutation ! non-basic}%
that is preserving a non-trivial
product structure (defined to be a bijection $\Omega \to
\Gamma^{\Delta}$ for $|\Gamma|, |\Delta| > 1$).  The domain for the
product action is the set of global sections
\index{global section}%
 (or transversals or functions
$\Delta \to \Gamma$) for the fibres.  Each fibre represents the set of
colourings on a particular edge, and a transversal represents a
complete colouring of the graph.  The bottom groups
\index{group ! bottom}%
$\Alt(m)$ are
primitive
\index{group ! permutation ! primitive}%
 and non-regular on $\Gamma$ if $m >3$, as is necessary for a
primitive $G^W_{m,n}$ product action on $\Omega$.  Notice that
$\Alt(m) \Wr \Sym(n) \leq G^W_{m,n} \leq \Sym(m) \Wr \Sym(n)$.  Now
$\Alt(m) \Wr \Sym(n)$ is primitive
\index{group ! permutation ! primitive}%
 if $m >3$, for by an argument in
the proof of the O'Nan--Scott Theorem:
\index{O'Nan--Scott Theorem}%
\index{O'Nan, M.}%
\index{Scott, L.}%
$A \Wr B$ is primitive
\index{group ! permutation ! primitive}%
 in the product action if and only if the action of $A$ on each coordinate or
fibre of $\Gamma$ is primitive and non-regular, and
$B$ is transitive~\cite{dixon}.

The cases $m = 3$ and $m \geq 4$ can be separately dealt with.
Firstly for $m \geq 4$ we have that $S_{m,n} \sd \Sym(n)$ has
a natural action on $S_{m,n}$ and it lies between $\Alt(m) \Wr
\Sym(n)$ and $\Sym(m) \Wr \Sym(n)$ in the product action.
\index{group ! action ! product}%
For $m > 3$ we have that $\Alt(m)$ is primitive
\index{group ! permutation ! primitive}%
 and not regular and also
$\Sym(n)$ is transitive on $n \choose 2$ edges.  Firstly this is equivalent
to saying that $\Alt(m) \Wr \Sym(n)$ is primitive, and secondly it
implies that $S_{m,n} \sd \Sym(n)$ is
primitive, as it lies above a primitive group~\cite{dixon}.

For $m=3$, note that $\Alt(3) = C_3$.  Now, a group of the form $G = C_p^l
\sd H \cong C_p \Wr H$ acts on $(\mathbb{Z}/p)^l$ (an
$l$-dimensional vector space $V$ over the finite field $\GF(p)$) with $H \leq
GL(l, p)$ being transitive of degree $l$.  The group $G$ acts transitively
because $C_p^l$ does.

Consider
\[C_3^{n \choose 2} \sd \Sym(n) \leq G^W_{3,n} \leq \Sym(3)^{n \choose 2} \sd \Sym(n),\]
where the symmetric group action has degree ${n \choose 2}$.  Here
$G^W_{3,n}$ is the semi-direct product $\Big(C_3^{n \choose 2}
. C_2^{(n-1)} \Big) \sd \Sym(n)$.  In our case $l = {n \choose 2}, p=3, H = C_2^{n-1} \sd \Sym(n)$.
We can parametrize the $3$ colours as $\{0, \pm 1\}$, so that the group $C_2$ acts to stabilize
$0$ and transpose $+$ and $-$.  We need to show that $H$ is
irreducible on $V$.  Let $W$ be an $H$-invariant subspace of $V$.  Take $w \in W$ to be a function
from edges to $\{0, \pm 1\}$ with fewest nonzero coordinates.  Apply
any vertex permutation or switching fixing the $0$ coordinates and
changing the signs of the others.  We claim that all edges on which $w
\neq 0$ pass through the same vertex.  

For let $i$ be any vertex such that $w \ne 0$ on some edge through
that vertex.  Switch $+$ and $-$ at edges through $i$ and subtract
from original to get $w - w^{\sigma} = x$, where $x = 0$ on all edges
not containing $i$.  Also $x(\{i, j\}) = - w(\{i, j\})$ $\forall j$,
and is nonzero on all edges containing $i$ where $w$ is nonzero, that is $x
\in W,\ x \ne 0$.  Also $x$ has fewer non-zero coordinates than $w$
 (when $w = 0$ on edges not containing $i$).  

Further $w \ne 0$ only on a single edge.  Suppose that $w(\{i, j\}) \ne
0$.  Switching at $j$ gives $w^{\sigma}(\{i, j\}) = -
w(\{i, j\})$ and if $k \ne i, j$ $w^{\sigma}(\{i, k\}) = w(\{i, k\})$.
So $w - w^{\sigma}$ is nonzero only on $i, j$.  So $W = V$. 

As $\Sym(n)$ acts transitively on edges so for each edge there exists
$w \in W$ which is non-zero only on that edge.  These span $V$.
Therefore $W = V$.  So $G^W_{3,n}$ is primitive
\index{group ! permutation ! primitive}%
 and finally we can state the result
\begin{theorem}
\label{extswgp}
The group $S^*_{m,n}:= S_{m,n} \sd \Sym(n)$
acts primitively on $\mathcal{G}_{m,n}$, for $m \geq 3$.
\end{theorem}
\index{group ! switching ! extended}%

This result can be seen as a strengthening of Theorem~\ref{strg}.  

Once again just as we saw above for the theory of switching classes we have a contrast between the $m=2$ random graph
\index{graph ! random}%
 where $S^*_{2,n}$ has non-trivial switching orbits for $n \ge 3$, and the
higher-adjacency graphs where $S^*_{m,n}$ acts primitively for $m
\ge3$.  In fact $S^*_{2,n}$ is not even transitive; its orbits are
isomorphism classes of switching classes.  For three vertex graphs there are
two switching classes of four graphs:

\vspace{50pt}

\begin{figure}[!h]
$$\xymatrix{
& {\bullet} &&& {\bullet} \ar@{-}[dr] \ar@{-}[dl]\\
{\bullet} && {\bullet} & {\bullet} && {\bullet}
}$$
$$\xymatrix{
& {\bullet} \ar@{-}[dl] &&& {\bullet} \ar@{-}[dr]\\
{\bullet} \ar@{-}[rr] && {\bullet} & {\bullet} \ar@{-}[rr] && {\bullet}
}$$\\

and\\

$$\xymatrix{
& {\bullet} \ar@{-}[dr] \ar@{-}[dl] &&& {\bullet}\\
{\bullet} \ar@{-}[rr] && {\bullet} & {\bullet} \ar@{-}[rr] && {\bullet}
}$$
$$\xymatrix{
& {\bullet} \ar@{-}[dl] &&& {\bullet} \ar@{-}[dr]\\
{\bullet} && {\bullet} & {\bullet} && {\bullet}
}$$\caption{Two switching classes for $S_{2,3}$}  
\end{figure}
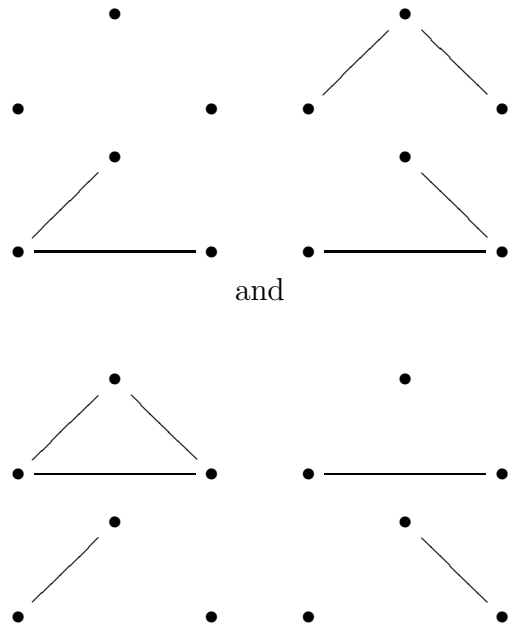 

\clearpage

\section{Transitive Extensions}

At this point, we can discern a property of the switching groups that
is alluded to by some of the results of the previous sections.

Let $H$ be a permutation group acting on a set $\Sigma$.  We
say that the permutation group $G$ is a
\emph{transitive extension}
\index{group ! permutation ! transitive extension}%
 of $H$ if it acts on
$\Omega = \Sigma \cup \{*\}$ such that $H$ is the stabilizer
of point ${*}$ acting on $\Sigma$.  A \emph{weak transitive
extension} is a transitive permutation group $G$ acting on
$\Omega = \Sigma \cup \{*\} \cup \Gamma$ such that $G_{\{*\}}$ setwise
stabilizes $\Sigma$ and acts on $\Sigma$ as $H$.
Here $\Gamma$ is a set of points such that $\Sigma \cap \Gamma =
{\emptyset}$.  Every multiply transitive
\index{group ! permutation ! multiply transitive}%
 group arises from some other
group by transitive extension. 

As an illustration of the difference between transitive extensions and
weak transitive extensions,
\index{transitive extension ! weak}%
we note that the Klein $4$-group
\index{group ! Klein}%
has no transitive extension~\cite[p.~232]{dixon}, so we would not get a transitive extension with
respect to a $5$-vertex set, but there is a weak one.  This is further
illustrated by the following two examples:

\bigskip

\head{Examples of Weak Transitive Extensions}

1.  The action of $\Sym(3) \times \Sym(3)$ on the nine points
$\{(i,j): 1 \le i,j \le 3 \}$.  The stabilizer of the point $* = (3,
3)$ is  $\Sym(2) \times \Sym(2) \cong V_4$ with $\Sigma = \{(1,1),
(1,2), (2,1), (2,2)\}$.  

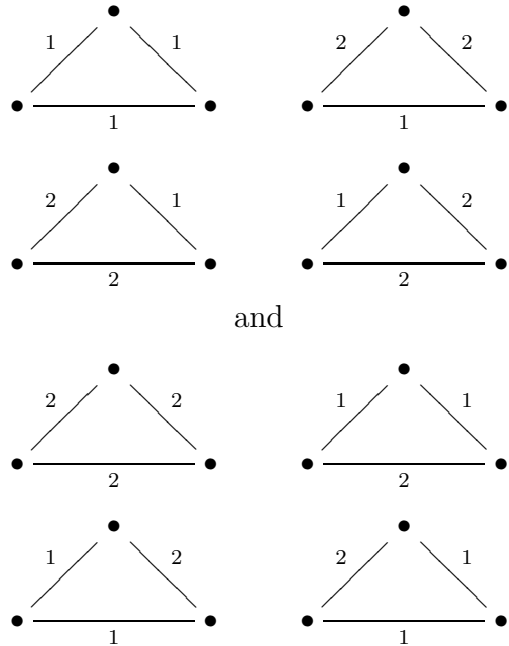
\begin{figure}[!h]$$\xymatrix{
& {\bullet} \ar@{-}[dr]^{1} \ar@{-}[dl]_{1} &&& {\bullet} \ar@{-}[dr]^{2} \ar@{-}[dl]_{2}\\
{\bullet} \ar@{-}[rr]_{1} && {\bullet} & {\bullet} \ar@{-}[rr]_{1} && {\bullet}
}$$
$$\xymatrix{
& {\bullet} \ar@{-}[dr]^{1} \ar@{-}[dl]_{2} &&& {\bullet} \ar@{-}[dr]^{2} \ar@{-}[dl]_{1}\\
{\bullet} \ar@{-}[rr]_{2} && {\bullet} & {\bullet} \ar@{-}[rr]_{2} && {\bullet}
}$$
and
$$\xymatrix{
& {\bullet} \ar@{-}[dr]^{2} \ar@{-}[dl]_{2} &&& {\bullet} \ar@{-}[dr]^{1} \ar@{-}[dl]_{1}\\
{\bullet} \ar@{-}[rr]_{2} && {\bullet} & {\bullet} \ar@{-}[rr]_{2} && {\bullet}
}$$
$$\xymatrix{
& {\bullet} \ar@{-}[dr]^{2} \ar@{-}[dl]_{1} &&& {\bullet} \ar@{-}[dr]^{1} \ar@{-}[dl]_{2}\\
{\bullet} \ar@{-}[rr]_{1} && {\bullet} & {\bullet} \ar@{-}[rr]_{1} && {\bullet}
}$$\caption{Example of weak transitive extension of $S_{3,3}$}
\label{wktrex}  
\end{figure}

2.  Consider the action of $S_{3,3}$ acting on ordered triples with
    elements $1, 2, 3$,$$
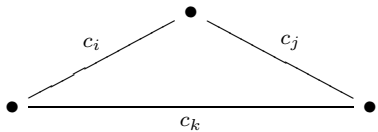
\xymatrix{
&&& {\bullet} \ar@{-}[drr]^{c_j} \ar@{-}[dll]_{c_i}\\
{} \ar@{}[r] & {\bullet} \ar@{-}[rrrr]_{c_k} &&&& {\bullet}
}$$Recall that $|S_{3,3}| = 27 . 4$.  Here $G = S_{3,3}$ and $$ * = 
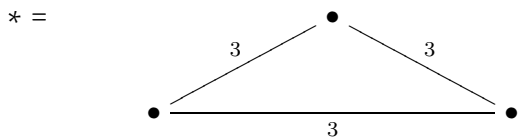
\xymatrix{
&&& {\bullet} \ar@{-}[drr]^{3} \ar@{-}[dll]_{3}\\
{} \ar@{}[r] & {\bullet} \ar@{-}[rrrr]_{3} &&&& {\bullet}
}$$  The pointwise stabilizer $G_{\{*\}} \cong V_4$ has two orbits of $4$ graphs each as shown in Figure~\ref{wktrex}.
    
Contrast this to the case where $m \ge 4$ where there is only one
switching class.

As a corollary to the high-transitivity of $S_{m,n}$ we find,

\begin{theorem}
\label{wte}
For $3 \le n \le \aleph_0$, $3 \le m < \aleph_0$ and for a fixed
value of $n$, the groups in the series of multiply transitive groups
given by 
\[ S_{3,n} \le S_{4,n} \le \ldots \le S_{m-1,n} \le S_{m,n} \le \ldots \]
each form a \emph{weak transitive extension} of the previous group
in the series, in their action on the set of all $m$-coloured graphs
on $n$ vertices.
\end{theorem}

\begin{proof}
The stabilizer in $S_{m,n}$ of a graph $X$ with all edges coloured $m$,
fixes the set of all graphs not using colour $m$ and acts as
$S_{m-1,n}$ on this set.  Clearly $S_{m-1,n} \le S_{m,n}$ and
fixes $X$.  So $S_{m-1,n} \le
(S_{m,n})_X$.  So by Theorem~\ref{swgpfm} for the first
equality and the Orbit-Stabilizer Theorem
\index{Orbit-Stabilizer Theorem}%
for the last, $|S_{m,n} : S_{m-1,n} | = m^{n \choose 2} = \sharp
\text{ coloured graphs on}\ m\ \text{colours} = |S_{m,n} : (S_{m,n})_X|$.  Therefore from
the containment and dimension arguments, $S_{m-1,n} = (S_{m,n})_X$.
\end{proof}

\smallskip

We end the section with an observation on generalized two-graphs.
\index{graph ! two-graph}%
\index{two-graph}%

Recall that two $3$-coloured complete graphs on $V$ are
\emph{P-equivalent} if they have the same red edges and each blue-green
triangle has the same parity of the number of green edges.  This can be extended to any number of colours by defining two $m$-coloured complete graphs as \emph{$P_m$-equivalent}~\label{P_m-equivalent}
\index{graph ! Pmeq@$P_m$-equivalent}%
 if they have the same edges on $m-2$ colours and the same parity of the number of green edges in each blue-green triangle.

\begin{figure}[!h]
$$\xymatrix{
{\ldots} \ar@{^{(}->}[r]^{\iota} & {\Aut(P_4(\mathfrak{R}_{4, \omega}))} \ar@{^{(}->}[r]^{\iota} & {\Aut(P_3(\mathfrak{R^{t}}))}
\ar@{^{(}->}[r]^{\iota} & {\Aut(\mathcal{T}(\mathfrak{R}))}\\
{\ldots} \ar@{^{(}->}[r]^{\iota} & {\Aut(\mathfrak{R}_{4, \omega})} \ar@{^{(}->}[u]^{e} \ar@{^{(}->}[r]^{\iota} & {\Aut(\mathfrak{R^{t}})} \ar@{^{(}->}[u]^{e} \ar@{^{(}->}[r]^{\iota} & {\Aut(\mathfrak{R})} \ar@{^{(}->}[u]^{e}
}$$ 
\caption{Relations between generalized two-graph and random graph automorphism groups}
\label{tgrg}
\end{figure}

Adding the point at infinity~\label{ptatinf} to the vertex set of an $m$-coloured random graph
\index{graph ! random ! $m$-coloured}%
 to get $V(\mathfrak{R}_{m, \omega}^{+}) = V(\mathfrak{R}_{m, \omega}) \cup \{\mathfrak{O}\}$ gives a new structure, $\mathcal{T}(\mathfrak{R}_{m, \omega})$ such that $\Aut(\mathcal{T}(\mathfrak{R}_{m, \omega}))$ is a transitive extension of $\Aut(\mathfrak{R}_{m, \omega})$.  For $m=2$, $\mathcal{T}(\mathfrak{R})$ is the unique
countable universal homogeneous two-graph.
\index{graph ! two-graph}%
\index{two-graph}%
  In fact $\Aut(\mathcal{T}(\mathfrak{R}))$ is a \emph{curious transitive
  extension}
\index{transitive extension ! curious}%
of $\Aut(\mathfrak{R})$~\cite{cam7}, that is it has a
transitive subgroup on $V(\mathfrak{R}^{+})$ which is isomorphic to $\Aut(\mathfrak{R})$.  This carries through for $m \ge 3$ and we obtain Figure~\ref{tgrg} as a result, where the $\iota$~\label{iota} maps are inclusions and the $e$ maps are curious transitive extensions.  The $\iota$ inclusions follow by considering the partition of the sets that make up the generalized two-graphs on different numbers of colours, and going colourblind in successive pairs of colours.

\chapter{Finitary Switching Groups}
\label{chap3}
\bigskip

To my own Gods I go.\\
It may be they shall give me greater ease\\
Than your \ldots tangled Trinities.
\begin{flushright}
Rudyard Kipling, \textit{Plain Tales from the Hills, `Lisbeth'}
\end{flushright}

\bigskip  

The next dream I want to present is an even more fantastic set of
theorems and conjectures.  Here I also have no theory and actually the
ideas form a kind of religion rather than mathematics.  The key
observation is that in mathematics one encounters many trinities.
\begin{flushright}
V. I. Arnold, \textit{The Arnoldfest, Fields Institute Communications
vol.24, p.32, eds. E. Bierstone \emph{et al.}, 1999}
\end{flushright}

\medskip

In general, whilst few permutation groups are finitary,
\index{group ! finitary}%
there are many infinitary ones such as any infinite group acting by the right regular action; in fact regular infinite groups are as cofinitary
\index{group ! cofinitary}%
as it is possible to achieve.
 
This chapter deals with variations on the switching group theme, beginning with a specific instance of finitary action, that is assuming that there are infinitely many colours, but allowing switchings
only on the first $m$ colours, the colours $m+1, m+2, \dots$
all being stabilized under group action.
\index{group ! action}%
  Whilst the switchings are finitary the action of the group on vertices will not be.  Label the infinitely many colours $1, 2, 3, \ldots$ and let $\mathcal{C}_m = \{1, \ldots, m\}$ denote the set of the first $m$ colours. 
Define the
\emph{finitary switching group}
\index{group ! switching ! finitary}%
acting on the relevant infinite graphs, to be the group
generated by switchings of pairs of the first $m$ colours, writing this as $S^{f}_{m,\omega} = \langle
\sigma_{c, d, X} : X \subseteq V \text{and}\ c, d \in \mathcal{C}_m \rangle$, where $V$ is the
vertex set of the graph, $|V| = \omega$.  The union $\bigcup_{m \geq 3} S^{f}_{m,\omega}$ is then the group of all finitary switchings.  The action of such groups is not transitive for by Lemmas~\ref{eqconsw} and~\ref{lmthmo} we cannot switch the complete red-coloured graph to the complete blue-coloured graph. 

\section{Basic Properties}

Our aim in this first section is to uncover results on the way to finding the structure of the \emph{full (or unrestricted) switching group},
\index{group ! switching ! full}%
 $S_{\omega,\omega}$ of all possible switchings, to be defined below.  First we pr\'ecis our findings on switching groups so far.  

The finite switching groups can be extended by defining $S^*_{m,n} :=
S_{m,n} \sd \Sym(n)$
\index{group ! switching ! extended}%
which is a group generated by the switchings and the vertex permutations.  It is
well defined since $\sigma_{c,d,X}^g=\sigma_{c,d,X^g}$ for $g \in \Sym(n)$ and so
$\Sym(n)$ normalizes $S_{m,n}$.  For any coloured graph $\Gamma$, the
group $\SAut(\Gamma)$ of switching automorphisms of
$\Gamma$ is the stabilizer of $\Gamma$ in $S^*_{m,n}$ -- equivalently, it is
the group of vertex permutations $g$ such that $\Gamma g=\Gamma\sigma$ for some
switching $\sigma$.  For $m=2$, $\SAut(\Gamma)$ is isomorphic to the
group of automorphisms of the two-graph
\index{graph ! two-graph}%
\index{two-graph}%
 derived from $\Gamma$, and so cannot be more than $2$-transitive
\index{group ! permutation ! $2$-transitive}%
 unless $\Gamma$ is trivial; only if $\Gamma$ is trivial, that is has just a single colour, can its switching group be $3$-transitive. On the other hand, if $m$ is finite but greater than $2$, then $\SAut(\Gamma)$ can be highly transitive, as it is for example
for $\Gamma = \mathfrak{R}_{m,\omega}$.  Whilst $S^*_{m,n}$ is primitive
\index{group ! permutation ! primitive}%
 on $\mathcal{G}_{m,n}$ if $m>2$
and $n < \omega$, the group $S^*_{m, \omega}$ for $m>2$ is not even transitive on the set of
all $m$-coloured graphs, as illustrated by the monochromatic graph example above.

The semidirect product structure of the finite switching groups derived in Theorem~\ref{swgpfm} is
retained when we take a countably infinite vertex set, so that
$S_{m,\omega} = N \sd K$, where $N$ is the normal subgroup of all elements of
$S_{m, \omega}$ which induce even
permutations of the colours on any edge, and $K$ is an infinite
product of copies of $C_2$.  The group $N$ lies strictly between the
direct and cartesian products of alternating groups of degree $m$ for
$m \ge 3$ colours and $K$ is an elementary abelian $2$-group
\index{group ! elementary abelian}%
isomorphic to the switching group $S_{2,\omega}$.  

\begin{figure}[!h]
$$\xymatrix{ 
{\mathcal{C}_m} & {}  & {} & {}  & {} & {}  & {}\\
& {\bullet~\bullet} \ar@{-}[u]  & {\bullet~\bullet} \ar@{-}[u] & {\bullet~\bullet}
  \ar@{-}[u] & {\bullet~\bullet} \ar@{-}[u] & {\bullet~\bullet} \ar@{-}[u] & {\bullet~\bullet}\ar@{-}[u] & {E}
}$$\caption{Fibre diagram for $n$--vertex $m$--coloured graph}
\label{fibdi}  
\end{figure}
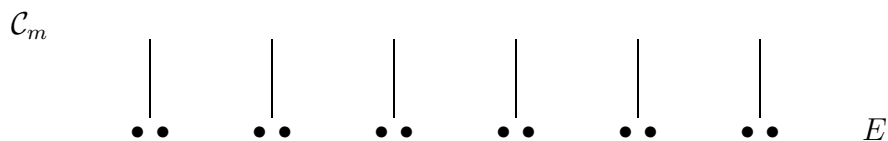

Now look at an extended version of the finitary switching groups.  We make similar constructions to that of the wreath products
\index{group ! permutation ! imprimitive}%
 of two permutation groups~\cite{cam1}.  Let $\Sym(V)$, be the `top group'
\index{group ! top}%
acting in its induced action on the set
$E$ of two-element subsets of the graph vertex set $V$ ($2 \le
|E| = {n \choose 2} \le \infty$ for an $n$-vertex set), and $H$ the `bottom group'
\index{group ! bottom}%
 (given by transpositions when it is finite), acting on the colour set $\mathcal{C}_m$ ($2
\le |\mathcal{C}_m| = m \le \infty$); that is, $H$ denotes a switching
group.  The operand is shown in Figure~\ref{fibdi}.

So $\Omega = \mathcal{C}_m \times E$ is a fibre bundle
\index{fibre bundle}%
over the set of edges $E$, with projection map $(c, e) \mapsto e$; each fibre takes
the form $\{(c, e) : c \in \mathcal{C}_m \}$ for fixed $e \in E$~\label{edgeofgraph} and represents a colour set on an edge.  Each of the heavy pairs of points in the diagram represents an edge, a point in a fibre is
a coloured edge and a global section
\index{global section}%
(a subset containing one point from each fibre) corresponds to a colouring of the graph.  
 \smallskip

We would like to go one step further and define a concept of unrestricted switchings for the countably infinite
random graph
\index{graph ! random}%
 $\mathfrak{R}_{\omega,\omega}$ with a countable infinity
of colours,  and we denote the group so generated by
$S_{\omega,\omega}$.  Note that $\bigcup_{m \geq 3} S^{f}_{m,\omega} < S_{\omega,\omega}$.  Let $X \subseteq V$ where $V$ is the vertex set and $\pi$~\label{piperm} be a permutation of $\Sym(m)$
with $|m| = \omega$.  For subgraph $\Gamma$ let the switching action be
$\sigma_{\pi, X} (\Gamma) = \Gamma'$ where the colour of edge $e$ in
graph $\Gamma'$ is

\begin{displaymath}
\mathcal{C}_{\Gamma'} (e) := \left\{ \begin{array}{ll} \mathcal{C}_{\Gamma'} (e)\pi & \textrm{if edge}\ e\ \textrm{goes between}\ X\
\textrm{and}\ X'\\
\mathcal{C}_{\Gamma} (e) & \textrm{otherwise,} \end{array} \right.
\end{displaymath}
where $X'$ is the complement of $X$ in $\Gamma$.  It can be checked straightforwardly that $\bigcup_{m \geq 3} S^{f}_{m,\omega}$ is obtained by restricting $\pi$ to be a finitary permutation.

Then define $S_{\omega,\omega} := \langle
\sigma_{X, \pi} | X \subseteq V\ \text{and}\ \pi \in \Sym(\mathcal{C}_m), m \subseteq \omega \rangle$.  In terms of the
above fibre diagram, those fibres representing edges both within and
outside $X$ are stabilized by the switching action, whilst for edges
between $X$ and its complement $X'$, $c \to c \pi$.  We have here a form
of subcartesian product of copies of $\Sym(\mathcal{C}_m)$ indexed by edges or
more properly an extended direct product which in a suitably defined
limit of action over all fibres or edges gives a cartesian product.  So
our product lies strictly between a direct and a cartesian product.  A
direct product is insufficient because a finite number of switchings
is insufficient to describe the switching operation on an infinite
colour set.  Nor can it be a cartesian product which gives a transitive group action on the relevant
function set where $S_{\omega, \omega} \le
\Cr_{E} \Sym(\mathcal{C}_{\omega})$~\label{Crprod} with the product action, if we identify a
coloured graph with a function $E \to \mathcal{C}_{\omega}$; this would
contradict the strict inequality $\SA(\Gamma) < \Sym(X)$ for a graph $\Gamma$ on a vertex set $X$,
from Theorem~\ref{sanesm}.  This new type of group product is worthy
of study in its own right in order to uncover its
properties.

\smallskip

Consider $S_{\omega, \omega} \sd \Sym(V)$ where $|V| = \omega$.  Using this group we can set up a wreath-like product under two different actions that are similar (but not identical) to the imprimitive
\index{group ! permutation ! imprimitive}%
 and product actions of a wreath product: 

(i) an imprimitive action 
\index{group ! permutation ! imprimitive}%
on the set $\Omega = \mathcal{C}_m \times E$ of coloured edges;  these sets are orbits of $S_{\omega, \omega}$ and blocks of imprimitivity of $S_{\omega, \omega} \sd \Sym(V)$, and the orbits are transitive and oligomorphic (see later in the chapter for the definition of a `parity' that provides the relevant countable first-order structure), and
\index{group ! oligomorphic}%
\index{first-order structure}%
(ii) a product action 
\index{group ! permutation ! product}%
on the set of all coloured graphs, which can be identified as the set of all
functions from $E$ to $\mathcal{C}_m$, or of global sections
\index{global section}%
of the fibred space $\mathcal{C}_m \times E \to E$.  This action is not oligomorphic since the number of edges on which two colourings agree is an invariant of pairs.  (The group $\SA(\Gamma)$ is the stabilizer of graph $\Gamma$ in the wreath product under the product action with $\Gamma$
identified with a global section of $\Omega = \mathcal{C}_m \times E$).  

To describe the imprimitive complete wreath product action 
\index{group ! permutation ! wreath product}%
let the cartesian
product $B$ of $|E|$ copies of the group $H$, one acting on
each fibre as $H$ acts on $\mathcal{C}_m$, be the set of functions from
$E$ to $H$ acting pointwise;  the action is given by $\mu
( (c, e), h) = ( c h(e), e )$, for $h \in
H$, with function $\mu : \Omega \times H \mapsto \Omega$.  For $|\mathcal{C}_m|, |E| > 1$, the fibres are
blocks of imprimitivity with congruence
 \index{congruence}%
  relation $\sim$ defined by
$(c, e) \sim (c', e')$\ $\emph{iff}$\
$e = e'$.  Now let $T$ be a copy of the group $\Sym(V)$ permuting the fibres: 
 $\mu ((\gamma, e), k) = (\gamma, e k)$ for $k \in
 \Sym(V)$.  Then $T$ normalizes $B$, and the semidirect product $G=B \sd T$
 is a wreath product of permutation groups.  We are interested in a
 subgroup of this wreath product where the full cartesian
 product of $H$ groups in $B$ is replaced by the switching group
 acting on only a proper subset of the fibres, although this subset
 can still be countably infinite.
 
The group $\SA\left(\mathfrak{R}_{m,n}\right) \sd \Sym(V)$ where $m \ge 3$ and $|V| = \omega$ is transitive on coloured edges.

Returning to finitary permutations, assuming the axiom of choice
\index{axiom of choice}%
define $\mathcal{S}_{\omega} := \cup_{m \ge 3}S^{f}_{m, \omega}$ whose elements induce finitary permutations on the colours of any edge, so $\mathcal{S}_{\omega}$ also has a semidirect product structure; a union of an ascending sequence of groups is always a group.  Another instructive way of writing $\mathcal{S}_{\omega}$ would be as $S^{f}_{\omega, \omega}$.

\begin{proposition}
$S_{\omega,\omega} \gneqq \mathcal{S}_{\omega}.$
\end{proposition}

This is clear because elements of the right-hand group can only change finitely many colours whereas the other group changes an arbitrarily large number of colours.

Let $m$ be finite.  Define another variation on the switching group theme, the \emph{enhanced switching group} by $ES_{m, n} :=
\Cr^{{\omega}}_{i=1} (S_{m, n})_{i=1} \sd \FSym(V)$,
\index{group ! switching ! enhanced}%
where $|V| = \omega$.

\begin{theorem}
 $ES_{m, n}$ is locally finite.
\index{group ! locally finite}%
\end{theorem}

\begin{proof}
This follows from the next two lemmas.
\end{proof}

\begin{lemma}
The cartesian power of a fixed finite group is locally finite.
\end{lemma}

\begin{proof}
Let $I$ be an infinite index set. Define $\Cr_{i \in I} (H_{i}) := \{
\emph{f}: I \longrightarrow \bigcup_{i \in I} H_{i}\ |\ \emph{f}(i)
\in H_{i}\ \forall i \in I\}$.  
If $\emph{f}_{1}$, $\emph{f}_{2}$ are two elements of the cartesian
power, the group operation is multiplication of components defined
by $(\emph{f}_{1} \emph{f}_{2})(i) = \emph{f}_{1}(i)\emph{f}_{2}(i)$,
where the right-hand side is a product in $H_{i}$.  Take $H_{i} = H$, a
fixed finite group.  Take $n$ elements of the cartesian product
$\emph{f}_{1}, \ldots, \emph{f}_{n}$.  There exists a partition $I =
I_{1} \cup I_{2} \ldots \cup I_{N}$ such that all $\emph{f}_{i}$ are constant on the components in each $I_{j}$.  We need to show that there is an embedding $
\langle \emph{f}_{1}, \ldots, \emph{f}_{n} \rangle \hookrightarrow
H^{N}$.  Define a map $\phi_i : (\emph{f}_{1}, \ldots, \emph{f}_{n}) \to
(\emph{f}_{1}(i), \ldots, \emph{f}_{n}(i))$ where $\phi_i = \phi_{i'}$ for $i, i' \in I_j$ $(1 \le j \le N)$.
Then this is an injection from $\langle \emph{f}_{1}, \ldots,
\emph{f}_{n} \rangle$ into $H^{N}$ which is a finite group.  Therefore
$\langle \emph{f}_{1}, \ldots, \emph{f}_{n} \rangle$ is a finite group.
\end{proof}

So the group $S_{m, \omega} = \langle \sigma_{c, d, X} | X \subseteq V,\ c, d \in \mathcal{C}_m \rangle (\le \Cr^{{\omega}}_{i=1} (\Sym(m))_{i})$ is locally finite.
\index{group ! locally finite}%
\begin{lemma}
\label{lgpitm}
The class of locally finite groups is extension-closed.
\index{group ! extension-closed}%
\end{lemma}

\begin{proof}
Let $N \vartriangleleft G$ and suppose that $N, G / N$ are locally finite.  Take $g_1, \ldots, g_n \in G$.  Their images under the natural homomorphism $G \to G / N$ generate a finite subgroup.  Let $H = \langle g_1, \ldots, g_n \rangle$.  The Reidemeister-Schreier algorithm
\index{Reidemeister-Schreier algorithm}%
\index{Schreier, O.}%
\index{Reidemeister, K}%
then gives us a finite generating set for $H \cap N$, so $H \cap N$ is finite.  Also $N / H \cap N \cong HN / N$ is finite by the first step.  So $H$ is finite.
\end{proof}

It then follows that,

\begin{theorem}
\label{lgpitm}
$\mathcal{S}_{\omega} \sd \FSym(V)$ where $|V| = \omega$ is locally finite.
\index{group ! locally finite}%
\end{theorem}

\begin{proof}
The class of locally finite groups is subgroup-closed.
\end{proof}

Note however that $S_{\omega,\omega}$ is not locally finite because it contains subgroups that are the product of an infinite number of elements of increasingly high order.  


We saw that the full group of switching automorphisms defines a reduct if $m=2$ but not if $m>2$. However we also have a proof that switching a specific pair of colours defines a 
reduct for $m>2$. It is probably true that the group of switching automorphisms
involving a specific proper subset of the colours always defines a reduct.

\section{Parity Equivalence of Coloured Graphs}
\label{profilsec}

In this section we continue to uncover further properties of switching groups with
finitely many colours.  It will transpire that the finitary switching groups are closed
\index{group ! closed}%
 in a sense that we will elucidate.

We differentiate between the group $S_{m, \omega} = \langle
\sigma_{c,d, X} : X \subseteq V, |V|=\omega, c,d \in \mathcal{C}_{m}
\rangle$ that is generated by switchings, and the group that is its closure in its action on the space $\mathcal{C}_m \times {\omega \choose 2}$ of coloured edges, which we denote $S_{m, \omega}^{cl}$.~\label{clswitchinggroup}
\index{group ! switching ! closure of}%
We will give a description of $S_{m, \omega}^{cl}$ as a closed permutation group
\index{group ! permutation ! closed}%
 whereby its elements act so as to preserve a parity condition.

First we observe a topological property of those switching groups $S_{m, n}$ for finite $m$ in the $n \to \omega$ limit -- that $S_{m, \omega}^{cl}$ is a
\emph{profinite group},
\index{group ! profinite}%
that is a compact totally disconnected topological group.
\index{topological group ! totally disconnected}%
\index{group ! topological ! totally disconnected}%
A metric space is \emph{connected}
\index{metric space ! connected}%
if it cannot be expressed as a union of two disjoint, nonempty subsets that are open in its induced topology.  A space is \emph{totally disconnected}
\index{space ! totally disconnected}%
if each of its components contains only a single point, which means that it has separation property 
\index{separation property}%
$T_1$,
\index{topological space ! $T_1$}%
 that is, every point is closed; more on this concept in Chapter~\ref{chap8}.  The groups in question are compact, being both closed and having only finite orbits~\cite{cam4} in an imprimitive
\index{group ! permutation ! imprimitive}%
  action to be defined below, and so can be interpreted as profinite groups, that is inverse limits of inverse systems of finite groups.  For each colour
set the inverse system is the sequence
\[ \ldots \to S_{m, n} \to S_{m, n-1} \to \ldots \to S_{m, \bar{n}} \] 
together with epimorphisms $\phi_{n,\bar{n}} :
S_{m,n} \to
S_{m,\bar{n}}$ (for $\bar{n} \le n$)
such that $\phi_{a,b}\phi_{b,c}=\phi_{a,c}$ whenever $c \le b \le
a$.  The maps are restriction homomorphisms regarding colourings of $K_{\bar{n}}$ as restrictions of colourings of  $K_n$.

The \emph{inverse\ limit} of the sequence is a universal topological group
\index{topological group ! universal}%
$G$ whose elements are all sequences $(g_n : n \in
\omega)$ with $g_{n} \in S_{m,n}$ such that $g_{n}\phi_{n,\bar{n}} =
g_{\bar{n}}$ for $\bar{n} \le n$, with group composition being
coordinatewise multiplication together with the homomorphisms $\theta_{n}$ such that
  $\theta_{n}\phi_{n,\bar{n}} = \theta_{\bar{n}}$, that is
$\theta_{n}$ acts as projection onto the $n$th coordinate, $\theta_{n} : G \to
  S_{m,n}$.  So $G$ induces the
finite groups $S_{m,n}$.  It is possible that the group
$S_{m, \omega}^{cl}$ also has non-profinite
proper subgroups but we do not concider these.  

All profinite groups are homomorphic images of some closed oligomorphic permutation group
\index{group ! permutation ! oligomorphic}%
with some prescribed degree of transitivity~\cite{evhe}.  There are several equivalences~\cite{waterhouse}~\cite[p.~19]{wilson} to a topological group
\index{topological group}%
being profinite, two of which are:

(i) $S_{m, \omega}^{cl}$ is isomorphic (as a topological group) to a closed subgroup
\index{group ! closed}%
 of a Cartesian product of finite groups. 

We can illustrate how the group $S_{m, \omega}^{cl}$ is closed in a countably infinite
cartesian product of finite symmetric groups, $S_{m, \omega}^{cl} \le
\Cr^{{\omega}}_{i=1} (\Sym(m))_{i}$.  Observing that the closure of the union of direct products $\Dr_{i < \infty} (G)_i$~\label{Drprod} is equal to the cartesian product $\Cr_{j = 1}^{\infty} (G)_j$ for a finite group $G$, we define $N := S_{m, \omega}^{cl} \cap \Cr^{{\omega}}_{i=1} (\Alt(m))_{i}$.  Clearly $N \vartriangleleft
S_{m, \omega}^{cl}$ and the quotient $S_{m, \omega}^{cl} / N$ is isomorphic to the
elementary abelian two-group $C_2^{\omega}$.  Now there exists a subgroup $S'$ of
index $2$ in $S_{m, \omega}^{cl}$ containing $N$ such that $S'$ induces
$\Sym(m)$ on each edge or fibre.  Note that there are uncountably many
subgroups of index $2$, but only countably many fail to induce $\Sym(m)$ on
some fibre.  We can choose $S'$ to contain the stabilizer  $S_{m,
\omega, (e_1, \ldots, e_r)}^{cl}$ of the
colours of a finite
tuple of edges $e_1, \ldots, e_r$.  Project onto the subgroup of
$S_{m, \omega}^{cl} / N$ which is trivial in coordinates $1, \ldots, r$, so that
$S' \ge S_{m, \omega, (e_1, \ldots, e_r)}^{cl}$.  So $S'$ has index $2$ in
$S_{m, \omega}^{cl}$ and is open
\index{group ! open}%
because it contains the pointwise stabilizer of a finite set.

(ii) $S_{m, \omega}^{cl}$
\index{group ! closed}%
 is compact and $\bigcap (N | N \vartriangleleft_O S_{m, \omega}^{cl}) = 1$, 
where $\vartriangleleft_O$ denotes an open normal subgroup.
\index{group ! open normal}%

We can illustrate this by taking for $N$ a sequence of fibrewise stabilizers of an increasing number of fibres, whose intersection is then trivial.

\bigskip

There are also various subgroup properties~\cite[p.~36]{wilson} such as,

(i) $S_{m, \omega}^{cl}$ has Sylow $p$-subgroups for some prime $p$,
\index{group ! Sylow subgroup}%

(ii)  Any two Sylow $p$-subgroups are conjugate in
$S_{m, \omega}^{cl}$, that is $\exists g \in S_{m, \omega}^{cl}$ such that $g^{-1} P_1 g =
P_2$.

There is a second equivalent description of $S_{m,\omega}^{cl}$, which explains the connection between profinite groups and groups with finite orbits via an equivalence of topologies.  

The switching group $S_{m,\omega}^{cl}$ is a closed subgroup of the cartesian
product of symmetric groups of degree $m \ge 3$, in the topology defined by
the imprimitive action
\index{group ! permutation ! imprimitive}%
 for the cartesian product where each
fibre comprises the colour set $\mathcal{C}_m = \{1, \ldots, m\}$ of an edge.  It is possible to firstly show that
on each fibre we have a convergent sequence of switchings using the
idea of pointwise convergence: the
distance $d(f, g)$ between two distinct automorphisms $f$ and $g$,
\[d(f, g) := 2^{-(min (n): f(n) \neq g(n)\ or\ f(n)^{-1} \neq g(n)^{-1})}\] 
defines the topology of pointwise convergence.
\index{topology ! of pointwise convergence}%
Then we can show that there is a convergence on a
cartesian product of the fibres.  By identifying the topology of pointwise
convergence with the profinite topology we can prove convergence of
products of switchings on the cartesian product of fibres.  

A profinite group having a countable base of open subgroups
\index{group ! open}%
(or equivalently a closed subgroup of $\Sym(\omega)$ with all orbits finite), is topologically isomorphic to the factor group of an automorphism group
of an $\aleph_{0}$-categorical structure
\index{aleph@$\aleph_0$-categorical}%
by a closed normal subgroup
$\Phi$~\cite[Lemma~3.1]{evhe}.  Furthermore, $\Phi$ can be chosen to be oligomorphic
\index{group ! oligomorphic}%
and
having no proper closed subgroups of finite index.  

Recall that the topology of pointwise convergence
\index{topology ! of pointwise convergence}%
refers to the limit of a sequence of
permutations on operand $X$ converging thus: $\lim_{n \to \infty} (g_n) = g$ if and
only if $\forall x_i \in X,\ \exists n_0 \in \omega$ such that $\forall
n > n_0,\ x_i g_n = x_i g$.  A basis
of open neighbourhoods of the identity consists of the pointwise
stabilizers of all finite sets.  But $S_{m,n} < \Sym(m)$ and the
topology of pointwise convergence on $\Sym(m)$, where $\mathcal{C}_m$ has the
discrete topology,
\index{topology ! discrete}%
is Hausdorff,
\index{Hausdorff space}%
so we have closure and convergence of
switchings on each fibre.  Now we identify the topology of pointwise
convergence (for some $n_0$ as above) with the profinite topology of
product switchings on the $n_0$ vertices.  This is done by taking
the inverse system with all epimorphisms $\phi_{n_0,n} : S_{m,n_0} \to
S_{m,n}$ $(\forall\ 1 \le n < n_0)$ such that $\phi_{n_0,n_1}
\phi_{n_1,n_2} = \phi_{n_0,n_2}$ whenever $n_2 \le n_1 \le n_0$.  

We omit the details but the conclusion is that there is convergence in the profinite topology
\index{topology ! discrete}%
 of product switchings induced from the discrete topology
\index{topology ! discrete}%
on each $S_{m,n}$.  The switching group $S_{m,\omega}^{cl}$ is an elementary type of profinite group, that is one contained in a cartesian product of finite groups.

Automorphism groups
\index{group ! automorphism}%
 of countably infinite structures $\mathcal{M}$ are Polish groups~\cite{cam10},
\index{group ! Polish}%
that is, complete separable metric spaces whose
underlying topology turns the group into a topological group.
\index{group ! topological}%
  The topology of pointwise convergence
\index{topology ! of pointwise convergence}%
turns $\Aut(\mathcal{M})$ into a complete metric space.
\index{metric space ! complete}%
The group $S_{m, \omega}^{cl}$ is a Polish group because it is compact and so complete.

A third equivalent description of $S_{m, \omega}^{cl}$ is as a group with commutator subgroup $\Cr_{i = 1}^{\infty} (\Alt(m))_i$.

\bigskip

We now derive a parity condition which acts as an invariant relation for the action of $S_{m, \omega}^{cl}$
\index{group ! closed}%
 on $3$-vertex subsets of $\mathfrak{R}_{m,\omega}$. 

\begin{figure}[!h]
$$\xymatrix{ 
& {\bullet~\bullet} \ar@{}[r] & {\bullet~\bullet} \ar@{}[r] & {\bullet~\bullet}\\
& {\bullet~\bullet} \ar@{-}[u] \ar@{.}[d] & {\bullet~\bullet} \ar@{-}[u] \ar@{.}[d] &
{\bullet~\bullet} \ar@{-}[u] \ar@{.}[d]\\
& {\bullet~\bullet} \ar@{}[r] & {\bullet~\bullet} \ar@{}[r] & {\bullet~\bullet}\\
& {\quad \bullet~\bullet}_{xy} \ar@{-}[u]  & {\quad \bullet~\bullet}_{xz} \ar@{-}[u] & {\quad \bullet~\bullet}_{yz} \ar@{-}[u]
}$$\caption{Fibre diagram for $\mathfrak{R}_{m,3}$}
\label{fdfthr}
\end{figure}
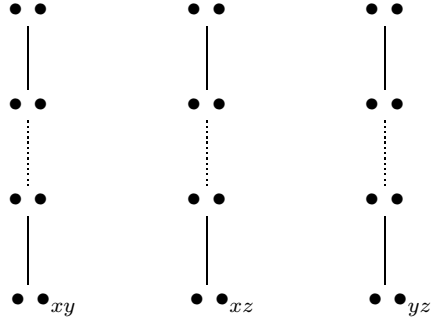

 Consider three vertices $x, y, z$ in a \emph{complete multicoloured graph} and $m$ possible colours
on the corresponding edges $xy, xz, yz$, so that we have an $m$-layered
fibre diagram as in Figure~\ref{fdfthr}.
Any switching changes the colours of an even number of edges.  It
also preserves an analogue of a two-graph,
\index{graph ! two-graph}%
\index{two-graph}%
 namely, the parity of the permutation of the colours on any $3$-set of vertices with edges
coloured from $\mathcal{C}_{m} = \{c_1, \ldots, c_m\}$.  Let $\mathsf{P}_m$ be a
$3m$-place relation which holds only if there exist vertex triples $x, y, z$ such
that the first $m$ arguments form the fibre above $xy$ (in some
order), the second $m$ form the fibre above $xz$, and
the last $m$ form the fibre above $yz$, and if $g_{xy},
g_{xz}$, and $g_{yz}$ are the permutations describing the order of
these arguments, then $g_{xy} g_{xz} g_{yz}$ is an even
permutation.  The element $g_{xy}$ maps the colour $c_i$ on edge $\{x, y\}$ to the $i$th
argument of $\mathsf{P}_m$.  For example taking $m=3$, $c_{12}$ to mean the colour on edges $\{1, 2\}$ and $(c_{12}, c_{13}, c_{23})
= (\mathfrak{r} , \mathfrak{b}
, \mathfrak{g})$, the permutation $(\mathfrak{r}, \mathfrak{b}, \mathfrak{g}) \to
(\mathfrak{g}, \mathfrak{b}, \mathfrak{r}) \to (\mathfrak{b}, \mathfrak{g},
\mathfrak{r}) \in \mathsf{P_m}$.  Notice that the definition of parity depends on the sequential composition of the individual colour changes on the edges, for if we were instead to take the following definition of composition: $(\mathfrak{r}, \mathfrak{b}, \mathfrak{g}) \to (\mathfrak{g}, \mathfrak{b}, \mathfrak{r})$ which is an
odd permutation followed by $(\mathfrak{r}, \mathfrak{b}, \mathfrak{g}) \to (\mathfrak{b}, \mathfrak{g}, \mathfrak{r})$ which is an even permutation, then composing these parities multiplicatively gives a product of permutations having odd parity.

We call two complete $m$-coloured graphs $\Gamma_1$ and $\Gamma_2$ on any vertex
set \emph{$\mathsf{P}_m$}-\emph{equivalent}~\label{mathsf{P}_m} 
\index{graph ! Pmequ@$\mathsf{P}_m$-equivalent}%
if each triple of vertices of $\Gamma_2$ is a parity preserving permutation of the colours of the corresponding vertex triple of $\Gamma_1$.  

\begin{lemma}
\label{rmlem}
$\mathsf{P}_m$ is preserved by any switching $\sigma_{c, d, X}$.
\end{lemma}

\begin{proof} 
Obvious if $X \cap \{x, y, z \} = \emptyset$ or $\{x, y, z\}$.  If $X
\cap \{x, y, z \} = \{x\}$ then $\sigma_{c, d, X}$ transposes colours
$c$ and $d$ on
edges $xy$ and $xz$ and fixes colours on $yz$.  So the parity of the
product is unchanged. 
\end{proof}

\begin{theorem}
\label{parthm}
The elements of $S_{m, \omega}^{cl}$ are precisely the permutations which fix all the fibres and preserve $\mathsf{P}_m$.
\end{theorem}

\begin{proof} 
Switching does not change the $\mathsf{P}_m$ parity of the edges within a $3$-set and as this is true for all $3$-sets so switching-equivalent graphs are $\mathsf{P}_m$-equivalent.  

For the converse we need to show that a $\mathsf{P}_m$-preserving permutation of edge colours is always achievable by a product of elements of $S_{m, \omega}^{cl}$.  We must first consider the vertex-permutation induced by a $\mathsf{P}_m$-preserving permutation of the whole graph and then show that this is achievable as a product of switching automorphisms.  

But as we have seen in the previous chapter, the group of switching automorphisms is highly transitive so its closure is the symmetric group.  So we can always undo a parity-preserving permutation by an element within $S_{m, \omega}^{cl}$, modulo a permutation of the elements of a fibre on each edge.  This can then be done one fibre at a time and then the limit can be taken, which lies within the wreath product $\lim_{X \to \infty} \Sym(\mathcal{C}_{m}) \Wr \Sym(X)$.
\index{wreath product}%

Now we deal with the permutation of colours/fibres on each edge.  There are $m !$ such permutations.  On a particular edge, we can permute the colours so as to preserve $\mathsf{P}_m$.  Working on more and more edges, we can switch to ensure that the fibres above an edge are left invariant after a sufficient point.  This will apply to an increasing number of edges.  Then pointwise convergence
\index{topology ! of pointwise convergence}%
 makes the limiting process work.
\end{proof}

It now follows that the group $S_{m, \omega}^{cl}$ is the automorphism group
\index{group ! automorphism}%
 of an equivalence relation whose classes are fibres on an edge together with the first-order $3m$-ary relation $\mathsf{P}_m$.  In the next chapter we give a slight extension of this result from complete graphs to more general graphs.

\medskip

To recap, we have introduced two different notions of parity of colours on random graph edges, but we note that a comparison of these is not straightforward as they seem to be somewhat different entities.

Firstly in Chapter~\ref{chap2} we defined two complete $3$-coloured graphs to be
\emph{P-equivalent}~\label{P-equivalent}
\index{graph ! P-equivalent}%
 if they have the same red edges and each blue-green triangle has the same parity of the number of green edges.  This is a \emph{global} concept used to compare two graphs, or two global sections of a fibre diagram, where each fibre represents the different colours on a specific edge.  

Then in this section, we defined two complete $m$-coloured graphs on any vertex set to be \emph{$\mathsf{P}_m$}-\emph{equivalent}~\label{mathsf{P}_m} 
\index{graph ! Pmequ@$\mathsf{P}_m$-equivalent}%
if each triple of vertices of one graph is a parity preserving permutation of the colours of the corresponding vertex triple of the other graph.  Consideration of triples of vertices at a time makes this a \emph{local} concept.

From our discussion in Chapter~\ref{chap2}, P-equivalence is a weaker condition that S-equivalence in general, and by Theorem~\ref{psequiv}, as the vertex set $n\to\infty$, the two become asymptotically equivalent in an appropriate set of measure 1, and they coincide for $\mathfrak{R^{t}}$.

This, together with Theorem~\ref{parthm}, appear to indicate that $S_{m, \omega}^{cl}$ is strictly an overgroup of $S_{m, \omega}$, and that the example we gave in our discussion of P-equivalence in Chapter~\ref{chap2} belongs to $S_{m, \omega}^{cl} \backslash S_{m, \omega}$.  Any switching which stabilises the blue-green parity (in triangles) but changes the red edges is another example.  

\head{Open Question}  Find an example of a graph for which $S_{m, \omega}$ is not closed.

\vspace{15pt}

We began the section by showing how to form an inverse limit of switching groups, and we end by showing how they can form a direct limit.
\index{direct limit}%
  Fix the number of colours $|\mathcal{C}_{m}| = m$.  Let $\mathsf{S}_{n}$~\label{mathsf{S}_{n}} be the group of switchings on $\mathfrak{R}_{m, n}$; any element of this group is the product of elementary switchings about a single vertex which we can denote $\sigma_{c, d, v}$.  We can extend $\mathsf{S}_{n}$ to act on $\mathfrak{R}_{m, n+1}$.  Denote by $\sigma^{(n)}_{c, d, v}$ the switching of colours $c$ and $d$ about $v \in \{1, \ldots, n\}$.  The map $\alpha^{n+1}_{n} : \sigma^{(n)}_{c, d, v} \mapsto \sigma^{(n+1)}_{c, d, v}$ is an injection from $\mathsf{S}_{n}$ to $\mathsf{S}_{n+1}$.  The map $\pi^{n+1}_{n} : \mathsf{S}_{n+1} \to \mathsf{S}_{n} : \sigma^{(n+1)}_{c, d, v} \mapsto \sigma^{(n)}_{c, d, v}$ is an epimorphism since we can define
\begin{displaymath}
\sigma^{(n+1)}_{c, d, v} \in \mathsf{S}_{n+1} \to \left\{ \begin{array}{ll}
\sigma^{(n)}_{c, d, v} \in \mathsf{S}_{n} & \text{if}\ v \in \{1, \ldots, n\}\\
1 & \text{otherwise}
 \end{array} \right.
\end{displaymath}
and
\begin{displaymath}
\sigma^{(n+1)}_{c, d, v} \sigma^{(n+1)}_{k, l, w} \in \mathsf{S}_{n+1} \to \left\{ \begin{array}{ll}
\sigma^{(n)}_{c, d, v} \sigma^{(n)}_{k, l, w} \in \mathsf{S}_{n} & \text{if}\ v, w \in \{1, \ldots, n\}\\
1 & \text{if}\ v, w \notin \{1, \ldots, n\}\\
\sigma^{(n)}_{c, d, v} \in \mathsf{S}_{n} & \text{if}\ v \in \{1, \ldots, n\},\\
& \text{and}\ w \notin \{1, \ldots, n\}.
 \end{array} \right.
\end{displaymath} Therefore there are homomorphisms both ways
$$
\xymatrix{{\mathsf{S}_{n}} \ar@/^/[rr]^{\iota_{n+1}} && {\mathsf{S}_{n+1}} \ar@/^/[ll]^{\pi_n}}$$
where $\ker(\pi_{n}) \vartriangleleft \mathsf{S}_{n+1}$ is a normal complement to $\im(\iota_{n+1})$,~\label{image} and also $\ker(\pi_{n}) \cap \im(\iota_{n+1}) = 1$.  The group $\mathsf{S}_{n}$ acts on the points above edges in $\{1, \ldots, n\}$, that is on coloured edges from this set.  So $\mathsf{S}_{n+1} \to \mathsf{S}_{n}$ is the restriction map, and its kernel consists of switchings fixing all coloured edges within $\{1, \ldots, n\}$.  The image map permutes the colours on the edges beween the $(n+1)$st  vertex and the other $n$ vertices.

Now we can define a direct system of groups comprising a family $(\mathsf{S}_{n})_{n \in \mathbb{N}}$ and homomorphisms $\alpha^{n_2}_{n_1} : \mathsf{S}_{n_1} \mapsto \mathsf{S}_{n_2}$ where $n_1 \le n_2$, satisfying (i) $\alpha_{n_1}^{n_1}$ is the identity on $\mathsf{S}_{n_1}$, (ii) $\alpha^{n_3}_{n_1} = \alpha^{n_2}_{n_1} \alpha^{n_3}_{n_2}$ whenever $n_1 \le n_2 \le n_3$.  We can define an equivalence relation $\sim$ on the set-theoretic union $\mathsf{S} = \bigcup_{n \in \mathbb{N}} \mathsf{S}_n$, defining $\sigma_{n_1} \sim \bar{\sigma}_{n_2}$ if there exists $n_3 \ge n_1, n_2$ such that  $\alpha^{n_3}_{n_1} (\sigma_{n_1}) = \alpha^{n_3}_{n_2} (\bar{\sigma}_{n_2})$, that is the images are equal.  Let $\bar{\mathsf{S}}_n$ be the image of the monomorphism $\theta_n: \mathsf{S}_n \to \mathsf{S}_D$ in $\mathsf{S}_D = \lim_{n \to \infty} \mathsf{S}_n$.  Now convert the set $\mathcal{S}_D$ of all equivalence classes $[\sigma_n] \in \bar{\mathsf{S}}_n$ of $\sigma_n \in \mathsf{S}_n$ into a group by defining the product $[\sigma_{n_1}][\sigma_{n_2}] = [\alpha^{n_3}_{n_1} (\sigma_{n_1}) \alpha^{n_3}_{n_2} (\sigma_{n_2})] = [\alpha^{n_3}_{n_1} (\bar{\sigma}_{n_1}) \alpha^{n_3}_{n_2} (\bar{\sigma}_{n_2})]$ and $1_D = [1_{\mathsf{S}_n}]$,  and $[\sigma_n]^{-1} = [\sigma_n^{-1}]$.  Then $\{\mathsf{S}_D,\ \theta_n\}$ is the \emph{direct limit} of the direct system.
\index{direct limit}%

It is not always the case that inverse and direct limits of homomorphisms accompany each other.  For example, whilst there is a direct limit of symmetric groups, this being the finitary symmetric group, there are no inverse homomorphisms. Where they do both exist,  as is the case for the switching group limits as defined in this section,   the direct colimit and inverse limit are categorically dual to each other, though this case is perhaps not so interesting from a category theory point of view.  

\chapter{Local Switchings \& Group Presentations}
\label{chap4}
\bigskip

All the mathematical sciences are founded on relations between physical laws and laws of numbers, so that the aim of exact science is to reduce the problems of nature to the determination of quantities by operations with numbers.
\begin{flushright}
James Clerk Maxwell, \textit{On Faraday's Lines of Force, 1856}
\end{flushright}

\bigskip

\section{Local Switching of Simple Graphs}

In this section we give a synopsis of the theory of \emph{local switching}
\index{switching ! local}%
 and its connection to signed graphs, root lattices
\index{lattice ! root}%
  and Coxeter groups as expounded by Cameron, Seidel and Tsaranov
\index{Cameron, P. J.}%
\index{Seidel, J. J.}%
\index{Tsaranov, S. V.}%
 in~\cite{camseits} for two-colour graphs.  In the next section we begin to extend this theory to multicoloured graphs.

A \emph{signed graph}
\index{graph ! signed}%
 $(\Gamma, f)$~\label{Gammaf} is a graph $\Gamma$ with a signing $f: E(\Gamma) \to \{+1, -1\}$ on the edges.  

Given a signing $f$ of the edge-set of $\Gamma$, the signing $f_X$ obtained from $f$ by reversing the sign of each edge which has one vertex in a subset $X \subseteq V(\Gamma)$ of the vertex set, then this defines on the set of signings an equivalence relation, called \emph{switching}, since 
\[ f_{\emptyset} = f,\quad\quad (f_X)_X = f,\quad\quad \text{and}\quad\quad (f_X)_Y = f_{X \symd Y}. \]

The \emph{signed switching classes}~\label{ssclass}
\index{switching ! signed class}%
 of $\Gamma$ are given by the equivalence classes $\{f_X : X \subseteq V \}$.
 
Observe firstly that the collection of ``odd'' signed cycles (those with an odd number of + signs) is invariant under switching.  Further observe that all signings of a tree are equivalent, because from any source we can inductively switch all signs into $+1$.
 
Signed graphs can be used to define root lattices
\index{lattice ! root}%
 $\mathbb{L}(\Gamma, f)$, Weyl 
\index{group ! Weyl}%
 $\W(\Gamma, f)$,~\label{weylgp} Coxeter
\index{group ! Coxeter}%
\index{Coxeter ! group}%
 $\Cox(\Gamma, f)$, and Tsaranov 
\index{group ! Tsaranov}%
 $\Ts(\Gamma)$ groups. 

The ordinary Coxeter group $\Cox(\Gamma)$ on a graph $\Gamma$ has a generator for each vertex of $\Gamma$ and product relations of order 3 or 2 corresponding respectively to adjacent or non-adjacent vertices.  Its factor group $\Cox(\Gamma, f)$, the Coxeter group of the signed graph $(\Gamma, f)$ is obtained by adding one relation $r_C$ for each odd signed cycle $C$ (whose edges carry an odd number of $+$ signs) in $(\Gamma, f)$, where for each induced cycle $r_C = \{x_1 x_2 \ldots x_{n-1} x_n x_{n-1} \ldots x_2\}^2$ of the Coxeter group $\Cox(\Gamma)$.  The term $r_C$ is called the \emph{cut element} for an induced cycle $C$ in $\Gamma$.  Note that the relator $r_C$ is unchanged (as an element in $\Cox(\Gamma)$) if we either start the cycle at a different point or if we go round the cycle in the other direction.  That is $\{x_1 x_2 \ldots x_{n-1} x_n x_{n-1} \ldots x_2\}^2 = 1$ in $\Cox(\Gamma) \Leftrightarrow \{x_2 x_3 \ldots x_n x_1 x_n \ldots x_3\}^2 = 1$  in $\Cox(\Gamma) \Leftrightarrow \{x_1 x_n \ldots x_3 x_2 x_3 \ldots x_n\}^2 = 1$ in $\Cox(\Gamma)$; to get from the first expression to the third, move $x_1$ from the front to the end and this is equivalent to reversing the cycle, since in Coxeter groups $x_i x_i^{-1} = 1$ and $(x_1 \ldots x_n) = (x_n \ldots x_1)^{-1}$.  Written in terms of the set $\mathcal{C}_f$ of induced odd cycles of $(\Gamma, f)$,
\[ \Cox(\Gamma, f) = \Cox(\Gamma) / \langle r_C : C \in \mathcal{C}_f \rangle^{\Cox(\Gamma)}. \]

The \emph{Witt representation}
\index{Witt representation}%
 of $\Cox(\Gamma)$ is defined as follows.  Let $\Gamma$ have $(1, 0)$-\emph{adjacency matrix}~\label{admx}
\index{adjacency matrix}%
$A$, and let $\mathbb{L}$ be the lattice having a basis $(e_1, \ldots, e_n)$ with Gram matrix
\index{Gram matrix}%
  $2 I - A$ (this is the \emph{intersection matrix}
\index{intersection matrix}%
 of the $-$ signing of $\Gamma$; see below).  The Weyl group $\W(\mathbb{L})$ is generated by the reflections $w_i$ in the hyperplanes perpendicular to $e_i$, for $i = 1, \ldots, n$.

That $\Cox(\Gamma, f)$ is invariant under switching, that is, the relation corresponding to an odd cycle is changed to an equivalent relation modulo $\Cox(\Gamma)$ on switching, follows because switching preserves odd cycles and graphs.

\begin{theorem}
If $f$ is any signing of the graph $\Gamma$ then the map $x_i \mapsto w_i$, $i = 1, \ldots, n$, induces a homomorphism $\Cox(\Gamma) \to \W(\Gamma, f)$ whose kernel contains the cut elements $r_C$ for all odd cycles $C$.  So $\W(\Gamma, f)$ is a homomorphic image of $\Cox(\Gamma, f)$.
\end{theorem}
\begin{proof}
To show that $\Cox(\Gamma) \to \W(\Gamma, f): x_i \mapsto w_i$ is a homomorphism, we must show that $w_i$ satisfy the defining relations of the Coxeter group.  Firstly note that $w_i$ has order 2.  If $i \sim j$, then $(e_i, e_j) = \pm 1$, and so the roots are at an angle $\pi/3$ or $2\pi/3$; so $w_i w_j$ is a rotation through the doubled angle $2\pi/3$ or $4\pi / 3$, and has order 3.  Similarly if $i \nsim j$ then $w_i w_j$ has order 2.

Consideration of the Weyl groups of odd cycles (~\cite[\S 3]{camseits}) then leads to the fact that the image of $r_C$ is the identity under the homomorphism.
\end{proof}

So $\im(r_C) = 1$.

Let $E$ be the adjacency matrix
\index{adjacency matrix}%
 of any connected signed graph $(\Gamma, f)$.  The \emph{intersection matrix}
\index{intersection matrix}%
 $2 I + E$ of $(\Gamma, f)$ is a symmetric $n \times n$ matrix with entries $2$ on the diagonal and $0, +1, -1$ elsewhere.  It can be interpreted as the Gram matrix
\index{Gram matrix}%
 of a root basis of the lattice $\mathbb{L}(\Gamma, f)$ (whose entries are the inner product of $n$ basis vectors).  That is, as the Gram matrix of the inner products of the basis vectors $e_1, \ldots, e_n$ of real $n$-space.  These vectors are roots having length $\sqrt{2}$ at angles $\frac{\pi}{2}, \frac{\pi}{3}$ or $\frac{2\pi}{3}$.  Then $\mathbb{L}(\Gamma, f)$ is the \emph{root lattice}, an even integral lattice spanned by vectors of norm 2.
 

Signed graphs with only the signing $``-"$ are called \emph{fundamental},
\index{fundamental signing}%
 because they correspond to fundamental bases for Coxeter groups.  If the graph $\Gamma$ has the fundamental signing $-$, then there are no odd cycles, and the homomorphism of the last theorem induces the isomorphism $\Cox(\Gamma)$ to $\W(\Gamma) : x_i \mapsto w_i$.  So we observe that $\Cox(\Gamma, -) = \Cox(\Gamma)$ where $-$ is the all-negative signing.

If $\Gamma$ is a tree,
\index{tree}%
 then  $\Cox(\Gamma, f) = \Cox(\Gamma)$ for any $f$, since all signings are equivalent to the fundamental signing $``-"$, and $\Cox(\Gamma, f) = \W(\Gamma, f)$, that is signed graphs correspond to fundamental bases for Coxeter groups.  For a cycle, there are just two switching classes, and so just two Coxeter groups, $\Cox(\Gamma)$ (the Coxeter group of type $\tilde{A}_{n-1}$) and the quotient $2^{n-1} \sd \Sym(n)$.

The determinant and signature of $2 I + E$ are invariants.  If $(\Gamma, f)$ is a signed complete graph, then its switching class is equivalent to a two-graph.
\index{graph ! two-graph}%
\index{two-graph}%
  Any such $(K_n, f)$ gives the $\pm1$ adjacency matrix
\index{adjacency matrix}%
 of a graph $\gamma$, and gives a Tsaranov group $\Ts(\gamma)$ which will be defined later.  We will say more about this special case at the end of the section.

Local switching, which we also define below, is a symmetric but not
transitive relation on the set $\Omega_n$ of switching classes of
signed graphs on $n$ vertices.  Its transitive closure partitions
$\Omega_n$ into clusters of connected classes, which can be used to
classify the lattices and groups.  Local switching leaves the root
lattice and Weyl group invariant, but does not in general leave
$\Cox(\Gamma, f)$ invariant (though it does so in the positive definite case).  Local switching has combinatorial, geometric and algebraic interpretations respectively for signed graphs, root
lattices and Weyl groups.  In the positive definite case, local
switching is equivalent to the Gabrielov transformation
\index{Gabrielov transformation}%
 in singularity theory~\cite{gabr}.
\index{singularity theory}%

Signed graphs arise all over mathematics.  T. Zaslavsky~\cite{zaslavsky}
\index{Zaslavsky, T.}%
has used them to study line graphs of finite loopless two-colour graphs.  It is instructive to recall one set of related results.  In this application, Cameron
\index{Cameron, P. J.}%
 has expressed generalizations of Mallows
\index{Mallows, C. L.}%
and Sloane's result~\cite{ms}
\index{Sloane, N. J. A.}%
on two-graphs
\index{graph ! two-graph}%
\index{two-graph}%
(see Chapter~\ref{chap2}) who treat the case of a complete underlying graph.  This was further extended by Wells~\cite{wells}
\index{Wells, A. L. jun.}%
who introduced the notion of \emph{even signings}
\index{signing ! even}%
as a generalization of even graphs and signed switching classes
\index{switching ! signed}%
which generalize ordinary switching classes.

A signing $f$ of a signed graph $(\Gamma, f)$ is \emph{even} if each vertex is incident with an even number of edges signed $``-"$.  Wells shows that the number of isomorphism types of even signings of a fixed simple graph $\Gamma$ equals the number of isomorphism types of signed switching classes of $\Gamma$.  He derives a formula for the number of these isomorphism types, and in addition, he finds a criterion for determining whether all signed switching classes fixed by a graph automorphism $\alpha$ actually contain signings fixed by $\alpha$.

Cameron
\index{Cameron, P. J.}%
 used cohomological algebra~\cite{cama} to explain properties of switching classes and two-graphs, and Cameron and Wells used it~\cite{camwel} to develop a theory of signings, signed switching classes and signatures, treating these as elements of certain distinguished vector spaces over $K = \{0, 1\}$, where the operations of switching and forming signatures are linear transformations on these binary vector spaces.  The binary vector spaces are of the form $V_i = \mathbb{F}_2^{X_i}$, where $X_{-1} = \{\emptyset\}$, $X_0$ denotes the set of vertices, $X_1$ denotes the set of edges (including vertices), and $X_2$ and $X_3$ are arbitrary sets of subgraphs.  The conditions to get a \emph{cochain complex} are to take $X_2$ as some cycles $C$ of the graph having all valencies even, whilst the edge of an $X_3$-graph should be in an even number of $X_2$-subgraphs.  Let $V_i$ be the binary vector space of formal sums of elements of $X_i$.  If $0 \leq i \leq 3$ define boundary maps $\partial_i : V_i \to V_{i-1}$ by taking the linear extension of the function mapping $X_i$-elements to the sum of its subgraphs in $X_{i-1}$.
 
Via coboundary maps $\delta^i$ this leads to the \emph{complex} sequence for dual vector spaces as in Figure~\ref{compseq}
\begin{figure}[!h]
$$\xymatrix{
{V^{-1}} \ar[r]^{\delta^{-1}} & {V^{0}} \ar[r]^{\delta^{0}} & {V^{-1}} \ar[r]^{\delta^{1}} & {V^{1}} \ar[r]^{\delta^{2}} & {V^{3}}
}$$ 
\caption{The Complex when conditions $\delta^{i} \delta^{i+1} = 0$ are satisfied for $-1 \leq i \leq 1$.}
\label{compseq}
\end{figure}

For $f \in V_0$ and $\{v, w\} \in X_1$, $\delta f (\{v, w\}) = f(w) - f(v) = f(w) + f(v)$, where the last equality follow as we are working$\mod 2$.  For $f \in V_1$ and $C \in X_2$, $\delta f (C) = \sum_{e \in C} f(e)$.  

Furthermore, $\delta^2 f (C) = \sum_{e \in C} \delta f(e) = 0$ if each $v$ is on 2 edges of $C$.  Also if $f \in V_1$ then $\delta^2 f (S) = \sum_{C \in S} \delta f(e) = 0$ if each edge is on 2 faces of $S \in X_3$.  So the even valency conditions on $X_2$ and $X_3$ are required to ensure $\delta^2 = 0$.

Then if $(\Gamma ; X_2, X_3)$ is a complex, it is possible to define binary vector spaces $H^i = \ker(\delta^i) / \im(\delta^{i-1})$, $0 \leq i \leq 2$.  The complex $(\Gamma ; X_2, X_3)$ is called \emph{exact} when all three binary vector spaces $H^0, H^1$ and $H^2$ vanish.  Cameron and Wells show that a sequence of dual binary vector spaces is exact precisely when the analogues of properties of switching-classes and two-graphs are satisfied, as encapsulated in, 
 
\begin{theorem}
$(\Gamma ; X_2, X_3)$ is exact if and only if all three of the following statements hold for all $X, Y \subseteq V(\Gamma)$, all $f, f_1, f_2 \in V'$ (where the elements of $V'$ correspond to signings of $\Gamma$), and $S \subseteq X_2$:

(i)  $f_X = f_Y$ if and only if $X = Y$ or $X = V(\Gamma) \backslash Y$;

(ii)  $f_1$ and $f_2$ are switching equivalent if and only if they have the same signature relative to $X_2$;

(iii)  A subset $S$ of $X_2$ is a signature if and only if each member of $X_3$ has an even number of elements of $S$ as subgraphs.
\end{theorem}

Zaslavsky
\index{Zaslavsky, T.}%
 has studied~\cite{zaslav} generalized two-graphs called \emph{togs},
\index{togs}%
 structures consisting of sets of polygons or cycles in a graph that satisfy relations that make them equivalent to switching classes of signings of a particular base graph.  
 
Further references to related studies of switching by Seidel and Taylor are~\cite{seideltaylor}~\cite{taylor} and~\cite{taylor1}.

The \emph{adjacency matrix}
\index{adjacency matrix}%
 $E = E(\Gamma, f)$ has entries 
 \begin{displaymath}
E_{ij} = \left\{ \begin{array}{ll}
f({i,j}) & \textrm{for $\{i, j\} \in E$}\\
0 & \textrm{otherwise.}
 \end{array} \right.
\end{displaymath}

Switching with respect to $X \subseteq V$ induces the transformation
$E \to DED$, where $D$ is the diagonal matrix with $-1$ (respectively
$+1$) at the positions of $X$ (respectively $V \backslash X$).  The
eigenvalues of $E$ and of $DED$ are the same and are called the
\emph{spectrum} \index{switching ! signed class ! spectrum}%
of the signed switching class.

As we have already noted, a tree has only one switching class because we can inductively switch all signs to $-1$, and a \emph{signed cycle}
\index{signed cycle}%
 has two switching classes.  The \emph{parity}
\index{signed cycle ! parity}%
 of a signed cycle is the parity of the number of its edges which carry a \emph{positive} sign.  The \emph{balance}
\index{signed cycle ! balance}%
 of a signed cycle, cf.~\cite{har} is the product of the signs on its edges.  If $\mathcal{C}$~\label{indcyc} is the set of induced cycles (those paths having no chords in the graph) in $(\Gamma, f)$, then in obvious notation $ \mathcal{C} = \mathcal{C}^{+} \cup \mathcal{C}^{-}$ (where $b = \pm$ is the balance).  Both parity and balance are switching class invariants.  For example, two signings of a connected graph are switching equivalent if and only if they have the same balance; equivalently stated, the space of cycles and cocycles are complementary~\cite{camwel}.  For signed complete graphs this reduces to the bijection between switching classes and two-graphs.
\index{graph ! two-graph}%
\index{two-graph}%

The \emph{root lattices}
\index{lattice ! root}%
 $\mathbb{L}(\Gamma, f)$ are formed from integral linear combinations of length $\sqrt{2}$ vectors (roots) which are at angles $\pi / 2, \pi / 3$ or $2 \pi / 3$.  Their Weyl groups $\W(\Gamma, f)$
\index{group ! Weyl}%
 are generated by the reflections $w_i$ in the hyperplanes orthogonal to the roots $e_i$, where

\[ w_i (x) = x - \frac{2(e_i, x)}{(e_i, e_i)} e_i = x - (x, e_i) e_i. \]

We review some salient points of this theory.  Assume that $\Gamma$ is connected.

Switching $(\Gamma, f)$ maps $e_i \to - e_i$, so leaves $\mathbb{L}(\Gamma, f)$ and $\W(\Gamma, f)$ invariant.

The group $\W(\Gamma, f)$ is finite if and only if $2 I + E$ is positive definite.  Then the root system is finite and $\Gamma$ is one of the \emph{spherical Coxeter-Dynkin diagrams}
\index{Coxeter-Dynkin diagram}%
 $A_n, D_n, E_6, E_7, E_8$, corresponding to integral bases of Euclidean
$\mathbb{L}(\Gamma, f)$
\index{lattice ! Euclidean}%
\index{Euclidean lattice}%
 of type $A, D, E$.  (Were $\Gamma$ to be disconnected, then it would be the direct sum of copies of these diagrams).  For if $2 I + E$ is positive definite then the orthogonal group
 \index{group ! orthogonal}%
 that it defines is compact, and its discrete subgroup $\W(\Gamma, f)$ is finite.  If $\Gamma$ is a tree, $\W(\Gamma, f)$ is signing-independent so we write $\mathbb{L}(\Gamma)$ and $\W(\Gamma)$.  The origin of the word spherical is that the Coxeter group $\Cox(\Xi)$ of a graph $\Xi$ in the so-called \emph{spherical case} is represented as a reflection group in spherical space.   As a footnote, we mention that at the heart of the ubiquity of the Coxeter-Dynkin diagram $ADE$-classification given by the Coxeter-Dynkin diagrams,
\index{Coxeter-Dynkin diagram}%
appears to be the ambiguity as to whether they describe either $3$-dimensional rotation groups and singularities (so that $E_6, E_7$, and $E_8$ represent the tetrahedron, octahedron and icosahedron), or root systems in higher dimensions (so that these three diagrams represent objects in Euclidean spaces of dimension 6, 7 and 8). 

The extended Coxeter-Dynkin graphs of $\tilde{D_n}, \tilde{E_6}, \tilde{E_7}, \tilde{E_8}$ are trees, but $\tilde{A_n}$ is a cycle with $n+1$ vertices, with $\W(A_n) \cong \Sym(n+1)$, and $\W(\tilde{A_n}) \cong \mathbb{Z}^n \sd \W(A_n)$.  Here $\mathbb{Z}^n$ is isomorphic to $\mathbb{L}(A_n)$ and the action of the quotient group by conjugation agrees with its natural action on the lattice.  Root lattices whose bases correspond to signed graphs with semi-definite or indefinite intersection matrices have infinite Weyl groups.  
Any connected signed graph with positive semi-definite intersection matrix
\index{intersection matrix}%
 is embeddable in a signed graph corresponding to a spherical root system.

We say that $2I + E$ (or the signed graph) has \emph{standard
 signature}
\index{standard signature}%
 if the signs of its eigenvalues are $\epsilon + + \ldots +$, where $\epsilon \in \{+, -, 0\}$.  If the intersection matrix
 \index{intersection matrix}%
  of a signed graph $(\Gamma, f)$ has standard signature, then the root lattice $\mathbb{L}(\Gamma, f)$ has a fundamental basis.  
(The intersection matrix of Slodowy's
\index{Slodowy, P.}%
 IM Lie algebra~\cite{slodowy} is related to the Cartan matrix of a Kac-Moody algebra
\index{Kac-Moody algebra}%
 by the maps
$-1 \to$
$ \begin{pmatrix} 
1&0\\
0&1
\end{pmatrix}$
and $+1 \to$
$ \begin{pmatrix} 
0&1\\
1&0
\end{pmatrix}$.)

\emph{Signed graphs may be considered to be global as compared to lattices which are local, because a signed graph corresponds to a particular root basis for the lattice; other signed graphs can give rise to the same lattice.  The localization of switching so as to relate signed graphs giving the same lattice motivated the following concept}:-

The \emph{local graph}
\index{graph ! local}%
 of $(\Gamma, f)$ at vertex $i$ has the set $N(i)$ of neighbours of
 $i$ as its vertex set, and as edges all edges $\{j, k\}$ of $\Gamma$
 for which $f(i, j) f(j, k) f(k, i) = -1$.  A \emph{rim}~\label{rim}
\index{graph ! local ! rim}%
 of $(\Gamma, f)$ at $i$ is any union of connected components of the local graph at $i$.

Let $J$ be any rim at $i$, and let $K = N(i) \backslash J$.
 The following defines \emph{local switching}
\index{switching ! local}%
 of $(\Gamma, f)$ with respect to $(i, J)$:

(i)  delete all edges of $\Gamma$ between $J$ and $K$; (any such edge $\{j, k\}$ satisfies $f(i, j) f(j, k) f(k, i) = +1$);

(ii)  $\forall j \in J, k \in K$ such that $j \nsim k$, introduce an edge $\{j ,k\}$ with sign chosen so that $f(i, j) f(j, k) f(k, i) = -1$;

(iii)  change the signs of all edges from $i$ to $J$;

(iv)  leave all other edges and signs unaltered.

Let $\mathbb{L}(\Gamma, f)$ be the root lattice corresponding to the signed graph $(\Gamma, f)$, and let $J$ be the rim at a vertex $i \in \Gamma$.  Local switching with respect to $(i, J)$ is the transformation $\prod_{j \in J} \sigma_{ij}$ of the root basis of $\mathbb{L}(\Gamma, f)$ corresponding to $(\Gamma, f)$, where if $w_i \in \W(\Gamma, f)$ then

\centerline{$\sigma_{ij} : (e_1, \ldots, e_j, \ldots, e_n) \mapsto (e_1, \ldots, w_i(e_j), \ldots, e_n).$}

Let $J$ be a rim at a vertex $i$.  The Weyl group $\W(\Gamma, f)$ is local switching invariant.  There are two equivalence relations associated with local switching on the class of signed graphs:

(i) at a vertex $i$, leaving the local graph at $i$ fixed,

(ii) the transitive closure of local switchings at a sequence of different vertices.

The equivalence classes of the transitive closure of the symmetric but not transitive relation on the set $\Omega_n$ of switching classes of $n$-vertex signed graphs, are called \emph{clusters} of order $n$.  Local switching preserves lattice invariants, such as the determinant and signature of $2I + E$.  If $(\Gamma, f)$ and $(\Gamma', f')$ are connected signed graphs such that the eigenvalues of $E(\Gamma, f)$ are greater than $-2$, then one can be locally switched into the other if and only if $\mathbb{L}(\Gamma, f)$ and $\mathbb{L}(\Gamma', f')$ have the same dimension \emph{and} type.

The braid group
\index{group ! braid}%
 $B_n = \langle b_i \rangle$ acts on the set $\mathcal{B}$ of ordered root bases:

\centerline{\quad\quad\quad\quad $b_i : (e_1, \ldots, e_i, e_{i+1}, \ldots, e_n) \mapsto (e_1, \ldots, e_{i+1}, w_{i+1}(e_i), \ldots, e_n).$}

A \emph{quasi-Coxeter element}
\index{quasi-Coxeter element}%
 $w = w_1 \ldots w_n \in \W(\Gamma, f)$ is associated with $(e_1, \ldots, e_n) \in \mathcal{B}$.  The element of $\W(\Gamma, f)$ corresponding to $b_i(e_1, \ldots, e_n)$ is $w_1 \ldots w_{i-1}w_{i+1}w_i^{w_{i+1}} \ldots w_n = w$.  The quasi-Coxeter element corresponding to a basis is constant on the orbits of $B_n$, for all root lattices.  In the positive definite $\mathbb{L}(\Gamma, f)$ case, the converse holds and $(e_1, \ldots, e_n) \mapsto (w_1 \ldots w_n)$ is a bijection between the set of $B_n$ orbits on $\mathcal{B}$ and the set of quasi-Coxeter elements of $\W(\Gamma, f)$.




If $\mathbb{L}$ is positive semi-definite, then $\W(\Gamma, f) = W_0
 \sd \mathbb{L}'$, where the abelian subgroup
\index{group ! abelian}%
 $W_0 \cong \bigoplus_{i = 1}^{n-d} (\mathbb{L}')_i$, $\mathbb{L}'$ is a positive definite root lattice with Weyl group \index{group ! Weyl}%
 $W'$, and the kernel of $\mathbb{L}$ has codimension $d$.  Here $W_0$ is free abelian
\index{group ! free abelian}%
 of rank $d (n - d)$ and $W'$ has type $A_d, D_d$ or $E_d$.  However if $n - d > 1$ then the homomorphism  from $\Cox(\Gamma) \to \W(\Gamma, f) : x_i \to w_i, (i = 1, \ldots, n)$, whose kernel contains the elements $r_C$ for all odd cycles $C$, is not necessarily faithful.

Let $\mathbb{L}(\Gamma, f)$ and $\mathbb{L}(\Gamma', f')$ be equivalent under local switching, and suppose that $2 I + E(\Gamma, f)$ is positive definite.  Then the following isomorphisms are true $\Cox(\Gamma, f) \cong \Cox(\Gamma', f') \cong \Cox(\Xi)$, for some Coxeter-Dynkin graph $\Xi$.  Note that positive definite means that in $(\Gamma, f)$ and $(\Gamma', f')$, all signed cycles have odd parity.  Also that any local switching is a composition of `elementary' local switchings, that is those having rims of order 1.

\head{Signed Complete Graphs} $(K_n, f)$.
\index{graph ! signed complete}%

Switching of signed complete graphs corresponds to Seidel
\index{Seidel, J. J.}%
 switching of graphs.  Off-diagonal entries of $E$ are either $-1$ (adjacency) or $+1$ (non-adjacency) in $\Gamma$.  In a complete graph the only induced cycles are the triangles.

Signed complete $n$-vertex graphs corresponding to ordinary graphs $\Gamma$ can be lifted to signed graphs $(\overline{\Gamma + 0}, +)$ on $n+1$ vertices, where the bar denotes complementation.

The signed graph denoted $(\overline{\Gamma + 0}, +)$, has $\Gamma + 0$ as the disjoint union of $\Gamma$ and an isolated vertex $0$.  The connected graph $(\overline{\Gamma + 0})$ yields a root lattice $\mathbb{L}(\overline{\Gamma + 0}, +)$.  Any signed complete graph can be represented as a bundle of stars (through any one root) in a root system.

\head{Tsaranov Groups}
\index{group ! Tsaranov}%

For an $n$-vertex graph $\Gamma$, the \emph{Tsaranov group} of $\Gamma$ is defined as,

$\Ts(\Gamma) : = \langle t_1, \ldots, t_n : t_i^3 = 1, (t_i t_j^{-1})^2 =1, \text{if } i \sim j, (t_i t_j)^2 =1, \text{if } i \nsim j \rangle$.

  For $(K_n, f)$, the last two relations can be written as $(t_i t_j^{f(\{i, j\})})^2 = 1$.  One of the key features is that switching equivalent signings give isomorphic Tsaranov groups. 

We can extend $\Ts(\Gamma)$ by an automorphism which inverts every generator:
\[ \Ts^{*}(\Gamma) = \langle \Ts(\Gamma), t : t^2 = 1,\ t t_i t = t_i^{-1},\ i = 1, \ldots, n \rangle. \]

With substitutions $x_0 = t,\ x_i = t u_i\ (i = 1, \ldots, n)$, 
\[ \Ts^{*}(\Gamma) \cong \Cox(\overline{\Gamma + 0}) / N \cong \Cox(\overline{\Gamma + 0}, +) \]
where $N$ is the normal closure of 
\[ \{ (x_0 x_i x_0 x_j)^2 : \{0, i, j\}\ is\ a\ triangle\ in\ \overline{\Gamma + 0}\},\] which are just the cut relations for the triangles if we impose the positive signing, and where the $+$ sign in the final term denotes a positive intersection matrix.  
\index{intersection matrix}%

\head{Two-graphs arising from Trees}
\index{graph ! two-graph}%
\index{two-graph}%

Let $T$ denote a tree which we assume has the fundamental signing, with its intersection matrix $2I - A(T)$ being the Gram matrix
\index{Gram matrix}%
 of a root basis $B = \{e_k\}$ of the root lattice $\mathbb{L} (T, -)$.  Take a fixed vertex $0$ as root.  For any other vertex $j$, let $(0, j) = (j_1 = 0, j_2, \ldots, j_s = j)$ be the unique path from $0$ to $j$.  Then $D = \{d_j\};\quad d_j = \Sigma_{i=0}^{s} e_{j_i}$ is a new basis for $\mathbb{L} (T, -)$.  It is possible to check that 
 \begin{displaymath}
(d_{j}, d_{k}) = \left\{ \begin{array}{ll}
2 & \textrm{if $j = k$}\\
1 & \textrm{if $(0, j) \subset (0, k)\ or\ (0, k) \subset (0, j)$}\\
0 & \textrm{otherwise.}
\end{array} \right.
\end{displaymath}

Thus the Gram matrix
\index{Gram matrix}%
 of $D$ has the form $2I + E(\Delta, +)$, where $\Delta$ is the graph on the same vertex set as $T$, in which $i \sim j$ if and only if $(0, i) \subset (0, j)$ or vice versa.  The graph $\Gamma(T)$ is said to \emph{arise from} the tree $T$ if the graph is obtained from the complement of $\Delta$ by deleting the vertex $0$ and all edges issuing from it.

In the tree $T$, we call $j$ a \emph{predecessor} of $k$, or $k$ a \emph{successor} of $j$, if $(0, j) \subset (0, k)$.  If in addition, $j \sim k$ (in $T$), then $k$ is an \emph{immediate successor} of $j$.

\begin{theorem}
Given a tree $T$ with vertex $0$, let $\Gamma(T)$ be defined as above.  Then
\[ \Ts^*(\Gamma(T)) \cong \Cox(T).\]
\end{theorem}

Different choices of root vertex $0$ yield switching-equivalent graphs.

\bigskip

To summarize:  in~\cite{camseits} the Coxeter group of a signed graph was defined to be the quotient of the usual Coxeter graph by relations from the presentation of the Weyl group, and in this case it is isomorphic to the Weyl group.  For each chordless cycle in the graph on vertices $x_1, \ldots, x_n$ carrying an odd number of $+$ signs there is such a relation of the form
\begin{equation}
  (x_1x_2 \ldots x_{n-1} x_n x_{n-1} \ldots, x_2)^2 = 1.
\end{equation}
This relation is independent of the starting point or direction of the cycle.  This paper also introduced the idea of \emph{local switching} which correspond to natural operations on bases for a given root system, as well as ``clusters'' -- partitions of signed graphs formed by repeated local switching.  The paper investigated the case of a non-positive definite adjacency matrix, which corresponds to a root system in a real vector space with indefinite inner product.  In the positive semidefinite case the authors used the theory of vanishing lattices from singularity theory to explain the positive semi-definite case where the Weyl group is a semi-direct product of the additive group of a lattice with a finite Weyl group. 

In~\cite{barot}, Michael Barot and Robert Marsh
\index{Marsh, R.}%
\index{Barot, M.}%
 discovered that the simply-laced ($ADE$) case of local switching corresponds to cluster mutation in the theory of \emph{cluster algebras}.
\index{cluster algebras}%
  We refer to the papers of Fomin and Zelevinsky~\cite{fomin}
\index{Fomin, S.}%
\index{Zelevinsky, A.}%
 for the definition and theory of these structures, but mention that a central theorem says that a mutation class of seeds is finite if and only if one of them has a graph which is a Coxeter-Dynkin diagram.  This is surprising partly because cluster algebras are associated with directed graphs, whereas Coxeter-Dynkin diagrams are undirected. In order to make the correspondence work, it is necessary to symmetrise the digraphs.  The Weyl group presentations in the $ADE$ case are the same in~\cite{camseits} as in~\cite{barot}, however the work in~\cite{barot} is more general in that they deal with the other Coxeter-Dynkin diagrams (of types $BFG$) as well, and that in~\cite{camseits} is more general in that they study the non-positive-definite case.

\section{Local Switching of Multicoloured Graphs}

In this section we begin a generalization of the notion of local switching from two-coloured to multicoloured graphs.  It will transpire that there is a slightly generalized closed switching group and accompanying conserved parity, analogous to those identified in the final section of the previous chapter.

Assume that $m \ge 3$.  Take as the basic object a complete graph with an ordering of the colours on each edge.  This object is equivalent to an $m$-tuple of edge-coloured graphs, say $\Gamma_1, \ldots, \Gamma_m$ such that for every edge $e$, all $m$ colours appear once on $e$ in $\Gamma_1, \ldots, \Gamma_m$.  We call this object, a \emph{set of complete multicoloured graphs with a list of edge-colours on each edge}, an \emph{SML graph}
\index{graph ! SML}%
$\Upsilon$~\label{smlgraph}, which comprises an underlying graph $\Gamma$ and an ordered tuple of colours on each edge.  

There is a dual way to view this definition.  Either (a) as a complete graph with an ordering of colours on each edge, or (b) as an $m$-tuple of edge-coloured complete graphs such that each colour occurs on a given edge in one of the graphs. 

Perhaps category theory
\index{category theory}%
  could provide a suitable terminology and techniques for its study, but we want to avoid being taken too far from the concepts of the two-coloured case, where the switching simply means colour transposition on an edge.

Take a standard reference SML graph with an ordered $m$-tuple $(c_1, \ldots, c_m)$ of colours on \emph{every} edge.  There is a map $\phi$ from $\Upsilon$ to signed complete graphs such that $\phi(\Upsilon)$ has sign $+$ (respectively $-$) on $e$ according to whether the list of colours on $e \in \Upsilon$ is an even (respectively odd) permutation of the starting colour configuration.

A switching $\sigma_{v, \pi}$ means apply permutation $\pi$ to the list of colours on all edges containing vertex $v$.

Now,  \begin{displaymath}
\phi(\sigma_{v, \pi} (\Upsilon)) = \left\{ \begin{array}{ll}
\phi(\Upsilon) & \textrm{if $\pi$ is an even switching of $\phi(\Upsilon)$ at $v$,}\\
\sigma_{v}(\phi(\Upsilon)) & \textrm{if $\pi$ is an odd switching of $\phi(\Upsilon)$ at $v$,}\\
& \textrm{where $\sigma_v$ is ordinary switching.}
\end{array} \right.
\end{displaymath}

Henceforth we write $(\Upsilon, \phi)$ for an SML graph, complete with its parity signings.  With the above assumption on the standard reference graph, for $m = 2$ colours $\phi$ would reduce to an isomorphism with: $(c_1, c_2) \to +$, and $(c_2, c_1) \to -$.

We can then generalize the $m = 2$ definitions to the case of $m \ge 3$:

\begin{itemize}

\item  The \emph{parity} (sign parity) of a signed cycle of multicoloured edges is the parity of the number of its edges which carry a \emph{positive} sign (derived from an even permutation of the colours on that edge.)

\item  The \emph{balance} of a signed cycle of multicoloured edges is the product of the signs of the colour permutations on its edges.
\end{itemize}

A signed cycle has two switching classes, as for the $2$-coloured case.  Both parity and balance are switching class invariants but are not independent graph invariances; they are essentially the same thing, for example, two signings of a connected graph are switching equivalent if and only if they have the same balance. 

The multicoloured version of some basic facts of two-coloured graphs are true, as in the next two results:

\begin{lemma}
\label{treeprop}
Two SML graphs based on a finite tree are switching equivalent.
\end{lemma}

\begin{proof}
Induct on the number of vertices.  Consider an SML graph $\Upsilon$ based on a tree $T$ with $n$ vertices and where $v$ is a leaf.

Let $T' = T \backslash \{v\}$, with $\Upsilon'$ the (restricted) SML graph based on $T'$.  By induction, switch $\Upsilon'$ such that all edge permutations are the identity.  The switchings involved will effect an induced switching about $v$, call it $\pi$, in the enlarged tree $T$.  Now a switch by $\sigma_{v, \pi^{-1}}$ at $v$ proves the proposition.
\end{proof}

We have stated the proof for finite trees, but because infinite trees have no leaves it would require some more formalism to prove the theorem still works.

\begin{proposition}
Two connected SML graphs are switching equivalent if and only if they have the same balance.
\end{proposition}

\begin{proof}
Consider first the spanning tree
\index{graph ! spanning}%
 of an arbitrary SML graph.  By Lemma~\ref{treeprop} we can switch so as to effect the identity permutation on all its edges.  Consider any other edge of the underlying graph and without loss of generality consider an even permutation of the colours on it.  Even permutations lie in the derived subgroup of the permutation group of colours.  So consider switchings about the two vertices $v, w$ on the chosen edge.  Take $g \in \sigma_{v, \pi_1}$ and $h \in \sigma_{w, \pi_2}$, where $\pi_1, \pi_2 \in \Sym(c_1, \ldots, c_m)$, we can see that the commutator $g^{-1} h^{-1} g h$ lies in $\Alt(c_1, \ldots, c_m)$, and so a graph switching is a product of commutators.  For the commutator is the identity on all edges except $\{v, w\}$.  This is clear for edges containing neither $v$ nor $w$.  If $\{v, u\}$ is an edge then $h$ acts trivially on colours of $\{v, u\}$ so the commutator acts as $g^{-1}g = 1$.  Similarly for edges containing $w$.

 Therefore SML graphs are switching equivalent if and only if the parities of the cycles are the same.

\end{proof}

This proposition proves that balance is an invariant with respect to switching.  

For each fixed number $m$ of colours, there is an inverse system of switching groups on graphs with $n$ vertices, which in the $n \to \omega$ limit gives a \emph{profinite group}.
\index{group ! profinite}%
  This group is isomorphic as a topological group
\index{topological group}%
 to a closed subgroup of a Cartesian product of
symmetric groups of degree $m \ge 3$, in the topology defined by
the imprimitive action
\index{group ! permutation ! imprimitive}%
 for the cartesian product where each
fibre comprises the colour set $(c_1, \ldots, c_m)$ of an edge.   There is a parity equivalence between two $m$-coloured SML graphs, if the ordered $m$ tuple of colours of each edge of one graph is a parity-preserving permutation of those on the second graph.   This parity is preserved by any switching $\sigma_{X, \pi} := \prod_{v \in X} \sigma_{v, \pi}$ which has an even parity, and the elements of the closed switching group are precisely the parity-preserving permutations fixing all the fibres on edges.  This slight generalization of Theorem~\ref{parthm}, embodied in Lemma~\ref{treeprop} and the subsequent proposition,  works because the cycle space is spanned by triangles.

\smallskip

What novel features arise in the local switching theory when we go from two to many colours?  

First we mention that we loose the original geometric motivation, which arises from the fact that if the symmetric adjacency matrix
\index{adjacency matrix}%
 of a signed graph $E$ (with entries $0, \pm 1$) has least eigenvalue $- \alpha$, then $\alpha I + E$ defines a set of lines whose angles $\theta$ satisfy $cos(\theta) \in \{0, \frac{1}{\alpha}\}$.  A signed complete graph defines a two-graph
\index{graph ! two-graph}%
\index{two-graph}%
 and a set of equiangular lines~\cite{seidel}.
\index{equiangular lines}%
 Seidel
\index{Seidel, J. J.}%
 switching reverses the directions of the lines.  But in the case of many colours, there is no clear geometric
interpretation; if there are as many angles between the lines as there are colours, then what does multicoloured switching do?

What we gain however is a \emph{duality} on SML graphs between graph number and edge colour, that is a duality between the $i$th graph having colour $j$ on an edge and the $j$th graph having colour $i$ on that edge.  Here colour refers to either the colour on a fixed reference edge or a reference ordered tuple of colours.  In the two-colour case this reduces to a \emph{self-duality}.  This fact provides us with $m^2$ matrices as a function of both colour and edge numbers.  For the inner product matrix of a set of vectors $v_i$ has $v_i \cdot v_i = 2$ along the diagonal and  $v_i \cdot v_j \in \{\pm 1, 0 \}$ $i \neq j$.  If the matrix is positive definite then these vectors span a root lattice.  Let $i, j \in \{1, \ldots, m\}$ and let $v, w$ be vertices of a graph.  Let $M_{ij}$ be the matrix with $v, w$ entry 0 if $v=w$ or if $v \nsim w$, and if $v \sim w$ then an entry of 1 in the case of the $i$th colour on the list on $\{v, w\}$ $(i \leq j)$ and 0 otherwise.  Each $M_{ij}$ is then a symmetric $0$--$1$ matrix such that
\[ \sum_i M_{ij} = \sum_j M_{ij} = A, \]
where $A$ is the adjacency matrix of the graph.
\index{adjacency matrix}%

We leave for future work, the study of a object that is a dual to an SML graph.  Define a \emph{dual SML graph} DSML
\index{graph ! DSML}%
 to be \emph{a colour with a list of graphs forming an SML graph}, which we can symbolise as $\Upsilon^D$.~\label{dsmlgraph}
 


\smallskip

Any symmetric is the inner product matrix of a set of vectors in some possibly indefinite (and degenerate) inner product space.  But it may happen that linear combinations of these vectors are dense in some subspace, so not a lattice. 
matrix with $0, \pm 1$ entries can be taken as the Gram matrix
\index{Gram matrix}%
 of a root basis giving rise to a lattice.  More specifically, given an SML graph $(\Upsilon, \phi)$, can we use it to define root lattices and Weyl, Coxeter and Tsaranov groups, and if so what do they look like?


In the two-colour case, the root lattice $\mathbb{L} (\Gamma, f)$ and Weyl group of $\mathbb{L}$ are local switching invariants, because switching $(\Gamma, f)$ reflects the root bases of the lattice, thus $e_i \mapsto - e_i$.  An SML graph $(\Upsilon, \phi)$ may be considered as a \emph{set of global} objects corresponding to a \emph{local} set of lattice root bases - the same set of $m$ lattices can give different ($m$-sets of) SML graphs.








But what would be the combinatorial origin of these? 


The signed $(\pm)$ switching that we have derived from the parity of permutations of colours on the graph edges simulates two-colour switching, and allows us to retrieve properties of random graphs
\index{graph ! random}%
 that are lost for graphs of three or more colours.  In particular there is a link between cohomology
\index{group ! cohomology}%
 and two-graphs
\index{graph ! two-graph}%
\index{two-graph}%
 which we now briefly recall, leaving the full details to references~\cite{cama}~\cite[p.57]{evans1}.

Recall from Chapter~\ref{chap2}, that the set of all graphs
on a given vertex set forms a $\mathbb{Z}_2$-vector space, where the
sum of two graphs is obtained by taking the symmetric difference of
their edge sets.  Switching a graph corresponds to adding to it a
complete bipartite graph.
\index{graph ! bipartite}%
 So the switching classes are the cosets of the subspace $V_0$ of all complete bipartite graphs, including the null and complete graphs; alternatively $V_0$ is the set of all
switching operations.
\index{switching ! operation}%

Let $\Gamma$ be a graph giving rise to a two-graph
$\mathcal{T}(\Gamma)$ with automorphism group
\index{group ! automorphism}%
 $G$.  Because any two such
graphs lie in the same switching class, we can define a derivation
function $d: G \to V_0 : g \mapsto (g \Gamma - \Gamma)$.  There is an element of the \emph{first cohomology group}
\index{group ! first cohomology}%
 of $G$ on $V_0$, $\gamma \in H^1(G, V_0)$,~\label{fstinv} called the \emph{first invariant}
\index{two-graph ! first invariant}%
 of $\mathcal{T}(\Gamma)$ with the property that $\gamma = 0$ if and only
if there is a graph $\Gamma' \in V_0 + \Gamma$ such that $G =
\Aut(\Gamma')$.  In other words, that there exist switching classes where the first cohomological invariant is non-zero, means that there are switching classes for which the group of automorphisms is strictly larger than the group of automorphisms of its graphs.  The invariant that we derive from the $\mathsf{LP}_m$ parity of the SML graphs is isomorphic to this one for simple graphs.  There is a double cover of
$\Aut(\mathcal{T}(\Gamma))$ which is strongly split
\index{group ! split extension}%
 if and only if $G$ preserves a graph in the same switching class as $\mathfrak{R}$.  The $2$-transitive
\index{group ! permutation ! $2$-transitive}%
automorphism group
\index{group ! automorphism}%
 of the two-graph corresponding to the random graph
\index{graph ! random}%
cannot preserve any graph in the same switching class as
$\mathfrak{R}$, so $\gamma$ cannot be zero and the double cover is not
strongly split.  As point stabilizers in $G$ are irreducible, the
double cover of the two-graph of the random graph is not split.
\index{graph ! two-graph}%
\index{two-graph}%

David Harries and Hans Liebeck~\cite{harries}
\index{Liebeck, H.}%
\index{Harries, D.}%
 define a permutation group $G$ on a finite set $X$ to be \emph{always exposable}
\index{group ! permutation ! always exposable}%
 if whenever $G$ stabilizes a switching class of graphs $\Gamma$ on $X$, $G$ fixes a graph in the switching class.  In terms of Cameron's graph invariants~\cite{cama}, this is equivalent to saying that $G$ is always exposable if the first invariant of $G$ and every two-graph $\mathcal{T}(\Gamma)$ on operand $X$ are zero.  The problem they solved was which permutation representations of $G$ are always exposable when $G$ is dihedral, noting that Mallows and Sloane had already solved the cyclic group case in~\cite{ms}.
\index{Mallows, C. L.}%
\index{Sloane, N. J. A.}%
 In~\cite{liebeckb}, Hans's son Martin
\index{Liebeck, M. W.}%
 solved the problem for several other types of groups.

\bigskip

We shall revisit the combination of graphs, groups and lattices in a
different context in Chapter~\ref{homcochap}, where we shall prove further results.

\bigskip

\head{Open Question. Coloured Coxeter and Tsaranov Groups}

In order to make contact with the theory expounded in~\cite{camseits} for two-colour graphs, elements of which were summarised in the first section of this chapter, and to fully develop multicoloured versions of the connections made through local switchings between lattices, Weyl groups and Lie algebras, we require multicoloured versions of Coxeter and Tsaranov groups, and one direction to pursue would be through presentations of multicoloured switching groups, which we outline in the two chapter appendices.


Recall that the \emph{Coxeter group}
\index{group ! Coxeter}%
\index{Coxeter ! group}%
$\Cox(\Gamma)$~\label{cox} of a simple (that is two-coloured) graph $\Gamma$ is defined to have a generator $s_v$ for each graph vertex $v$ such that $s_v^{2}=1, (s_v s_w)^3 = 1$
if $v \sim w$, and $(s_v s_w)^2 = 1$ if $v \nsim w$ because all the relations have even length.  The set of all products of even length in the generators forms a normal subgroup of index $2$, called the \emph{even part} of the Coxeter group,
$\Cox^{+}(\Gamma)$~\label{cox+}.  Analogously the \emph{Tsaranov group},
\index{group ! Tsaranov}%
\index{Tsaranov, S. V.}%
$\Ts(X, \Omega)$~\label{Ts} has a generator $t_v$ for each vertex such that $t_v^{3}=1,
(t_v t^{-1}_w)^2 = 1$ if $v \sim w$, and $(t_v t_w)^2 = 1$ if $v \nsim
w$.  Switching a graph with respect to a vertex set corresponds to
replacing Tsaranov's generators for vertices in this set with their
inverses; so $\Ts(X, \Omega)$ is an invariant of the switching class,
and so therefore, of the two-graph.  
 
Find a generalization of the theory of this chapter, and in particular of multicoloured Coxeter and Tsaranov groups that reduces to these definitions in the two-coloured case.  Could the derived coloured Coxeter and Tsaranov groups have the SML graphs
\index{graph ! SML}%
 as operands?

\bigskip
\bigskip

\head{Some Results on Trees}

A \emph{tree}
\index{tree}%
is a connected graph without cycles.  There are maps from trees to two-graphs with a Tsaranov group $\Ts(X, \Omega)$ being represented by the two-graph~\cite{cam3a}.  A two-graph contains a graph $\Gamma$ as an \emph{induced subgraph}
\index{graph ! induced subgraph}%
 if $\Gamma$ is a subgraph of some graph in the switching class of the two-graph.  The Tsaranov group is the even subgroup of the Coxeter group of the tree and if $(X, \Omega)$ arises from a tree
\index{tree}%
 $T$ by a certain map, then $\Ts(X, \Omega) \cong \Cox^{+}(T)$~\cite{camseits}, which is proved in~\cite{seideltsar} as follows.  
 
 \begin{theorem}[Seidel and Tsaranov]
\index{Seidel, J. J.}%
\index{Tsaranov, S. V.}%
If $(X, \Omega)$ is a two-graph arising from a tree $T$, then $\Ts(X, \Omega) \cong \Cox^{+}(T)$.
\end{theorem}

\begin{proof}
Assume the Coxeter group $\Cox(T)$ is generated by involutions representing the vertices of the tree $T$.

We require the following construction~\cite{cam3a}.  Let $E(T)$ be the edge set of the connected tree $T$.  The triples of edges of $T$ are of two kinds:  one edge is situated between the other two, or none is between the other two.  If $\Omega$ is the set of triples of edges in $T$ of the type none between the others, then consideration of the possibilities for quadruples of edges leads us to $(X, \Omega)$ being a two-graph.
\index{graph ! two-graph}%
\index{two-graph}%

We will use an equivalent way of phrasing this construction from $(n+1)$-vertex trees $T$ to two-graphs $(X, \Omega)$ on $n$ vertices, which has been given by Ivanov~\cite{seideltsar}.
\index{Ivanov, A. A.}%
 Let $v$ be a terminal vertex of $T$. Define a graph $\Gamma_v$ on the vertex set $V(T) \backslash \{v\}$, where two vertices $x, y$ are adjacent if and only if neither one of the shortest paths, $x, v$ and $y, v$ in $T$ contains the other.  Distinct vertices $u, v \in T$ yield switching equivalent graphs $\Gamma_u = \Gamma_v$.  This gives a two-graph $(X, \Omega)$ from the tree $T$.

Now, let $\Ts(\Gamma_v)$ be generated by order-3 elements $x_w$, indexed by $w \in V(\Gamma_v) = V(T) \backslash \{v\}$.  Let $(v, s, t, \ldots, u, w)$ denote the shortest path connecting $v$ and $w$.  Associate with each vertex $w$ the element $y_w := v \cdot v^{st \ldots uw} \in \Cox(T)$.

It can be shown that the following two maps $\phi, \phi'$ are both epimorphisms
\[ \phi : \Ts^{*}(\Gamma_v) \to \Cox(T) : x_w \mapsto y_w,\quad\quad x_0 \mapsto v \]
\[ \phi' : \Cox(T) \to \Ts^{*}(\Gamma_v) : w \mapsto (x_0 x_w)^{x_u^{-1}} = x_u x_0 x_w x_u^{-1},\quad\quad v \mapsto x_0 \]

The composite map $\phi \circ \phi'$ acts on the generators of $\Cox(T)$ by
\[ (\phi \circ \phi')(w) = \phi(\phi'(w)) = \phi(x_u x_0 x_w x_u^{-1}) = y_u y_0 y_w y_u^{-1}. \]

This term can be shown to equal $w$ and so $\phi \circ \phi'$ is the identity, and both $\phi$ and $\phi'$ are isomorphisms.
\end{proof}

If $T$ is a tree with edge set $E$ and $\Omega$ is the set of $3$-subsets of $E$ not contained in paths in $T$, then $(E, \Omega)$ is a two-graph.
\index{graph ! two-graph}%
\index{two-graph}%
 There is a map~\cite{cam3b} from trees to two-graphs characterising its image by excluded substructures.  It was proved that a two-graph $(E, \Omega)$, arises from a tree by this construction if and only if it contains neither the pentagon nor the hexagon as induced substructures.  Further if the two-graphs arising from trees $T_1, T_2$ are isomorphic then $T_1 \cong T_2$.  With a rooted tree $(T, r)$,~\label{roottree} there is a partial order on the vertices with $r$ as the least element.  The graph $\Gamma(T, r)$ of this partial order, with $r$ deleted has two vertices, $x, y$ nonadjacent if and only if the path from $r$ to $x$ contains $y$ or vice versa.  A graph $\Gamma$ is isomorphic to $\Gamma(T, r)$ for some rooted tree if and only if $\Gamma$ contains neither a path of length 3 nor 2 disjoint edges as induced subgraphs.  Furthermore if $\Gamma(T_1, r_1)$ and $\Gamma(T_2, r_2)$ are isomorphic graphs then $(T_1, r_1)$ and $(T_2, r_2)$ are isomorphic rooted trees.  Also under this map a graph $\Gamma \cong \Gamma(T, r)$ is associated to some rooted tree $(T, r)$ if and only if $\Gamma$ is \emph{N-free}
 \index{graph ! N-free@$N$-free}%
(where an \emph{N} is a path of length $3$) and contains no two disjoint edges if and only if $\Gamma_x$ is \emph{N-free}, where $\Gamma_x$ is the graph in its switching class for which $x$ is isolated, for some or all $x \in X$~\cite{cam3b}.

Let $\Gamma_x$, $x \in X$ be the unique graph in the switching class for which $x$ is an isolated vertex.  A pentagon can be switched into a path of length 3 and an isolated vertex.  The following are equivalent:

(i) $\Gamma_x$ is N-free for some $x \in X$;

(ii) $\Gamma_x$ is N-free for all $x \in X$;

(iii) The two-graph $(X, \Omega)$ is pentagon-free.
\index{graph ! two-graph}%
\index{two-graph}%

The amalgamation property can be restored to the class of pentagon-free two-graphs by adding a quaternary relation.

By a second construction~\cite{cam3a}, a two-graph arises from a tree if and only if the two-graph does not contain the pentagon as an induced substructure.  Here, only isomorphic coloured reduced trees yield isomorphic two-graphs.  There are also formulae~\cite{cam3b} for the number of labeled two-graphs obtained by the above constructions.

\head{Open Questions.}

1.  Can we generalize the two-colour results of Cameron~\cite{cam3a}~\cite{cam3b}
\index{Cameron, P. J.}%
perhaps using a multicoloured parity, for example $P_m$-equivalence
\index{graph ! Pmeq@$P_m$-equivalent}%
of Chapter~\ref{chap2} or the $\mathsf{P}_m$ of Chapter~\ref{chap3} or some other version, and find maps from trees with multicoloured edges to $\Ts(\mathfrak{R}_{m,n})$?  Is there a condition, perhaps related to monochromatic N-freeness Chapter~\ref{chapFD}), under which the coloured graphs are associated to rooted trees?

2.  Homomorphisms between edge-coloured graphs
\index{graph ! homomorphism}%
 and Coxeter groups arise naturally and have been studied before, for example in~\cite{alonmarshall}, and can also lead to a fruitful avenue for research.  For two edge-coloured graphs $\Gamma_1 = (V_1, E_1), \Gamma_2 = (V_2, E_2)$ with no multiple edges or loops, there is a \emph{homomorphism} $\phi : V_1 \mapsto V_2$ if for every pair of adjacent vertices $u, v \in \Gamma_1$, $\phi(u)$ and $\phi(v)$ are adjacent in $\Gamma_2$ and the colour of edge $\phi(u)\phi(v)$ is the same as that of the edge $uv$.  

Form an edge-coloured graph $\gamma(\mathfrak{R}_{m,n})$ from the presentation of $\Ts(\mathfrak{R}_{m,n})$ taking the generators $a_{ij}^{k}$ as vertices and joining $a_{ij}^{k}$ and $a_{ij}^{l}$ whenever $(a_{ij}^{k} a_{ij}^{l})^2 = 1$.  The index two subgroup $\Ts^{0}(\mathfrak{R}_{m,n})$ consists of products of an even number of generators of the form $\{\gamma_{kl} = a_{ij}^{k} a_{ij}^{l}\}$ and relations $a_{ij}^{3}=1$.  The $\gamma_{kl}$ represent directed edges in $\gamma(\mathfrak{R}_{m,n})$.  The cubed relations correspond to directed triangles and relations of the form $\gamma_{kl} \gamma_{lm} \ldots \gamma_{qk}$ to directed circuits.  A free product
\index{group ! free product}%
 of Tsaranov groups leads to a disjoint union of graphs.  Can this formalism be developed?

\section{Appendix 1: A First Presentation of $S_{m,n}$}

A \emph{Coxeter group} 
\index{Coxeter ! group}%
\index{group ! Coxeter}%
has generators $s_1, \ldots, s_k$ and relations $s_i^2 = 1$ and $(s_i s_j)^{m_{ij}} = 1$, where $m_{ij} = \{2, 3, \ldots, \infty\}$.  The finite Coxeter groups are the finite real reflection groups~\cite{hum}.  
\index{group ! reflection}%

The Coxeter group $\Cox(\Gamma)$ may be defined on a graph $\Gamma$ and is generated by vertices $s_i \in \Gamma$ $(i \in I)$, such that
\begin{displaymath}
\begin{array}{ll}
s_i^{2}=1 & \textrm{$\forall s_i \in \Gamma$};\\
(s_i s_j)^2 = 1 \Leftrightarrow s_i s_j = s_j s_i & \textrm{if $s_i, s_j$ are disjoint};\\
(s_i s_j)^3 = 1 \Leftrightarrow s_i s_j s_i = s_j s_i s_j & \textrm{if $s_i, s_j$ are joined by a single edge};\\
(s_i s_j)^{m_{ij}} = 1& \textrm{where an edge is labelled $m_{ij}$ or}\\
{} & \textrm{there are $m_{ij} - 2$ parallel edges.}
\end{array}
\end{displaymath}\\
For example, $\Sym(n)$ is a Coxeter group represented by the following diagram:
$$\xymatrix{
{\circ} \ar@{-}[r] & {\circ} \ar@{-}[r] & {\circ} \ar@{-}[r] & {} {\ldots} {\ldots} \ar@{-}[r] & {\circ}
}$$
which has $n-1$ nodes indexed by $(i \in I)$, and where $s_i \mapsto (i, i+1)$. 

In this chapter we derive two presentations of the
switching groups $S_{m,n}$ for graphs on $n$ vertices and $m$ colours, the second of which leads us to define and study a type of generalized Coxeter group, with the following as basic motivation.
\index{Coxeter ! group}%
\index{group ! Coxeter}%

Any group is isomorphic to a quotient of a free group.
\index{group ! free}%
  If the map $\phi: F \to G$ is a homomorphism from an $n$-generator free group $F$ onto an $n$-generator group $G$, then $G$ is presented by the generators of $\ker(\phi)$ as normal subgroup of $F$ and the relations $\ker(\phi)$.   To show that there is a homomorphism from a group $G = \langle g_1, \ldots, g_n | r_1, \ldots r_k \rangle$  to a group $H$, we need to find $h_1, \ldots, h_n \in H$ such that $r_i(h_1, \ldots, h_n) = 1$ in $H$.  Then $g_i \mapsto h_i$ $(i = 1, \ldots, n)$ extends to a unique homomorphism $G \to H$.  The following theorem is instructive,

\begin{theorem}[Dyck's theorem]
Let $\phi: F \to G$ and $\phi': F \to H$ denote homomorphisms from the free group onto groups $G$ and $H$ with $\ker(\phi) \le \ker(\phi')$.  Then there is a homomorphism $\phi'': G \to H$ such that $\ker(\phi'') = \ker(\phi') / \ker(\phi)$.  
\end{theorem}

\bigskip

The first presentation of $S_{m,n}$ is given in terms of switchings.  To build an intuition for what the relations satisfied look like for a general $m$ and $n$, we work out the simplest $m=3=n$ and $m=3, n=4$ cases using
$\mathsf{GAP}$ 4
\index{gap@$\mathsf{GAP}$}%
to verify our results.

As we saw in a previous chapter, the group $S_{3,3}$ is $4$-generated
by $\sigma_{i,j,1}$, $\sigma_{j,k,1}$, $\sigma_{i,j,2}$,
$\sigma_{j,k,2}$.  It has four types of relations: 

(i) those corresponding to the generators being involutions, for example
$\sigma_{i,j,1}^{2} = 1$,

(ii) Klein $4$-group `row' relations of the form $\sigma_{i,j,1}
\sigma_{i,j,2} \sigma_{i,j,1} \sigma_{i,j,2} = 1$,
\index{group ! Klein}%

(iii) $\Sym(3)$-group `column' relations of the form $\sigma_{i,j,1}
\sigma_{j,k,1} \sigma_{i,j,1} = \sigma_{j,k,1} \sigma_{i,j,1}
\sigma_{j,k,1}$,

(iv) two extra `diagonal' relations required to get a finite group from our starting
infinite free group,
$(\sigma_{i,j,1} \sigma_{j,k,2})^2 = (\sigma_{i,j,1} \sigma_{j,k,1}
\sigma_{i,j,1} \sigma_{i,j,2})^2$, and $(\sigma_{i,j,1}
\sigma_{j,k,2})^2 = (\sigma_{j,k,1} \sigma_{i,j,2} \sigma_{j,k,2}
\sigma_{i,j,2})^2$.  The origin of these relations can be demystified 
by noticing that $(\sigma_{i,j,1} \sigma_{j,k,2})^2$ is a $3$-cycle
$(i j k)$ on edge $\{1,2\}$ and the identity elsewhere, and that
$(\sigma_{i,j,1} \sigma_{j,k,1} \sigma_{i,j,1} \sigma_{i,j,2})^2 =
(\sigma_{i,k,1} \sigma_{i,j,2})^2$ is the same.  Similarly for the
second relation.

Then starting with a free group on these four generators, we can
factor out the relations to get a finitely presented group of the same
order as $S_{3,3}$, namely $108$, as follows:

With \texttt{f.1, f.2, f.3, f.4} replacing  $\sigma_{i,j,1}, \sigma_{j,k,1},
\sigma_{i,j,2}, \sigma_{j,k,2}$ respectively this translates to
$\mathsf{GAP}$ 4 commands as follows:\\
\texttt{ gap> f := FreeGroup("R1", "R2", "R3", "R4");; }\\
\texttt{ gap> g := f / [ f.1*f.1, f.2*f.2, f.3*f.3, f.4*f.4, }\\ 
\texttt{\                f.1*f.2*f.1*f.2*f.1*f.2, f.1*f.3*f.1*f.3, }\\
\texttt{\                f.3*f.4*f.3*f.4*f.3*f.4, f.2*f.4*f.2*f.4, }\\
\texttt{\                f.1*f.4*f.1*f.4*f.3*f.1*f.2*f.1*f.3*f.1*f.2*f.1, }\\ 
\texttt{\                f.1*f.4*f.1*f.4*f.3*f.4*f.3*f.2*f.3*f.4*f.3*f.2 ]; }\\
\texttt{ gap> Size(f); infinity }\\
\texttt{ gap> Size(g); 108 }\\


That the relations in (iv) above are satisfied in $S_{3,3}$ is easily
checked in $\mathsf{GAP}$ 4 through\\
\texttt{ gap> R1*R4*R1*R4=R1*R2*R1*R3*R1*R2*R1*R3; true  }\\
\texttt{ gap> R1*R4*R1*R4=R2*R3*R4*R3*R2*R3*R4*R3; true  }

We need to uncover the structure of the group \texttt{g}.\\

That $\texttt{g}$ obeys the same relations as the group $\texttt{G}$
from $\S7$ of Chapter~\ref{chap2} means that by Dyck's theorem there is a homomorphism from one to the other.  Then given that the two groups are of the same order, this is actually an isomorphism.  So $\texttt{g} \cong S_{3,3}$.

\smallskip

The next case to consider is that of $4$ vertices and
$3$ edge-colours.  The non-trivial switchings of colours $c, d, e$ about vertices $1, 2$ or $3$ are given by
\begin{figure}[!h]
$$\xymatrix{
{\sigma_{c,d,1}\ (\texttt{f.1})} & {\sigma_{c,d,2}\ (\texttt{f.3})} & {\sigma_{c,d,3}\ (\texttt{f.5})} \\
{\sigma_{d,e,1}\ (\texttt{f.2})} &
{\sigma_{d,e,2}\ (\texttt{f.4})} & {\sigma_{d,e,3}\ (\texttt{f.6})}
}$$
\caption{Generators of $S_{3,4}$}  
\end{figure}

Here $\sigma_{c,d,1}$ means a switch of colours $c$ and $d$ whenever
they occur on an edge about vertex $1$.  With \texttt{f.1, \ldots, f.6} replacing  $\sigma_{c,d,1}, \ldots,
\sigma_{d,e,3}$, this translates to $\mathsf{GAP}$ 4 commands that are similar to those above for $S_{3,3}$ and lead to:\\
\texttt{ gap> Size(g1); 5832 }\\

This order of $\texttt{g1}$ is the same as that of the group $S_{3,4}$.  It can be checked using $\mathsf{GAP}$ 4 that the generators of \texttt{g1} satisfy the same relations as those of $S_{3,4}$, so there is a homomorphism from \texttt{g1} to $S_{3,4}$.   Therefore the groups are isomorphic.  The 3 column relations and 6 row relations are followed by 12 `diagonal' relations which are motivated for example by the following diagram:\begin{figure}[!h]
$$
\xymatrix{
{1} \ar@{-}[dr]^{(cd)} \ar@{-}[dd]_{(cd)(de)=(ced)} &&&& {1} \ar@{-}[dr]^{(ce)} \ar@{-}[dd]_{(ce)(cd)=(ced)}\\
& {\Gamma \backslash \{1, 2\}} &&&& {\Gamma \backslash \{1, 2\}}\\
{2} \ar@{-}[ur]_{(de)} &&&& {2} \ar@{-}[ur]_{(cd)}
}    
$$
\caption{Motivating the diagonal relations of $S_{m,n}$}  
\end{figure}
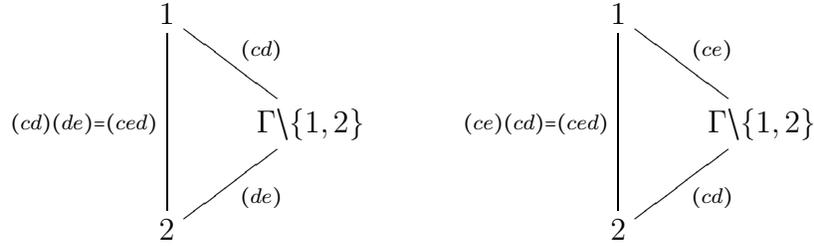

This gives $(\sigma_{c,d,1} \sigma_{d,e,2})^2 = (\sigma_{c,e,1}
\sigma_{c,d,2})^2$ which becomes 

$(\texttt{f.1*f.4})^2 = (\texttt{f.1*f.2*f.1*f.3})^{2}$.  

Other examples of composite relations are:\\
$\sigma_{c,e,1} = \texttt{f.1*f.2*f.1}$, $\sigma_{c,e,2} =
\texttt{f.3*f.4*f.3}$, $\sigma_{c,e,3} = \texttt{f.5*f.6*f.5}$,
$\sigma_{c,d,4} = \texttt{f.1*f.3*f.5}$, $\sigma_{d,e,4} =
\texttt{f.2*f.4*f.6}$.  

The two cases above having been identified, we use the Reidemeister-Schreier algorithm
\index{Reidemeister-Schreier algorithm}%
\index{Schreier, O.}%
\index{Reidemeister, K}%
to give us a presentation of $S_{m,n}$ as follows.  From the proof of Theorem~\ref{swgpfm} there is a homomorphism $S_{m,n} \to C_2^{n-1}$, the kernel of which is a group in which we identify a (right) Schreier transversal with the Schreier property.  The Schreier transversal has  $|S_{m,n} / S_{m,n-1}| = n-1$ elements in it, which we can take to be $\sigma_{c,d,i}$ for $(1 \le i \le n-1)$.

Using the Reidemeister-Schreier algorithm, we can find a presentation for this subgroup via coset representatives of it in $S_{m,n}$.  This is laid out in the table below.

\vspace{2cm}

\begin{displaymath}
\begin{array}{|c||c|c|c|}
\hline
{}&\sigma_{c,d,j}&\sigma_{d,c',j} \sigma_{c,d,i} \sigma_{d,c',j}&\sigma_{c',d',i} \sigma_{c,d,i} \\
{}&{(1 \le j \le n)}&(1 <  j \le n)&(1 \le i \le n)\\
{}&{(i \ne j)}&(\forall c \ne d \ne c')&(\{c,d\} \ne \{c',d'\})\\
\hline
\hline
\sigma_{c,d,i}&\sigma_{c,d,i} \cdot \sigma_{c,d,j}&\sigma_{c,d,i} \cdot (\sigma_{d,c',j}\sigma_{c,d,i}.\sigma_{d,c',j})&\sigma_{c,d,i} \cdot (\sigma_{c',d',i} \cdot\sigma_{c,d,i}) \\
{}&{\cdot (\sigma_{c,d,i} \sigma_{c,d,j})}&\cdot (\sigma_{c,d,j}\sigma_{c,c',i})^2&\cdot (\sigma_{c',d',i}\sigma_{c,d,i} \sigma_{c',d',i})\\
\hline
\end{array}
\end{displaymath}

\vspace{2cm}

Finally we arrive at:

$S_{m,n} = \langle \sigma_{c,d,i}, 1 \le i \le n, 1 \le c < d \le m,\ |\
\sigma_{c,d,i}^2 = 1, \forall c, d, i,\\ (\sigma_{c,d,i}
\sigma_{c,d,j})^2=1, \forall c, d, (i \ne j),\ (\sigma_{c,d,i}
\sigma_{c',d',i})^3=1, \forall i\ \forall \{c,d\} \ne \{c',d'\},\\
(\sigma_{c,d,i} \sigma_{d,c',j})^2=(\sigma_{c,c',i} \sigma_{c,d,j})^2
\ ( = (\sigma_{c,d,i} \sigma_{d,c',j} \sigma_{c,d,i}
\sigma_{c,d,j})^2),\ (\forall c \ne d \ne c')\\ (\forall i<j) \rangle$.

\bigskip

The \emph{$m$-strand braid group}~\label{braidgroup}
\index{group ! mstrand@$m$-strand braid}%
is defined by $B_m := \langle b_1, \ldots, b_{m-1}\ |\ b_i b_{i+1} b_i\\ = b_{i+1} b_i
b_{i+1},\ i=1, 2, \ldots, m-2,\ b_i b_j = b_j b_i,\ |i-j| \ge 2
\rangle$.  The braid group has the symmetric group as a quotient.  There is a homomorphism from the direct product of two braid groups into
$S_{3,3}$ via symmetric groups.  To see this write the generators and relations of the $3$-braid group suggestively as
$B_3 \times B_3 = \langle b_1^{(1)}, b_2^{(1)} \rangle \times \langle
b_1^{(2)}, b_2^{(2)} \rangle = \langle a, b \rangle \times \langle
c, d \rangle$ $$\xymatrix{
& a\ (b_1^{(1)}) & c\ (b_1^{(2)}) \\
& b\ (b_2^{(1)}) & d\ (b_2^{(2)}) 
}$$ satisfying $aba = bab, cdc=dcd, ac=ca, bd=db, ad=da, bc=cb$
and that of the switching group as 
$S_{3,3} =$ $$\xymatrix{
& A\ ({\sigma_{i,j,1}}) & C\ ({\sigma_{i,j,2}}) \\
& B\ ({\sigma_{j,k,1}}) & D\ ({\sigma_{j,k,2}}) 
}$$ satisfying $A^2=1$, $B^2=1$, $C^2=1$, $D^2=1$, $ABA = BAB,
CDC=DCD, AC=CA, BD=DB, (ABAC)^2=(AD)^2, (BCDC)^2=(AD)^2$.

To show that there is a homomorphism from $B_3 \times B_3 \to S_{3,3}$ we need to verify that the $S_{3,3}$ relations imply the $B_3 \times B_3$ relations.  This is readily confirmed using $\mathsf{GAP}$ 4.

Consider the transformations:\ $a \to ADA,\ d \to DAD,\ b \to BCB,\\ c \to
CBC$.  Using $\mathsf{GAP}$ 4 we can verify that $ADA, DAD, BCB, CBC$
generate a group of order $36$ using the defining relations of
$S_{3,3}$ we gave in Chapter~\ref{chap2} above.  In the $\mathsf{GAP}$ 4 code that follows we denote by $\texttt{Q1}$ the homomorphic image of the mapping $B_3 \times B_3 \to S_{3,3}$, and we have to show that $\texttt{Q1} = \Sym(3) \times \Sym(3)$.

\texttt{ gap> S := SymmetricGroup(27); }

\texttt{ gap> P1 := (1,7)(2,4)(3,8)(5,6)(9,18)(10,25)(11,20) }\\
\texttt{\                (12,22)(13,24)(16,19)(17,26)(21,27); }\\
\texttt{ gap> P2 := (1,14)(2,20)(3,23)(4,18)(5,26)(6,27)(9,11)}\\
\texttt{\                 (10,15)(12,13)(16,24)(17,21)(19,22); }\\
\texttt{ gap> P3 := (1,6)(2,8)(3,4)(5,7)(9,27)(10,19)(11,17)}\\
\texttt{\                  (12,22)(14,23)(16,25)(18,21)(20,26); }\\
\texttt{ gap> P4 := (1,13)(2,24)(3,17)(4,19)(5,26)(7,25)(9,15) }\\
\texttt{\                 (10,11)(12,14)(16,20)(18,22)(21,23); }\\
\texttt{ gap> G := Subgroup(S,[P1,P2,P3,P4]); }\\
\texttt{ gap> gen1 := P1*P4*P1; (1,10)(2,16)(4,13)(6,17)(7,24)
}\\
\texttt{\                 (8,26)(9,12)(11,19)(14,22)(15,18)(20,25)(23,27) }\\
\texttt{ gap> gen2 := P4*P1*P4; (1,2)(3,5)(4,20)(6,26)(7,11)(8,17) }\\
\texttt{\                (10,16)(13,25)(14,18)(15,22)(19,24)(23,27) }\\
\texttt{ gap> gen3 := P2*P3*P2; (1,3)(2,5)(4,17)(6,11)(7,26)(8,20) }\\
\texttt{\                (9,21)(13,19)(14,27)(15,22)(18,23)(24,25) }\\
\texttt{ gap> gen4 := P3*P2*P3; (1,9)(3,21)(4,14)(6,23)(7,20) }\\
\texttt{\                (8,26)(10,12)(11,18)(13,22)(15,19)(17,27)(24,25)}\\
\texttt{ gap> Q1 := Subgroup(S,[gen1, gen2, gen3, gen4]);  }\\
\texttt{ gap> Size(Q1); 36 }\\
\texttt{ gap> IsSubgroup(G,Q1); true }\\
\texttt{ gap> IsNormal(G,Q1); false }\\

So certainly $\texttt{Q1}$ has the correct order, and it remains to check that it has the required structure.  We can do that in  $\mathsf{GAP}$ 4 through commands such as\\
\texttt{ gap> gen1*gen3*gen1=gen3*gen1*gen3; true  }\\
\texttt{ gap> gen1*gen2=gen2*gen1; true  }

The homomorphism $B_3 \times B_3 \to S_{3,3}$ is not a surjection.  Rather we have that $B_3 \times
B_3 \onto \Sym(3) \times \Sym(3) \hookrightarrow S_{3,3}$.  This says
that in a sense a $3$-switching is the commuting pair of two
$3$-braidings acting independently on two different triples of strings, modulo a $3$-cycle.  To see that this generalizes to $m, n > 3$ consider the following graph
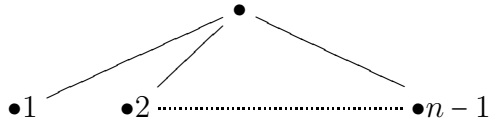
\begin{figure}[!ht]
$$\xymatrix{ 
&& {\bullet} \ar@{-}[dll] \ar@{-}[dl]
\ar@{-}[drr] \\
{\bullet}1  & {\bullet}2 \ar@{.}[rrr] &&& {\bullet}n-1
}$$
\end{figure}

This represents the product of $n-1$ copies of $\Alt(m)$, one on each edge.  Switchings at the end vertices generate $C_2^{n-1}$.  One switching at the apex vertex acts on each of a direct product of alternating groups.  So there is an embedding of the direct product $\Dr_{1 \le i \le n-1} \Sym(m)^i$ of groups into $S_{m,n}$.

\smallskip

That switching and hyperoctahedral groups
\index{group ! hyperoctahedral}%
have the same direct factors, but complementary normal subgroups and that switching and 
braid groups
\index{group ! braid}%
 are homomorphically related, suggests the possible existence of representations for the switching groups related to those of type A or B Hecke algebras;
\index{Hecke algebra}%
 see Chapter~\ref{chapFD}.

\section{Appendix 2: A Second Presentation of $S_{m,n}$}

A second type of $S_{m,n}$ presentation is given in terms of
involutions about vertices and $3$-cycles about edges, and can be recovered
by inspection of the form of the switching groups.  For example $S_{3,3} \cong
(\Alt(3))^{3}\sd (C_2)^{2}$ consists of $3$ generators $a_1, a_2, a_3$ that are intended to
represent three $3$-cycles and two involutions $b_1$ and $b_2$ whose
actions are shown by the following picture
$$
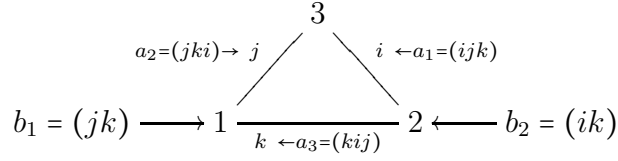
\xymatrix{ 
&& {3} \ar@{-}[dl]_{a_2=(j k i) \rightarrow \ j}  \ar@{-}[dr]^{i\ \leftarrow a_1=(i j k)}\\
{b_1=(jk)} \ar@{->}[r] & {1} \ar@{-}[rr]_{k\ \leftarrow a_3 = (k i j)} && {2} & {b_2=(ik)} \ar@{->}[l]
}$$
The presentation is given by

$S_{3,3} = \langle a_1, a_2, a_3, b_1, b_2\ |\ a_i^3 = 1,\ [a_i, a_j] =
1,\ b_i^2 = 1,\ [b_i, b_j] = 1,\\ b_1^{-1} a_1 b_1 = a_1,\ b_1^{-1} a_2 b_1
= a_2^{-1},\ b_1^{-1} a_3 b_1 = a_3^{-1},\ b_2^{-1} a_1 b_2 = a_1^{-1},\
b_2^{-1} a_2 b_2\\ = a_2,\ b_2^{-1} a_3 b_2 = a_3^{-1} \rangle$.
It is easily verified that this gives $S_{3,3}$ by\\
\texttt{ gap> f := FreeGroup("R1", "R2", "R3", "R4", "R5");; }\\
\texttt{ gap> g := f / [ f.1*f.1*f.1, f.2*f.2*f.2, f.3*f.3*f.3, }\\
\texttt{\                f.4*f.4, f.5*f.5, f.1*f.2*f.1*f.2,}\\ 
\texttt{\                f.1*f.3*f.1*f.3, f.2*f.3*f.2*f.3, f.4*f.5*f.4*f.5,}\\
\texttt{\                f.4*f.1*f.4*f.1$^-1$, f.4*f.2*f.4*f.2$^-1$, f.4*f.3*f.4*f.3$^-1$,}\\
\texttt{\                f.5*f.1*f.5*f.1, f.5*f.2*f.5*f.2, f.5*f.3*f.5*f.3 ];}\\
\texttt{ gap> Size(g); 108 }\\

We can use $\mathsf{GAP}$ 4 to show that this $S_{3,3}$ presentation is isomorphic to the one in the previous section.

This presentation generalizes to $n > 3$.  First for convenience we change notation, for example letting $a_1$ above become $a_{23}$, etc.

Extrapolating from $3$ to a finite number $n$ of vertices, we can give
generators for $S_{3,n}$ as $b_1, \ldots, b_{n-1}, a_{ij}, 1 \leq i <
j \leq n$, where the $b_i$ are transpositions of $2$ colours on edges
$i, j$ for all $j \ne i$ and the $a_{ij}$ are $3$-cycles on colours on
edge $\{i, j\}$.  The relations are:

$a_{ij}^{3}=1$,\quad $[a_{ij},a_{kl}]=1$ $(i,j \neq k, l)$,\quad $b_i^{2}=1$,\quad $[b_i,b_{j}]=1$ $(i \ne j)$,
\begin{displaymath}
b_i^{-1} a_{jk} b_i = \left\{ \begin{array}{ll}
a_{jk}^{-1} & \textrm{if $j=i$ or $k=i$}\\
a_{jk} & \textrm{if $i \notin \{j,k\}$.}
\end{array} \right.
\end{displaymath}

Certainly $S_{3,n}$ satisfies this presentation.  But we must show
that this gives a defining set of relations for $S_{3,n}$.  If
$G$ is a group defined by this presentation then $|G| \leq |S_{3,n}|$
and $G \onto S_{3,n}$.  Therefore we have equality and a verification
that this is a presentation for $S_{3,n}$.

\bigskip

We can extrapolating further from $m=3$ to any finite number $m$ of colours, and
use the presentation

\centerline{$\Alt(m)= \langle a_1, \ldots, a_{m-2} |  a_k^{3}=1,\ (a_k a_l)^2=1\ (k \neq l) \rangle$.}

This gives the generators for $S_{m,n}$ as: $b_i$ where $b_i$ is a $2$-cycle $(m-1, m)$ on colours on all edges $\{i, j\}$ for all $1 \le i < j \le n$, and $a_{ij}^{k}$ for $1
\le k \le m-2$ which are $3$-cycles $(k, m-1, m)$ of colours on all edges $\{i, j\}$.  Then
$(a_{ij}^{k} a_{ij}^{l})^2 = ( (k, m-1, m)(l, m-1, m) )^2
= ( (k\ m)(l\ m-1) )^2 = 1$.  The relations are:

$(a_{ij}^{k})^{3}=1,\quad (a_{ij}^{k}\ a_{ij}^{l})^{2} = 1\ (l \neq k),\quad [a_{ij}^{k},\ a_{i'j'}^{l}] = 1\ (\text{if}\ \{i',j'\} \neq \{i,j\}),\quad b_i^{2}=1,\ [b_i, b_{i'}] = 1\ (\text{if}\ i \ne i')$,
\begin{displaymath}
b_h^{-1} a_{ij}^{k} b_h = \left\{ \begin{array}{ll}
(a_{ij}^{k})^{-1} & \textrm{if $h=i$ or $h=j$}\\
a_{ij}^{k} & \textrm{if $h \notin \{i,j\}$.}
 \end{array} \right.
\end{displaymath}

Each edge $\{i, j\}$ represents one alternating group in the direct product expansion of $S_{m,n}$.  To prove that this presentation defines precisely $S_{m,n}$, first recall the presentation 

$\Sym(m) = \langle a_{ij}^1, \ldots, a_{ij}^{m-2},\ b_{i} : (a_{ij}^k)^3 = 1 = (a_{ij}^k a_{ij}^l)^2\ (i \ne j),\\ b_{i}^2 = 1,\ b_{i} a_{ij}^k b_{i} = (a_{ij}^k)^{-1} \rangle$ for fixed $i, j$ and for all
$k$.  

Making the substitutions $a_{ij}^{k} \mapsto (k\ m-1\
m),\ b_{i} \mapsto (m-1\ m)$ we see that the $a_{ij}^{k}$ generate $\Alt(m)$. Then
proceeding as for $S_{3,n}$ proves the presentation.






\chapter{Some Ring Theory around Random Graphs}
~\label{ringchap}

\bigskip

What is now proved was once only imagined
\begin{flushright}
William Blake, \textit{The Marriage of Heaven and Hell, (1793)} 
\end{flushright}

\bigskip

This chapter takes as its starting point the work of Cameron
\index{Cameron, P. J.}%
 on the algebra of $\mathbb{Q}$-valued functions on $n$-subsets of an infinite set $\Omega$, as outlined in Appendix~\ref{Polynomials}.  The algebra carries information on group actions
\index{group ! action}%
 on finite subsets of $\Omega$. We will uncover further results on the structure of this algebra in specific cases.

\section{Polynomial Algebras for Random Graphs}

In this first section we build on the application of this theory of permutation group algebras to the case of multicoloured random graphs~\cite{cam7}~\cite{cam8}.
\index{graph ! random ! $m$-coloured}%

\medskip

Let $\mathfrak{R}_{m,n} \in \mathcal{G}_{m,n}$ denote the countable
$m$-coloured random graph,
\index{graph ! random ! $m$-coloured}%
with $n$ vertices and $m$ colours in the set of all such simple graphs.  Fix $m$ and let $\mathcal{G}$ be the class of all
finite $m$-coloured graphs, $V_{n}$ the vector space of functions from
the isomorphism types of $n$-element structures in $\mathcal{G}$ to
$\mathbb{Q}$, and  $\mathcal{A} := \oplus_{n \ge 0}
V_{n}$, the \emph{age}
\index{age}%
 of $\mathfrak{R}_{m,n}$ being the class of
all finite structures embeddable in $\mathfrak{R}_{m,n}$ as induced
substructures.  We will determine the structure of
$\mathcal{A} = \mathcal{A}^{Aut(\mathfrak{R}_{m,n})}$, working with
a countably infinite set, $|n| = \aleph_{0}$.

Consider the following two conditions:

(C1) If a structure $S$ is partitioned into disjoint induced
substructures $S_{1}, S_{2}, \ldots$ then $S_{1} \circ S_{2}
\circ \ldots \le S$.  (Here $\circ$ is a binary, commutative and associative law of
composition, such that if $X$ and $Y$ are finite subgraphs, $X \circ Y$ is a subgraph with $|X| + |Y|$ vertices.  The composition referred to here is that of structures that are not themselves compositions.  Also the partial order $\leq$ indicates ``involvement'' on the class of $n$-element structures for each $n$).

(C2)  Any structure has a unique representation as a composition of
connected structures.  We assume $\mathcal{G}$ to be a class of finite
relational subgraphs closed under isomorphism and taking induced
subgraphs.

If an age $\mathcal{A}$ has an isomorphism-closed subclass of ``connected'' structures, a ``composition'', and a partial order on $n$-element structures in the form of ``involvement'' such that (C1) and (C2) hold, then the algebra of $\mathcal{A}$ is freely generated by the characteristic functions of ``connected'' structures.  The words in quotes do not have their usual meaning here.

If we take the age
\index{age}%
 consisting of all finite graphs, let ``connected'' have its usual meaning, ``involvement'' be spanning subgraph,
\index{graph ! spanning}%
 and ``composition'' be disjoint union, then (C1) and (C2) hold.

More concretely our generalization is to take finite $3$-coloured complete graphs, where:

$A$ is a ``subgraph'' of $B$ if by changing some of the $\mathfrak{b}$ or
$\mathfrak{r}$-adjacencies of $B$ to $\mathfrak{g}$, that is `throwing away' some $\mathfrak{b}$ or $\mathfrak{r}$ adjacencies, we get $A$;  ``connected'' means $\mathfrak{r} \mathfrak{b}$-connected, that is any two vertices are joined by a path using only red and blue adjacencies, taking the green colour $\mathfrak{g}$ to represent a non-adjacency;  ``involvement'' meaning changing some green edges to red or blue;  ``composition'' meaning a disjoint union of finite $\mathfrak{r}\mathfrak{b}\mathfrak{g}$-subgraphs, with only $\mathfrak{g}$ edges between the parts.  If we adopt these meanings then again conditions (C1) and (C2) hold.

We proceed as in~\cite[Theorem 2.1]{cam7}.

\begin{theorem}
\label{polalgcon}
Assuming (C1) and (C2), $\mathcal{A}$ is a polynomial
algebra generated by the characteristic functions of the connected
structures.
\end{theorem}

\begin{proof}
Let $S = \bigsqcup S_{i}$, where $|S| = n$ the $S_{i}$ are connected
structures and the union is disjoint.
Let $V_{n}$ and $\mathcal{A} := \oplus_{n \ge 0}
V_{n}$, be as defined above.  By definition the
characteristic functions $\chi_{S}$ for $|S|=n$ form a basis for the
vector space $V_{n}$.  There is a bijection between
$\chi_{S}$ and the monomials $\phi_{S} = \chi_{S_{1}} \chi_{S_{2}}
\ldots$ of total weight $n$, where the $S_{i}$ are connected and
$S_{1} \circ S_{2} \circ \ldots = S$.  The matrix expressing $\phi_{S}$ in
terms of the basis vectors $\chi_{S'}$, has non-zero coefficient for
$\chi_{S}$ in the row corresponding to $\phi_{S}$.  If $\chi_{S'}$
also has non-zero coefficients then $S'$ has a partition into induced
subgraphs isomorphic to the $S_{i}$.  Thus $S = \bigsqcup S_{i} \le
S'$.  The matrix corresponding to this is upper triangular where the ordering extends the partial ordering of involvement, with non-zero diagonal; so it is invertible.  So the monomials of weight
$n$ themselves form a basis for $V_{n}$.
\end{proof}

To clarify the origin of the form of the matrix note
that there is a partial order inherent in the choice of green edge
regarded as a non-adjacency:
$$\xymatrix{ {\bullet} \ar@{-}[r]^{\mathfrak{g}} & {\bullet}}$$
followed by either a red edge or a blue edge
$$\xymatrix{ {\bullet} \ar@{-}[r]^{\mathfrak{r}} & {\bullet} & {or} &
  {\bullet} \ar@{-}[r]^{\mathfrak{b}} & {\bullet} }$$    
That we can extend this partial order on $n$-element structures to a total
order implies upper triangular transition matrices between bases, where the partial order on
the coloured edges is induced onto the $\chi$s.

The theorem clearly holds for graphs on any fixed number $m$ of colours.

\begin{corollary}
For $H= \Aut(\mathfrak{R^{t}})$, $\mathcal{A}^{H}$ is a
polynomial algebra,
\index{polynomial algebra}%
and therefore an integral domain.
\index{integral domain}%
\end{corollary}

\begin{proof}
With the identifications as above, applying the theorem gives that
$\mathcal{A}(\mathcal{G}) =
\mathcal{A}^{\Aut(\mathfrak{R^{t}})}$ is a polynomial
algebra, whose generators correspond to the finite connected graphs.
\end{proof}
 
\begin{corollary}
For $H= \Aut(\mathfrak{R}_{m,\omega})$, $4 \le m \le \aleph_{0}$,
$\mathcal{A}^{H}$ is a polynomial algebra.
\index{polynomial algebra}%
\end{corollary}

\begin{proof}
The same method as above works for $m$ colours if we take connected to
mean connected in the first $m-1$ colours and take the final colour
$m$ to represent non-adjacency.
\end{proof}

Because $S_{m, \omega}$ is highly transitive
\index{group ! permutation ! highly transitive}%
for $m \ge 3$, $\mathcal{A}^{S_{m, \omega}}$ is a polynomial ring in
one variable.  Another way to see this is as a corollary of the theorem.  

\begin{corollary}
For $H= S_{m, \omega}, D_{m, \omega}\ or\ B_{m, \omega}$, $4 \le m \le \aleph_{0}$,
$\mathcal{A}^{H}$ is an integral domain.
\end{corollary}

\begin{proof}
This is clear since $G_{1} \le G_{2}$ implies that $\mathcal{A}^{G_{2}} \le \mathcal{A}^{G_{1}}$ and a subring of an integral domain is an integral domain.  Finally of course a polynomial ring is an integral domain.
\end{proof}

Here $D_{m, \omega}$ is understood to mean take any
subgroup $H \le \Sym(m)$ and take all permutations of the vertices
which induce a permutation from $H$ on the colours, whilst for $S_{m, \omega}$ we
take any partition of the set of colours and allow switching of colours between the parts.

If $\mathcal{A}^{G}$ is an integral domain, we say that $G$ is \emph{entire}.
\index{group ! entire}%

Among the $m$ distinct colours of the $m$-coloured random graph,
\index{graph ! random ! $m$-coloured}%
$\mathfrak{R}_{m, \omega}$ ($m \ge 2$), choose one distinguished
``transparent'' colour $c_{0}$.  Taking connectedness to mean
connected in the other $m-1$ colours yields, as we have seen, a
polynomial algebra,
\index{polynomial algebra}%
 the stabilized colour being the
``non-adjacency'' colour.  As the action on the colour set
$\mathcal{C}_m$ is transitive, $\Sym(m-1)$ stabilizes $c_{0}$ and permutes the other colours forming a single conjugacy class of subgroups of $\Sym(m)$.  This translates into the fact that $\Sym(m-1)$ acts on
the polynomial algebra so as to permute the polynomial generators, the
monomials $\phi_{S}$.  (The universality of the graph
$\mathfrak{R}_{\omega, \omega}$ gives rise to a universal polynomial ring
associated to it).  So there is an invariant subalgebra generated by
symmetric polynomials, but is it a polynomial algebra?

Furthermore consider the graph $\mathfrak{R}_{m, \omega}$ on $m$
colours $c_0, \ldots, c_{m-1}$ with $c_0$ transparent and $ G =
\Aut(\mathfrak{R}_{m, \omega})$.  Then $H \le \Sym(m-1)$
leads to $ D^H_{m, \omega}$ the group of all
permutations of vertices which fix $c_0$ and induce an element of $H$
on the remaining colours.  This gives rise to the \emph{Galois
  Correspondence}
 \index{Galois ! correspondence}%
\[ H \le \Sym(m-1) \leftrightarrow \mathcal{A}^{D^H_{m,
    \omega}}. \]
For all subgroups of $\Sym(m)$ of degree $m$ the algebra is an
integral domain because in this case the automorphism groups
\index{group ! automorphism}%
 are reducts. 
The correspondence arises because the larger the group the smaller the algebra.
 If $H = 1$ then $G = D^{\{1\}}_{m, \omega}
\leftrightarrow \mathcal{A}^{G}$ where $\mathcal{A}^{G}$ is a polynomial ring generated by
connected graphs.  If $H = \Sym(m-1)$ then $\mathcal{A}^{{D}^{H}_{m, \omega}}$ is a polynomial ring if $m=3$ but we conjecture that it is not if $m > 3$.  

\head{Open Question}  Is it true that $\mathcal{A}^{{D}^{H}_{m, \omega}}$ is a polynomial
ring $\Leftrightarrow H$ is a Young subgroup
\index{group ! Young subgroup}%
of $\Sym(m-1)$ corresponding to a partition with parts of size $1$ except at most one
of size $2$.

This would be false without the condition on the partitions, for
$\Sym(m-1)$ itself is a Young subgroup.  Further, even a partition of
four colours into two pairs, say $\{c_1, c_2\}$ $\{c_3, c_4\}$
would allow a regular orbit under the Klein $4$-group,
\index{group ! Klein}%
giving the orbit
$$
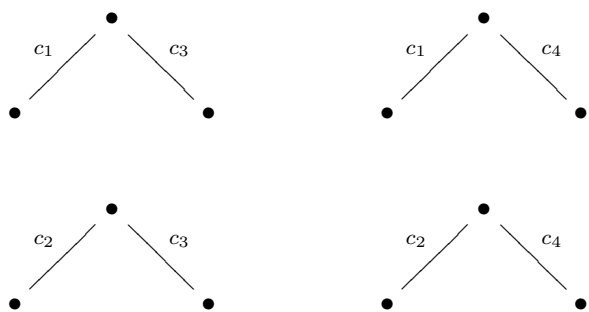
\xymatrix{
& {\bullet} \ar@{-}[dr]^{c_3} \ar@{-}[dl]_{c_1} &&&& {\bullet} \ar@{-}[dl]_{c_1} \ar@{-}[dr]^{c_4}\\ 
{\bullet} && {\bullet} && {\bullet} && {\bullet}\\
& {\bullet} \ar@{-}[dr]^{c_3} \ar@{-}[dl]_{c_2} &&&& {\bullet} \ar@{-}[dl]_{c_2} \ar@{-}[dr]^{c_4}\\ 
{\bullet} && {\bullet} && {\bullet} && {\bullet}
}$$
where the transparent edge is $c_0$.  However the action of the group
$C_2$ is always regular, so the open question is probably true as formulated.

More generally we could ask why do we get a polynomial ring if the orbit is \emph{regular}?  If $G$ is a finite group acting on $\mathbb{C}[X]$, where $X$ is any set of indeterminates, and there is a $G$-orbit $Y$ such that $G \neq \Sym(Y)$ then is $\mathbb{C}[X]^G$ a polynomial ring?

We can get a stronger result as follows.  Imagine $\Sym(n)$ to act on the
ordered partition $S = (S_{1}, \ldots, S_{j})$ a natural way, that is as $S^{g} = (S_{1}^{g}, \ldots, S_{j}^{g})$.  Then the stabilizer of an element of the partition contains a
\emph{Young subgroup}
\index{group ! Young subgroup}%
of $\Sym(n)$ which is isomorphic to
$G := \Sym(n_{1}) \times \Sym(n_{2}) \times \dots \times \Sym(n_{j})$,
where $|S_{i}| = n_{i}$ and $n = n_{1} + \ldots + n_{j}$, $n_{1} \ge
\ldots \ge n_{j} > 0$ is a partition of $n$.  Now it is known that if $G$ is a Young subgroup
of $\Sym(n)$ acting on a polynomial ring in $n$ variables, then the ring of invariants of $G$ is a polynomial ring in the elementary symmetric functions of degree $1, \ldots, n$.  So here, $\mathcal{A}^{\Sym(n_{1})} \otimes \ldots \otimes \mathcal{A}^{\Sym(n_{j})}$ is isomorphic to a ring of invariants on $n_{i}$ homogeneous generators of degrees $1, \ldots, n_{i}$ which is a tensor product of invariants of the $n_{i}$ orbits, for $1 \le i \le j$~\cite[p.~79]{cam6}.

From the character theory of the symmetric group~\cite{macd}, if a
partition~\label{partition} $\pi = (\pi_1, \pi_2,
\ldots)$ determines a monomial $x^{\pi} =
x^{\pi_1}x^{\pi_2} \ldots$, then the \emph{monomial symmetric
  function} $m_{\pi}$
\index{symmetric function ! monomial}%
is the sum of all distinct monomials
obtainable from $x^{\pi}$ by permutations of the $x's$, for
example, $m_{(21)} = \sum x_i^{2}x_j$ summed over all pairs $(ij)$ such that
$i \neq j$. When $\pi = (1^r)$, $m_{\pi}$ is the \emph{r}th
\emph{elementary symmetric function}
\index{symmetric function ! elementary}%
\[e_r = m_{(1^r)} = \sum_{i_1 <
  \ldots < i_r} x_{i_1} \ldots x_{i_r}\] with $e_0 = 1$.  The equations
$e_{\pi'} = m_{\pi} + \sum_{\mu < \pi} a_{\pi\mu} m_{\mu}$ with nonnegative integer coefficients $a_{\pi\mu}$ are
of triangular form so can be solved for $m_{\pi}$ in terms of the
$e_{\mu}$.  Now the $\chi$ and $\phi$ functions in Theorem~\ref{polalgcon} form two different bases
for the algebra of symmetric functions, the transition matrix between
these being triangular.  By analogy, writing our polynomials homogeneously, sums of
the $\chi_{S}$ functions could be said to be \emph{elementary
  symmetric functions}, with the monomials $\phi_{S}$ forming the
so-called \emph{monomial symmetric functions}, and the upper
triangular matrix in the above theorem being the character table
\index{character ! table}%
of $\Sym(n)$.  Inverting upper triangular matrices is done with the aid
of the M\"obius function~\cite{cameron1}.
\index{mobius@M\"obius function}%

We note two points: (1) the age
\index{age}%
 may be an integral domain,
without being a polynomial algebra.
\index{polynomial algebra}%
  (2)  If it is a polynomial
algebra, we are expecting that $H$
be a Weyl group
\index{group ! Weyl}%
of some type, because only for the symmetric group in
its natural action do the invariants form a polynomial algebra (which
is generated by symmetric functions).

Consider the induced action of $G = \Sym(m-1)$ on the
field of fractions $L$ of $\mathcal{A}^{G}$, which is isomorphic to
$\mathbb{Q}(x_{1}, \ldots, x_{n})$ a number field
\index{number field}%
which is a purely
transcendental extension of $\mathbb{Q}$ of transcendence degree $n$,
with $G$ being the group of field automorphisms.  Here $G$ is a Galois
group
\index{group ! Galois}%
and $L$ is a finite Galois extension
\index{Galois ! extension}%
 of the fixed field $L^{G}$
which also has transcendental degree $n$ over $\mathbb{Q}$ because
the $G$-action preserves the grading of $\mathcal{A}^{G}$ by degree.

We have the following links, by expansion of a known schema~\cite[p.~181]{sita}.

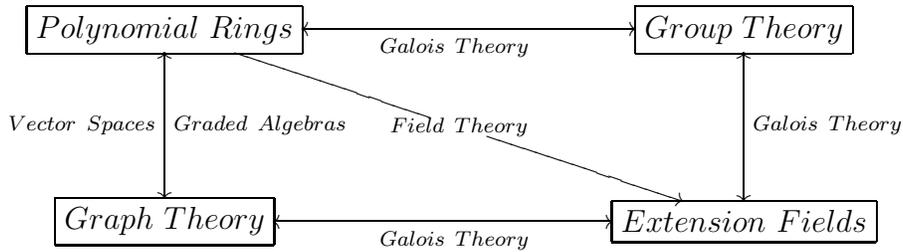
$$\xymatrix{ 
*+[F]{Polynomial\ Rings}  \ar[rrrr]_{Galois\ Theory} \ar[dd]_{Vector\
Spaces} \ar[ddrrrr]|-{Field\ Theory}
&&&& *+[F]{Group\ Theory} \ar[llll] \ar[dd]^{Galois\ Theory}\\
&&\\
*+[F]{Graph\ Theory} \ar[rrrr] \ar[uu]_{Graded\ Algebras} &&&&
*+[F]{Extension\ Fields} \ar[llll]^{Galois\ Theory} \ar[uu]\\
}$$

The arrows in this diagram are not inteded to signify
functors but our expectation of interesting connections between the
fields via theorems extending the results of this section, resulting in a dictionary which would enable us
to translate the terminology of one area into those of another.

Recall that $K$ is a number field if the dimension $[K : \mathbb{Q}]$
of $K$ as a vector space over $\mathbb{Q}$ is finite, and $K$ is a
Galois extension of $\mathbb{Q}$ giving a Galois group of $K/\mathbb{Q}$ if and only if:

$|\Gal(K/\mathbb{Q})| = |\Aut(K)| = [K : \mathbb{Q}] = \sharp$ of field
automorphisms $K \hookrightarrow K$.

Fields have no non-trivial ideals so these homomorphisms are always
 $1$--$1$.  To find number fields that are Galois over $\mathbb{Q}$,
factor any $f(X) \in \mathbb{Q}[X]$ over $\mathbb{C}$, $f(X) = a
(X-\alpha_{1})(X-\alpha_{2}) \ldots (X-\alpha_{n})$, letting
$K=\mathbb{Q}(\alpha_{1}, \ldots,\alpha_{n})$ be the smallest subfield
of $\mathbb{C}$ containing these $\alpha_{i}$s.  For $\sigma: K
\hookrightarrow \mathbb{C}$, we firstly have that $\sigma(K)=K$
(showing that $K$ is Galois over $\mathbb{Q}$) and also the
$\sigma(\alpha_{i})$ are roots of $f(X)$.  The field $K$ is the
\emph{splitting\ field}
\index{splitting field}%
of $f(X)$ over $\mathbb{Q}$ if and only if it is a
Galois extension, that is $\sigma(\alpha) \in K, \forall \sigma$.
There is a Galois group action
\index{group ! action}%
 $\Gal(K/\mathbb{Q})$ on the abelian group
\index{group ! abelian}%
of $K$-rational points on a cubic curve, given by $C(K) = \{(x, y): x,y \in K,
and\ y^{2} = x^{3} + ax^{2} + bx + c\} \cup {\mathfrak{O}}$, where $\mathfrak{O}$ is the point at infinity.~\label{ptatinf}

Natural open questions arise; is there any \emph{graph-theoretical} significance to the reducts being Galois groups?  For example, is $S_{m,n}$ a Galois group over $\mathbb{Q}$ in a way that is related to its definition?  

In our work, the interesting groups are often infinite and a description of the relevant Galois theory
\index{Galois ! theory}%
appears in an appendix to~\cite{swin}.  The groups $S_{m, \omega}^{cl}$ are profinite
\index{group ! profinite}%
 and therefore Galois
\index{group ! Galois}%
 for each $m$ but they are infinite so realizing their Galois action is not straightforward.  One way to proceed in principle is to take $[\Gamma_n]$ as an equivalence class of $n$-vertex graphs on a fixed number $m \ge 3$ of colours and to define $\mathcal{S} : = \bigsqcup^{\infty}_{n=0} [\Gamma_n]$ as a ring of scalars in which the $\Gamma_n$ are formal symbols.  Then $\mathcal{S} [x, \ldots] = \bigsqcup^{\infty}_{n=0} [\Gamma_n] [x, \ldots]$ is a ring of polynomials in as many variables as we wish.  As each element of a \emph{finite subgroup} of $S_{m, \omega}^{cl}$ acts on a set of graphs, dually it permutes the roots of the polynomials.

We will stick with this theme and make some further comments to end the section.

From Theorem~\ref{strg}, for $m > 2$, the group $S_{m,n}$ acts transitively on the set
$\mathcal{G}_{m,n}$ of $m$-edge-coloured graphs on $n$ vertices, and it follows that all coloured graphs on a fixed number $n$ of vertices are equivalent from the viewpoint of switchings.  This leads to the possiblity of a `Switching Ring' $\mathcal{SR}$ whose elements are the equivalence classes of switching-equivalent graphs which we denote $[\Gamma]$.

Firstly we can form a \emph{semiring} using addition of sets of equivalence classes in the semiring to be the natural addition of commutative monoid
\index{monoid}%
 coefficients
\[ \sum_{i} \alpha_i [\Gamma_i] +  \sum_{i} \beta_i [\Gamma_i] =  \sum_{i} (\alpha_i + \beta_i) [\Gamma_i],\]
where $\Gamma_i$ are the multicoloured graph basis elements.  We use addition of vertex sets for multiplication, inducing a map $[\Gamma_i] \times [\Gamma_j] \to [\Gamma_{i + j}]$, where we do not have to define a colouring rule for the edges of the product graph because all graphs on a given number of vertices are equivalent and so $[\Gamma_{i + j}]$ is independent of the choice of representatives of $[\Gamma_{i}]$ and $[\Gamma_{j}]$.  (Other examples of semirings are: (1) $\mathbb{N}[x]$ = polynomials on a generator $x$ with natural number coefficients forming a free commutative semiring; (2) the ideals of a ring under addition and multiplication of ideals).  Each element of our semiring is a formal sum of $n$-vertex multicoloured graphs, where two graphs with different numbers of vertices are taken to be linearly independent in the basis.  

The zero coefficient gives the additive and multiplicative semiring zero and the \emph{empty graph} 
\index{graph ! empty}%
$\Gamma_{0} = (V(\Gamma)= \emptyset, E(\Gamma) = \emptyset)$ on no vertices and (therefore) no edges is defined to be the multiplicative identity.  From this follows additive associativity and commutativity, multiplicative associativity and left and right distributivity of multiplication over addition.  We further assume multiplication to be commutative giving us a commutative semiring of graph equivalence classes.

In order to define the stronger concept of rings of graphs, we turn the commutative monoid (on +) that we have thus far formed into an abelian group by allowing additive inverses of the graphs.  There are many ways to do this, for example using the inverses in the coefficient ring (or field).  This approach merely reproduces the theories of polynomial rings (or algebras) over one variable.

\bigskip
\bigskip

\begin{proposition}
The ring $\bigsqcup^{\infty}_{n=0} \mathbb{Z} [\Gamma_n] \cong \mathbb{Z} [x]$
\end{proposition}

\begin{proof}
With the homomorphism $[\Gamma_i] \mapsto x^i$ addition is clear and multiplication follows from $x^i \times x^j = x^{i+j}$.  The map reverses, giving the result.
\end{proof}

\begin{proposition}
\label{grpheqalg}
The ring $\bigsqcup^{\infty}_{n=0} \mathbb{Q} [\Gamma_n] \cong \mathbb{Q} [x]$
\end{proposition}

\begin{proof}
Defining multiplication by $[\Gamma_i] \times [\Gamma_j] =$ $i+j \choose i$ $[\Gamma_{i+j}]$ in the domain algebra.  Then the map $[\Gamma_i] \mapsto x^i / i!$ is an isomorphism because $x^i / i! \cdot x^j / j! =$ $i+j \choose i$ $x^{i+j} / {(i+j)!}$
\end{proof}

The following result can be deduced from the last proposition.

\begin{corollary}
\label{algcor}
Let $G$ be the automorphism group
\index{group ! automorphism}%
 of the countable homogeneous
switching class, and let $\mathcal{A}^{G}$ denote the algebra which
encodes information about the action of $G$ on finite subsets of
equivalence classes.  Then

(a) $G$ is highly transitive,
\index{group ! permutation ! highly transitive}%
 so $\mathcal{A}^{G} \cong \mathcal{A}^{\Sym(x)}$ 

(b) $\mathcal{A}^{G}$ is the graph equivalence class algebra with multiplication defined in Proposition \ref{grpheqalg}.
\end{corollary}

This suggests the following~\cite[p.~120]{cam6a}:

\head{Open Question}  Given that for graphs on two colours, the finite two-graphs
\index{graph ! two-graph}%
\index{two-graph}%
form a Fra\"{\i}ss\'e class,
\index{Fra\"{\i}ss\'e class}%
 denoting the automorphism group
 \index{group ! automorphism}%
  of the countable homogeneous two-graph by $\Aut(\mathcal{T}(\mathfrak{R}))$, is $\mathcal{A}^{\Aut(\mathcal{T}(\mathfrak{R}))}$ a polynomial algebra?
\index{polynomial algebra}%


\section{An Isomorphism of Two Algebras}
 
Two relational structures
\index{relational structure}%
 with the same age have the same algebra up to isomorphism of graded algebras, hence the name \emph{age algebra}.
\index{age algebra}%
In fact the homogeneous components are indexed by an age of the appropriate size and so algebras built out of the age are equal and not just isomorphic.  The converse is not true.

In this section we will find an isomorphism between two algebras with very different structures in the finite algebra case, one graded and the other semisimple, and see how an isomorphism arises between them in the infinite case.

Let $\mathcal{M}$ be a relational structure
\index{relational structure}%
 on a set $X$, finite or infinite. For any
natural number $n$, let $V_n$ denote the vector space of all functions from
$n$-element subsets of $X$ to $\mathbb{C}$ which are constant on isomorphism
classes of $n$-element substructures of $\mathcal{M}$, and let
\[\mathcal{A}=\bigoplus_{n\ge0}V_n.\]

There are two ways to define multiplication on $\mathcal{A}$; each turns it
into a commutative and associative algebra.
\begin{itemize}
\item{Cameron:} For $f\in V_n$, $g\in V_m$, put
\[(f\circ g)(Y)=\sum\{f(U)g(V):U\cup V=Y,U\cap V=\emptyset\}.\]
\item{Glynn:} For $f\in V_n$, $g\in V_m$, put
\[(f*g)(Y)=\sum\{f(U)g(V):U\cup V=Y\}.\]
\end{itemize}
\index{Glynn, D. G.}%

In the case where $X$ is finite, these algebras are very different. Cameron's
is graded (since the product of elements of $V_n$ and $V_m$ is non-zero only
on sets of cardinality $m+n$, and so belongs to $V_{n+m}$), but any element
with $V_0$ component zero is nilpotent
(indeed $f^N=0$ if $N>|X|$).

 Glynn's algebra is not graded (the displayed product $fg$ can have non-zero components in $V_k$ for $\max(n,m)\le k\le n+m$), but is semisimple and commutative.  Glynn's construction~\cite{glynn} of rings of geometries
\index{ring of geometries}%
 (the term geometry includes graphs, matrices, projective planes, and many other combinatorial structures), in which addition is disjoint union and multiplication of graphs $\Gamma_i, \Gamma_j$ on vertices $n_i, n_j$ is defined by $\Gamma_i \times \Gamma_j = \sum \mu^k_{ij} \Gamma_k$ where $\mu^k_{ij}$ is the number of ways of writing $\Gamma_k = \gamma_i \cup \gamma_j$ where $\gamma_i \cong \Gamma_i$ and $\gamma_j \cong \Gamma_j$.  Here $\mu^k_{ij} \ne 0 \Rightarrow \max(n_i, n_j) \le | \Gamma_k| \le n_i + n_j$.  His ring is defined on any closed set
\index{closed set}%
  of geometries with coefficients in $\mathbb{Z}$ and the geometries as basis elements; coefficients in a field would instead give an algebra.  The semisimplicity arises from the fact that each geometry in a finite closed set of them corresponds to a homomorphism to $\mathbb{Z}$ and to a principal idempotent of the ring, and these idempotents are pairwise orthogonal~\cite[Theorem~2.5]{glynn}.  A closed set of $n$ geometries is isomorphic to the direct sum of $\mathbb{Z}$, $n$ times.

We begin by investigating the very simplest infinite case,
that where $\mathcal{M}$ has no structure at all, consisting simply of the set $X$.
In this case, $V_n$ is $1$-dimensional, spanned by the constant function 
$f_n$ on $n$-sets with value $1$. In either algebra, we have
\[f_1^n=n!\,f_n+\hbox{lower terms},\]
(the lower terms are absent in Cameron's case), so the algebra is a polynomial
algebra in one variable, generated by $f_1$. 

We calculate explicitly an isomorphism from Glynn's algebra to Cameron's in
this case. Such an isomorphism $\theta$ will be of the form
\[\theta(f_n)=\sum_{i=1}^na_{n,i}f_i,\]
with $a_{n,n}=1$ for all $n$. We have
\[f_n*f_1=(n+1)f_{n+1}+nf_n,\qquad f_n\circ f_1=(n+1)f_{n+1}.\]
Hence we have
\[(n+1)\theta(f_{n+1})+n\theta(f_n) = \theta(f_n*f_1)
= \theta(f_n)\circ\theta(f_1) = \theta(f_n)\circ f_1.\]
Thus, we find the recurrence relation
\[(n+1)a_{n+1,i}=ia_{n,i-1}-na_{n,i},\]
with the convention that $a_{n,0}=a_{n,n+1}=0$.

We find easily that $a_{n,n}=1$ and $a_{n,1}=(-1)^{n-1}/n$. For $n\le6$,
computation gives the following results:

\centerline{$ \begin{pmatrix} 
1\cr
-\frac{1}{2}&1\cr
\frac{1}{3}&-1&1\cr
-\frac{1}{4}&\frac{11}{12}&-\frac{3}{2}&1\cr
\frac{1}{5}&-\frac{5}{6}&\frac{7}{4}&-2&1\cr
-\frac{1}{6}&\frac{137}{180}&-\frac{15}{8}&\frac{17}{6}&-\frac{5}{2}&1\cr
\end{pmatrix}$.}

The inverse of this matrix (giving the reverse isomorphism) is

\centerline{$ \begin{pmatrix} 
1\cr
\frac{1}{2}&1\cr
\frac{1}{6}&1&1\cr
\frac{1}{24}&\frac{7}{12}&\frac{3}{2}&1\cr
\frac{1}{120}&\frac{1}{4}&\frac{5}{4}&2&1\cr
\frac{1}{720}&\frac{31}{360}&\frac{3}{4}&\frac{13}{6}&\frac{5}{2}&1\cr
\end{pmatrix}$.}

Inspection of the numbers reveals that in the second matrix, modulo a small shift, the pattern is of Stirling numbers of the second kind, and in the first matrix Stirling numbers of the first kind, the two families of Stirling numbers being mutual inverses.  In fact normalizing the numbers by multiplying by the highest common multiple, we arrive at the Triangle of Numbers $T(n,k) = k! S(n,k)$ read by rows $(n \geq 1, 1 \leq k < n)$ (sequence A019538  in~\cite{sloane}), where $S(n, k)$ is the Stirling number of the second kind.
\index{Stirling number ! of second kind}%

The reason for the occurrence of the $S(n, k)$
\index{Stirling number ! of second kind}%
 is that the overlapping terms in Glynn's algebra (that is the terms $a_{n,k}$ in $\theta(f_n)=\sum_{k=1}^n a_{n,k} f_k$) are combinatorially precisely the same as having indistinct elements in $G$-orbits on $k$-tuples in an implied wreath product
\index{wreath product}%
 action and resulting fibre bundle, where each $a_{n,k}$ term is the number of partitions of the $n$-th fibre with $k$ parts.  The group action,
\index{group ! action}%
 say $G \Wr H$
\index{wreath product}%
 acting on a set $X \times Y$ producing the partitions is implied, and the $k$ parts in a fibre are the ``connected'' subsets, i.e. connected graphs.  That the Stirling numbers occur in a lower triangular matrix is obviously because $k \leq n$.  Since orbits and partitions are equivalent, the equation $\theta(f_n)=\sum_{k=1}^na_{n,k}f_k$ is equivalent to the appearance of the Stirling numbers of the second kind that are the entries of the transition matrix in $F_{n}^{*} = \sum^{n}_{k=1} S(n, k) F_{k}$; (see Appendix~\ref{EnumerationandReconstruction}).

\bigskip

More generally, assume as we found above, that $f_n$ is a linear function of $f_1 \ldots f_{n-1}$ for all $n$.  By induction,

\begin{align*}
(\bullet) \times f_n &= (\bullet) \times \phi_1(f_1, \ldots, f_{n}) \\
     &= \phi_0(f_{n+1}, \phi_1(f_1, \ldots, f_{n})) \\
     &= \phi_0(f_{n+1}, \phi_1(f_1, \ldots, f_{n}), \phi_2(f_1, \ldots, f_{n-1})) \\
     &= \ldots \\
     &= \phi_0(f_{n+1}, \phi_1(f_1, \ldots, f_{n}), \phi_2(f_1, \ldots, f_{n-1}), \ldots, \phi_{n-1}(f_1,  f_2)),
\end{align*}
where each of $\phi_1, \ldots, \phi_{n-1}$ are previously determined functions with known coefficients as in the above matrices, and the leading $f_{n+1}$ terms in the $\phi_0$ expression is identical for the Cameron and Glynn algebras.

We conclude that 
\begin{theorem}
\label{camglyn}
In the limit of infinite pure point sets without any structure, the Cameron algebra and the Glynn algebra are isomorphic.
\end{theorem}

Can we find an isomorphism between the Cameron and Glynn algebras for more general structures?  In the setting of vector spaces over finite fields, the Gaussian (or $q$-binomial) coefficient plays the role of Stirling numbers.  Thus there may be a vector space version of the algebra isomorphism.  

\medskip

The next case to study is that of two infinite sets with a direct product of symmetric groups $\Sym(\infty) \times \Sym(\infty)$,
\index{group ! symmetric}%
 but we leave this for the interested reader.

\medskip

The question as to whether or not there is an algebra isomorphism in the case of graphs is open.  We have not attempted a full proof in the case of graphs but we realize the local isomorphism for graphs of up to three vertices in the appendix to this chapter.  Either algebra can be defined for graphs in general or for random graphs in particular.  We know that the Cameron algebra is a polynomial ring in infinitely many variables (corresponding to connected graphs for $\mathfrak{R}$, or graphs which are connected when a distinguished colour is made transparent for $\mathfrak{R^{t}}$.)  The former has been studied, see for example the first section of this chapter and references therein; we do not think that Glynn's algebra has been looked at in this case.

As we saw in the first section of this chapter, homogeneous structures have invariant rings that are polynomial rings, where the homogeneous generators are the characteristic functions of the ``connected'' structures in the age; therefore they are integral domains.  

\bigskip

\emph{Inexhaustibility}

As we mentioned earlier in the book, there are constructions of non-isomorphic structures with equal ages~\cite{drostemac}.  However two homogeneous structures with isomorphic ages are isomorphic~\cite{cam6}; (note that this is no longer true if the structures are just $\aleph_0$-categorical). 

An \emph{inexhaustible homogeneous structure}
\index{structure ! inexhaustible}%
 $\mathcal{M}$ is one for which $\mathcal{M} \backslash A \cong \mathcal{M}$, for all finite $A \subseteq \mathcal{M}$, or equivalently, $\Aut(\mathcal{M})_{A}$ has only finite orbits.  The notion of \emph{inexhaustible} relational structures originates in Fra\"{\i}ss\'e's work~\cite{frai}.
\index{Fra\"{\i}ss\'e, R.}%
 
 El-Zahar and Sauer
\index{El-Zahar, M.}%
\index{Sauer, N.}%
 studied~\cite{elzahar} Ramsey-type
\index{Ramsey theory}%
 properties of infinite relational structures, connecting different types of indivisibility and inexhaustibility.  B\"or\"oczky, Sauer and Zhu 
\index{B\"or\"oczky, K.}%
\index{Sauer, N.}%
\index{Zhu, X.}%
 further studied~\cite{boroczky} variations on the theme of inexhaustibility.  They show that \emph{inexhaustible} $\Rightarrow$ \emph{weakly inexhaustible} $\Rightarrow$ \emph{age-inexhaustible} $\Rightarrow$ \emph{closure-inexhaustible}; see their paper for definitions of these and other related concepts.  For an infinite structure on an infinite domain, the age
\index{age}%
 of $\mathcal{M}$ is \emph{inexhaustible} if given any two finite substructures of $\mathcal{M}$, is it possible to find disjoint copies of them within $\mathcal{M}$.  A homogeneous structure is age-inexhaustible if and only if all orbits of the empty set $A = \emptyset$ are infinite, such as for the Cameron algebra.  For finite languages weakly inexhaustible and age-inexhaustible are equivalent, and in the language of graphs, closure-inexhaustible is also equivalent.  For infinite languages the hierarchy is strict.
 
If $\mathcal{M}$ encodes $G = \Sym(X)$, then $\ker(\mathcal{M})$ is the union of the finite $G$-orbits of the 1-element sets, so if $G$ has no finite orbit then $\ker(\mathcal{M}) = \emptyset$, and $\mathcal{A}^G$ is an integral domain.  Pouzet showed~\cite{pouzet1} that the age algebra
\index{age algebra}%
 $\mathbb{C} \mathcal{A}(\mathcal{M})$ of $\mathcal{M}$ is an integral domain if and only if $\mathcal{M}$ is age-inexhaustible.  The amalgamation condition for weakly inexhaustible structures is strong amalgamation,
\index{amalgamation property ! strong}%
 where the ages on labelled and unlabelled structures are the same (see Appendix~\ref{TheoryofRelationalStructures} and~\cite{cam2cc}).

We can now state a conjecture:

\begin{conjecture}
If $\mathcal{M}$ is an age-inexhaustible relational structure
\index{relational structure}%
 then the Cameron and Glynn algebras are isomorphic.
\end{conjecture}

Earlier in the section we proved the conjecture for sets with no structure (or equivalently for $\Sym(X)$) using Stirling numbers.  The case of groups, for example $G = \Sym(X)$, is a special case of that of relational structures.   We shall prove the converse of this conjecture, finding it instructive nevertheless to give the argument for the group case separately.  For we shall show in the case of both groups and relational structures, that if $\mathcal{M}$ is not age-inexhaustible then there is a finite set $S$ consisting of all points of some given type.  So if $f$ is a characteristic function of $S$ then$f^2 = 0$ in the Cameron algebra, but the Glynn algebra for finite sets, being semisimple has no nonzero nilpotent elements.

\emph{Group Case:}

A statement of the conjecture in the case of a group $G$ is:  if $G$ has no finite orbits then the Cameron and Glynn algebras are isomorphic.

To prove that the converse is true, let $S$ be a finite orbit and let $f$ its characteristic function so that $f^2 = 0$ in the Cameron algebra, as is proved in previous work on this subject, see~\cite{cam8} and references therein.  We need to show that the Glynn algebra never contains such elements.  Let $G = \Sym(X)$ so that the algebra is spanned by $f_0, f_1, f_2, \ldots$, where $f_k$ is the characteristic function of the set of all $k$-sets.  Then
\[ \big( \sum_{k = 1}^{n} a_k f_k \big)  \big( \sum_{k' = 1}^{m} b_{k'} f_{k'} \big), \quad a_n, b_m \neq 0    \]
has a term $a_n b_m f_{n+m} \neq 0$ and lower terms.  For $f \neq 0$, write $f = \sum a_i f_{(k, i)}$, where $f_{(k, i)}$ is the characteristic function of some orbit of $G$ on $k$-sets.  Then $f^2 = \sum a_i b_j f_{(k, i)} f_{(l, j)}$, and the argument goes through.

\emph{Relational Structure Case:}

A statement of the conjecture in the case of a group $G$ is:  if  $\mathcal{M}$ is not an age-inexhaustible relational structure then the Cameron and Glynn algebras are not isomorphic.

For each $k \in \mathbb{N}$ and for each isomorphism type $T_i$ of $k$-element substructures of $\mathcal{M}$, let $f_{k, i}$ be the characteristic function of the set of $k$-element substructures isomorphic to $T_i$.  Then 
\[ (f_{k, i} \star f_{l, j}) (S) = \sharp (S_1, S_2) : S_1 \cup S_2 = S, \]
where $f_{k, i}(T_1) \neq 0, f_{k, j}(T_2) \neq 0, S_1 \in T_1, S_2 \in T_2$.  (In the case of the Cameron algebra there would be the the additional condition that $S_1 \cap S_2 = \emptyset$).  Extend linearly to the span of all such terms.  

We claim that if $U \subseteq X = \dom(\mathcal{M})$, then we find that the Glynn algebra of $\mathcal{M}_{| U}$, not necessarily a subalgebra, but $(f_{k, i} \star f_{l, j}) (U)$ is the same in this smaller algebra than in the whole algebra.  

Suppose that $f^2 = 0$.  Then take any subset $U$ for which $f_{| U} \neq 0$,  $f^{2}_{| U} = 0$.  This meets the hypotheses of Glynn's algebra, but gives us a contradiction as in the finite case this algebra is semisimple and has no nilpotent elements. 


\bigskip


\emph{Growth Rates and Reconstruction}

The algebraic considerations of this section are related to vertex reconstruction~\cite{cam6a}~\cite{cam6b}~\cite{pouzetthieryIII};
\index{reconstruction}%
\index{growth rates}%
(see Appendix~\ref{EnumerationandReconstruction}).  We end the chapter by stating and motivating a conjecture that relates growth rates and reconstruction.  

\begin{conjecture}
Let $\mathcal{M}$ be a relational structure
\index{relational structure}%
 and let $G = \Aut(\mathcal{M})$.  If $F_n(G)$ is known together with a bound on its growth, then is $\mathcal{M}$ reconstructible.
\end{conjecture}

What evidence do we have that such a statement might be true?

It is plausible that if growth is not too fast, for example is slower than $n^2$, then we may be able to prove reconstruction.  In the case of polynomial sequence growth, M. Pouzet
\index{Pouzet, M.}%
 has a structural characterization of the growth and it is slow enough for the structures to be reconstructible.
 
If the growth rate is too fast then reconstruction is not possible, for example in the case of $n$-element structures, reconstruction fails if $f_n > (f_{n-1})^n$, so this represents an upper bound for $(n-1)$-element substructures.  Take for instance $f_n = (f_{n-1})^n$ for $n \geq 2$ in the case that $f_1 = 2$ giving $f_n = 2^{n!}$ which is an exponential of an exponential (by Stirling's formula); this is faster than the growth rate for graphs which is merely exponential.  

There is probably a band in the middle of these two possibilities where we can ask questions such as 
\head{Open Question}  

Is it true that there are two ages $\mathcal{A}$ and $\mathcal{B}$ such that $f_n(\mathcal{A}) = f_n(\mathcal{B})$ $\forall n$ whilst $\mathcal{A}$ is reconstructible but  $\mathcal{B}$ is not.


Age reconstruction of graphs is equivalent to the reconstruction of the random graph~\cite{cammar}.  
\index{graph ! random}%

\smallskip

Finally we mention a reference to the Hopf algebras with bases labelled by graphs and hypergraphs, provided with natural embeddings into a polynomial algebra
\index{polynomial algebra}%
 in infinitely many variables, and graded by the number of edges; these algebras can be considered as generalizations of symmetric functions~\cite{novelli}.
\index{symmetric function}%

\section{Appendix: Cameron $\&$ Glynn Algebra Graph Isomorphisms}

We will demonstrate an isomorphism between the Cameron and Glynn algebras in Theorem~\ref{camglyn} for  graphs on $1, 2$, and $3$ vertices.

The Cameron algebra expression is given first, followed by the `leads to' symbol $\leadsto$ and then the equivalent Glynn algebra term.  The latter includes smaller terms corresponding to overlapping graphs in Glynn's algebra, with coefficients to be determined.  The appropriate isomorphism is defined by mapping an element of the Cameron algebra onto the corresponding Glynn algebra term, by showing that the Glynn coefficients can be assigned unambiguous values.  

\textbf{One and Two vertex graphs}

$\left( 
\def\objectstyle{\scriptstyle} 
\def\labelstyle{\scriptstyle} 
\vcenter{\xymatrix @-1.2pc { 
\emptyset}} 
\right)$ $\quad \leadsto \quad$ $\left( 
\def\objectstyle{\scriptstyle} 
\def\labelstyle{\scriptstyle} 
\vcenter{\xymatrix @-1.2pc @ur { 
\emptyset  }} 
\right)$

$\left( 
\def\objectstyle{\scriptstyle} 
\def\labelstyle{\scriptstyle} 
\vcenter{\xymatrix @-1.2pc { 
\bullet }} 
\right)$ $\quad \leadsto \quad$ $\left( 
\def\objectstyle{\scriptstyle} 
\def\labelstyle{\scriptstyle} 
\vcenter{\xymatrix @-1.2pc @ur { 
\bullet\ +\ a\ \emptyset  }} 
\right)$

$\left( 
\def\objectstyle{\scriptstyle} 
\def\labelstyle{\scriptstyle} 
\vcenter{\xymatrix @-1.2pc { 
(\bullet \ar@{-}[r] & \bullet) }} 
\right)$ $\quad \leadsto \quad$ $\left( 
\def\objectstyle{\scriptstyle} 
\def\labelstyle{\scriptstyle} 
\vcenter{\xymatrix @-1.2pc { 
(\bullet \ar@{-}[r] & \bullet)\ +\ b\ \bullet\ +\ c\ \emptyset  }} 
\right)$

$\left( 
\def\objectstyle{\scriptstyle} 
\def\labelstyle{\scriptstyle} 
\vcenter{\xymatrix @-1.2pc { 
(\bullet \ar@{} & \bullet) }} 
\right)$ $\quad \leadsto \quad$ $\left( 
\def\objectstyle{\scriptstyle} 
\def\labelstyle{\scriptstyle} 
\vcenter{\xymatrix @-1.2pc { 
(\bullet \ar@{} & \bullet)\ +\ d\ \bullet\ +\ e\ \emptyset  }} 
\right)$

Together with multiplication rules,

$\def\objectstyle{\scriptstyle} 
\def\labelstyle{\scriptstyle} 
\vcenter{\xymatrix @-1.2pc { 
\bullet\ \times\ \emptyset\ =\ \bullet }}
$\quad   and\quad $\def\objectstyle{\scriptstyle} 
\def\labelstyle{\scriptstyle} 
\vcenter{\xymatrix @-1.2pc { 
(\bullet \ar@{}[r] & \bullet)\ \times\ \emptyset\ =\ 2\ (\bullet \ar@{}[r] & \bullet) }}
$\quad  and\quad $\def\objectstyle{\scriptstyle} 
\def\labelstyle{\scriptstyle} 
\vcenter{\xymatrix @-1.2pc { 
(\bullet \ar@{-}[r] & \bullet)\ \times\ \emptyset\ =\ 2 (\bullet \ar@{-}[r] & \bullet)}}
$  

we have that for

Cameron: 

$\def\objectstyle{\scriptstyle} 
\def\labelstyle{\scriptstyle} 
\vcenter{\xymatrix @-1.2pc { 
\bullet\ \times\ \bullet\ =\ 2 (\bullet \ar@{-}[r] & \bullet)\ +\ 2 (\bullet \ar@{}[r] & \bullet) }}
$ 

$\quad \leadsto \quad$ 

Glynn: 

$\left( 
\def\objectstyle{\scriptstyle} 
\def\labelstyle{\scriptstyle} 
\vcenter{\xymatrix @-1.2pc { 
\bullet\ +\ a\ \emptyset }} 
\right)$ $\times$ $\left( 
\def\objectstyle{\scriptstyle} 
\def\labelstyle{\scriptstyle} 
\vcenter{\xymatrix @-1.2pc @ur { 
\bullet\ +\ a\ \emptyset }} 
\right)$ $\quad=\quad$ $\def\objectstyle{\scriptstyle} 
\def\labelstyle{\scriptstyle} 
\vcenter{\xymatrix @-1.2pc { 
2 (\bullet \ar@{-}[r] & \bullet)\ +\ 2(\bullet\ \ar@{}[r] & \bullet)\ +\ \bullet\ +\ 2\ a\ \bullet\ + a^2\ \emptyset }} 
$ 

(The single vertex in the last expression is from the superposition of one vertex on top of another.)  But also

$\def\objectstyle{\scriptstyle} 
\def\labelstyle{\scriptstyle} 
\vcenter{\xymatrix @-1.2pc { 
\bullet\ \times\ \bullet\ =\ 2 (\bullet \ar@{-}[r] & \bullet)\ +\ 2 (\bullet \ar@{}[r] & \bullet)\quad \leadsto \quad}}
$

$\def\objectstyle{\scriptstyle} 
\def\labelstyle{\scriptstyle} 
\vcenter{\xymatrix @-1.2pc { 
\quad\quad\quad 2 (\bullet \ar@{-}[r] & \bullet)\ +\ 2b\ \bullet\ +\ 2c \emptyset\quad  +\quad 2 (\bullet \ar@{} & \bullet)\ +\ 2d\ \bullet\ +\ 2e \emptyset }}$.

So
\begin{eqnarray*}
2a + 1 =  2b + 2d\\
a^2 = 2c + 2e.
\end{eqnarray*}

One possible isomorphism is found by setting $a = -\frac{1}{2}$ in order to eliminate the single-vertex terms on the right-hand side, giving

$\def\objectstyle{\scriptstyle} 
\def\labelstyle{\scriptstyle} 
\vcenter{\xymatrix @-1.2pc { 
\bullet }}$ $\quad \leadsto \quad$ 
$\def\objectstyle{\scriptstyle} 
\def\labelstyle{\scriptstyle} 
\vcenter{\xymatrix @-1.2pc @ur { 
\bullet\ -\ \frac{1}{2}\ \emptyset  }}$

$\def\objectstyle{\scriptstyle} 
\def\labelstyle{\scriptstyle} 
\vcenter{\xymatrix @-1.2pc { 
(\bullet \ar@{-}[r] & \bullet) }}$ $\quad \leadsto \quad$
$\def\objectstyle{\scriptstyle} 
\def\labelstyle{\scriptstyle} 
\vcenter{\xymatrix @-1.2pc { 
(\bullet \ar@{-}[r] & \bullet)\ +\ b\ \bullet\ +\ c\ \emptyset  }}$

$\def\objectstyle{\scriptstyle} 
\def\labelstyle{\scriptstyle} 
\vcenter{\xymatrix @-1.2pc { 
(\bullet \ar@{} & \bullet) }}$ $\quad \leadsto \quad$
$\def\objectstyle{\scriptstyle} 
\def\labelstyle{\scriptstyle} 
\vcenter{\xymatrix @-1.2pc { 
(\bullet \ar@{} & \bullet)\ -\ b\ \bullet\ +\ e\ \emptyset  }}$

and thus the isomorphism

$\def\objectstyle{\scriptstyle} 
\def\labelstyle{\scriptstyle} 
\vcenter{\xymatrix @-1.2pc { 
\bullet\ \times\ \bullet\ =\ 2 (\bullet \ar@{-}[r] & \bullet)\ +\ 2 (\bullet \ar@{}[r] & \bullet)\quad {\longrightarrow} \quad 2 (\bullet \ar@{-}[r] & \bullet)\ +\ 2 (\bullet \ar@{} & \bullet)\ +\ \frac{1}{4} \emptyset }}$.

A second possible isomorphism retains the single vertex term and eliminates the no-vertex term, as follows.  Set
\begin{eqnarray*}
a =  0;\quad c = e = 0.
\end{eqnarray*}
Then \begin{eqnarray*}
d = \frac{1}{2} - b.
\end{eqnarray*}
So the isomorphism is

$\quad\quad\quad\quad \def\objectstyle{\scriptstyle} 
\def\labelstyle{\scriptstyle} 
\vcenter{\xymatrix @-1.2pc { 
\bullet }}$ $\quad \leadsto \quad$ 
$\def\objectstyle{\scriptstyle} 
\def\labelstyle{\scriptstyle} 
\vcenter{\xymatrix @-1.2pc @ur { 
\bullet }}$

$\quad\quad\quad\quad  \def\objectstyle{\scriptstyle} 
\def\labelstyle{\scriptstyle} 
\vcenter{\xymatrix @-1.2pc { 
(\bullet \ar@{-}[r] & \bullet) }}$ $\quad \leadsto \quad$
$\def\objectstyle{\scriptstyle} 
\def\labelstyle{\scriptstyle} 
\vcenter{\xymatrix @-1.2pc { 
(\bullet \ar@{-}[r] & \bullet)\ +\ b\ \bullet }}$

$\quad\quad\quad\quad  \def\objectstyle{\scriptstyle} 
\def\labelstyle{\scriptstyle} 
\vcenter{\xymatrix @-1.2pc { 
(\bullet \ar@{} & \bullet) }}$ $\quad \leadsto \quad$
$\def\objectstyle{\scriptstyle} 
\def\labelstyle{\scriptstyle} 
\vcenter{\xymatrix @-1.2pc { 
(\bullet \ar@{} & \bullet)\ +\ (\frac{1}{2}- b)\ \bullet  }}$

in turn giving the isomorphism

$\def\objectstyle{\scriptstyle} 
\def\labelstyle{\scriptstyle} 
\vcenter{\xymatrix @-1.2pc { 
\bullet\ \times\ \bullet\ =\ 2 (\bullet \ar@{-}[r] & \bullet)\ +\ 2 (\bullet \ar@{}[r] & \bullet)\quad {\longrightarrow} \quad 2 (\bullet \ar@{-}[r] & \bullet)\ +\ 2 (\bullet \ar@{} & \bullet)\ +\ \bullet }}$.

It is this second isomorphism that we will use.

\textbf{Three vertex graphs}

(i)  Cameron: 

$\def\objectstyle{\scriptstyle} 
\def\labelstyle{\scriptstyle} 
\vcenter{\xymatrix @-1.2pc { 
\bullet\ \times\ (\bullet \ar@{}[r] & \bullet) }}
$ 

$\quad\quad\quad=\quad$  $3$ $\left( \def\objectstyle{\scriptstyle} 
\def\labelstyle{\scriptstyle} 
\vcenter{\xymatrix @-1.2pc { 
& \bullet  \ar@{}[dr]  \ar@{}[dl] \\
\bullet  \ar@{}[rr] && \bullet }}
\right)$  $+$ 2 $\left( \def\objectstyle{\scriptstyle} 
\def\labelstyle{\scriptstyle} 
\vcenter{\xymatrix @-1.2pc { 
& \bullet  \ar@{}[dr]  \ar@{}[dl] \\
\bullet  \ar@{-}[rr] && \bullet }}
\right)$ $+ $ $\left( \def\objectstyle{\scriptstyle} 
\def\labelstyle{\scriptstyle} 
\vcenter{\xymatrix @-1.2pc { 
& \bullet  \ar@{-}[dr]  \ar@{-}[dl] \\
\bullet  \ar@{}[rr] && \bullet }}
\right)$ 

$\quad \leadsto \quad$ 

Glynn: 

$\left( 
\def\objectstyle{\scriptstyle} 
\def\labelstyle{\scriptstyle} 
\vcenter{\xymatrix @-1.2pc { 
\bullet\ +\ a\ \emptyset }}
\right)$ $\times$ $\left( 
\def\objectstyle{\scriptstyle} 
\def\labelstyle{\scriptstyle} 
\vcenter{\xymatrix @-1.2pc { 
(\bullet \ar@{}[r] & \bullet)\ +\ d\ \bullet\ +\ e\ \emptyset }} 
\right)$ 

$\quad\quad\quad=\quad$  $3$ $\left( \def\objectstyle{\scriptstyle} 
\def\labelstyle{\scriptstyle} 
\vcenter{\xymatrix @-1.2pc { 
& \bullet  \ar@{}[dr]  \ar@{}[dl] \\
\bullet  \ar@{}[rr] && \bullet }}
\right)$  $+$ 2 $\left( \def\objectstyle{\scriptstyle} 
\def\labelstyle{\scriptstyle} 
\vcenter{\xymatrix @-1.2pc { 
& \bullet  \ar@{}[dr]  \ar@{}[dl] \\
\bullet  \ar@{-}[rr] && \bullet }}
\right)$ $+ $ $\left( \def\objectstyle{\scriptstyle} 
\def\labelstyle{\scriptstyle} 
\vcenter{\xymatrix @-1.2pc { 
& \bullet  \ar@{-}[dr]  \ar@{-}[dl] \\
\bullet  \ar@{}[rr] && \bullet }}
\right)$ 

$\quad\quad\quad\quad+$ $(ad + e) \emptyset$ $\left( \def\objectstyle{\scriptstyle} 
\def\labelstyle{\scriptstyle} 
\vcenter{\xymatrix @-1.2pc { 
\bullet }}
\right)$ $+$  $(ae \emptyset^{2})$ 

$\quad\quad\quad\quad+$ $\def\objectstyle{\scriptstyle} 
\def\labelstyle{\scriptstyle} 
\vcenter{\xymatrix @-1.2pc { 
2\ d\ (\bullet \ar@{-}[r] & \bullet)\ +\ (2a + 2d)\ (\bullet \ar@{}[r] & \bullet)\ +\ d\ (\bullet)}}$

Let\\

$\left( \def\objectstyle{\scriptstyle} 
\def\labelstyle{\scriptstyle} 
\vcenter{\xymatrix @-1.2pc { 
& \bullet  \ar@{}[dr]  \ar@{}[dl] \\
\bullet  \ar@{}[rr] && \bullet }}
\right)$ $\mapsto$ $\left( \def\objectstyle{\scriptstyle} 
\def\labelstyle{\scriptstyle} 
\vcenter{\xymatrix @-1.2pc { 
& \bullet  \ar@{}[dr]  \ar@{}[dl] \\
\bullet  \ar@{}[rr] && \bullet }}
\right)$ $+$ $\alpha$ $\left( \def\objectstyle{\scriptstyle} 
\def\labelstyle{\scriptstyle} 
\vcenter{\xymatrix @-1.2pc {
\bullet  \ar@{-}[rr] && \bullet }}
\right)$  $+$ $\alpha'$ $\left( \def\objectstyle{\scriptstyle} 
\def\labelstyle{\scriptstyle} 
\vcenter{\xymatrix @-1.2pc {
\bullet  \ar@{}[rr] && \bullet }}
\right)$ 

$\quad\quad\quad\quad\quad\quad\quad\quad\quad\quad\quad\quad\quad+\ \alpha''\ (\bullet)$

$\left( \def\objectstyle{\scriptstyle} 
\def\labelstyle{\scriptstyle} 
\vcenter{\xymatrix @-1.2pc { 
& \bullet  \ar@{-}[dr]  \ar@{-}[dl] \\
\bullet  \ar@{-}[rr] && \bullet }}
\right)$ $\mapsto$ $\left( \def\objectstyle{\scriptstyle} 
\def\labelstyle{\scriptstyle} 
\vcenter{\xymatrix @-1.2pc { 
& \bullet  \ar@{-}[dr]  \ar@{-}[dl] \\
\bullet  \ar@{-}[rr] && \bullet }}
\right)$ $+$ $\beta$ $\left( \def\objectstyle{\scriptstyle} 
\def\labelstyle{\scriptstyle} 
\vcenter{\xymatrix @-1.2pc {
\bullet  \ar@{-}[rr] && \bullet }}
\right)$   $+$ $\beta'$ $\left( \def\objectstyle{\scriptstyle} 
\def\labelstyle{\scriptstyle} 
\vcenter{\xymatrix @-1.2pc {
\bullet  \ar@{}[rr] && \bullet }}
\right)$ 

$\quad\quad\quad\quad\quad\quad\quad\quad\quad\quad\quad\quad\quad+\ \beta''\ (\bullet)$

$\left( \def\objectstyle{\scriptstyle} 
\def\labelstyle{\scriptstyle} 
\vcenter{\xymatrix @-1.2pc { 
& \bullet  \ar@{-}[dr]  \ar@{-}[dl] \\
\bullet  \ar@{}[rr] && \bullet }}
\right)$ $\mapsto$ $\left( \def\objectstyle{\scriptstyle} 
\def\labelstyle{\scriptstyle} 
\vcenter{\xymatrix @-1.2pc { 
& \bullet  \ar@{-}[dr]  \ar@{-}[dl] \\
\bullet  \ar@{}[rr] && \bullet }}
\right)$ $+$ $\gamma$ $\left( \def\objectstyle{\scriptstyle} 
\def\labelstyle{\scriptstyle} 
\vcenter{\xymatrix @-1.2pc {
\bullet  \ar@{-}[rr] && \bullet }}
\right)$   $+$ $\gamma'$ $\left( \def\objectstyle{\scriptstyle} 
\def\labelstyle{\scriptstyle} 
\vcenter{\xymatrix @-1.2pc {
\bullet  \ar@{}[rr] && \bullet }}
\right)$ 

$\quad\quad\quad\quad\quad\quad\quad\quad\quad\quad\quad\quad\quad+\ \gamma''\ (\bullet)$

$\left( \def\objectstyle{\scriptstyle} 
\def\labelstyle{\scriptstyle} 
\vcenter{\xymatrix @-1.2pc { 
& \bullet  \ar@{}[dr]  \ar@{}[dl] \\
\bullet  \ar@{-}[rr] && \bullet }}
\right)$ $\mapsto$ $\left( \def\objectstyle{\scriptstyle} 
\def\labelstyle{\scriptstyle} 
\vcenter{\xymatrix @-1.2pc { 
& \bullet  \ar@{}[dr]  \ar@{}[dl] \\
\bullet  \ar@{-}[rr] && \bullet }}
\right)$ $+$ $\delta$ $\left( \def\objectstyle{\scriptstyle} 
\def\labelstyle{\scriptstyle} 
\vcenter{\xymatrix @-1.2pc {
\bullet  \ar@{-}[rr] && \bullet }}
\right)$   $+$ $\delta'$ $\left( \def\objectstyle{\scriptstyle} 
\def\labelstyle{\scriptstyle} 
\vcenter{\xymatrix @-1.2pc {
\bullet  \ar@{}[rr] && \bullet }}
\right)$ 

$\quad\quad\quad\quad\quad\quad\quad\quad\quad\quad\quad\quad\quad+\ \delta''\ (\bullet)$

So 

$\def\objectstyle{\scriptstyle} 
\def\labelstyle{\scriptstyle} 
\vcenter{\xymatrix @-1.2pc { 
\bullet\ \times\ (\bullet \ar@{}[r] & \bullet) }}
$ 

$\quad\quad\quad=\quad$  $3$ $\left( \def\objectstyle{\scriptstyle} 
\def\labelstyle{\scriptstyle} 
\vcenter{\xymatrix @-1.2pc { 
& \bullet  \ar@{}[dr]  \ar@{}[dl] \\
\bullet  \ar@{}[rr] && \bullet }}
\right)$  $+$ 2 $\left( \def\objectstyle{\scriptstyle} 
\def\labelstyle{\scriptstyle} 
\vcenter{\xymatrix @-1.2pc { 
& \bullet  \ar@{}[dr]  \ar@{}[dl] \\
\bullet  \ar@{-}[rr] && \bullet }}
\right)$ $+ $ $\left( \def\objectstyle{\scriptstyle} 
\def\labelstyle{\scriptstyle} 
\vcenter{\xymatrix @-1.2pc { 
& \bullet  \ar@{-}[dr]  \ar@{-}[dl] \\
\bullet  \ar@{}[rr] && \bullet }}
\right)$ 

$\quad \leadsto \quad$

$\quad \quad 3 \left( \def\objectstyle{\scriptstyle} 
\def\labelstyle{\scriptstyle} 
\vcenter{\xymatrix @-1.2pc { 
& \bullet  \ar@{}[dr]  \ar@{}[dl] \\
\bullet  \ar@{}[rr] && \bullet }}
\right)$ $+$ $\alpha$ $\left( \def\objectstyle{\scriptstyle} 
\def\labelstyle{\scriptstyle} 
\vcenter{\xymatrix @-1.2pc {
\bullet  \ar@{-}[rr] && \bullet }}
\right)$  $+$ $\alpha'$ $\left( \def\objectstyle{\scriptstyle} 
\def\labelstyle{\scriptstyle} 
\vcenter{\xymatrix @-1.2pc {
\bullet  \ar@{}[rr] && \bullet }}
\right)$ $\alpha''\ (\bullet)$

$+$ $\left( \def\objectstyle{\scriptstyle} 
\def\labelstyle{\scriptstyle} 
\vcenter{\xymatrix @-1.2pc { 
& \bullet  \ar@{-}[dr]  \ar@{-}[dl] \\
\bullet  \ar@{}[rr] && \bullet }}
\right)$ $+$ $6 \gamma$ $\left( \def\objectstyle{\scriptstyle} 
\def\labelstyle{\scriptstyle} 
\vcenter{\xymatrix @-1.2pc {
\bullet  \ar@{-}[rr] && \bullet }}
\right)$   $+$ $6 \gamma'$ $\left( \def\objectstyle{\scriptstyle} 
\def\labelstyle{\scriptstyle} 
\vcenter{\xymatrix @-1.2pc {
\bullet  \ar@{}[rr] && \bullet }}
\right)$ $+\ 6 \gamma''\ (\bullet)$

$+ 2$ $\left( \def\objectstyle{\scriptstyle} 
\def\labelstyle{\scriptstyle} 
\vcenter{\xymatrix @-1.2pc { 
& \bullet  \ar@{}[dr]  \ar@{}[dl] \\
\bullet  \ar@{-}[rr] && \bullet }}
\right)$ $+$ $3 \delta$ $\left( \def\objectstyle{\scriptstyle} 
\def\labelstyle{\scriptstyle} 
\vcenter{\xymatrix @-1.2pc {
\bullet  \ar@{-}[rr] && \bullet }}
\right)$   $+$ $3 \delta'$ $\left( \def\objectstyle{\scriptstyle} 
\def\labelstyle{\scriptstyle} 
\vcenter{\xymatrix @-1.2pc {
\bullet  \ar@{}[rr] && \bullet }}
\right)$ $+\ 3 \delta''\ (\bullet)$

But we have already set $a = e = 0$.   Matching up coefficients leaves us with only three non-zero ones, $\alpha = 2 d$,  $\alpha' = 2 d +2 a$, $\alpha'' = d$, where as previously determined
\begin{eqnarray*}
d = \frac{1}{2} - b.
\end{eqnarray*}

(In the next step we shall see that $b= 0,\ d = \frac{1}{2}$).

\bigskip

(ii)  Cameron: 

$\def\objectstyle{\scriptstyle} 
\def\labelstyle{\scriptstyle} 
\vcenter{\xymatrix @-1.2pc { 
\bullet\ \times\ (\bullet \ar@{-}[r] & \bullet) }}
$ 

$\quad\quad\quad=\quad$  $\left( \def\objectstyle{\scriptstyle} 
\def\labelstyle{\scriptstyle} 
\vcenter{\xymatrix @-1.2pc { 
& \bullet  \ar@{-}[dr]  \ar@{-}[dl] \\
\bullet  \ar@{-}[rr] && \bullet }}
\right)$ $+$   $\left( \def\objectstyle{\scriptstyle} 
\def\labelstyle{\scriptstyle} 
\vcenter{\xymatrix @-1.2pc { 
& \bullet  \ar@{}[dr]  \ar@{}[dl] \\
\bullet  \ar@{-}[rr] && \bullet }}
\right)$ $+$  2 $\left( \def\objectstyle{\scriptstyle} 
\def\labelstyle{\scriptstyle} 
\vcenter{\xymatrix @-1.2pc { 
& \bullet  \ar@{-}[dr]  \ar@{-}[dl] \\
\bullet  \ar@{}[rr] && \bullet }}
\right)$ 

$\quad \leadsto \quad$ 

Glynn: 

$\left( 
\def\objectstyle{\scriptstyle} 
\def\labelstyle{\scriptstyle} 
\vcenter{\xymatrix @-1.2pc { 
\bullet\ +\ a\ \emptyset }}
\right)$ $\times$ $\left( 
\def\objectstyle{\scriptstyle} 
\def\labelstyle{\scriptstyle} 
\vcenter{\xymatrix @-1.2pc { 
(\bullet \ar@{-}[r] & \bullet)\ +\ b\ \bullet\ +\ c\ \emptyset }} 
\right)$ 

$\def\objectstyle{\scriptstyle} 
\def\labelstyle{\scriptstyle} 
\vcenter{\xymatrix @-1.2pc { 
\bullet\ \times\ (\bullet \ar@{-}[r] & \bullet) }}
$ 

$\quad\quad\quad=\quad$  $\left( \def\objectstyle{\scriptstyle} 
\def\labelstyle{\scriptstyle} 
\vcenter{\xymatrix @-1.2pc { 
& \bullet  \ar@{-}[dr]  \ar@{-}[dl] \\
\bullet  \ar@{-}[rr] && \bullet }}
\right)$ $+$  $\left( \def\objectstyle{\scriptstyle} 
\def\labelstyle{\scriptstyle} 
\vcenter{\xymatrix @-1.2pc { 
& \bullet  \ar@{}[dr]  \ar@{}[dl] \\
\bullet  \ar@{-}[rr] && \bullet }}
\right)$ $+$  2 $\left( \def\objectstyle{\scriptstyle} 
\def\labelstyle{\scriptstyle} 
\vcenter{\xymatrix @-1.2pc { 
& \bullet  \ar@{-}[dr]  \ar@{-}[dl] \\
\bullet  \ar@{}[rr] && \bullet }}
\right)$ 

$\quad\quad\quad\quad\quad$ +  b $\left( \def\objectstyle{\scriptstyle} 
\def\labelstyle{\scriptstyle} 
\vcenter{\xymatrix @-1.2pc { 
\bullet\ \times\ \bullet }}
\right)$ $+$ $c \emptyset$ $\left( \def\objectstyle{\scriptstyle} 
\def\labelstyle{\scriptstyle} 
\vcenter{\xymatrix @-1.2pc { 
 \bullet }}
\right)$ $+$ $\def\objectstyle{\scriptstyle} 
\def\labelstyle{\scriptstyle} 
\vcenter{\xymatrix @-1.2pc { 
a\ \emptyset\ \times (\bullet \ar@{-}[r] & \bullet)\ +\ a\ b\ \emptyset\ (\bullet) }}
$ 

$\quad\quad\quad\quad\quad$ $+$  $\left( \def\objectstyle{\scriptstyle} 
\def\labelstyle{\scriptstyle} 
\vcenter{\xymatrix @-1.2pc { 
a\ c\ \emptyset^{2} }}
\right)$

$\quad\quad\quad=\quad$ $\left( \def\objectstyle{\scriptstyle} 
\def\labelstyle{\scriptstyle} 
\vcenter{\xymatrix @-1.2pc { 
& \bullet  \ar@{-}[dr]  \ar@{-}[dl] \\
\bullet  \ar@{-}[rr] && \bullet }}
\right)$ $+$  $\left( \def\objectstyle{\scriptstyle} 
\def\labelstyle{\scriptstyle} 
\vcenter{\xymatrix @-1.2pc { 
& \bullet  \ar@{}[dr]  \ar@{}[dl] \\
\bullet  \ar@{-}[rr] && \bullet }}
\right)$ $+$  2 $\left( \def\objectstyle{\scriptstyle} 
\def\labelstyle{\scriptstyle} 
\vcenter{\xymatrix @-1.2pc { 
& \bullet  \ar@{-}[dr]  \ar@{-}[dl] \\
\bullet  \ar@{}[rr] && \bullet }}
\right)$ 

$\quad\quad\quad\quad+$ $\def\objectstyle{\scriptstyle} 
\def\labelstyle{\scriptstyle} 
\vcenter{\xymatrix @-1.2pc { 
2\ b\ (\bullet \ar@{-}[r] & \bullet)\ +\ 2\ b\ (\bullet \ar@{}[r] & \bullet) }}$ $+$ $b$ $\left( \def\objectstyle{\scriptstyle} 
\def\labelstyle{\scriptstyle} 
\vcenter{\xymatrix @-1.2pc { 
 \bullet }}
\right)$

since $a = 0 = c$.  The single vertex term comes from the overlap of the two single vertices.  






We can eliminate the single vertex term by setting $b = 0$.  Then $d = \frac{1}{2}$ in the previous calculation.

\smallskip

As a check, notice that if we define the local Glynn maps by

$\quad\quad\quad\quad$  $\left( \def\objectstyle{\scriptstyle} 
\def\labelstyle{\scriptstyle} 
\vcenter{\xymatrix @-1.2pc { 
& \bullet  \ar@{-}[dr]  \ar@{-}[dl] \\
\bullet  \ar@{-}[rr] && \bullet }}
\right)$ $\quad \mapsto \quad$ $\left( \def\objectstyle{\scriptstyle} 
\def\labelstyle{\scriptstyle} 
\vcenter{\xymatrix @-1.2pc { 
& \bullet  \ar@{-}[dr]  \ar@{-}[dl] \\
\bullet  \ar@{-}[rr] && \bullet }}
\right)$ $+$ $x$ $\left( \def\objectstyle{\scriptstyle} 
\def\labelstyle{\scriptstyle} 
\vcenter{\xymatrix @-1.2pc {
\bullet  \ar@{-}[rr] && \bullet }}
\right)$ 

$\quad\quad\quad\quad$  $\left( \def\objectstyle{\scriptstyle} 
\def\labelstyle{\scriptstyle} 
\vcenter{\xymatrix @-1.2pc { 
& \bullet  \ar@{-}[dr]  \ar@{-}[dl] \\
\bullet  \ar@{}[rr] && \bullet }}
\right)$ $\quad \mapsto \quad$ $\left( \def\objectstyle{\scriptstyle} 
\def\labelstyle{\scriptstyle} 
\vcenter{\xymatrix @-1.2pc { 
& \bullet  \ar@{-}[dr]  \ar@{-}[dl] \\
\bullet  \ar@{}[rr] && \bullet }}
\right)$ $+$ $y$ $\left( \def\objectstyle{\scriptstyle} 
\def\labelstyle{\scriptstyle} 
\vcenter{\xymatrix @-1.2pc {
\bullet  \ar@{-}[rr] && \bullet }}
\right)$ 

$\quad\quad\quad\quad$  $\left( \def\objectstyle{\scriptstyle} 
\def\labelstyle{\scriptstyle} 
\vcenter{\xymatrix @-1.2pc { 
& \bullet  \ar@{}[dr]  \ar@{}[dl] \\
\bullet  \ar@{-}[rr] && \bullet }}
\right)$ $\quad \mapsto \quad$ $\left( \def\objectstyle{\scriptstyle} 
\def\labelstyle{\scriptstyle} 
\vcenter{\xymatrix @-1.2pc { 
& \bullet  \ar@{}[dr]  \ar@{}[dl] \\
\bullet  \ar@{-}[rr] && \bullet }}
\right)$ $+$ $z$ $\left( \def\objectstyle{\scriptstyle} 
\def\labelstyle{\scriptstyle} 
\vcenter{\xymatrix @-1.2pc {
\bullet  \ar@{-}[rr] && \bullet }}
\right)$ 

then on the one hand

$\def\objectstyle{\scriptstyle} 
\def\labelstyle{\scriptstyle} 
\vcenter{\xymatrix @-1.2pc { 
\bullet\ \times\ (\bullet \ar@{-}[r] & \bullet) }}
$  $\quad \leadsto \quad$
$\def\objectstyle{\scriptstyle} 
\def\labelstyle{\scriptstyle} 
\vcenter{\xymatrix @-1.2pc { 
\bullet\ \times\  ((\bullet \ar@{-}[r] & \bullet)\ +\ b\ \bullet\ +\ c\ \emptyset )  }}$

$\quad\quad\quad\quad\quad\quad=\quad$  $\left( \def\objectstyle{\scriptstyle} 
\def\labelstyle{\scriptstyle} 
\vcenter{\xymatrix @-1.2pc { 
& \bullet  \ar@{-}[dr]  \ar@{-}[dl] \\
\bullet  \ar@{-}[rr] && \bullet }}
\right)$ $+$ 2 $\left( \def\objectstyle{\scriptstyle} 
\def\labelstyle{\scriptstyle} 
\vcenter{\xymatrix @-1.2pc { 
& \bullet  \ar@{-}[dr]  \ar@{-}[dl] \\
\bullet  \ar@{}[rr] && \bullet }}
\right)$  $+$ $\left( \def\objectstyle{\scriptstyle} 
\def\labelstyle{\scriptstyle} 
\vcenter{\xymatrix @-1.2pc { 
& \bullet  \ar@{}[dr]  \ar@{}[dl] \\
\bullet  \ar@{-}[rr] && \bullet }}
\right)$ 

$\quad\quad\quad\quad\quad\quad\quad$ $+$  $\def\objectstyle{\scriptstyle} 
\def\labelstyle{\scriptstyle} 
\vcenter{\xymatrix @-1.2pc { 
2\ b\ (\bullet \ar@{-}[r] & \bullet)\ +\ 2\ b\ (\bullet \ar@{}[r] & \bullet) }}$ $+$ $b$ $\left( \def\objectstyle{\scriptstyle} 
\def\labelstyle{\scriptstyle} 
\vcenter{\xymatrix @-1.2pc { 
 \bullet }}
\right)$

and on the other hand\\

$\quad\quad$ $\left( \def\objectstyle{\scriptstyle} 
\def\labelstyle{\scriptstyle} 
\vcenter{\xymatrix @-1.2pc { 
& \bullet  \ar@{-}[dr]  \ar@{-}[dl] \\
\bullet  \ar@{-}[rr] && \bullet }}
\right)$ $+$ 2 $\left( \def\objectstyle{\scriptstyle} 
\def\labelstyle{\scriptstyle} 
\vcenter{\xymatrix @-1.2pc { 
& \bullet  \ar@{-}[dr]  \ar@{-}[dl] \\
\bullet  \ar@{}[rr] && \bullet }}
\right)$  $+$ $\left( \def\objectstyle{\scriptstyle} 
\def\labelstyle{\scriptstyle} 
\vcenter{\xymatrix @-1.2pc { 
& \bullet  \ar@{}[dr]  \ar@{}[dl] \\
\bullet  \ar@{-}[rr] && \bullet }}
\right)$  $\quad \leadsto \quad$  

$\quad\quad\quad\quad\quad\quad\quad$   $\left( \def\objectstyle{\scriptstyle} 
\def\labelstyle{\scriptstyle} 
\vcenter{\xymatrix @-1.2pc { 
& \bullet  \ar@{-}[dr]  \ar@{-}[dl] \\
\bullet  \ar@{-}[rr] && \bullet }}
\right)$ $+$ $x$ $\left( \def\objectstyle{\scriptstyle} 
\def\labelstyle{\scriptstyle} 
\vcenter{\xymatrix @-1.2pc {
\bullet  \ar@{-}[rr] && \bullet }}
\right)$ 

$\quad\quad\quad\quad\quad+\quad$   2 $\left( \def\objectstyle{\scriptstyle} 
\def\labelstyle{\scriptstyle} 
\vcenter{\xymatrix @-1.2pc { 
& \bullet  \ar@{-}[dr]  \ar@{-}[dl] \\
\bullet  \ar@{}[rr] && \bullet }}
\right)$  $+$ $2 y$ $\left( \def\objectstyle{\scriptstyle} 
\def\labelstyle{\scriptstyle} 
\vcenter{\xymatrix @-1.2pc {
\bullet  \ar@{-}[rr] && \bullet }}
\right)$ 

$\quad\quad\quad\quad\quad+\quad\quad$   $\left( \def\objectstyle{\scriptstyle} 
\def\labelstyle{\scriptstyle} 
\vcenter{\xymatrix @-1.2pc { 
& \bullet  \ar@{}[dr]  \ar@{}[dl] \\
\bullet  \ar@{-}[rr] && \bullet }}
\right)$ $+$ $z$ $\left( \def\objectstyle{\scriptstyle} 
\def\labelstyle{\scriptstyle} 
\vcenter{\xymatrix @-1.2pc {
\bullet  \ar@{-}[rr] && \bullet }}
\right)$ 

In order to match up the two expressions, we must take $b = 0$, and then we must set $x = y = z = 0$.

(iii)  Cameron:

$\def\objectstyle{\scriptstyle} 
\def\labelstyle{\scriptstyle} 
\vcenter{\xymatrix @-1.2pc { 
\bullet\ \times\ \bullet\ \times\ \bullet }}
$ 

$\quad\quad\quad=\quad$  $\left( \def\objectstyle{\scriptstyle} 
\def\labelstyle{\scriptstyle} 
\vcenter{\xymatrix @-1.2pc { 
& \bullet  \ar@{-}[dr]  \ar@{-}[dl] \\
\bullet  \ar@{}[rr] && \bullet }}
\right)$  $+$ 6 $\left( \def\objectstyle{\scriptstyle} 
\def\labelstyle{\scriptstyle} 
\vcenter{\xymatrix @-1.2pc { 
& \bullet  \ar@{}[dr]  \ar@{}[dl] \\
\bullet  \ar@{}[rr] && \bullet }}
\right)$ +  2 $\left( \def\objectstyle{\scriptstyle} 
\def\labelstyle{\scriptstyle} 
\vcenter{\xymatrix @-1.2pc { 
& \bullet  \ar@{}[dr]  \ar@{}[dl] \\
\bullet  \ar@{-}[rr] && \bullet }}
\right)$

$\quad \leadsto \quad$ 

Glynn: 

$\left( 
\def\objectstyle{\scriptstyle} 
\def\labelstyle{\scriptstyle} 
\vcenter{\xymatrix @-1.2pc { 
\bullet\ +\ a\ \emptyset }}
\right)$ $\times$ $\left( 
\def\objectstyle{\scriptstyle} 
\def\labelstyle{\scriptstyle} 
\vcenter{\xymatrix @-1.2pc { 
\bullet\ +\ a\ \emptyset }} 
\right)$  $\times$ $\left( 
\def\objectstyle{\scriptstyle} 
\def\labelstyle{\scriptstyle} 
\vcenter{\xymatrix @-1.2pc { 
\bullet\ +\ a\ \emptyset }} 
\right)$ 

$\quad\quad\quad=\quad$  $\left( \def\objectstyle{\scriptstyle} 
\def\labelstyle{\scriptstyle} 
\vcenter{\xymatrix @-1.2pc { 
& \bullet  \ar@{-}[dr]  \ar@{-}[dl] \\
\bullet  \ar@{}[rr] && \bullet }}
\right)$  $+$ 6 $\left( \def\objectstyle{\scriptstyle} 
\def\labelstyle{\scriptstyle} 
\vcenter{\xymatrix @-1.2pc { 
& \bullet  \ar@{}[dr]  \ar@{}[dl] \\
\bullet  \ar@{}[rr] && \bullet }}
\right)$ +  2 $\left( \def\objectstyle{\scriptstyle} 
\def\labelstyle{\scriptstyle} 
\vcenter{\xymatrix @-1.2pc { 
& \bullet  \ar@{}[dr]  \ar@{}[dl] \\
\bullet  \ar@{-}[rr] && \bullet }}
\right)$

where the subsidiary terms vanish as we have already set $a = 0$.

\smallskip

Note that in all three above cases, the no-vertex and single vertex terms have been eliminated, as we required.

\bigskip

The above method can in principle be continued to obtain further isomorphisms.   There is a separate combinatorial  problem of how many graphs can be fitted into a larger graph on a given vertex set, and one approach to its resolution is as follows.  

Denoting the number of $i$-vertex graphs that can be embedded as induced subgraphs of a $j$-vertex graph ($j \geq i$) by a double angled bracket, we have that

\begin{equation}
\quad \llangle \Gamma_j : \Gamma_i \rrangle = \sum  \llangle \Gamma_j : \Gamma_{j - 1}  \rrangle \llangle \Gamma_{j-1} : \Gamma_i \rrangle
\end{equation}

\begin{equation}
= \left[ {j-1 \choose j-1} + \ldots + {j-1 \choose 2} + {j-1 \choose 1} \right] \llangle \Gamma_{j-1} : \Gamma_i \rrangle
\end{equation}\\





where the sum is over all $j-1$ vertex graphs.  This gives a recursive combinatorial expression, with the possible use of the following known identities:-

(1)  the recurrence relation for Bell numbers:\\
\index{Bell number}%
$\Bn_n = \sum_{k=1}^{n} {n-1 \choose k-1} \Bn_{n-k} = \sum_{k=1}^{n} S(n, k)$,
\index{Stirling number ! of second kind}%

(2)  ${n \choose k} = {n-1 \choose k-1} + {n-1 \choose k}$,

(3)  $L(n, k) = \left( \frac{n-1}{k-1} \right) \frac{n!}{k!}$;\quad $L(n, k+1) = \frac{n-k}{k(k+1)} L(n, k)$;\quad and\\ $(-1)^n L(n, k) = \sum_z (-1)^z s(n, k) S(n, k)$,

where the unsigned Lah number $L(n, k)$ counts the number of ways a set of $n$ elements can be partitioned into $k$ nonempty linearly ordered subsets, and $s(n, k)$ is the Stirling number of the first kind.
\index{Stirling number ! of first kind}%

\chapter{Random Graphs as Homogeneous Cayley Objects}
~\label{homcochap}
\index{Cayley, A.}%
\bigskip  
\bigskip  

It has been observed frequently in mathematics that apparently
unrelated theories can lead to the same special objects.  This often
indicates a potential for new relationships yielding further
development of these theories.
\begin{flushright}
Frenkel, Lepowsky and Meurman, from the Introduction to 
\end{flushright}
\begin{flushright}
\textit{Vertex Operator Algebras and the Monster}
\end{flushright}

\medskip

The Eight Immortals of Chinese mythology $\ldots$ are not
strictly speaking, gods.  They are legendary personages who became
immortal $\ldots$.  These eight characters have nothing in common, and
it is hard to say how the Taoists came to make them into an almost
inseparable group.
\begin{flushright}
\textit{New Larousse Encyclopedia of Mythology}
\end{flushright}

\medskip

The flask of oil lasted for eight days.  Eight is a significant
number.  If six stands for the physical creation and seven for the
spiritual in the midst of the physical (Shabbat), then eight stands
for that which is completely outside our world.  It stands for the
World to Come.
\begin{flushright}
\textit{E. E. Dessler, Sanctuaries in Time, 1994}
\end{flushright}

\medskip

The Eightfold Path is the Buddhist Manual to Self-Enlightenment.  Much of it is to be found in other religions and philosophies $\ldots$.
\begin{flushright}
\textit{C. Humphreys, Buddhism, 1954, p.108}
\end{flushright}

\medskip

The Fathers of the Mishnah $\ldots$ were not fantastic fools, but
subtle philosophers, discovering the reign of universal law through
the exceptions, the miracles that had to be created specially and were
still a part of the order of the world, bound to appear in due time
much as apparently erratic comets are.
\begin{flushright}
\textit{I. Zangwill}
\end{flushright}

\bigskip

\section{$\mathfrak{R}$ as a Cayley Graph}
\index{graph ! Cayley}%
A graph admits a regular action of a group if and only if it is a Cayley graph for the group.  In this first section we turn our attention to the two-coloured random graph,
\index{graph ! random}%
 recalling Cameron's construction
\index{Cameron, P. J.}%
of $\mathfrak{R}$ as a Cayley Graph and deriving results allied to the following~\cite{cameron}: 

\begin{theorem} 
\label{pnca}
$\mathfrak{R}$ admits $2^{\aleph_{0}}$ pairwise non-conjugate cyclic
automorphisms.
\index{automorphism ! cyclic}%
\end{theorem}

It is a corollary of this result that the random graph $\mathfrak{R}$ has uncountably many non-conjugate cyclic automorphisms that the automorphism group $\Aut(\mathfrak{R})$ of the random graph has uncountable order.  It is also possible to reach this conclusion by viewing $\Aut(\mathfrak{R})$ as a Polish group.

The proof of this theorem uses the construction of $\mathfrak{R}$ as a Cayley graph $\Gamma(S)$, $S \subset \mathbb{N}$,~\label{SsubsetN} whose vertices are integers, and an edge $x \sim
y$ is formed for $x, y \in \Gamma(S) \Leftrightarrow |x-y| \in S$.  The
cyclic shift $x \mapsto x+1$ is an automorphism of $\Gamma(S)$.  For each
$S\ (= \{x > 0 : x \sim 0 \})$ form a graph $\Gamma(S)$, and if $S_1 \neq
S_2$ but $\Gamma(S_1) \cong \Gamma(S_2)$, then the corresponding cyclic shifts
are non-conjugate graph automorphisms for a countable graph $\Gamma$; (see Fact 2 below).
There is a natural bijection between countable zero-one sequences
\index{zero-one sequence}%
 $s = s_1 s_2 \ldots$ and $\{S\}$.  Furthermore $\Gamma(S) \cong \mathfrak{R}$ if and
only if this characteristic function of $S$ is \emph{universal}, that
is, any finite zero-one sequence
\index{zero-one sequence}%
 occurs as a consecutive subsequence.  A random zero-one sequence is universal with probability $1$ because a given sequence $\sigma$ of length $k$ has positive probability $2^{-k}$ of occurring in $k$ given consecutive positions, and we can find infinitely many disjoint sets of consecutive positions.

That $\mathfrak{R}$ has cyclic automorphisms
\index{automorphism ! cyclic}%
 shows that it is a Cayley graph for the infinite cyclic group.
\index{group ! cyclic}%
It follows from the theorem that $\Aut(\mathfrak{R})$ has $2^{\aleph_{0}}$ conjugacy classes of regular subgroups.
\index{group ! regular}%
 An $m$-valued characteristic function in the above construction would give $\mathfrak{R}_{m,\omega}$.

This method can be extended to find all the possible cycle structures for automorphisms of $\mathfrak{R}$. For example, there cannot be an automorphism which has one fixed point and permutes the others in a single cycle, since this would force the fixed point to be joined to all or no other vertices. However, there is an automorphism which fixes a point and permutes the remaining vertices in two cycles (one for neighbours, and one for non-neighbours, of the fixed point).

From the construction of cyclic automorphisms, we see that $\mathfrak{R}$ has a two-way infinite Hamiltonian path.
\index{graph ! Hamiltonian path}%
 (Choose a universal set containing 1; then $v_n$ is joined to $v_{n+1}$ for all $n$.)  

\bigskip

We will now use the following way of stating the \emph{universality} of the set $S$, namely that given any finite zero-one sequence
\index{zero-one sequence}%
 $(e_0, \ldots, e_{n-1})$, there exists a natural number $N$ such that, for $0 \le i \le n-1$, $N + i \in S$ if and only if $e_i = 1$.

\head{Fact 1}  \emph{$\Gamma(S) \cong \mathfrak{R}$ if and only if $S$ is universal.}

\begin{proof}
Suppose that $S$ is universal.  Let $U$ and $V$ be finite disjoint sets of vertices.  Let $n = \max(U \cup V) - \min(U \cup V) + 1, r = \max(U \cup V)$.  For $0 \le i \le n-1$, let $e_i = 1$ if $r - i \in U$, $e_i = 0$ if $r - i \in V$, and choose arbitrarily whether $e_i$ is 0 or 1 if $r - i \notin U \cup V$.  Choose $N$ such that $N + i \in S$ if and only if $e_i = 1$ (by universality).  Now $N + i = |(N + r) - (r - i)|$, so $N + r \sim r - i$ for $r - i \in U$ and $N + r \nsim r - i$ if $r - i \in V$.  Thus $N + r$ witnesses property $(*)$ for $(U, V)$.

Conversely,  suppose that $\Gamma(S) \cong \mathfrak{R}$; we have to show that $S$ is universal, so suppose that $s$ is a finite zero-one sequence
\index{zero-one sequence}%
 of length $n$, say.  Let $s^{\dagger}$ denote the reversal of $s$, and $t$ the concatenation of $s$ with $s^{\dagger}$.  Now let $U = \{i: t_i = 1\}$ and $V = \{i: t_i = 0\}$.  There exists a witness $z$ for $(*)$ for the pair $(U, V)$, and by reversing the previous argument we see that $z$ gives rise to a subsequence $t$ in the characteristic function of $S$ if it is negative, and $t^{\dagger}$ if is at least $2n$.  In either case, $s$ is realized.
\end{proof}

\head{Fact 2}  \emph{Let $\Gamma_1$ and $\Gamma_2$ be countable graphs admitting cyclic automorphisms
\index{automorphism ! cyclic}%
 $\sigma_1$ and $\sigma_2$.  Choose any vertex $v_i$ of $\Gamma_i$, and let
\[ S_i = \{ m > 0 : v_i^{\sigma_i^{m}} \sim v_i \} \]
for $i = 1,2$.  Then the set $S_i$ does not depend on the choice of $v_i$.  Moreover, $S_1 = S_2$ if and only $\Gamma_1 \cong \Gamma_2$ and these graphs can be identified so that $\sigma_1$ and $\sigma_2$ are conjugate in the automorphism group.}
\index{group ! automorphism}%

\begin{proof}
Since every vertex $w_i$ of $\Gamma_i$ is of the form $v_i^{\sigma_i^{r}}$ for some $r$, we see using the fact that $\sigma_i$ is an automorphism that
\[ w_i \sim w_i^{\sigma_i^{m}} \Leftrightarrow v_i^{\sigma_i^{r}} \sim v_i^{\sigma_i^{m+r}} \Leftrightarrow v_i \sim v_i^{\sigma_i^{m}}. \]

So the set $S_i$ is independent of the choice of $v_i$.

A very similar argument shows that if $\sigma_i$ is replaced by a conjugate automorphism then the same set $S_i$ is obtained.

Now clearly $\Gamma_i \cong \Gamma(S_i)$ (with the earlier notation).  If $S_1 \cong S_2$, then clearly $\Gamma_1 \cong \Gamma_2$.  Moreover, the isomorphism between them conjugates $\sigma_1$ to $\sigma_2$.  The converse is clear from the preceding remark.
\end{proof}

\head{Fact 3}  \emph{The set of universal subsets of $\mathbb{N}$ is of full measure and residual.}

\begin{proof}
The occurence of a fixed zero-one sequence
\index{zero-one sequence}%
 in the characteristic function of $S$ is clearly of full measure, open and dense.
\index{dense open set}%
  There are only countably many such finite sequences.
\end{proof}
\index{measure}%
The phenomenon of chaos
\index{chaos}%
has several characteristic features.  One is
sensitivity to initial conditions and another is randomness.  In the
study of the symbolic dynamics of the horseshoe function, there arises the
set of all sequences as maps to a $2$-element set $\mathbb{Z} \to \{0,
1\}$~\cite{batte}.  If $\sigma$ is a function on the set of sequences
that shifts each entry to the left by one position, then the sequence
$\ldots 00000000 \ldots$ will be a fixed point for $\sigma$, and
$\ldots 01010101 \ldots$ will have period $2$, and $\ldots 10010001
\ldots$ is not periodic.  S. Smale
\index{Smale, S.}%
has shown that all the dynamics of $\sigma$ are embedded in
the horseshoe, so shares the latter's infinite number of targets (and
also periodic
points though this point is irrelevant for us), and the existence of initial states with random
orbits.  The randomness of the horseshoe is thus at least to some
small extent encoded in the randomness of the Cayley graph
\index{graph ! Cayley}%
 of $\mathfrak{R}$.  We must qualify this statement by noting that because
\emph{almost all} sequences give a graph isomorphic to
$\mathfrak{R}$ there will be $2^{\aleph_{0}}$ exceptions, but these
form a null (respectively meagre) set in the sense of measure
\index{measure}%
(respectively category).  Some examples of exceptions are (i) the sequence $10101010 \ldots$
with the edge-formation rule involves joining $x$ to $y$ if $|x - y|$ is odd yields the
complete bipartite graph,
\index{graph ! random ! bipartite}%
 the two blocks being the odd numbers and the
even numbers; (ii) the sequence $01010101 \ldots$ with the joining rule
join $x$ to $y$ if $|x - y|$ is even yields the disjoint union of two
complete graphs, one on the odd numbers and the other on the even numbers; (iii) sequences
$\sigma_n$ which are zero when $n$ is even and $0$ or $1$ randomly if $n$
is odd give almost surely the random bipartite graph.
\index{graph ! random ! bipartite}%
  Placing a $1$ at the zeroth position of the above sequences corresponds to a vertex
with a loop.

\smallskip

The construction of a graph $\Gamma$ as a Cayley graph given
 in~\cite{cameron} has that $\Gamma = \Cay(\mathbb{Z}, S \cup
 (-S))$, where the vertex set is $\mathbb{Z}$ and $x \sim y
 \Leftrightarrow |x - y| \in S$ for all subsets $S \subseteq
 \mathbb{Z}$.  Then $x \sim y \Leftrightarrow (x+1)
 \sim (y+1)$ and so $\mathfrak{R}$ admits a cyclic shift as an
 automorphism.  Let $\sigma_S$ be the automorphism of $\mathfrak{R}$
 obtained by identifying $\mathfrak{R}$ with $\Cay(\mathbb{Z}, S \cup
 (-S))$.  Suppose $\sigma_S \sim \sigma_{S'}$ in
 $\Aut(\mathfrak{R})$.  Then there is an automorphism $\beta$ of 
 $\mathfrak{R}$ such that $\sigma_{S'} = \beta^{-1} \sigma_S \beta$; (see Fact 2 above).
 Pictorially we have that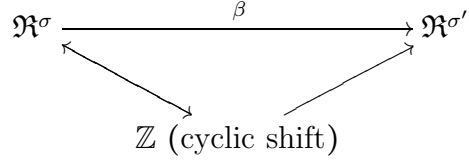
\begin{figure}[!h]$$\xymatrix{
{\mathfrak{R}^{\sigma}} \ar[rr]^{\beta} \ar[dr] && {\mathfrak{R}^{\sigma'}} \ar@{}[dl]\\
& {\mathbb{Z}\ (\text{cyclic shift})}  \ar[ur] \ar[ul]
}$$\caption{Automorphisms of $\mathfrak{R}$}  
\end{figure} 

\smallskip

Recall that a permutation $\sigma$ of $\mathfrak{R}$ is
an \emph{almost automorphism}
\index{group ! almost automorphism}%
if the set of unordered pairs of vertices $\{\{x, y\}
\in [\mathfrak{R}]^2 : f\{\sigma x, \sigma y\} \neq f\{x, y\}\}$ where $f :
[\mathfrak{R}]^2 \to \{\text{edge, non-edge}\}$ is finite. 

The set of almost automorphisms form a group.  Two conjugate cyclic automorphisms in $\Aut(\mathfrak{R})$ stay conjugate in $\AAut(\mathfrak{R})$.  Now we are going to prove the converse, that is that two almost automorphisms never become conjugate,

\begin{theorem} 
\label{pncatwo}
\index{automorphism ! cyclic}%
The $2^{\aleph_{0}}$ cyclic automorphisms
\index{automorphism ! cyclic}%
 of $\mathfrak{R}$ that are non-conjugate in the group $\Aut(\mathfrak{R})$, remain non-conjugate in the group $\AAut(\mathfrak{R})$ of almost automorphisms. 
\end{theorem}

\begin{proof}
Suppose $\sigma' = \beta^{-1} \sigma \beta$ where $\beta \in
\AAut(\mathfrak{R})$ (as opposed to $\beta \in \Aut(\mathfrak{R})$ thus
far).  Put $y_n = x_n \beta$ and $\sigma' : y_n \mapsto
y_{n+1}$.  Take any $s \in
S$.  There are infinitely many pairs $m, n$ with $|m - n| \in S$, that
is $x_m \sim x_n$.  Only finitely many edges are destroyed by $\beta$.
So $\exists m, n$ such that $|m - n| \in S$, with $y_m \sim y_n$.  So $s = |m - n| \in S'$.

A similar argument goes through for non-edges showing that $s \notin S
\Rightarrow  s \notin S'$.  Therefore $S = S'$.
\end{proof}

The essence of the above argument is that there are infinitely many pairs $m, n$ for which $|m - n|$ is any given positive integer $r$.  If $r \in S$ then $m \sim n$ for all such pairs; if $r \notin S$ then $m \sim n$ for none.  We cannot change one case to the other by adjusting only finitely many edges.

\begin{corollary} 
$\mathfrak{R}$ admits cyclic almost automorphisms
\index{group ! almost automorphism}%
which are not automorphisms.
\end{corollary}

\begin{proof}
Let $\sigma$ to be a cyclic automorphism
\index{automorphism ! cyclic}%
 of $\mathfrak{R}$.  Take a
suitable almost automorphism $\beta$ and calculate $\beta^{-1} \sigma
\beta$.  This is a cyclic almost automorphism, because $\sigma$ is
cyclic and it is also the product of almost automorphisms.  We need to show for suitable choice of $\beta$, that it is not an automorphism.  Suppose $\{1,2\}$ is an edge and $\beta :
\{1,2\} \to \{1^{\beta},2^{\beta}\}$, where $\{1^{\beta},2^{\beta}\}$ is a
non-edge and preserves all other adjacencies.  Then, $\{1^{\beta},2^{\beta}\}^{\beta^{-1} \sigma \beta} = \{2^{\beta},3^{\beta}\}$.  But $\{2,3\}$ is an edge, being the image of
$\{1,2\}$ under $\sigma$.  Therefore $\{2^{\beta},3^{\beta}\}$ is an
edge, since $\beta$ preserves all adjacencies other than $\{1,2\}$.  Therefore $\beta^{-1} \sigma \beta$ is not an automorphism.
\end{proof}

We have an example of an even stronger proposition - that
$\mathfrak{R}$ admits cyclic almost automorphisms which are not
conjugate to automorphisms in $\AAut(\mathfrak{R})$.

Let $S$ be a subset of the natural numbers $\mathbb{N}$, so that
we have an isomorphism between $\Cay(\mathbb{Z}, S \cup (-S))$ and
$\mathfrak{R}$.  As
we noted in the paragraph after the statement of Theorem~\ref{pnca}, this
is true if and only if the characteristic function corresponding to $S$ is universal.  Then
$S\setminus\{1\}$ is also universal. Define a graph $\Gamma$ with vertex set
$\mathbb{Z}$ such that $m\sim n$ if and only if either
\begin{itemize}
\item $|m-n|\in S$, $|m-n|>1$; or
\item $|m-n|=1$, $m,n\ge0$.
\end{itemize}

This is the Cayley graph of $\mathbb{Z}$ with respect to $S$, with the ``negative
half'' of one cycle of edges removed.

Claim 1: $\Gamma\cong \mathfrak{R}$. Given $U,V$, there are infinitely many $z$
correctly joined in $\Cay(\mathbb{Z},S)$, and certainly only finitely many of them can
differ from an element of $U\cup V$ by one.

Claim 2: The shift $g:n\to n+1$ is almost an automorphism. For the
non-edge $\{-1,0\}$ is mapped to the edge $\{0,1\}$; but all other
adjacencies are preserved.

Claim 3: $g$ is not conjugate to an automorphism. For consider $h^{-1}gh$,
where $h$ is an arbitrary almost automorphism.
\index{group ! almost automorphism}%
This maps $nh$ to $(n+1)h$. If it were an automorphism then $mh\sim nh$ if and only if $|m-n|\in T$ for some set $T$. But since $h$ is an almost automorphism we find $S=T=S\setminus\{1\}$, a contradiction.

An alternative argument for Claim 3: Let $a$ be an automorphism, and $h$ an
almost automorphism. Let $C(g)$ be the set of adjacencies changed by the almost automorphism $g$. 
Then $C(h^{-1}ah)$ is the symmetric difference of $C(h^{-1})$ and
$C(h)a^{-1}h$.
(For if $e$ is any pair of vertices and if $e\in C(h^{-1})$ then $e$ is
changed by $h^{-1}$, and if $eh^{-1}a\in C(h)$ then the image of $e$ under
$h^{-1}a$ is changed by $h$; if it is changed twice, it goes back to its
original state).  Since these two sets have the same size, $C(h^{-1}ah)$ has 
even cardinality. But since $C(g)$ is a singleton, $g$ cannot be conjugate 
to an automorphism.

Cyclic automorphisms
\index{automorphism ! cyclic}%
 of the random graph are investigated using methods of Baire category
\index{Baire category theorem}%
 and measure theory in~\cite{camd}.

\medskip

Recall that Henson's graph $\mathfrak{H}_k$ 
\index{graph ! Henson}%
 contains no $K_k$ but has all finite $K_k$-free graphs as induced subgraphs.  Henson
\index{Henson, C. W.}%
 observed that whilst $\mathfrak{H}_3$ admits cyclic automorphisms, the graphs $\mathfrak{H}_k$ $(k \geq 4)$ do not.  For these cyclic automorphisms, observe that $\Gamma(S)$ is triangle-free
\index{graph ! triangle-free}%
 if and only if $S$ is sum-free.
A set $S$ of natural numbers is said to be \emph{sum-free},
\index{sum-free set}%
if, for all $x,y\in S$, we have $x+y\not\in S$. The obvious example is the set of odd numbers.  Schur~\cite{schur}
\index{Schur, I}%
proved that the natural numbers cannot be partitioned into finitely many sum-free sets; this is one of the earliest results in Ramsey theory. 
\index{Ramsey theory}%
 A sum-free set $S$ cannot be universal because if $k \in S$ then no two elements of $S$ can differ by $k$, and so no finite sequence with ones in positions $k$ apart can occur in $S$.  A sum-free set with characteristic function $S$ is \emph{sf-universal}~\cite{cam2c}
\index{sf-universal}%
 if given any finite zero-one sequence
\index{zero-one sequence}%
  $\sigma = (e_1, \ldots, e_n)$, either (i) $\exists i, j (1 \leq i \le j \leq n) \text{ with } e_i = e_j = 1 \text{ and } j - i \in S$, or (ii) $\sigma$ is a subsequence of $S$.  It can be shown that for any sum-free set $S$, $\Gamma(S) \cong \mathfrak{H}_3$ if and only if $S$ is sf-universal.  Finally the following was proved in ~\cite{cam2c}.
 
\begin{proposition}[P. J. Cameron]
\index{Cameron, P. J.}%
Let $v_0, \ldots, v_n$ be vertices of $\mathfrak{H}_3$, and suppose that the map $e : v_i \to v_{i+1} (i = 0, \ldots, n-1)$ is an isomorphism of induced subgraphs.  Then there is a cyclic automorphism of $\mathfrak{H}_k$ extending $e$.
\end{proposition}

\medskip

We end the section with a result of Macpherson.
\index{Macpherson, H. D.}%
  A Polish group
\index{group ! Polish}%
   $G$ has a \emph{cyclically dense conjugacy class}
\index{cyclically dense conjugacy class}%
if $\exists g, h \in G$ so that $\{g^{n} h g^{-n}\}_{n \in \mathbb{Z}}$ is dense in $G$.  In this case $G$ is topologically $2$-generated, that is it has a dense $2$-generated subgroup.  In~\cite{macph} it is shown that $\Aut(\mathfrak{R})$ has this property.

\section{A Construction of $\mathfrak{R^{t}}$ as a Homogeneous Cayley Object}
\label{sec2ch6}
\index{Cayley object}%

We use a two-step approach to achieving the stated aim of the section title, using the modular group
\index{group ! modular}%
\[ \gimel = \PSL(2, \mathbb{Z}) = \langle \sigma, \rho : \sigma^{2} =
\rho^{3} = 1 \rangle. \]~\label{gimel}~\label{PSL}
We will use the Hebrew character $\gimel$ (`gimel') to
denote the modular group because we have already reserved the usual
symbol $\Gamma$ for
graphs.  That $\gimel \cong C_{2} \ast C_{3}$ so that it has so few
defining relations, implies that it has many
epimorphic images, one of which is $\Sym(3)$.  For example $\PSL(2, \mathbb{Z}) \cong B_{3} / \zeta(B_3)$ where the braid group
\index{group ! braid}%
 $B_{3} = \langle x, y\ |\ xyx = yxy \rangle$ maps surjectively
onto $\PSL(2, \mathbb{Z})$
\index{group ! modular}%
 via
$x \mapsto$
$ \begin{pmatrix} 
1&1\\
0&1\\
\end{pmatrix}$, 
$y \mapsto$
$ \begin{pmatrix} 
1&0\\
-1&1\\
\end{pmatrix}$, and the centre of $B_{3}$ is generated by $(xy)^3$. 
The Kurosh subgroup theorem
\index{Kurosh subgroup theorem}%
 dictates that subgroups of $\gimel = C_{2} \ast C_{3}$ are free products
\index{group ! free product}%
  of cyclic groups of orders $2, 3$ or $\infty$, so ruling out almost all
other possible subgroups.  The determination of the finite simple
quotients of $\PSL(2, \mathbb{Z})$
\index{group ! modular}%
 is a well-studied
area~\cite{shalev}.  Its above free product
presentation implies that such quotients, said to be
$(2, 3)$-$\emph{generated}$, are generated by elements $x, y$ such that
  $x^{2} = 1 = y^{3}$.  As in random graph theory, this is often proved
probabilistically by demonstrating that almost all groups in a certain
class have this property.

First we will show that there is an extension of a certain free product
\index{group ! free product}%
 $F$ by a group isomorphic to $C_{3}$, the extended group being isomorphic to $\gimel$, and the extending group cyclically
permuting three colour classes.  Then we show that the group $F =
\ker(\gimel \to C_{3})$ has $\mathfrak{R^{t}}$ as a Cayley graph by assigning its inverse pairs of elements randomly to three appropriately chosen classes.  

The main result of this section will be that $\mathfrak{R^{t}}$ is a Cayley object
\index{Cayley object}%
 for the group $F \cong C_{2} \ast C_{2} \ast C_{2}$ in such a way that $C_3$ permutes the $3$ colours cyclically, and $\gimel \cong F . C_{3}$.

The kernel of the homomorphism $\gimel \to C_{3} = \langle \tau
\rangle, \sigma \mapsto 1, \rho \mapsto \tau$,~\label{tau} which we have called $F$, can be obtained~\cite{johns} using Schreier's method, which we now pr\'ecis.
\index{Schreier's method}%
\index{Schreier, O.}%

Consider the group $F_{1} = \langle x, y, \ldots \rangle$, with $H_{1}$ a fixed
subgroup of $F_{1}$.  A set $\emptyset \neq S \subseteq F_{1}$ has the
Schreier property if $w = x_{1} \ldots x_{n} \in S \Rightarrow x_{1}
\ldots x_{n-1} \in S$, where $l(w) = n \geq 1$.  A Schreier
transversal for $H_{1}$ in $F_{1}$ is a (right) transversal for
$H_{1}$ with the Schreier property.  The algorithm for getting the
generators of the kernel is given by the following theorem~\cite[p.~17]{cam1}:

\bigskip

\begin{theorem} 
\label{scrtm}
Let $H$ be a subgroup of index $n$ in $G = \langle
  g_{1} \ldots g_{r} \rangle$, with coset representatives $x_{1}
  \ldots x_{n}$, where $x_{1} = 1$ is the representative of $H$.  Let
  $\bar{g}$ be the representative of the coset $Hg$.  Then 
\[ H = \langle x_{i} g_{j} (\overline{x_{i}g_{j}})^{-1} : 1 \leq i \leq n,
  1 \leq j \leq r \rangle. \]
\end{theorem}

By inspection, our Schreier transversal is $\{1, \rho, \rho^{2}\}$.
From the above theorem, the generators of $F$ are given by
\begin{figure}[!h]
$$\xymatrix{
& {} & {\sigma} & {\rho} \\
& {1} & {1.\sigma.1^{-1}} & {1.\rho.\rho^{-1}} \\
& {\rho} & {\rho.\sigma.\rho^{-1}} & {1.\rho^{2}.\rho^{-2}} \\
& {\rho^{2}} & {\rho^{2}.\sigma.\rho^{-2}} & {1.\rho^{3}.1^{-1}} \\
}$$\caption{Generators of $F$}  
\end{figure} 

Two of the generators are trivial because we have taken a Schreier
transversal and the third is $\rho^{3} = 1$.  So $F = \langle a, b, c \rangle$, $a = \sigma, b = \rho \sigma
\rho^{-1}, c = \rho^{2} \sigma \rho^{-2}$. 

The group $F$ is generated by three
involutions.  It is known~\cite{iver} that $C_{2} \ast C_{2}
\ast C_{2}$ is the unique index $3$ subgroup of $\gimel$, so $F$ is
the free product
\index{group ! free product}%
 $C_{2} \ast C_{2} \ast C_{2}$.  (We mention in passing that
$\SL(2, \mathbb{Z})$ has two index $3$ subgroups which are conjugate
but not identical and so are not normal subgroups~\cite[p.~79]{hirzebruch}).  This can also be shown as follows.  Call the
quotient $H$, that is $|\gimel : F| = |H| = |\langle \tau
\rangle| = 3$, where $\tau$ is in the image of the homomorphism $\gimel
\to C_{3}: \rho \mapsto \tau$, with kernel $F$.  Certainly
$F$ contains $\sigma$ and its conjugates $\tau \sigma \tau^{-1}$
and $\tau^{-1} \sigma \tau$.  Setting $\rho = \tau
\sigma$ gives that $\gimel / F$ is cyclic of order $3$ generated by
$\tau$.  Schreier's method can also be used to give the relations satisfied by the
generators as well as the generators themselves.  Yet another way to find the structure of $F$ is geometric.  Recall that $F= C_{2} \ast C_{2} \ast
C_{2}$ is a triangle group
\index{group ! triangle}%
 with all three triangle vertices at infinity~\cite{lynd}.  The modular and triangle groups are examples of
Fuchsian groups,
\index{group ! Fuchsian}%
 being discrete subgroups of $\PSL(2, \mathbb{R})$.
Thus if $\Delta$ denotes the hyperbolic triangle with vertices $e^{\pi i / 3},
e^{2 \pi i / 3}, \infty$, then $\tau^{-1} \Delta \cup \Delta \cup \tau
\Delta$ represents a fundamental domain~\cite[Lemma V.1.4]{iver}.  So there are side transformations $\sigma, \tau \sigma \tau^{-1},
\tau^{-1} \sigma \tau, \tau^{3}$ and one cycle relation: $(\tau^{-1} \sigma
\tau) \sigma (\tau \sigma \tau^{-1})\tau^{3} = 1$; the latter
eliminates the generator $\tau^{3}$.  Since all three transformations
are involutions, we conclude that $F$ is the free product
\index{group ! free product}%
 of three involutions.

To summarize $H \cong C_{3} \vartriangleleft \Aut^{*}(\mathfrak{R^{t}})
/ \Aut(\mathfrak{R^{t}}) \leq \Sym(3)$.  Because $\gimel$ has a normal
subgroup $F$ and $\langle  \tau \rangle \cong C_{3}$ is a complement for $F$ the
extension of $F$ by $H$ splits.
\index{group ! split extension}%
  There is then an action of $\gimel$ on cosets of $\langle \tau \rangle$ with $F$ a
regular normal subgroup.
\index{group ! regular normal subgroup}%
We need to construct a copy of
$\mathfrak{R^{t}}$ on this set such that $F \leq \Aut(\mathfrak{R^{t}})$ and $\tau$ cyclically permutes the colour classes.

A mathematical structure such as a graph is a \emph{Cayley Object}
\index{Cayley object}%
~\cite{cam11},
$O$, for the group $G$ if its point set are the elements of $G$ and right
multiplication by any element of $G$ is an automorphism of $O$.  A Cayley graph
\index{graph ! Cayley}%
for a group $F$ takes the form $\Cay(F,
S)$, where $S$ is an inverse-closed subset of $F \backslash \{1\}$ with vertex
set $F$ and edge set $\{\{g, sg\}: g \in F, s \in S\}$.  As $F$ is countable, we can enumerate the inverse pairs of non-identity
elements of $F$ as $\{g_{1}, g_{1}^{-1}\}, \{g_{2}, g_{2}^{-1}\}, \ldots$.  Baire
category theory
\index{Baire category theorem}%
 can sometimes be used for building a homogeneous
Cayley object
\index{Cayley object}%
for a group $G$, by showing that almost all
$G$-invariant objects of the required
type, that is a residual set of them, are homogeneous.  In order to make sense of the notion of a residual set of Cayley graphs for $F$, specify Cayley graphs by paths in a ternary tree,
    whereby the three descendents of a node or vertex at level $n$
    correspond to including or excluding the inverse pair $\{g_{n+1}, g_{n+1}^{-1}\}$
      in one of the colour classes of $F$; call these colour classes
$X, Y, Z$.
\begin{figure}[!h]
$$\xymatrix{
& {\bullet}^{\{g_{n+1}, g_{n+1}^{-1}\} \in X } \ar@{-}[dr] & {\bullet}^{\{g_{n+1},
g_{n+1}^{-1}\} \in Y } \ar@{-}[d] & {\bullet}^{\{g_{n+1}, g_{n+1}^{-1}\} \in Z }
\ar@{-}[dl]\\
&& {\bullet} \ar@{-}[d]\\
&& {\bullet}
}$$\caption{Colour classes of $F$}
\label{ccif}
\end{figure}
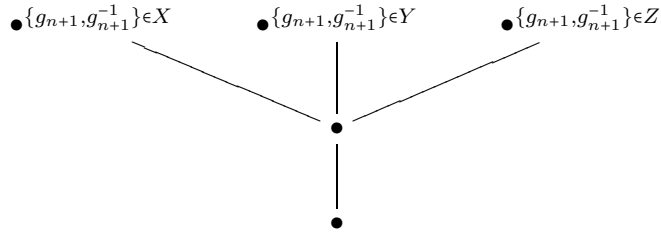 

The set of three-coloured Cayley graphs is identified with the set of paths in the ternary tree; see remark (b) after the next theorem. 
In Baire category theory
\index{Baire category theorem}%
 if an object is specified by a
countable sequence of choices, then the existence of one such object
with a given property $\mathcal{P}$ can be proved by showing that
$\mathcal{P}$ holds for `almost all' choices.  Let $\mathcal{P}(T)$
denote the set of paths of countable length starting at the root of a
tree $T$, whose nodes at height $n$ are labelled by structures on
$\{1, \ldots, n\}$.  We define the distance between distinct paths $p$
and $p'$ to be $f(n)$, where $n$ is the height of the last node at
which $p$ and $p'$ agree, and $f$ is any strictly decreasing function
tending to zero.  The complete metric space 
\index{metric space ! complete}%
 to which the Baire
category theorem will be applied arises from paths in rooted trees of
countable height.  A Cauchy sequence
\index{Cauchy sequence}%
in this space is a sequence of
paths agreeing on increasingly longer initial segments, and so has a
unique limiting path.  In this way we achieve a complete metric space.  We
require an interpretation of openness and denseness to formulate
residual sets in this space.  An \emph{open ball}
\index{open ball}%
consists of all paths in the tree containing a given vertex.  A set $S$ of paths is open (or
\emph{finitely determined}) if each path in $S$ has a vertex such that
every path through this vertex is in $S$.  A set $S$ is \emph{dense}
(or \emph{always reachable}) if it meets every open ball, i.e. if all
vertices lie on some path in $S$.  The triality graph
\index{graph ! triality}%
is the Fra\"{\i}ss\'e limit
\index{Fra\"{\i}ss\'e limit}%
of the class of all appropriately
defined finite $3$-coloured graphs, so with a countable vertex set the
isomorphism class of $\mathfrak{R^{t}}$ is residual in the set of $3$-coloured
graphs on $\mathbb{N}$.  Thus it makes sense to talk of a residual set of
Cayley graphs for $F$.  The metric space has an underlying
topological space,
\index{topological space}%
 where the topology is on the collection of inverse
closed subsets of $F$.

The strategy of the construction is as follows:  we partition the
non-identity inverse pairs of elements of $F$ into $3$ classes $X, Y,
Z$ such that (a) $g \in X \Rightarrow g^{-1} \in X$ (b) $g \in X
\Rightarrow g^{\tau} \in Y, g^{\tau^2}, \in Z$, etc.  Give $\{g, h \}$ colour red if $gh^{-1} \in
X$, green if $gh^{-1} \in Y$ and blue if $gh^{-1} \in Z$.  In other words
$x \in X \Rightarrow \rho x \rho^{-1} \in Y, \rho^2 x \rho^{-2} \in
Z$, or in terms of elements of $F$, $word(x, y, z) \in X \Rightarrow word(y, z, x)
\in Y \Rightarrow word(z, x, y) \in Z$.  Choose $X, Y, Z$ in such a
way as to ensure that the resulting object is isomorphic to
$\mathfrak{R^{t}}$.  Finally we can state our result.

\begin{theorem} 
\label{mdgptm}
Let $F = C_{2} \ast C_{2} \ast C_{2}$.  The set of Cayley graphs
for $F$ which are isomorphic to $\mathfrak{R^{t}}$ is residual in the
set satisfying (a) and (b) above.
\end{theorem}
\index{graph ! Cayley}%

\begin{proof}
Let $\chi(A, B, C)$ be the set of Cayley graphs satisfying the
$1$-point extension property (equivalently the I-property
\index{I-property}%
 or
homogeneity) of $\mathfrak{R^{t}}$:  for all finite graphs $A, B$ with $A \subseteq B$ and $|B| = |A| + 1$,
and any $C \subset F$ with $|C| = |A|,\ (C \cong A) \Rightarrow (\exists
z) (C \cup \{z\} \cong B).$  There are countably many such conditions, so by proving that
$\chi(A, B, C)$ is open and dense, we demonstrate that our set is a
countable intersection of dense open sets, and therefore residual.
\index{dense open set}%

We must verify that $\chi(A, B, C)$ is open and dense.  That $z$ witnesses $\chi(A, B, C)$ depends only on the colours assigned to $x y^{-1}$ for
$x, y \in C\cup\{z\}$, that is a finite number of choices.  So
$\chi(A, B, C)$ is the finite union of open balls,
\index{open ball}%
so is open.  To prove $\chi(A, B, C)$ is dense, assume that the
first $n$ inverse pairs comprising a finite subset $S_{0}$ of $S$ have
been chosen, including all $x y^{-1}$ for $x, y \in C$, thereby determining the
structure of $C$.  If $C \ncong A$, the argument is finished.
Otherwise, we must next choose the above implication.  At each step of
the $(*_t)$-condition, we have three finite disjoint subsets of vertices $U, V, W \subset F$, with $C
= U \cup V \cup W$, and we must find some $z \in F \backslash C$ and
extend the colouring so that $zu^{-1} \in X$ for every vertex $u$ in the finite set U, $zv^{-1} \in Y$
for every vertex $v$ in the finite set $V$, $zw^{-1} \in Z$ for every
vertex $w$ in the finite set $W$.  That is we colour the edges
$zu^{-1}$ red, the edges $zv^{-1}$ green and the edges $zw^{-1}$ blue.
We must eliminate certain elements that could not satisfy the ($*_{t}$) condition:

(i)  We claim that $\exists z$ such that no $zc^{-1} (c \in C)$ has yet
been assigned one of the three colours.  For denoting by $\Phi$ the
finite set whose elements have
been assigned colours, then $zc^{-1} \in \Phi \Rightarrow
zc^{-1}=\phi \in \Phi \Rightarrow z = \phi c$.  Since there are only
finitely many $\phi$ and $c$, only finitely many $z$ are excluded.

(ii)  We must also disqualify any $z$ which would be forced to have the same adjacency to any
two of $u \in U, v \in V$ and $w \in W$.  These have the form $zu^{-1}
= (zv^{-1})^{-1}$ with $u, v \in C$, that is $(zu^{-1})^{2} =
vu^{-1}$.  For a given $u$ and $v$ these have the $z$ in question as
translates of non-principal square-root sets~\cite{cam11},
\index{square-root set}%
that is $\{z : (zu^{-1})^{2} = vu^{-1}\} = (vu^{-1})^{\frac{1}{2}} u$.  But elements of
free products
\index{group ! free product}%
 are reduced words in the generators and we want each
word to have at most one square root; non-cyclically reduced words are conjugates of reduced words.  By conjugation if necessary, we
can assume that the word does not begin and end with the same letter.
Note that conjugation does not change the number of square roots.  It
then follows that each word has at most
one square root, for if $word(u, v, w) = word'(u, v, w)\ word'(u, v,
w)$ then $(word(u, v, w))^{\frac{1}{2}} = word'(u, v, w)$.  Any cancellation
here would require the last letter of the first $word'$ to be the same
as the first letter of the second $word'$, contrary to assumption.
Again only finitely many such $z$ are excluded.

(iii)  Next consider $z$ for which $(zu^{-1})^{\tau} = zv^{-1}$, and
$u, v \in C$.  Suppose that $z$ begins with a certain letter $a$.  If not all of $z$
cancels into $u^{-1}$ or $v^{-1}$ then $zu^{-1}$ will also begin with
$a$, as will $zv^{-1}$.  However under the action of $\tau$,
$(zu^{-1})^{\tau}$ will be a word whose first letter is $b$, contradiction.  So in order to
avoid getting $(zu^{-1})^{\tau} = zv^{-1}$, we simply choose all our
$z$'s to have reduced length greater than any element so far put into $C$; so only
finitely many such $z$ are excluded.

(iv)  Finally we must show that given a finite set $C$,
there is $z$ such that $(zu^{-1})^{\tau} \neq (zv^{-1})^{-1}$ for $u,
v \in C$.  So we ask
when does $(zu^{-1})^{\tau} = vz^{-1}$, i.e. $z^{\tau} = vz^{-1}u^{\tau}$?
Take $z$ to be a repeated series of letters $abc$, such that it is at
least twice as long as any element of $C$.  So $z = \ldots abcabcabc
\ldots$ and $z^{-1} = \ldots cbacbacba \ldots$.  Then clearly, with $z
= (abc)^{n}$ and $n$ sufficiently large, the required inequality
holds.  There are infinitely many $z$.

So for any three finite disjoint sets $U, V$ and $W \subset
F$, there exists $z \in F \backslash C$ such that, for all $u \in U, v \in V, w
\in W$, we have $(zu^{-1})^{\tau} \neq (zv^{-1})^{-1}$,
$(zu^{-1})^{\tau^{2}} \neq (zv^{-1})^{- \tau}$, and $(zu^{-1}) \neq
(zv^{-1})^{- \tau^{2}}$ and also $zu^{-1}, zv^{-1}, zw^{-1}$ have not
yet been assigned to $X, Y, Z$.  Now assign all $zu^{-1}$ to $X$, all
$zv^{-1}$ to $Y$ and all $zw^{-1}$ to $Z$; then the induced subgraph
on $C \cup \{z\}$ is isomorphic to $B$, no matter what further
choices we make.

Our argument that $(zu^{-1})^{\tau} \neq (zv^{-1})^{-1}$ provides
another reason why any group larger than $C_{3}$ acting on the colours
would cause the construction to
fail, for even an extra $C_{2}$ involution invalidates the argument in (iv).

In total, we have excluded only finitely many $z$, so infinitely many
choices remain for us to satisfy the $1$-point extension hypothesis of
the theorem.  We have proved that the set $\chi(A, B, C)$ of Cayley
graphs is dense.  The Baire category theorem
\index{Baire category theorem}%
 implies that the
intersection of all the sets $\chi(A, B, C)$ is residual and hence
non-empty.  It follows that $\gimel$
acts on the vertices of $\mathfrak{R^{t}}$ with $F$ as a regular
subgroup because it is transitive, and the stabilizer of a vertex is
the identity.  Hence result.
\end{proof}

It follows that $\PSL(2, \mathbb{Z}) < \Aut(\mathfrak{R^{t}})$.  

It is open whether or not $\PGL(2, \mathbb{Z}) < \Aut(\mathfrak{R^{t}})$, where the \emph{extended modular group} $\PGL(2, \mathbb{Z}) := \GL(2, \mathbb{Z}) / \{ \pm I\}$ and $\PSL(2, \mathbb{Z})$ is a subgroup of index 2 in $\PGL(2, \mathbb{Z})$~\cite{jonesth}.

The group $C_{2} \ast C_{2} \ast C_{2}$ has been sufficiently important in exhibiting the structure of $\mathfrak{R^{t}}$ as a Cayley graph that in the interest of potential future extensions of this section we make some tangential remarks.

\head{Remarks}

\begin{itemize}
\item[(a)] Exhibiting a multicoloured random graph as the Cayley graph of any countable group requires the satisfaction of the square root condition, pairwise over all the colours.  This condition is derived in~\cite{cameron}, but see also Theorem~\ref{sqrootthm} later on in this chapter.

\item[(b)] The \emph{modular group of first level}
\index{group ! modular}%
 given
by $\gimel = \gimel(1) =
 \PSL(2, \mathbb{Z})$ $= \SL(2, \mathbb{Z}) / \{\pm I\}$, where $\SL(2,
\mathbb{Z}) := \{ \frac{az + b}{cz + d}\ |\ a, b, c, d \in \mathbb{Z},
ad-bc=1 \}$ has subgroups arising as follows.  The
\emph{modular group of level n}, also called the \emph{principal modular subgroup of level n}, $\gimel(n) \vartriangleleft
\gimel(1)$ is defined by the requirement that $[ab/cd]$ be congruent
to the identity $[10/01] \pmod n$.  If $\gimel(1)$ is generated by the
 substitutions $A: z \to z + 1$ and $B: z \to -1/z$, then $\gimel(2)$
 is the free group on $A^{2} = [12/01]$ and $BA^{-2}B = [10/21]$.  It
 can be shown~\cite{mcmo} that $\gimel(1) / \gimel(2) \cong
\Sym(3)$.  Therefore we see that an index $6$ subgroup of $\gimel$
exists, but for Theorem~\ref{mdgptm} to work we require the action of $C_{2} \ast C_{2} \ast
C_{2}$ rather than $\gimel(2)$.

\item[(c)] J. H. Conway
\index{Conway, J. H.}%
 found a way of determining whether or not the plane can be tesselated by given tiles.  A relation in a group is a word in the generators which represents the identity element 1.  A relation in an infinite finitely presented group is a tile and a planar region can be tiled only if the group element describing the boundary of the region is 1.  The Cayley graph
\index{graph ! Cayley}%
$\Gamma(G)$ of a group $G$ extends to a 2-complex $[\Gamma(G)]^2$~\label{Gamma(G)]^2} 
\index{two-complex@$2$-complex}%
where at each relevant vertex we attach as many discs as there are relations $r_i$, so that the disc boundary traces out the words $r_i$.  In~\cite{thurston} Thurston
\index{Thurston, W. P.}%
shows how an element $g(\pi) \in F = C_{2} \ast C_{2} \ast C_{2}$ is determined by a path $\pi$ in the 1-skeleton of a hexagonal grid.

\item[(d)] A. W. M. Dress
\index{Dress, A. W. M.}%
defined~\cite{dress} a \emph{(thin) chamber system} (of rank 2)
\index{chamber system}%
to be just a set $\mathsf{C}$ on which the above free Coxeter group $F = \langle u, v, w | u^2 = v^2 = w^2 = 1 \rangle$ acts from the right.  An isomorphism between two chamber systems $\mathsf{C}$ and $\mathsf{C}'$ is a bijection $f : \mathsf{C} \to \mathsf{C}'$ such that $f (\mathsf{C} g) = f(\mathsf{C}) g$ for all $C \in \mathsf{C}$ and $g \in F$.  The relation between tilings and chamber systems is established and it is proved that there is a $1$--$1$ correspondence between isomorphism classes of non-degenerate tilings of the Euclidean plane
\index{Euclidean plane}%
and isomorphism classes of infinite chamber systems satisfying given conditions.  A metric can be defined on $\mathsf{C}$ by $d(C, C') := \inf(l(g) | g \in F,\ Cg = C')$ $\forall\ C, C' \in \mathsf{C}$, and $d(C, C') = \infty$ if and only if $C, C'$ lie in different $F$-orbits.
\end{itemize}

Recall that to specify a transitive action of $\gimel$ on a finite vertex set $V$ with pointwise stabilizer $\gimel_0$ of one distinguished vertex $v_0$ is equivalent to specifying a finite index subgroup $\gimel_0$ of $\gimel$, where the vertices of $V$ can be identified with the cosets of $\gimel_0 \backslash \gimel$, and $v_0$ with the coset $\gimel_0$.  There is another work~\cite{mckayseb} which studies the regular action of the
    whole modular group as opposed to the index $3$ subgroup in our
    previous theorem.  A coset multigraph consisting of triangles and single edges
    is given for a permutation
    representation of $\gimel = \PSL(2, \mathbb{Z})$
\index{group ! modular}%
     on the cosets of a point
    stabilizer $\gimel^{\mu}$, for every torsion-free
 \index{group ! torsion-free}%
     genus zero
    congruence subgroup
    \index{group ! congruence subgroup}%
 $\gimel^{\mu}$ of index $\mu$ in $\gimel$.  The Cayley graph of $\PSL(2, \mathbb{Z}) = \langle
    \sigma, \rho : \sigma^{2} = 1 = \rho^{3} \rangle$ is the infinite
    free trivalent tree with each node replaced by a positively
    oriented triangle corresponding to $\rho$.  Collapsing the triangles
    to points gives the set of vertices of our above constructed Cayley graph.  In a sense this is a special case of the known result that if $H \le G$, $g \in G, g^2 = 1$ and $|H : H \cap H^{g}| = m$ then $G$ acts vertex and edge-transitively on an $m$-valent graph.  (More is known from applications of group amalgams to algebraic graph theory~\cite[p.~417]{ivanovs}.  Let $x$ and $\{x, y\}$ be a vertex and an edge of the graph $\Gamma$ with respective stabilizers $H$ and $K$ in $G = \Aut(\Gamma)$.  If $G$ is vertex-transitive then $\{x, y\}$ is invertible so $|K : H \cap K| = 2$.  If $\Gamma$ is trivalent then $|H : H \cap K| = 3$.  If $\Gamma$ is connected it follows that $G = \langle H, K \rangle$.  The $G$-action on $\Gamma$ is faithful, so $H$ and $K$ have no common nontrivial normal subgroups).

\smallskip

Recall that the \emph{principal congruence subgroup}
\index{group ! principal congruence subgroup}%
 of level $n \ge 0$ is a subgroup of $\SL(2, \mathbb{Z})$
\index{group ! SL()@$\SL(2, \mathbb{Z})$}%
defined by, $\gimel(n) = \Bigg\{\begin{pmatrix} 
a&b\\
c&d\\
\end{pmatrix} \in \gimel\ |\ a \equiv 1, b\ \equiv 0, c\ \equiv\ 0, d\ \equiv 1 \pmod{n} \Bigg\}$.  

A \emph{congruence subgroup}
\index{group ! congruence subgroup}%
 is a subgroup of $\SL(2, \mathbb{Z})$
containing $\gimel(n)$ for some $n$, for example, $\gimel_{0}(n) = \Bigg\{\begin{pmatrix} 
a&b\\
c&d\\
\end{pmatrix} \in \gimel\ |\ c\ \equiv\ 0 \pmod{n} \Bigg\}, n \in \mathbb{N}$.  By
definition the sequence
\[ 1 \to \gimel(n) \to \SL(2, \mathbb{Z}) \to \SL(2, \mathbb{Z}/n\mathbb{Z}) \]
is exact.

\section{Lattice Constructions of Random Graphs}

A lattice $\mathbb{L}$~\label{mathbb{L}}
\index{lattice}%
 in a real finite-dimensional vector space $V$
is a subgroup of $V$ satisfying one of the following equivalent
conditions:

(i)  $\mathbb{L}$ is discrete and $V / \mathbb{L}$ is compact;

(ii)  $\mathbb{L}$ is discrete and generates the $\mathbb{R}$-vector
space $V$;

(iii)  $V$ has an $\mathbb{R}$-basis $(e_1, \ldots, e_n)$ which is a
$\mathbb{Z}$-basis of $\mathbb{L}$ (that is $\mathbb{L} = \mathbb{Z}
e_1 \oplus \ldots \oplus \mathbb{Z} e_n)$.

So a lattice in $\mathbb{R}^{n}$ is the set of integral linear
combinations of $n$ linearly independent vectors.  The lattices generated by vectors of norm $1$ are $\mathbb{Z}^n$ with the standard bilinear form.  (This norm is the squared norm as we shall be dealing with positive-definite lattices).  Those generated by vectors of norm $2$, which are called \emph{roots}, are the root lattices
\index{lattice ! root}%
have a classification theorem, whilst those of norm $3$ do not as yet.  Even unimodular
lattices, those having one
point per unit volume and in which every squared length is even, exist only when the dimension is a multiple of $8$, and the
unique $8$-dimensional one is denoted $E_8$ because it is the root
lattice of the corresponding Lie algebra.  In any even
positive-definite lattice,
\index{lattice ! positive-definite}%
 the \emph{roots}
\index{lattice ! roots}%
 are those vectors with
minimal non-zero norm of $2$; they are special because whereas
reflecting through
vectors will not in general map the lattice to itself, reflecting
through roots is a lattice automorphism.  Those lattices that are
spanned by their roots are precisely the orthogonal direct sums of the so-called
$ADE$ lattices, that is \emph{Dynkin diagrams}.
\index{lattice ! ADE}%
\index{Dynkin diagram}%
There is here a direct link between lattices and graphs.  For example,  if $V = \mathbb{L} \otimes_{\mathbb{Z}} \mathbb{R}$ is a real vector space, $\{\alpha\} \in \mathbb{L}$ the set of roots of $\mathbb{L}$, $\Delta$ the set of roots corresponding to the walls of a particular chamber (connected component of $V$), then we can build a graph $\Gamma_{\mathbb{L}}$ whose vertices are the elements of $\Delta$ and edges between those vertices $\alpha, \beta$ for which $(\alpha, \beta) = -1$.  

\medskip

We will prove a general theorem connecting random graphs of any number of colours
and integral lattices.  Our main aim however is to find links to the
primary even unimodular lattices and with this in mind we prove a
$1$--$1$
correspondence between $\mathfrak{R^{t}}$ and the Leech lattice
\index{lattice ! Leech}%
which
is the only even unimodular lattice in $24$ dimensions with no
vectors of squared length $2$.  This is one of the $24$ even unimodular
lattices in $24$ dimensions classified by Niemeier;
\index{Niemeier, H-V}%
 it represents the densest packings in $24$ dimensions.  We shall refer to a spherical layer of vectors in the lattice that are
all the same distance from the origin as a \emph{shell}.
\index{lattice ! shell}%

Now for the constructions, firstly working with a general colour set
of size $m$.

A Cayley graph
\index{graph ! Cayley}%
for the lattice $\mathbb{L}$, here regarded as a
subgroup of a finite-dimensional vector space, takes the form
$\Cay(\mathbb{L}, S)$, where the vertices are $l \in \mathbb{L}$ and
the inverse-closed subset
$S \subseteq \mathbb{L} \backslash \{1\}$ gives the edge set $\{\{l,
l + s\}: l \in \mathbb{L}, s \in S\}$ written additively.  

We not only want to produce a Cayley graph for the lattice, but more
specifically one satisfying the extra condition that it is invariant under the automorphism
group of the lattice, for then the lattice isometries will be Cayley
graph automorphisms.  That is, denoting the
additive group of the lattice by $\mathbb{L}^{+}$, we are embedding $\mathbb{L}^{+} \sd
\Aut(\mathbb{L})$ into $\Aut(\mathfrak{R}_{m,\omega})$ as a
permutation group.  The automorphism groups
\index{group ! automorphism}%
 of the random graphs for any number of colours are uncountable.

Partition the elements of $\mathbb{L}$ into $m$ disjoint classes $U_1,
\ldots, U_m$ corresponding to colours $c_1, \ldots, c_m$.  One way to
do this which yields the extra condition is to partition $\mathbb{N}$ into $m$ disjoint classes $U_1, \ldots, U_m$,
and letting $\tilde{U_i}$ be the set of lattice vectors $v$ with
$\parallel v \parallel \in U_i$ give $\{x, y\}$ colour $c_i$ if $y - x
\in \tilde{U_i}$.  We describe this procedure as ``assigning shells in
$\mathbb{L}$ to $m$ boxes.''

We must choose our $\{U_i\}$ so as to ensure that in the sense of Baire category,
\index{Baire category theorem}%
 almost all colourings result in an object isomorphic to $\mathfrak{R}_{m,\omega}$.

\begin{theorem} 
\label{rglttm}
Let $\mathbb{L}$ be a lattice in $\mathbb{R}^d$, for $d \geq 2$ and
take an integer $m \geq 2$.  Assigning the shells in $\mathbb{L}$ to $m$ boxes, the set of Cayley
graphs
\index{graph ! Cayley}%
for $\mathbb{L}$ with lattice vectors as vertices and vector
pairs as edges, which are isomorphic to the random $m$-coloured graph
is residual.
\end{theorem}

\begin{proof}
The topological space in which the Baire category argument takes place
is the space of assignments of shells to boxes.  Construct a
tree whose vertices at level $n$ are assignments of the first $n$
shells to boxes;  the descendents of the vertices of level $n$ are
obtained by assigning the $(n+1)^{st}$ shell to a box.  Paths in this
tree now correspond to assignments of all shells to boxes.  Our approach is to partition the distances between the lattice vectors
into $m$ colour (equivalence) classes.  We throw vectors, one shell at
a time, into a coloured box.  For all finite graphs $A, B$ and $C$, let $\chi(A, B, C)$ be the
set of Cayley graphs satisfying the $1$-point extension property for countable
$m$-coloured random graphs
\index{graph ! random ! $m$-coloured}%
 $\mathfrak{R}_{m,\omega}$: for all finite graphs $A, B$ with $A \subseteq B$ and $|B| = |A| + 1$,
and any $C \subset \mathbb{L}$ with $|C| = |A|,\ (C \cong A) \Rightarrow (\exists
z) (C \cup \{z\} \cong B).$

The set $\chi(A, B, C)$ is open because only those finite number of
vertices in $C$ whose distances fall within one of the $m$
colour classes must be witnessed by any chosen $z$, so the set of
distances is finitely determined.  To prove that
$\chi(A, B, C)$ is dense assume that the first $n$ vertex pairs
comprising a finite subset of $S$ have been chosen and that all
edges between vertex pairs of $C$ have been coloured.  If $C \ncong A$
we are done, otherwise we must demonstrate that the $1$-point extension
property holds.  Choose $m$ finite disjoint subsets of
vertices $\{U_i\}^m_{i=1} \subset \mathbb{L}$, with $C \supseteq \bigcup_i U_i$ and
$z \in \mathbb{L} \backslash C$ such that $d(z, u_i) = c_i, \forall u_i \in
U_i$.  Without loss of generality assume $C = \bigcup_i U_i$.
Defining $H_U = \{z : d(z, u_i) = d(z, u_j),\ 1 \leq i \neq j \leq m,\ u_i \in U_i\}$, we observe that $H_U$ is a finite union of hyperplanes.

Given that the lattice $\mathbb{L}$ is infinite in $d$ independent directions and that $C$ is finite, any $z \in \mathbb{L} \backslash C$ with $||z||$ sufficiently large will have a distance from any element of $C$ that has not yet been assigned to a colour class.
We must also disqualify any $z$ which has the same adjacency to, i.e.
is equidistant from any two of $u_i \in U_i, u_j \in U_j, 1 \leq
i \neq j \leq m$; in other words those lying on a hyperplane.  This is equivalent to showing that we
cannot cover a lattice $\mathbb{L}$ in $\mathbb{R}^d$ by finitely many
affine hyperplanes and a finite set (of vertices chosen up to that point).  Let $H_1, H_2, \ldots, H_r$ be a
finite collection of hyperplanes and $S$ a finite set of lattice points of size $|S| =
s$.  Choose a large prime $p$ and look at $\mathbb{L} / p\mathbb{L}$.
Now $| \mathbb{L} / p\mathbb{L} | = p^n$ for some integer $1 \leq n
\leq d$ (compare $\mathbb{L}$ to $\mathbb{Z}$ and $ \mathbb{L} / p\mathbb{L} $ to $\mathbb{Z}^{n}_{p}$).  But the image under the
reduction map of each hyperplane has size $p^{n-1}$ as the hyperplane
has codimension $1$ in the space, so the image of
the hyperplanes and $S$ contain $\leq r \cdot p^{n-1} + s$ points.  As
long as $p^n > r \cdot p^{n-1} + s$ we are done.  This is true for
$p$ sufficiently large and so the result follows by Euclid's theorem that there exists an infinity of primes.
\end{proof}

It is immediate that $\Aut(\mathbb{L}) \leq \Aut(\mathfrak{R}_{m,\omega})$, since automorphisms of $\mathbb{L}$ are isometries.

\bigskip

\begin{corollary}
The free abelian groups
\index{group ! free abelian}%
 $\mathbb{Z}^d$ act on random graphs $\mathfrak{R}_{m,\omega}$ for integers $d, m \geq 2$, so as to preserve colours.
\index{graph ! random ! $m$-coloured}%
\end{corollary}

\begin{proof}
$\mathbb{Z}^d$ is an integer lattice in $\mathbb{R}^d$.
\end{proof}

A $d$--dimensional \emph{crystallographic (or Bieberbach) group}
\index{group ! crystallographic (Bieberbach)}%
 is a cocompact discrete group of isometries of $d$--dimensional Euclidean space.
\index{Euclidean space}%

\begin{corollary}
The Bieberbach groups are subgroups of\\ $\Aut^{*}(\mathfrak{R}_{m,\omega})$.
\end{corollary}

\begin{proof}
This is immediate from the structure of Bieberbach groups as free abelian by finite.
\end{proof}

For the next theorem, we require the complexified version of the Leech
lattice $\mathbb{L}_L$~\cite[p.293]{conwayslo}.
\index{lattice ! Leech}%

For every vector $x \in
\mathbb{L}_L$, defining a vector in $\mathbb{L}_L^{\mathbb{C}}$ by $x
(a + b \omega) = ax + b(x \omega)$, where $\omega$ is a primitive cube
root of $1$, gives a $12$-dimensional lattice
(or module) over the Eisenstein integers
\index{Eisenstein integers}%
 $\mathbb{Z}[\omega]$ (see Appendix~\ref{NumberTheory}) whose automorphism group
 \index{group ! automorphism}%
  is the central extension $\Aut(\mathbb{L}_L^{\mathbb{C}}) = 6 \cdot \Suz$, where
$\Suz$~\label{Suz} denotes the Suzuki group.
\index{group ! Suzuki}%
Vectors
$x_i, y_i, z_i$ $(i \in PL(11))$ have a definition that includes
triplet formation $x_i + y_i + z_i = 0$ with Conway group
\index{group ! Conway}%
 element $\omega \in \Co_{0}$~\label{Conway} having
the transitive action $\omega: x_i \to y_i \to z_i \to x_i$.  Vectors
in $\mathbb{L}_L^{\mathbb{C}}$ can be written in terms of the ternary
Golay code.
\index{Golay code ! ternary}%
Also $\mathbb{L}_L^{\mathbb{C}}$ is the union of a
lattice and two translates.  The underlying real lattice $\mathbb{L}_L^{\mathbb{R}}$ of $\mathbb{L}_L^{\mathbb{C}}$ is scaled to have minimal norm $6$.  More on the sporadic simple groups of Suzuki and Conway can be found in~\cite{conwaya}.

\begin{theorem}
\label{cplltm}
There is a $1$--$1$ correspondence between vectors in
$\mathbb{L}_L^{\mathbb{C}}$ and vertices in $\mathfrak{R^{t}}$, and
$\mathfrak{R^{t}}$ is a Cayley object
\index{Cayley object}%
for $\mathbb{L}_L^{\mathbb{C}}$
such that multiplication by $\omega$ permutes the $3$ colours cyclically.
\end{theorem}

\begin{proof}
The first part of the statement follows from a similar type of
argument to that of Theorem~\ref{mdgptm}, given that $\mathbb{L}_L$ and
$\mathbb{L}_L^{\mathbb{C}}$ are different descriptions of the same
object.  The topological space
\index{topological space}%
 in which the Baire category argument takes place
is the space of inverse pairs of vectors in $\mathbb{L}_L^{\mathbb{C}}$, with
a partitioning of the vector pairs into $3$
colour classes.  The argument must satisfy the action of
$\omega$ which imitates the $C_{3}$ action on the colours that arises in Theorem~\ref{mdgptm}, so it remains to verify this.  Notice that the Eisenstein integers are simply an index, so the
$\omega$-action is effectively $(a + b \omega) x\ (=x_i) \to (a \omega + b \omega^2)x\ (=y_i) \to (a
\omega^2 + b) x\ (=z_i) \to (a + b \omega) x\ (=x_i)$
that is $(a, b) \to (-b, a-b) \to (b-a, -a) \to (a, b)$.  Assign
vectors in $\mathbb{L}_L^{\mathbb{C}}$ randomly to the three `colour'
classes $x_i, y_i, z_i$.  To show that this gives $\mathfrak{R^{t}}$,
we must show that we can still satisfy the $1$-point extension
property once we have thrown away the errant vectors that are
equidistant from two differently coloured vectors.  That is we want
vectors $z_i = (p + q \omega) z$ such that no two of $d \big( (p,q) -
(a, b) \big)$, $d \big( (p,q) - (-b, a-b) \big)$ and $d \big( (p,q) -
(b-a, -a) \big)$ are equal, for Euclidean metric distance function
$d$. 
\index{Euclidean metric}%
 Taken pairwise, these conditions give us at worst three
quadratic equations (in $p$ and $q$) defining three quadratic curves
lying in three planes, to be avoided.  By the same type of argument
employed in Theorem~\ref{rglttm} we can
evade any finite number of such planes and still retain a countable
infinity of such vertices $z_i$ to build $\mathfrak{R^{t}}$.  So
$\mathfrak{R^{t}}$ is a Cayley object
\index{Cayley object}%
 for $\mathbb{L}_L^{\mathbb{C}}$.  If a vector $v$ is given a certain
colour then $-v$ acquires the same colour.  Also since $\omega$ is an
isometry, $v, \omega v$ and $\omega^{2} v$ are isometric so once a
colour is chosen for $v$, $\omega v$ and $\omega^{2} v$ get assigned
colours unambiguously according to the $3$-cycle of the colours.  Hence result.
\end{proof} 

\head{Remarks}

1.  Note that $x$ and $-x$ have the same colour, as do $y$ and $-y$,
   and $z$ and $-z$, so we actually have two mutually dependent $3$-cycles.

2.  The fixed-point-free action
\index{group ! action ! fixed-point-free}%
 of $\omega$ ensures that there is never ambiguity in assigning vectors a colour.

3.  This last theorem, together with Theorem~\ref{mdgptm} demonstrate that a Cayley
    object is not associated with a unique group or structure.  We already have a much better example of this: $\mathfrak{R}$ is a Cayley graph for any countable group satisfying the square-root condition (see Appendix~\ref{PermutationGroups} and~\cite{cameron}), and in particular, any abelian group $A$ 
\index{group ! abelian}%
 in which the subgroup $\{a : 2a = 0\}$ has infinite index.

4.  From Euclidean geometry
\index{Euclidean geometry}%
 it is well-known that the union of a
    finite number of hyperplanes in a real vector space does not
    exhaust the space.  This is used, for example, in proving that a
    root system for a Lie algebra has a base~\cite[p.~48]{hum1}.
\index{Lie algebra}%

5.  Theorem~\ref{cplltm} uses little other than the existence of a fixed-point free automorphism of order 3.  The result for a more general lattice is given below in Theorem~\ref{centthree}.

\bigskip

That $\Aut(\mathbb{L}_L) \leq  \Aut(\mathfrak{R^{t}})$ follows from Theorem~\ref{rglttm}.  Furthermore $\Aut(\mathbb{L}_L^{\mathbb{C}}) \leq
\Aut(\mathfrak{R^{t}}) \sd C_3$, the group on the right-hand side
being the reduct of $\mathfrak{R^{t}}$ that permutes the $3$ colours cyclically.
In fact we have that for the Conway groups,
\index{group ! Conway}%
\[ \Co_{0} <  \Co_{\infty} < \Aut(\mathfrak{R^{t}}), \]
where $\Co_{\infty} := \mathbb{L}_{L} \sd \Aut(\mathbb{L}_{L})$ is the
affine automorphism group
\index{group ! automorphism}%
 of $\mathbb{L}_{L}$ arrived at when
$\mathbb{L}_{L}$ acts on itself by translations.

\section{Difficulties in Directly Proving Cartan Triality of $\mathfrak{R^{t}}$}

In this section we find more links between graphs, groups and lattices.  We shall conclude that any hope of finding a direct link using the Cayley object method, between the two occurrences of triality, the normalizing action of $\Sym(\mathfrak{r}, \mathfrak{b}, \mathfrak{g} )$ on $\Aut(\mathfrak{R^{t}})$ and that of Cartan is likely to fail, because of the absence of a fixed-point-free group action
\index{group ! action ! fixed-point-free}%
 which we would require to achieve sufficient randomness to build $\mathfrak{R^{t}}$.

In the preface we conjectured that there exists a connection between the colour triality of $\mathfrak{R^{t}}$ as given by the outer automorphic action of $T(\mathfrak{R^{t}}) \cong \Sym(\mathfrak{r} , \mathfrak{b} , \mathfrak{g})$ on $\Aut(\mathfrak{R^{t}})$ 
\index{Cartan triality}%
and a generalization of the algebraic triality outer automorphisms
\index{automorphism ! outer}%
 of the $8$-dimensional projective reduced orthogonal group.
\index{group ! orthogonal}%
Demonstrating such a connection would involve showing that the $\Sym({\mathfrak{r}, \mathfrak{b},\mathfrak{g}})$ group action
\index{group ! action ! $\Sym({\mathfrak{r}, \mathfrak{b},\mathfrak{g}})$}%
normalizing $\Aut(\mathfrak{R^{t}})$ can be traced back to, or in fact be induced by a generalization of the Cartan triality
\index{Cartan triality}%
\index{Cartan, \'E.}%
$\Sym(3)$ action on the vector and two spinor representations
\index{spinor (spin representation)}%
 of $\Spin(8)$.
\index{group ! spin}%
\index{group ! Spin(8)@$\Spin(8)$}%
In this section we will indicate why there are difficulties in proving a direct connection.  

\begin{lemma}
\label{reglem}
If $\mathbb{L}$ is a $d$-dimensional lattice for $d \ge 2$, the proportion of regular orbits to
the total number of orbits of $\Aut({\mathbb{L}})$ acting on $\mathbb{L}$ tends to $1$, as the
lengths of vectors approach infinity.
\end{lemma}

\begin{proof} (Proof from first principles).
The number of lattice vertices at distance at most
$n$ from the origin is $O(n^d)$.  For any $g \in
\Aut(\mathbb{L})$, at most $O(n^{d-1})$ points are fixed by $g$,
since $xg=x$ is a set of linear equations with non-zero rank. The
number of automorphisms is finite, so as $n \to \infty$, the
proportion of regular orbits to the total number of orbits approaches $1$.
\end{proof}

In other words almost all vectors lie in regular orbits of $\Aut(\mathbb{L})$.

\begin{proof} (Second proof)
The theta-function~\label{Thetafn}
\index{theta-function}%
\[\Theta(\mathbb{L}) = \sum_{v\in \mathbb{L}}x^{|v|^2},\]
is the generating function for lattice vectors by length. Now, for any
automorphism $g \in G = \Aut(\mathbb{L})$ of $\mathbb{\mathbb{L}}$, the fixed points of $g$ form a
sublattice $\mathbb{L}^g$, with theta-function $\Theta(\mathbb{L}^g)$. By summing these
and dividing by $|G|$, we get the
generating function for the number of orbits; call it $\Lambda(\mathbb{L},G)$. Similarly one can compute $\Lambda(\mathbb{L},H)$ for any
subgroup $H$ of $G$. Now a M\"obius inversion
\index{mobius@M\"obius inversion}%
 over the lattice of subgroups of $G$ gives a generating function for regular orbits. 

If $\Theta(\mathbb{L}) = \sum_{r=0}^{\infty} a_r x^{r}$ where $a_r =
|\{ v \in \Gamma | v . v = 2r\}|$ then $a_r$ counts the lattice vertices
lying on a sphere of radius $\sqrt{2r}$ around the origin.  So $a_r$ grows like
the area of the sphere, that is $(\sqrt{2r})^{d-1}$.  But $\lim_{r \rightarrow
\infty} \sqrt[r]{(\sqrt{2r})^{d-1}} =1$~\cite{ebeling} and so $\lim_{r \to \infty} \frac{\sqrt[r]{(\sqrt{2r})^{d-1}}}{\sqrt[r]{(\sqrt{2r})^{d-2}}} =1$.

\end{proof}

This result for vectors in lattices mirrors the one for regular orbits of permutation groups on the power set~\cite{cam2b}.  

We next work towards identifying different trialities associated with $\mathfrak{R^{t}}$, making
use of the $\mathbb{ATLAS}$~\cite{conwaya} in the initial discussion.
\index{atlas@$\mathbb{ATLAS}$}%
We denote the cube root of unity by $\omega$ and consider the elements $g$ of order 3 in the Conway group $\Co_{0} = 2 \ns \Co_{1}$, a non-split extension.
\index{group ! Conway}%
\index{group ! split extension}%

Let $g$ be an element of order 3, and let $f= \fix(g)$ be the number of its fixed points.  Then the $24$--dimensional matrix representation of $g$ has the following breakdown of eigenvalues:  $\sharp$ eigenvalues of value $1 = f$,  $\sharp$ eigenvalues of value $\omega = 12 - f/2$,  $\sharp$ eigenvalues of value $\omega^2 = 12 - f/2$.  So the character $\chi(g)$ of $g$ is $\chi(g) = f - (12 - f/2) = 3f/2 -12$.

\begin{figure}[ht]
\begin{tabular}{|l||c|c|c|c|}
\hline
&3A&3B&3C&3D\\
\hline
\hline
$\chi(g)$&-12&6&-3&0\\
\hline
$f$&0&12&6&8\\
\hline
\end{tabular}
   \caption{Permutation characters and fixed points in conjugacy
   classes of $\Co_{1}$ 3-cycles}
\label{pchartab}
\end{figure}

The character in Table~\ref{pchartab} is $\chi_{102}$~\cite[pp.183, 186]{conwaya}.  By inspection of the table, the 3A fixed-point-free
\index{group ! action ! fixed-point-free}%
element is the central element $\omega \in 6 \cdot \Suz$ (see the prelude to Theorem~\ref{cplltm}).
\index{group ! Suzuki}%
 We are justified in keeping the same notation for the corresponding element
of $\Co_{1}$ because from~\cite[pp.131]{conwaya} we know that the unit
scalars which can be regarded as a fixed-point-free subgroup $2A_3$ of
$2\Co_{1}$ can be extended to fixed-point-free groups $2A_4$ and $2A_5$
which yield the quaternionic and icosian Leech lattices.
\index{lattice ! Leech ! icosian}%
\index{lattice ! Leech ! quaternionic}%
The quaternions $\pm 1, \pm i, \pm j, \pm k, \frac{-1+i+j+k}{2}$
multiplicatively generate $2A_4$~\cite[pp.97]{conwaya}.
J. Tits~\cite{tits}
\index{Tits, J.}%
 has given an account of how centralizers in $\Co_{0}$ of double coverings of alternating groups form an important series of groups.  

The 3D element $\tau$ fixes each $\mathbb{L}_{E_8}$ and permutes the three $\mathbb{L}_{E_8}$ cyclically, which given that $\chi(g)=0$ therefore has the matrix form
$ \begin{pmatrix} 
0&*&0\\
0&0&*\\
*&0&0
\end{pmatrix}$, with each entry being 8--dimensional.  If we identify
the icosian ring
\index{icosian ring}%
 with $\mathbb{L}_{E_8}$, then each vector in the Leech lattice over the icosians can be written as an element of $\oplus_i \mathbb{L}^{(i)}_{E_8}$.  The double cover of the Hall-Janko group $2 \cdot J_2$
\index{group ! Hall-Janko}%
  can be realized as the automorphism group
\index{group ! automorphism}%
 of the Leech lattice over the icosians, and contains an element that permutes these three $\mathbb{L}_{E_8}$~\cite{wilsona}.  Griess
\index{Griess Jnr., R. L.}%
 has shown~\cite[p.83]{grie1} that $\Aut(\mathbb{L}_{L})$ acts transitively on the set of sublattices isometric to $\sqrt{2}\ \mathbb{L}_{E_8}$.

The $3B$ element also has an interpretation.  The $3B$ element fixes half of the 24 dimensions, so in an icosian description of $\mathbb{L}_{L}$ where the 3 imaginary units of the quaternion algebra are $e_1, e_2, e_3$, this fixes the real unit $1$ and the combination $e_1 + e_2 + e_3$.  So $3B$ acts as the 3-cycle of the 3 constituent $\mathbb{L}_{E_8}$s inside the icosian Leech lattice
\index{lattice ! Leech ! icosian}%
(see the paper by Wilson~\cite{wilsona}).
\index{Wilson, R. A.}%
  According to~\cite{kocaozdes} this 3-cycle of the simple roots   $e_1, e_2, e_3$, which is an element of the binary tetrahedral automorphism group
\index{group ! automorphism}%
 of the $\SO(8)$
\index{group ! SO(8)@$\SO(8)$}%
 root lattice, is also the outer automorphic 3-cycle of this lattice.
\index{automorphism ! outer}%

Notwithstanding the general result Theorem~\ref{rglttm} connecting coloured random graphs and lattices, we would like to explore a little further what happens in the specific case of the exceptional lattices, $\mathbb{L}_{L}$ and $\mathbb{L}_{E_8}$.   First we state a theorem that connects these two lattices.

\begin{theorem}[Lepowski and Meurman]
\index{Lepowski, J.}%
\index{Meurman, A.}%
\label{lmlemma}
 If we renormalize the $2$-dimensional lattice so that $\mathbb{L}_{L} / \sqrt{2}$ is even unimodular without
vectors of length $2$, then~\cite{lepmeu}
\[ 2 ( \mathbb{L}_{E_8} \oplus \mathbb{L}_{E_8} \oplus
\mathbb{L}_{E_8} ) \subset \mathbb{L}_{L} \subset ( \mathbb{L}_{E_8})^{3}. \]
\end{theorem}

In the next proposition, the element $\omega$ acts as an automorphism by scalar multiplication by $e^{\frac{2 \pi i}{3}}$ on the complex Leech lattice $\mathbb{L}_L^{\mathbb{C}}$,
\index{lattice ! Leech ! complex}%
 permutes the three $\mathbb{L}_{E_8}$ lattices and lies in the centre of the group $\Co_{1}$.   Acting with $\omega$ turns the 24-dimensional vector space $E_8^3 / (2 E_8)^3$ over $\mathbb{F}_2$ into a vector space over $\mathbb{F}_4 = \mathbb{F}_2(\omega)$.  The vector space $\mathbb{L}_L / (2 E_8)^3$ is then $12$-dimensional.  The existence of a fixed-point free
 \index{group ! action ! fixed-point-free}%
  automorphism of order 3 was also used in Theorem~\ref{cplltm}.    The other $3$-cycle to consider is $\tau$, an automorphism of $\mathfrak{R^{t}}$.   The fixed points of $\tau$ lie in 3-cycles of $\omega$, so the colouring of $\omega$-cycles is fixed by $\tau$, and we need to show that the total action has the form $C_3 \times C_3$.  The most natural way to show that $\tau$ commutes with $\omega$ is to demonstrate that $\omega$ preserves $\mathbb{L}_L^{\mathbb{C}}$ (or just that it preserves the complex structure on $\mathbb{R}^{24} \geq 2 \mathbb{L}_{E_8} \oplus 2 \mathbb{L}_{E_8} \oplus 2 \mathbb{L}_{E_8}$).  This latter property follows since as we noted above the $\omega$ operation takes $x_i \to y_i \to z_i \to x_i$ for each $i \in PL(11)$.  The $\tau$-action is $\tau : (E_8)_1 \to (E_8)_2 \to (E_8)_3 \to (E_8)_1 : a \mapsto a \tau \mapsto a \tau^2 \mapsto a$, where a triple $(a_1, a_2, a_3) \to (a_3 \tau, a_1 \tau, a_2 \tau)$.  If a vector is fixed by $\tau$ then so is its image under $\omega$, so the vectors $v, v \omega, v \omega^2$ are coloured the same.

We will see that we can choose a colouring from $\mathfrak{r}, \mathfrak{b}, \mathfrak{g}$ of the 3-cycles of both $\tau$ and $\omega$.  
  
\begin{proposition}
\label{24cor}
\item{(i)} The colouring of the vectors has a $C_3 \times C_3$ symmetry, this being the 3-cycle of $\mathbb{L}_{L}$--vectors in $\Co_{1}$ generated by the order 3 element $\omega$ and the 3-cycle $\tau$ on the edge colours of $\mathfrak{R^{t}}$ that is a coset of $\Aut(\mathfrak{R^{t}})$.
\item{(ii)} $\mathfrak{R^{t}}$ is a Cayley graph
\index{graph ! Cayley}%
for $2 ( \mathbb{L}_{E_8} \oplus \mathbb{L}_{E_8} \oplus \mathbb{L}_{E_8} )$ as a sublattice of $\mathbb{L}_{L}$.
\end{proposition}

\begin{proof}
(i)  The topological space
\index{topological space}%
 in which a Baire category
\index{Baire category theorem}%
 argument can be fashioned is the space of assignments of 3-cycle orbits into 3 colour classes.  We can ignore the factor of 2 outside the direct sum in the previous lemma because our colouring is of orbits of lattice vectors rather than their lengths; the origin of this factor is given in the statement of Theorem~\ref{lmlemma}.  

A typical vector is of the form $(a_1, a_2, a_3)$, where $a_j \in \mathbb{L}^{(j)}_{E_8}$.  Vectors are coloured as $(a, a \tau, a \tau^2)$-triples, and fixed points are of the form $a \oplus a \tau \oplus a \tau^2$.  Assign vectors of the form $(a, 0, 0)$ (respectively $(0, b, 0)$, $(0, 0, c)$) the colour $\mathfrak{r}$ (respectively $\mathfrak{b}, \mathfrak{g}$).  Any vector $(a, b, c)$ not of this form and not fixed by $\tau$ can be randomly assigned a colour whilst its images under $\tau$ can be coloured according to a 3-cycle of colours.  This applies to the fixed points of $\tau$ which lie in 3-cycles of $\omega$.

The finite set $S$ of points that are coloured is a union of orbits of $\langle G(2, 2), \tau, -1 \rangle$, 
\index{group ! G(2, 2)@$G(2, 2)$}%
where the last element arises because we want to colour a vector and its negative the same.  Here the group $G(2, 2)$ is the automorphism group
\index{group ! automorphism}%
 of the octonions
\index{octonions}%
 or, equivalently, the subgroup of $\SO(7)$
 \index{group ! SO(7)@$\SO(7)$}%
  that stabilizes a singular vector in its 8-dimensional real spinor representation.
\index{spinor (spin representation)}%
  The fixed points of $\tau$ in the commutator subgroup $\Omega^+(8, 2)$ of $O^+(8, 2)$ form the group $G(2, 2)$.  
  
(ii)  The vectors in $\oplus_i \mathbb{L}^{(i)}_{E_8}$ $(i = 1, 2, 3)$ are randomly coloured either according to 3-cycles of $\tau$ on the three $\mathbb{L}_{E_8}$ lattices (if the vectors do not belong to $\fix(\tau)$), or from the discussion in the preamble to this proposition, according to 3-cycles of $\omega$ (if they do belong to $\fix(\tau)$).  The two sets of vectors involved are disjoint and as we noted the two types of colouring are done independently.  The result goes through as a special case of the analysis of the next theorem.
\end{proof}

The argument for a more general lattice is the following,

\begin{theorem}
\label{centthree}
Let $\mathbb{L}$ be a $d$-dimensional lattice for $d \ge 2$, and let $\langle \tau \rangle \le \Aut(\mathbb{L})$, where $\tau$ is of order $3$ fixing only the origin $0$.  Assume that $\tau$ fixes no vector (see the special case in the previous proposition).  Then each $\langle \tau \rangle$-orbit is made up of three vectors permuted by $\tau$.  If we choose one point in each orbit, with the further requirement that if $v$ is chosen then so is $-v$ (the colour of a lattice vector and its negative must be the same), colour chosen points randomly, so that $\col(\{u, v\}) = \col(u - v)$ where the left-hand colour is that of an edge and the right-hand colour is that of a lattice vector, and extend to a  colouring of the other two vectors cyclically by $\mathfrak{r} \to \mathfrak{b} \to \mathfrak{g} \to \mathfrak{r}$, then with probability 1 we obtain a coloured Cayley graph for the additive group of the lattice that is isomorphic to $\mathfrak{R^{t}}$.  Furthermore, this can be extended to the action of both $\tau$ and $\omega$ (which is as above).
\end{theorem}

\begin{proof}

Suppose that a finite set $S = U \cup V \cup W$ of vectors is given.  To satisfy the one-point extension 
\index{one-point extension property}%
 property $(*_t)$ required to identify $\mathfrak{R^{t}}$ we require a vertex $z$ such that

(i)  $z \notin S$,

(ii)  $z - x \notin S$ for $x \in S$,

(iii)  $(z - x_1) g \neq z - x_2$ for $x_1, x_2 \in S$, $g \in \langle -I, \tau \rangle$.

Both (i) and (ii) exclude finitely many vectors.  If $(z - x_1) g = z - x_2$ we require that $g \neq 1$ as $x_1 \neq x_2$.  Then $z g - z =  x_1 g - x_2$, where the elements on the right-hand side are only finite in number.  So as $z (g-1) = a$, say, $z$ lies in a coset of $\ker(g-1)$, which given that $g \neq 1$ is a proper subspace.  So we have to exclude the finite union of finitely many proper subspaces.  So there are still infinitely many $z$ left to satisfy $(*_t)$.  For any further point $z$ the colours of $z - x$ are independent for $x \in X$.  

If $s=|S|=|U| \cup |V| \cup |W|$, then the probability that $z$ is not correctly joined is $1 - \frac{1}{3^{s}}$.  Choosing $z_1, \ldots, z_i$ the probability $p(z_1, \ldots, z_i)$ that these vertices are not correctly joined is 
\begin{displaymath}
 \lim_{i \to \infty} \left(1 - \frac{1}{3^{s}}\right)^{i} = 0.
\end{displaymath}

Consider triples of orbits.  Colour the points in one orbit randomly and colour all their $\tau$-images and $\tau^2$-images with the other two colours.  Let $g_l \in \Aut(\mathbb{L})$ and let the corresponding image under subgroup inclusion $\Aut(\mathbb{L}) \le \Aut(\mathfrak{R^{t}})$ be $g_r \in \Aut(\mathfrak{R^{t}})$.  Then it follows that $g_r \in C_{\Aut(\mathfrak{R^{t}})} (\tau)$.  It is clear that $\tau \in \Aut^{*}(\mathfrak{R^{t}})$ is inducing the 3-cycle on the colours.

We can extend this to the situation with $g \in \langle -I, \tau, \omega \rangle$ by a similar argument.
\end{proof}

By Theorem \ref{rglttm} almost all the lattice vector pair colourings can be put
 into $1$--$1$ correspondence with $\mathfrak{R^{t}}$ edge colourings.
 According to~\cite{kocaozdes} the 3-cycle on the three
 $\mathbb{L}_{E_8}$ lattices is precisely the Cartan triality
\index{Cartan triality}%
3-cycle which we denote $\tau_C$~\label{tauC}.  Result (ii) of Proposition~\ref{24cor} means that the Cayley graph construction of Theorem~\ref{centthree} does not suffice for us to claim that $\Aut(\mathfrak{R^{t}})$ derives at least its 3-cycle outer automorphism
\index{automorphism ! outer}%
 from the exceptional 8--dimensional Cartan triality of $E_8$ representations.  The presence of the fixed points of the $\tau$ orbit restricts the otherwise free choice of colourings and thus precludes full randomness; this forbids us from making a $1$--$1$ correspondence of the 3-cycles.  The fixed points of $\tau$ are not an artifact of the particular colouring scheme that we employed but rather of the overall method in conjunction with the use of 8-dimensional spaces, for they would still be present had we assigned the orbits of pairs of lattice vectors of the Cayley right regular action of $\mathbb{L}_{E_8}$ on itself, as we show next.  
 
 If we were colouring regular orbits rather than 3-cycles then \emph{every} regular orbit has a subset of size $|G(2, 2)| = 12096$~\cite[p.125]{conwaysm}~\cite[p.14]{conwaya}
\index{group ! G(2)@$G(2)$}%
 of coloured vectors representing lattice points fixed by $\tau$, in other words group elements commuting with the cycle action.  Let $G_1 = \Aut(\mathbb{L}_{E_8})/\{\pm 1\} \cong
SO^{+}(8, \mathbb{Z} / 2 \mathbb{Z})$.
\index{group ! SO(8)@$\SO(8)$}%
  Notice that we can identify a vector $v \in \mathbb{L}_{E_8}$ with
its translate $t_v \in \mathbb{L}_{E_8}^{+}$, where $\mathbb{L}_{E_8}^{+} = \{(x \mapsto x
+ v): x, v \in \mathbb{L}_{E_8} \}$ is the additive group of the
lattice; $\mathbb{L}_{E_8}^{+}$
is a permutation of $\mathbb{L}_{E_8}$.  Identify the elements in the regular orbits of $G_1$
with elements of $G_1$, so that an infinite number of
copies of the finite group $G_1$ are being
considered.  Because the group $\Aut(\mathbb{L}_{E_8})$ acts as an additive
general linear group of the lattice, its elements act
like lattice isometries.  But by definition of a linear map,
stabilizing lengths requires stabilization of the origin, and this is
conjugate to stabilizing another point of the lattice, call it $z$.

Let $S$ be a finite set of points in the lattice $\mathbb{L}_{E_8}$.
We want to find a point $z$ such that 

(i)  all points of $S$ lie in different regular orbits of the
pointwise stabilizer $G_z$ of $z$ which is conjugate to
$\Aut(\mathbb{L}_{E_8})$, and

(ii)  the length of $z$ is arbitrarily large. 

Property (i) is required to ensure that every lattice point can be coloured unambiguously.

The point $c \in S$ is in a regular orbit of $G_z$ if and only if the
point $z$ is in a regular orbit of $G_c$.  For each $c \in S$ and for
each $g \in G_c, g \ne 1$, the set of fixed points of $g$ is contained
in an affine subspace.  We must choose $z$ from outside all such
subspaces and so outside of their union.  This is always possible
because the affine subspace in question is a lattice translation of a
finite proper sublattice, so that the union of the translated points
is a translation of a finite set.  Therefore the colouring of $\mathbb{L}_{E_8}$ can be extended to a colouring of all the pairs of vectors in $\oplus_i \mathbb{L}^{(i)}_{E_8}$ in a way that is invariant under translations.

Property (ii) is easily achieved by choosing the length of $z$ to exceed the sum of the lengths of all vectors chosen up to that stage.

In conclusion, $\mathfrak{R^{t}}$ can be constructed by colouring regular orbits of $\mathbb{L}_{E_8}$ acting on itself.  However, any approach to demonstrating links between $\mathfrak{R^{t}}$ and Cartan triality
\index{Cartan triality}%
 based on $\mathbb{O}$ that requires a fixed-point-free action,
\index{group ! action ! fixed-point-free}%
  is likely to be thwarted because $G(2) = \Aut(\mathbb{O})$ $(= \Spin(7) \cap \Spin(7) \subset \Spin(8)$)
\index{group ! Spin(8)@$\Spin(8)$}%
 and thus there will always be stabilized points; the Lie algebra of $G_2$
\index{group ! G(2)@$G(2)$}%
 is the triality-invariant part of the Lie algebra of $\SO(8)$.
\index{group ! SO(8)@$\SO(8)$}%
  Another way to see this is to note that if $e_0 = 1, e_1, \ldots, e_7$ are the units of the octonion algebra then $e_i e_j = - \delta_{ij} + \sigma_{ijk} e^k$ and $G(2)$ stabilizes the coefficients in the completely antisymmetric $\sigma_{ijk}$.  (We should note that the complex Lie group $G(2)$ has two real forms, only one of which is compact,  and it is the compact form that can be lifted to the universal covering $\Spin(7)$ of $\SO(7)$; the compact group $G(2)$ decomposes the 8-dimensional irreducible real \emph{spin representation} into the trivial and the 7-dimensional representation.

Any condition such as the existence of fixed points is a constraint on the randomness of the edge-colouring.  By the above results and by Theorem~\ref{nagyvothm} of Appendix 9, the existence of fixed points creates an ambiguity in colouring and thus spoils our attempt to link the colour triality of $\mathfrak{R^{t}}$ with Cartan triality,
\index{Cartan triality}%
 and this is likely to be equally true of different attempts to create links via lattices and 3-nets.

\vspace{25pt}

For more graphs represented by roots of $\mathbb{L}_{E_8}$, but in a different context, see~\cite[p.~103]{brco}.

\vspace{25pt}

In the next two sections we will take a somewhat different approach to proving the Cartan Triality of $\mathfrak{R^{t}}$.

\section{Groups with Triality}

A \emph{group with triality}
\index{group ! with triality}%
 is a group $G$ with automorphisms $\rho$ and $\sigma$ satisfying
\[\rho^3=\sigma^2=(\rho\sigma)^2=1\]
and the \emph{triality identity}
\index{triality identity}%
 \[[g,\sigma][g,\sigma]^\rho[g,\sigma]^{\rho^2}=1\]
for all $g \in G$, where $[g,\sigma]=g^{-1}g^\sigma$.  

Obviously, $\langle\rho,\sigma\rangle$ is a homomorphic image of $\Sym(3)$.  More precisely, fix a homomorphism $h: \Sym(3) \to \Aut(G)$ and identify $\rho, \sigma$ with their images under $h$; then $G$ is a \emph{group with operators $\Sym(3)$}.
\index{group ! with operators}%
The triality identity is independent of a particular choice of $\rho$ and $\sigma$~\cite{doro}.

The study of abstract groups with triality was begun by Doro~\cite{doro}
\index{Doro, S.}%
following Glauberman~\cite{glauberman}.
\index{Glauberman, G.}%
 A discussion this concept is given in~\cite{hallnagy}; see also Appendix~\ref{LoopTheory}.
\index{triality identity}%

\emph{The only finite simple groups with triality are $D_4(q)$~\cite{liebeck}.
\index{group ! D(4)@$D(4)$}%
  So there is a sense in which the concept of a group with triality is a generalization of Cartan triality.
\index{Cartan triality}%
  In fact the origin of the name ``group with triality''   as well as the motivating example is Cartan's triality group (Appendix~\ref{carsec})}.

\head{Examples}

\begin{enumerate}

\item

Any group is a group with triality in a trivial way, taking $\rho$ and $\sigma$ to be the identity automorphism.

\item

If $G$ is a group with triality, and $H$ a subgroup invariant under $\rho$ and $\sigma$, then $H$ is also a group with triality.

\item

Let $H$ be an abelian group of exponent~$3$,
\index{group ! abelian of exponent~$3$}%
 and $G$ the semidirect product of $H$ with an involution $t$ inverting $H$. Take $\rho$ to be the inner automorphism
 \index{automorphism ! inner}%
 induced by an element $h\in H$, and $\sigma$ the inner automorphism induced by $t$. The relations satisfied by $\rho$ and $\sigma$ are clear. We have $[g,\sigma] \in H$, and so this element has order~$3$ and is fixed by $\rho$. So $G$ is a group with triality. By the first example, $H$ is also a group with triality. (Note that $\rho$ acts trivially on $H$.)

\item

In particular, $\Sym(3)$, acting on itself by inner automorphisms, is a group with triality.

\item

Let $A$ be any group, and define permutations $\rho$ and $\sigma$ of $G=A \times A$ by
\[(x,y)^{\rho}=(y,y^{-1}x^{-1}),\qquad (x,y)^{\sigma}=(y,x).\]
It is easily checked that $\sigma$ is an automorphism of order $2$. Also $\rho$ has order~$3$, since 
\[(x,y)^{\rho^2}=(y^{-1}x^{-1},xyy^{-1})=(y^{-1}x^{-1},x),\] 
and so 
\[(x,y)^{\rho^3}=(x,x^{-1}xy)=(x,y).\] 

The triality identity
\index{triality identity}%
 is satisfied, since we have
\[(x,y)(y,y^{-1}x^{-1})(y^{-1}x^{-1},x)=(1,1).\] 
If $A$ is abelian, then $\rho$ is also an automorphism and $(\rho\sigma)^2=1$. So, if $A$ is abelian, then $G=A \times A$ is a group with triality.

\item

In particular, the Klein $4$-group
\index{group ! Klein}%
 (whose automorphism group is $\Sym(3)$) is a group with triality.
\index{group ! automorphism}%

\item

$\Sym(3)$ is a subgroup of $\Sym(4)$, and so acts on it by conjugation. But $\Sym(4)$ is not a group with triality. For we may take $\rho=(1,2,3)$ and $\sigma=(1,2)$. If $g=(2,4)$, then $[g,\sigma]=(2,1,4)$, and so $[g,\sigma][g,\sigma]^\rho[g,\sigma]^{\rho^2} = (2,1,4)(3,2,4)(1,3,4) = (2,3)(1,4)$.
This example also shows that, if $G$ is a group with an action of $\Sym(3)$, then the set of elements of $G$ satisfying the triality identity is not necessarily a subgroup of $G$. For both the Klein group and $\Sym(3)$ satisfy the identity, and their product is $\Sym(4)$.

\item

If $S$ is a nonabelian simple group then the wreath product $S \Wr \Sym(3) = (S_1 \times S_2 \times S_3) \sd \Sym(3)$ is a group with triality~\cite{hallji}, where the base group $S_1 \times S_2 \times S_3$ has automorphism group
\index{group ! automorphism}%
  $\Aut(S) \Wr \Sym(3) = (\Aut(S_1) \times \Aut(S_2) \times \Aut(S_3)) \sd \Sym(3)$. See also~\cite{doro} and~\cite{nagyval}.

\end{enumerate}

A good reference to the triality of the $8$-dimensional orthogonal group is~\cite{fultha}.
\index{group ! orthogonal}%

Let $G$ be a group with automorphisms $\rho$ and $\sigma$ satisfying $\rho^3=\sigma^2=(\rho\sigma)^2=1$.  For $g \in G$ , $\sigma, \rho \in \Sym(3)$, let $\Phi(g) = [g, \sigma] [g, \sigma]^{\rho} [g, \sigma]^{\rho^2}$~\label{trid} be the triality identity.  So $G$ is a group with triality if and only if $\Phi(G) = 1$
\index{triality identity}%
 for all group elements $g \in G$.  If $\Phi(g) = 1\ \forall g \in G$\
and $H, K < G$ are both groups with triality with the same $\sigma$ and $\rho$, then $\Phi(hk) = 1\ \forall h \in H, k \in K$.

\bigskip

The following simple propositions give a taste of the group theory.

\begin{proposition}
A subgroup $K \le G$ is $\langle \rho, \sigma \rangle$--invariant if for $k \in K$, $[\sigma, k] = 1$.
\end{proposition}
\begin{proof}
An easy calculation shows $\Phi(k) = (\sigma^{k} (\sigma \rho^{2}))^3 = ([k, \sigma]
\rho^2)^3$ from which we have that 
\begin{align*} 
\Phi(k^{\sigma}) &= ([k^{\sigma}, \sigma] \rho^2)^3\\
     &= ((\sigma^{-1} k \sigma)^{-1} \sigma^{-1} (\sigma^{-1} k \sigma) \sigma \rho^{2})^3\\
     &= (\sigma^{-1} k^{-1} \sigma k \rho^{2})^3\\
     &= ([\sigma, k] \rho^2)^3.
\end{align*}
So if $[\sigma, k] = 1$ then $\Phi(k^{\sigma}) = 1$.

\smallskip

Also,
\begin{align*} 
\Phi(k^{\rho}) &= ([k^{\rho}, \sigma] \rho^2)^3\\
     &= (\rho^{-1} k^{-1} \rho \sigma^{-1} \rho^{-1} k \rho \sigma \rho^{2})^3\\
     &= (\rho^{-1} k^{-1} \sigma \rho k \sigma \rho)^3
\end{align*}
because $\rho \sigma = \sigma \rho^{2} \Leftrightarrow \rho \sigma \rho^{-1} = \sigma \rho$.  If $[\sigma, k] = 1$ then
\begin{align*} 
\Phi(k^{\rho}) &= (\rho^{-1} k^{-1} \sigma \rho \sigma k \rho)^3\\
     &= (\rho^{-1} k^{-1} \rho^{-1} k \rho)^3\\
     &= 1.
\end{align*}
\end{proof}

\begin{proposition}
Let $G$ be a group with triality.  If $H, K \le G$, $[\sigma, h]=1$, and $K$ and $G$ are triality groups with the same $\rho, \sigma$ then $\Phi(hk) = 1\ \forall h \in H, k \in K$.
\end{proposition}
\begin{proof}
 $\Phi(hk) = ([(hk),\sigma]\rho^{-1})^3$.  If $[\sigma, h] = 1$ then
$\Phi(hk) = \Phi(k) = 1$.
\end{proof} 

In particular if $h, k \in C_{G}(\sigma)$ then $\langle h k \rangle
\in C_{G}(\sigma)$.  This potentially gives us a way of constructing a subgroup of
$\Aut(\mathfrak{R^{t}})$ which is generated by those elements in
$\Aut(\mathfrak{R^{t}})$ that centralize the transpositions making up
$\sigma$, and which is also a group with triality.  One way to proceed
then is to define $S : = \{g \in \Aut(\mathfrak{R^{t}}): \Phi(g) = 1\}$ and to search for a large subgroup contained in $S$ and invariant under $\sigma$ and $\rho$.  In the next section we turn to a specific realization of this.


\section{Split Extension of $\Aut(\mathfrak{R^{t}})$}
\index{group ! split extension}%

Our aim is firstly to construct $\mathfrak{R^{t}}$ in such a way as to show that there is a group isomorphic to $\Sym(\mathfrak{r} , \mathfrak{b} , \mathfrak{g})$ acting as a group of automorphisms of the graph, and secondly to prove the semidirect product $\Aut^{*}(\mathfrak{R^{t}}) := \Aut(\mathfrak{R^{t}}) \sd T(\mathfrak{R^{t}})$ where $T(\mathfrak{R^{t}}) \cong \Sym(\mathfrak{r} , \mathfrak{b} , \mathfrak{g})$.  


\begin{theorem}
\label{outgpwttr}
There is a construction of $\mathfrak{R^{t}}$ having a group isomorphic to $\Sym(\mathfrak{r}, \mathfrak{b}, \mathfrak{g})$ acting as an automorphism group.
\end{theorem}
\begin{proof} 
The vertex set of $\mathfrak{R^{t}}$ is countably infinite so we can label its elements by $\mathbb{N}$.

First we prove the existence of an edge-colour transposition $\sigma$ and its normalizing action on $\Aut(\mathfrak{R^{t}})$.  We do not strictly need this step and can go straight to demonstrating the semiregular
\index{group ! permutation ! semiregular}%
 action of a group isomorphic to $\Sym(3)$ on 6-sets of vertices, but we do so in order to illustrate the method which would not work if there were only 2 colours.  

For concreteness, let $\sigma = (0\ 1)\ (2\ 3)\
\ldots$ represent vertex permutations.  Colour edges so that we get $\mathfrak{R^{t}}$, take $\sigma \in \Aut(\mathfrak{R^{t}})$ to be the vertex permutation that induces a transposition $(\mathfrak{b} \mathfrak{g})$ of edge colours.  Let the edge $\{2i, 2i+1\}$ have colour $\mathfrak{r}$.  Choose the colour of edge $\{2i, 2j\}$ randomly from respectively
$\mathfrak{r}, \mathfrak{b}$ or $\mathfrak{g}$ and then the edge $\{2i+1,
2j+1\}$ will be coloured respectively $\mathfrak{r}, \mathfrak{g}$ or
$\mathfrak{b}$.  The edge colours on $\{2i, 2j+1\}$ and $\{2i+1,
2j\}$ will be similarly transposed with each other.  We need to show that the resulting graph is isomorphic to
$\mathfrak{R^{t}}$.  Let $U, V, W$ be finite disjoint vertex sets and
choose a vertex $z$, say $z=2p$, not lying in any $2$-cycle containing
points of $U, V$ or $W$ ($2$-cycle as we are considering $\sigma$).  Attach $z$ to $2i$ and $2i+1$ with edges
whose colour is arbitrarily chosen.  Colour the edges
joining $z$ to the vertices in $U \cup V \cup W$ at random.  
$$
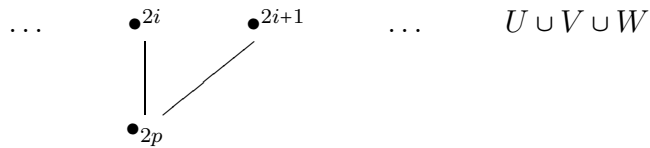
\xymatrix{ 
{\ldots} & \ar@{-}[d] {\bullet}^{2i} & {\bullet}^{2i+1} \ar@{-}[dl] & {\ldots} &
{U \cup V \cup W}\\
{} & {\bullet}_{2p}
}
$$
The usual construction of $\mathfrak{R^{t}}$ goes through and by
homogeneity and universality of $\mathfrak{R^{t}}$, $\sigma$ acts to normalize $\Aut(\mathfrak{R^{t}})$.

A 
 treatment for $\rho \in \Sym(\mathfrak{r}, \mathfrak{b},
\mathfrak{g}) \cong \langle \sigma, \rho\ |\ \sigma^2 = \rho^3 =
(\sigma \rho)^3 =1 \rangle$ to the one given for $\sigma$ will show
that $\rho$ normalizes $\Aut(\mathfrak{R^{t}})$.  

Consider the 6--vertex set $\{0, \ldots, 5\}$ and parametrise the $\Sym(3)$-action by:
 
\quad\quad\quad\quad\quad\quad\quad\begin{tabular}{cccccc}
1 & $\rho$ & $\rho^2$  & $\sigma$ & $\sigma \rho$ & $\sigma \rho^2$\\
0 & 1 & 2 & 3 & 4 & 5.
\end{tabular}

We colour edges in a $\Sym(3)$-invariant way, in particular with $\Sym(3)$ acting on $2$-element subsets of itself.  Take $\sigma =
(\mathfrak{r})(\mathfrak{b}\mathfrak{g})$ to fix the colour $\mathfrak{r}$, and take $\rho =
(\mathfrak{r}\mathfrak{b}\mathfrak{g})$.  We give edge $\{0, 1\}$ the fiducial colour $\mathfrak{r}$.  If edge $\{0, 1\}$ is $\mathfrak{r}$ then $\{1, 2\}$ becomes $\mathfrak{b}$ and $\{2, 0\}$ becomes $\mathfrak{g}$.  Acting by $\sigma$ and $\rho$ gives the other twelve edge colours.  For example $\{0, 3\}\sigma = \{3, 0\}$ so the colour of edge $\{0, 3\}$ must be fixed by $\sigma$, and therefore must be red.  Once edge $\{0, 1\}$ is specified as $\mathfrak{r}$ the remainder of the edges must be coloured as given in the table in Figure~\ref{dettab}:
\begin{figure}[ht]
\begin{tabular}{|l||c|c|c|c|c|c|}
\hline
&0&1&2&3&4&5\\
\hline
\hline
0&&$\mathfrak{r}$&$\mathfrak{g}$&$\mathfrak{r}$&$\mathfrak{g}$&$\mathfrak{b}$\\
\hline
1&$\mathfrak{r}$&&$\mathfrak{b}$&$\mathfrak{g}$&$\mathfrak{b}$&$\mathfrak{r}$\\
\hline
2&$\mathfrak{g}$&$\mathfrak{b}$&&$\mathfrak{b}$&$\mathfrak{r}$&$\mathfrak{g}$\\
\hline
3&$\mathfrak{r}$&$\mathfrak{g}$&$\mathfrak{b}$&&$\mathfrak{b}$&$\mathfrak{r}$\\
\hline
4&$\mathfrak{g}$&$\mathfrak{b}$&$\mathfrak{r}$&$\mathfrak{b}$&&$\mathfrak{g}$\\
\hline
5&$\mathfrak{b}$&$\mathfrak{r}$&$\mathfrak{g}$&$\mathfrak{r}$&$\mathfrak{g}$&\\
\hline
\end{tabular}
   \caption{Deterministic $\Sym(3)$-invariant colouring of 6-vertex clique}
\label{dettab}
\end{figure}

To reiterate, if $\{i, j\}$ denotes an edge with colour $c$ and $g \in \Sym(3)$ then we require that edge $\{i^g, j^g\}$ in a semiregular
\index{group ! permutation ! semiregular}%
 action, has colour $c^g$ in its natural action on $\{\mathfrak{r}, \mathfrak{b}, \mathfrak{g}\}$.  Consider different orbits $O_k = \{6k, \ldots, 6k + 5\}$ for $\Sym(3)$; the edges within each such set are \emph{deterministically} coloured so as to respect the $\Sym(3)$-action.  The randomness comes in because the first vertex of one $6$-vertex orbit is \emph{randomly} joined to one vertex of another $6$-vertex orbit.  This scheme is illustrated in Figure~\ref{dettab}, where edge $\{0, 1\}$ has a free choice of colours and all other colours are forced, for example the colour of edge $\{0, 3\}$ is fixed by $\sigma$ and is therefore $\mathfrak{r}$.

(There may be other colouring schemes, for example for $l > k$ we might colour the edges $\{6k + j, 6l\}$ for $j = 0, \ldots, 5$ \emph{randomly}, whence the colours of the other edges from $O_k$ to $O_l$ are determined by the proper action of $\Sym(3)$.  In particular if the edge $\{6k + j, 6l\}$ $(j \epsilon \{0, \ldots, 5\})$ is coloured randomly, then the edge $\{6k + j, 6l + i\}$ is the image of $\{6k + j', 6l\}$ $(j' \epsilon \{0, \ldots, 5\})$ under an element of $\Sym(3)$, this being the element corresponding to $i$.)

It remains to check that  $(*_t)$ holds with probability $1$, that is with probability $1$ some point $6n$ with $n$ large enough has correct joins to the given finite subset.  The usual limiting argument on the product of the orbits, each orbit $O_k$ being considered independently, gives that the probability of not getting $\mathfrak{R^{t}}$ is
\begin{displaymath}
 \lim_{k \to \infty}  \left(1 - \frac{1}{3^{s}}\right)^{k}  = 0,
\end{displaymath}
where $s=|S|=|U| \cup |V| \cup |W|$.
\end{proof}

\bigskip

A direct proof of the next result using the previous theorem would require the involved step-by-step verification of the triality identity at every stage of the graph extension process.  We can avoid this by noting that it follows as a corollary of the loop construction in the next section.

\begin{theorem}
\label{outgpwttr1}
The group $\Aut(\mathfrak{R^{t}})$ has a subgroup $G$ such that $(G,
T(\mathfrak{R^{t}}))$ is a group with triality, with $T(\mathfrak{R^{t}})) \cong \Sym(\mathfrak{r}, \mathfrak{b}, \mathfrak{g})$.
\end{theorem}

\bigskip

With the topology of pointwise convergence
\index{topology ! of pointwise convergence}%
 on $\Aut(\mathfrak{R})$, a subgroup is dense
\index{group ! dense}%
 in this group if it has the same orbits as $\Aut(\mathfrak{R})$ in the induced action on the set of vertex $n$-tuples, for all $n \geq 1$.  In~\cite{bhatmac}, Bhattacharjee and Macpherson,
\index{Bhattacharjee, M.}%
\index{Macpherson, H. D.}%
give a construction of a dense locally finite subgroup of the group $\Aut(\mathfrak{R})$.  They build finite vertex induced subgraphs $\Delta_i$ of $\Aut(\mathfrak{R}) =  \bigcup_{i < \omega} \Delta_i$ and build finite groups $\hat{G_i} \leq \Aut(\Delta_i)$ such that $\hat{G_i}  =  G_i^{\Delta_i}$ (the group induced by $G_i$ on $\Delta_i$) for each $i < \omega$.  Their method is based on arguments of Herwig and Lascar
\index{Herwig, B.}%
\index{Lascar, D.}%
 in~\cite{herwiglas}, these being adaptations of Hrushovski's~\cite{hrush}),
\index{Hrushovski, E.}%
 building automorphisms as unions of finite partial automorphisms.
\index{automorphism ! partial}%
Is it possible to interweave this theory together with Theorem~\ref{outgpwttr} to a sufficiently random colouring of any finite graph extensions to be executed with repeated verifications of the triality identity to strengthen the previous theorem? 

\head{Open Question}   Does the group $\Aut(\mathfrak{R^{t}})$ have a dense locally finite subgroup $H$ such that $(H, T(\mathfrak{R^{t}}))$ is a group with triality, with $T(\mathfrak{R^{t}})) \cong \Sym(\mathfrak{r}, \mathfrak{b}, \mathfrak{g})$?

The result can be further strengthened; by Theorem~\ref{toplist}, the pointwise stabilizer of a finite vertex set in the group $H$ if it exists, is both open and dense 
\index{dense open set}%
in the group $\Aut(\mathfrak{R^{t}})$.  This strengthening may be interesting for future applications where density is insufficient by itself; recall that both $\mathbb{Q}$ and $\sqrt{2} \mathbb{Q}$ are dense in $\mathbb{R}$, but $\mathbb{Q} \cap \sqrt{2} \mathbb{Q} = \emptyset$.

\medskip

We end the section with a proposition and an application of it.

\begin{proposition}
If $G$ is a group for which $\zeta(G) = 1$, $S \le \Aut(G)$, $S \cap \Inn(G) = 1$, and $S$ acts via its embedding in $\Aut(G)$, then $C_{G \sd S}(G) = 1$.
\end{proposition}
\begin{proof}
Let $c = gs \in C_{G \sd S}(G)$.  Then the automorphism of $G$ induced by $c$ is trivial, so the outer automorphism
\index{automorphism ! outer}%
 induced by $s$ is trivial, where by assumption $s = 1$.  Then $c = g \in \zeta(G)$, so $c = 1$.
\end{proof}

Now $G = \Aut(\mathfrak{R^{t}})$ is simple and so centreless, and if $\Aut^{*}(\mathfrak{R^{t}})) := \Aut(\mathfrak{R^{t}}) \sd \Sym(3)$ denotes the split extension
\index{group ! split extension}%
 of $\Aut(\mathfrak{R^{t}})$ by $\Sym(3)$ then the conditions of the above proposition are satisfied so that we have $C_{\Aut^{*}(\mathfrak{R^{t}})}(\Aut(\mathfrak{R^{t}})) = 1$.  This means that $\Aut^{*}(\mathfrak{R^{t}})$ acts faithfully on $\Aut(\mathfrak{R^{t}})$, so $\Aut^{*}(\mathfrak{R^{t}}) \le \Aut(\Aut(\mathfrak{R^{t}}))$.  To show that there is equality here we trace the steps of the previous proposition, as follows.

Let $\alpha$ be any automorphism of $G = \Aut(\mathfrak{R^{t}})$, and $v$ a vertex.  Then $(G_v)^{\alpha}$ is a subgroup of countable index, so by the small index property 
\index{small index property (SIP)}%
contains the pointwise stabilizer of a finite set; (automorphisms of $G$ map one subgroup of small index in the group to another).  But $G$ has just four $G_v - G_v$ double cosets (corresponding to equality and the three edge colours), and one can show that the rank of any other small index subgroup is more than four.  For example on red edges, two red edges could be equal, or they could meet at a vertex giving a rank of three corresponding to the three colours on the third edge, or they could be disjoint so giving a further rank of $3^4$.  So $G_v^{\alpha}$ is a vertex stabilizer.  Thus $\alpha$ is induced by an element of $\Sym(V(\mathfrak{R^{t}}))$ normalizing $G$.  Such an element must permute the three edge colours, so lies in $\Aut^{*}(\mathfrak{R^{t}})$.  So we have proved that

\begin{theorem}
\label{autsymaut}
$\Aut^{*}(\mathfrak{R^{t}}) = \Aut(\mathfrak{R^{t}}) \sd \Sym(\mathfrak{r}, \mathfrak{b}, \mathfrak{g})$.
\end{theorem}

Therefore $\Aut^{*}(\mathfrak{R^{t}})$ is the full automorphism group
\index{group ! automorphism}%
 of $\Aut(\mathfrak{R^{t}})$.

In the next chapter we return to the subject of extensions of automorphism groups
\index{group ! automorphism}%
 of random graphs of any number of colours.

\section{A Loop Construction of $\mathfrak{R^{t}}$}
\label{bolsec}

We have seen that generalizations of the concept of Cayley digraphs
\index{digraph}%
 exist for a wider range of objects than just groups.  Here is another example.

Let $Q$ be a quasigroup
\index{quasigroup}%
 and let $E \subseteq Q$.  The \emph{quasigroup digraph} $\Gamma = QG(Q, E)$ is defined by vertices $V(\Gamma) = Q$ and edges $A(\Gamma) = \{(u, ue) | u \in Q, e \in E \}$.  If $E \subseteq Q$ is right associative that is satisfies $a(bE) = (ab)E$ for all $a, b \in E$, then $\Gamma$ is said to be a \emph{quasi-Cayley digraph}.
\index{digraph ! quasi-Cayley}%

The following result was proved by D\"{o}rfler~\cite{dorfler}
\index{D\"{o}rfler, W.}%
\begin{theorem} 
A digraph $\Gamma$ is regular if and only if it is a quasigroup digraph.  Moreover, if $\Gamma = QG(Q, E)$ with $E$ right-associative, then $\Gamma$ is vertex-transitive.
\end{theorem}

Sabidussi
\index{Sabisussi, G.}%
 showed~\cite{sabisussi1} that every vertex-transitive graph is a retract of a Cayley graph; a retract is defined later in the book.  Gauyacq
\index{Gauyacq, G.}%
 found examples of vertex-transitive graphs that are not quasi-Cayley, including an infinite family of Kneser graphs.
\index{graph ! Kneser}%
She proved a theorem that gave a condition on the vertex set of a connected graph to be quasi-Cayley~\cite{gauyacq}\cite{laurisc}.

In~\cite{kelarev}, Kelarev
\index{Kelarev, A.}%
 and Praeger
\index{Praeger, C. E.}%
 study analogues of Cayley graphs for semigroups
\index{semigroup}%
  and show that they sometimes enjoy properties analogous to those of the Cayley graphs of groups.

\medskip

We now undergo our own study pertaining to infinite random graphs.

A Cayley graph for a group $G$ can be defined either as

(a)  a graph on the vertex set $G$ admitting the action of $G$ by right multiplication as a group of automorphisms;

or

(b)  a graph of the form $\Cay(G, S)$ with vertex set $G$ and edges $\{g, sg\}$ for $g \in G, s \in S$, where $S$ is an inverse-closed subset of $G \backslash \{1\}$.

Clearly (a) is the ``natural'' definition but there are problems applying it if $G$ is a loop.
\index{loop}%
  Indeed, if $G$ is nonassociative then (b) does not imply (a), since right multiplication by $h$ would take the edge $\{g, sg\}$ to $\{gh, (sg)h\}$, which is not necessarily the same as $\{gh, s(gh)\}$.  So for loops we have to abandon (a) and use (b) as a definition.

\medskip

Before continuing with Cayley graphs of loops, we specify one more concept.

A \emph{normal Cayley graph}
\index{graph ! Cayley ! normal}%
 of a group is one that admits both left and right multiplications by its elements as automorphisms.  The condition for $\Cay(X, S)$ to be a normal Cayley graph for a group $X$ is that $x y^{-1} \in S \Rightarrow (ax) (ay)^{-1} \in S = a (x y^{-1}) a^{-1} \in S$, that is $S$ is closed under conjugation.  Therefore an equivalent definition~\cite{cam11} is that a Cayley graph $\Cay(X, S)$ on a group $X$ is \emph{normal} if $S \subseteq X$ is a normal subset (one fixed by conjugation by all elements of $X$); for example Cayley graphs for abelian groups are normal.  So if a group $G$ has $\mathfrak{R^{t}}$ as a normal Cayley graph then $\{G\} \backslash \{1\} = S_{\mathfrak{r}} \bigsqcup S_{\mathfrak{b}} \bigsqcup S_{\mathfrak{g}}$, where each of $S_{\mathfrak{r}}, S_{\mathfrak{b}}, S_{\mathfrak{g}}$ are both inverse-closed and \emph{conjugation-closed} sets of elements.  More on conjugation-closure is given in Lemma~\ref{moufcond}(b).  Note that because a graph can be isomorphic to the Cayley graph of more than one group, and the condition of normality depends on the group rather than being a graph-theoretic property, then normality is effectively a joint property of the graph and the group.  

In a Cayley graph for a general loop
\index{loop}%
 the above two definitions are not equivalent.  Because of the nonassociativity of loop multiplication, to claim that a Cayley graph for Moufang loops is invariant under right or left multiplication means that extra conditions must be imposed.

We need to find conditions on countably infinite loops that ensure that the Cayley graphs representing these loops are isomorphic to $m$-coloured random graphs.  The condition in question turns out to be the same as the one for groups.  Recall from~\cite{cameron} the following result for countable groups.

\begin{theorem}
\label{sqrootthm}
The following three conditions on the countable group $G$ are equivalent:

(a)  some Cayley graph of $G$ is isomorphic to the random graph;
\index{graph ! Cayley}%
\index{graph ! random}%

(b)  with probability $1$, a random Cayley graph of $G$ is isomorphic to the random graph;

(c)  for any two finite disjoint subsets $U$ and $V$ of $G$, there exists $z \in G$ such that, for all $u \in U, v \in V$, we have $(z u^{-1})^2 \ne v u^{-1}$.
\end{theorem}

This theorem shows that if $G$ is any regular permutation group, then either no $G$-invariant graph is isomorphic to $\mathfrak{R}$ or almost all are.  It can be shown using a combination of the square root condition times $C_{\infty}$ that every countable group is a subgroup of $\Aut(\mathfrak{R})$.

The theorem readily generalizes to more than two colours by considering the colours pairwise.  In generalizing the two-colour case we say that a group is \emph{$\mathfrak{R}_{m}$-genic}
\index{graph ! Rmgenic@$\mathfrak{R}_{m}$-genic}%
if it satisfies this theorem with the following condition:\quad for all $m$-tuples $(U_1, \ldots, U_m)$ of finite pairwise disjoint subsets, $\exists z$ such that if $u, v$ belong to distinct sets $U_i$, $U_j$ then $(z u^{-1})^2 \ne v u^{-1}$.  This condition appears to get stronger as $m$ increases, but does it really?  That is, is there a group which satisfies it for $m$ but not for $m-1$?  More concretely, consider the following three conditions:

$(*)_2$\quad for finite disjoint $U, V$ $\exists z$ such that $\forall u, v$ $(z u^{-1})^2 \ne v u^{-1}$.

$(*)_m$\quad for finite disjoint $U_1, \ldots, U_m$ $\exists z$ such that $\forall u \in U_i, v \in U_j (i \ne j)$ $(z u^{-1})^2 \ne v u^{-1}$.

$(*)_{\dagger}$\quad for a finite set $U$ $\exists z$ such that $\forall u, v \in U, u\ne v$ $(z u^{-1})^2 \ne v u^{-1}$.

Certainly the implications $(*)_{\dagger} \Rightarrow (*)_m \Rightarrow (*)_2$ hold, but whether or not they reverse is an open question.

In terms of Baire category, assuming that $G$ satisfies the $\mathfrak{R}_{m}$-genic condition, the set of Cayley graphs for $G$ which are isomorphic to $\mathfrak{R}_{m, \omega}$ is residual.  

We use the same nomenclature for loops
\index{loop}%
 as for groups, and refer to Appendix~\ref{LoopTheory} for the definition of a Moufang loop.
\index{Moufang ! loop}%
So a Cayley graph for the Moufang loop $Q$
takes the form $\Cay(Q, S)$, where the vertices are $q \in Q$ and
the inverse-closed subset $S \subseteq Q \backslash \{1\}$ gives the
edge set $\{\{q, qs^{-1}\}\ : q \in Q, s \in S\}$.

\begin{theorem}
\label{lpcond}
The following three conditions on the countable Moufang loop $Q$ are equivalent:

(a)  some Cayley graph of $Q$ is isomorphic to the random graph;
\index{graph ! random}%

(b)  with probability $1$, a random Cayley graph of $Q$ is isomorphic to the random graph;

(c)  for any two finite disjoint subsets $U$ and $V$ of $Q$, there exists $z \in Q \backslash \{U \cup V\}$ such that, for all $u \in U, v \in V$, we have  $(z u^{-1})^2 \ne v u^{-1}$.
\end{theorem}

\begin{proof}
Clearly (b) $\Rightarrow$ (a).

To show that (a) $\Rightarrow$ (c), let $Q$ be a Moufang loop.  Firstly, the identity $a(xy^{-1}) = a (xa) \cdot a^{-1} y^{-1}$ for elements of $Q$ follows because $a(xy^{-1}) = a [x(a (a^{-1} y^{-1}))] = [((ax)a) \cdot (a^{-1} y^{-1})] = a (xa) \cdot a^{-1} y^{-1}$.  

Let $S$ be an inverse-closed subset of $Q$ which does not contain the identity.  Suppose that $\Cay(Q, S) \cong \mathfrak{R}$, and let $U$ and $V$ be finite disjoint subsets of $Q$.  Choose $z$ according to the property $(*)$.  For all $u \in U, v \in V$, if $(z u^{-1})^2 = v u^{-1}$ then $(z u^{-1}) (z u^{-1}) u = (v u^{-1}) u$, so $(z u^{-1}) z = v$ giving $z u^{-1} = v z^{-1}$ which leads to $\{z, u\}$ and $\{z, v\}$ being both edges or both non-edges depending on whether or not $zu^{-1} \in S$.  Then we deduce that $zu^{-1} \ne vz^{-1}$, so $(z u^{-1})^2 = v u^{-1}$ as required.

To show that (c) $\Rightarrow$ (b), first note that given any $U$ we can always enlarge it by adding a finite number of elements, so there are an infinite number of elements $z$ satisfying (c).  For such $z$, we always have $(z u^{-1})^2 \neq v u^{-1}$.  So the probability that $z$ is not correctly joined in a random Cayley graph is at most $1 - \frac{1}{2^m}$, where $m = |U \cup V|$.  (An example of a smaller probability would be if, say $(z u_1^{-1})^2 = u_2 u_1^{-1}$, so that $\{z, u_1\}$ and $\{z, u_2\}$ are both edges or non-edges).  Given finitely many such $z$, say $z_1, \ldots, z_n$, only finitely many more $z$ for which the joins of $z$ to $U \cup V$ are dependent on those of $z_1, \ldots, z_n$, for such dependence would require that $z x^{-1}= (z_i y^{-1})^{\pm 1}$ for some $i$ and some $x, y \in U \cup V$, and each such equation has a unique solution in $z$.  So we can construct an infinite sequence of $z$ for which these events are independent, and conclude that with probability 1 some $z$ is correctly joined.  There are only countably many choices of $U$ and $V$ so (b) follows.

\end{proof}

This theorem shows that either no Cayley graph for $Q$ is isomorphic to $\mathfrak{R}$ or almost all are,  and the proof is a simple variant to the argument given in~\cite{cameron} for the case of groups.  The theorem for Moufang loops readily generalizes to more than two colours by considering the colours pairwise.

We work with Moufang loops whose elements are assumed to have square roots~\cite[p.121]{bruck}.  For $l \in L$, define the \emph{square-root set} $\sqrt l : = \{x \in L: x^2 = l\}$.  If $l \ne 1$ it is \emph{non-principal}.  A \emph{translate} of a square-root set is obtained by multiplying it on the right by a fixed loop element.  The following sufficient condition to be found in~\cite{cameron} for groups, follows as a corollary to the above theorem applied to Moufang loops.

\begin{corollary}
If the countable Moufang loop $L$ cannot be expressed as the finite union of translates of non-principal square-root sets together with a finite set, then $L$ is $\mathfrak{R}_{m}$-genic.
\end{corollary}
\index{graph ! Rmgenic@$\mathfrak{R}_{m}$-genic}%
 
Note that left and right multiplications are bijections in a loop, for example $xa = ya \Rightarrow x = y$.

\begin{lemma} 
\label{moufcond}
\item{(a)} If $\RMlt(Q)$~\label{RMlt} is the group generated by all
  the right multiplications on $Q$, then $\RMlt(Q) \leq
  \Aut(\Gamma)$, where $\Gamma$ is any Cayley graph for $Q$.
\item{(b)} If $\LMlt(Q)$~\label{LMlt} is the group generated by all
  the left multiplications on $Q$, and if in $\Cay(Q, S)$ we have
  that $x y^{-1} \in S \Rightarrow [ (ax) y^{-1}] a^{-1} \in S$ then
  $\Cay(Q, S)$ is a normal Cayley graph admitting $\LMlt(Q)$.
\end{lemma}

\begin{proof}
(a)  For the set of elements $S \subseteq Q$, join two elements $x, y \in S$ to form an edge denoted $x \sim y$ if and only if $xy^{-1} \in S$.  In order to show that $\RMlt(a) \in \Aut(Q)$ we must prove that $(xa)
  (ya)^{-1} = xy^{-1}$ for elements of $Q$.  Using the identity $a(xy^{-1}) = a (xa) \cdot a^{-1} y^{-1}$, we find that,
\[ a^{-1} [a(xy^{-1})] = a^{-1} [a (xa) \cdot a^{-1} y^{-1}] \Rightarrow\] \[xy^{-1} = a^{-1} [(a (xa)) \cdot (a^{-1} y^{-1})] = a^{-1} [(a (xa)) \cdot (ya)^{-1} ] =\] \[a^{-1} [(ax)a \cdot (ya)^{-1}) ] = a^{-1} [a (x (a \cdot (ya)^{-1}))] = x ( a \cdot (ya)^{-1} ) =\] \[ (x a a^{-1}) [a (ya)^{-1}]  = (x a) a^{-1} (a a^{-1}) [a (ya)^{-1}]=\] \[ (x a) (a^{-1} a) a^{-1} [a (ya)^{-1}]= (xa) (ya)^{-1}.\] 

(b)  We can assume that $a \ne x \ne y^{-1}$, because if any two are
 equal then their composition associates, and we revert to the group
 case.  More precisely, by Moufang's theorem
 \index{Moufang's theorem}%
in a Moufang loop any three elements which associate generate a group, and by corollary every Moufang loop is diassociative (that is, any pair of elements whatsoever generates a group).  It also follows from Artin's diassociativity theorem
\index{Artin diassociativity theorem}%
 any two elements of a Moufang loop generate a subgroup.  Then
\begin{align*}
(ax) (ay)^{-1} &= [(ax)a] (a^{-1} a) [a^{-1} (ay)^{-1}] \\
     &= [a(xa)] (a^{-1} (ay)^{-1} ) \\
     &= ( ( [a(xa)] a^{-1} ) y^{-1} ) a^{-1} \\
     &= [ (ax) y^{-1}] a^{-1}.
\end{align*}
thus giving the condition that $x y^{-1} \in S \Rightarrow [ (ax) y^{-1}] a^{-1} \in S$ for $\LMlt(Q) \leq  \Aut(\Gamma)$.  If this condition, which we shall describe as \emph{conjugation-closure}, is satisfied together with (a) then it follows that $ \Mlt(Q) = \langle \LMlt(Q), \RMlt(Q) \rangle \leq \Aut(\Gamma)$.  
\end{proof}

The condition for a Cayley graph of a loop to be normal given in Lemma~\ref{moufcond}(b) is the loop equivalent to the one for groups given in the paragraph before the Lemma.
For a Moufang loop $Q$, let $\zeta(Q) := C(Q) \cap \Nuc(Q)$ be its centre,
\index{Moufang ! centre}%
where $C(Q)$ is the \emph{Moufang centre}
\index{Moufang ! Moufang centre}%
 and $\Nuc(Q)$~\label{Nuc} is the \emph{Moufang nucleus};
\index{Moufang ! Moufang nucleus}%
  (see Appendix~\ref{LoopTheory} for the definitions).

Let $K(Q)$ be the class of Moufang loops $Q$, such that 
\[ 1 < \zeta(Q) \le C(Q) < Q,\]
\[ exp(C(Q)) = 3,\]
\leftline{and,}
\[ Q \ncong \zeta(Q) \times L, \text{for some subloop}\ L.\]

Let $K'(Q)$ be the (large and general) class of Moufang loops that are \emph{not} members of $K(Q)$.  J. D. Phillips
\index{Phillips, J. D.}%
 proved the following theorem~\cite{phillips}~\cite{phillips1},

\begin{theorem}
\label{phillipsthm}
For a Moufang loop $Q$ belonging to $K'(Q)$, $\Mlt(Q)$ is a group with triality
\index{group ! with triality}%
 if and only if
\[ \zeta(Q) = \Nuc(Q),\]
\[ exp(C(Q)) = 3, and\]
\[ Q \cong \zeta(Q) \times L, \text{for some subloop}\ L\ \text{with}\ \Nuc(L) = 1; \text{or}\ C(Q) = Q.\]
\end{theorem} 

It is also known that there is an outer automorphism
\index{automorphism ! outer}%
 $\rho$ of order 3 of the loop's multiplication group~\label{Mlt} (see Appendix~\ref{LoopTheory}) acting as
\[ L(q)^{\rho} = R(q),\quad R(q)^{\rho} = P(q),\quad P(q)^{\rho} = L(q)\quad\quad\quad (\dagger) \]
where $q \in Q$ and $P(q) = R(q) L(q) = L(q) R(q)$, and an involutory automorphism acting as
\[  R(q)^{\sigma} = L(q^{-1}),\quad L(q)^{\sigma} = R(q^{-1}), \quad P(q)^{\sigma} = R(q)^{-1} L(q)^{-1},\quad\quad (\dagger\dagger) \]
where $L(q)^{-1} = L(q^{-1}), R(q)^{-1} = R(q^{-1}), P(q)^{-1} = P(q^{-1})$.  

Also $\langle \rho, \sigma \rangle \cong \Sym(3)$.

The action of $\rho$ in the next result, which is equivalent to the edge-colour permuting 3-cycle $\tau$ in Theorem~\ref{mdgptm}, is assumed to have no fixed points except the identity.  That we can make this assumption is shown in the dialogue after the proof of the theorem.

Doro proved~\cite{doro}
\index{Doro, S.}%
 that each Moufang loop with trivial nucleus has a multiplication group with triality.  We will use this fact in the next theorem.  We aim to produce a construction of $\mathfrak{R^{t}}$ in which the colour $\Sym(\mathfrak{r} , \mathfrak{b}, \mathfrak{g})$ as well as the generalized Cartan $\Sym(3)$ are present, in order that further research can potentially unveil the relationship between them.  

\begin{theorem}
\label{moulthm}
Let $(Q, \cdot)$ be a countably infinite Moufang loop with trivial nucleus and trivial Moufang centre, 
and satisfying the conditions of Theorem~\ref{lpcond} and Theorem~\ref{phillipsthm}.  


We can find $\mathfrak{R^{t}}$ as a Cayley graph for $(Q, \cdot)$ such that $\Sym(\mathfrak{r} , \mathfrak{b}, \mathfrak{g}) : = \langle \rho, \sigma \rangle$ acts so as to permute the colours.
\end{theorem}

\begin{proof}
The topological space
\index{topological space}%
 in which the Baire category argument takes place is the space of paths in a ternary tree representing a countable sequence of choices.  The complete metric space which is required for the application of this theory arises from paths in rooted trees of countable height.  This space is effectively isomorphic to that used to show that the set of Cayley graphs for $C_{2} \ast C_{2} \ast C_{2}$ which are isomorphic to $\mathfrak{R^{t}}$ is residual in a conditional set, so we do not need to repeat those arguments.  

Partition the non-identity inverse pairs of elements of $Q$ into $3$ inverse-closed, conjugation-closed classes $S_1, S_2, S_3$ (corresponding to edges joined to the identity by red, blue and green colours), such that
\[ q \in S_i \Rightarrow q^{-1} \in S_i\quad\quad (i = 1, 2, 3)\quad\quad\quad\quad (a) \] 


If (A) we assign pairs of loop elements to one of three colours, red if $ql^{-1} \in S_1$, blue if $ql^{-1} \in S_2$ and green if $ql^{-1} \in S_3$, and (B) there is an order 3 group element $\rho$ such that
\[ q \in S_1 \Rightarrow q^{\rho} \in S_2, q^{\rho^2} \in S_3,\quad\quad\quad\quad\quad\quad\quad (b) \]
and there is an order 2 group element $\sigma$ such that
\[ \sigma: S_1 \to S_2 \to S_1,\ S_3 \to S_3,\quad\quad\quad\quad\quad\quad (c) \]
then with the loop elements as graph vertices and pairs of elements as graph edges randomly coloured, we must show that the set of normal Cayley graphs representing the loop $(Q, \cdot)$ which are isomorphic to $\mathfrak{R^{t}}$ is residual in the set satisfying condition $(a), (b)$ and $(c)$.

We must choose $S_1, S_2, S_3$ in such a way as to ensure that the resulting object is isomorphic to $\mathfrak{R^{t}}$. The existence of $\mathfrak{R^{t}}$ is equivalent to the satisfaction of the following condition:

($*_{t}$)  If $U$, $V$ and $W$ are finite
disjoint sets of graph vertices, then there exists a vertex $z$,
joined to every vertex in $U$ with a red edge, to every vertex in $V$ with a blue edge, and joined to every vertex in $W$ with a green edge.

(i)  We claim that $\exists z$ such that no $sz^{-1} (s \in S)$ has yet
been assigned one of the three colours.  For if $\Phi$ is the
finite set whose elements have already been so assigned, then $sz^{-1} \in \Phi$ so only the finitely many $z$, those of the form $\phi^{-1} s, \phi^{-1} \in \Phi$ need to be excluded.  Nonassociativity does not arise in compositions of pairs of Moufang elements.

(ii)  We want to eliminate $z$ for which $uz^{-1} = (vz^{-1})^{-1}$ for $u \in U, v \in V$.  For such elements then $zu^{-1} = vz^{-1}$ (because Moufang loops have the inverse property and their elements have two-sided inverses).  Substituting $z = z' u$, gives 

\[ (z' u) u^{-1} = v (u^{-1} z'^{-1}) \Rightarrow z' = v (u^{-1} z'^{-1}) \Rightarrow u^{-1} z' = u^{-1} (v (u^{-1} z'^{-1})) \] 
\[ \Rightarrow (u^{-1} z') z' = (u^{-1} v) u^{-1} \Rightarrow u^{-1} (z')^{2} = u^{-1} (v u^{-1}) \Rightarrow \]
\[ \Rightarrow z' =  (vu^{-1})^{\frac{1}{2}} \Rightarrow z =  (vu^{-1})^{\frac{1}{2}} u. \]

By Theorem~\ref{lpcond}, this is precisely the form of $z$ that we must eliminate in order to be able to form any Cayley graph from elements of a loop.


Let $u \in S_1, v \in S_2, w \in S_3$.  We need to assume that the images of all $uz^{-1}, vz^{-1}, wz^{-1}$ under $\sigma$, $\rho$ and $\rho^{2}$ have not yet had colours assigned to them, so that we can assign them colours appropriately; to achieve randomness we must have a free choice of colours.  

\emph{Note: the elements $\sigma$, $\rho$ and $\rho^{2}$ act on $\Mlt(Q)$ and not on $Q$}.  

This leads us to consider the following cases.

(iii)  We need to exclude those $z$ for which $(zu^{-1})^{\rho} = zv^{-1}$.  Now $z^{\rho} u^{- \rho} = zv^{-1} \Rightarrow z^{-1} (z^{\rho} u^{- \rho}) = z^{-1} (zv^{-1}) = v^{-1} \Rightarrow z^{\rho} [z^{-1} (z^{\rho} u^{- \rho})] = z^{\rho} v^{-1} \Rightarrow [(z^{\rho} z^{-1}) z^{\rho}] u^{- \rho} = z^{\rho} v^{-1} \Rightarrow z^{2 \rho - 1} u^{- \rho} = z^{\rho} v^{-1} \Rightarrow z^{\rho - 1} u^{- \rho} = v^{-1} \Rightarrow$\[ z^{\rho - 1} = v^{-1} u^{- \rho^{-1}}\quad\quad (\alpha)\]
As a check, substituting $z^{\rho - 1}$ from $(\alpha)$ in the original equation for $(zu^{-1})^{\rho}$ gives $z^{\rho} u^{- \rho} = (z z^{\rho - 1}) u^{- \rho} = (z v^{-1} u^{- \rho^{-1}}) u^{- \rho} = zv^{-1}.$

By disqualifying the $z$ of the form $z = (v^{-1} u^{- \rho^{-1}})^{(\rho - 1)^{-1}}$, we can achieve the required inequality.

(iv)  We need to exclude those $z$ for which $(zu^{-1})^{\rho^{2}} = zw^{-1}$.  But this is just  $(z u^{-1})^{\rho^{-1}} = zw^{-1}$.  By comparison with (iii) we see that we need to disqualify $z$ of the form $z = (w^{-1} u^{- \rho})^{(\rho^{-1} - 1)^{-1}}$ in order to satisfy the required inequality.

(v)  We need to exclude those $z$ for which $(zu^{-1})^{\rho} = (zv^{-1})^{-1}$.  This is simply $(z u^{-1})^{\rho^{-1}} = zv^{-1}$, so we require $z \neq (v^{-1} u^{- \rho})^{(\rho^{-1} - 1)^{-1}}.$

(vi)  We need to exclude those $z$ for which $(zu^{-1})^{\rho^2} = (zw^{-1})^{-1}$.  This implies $(z u^{-1})^{\rho} = zw^{-1}$, so we eliminate $z = (w^{-1} u^{- \rho^{-1}})^{(\rho - 1)^{-1}}.$

(vii)  We need to exclude those $z$ for which $((zv^{-1})^{-1})^{\rho} = zw^{-1}$.  Thus $(z v^{-1})^{- \rho} = zw^{-1}$, so we exclude $z = (w^{-1} v^{\rho^{-1}})^{(- \rho - 1)^{-1}}.$

(viii)  We need to exclude those $z$ for which $(zu^{-1})^{\rho^2} = (vw^{-1})$, that is $z = (vw^{-1})^{\rho} u.$

(ix)  We need to exclude those $z$ for which $((zv^{-1})^{-1})^{\rho} = (vw^{-1})$, that is $z = (vw^{-1})^{- \rho^{-1}} v.$

(x)  Consider $(uw^{-1})^{\rho} = (u_1 w_1^{-1})$.  This relation has already been determined and is only possible if $uw^{-1} \in S_1$ (respectively $S_2, S_3$) implies $u_1 w_1^{-1} \in S_2$ (respectively $S_3, S_1$).

Next we need to show that the involution $\sigma: S_1 \to S_2 \to S_1$, $S_3 \to S_3$ does not affect the satisfaction of ($*_{t}$).

(xi)  We need to exclude those $z$ for which $(zu^{-1})^{\sigma} = zv^{-1}$.  By comparison with (iii) we see that we need to disqualify $z$ of the form $z = (v^{-1} u^{- \sigma^{-1}})^{(\sigma - 1)^{-1}}$.

(xii)  We need to exclude those $z$ for which $(zu^{-1})^{\sigma} = (zv^{-1})^{-1}$, that is for which $(zu^{-1})^{\sigma^{-1}} = zv^{-1}$.  Thus we need to disqualify $z$ of the form $z = (v^{-1} u^{- \sigma})^{(\sigma^{-1} - 1)^{-1}}$.

(xiii)  We need to exclude those $z$ for which $(zw_1^{-1})^{\sigma} = zw_2^{-1}$.  By comparison with (iii) we see that we need to disqualify $z$ of the form $z = (w_2^{-1} w_1^{- \sigma^{-1}})^{(\sigma - 1)^{-1}}$.

(xiv)  We need to exclude those $z$ for which $(zu^{-1})^{\sigma} = vw^{-1}$.  It is immediate that we need to disqualify any $z$ of the form $z = (vw^{-1})^{\sigma^{-1}} u$.

(xv)  Consider $(u_1 v_1^{-1})^{\sigma} = (u_2 v_2^{-1})$.  This relation has already been determined and is only possible if $u_1 v_1^{-1} \in S_1$ (respectively $S_2$) implies $u_2 v_2^{-1} \in S_2$ (respectively $S_1$), and conversely.

All other possibilities fit into one of these categories.

To summarize, for any three finite disjoint sets $U, V$ and $W \subset Q$, there exists $z \in Q \backslash S$ such that, for all $u \in U, v \in V, w \in W$, we can find $uz^{-1}, vz^{-1}, wz^{-1}$ which have not yet been assigned to $S_1, S_2, S_3$ and for which $(uz^{-1})^{\rho} \neq (vz^{-1})^{\pm 1}$, $(vz^{-1})^{\rho} \neq (wz^{-1})^{\pm 1}$, and $(wz^{-1})^{\rho} \neq (uz^{-1})^{\pm 1}$, so that the images of $uz^{-1}, vz^{-1}, wz^{-1}$  under $\rho, \rho^{2}$ have not been pre-assigned.  Now assign all $uz^{-1}$ to $S_1$, all $vz^{-1}$ to $S_2$ and all $wz^{-1}$ to $S_3$.

In total, we have excluded only finitely many $z$, so infinitely many choices remain for us to satisfy the required $1$-point extension.  The set $\chi(A, B, C)$ of Cayley graphs
\index{graph ! Cayley}%
is dense.  The Baire category theorem implies that the intersection of all the sets $\chi(A, B, C)$ is residual and hence non-empty.  Hence we have built a copy of $\mathfrak{R^{t}}$.
\end{proof}

Assume that the relevant criteria are satisfied and that we are
working in a normal Cayley graph and so all three types of
multiplication of $\Mlt(Q)$, left, right and bi-multiplication
\index{Moufang ! bi-multiplication}%
are automorphisms of $Q$.  Then it is immediate from the last theorem that
\[ \Mlt(Q) \le \Aut(\mathfrak{R^{t}}). \]

\bigskip

\head{Open Question}  

Can we find an action on the vertices of $\mathfrak{R^{t}}$ that induces a normalising action on $\Mlt(Q)$.


\head{Open Question}  

By~\cite[Theorem~III.2.7]{pflug}, if two loops are isotopic then all of their respective multiplication groups, as well as their inner mapping groups and nuclei are isomorphic.  Does the multiplication group determine the centre, nucleus and other loop properties?  So is it possible to conclude that if $Q$ satisfies the conditions of the theorem then so does any loop isotopic to $Q$?

\bigskip

The associative law holds for mappings so $L(q)$ and $R(q)$ are associative, even though the loop composition is not.  Because they are permutations, $L(q), R(q)$ lie in $\Sym(Q)$.  For elements lying outside the Moufang centre, $L(q), R(q)$ and $P(q)$ are mappings (and thus always single-valued) from $Q \backslash C(Q) \to \Sym(Q)$; $L, P, R$ are also $1$--$1$ and are pairwise disjoint.  So by $(\dagger)$ above, $\rho$ acts fixed-point-freely
\index{group ! action ! fixed-point-free}%
  on the Moufang loop multiplication groups, other than on the identity.  (However~\cite{phillips2} $\rho$ fixes both
 $R(c)$ for an order 3 element $c$ of the Moufang centre, and automorphic \emph{inner mappings}
\index{Moufang ! inner mapping group}%
\index{group ! inner mapping}%
 of a loop with trivial nucleus.  The inner mapping group $I(Q)$ comprises those elements of $\Mlt(Q)$ that stabilize the identity and form a subgroup of $\Mlt(Q)$ generated by $R(p, q) = R(p) R(q) R(pq)^{-1}, L(p, q) = L(p) L(q) L(qp)^{-1}$, and $T(p) = R(p) L(p)^{-1}$).  
 
Whilst inner mappings are not always loop automorphisms, the collection of all inner automorphisms
\index{automorphism ! inner}%
 of a loop $Q$ form a subgroup of $I(Q)$~\cite[pp.~27, 28]{pflug}).

Not much is known about infinite Moufang loops, but notice that for the previous theorem we chosen a loop that is totally non-commutative and nonassociative, and therefore a loop that is far away from having a group structure as possible.  What if $(Q, \cdot)$ was a group?  Then $\Nuc(Q) = Q$ as all the elements associate.  So by Theorem~\ref{phillipsthm}, $\Nuc(Q) = \zeta(Q) = C(Q) = Q$.  So the group must be abelian,
\index{group ! abelian}%
 and so $\Mlt(Q) = Q$; of necessity the group is of exponent 3.  As we noted in example (5) in the section on groups with triality,
\index{group ! with triality}%
 for an abelian group $A$, the group with triality is, using additive notation, $A \oplus A$.  Then the action of $\rho$ is $\rho : (x, 0) \mapsto (0, x) \mapsto (-x, -x) \mapsto (0, -x) + (x, x) = (x, 0)$.  Also $\rho : (x, y) \mapsto (-y, x-y)$ which is fixed if and only if $x = -y$ and $x-y = y \Rightarrow x = 2y = -y$.  Also $(x, -x) \mapsto (x, x+x) = (x, -x)$.  So $\rho$ has no fixed points outside the Moufang centre.
 
If instead of removing nonassociativity from Theorem~\ref{moulthm} we keep nonassociativity and remove non-commutativity, the theorem still fails.  A commutative Moufang loop (CML) (see Appendix~\ref{LoopTheory}) would fail to satisfy the last theorem.  If $Q$ were a group, then $\Mlt(Q)$ admits $\rho$ only if $Q$ is abelian of exponent 3~\cite[p.~1485]{phillips}.  However if $Q$ were a CML, necessarily of exponent 3~\cite[Lemma~1]{phillips}, then the triality of $\Mlt(Q)$ is trivial~\cite[Theorem~3]{phillips} (because then $L(q) = R(q) = P(q)).$  So we take $Q \in K'(Q)$ to be such that $\zeta(Q) = \Nuc(Q) = C(Q) = 1$, and we can further assume if necessary that there are an infinite number of elements of unbounded exponent in $Q$.  That Moufang loops with trivial nucleus have multiplication groups with triality was shown in~\cite{doro}.
\index{group ! with triality}%

It is immediate that the multiplication group of the Moufang loop of the last theorem is a group with triality, that is non-trivial on all elements.

\bigskip

Let us recap where we are and where we believe further work could produce elegant results.

\head{Open Question}  

\bigskip

\textit{We have already constructed $\mathfrak{R^{t}}$ so that $\Mlt(Q) \le \Aut(\mathfrak{R^{t}})$.  We would like to find an action on the vertices of $\mathfrak{R^{t}}$ that induces a normalizing action on $\Mlt(Q)$.}

Let us recap where we are and what remains to be achieved.

(i)  We know that colour $\Sym(\mathfrak{r} , \mathfrak{b}, \mathfrak{g})$ acts as outer automorphisms of $\Aut(\mathfrak{R^{t}})$.  

(ii)  We have a construction that shows that the extension of $\Aut(\mathfrak{R^{t}})$ by colour $\Sym(3)$ splits,
\index{group ! split extension}%
 that is we can realize this $\Sym(3)$ as a subgroup of $\Sym(\mathfrak{R^{t}})$ which normalizes $\Aut(\mathfrak{R^{t}})$.

(iii)  We can choose a Moufang loop $Q$ such that it admits triality and has $\mathfrak{R^{t}}$ as a Cayley graph.  Then the Cartan-like triality $\Sym(3)$ acts as outer automorphisms of $\Mlt(Q)$.  

(iv)  Can we find a $\Sym(3)$ in $\Sym(Q) = \Sym(\mathfrak{R^{t}})$ normalizing $\Mlt(Q)$ and inducing the triality $\Sym(3)$ on it?  If so then this extension splits.  (It is possible to build extensions, some non-split
\index{group ! non-split extension}%
 and this unique one that splits.)
\index{group ! split extension}%
Can we demonstrate that the split extension that normalizes $\Mlt(Q)$ simultaneously acts on the vertices of $\mathfrak{R^{t}}$?  

(v)  Can we by careful choices in (i) and (ii) find a nice relationship between the two $\Sym(3)$s ($\Sym(\mathfrak{r} , \mathfrak{b}, \mathfrak{g})$ and the triality $\Sym(3)$).  For example can we make them equal? commute?

\bigskip

\smallskip

\begin{figure}[!ht]
$$\xymatrix{
{\Sym(Q)} \ar[d]\\
{\mbox{Cartan-like\ Triality}} \ar[d]^{acts\ on}  && {\Sym(\mathfrak{r} , \mathfrak{b}
  , \mathfrak{g})} \ar[d]^{acts\ on}\\
{\Mlt(Q)} \ar[d]^{acts\ on} \ar[rr]^{embedding} && {\Aut(\mathfrak{R^{t}})}
\ar[dll]^{acts\ on}
\\
{Q}
}$$
\caption{Generalized Cartan triality \& $\mathfrak{R^{t}}$ triality}  
\label{trifig}
\index{Cartan triality}%
\end{figure}
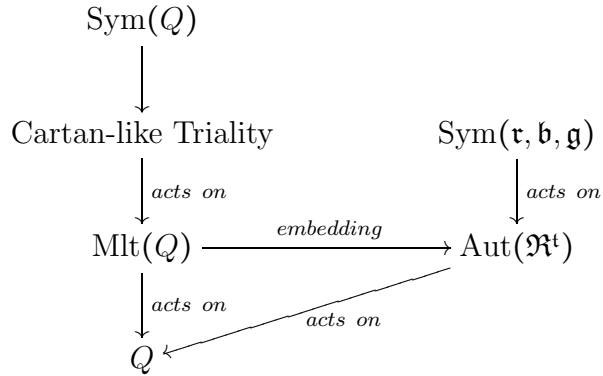


Note that in the right-hand side of Figure~\ref{trifig}, $\Sym(\mathfrak{r} , \mathfrak{b}, \mathfrak{g})$ can act directly on $Q$ because $\Sym(Q)$ contains both $\Sym(\mathfrak{r} , \mathfrak{b}, \mathfrak{g})$ and $\Aut(\mathfrak{R^{t}})$, and by Theorem~\ref{autsymaut} the first of these groups normalizes the second.  There is however no equivalent action on the left-hand side of the diagram.  Whilst the generalized Cartan-like group acts on $\Mlt(Q)$ it cannot act directly on $Q$, for if we could get $\langle \tilde{\rho}, \tilde{\sigma} \rangle$ to act on $Q$ then $L(q^{\tilde{\rho}}) = R(q)$, and for example applying this to the identity, $q^{\tilde{\rho}} = q$ for all $q \in Q$, so $\tilde{\rho} = 1$; or applying this to any $x \in Q$, $(x q)^{\tilde{\rho}} = q x$ so that in the particular case of $x = q$, $q q^{\tilde{\rho}} = q q$ and we reach the same contradiction, since $\tilde{\rho}$ has no fixed points.  Thus
the two left-hand actions are not transitive and so we cannot have either a direct triality action on the loop nor is there an equivalent of the split extension
\index{group ! split extension}%
 of the groups on the right.


 Although we are asking for a relationship between the colour and Cartan-type trialities, we should note the unusual and unexpected nature of this connection were it to be found.  A priori, how do we begin to connect a graph that is countably infinite and whose automorphism group is uncountably infinite, with a vector space that has dimension 8, and whose quadratic-form-preserving automorphism groups (the orthogonal groups $\Spin(8)$, and $P\Omega^{+}(8, K)$, over \emph{any} field $K$ including a finite field, making this group a finite group)?  On another level, picking out the triality graph from the class of all three-edge coloured countable universal homogeneous graphs has probability one, and so this graph is the norm in its class, whereas the 8-dimensional phenomenon is unique to that dimension.  It is possible that such a purported connection can be found through the generalized notion of groups with triality.  The tantalizing prospect is that in some sense this makes $\mathfrak{R^{t}}$ an exceptional object, inheriting its exceptionality from Cartan triality
\index{Cartan triality}%
\index{Cartan, \'E.}%
 and perhaps originating in the octonions $\mathbb{O}$~\cite{baez} \cite{stil}.
\index{octonions}%

\emph{This is a good place to recall Doro's result, further developed by Hall, that groups with triality and Moufang loops are essentially the same thing, in that nonassociative Moufang loops correspond to simple groups with triality; the group with triality is the multiplicative group of the Moufang loop, though see Appendix~\ref{LoopTheory} for a more detailed explanation of this.  Furthermore, as we have already noted in a previous section, the only finite simple groups with triality are $D_4(q)$~\cite{liebeck}.
\index{group ! D(4)@$D(4)$}%
  Thus with the construction of the triality-supporting Moufang loop $Q$ for which $\Mlt(Q) \le \Aut(\mathfrak{R^{t}})$ (Theorem~\ref{moulthm}), and with the subsequent implication of a subgroup of $\Aut(\mathfrak{R^{t}})$ as a group with triality (Theorem~\ref{outgpwttr}), we are able to give some justification to the name ``triality graph''
\index{graph ! triality}%
 that we have given to the graph that we have symbolised $\mathfrak{R^{t}}$, in the hope that in the future some connection to (a generalized version of) the outer automorphism triality discovered by \'E. Cartan, and the colour triality may be demonstrated.}
\index{automorphism ! outer}%
\index{Cartan, \'E.}%
\index{Doro, S.}%
\index{Hall, J. I.}%

\bigskip

\head{Open Question}  Is the multiplication group of our loop in theorem~\ref{moulthm} simple?  (The overgroup $\Aut(\mathfrak{R^{t}})$ is simple).

\bigskip

Finally, it is not impossible to connect the number 8 and countable infinity.  The Lie algebra $\mathfrak{sl}_3$
\index{Lie algebra}%
 has dimension 8 whilst its universal enveloping algebra has an infinite dimension; though we are not claiming that this link has anything to do with our results.

\chapter{Random Graph Outer Automorphisms \& Variations}
\label{chap7}
\smallskip

In a dream in the Spring of 1951 the word `automorphism' (taken from
mathematics) came flying towards me.  This is a word for the picturing
of a system on itself, a reflection of the system into itself, for a
process, that is, in which the inner symmetry, the connexional
richness (relations) of a system reveals itself.
\begin{flushright}
Wolfgang Pauli, \textit{Letter to Jung, 27 February 1952~\cite[p.~503]{enz}}
\end{flushright}

\smallskip

We found that objectivity means invariance with respect to the group of automorphisms.  Reality may not always give a clear answer to the question what the actual group of an automorphism is, and for the purpose of some investigations it may be quite useful to replace it by a wider group.
\begin{flushright}
H. Weyl, \textit{Symmetry, Princeton U. P. 1952, p.132}
\end{flushright}

What we learn from our whole discussion and what has indeed become a guiding principle in modern mathematics is this lesson: \emph{Whatever you have to do with a structure-endowed entity $\Sigma$ try to determine its group of automorphisms}, the group of those element-wise transformations which leave all structural relations undisturbed.   You can expect to gain a deep insight into the constitution of $\Sigma$ in this way. 
\begin{flushright}
H. Weyl, \textit{Symmetry, Princeton U. P. 1952, p.144}
\end{flushright}

\smallskip

In standard terminology, an ``outer automorphism'' is an automorphism which is not inner; and the ``outer automorphism group''
\index{group ! automorphism}%
 is the quotient of an automorphism group by the inner automorphism group.
\index{automorphism ! inner}%
\index{automorphism ! outer}%

\section{On Extensions of $\Aut(\mathfrak{R}_{m,\omega})$}

\subsection{Introduction}

In previous chapters we already looked at the outer automorphism group
\index{group ! automorphism}%
 of $\mathfrak{R^{t}}$ and now we will begin to develop a systematic theory.

We define an \emph{outer automorphism of an edge-coloured graph}
\index{graph ! outer automorphism}%
 $\Gamma$ to be a
permutation of its edge-colours modulo a permutation preserving colours acting naturally on the set of colours, which is an element of the quotient
group (denoted $\TAut(\Gamma)$)~\label{colrset}
\index{graph ! all automorphisms}%
 of the group of all graph automorphisms which preserve the partition of $E(\Gamma)$ into colour classes, by the
group of all \emph{inner} automorphisms
\index{graph ! outer automorphism}%
 fixing each colour class (denoted $\Aut(\Gamma)$).  The group of such permutations is the outer automorphism group
\index{group ! automorphism}%
  of the graph,
\index{graph ! outer automorphism}%
  and $\Aut(\Gamma) \lhd \TAut(\Gamma)$.  Therefore the inner automorphism group is the kernel of the action on the set $\mathcal{C}$~\label{colrset} of colours, so $\TAut(\Gamma) / \Aut(\Gamma)$ is isomorphic to a subgroup of $\Sym(\mathcal{C})$.

We shall further clarify this definition below but first we illustrate it with some examples. 
\head{Examples} 
\begin{enumerate}
\item  If $\Gamma$ is a graph with possible edge colours $\mathfrak{r}$ and
$\mathfrak{b}$ on a vertex set $V(\Gamma$) and if $(V(\Gamma),
\mathfrak{r}, \mathfrak{b}) \cong (V(\Gamma), \mathfrak{b}, \mathfrak{r})$ then there is an
automorphism of $V(\Gamma)$ interchanging $\mathfrak{r}$ and
$\mathfrak{b}$ which is outer.  An example of this is the
labelled graph
\index{graph ! labelled}%
 in Figure 1 where the solid and dotted lines represent different colours together with the vertex permutation $x \mapsto 2x \pmod 5$, where $x \in \{0,
\ldots, 4\}$.
\begin{figure}[!h]
$$\xymatrix{ 
&& {\bullet}^{0} \ar@{.}[ddl] \ar@{.}[ddr] \ar@{-}[dl] \ar@{-}[dr]\\
& {\bullet}^{4} \ar@{.}[rr] \ar@{-}[d] \ar@{.}[drr] && {\bullet}^{1} \ar@{-}[d] \ar@{.}[dll]\\
& {\bullet}^{3} \ar@{-}[rr] &&  {\bullet}_{2}
}$$ 
\caption{Labelled graph with outer automorphisms}  
\end{figure}
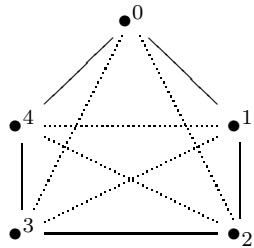

\item  
An example of an outer automorphism
\index{automorphism ! outer}%
 of a graph, is the $C_2$ vertex permutations that induce transpositions of edges and non-edges in the graph $\mathfrak{R}$.  The automorphism groups
\index{group ! automorphism}%
  $\Aut(\mathfrak{R}_{m,\omega})$ of the $m$-coloured random graphs
\index{graph ! random ! $m$-coloured}%
   are simple,
\index{group ! simple}%
 so are generated by any
non-trivial conjugacy class.  By analogy with the
two-coloured random
graph, $\Aut^{*}(\mathfrak{R}_{m,\omega}) / \Aut(\mathfrak{R}_{m,\omega})
\cong \Sym(m)$~\label{allaut} so that $\Sym(m)$ normalizes
$\Aut(\mathfrak{R}_{m,\omega})$ and so is outer; this is essentially
because of the universal nature of the injectivity property defining
the graphs.  Contrast this with the result~\cite[p.~133]{dixon} that
the group of outer automorphisms of a \emph{finite simple} group
\index{group ! finite simple}%
has a normal series of the form: $A \lhd B \lhd C$ where
$A$ is abelian,
\index{group ! abelian}%
$B / A$ is cyclic
\index{group ! cyclic}%
and $C / B \cong 1, \Sym(2)$ or $\Sym(3)$.
\end{enumerate}

Rubin
\index{Rubin, M.}%
studied~\cite{rubin} the reconstruction of $\aleph_0$-categorical
\index{aleph@$\aleph_0$-categorical}%
structures from their automorphism groups,
\index{group ! automorphism}%
 that is given two such structures $M_1, M_2$ for which $\Aut(M_1) \cong \Aut(M_2)$ does it follow that the two permutation groups $\langle \Aut(M_1), |M_1| \rangle$ and $\langle \Aut(M_2), |M_2| \rangle$ are isomorphic?  This is not the case in general, so a strengthening condition is required.  He proved that the isomorphism of the two permutation groups is equivalent to the existence of a bijection $f: |M_1| \to |M_2|$ such that for every $A \subseteq |M_1|^n,$ $A$ is $\emptyset$-definable
\index{definable ! $\emptyset$-definable}%
in $M_1$ if and only if $f(A)$ is $\emptyset$-definable in $M_2$.  A structure $M$ has \emph{no algebraicity}
\index{structure ! with no algebraicity}%
if, for every finite $A \subseteq |M|$ and $a \in |M| - A$, $a$ is not algebraic over $A$.  Let $C$ be the class of $\aleph_0$-categorical structures with no algebraicity, and $M \in C$.  Rubin
\index{Rubin, M.}%
 calls the structure $M$ \emph{group-categorical in C}  
\index{structure ! group-categorical}%
if for every $N \in C$ it is the case that $\Aut(M) \cong \Aut(N) \Rightarrow \langle \Aut(M), |M| \rangle \cong \langle \Aut(N), |N| \rangle$ as permutation groups.  The outer automorphism
\index{automorphism ! outer}%
 group of the automorphism group
\index{group ! automorphism}%
  of a group-categorical structure is easily determined because then every automorphism of $\Aut(M)$ is induced by a permutation of $|M|$.  Rubin cites $\Out(\mathfrak{R}_{m,\omega}) \cong \Sym(m)$~\label{outaut} as an example.  For a finite $A \subseteq M$, $a \in M$ is algebraic over $A$ if $a$ lies in a finite orbit of $\Aut(M)_{A}$.  So a structure that has no algebraicity has no finite orbits, for example the random graph.
\index{graph ! random}%
  However its line graph does have algebraicity.  By a corollary to Whitney's Theorem~\cite[p.~13]{laurisc},
\index{Whitney's Theorem}%
 which follows also for infinite graphs, the random graph $\mathfrak{R}$ and its line graph
\index{graph ! random ! line graph}%
which have the same automorphism group,
\index{group ! automorphism}%
 though the action on vertices is quite different.  That the line graph of the random graph has algebraicity is exemplified as follows.  Consider a triangle $xyz$ and let $M$ be the line graph of $\mathfrak{R}$, $A = \{ \{x, y\}, \{x, z\}\}$, $a = \{y, z\}$.  Fixing $x$, $\{x, y\}$ and $\{x, z\}$ means that an automorphism can either fix $y$ and $z$ or transpose them.  So $\Aut(M)_{A}$ fixes $a$.  

\smallskip

If $\zeta(G)$ denotes the centre of a group $G$ and $\Inn(G)$~\label{Inn} denotes its group of inner automorphisms,
\index{automorphism ! inner}%
 then as is well-known $G / \zeta(G) \cong \Inn(G)$.  We can make the identification $\Aut(\mathfrak{R}_{m,\omega}) \cong
\Inn(\Aut(\mathfrak{R}_{m,\omega}))$ because
$\Aut(\mathfrak{R}_{m,\omega})$ is simple and so centreless, (otherwise
its centre would be a non-trivial normal subgroup).  Furthermore, $\Aut(\mathfrak{R}_{m,\omega})$ is nonabelian, for otherwise given its transitive action on the vertices, this action would in fact be regular, which it certainly is not, for example it is possible to stabilize a vertex and move others; (the possible cycle types are given in Appendix~\ref{PreviousResults}).  A group $G$ is \emph{complete}
\index{group ! complete}%
 if both its centre and outer automorphism
 \index{automorphism ! outer}%
 groups are trivial, or equivalently the conjugation map $G \to \Aut(G)$ is an isomorphism.  Whilst $\Aut(\mathfrak{R^{t}})$ is not complete because $\Aut(\mathfrak{R}_{m,\omega}) \cong \Inn(\Aut(\mathfrak{R}_{m,\omega}))$ the group $\Aut(\Aut(\mathfrak{R^{t}})) \cong \Aut^{*}(\mathfrak{R^{t}}) = \Aut(\mathfrak{R^{t}}) \sd \Sym(\mathfrak{r} , \mathfrak{b} , \mathfrak{g} )$ is complete by a theorem of Burnside~\cite[p.414]{rob} 
\index{Burnside, W.}%
 which states that if $G$ is a nonabelian simple (finite or infinite) group
\index{group ! simple}%
then $\Aut(G)$ is complete. 

\smallskip

The isomorphism class of $\mathfrak{R}_{m,\omega}$ is residual in the set of all $m$-coloured complete graphs on a fixed countable vertex set. (See~\cite{cam6} for discussion).  Truss~\cite{truss1} showed that the group of permutations of the vertex set fixing all the colours, $\Aut(\mathfrak{R}_{m,\omega})$, is a simple group. 

As we stated at the beginning of the section, for any permutation $\pi$ of the set of colours, let $\mathfrak{R}_{m,\omega}^\pi$ be the graph obtained by applying $\pi$ to the colours. Then $\mathfrak{R}_{m,\omega}^\pi$ is
countable, universal and homogeneous, and hence isomorphic to $\mathfrak{R}_{m,\omega}$. This means that,
if $\Aut^*(\mathfrak{R}_{m,\omega})$ is the group of permutations of
the vertex set which induce permutations of the colours, then 
$\Aut^*(\mathfrak{R}_{m,\omega})$ induces the symmetric group $\Sym(m)$ on the colours; 
so $\Aut^*(\mathfrak{R}_{m,\omega})$ is an extension of $\Aut(\mathfrak{R}_{m,\omega})$ by $\Sym(m)$. 

The first question we consider here is: when does this extension split?
\index{group ! split extension}%
 In other words, when is there a complement for $\Aut(\mathfrak{R}_{m,\omega})$ in $\Aut^*(\mathfrak{R}_{m,\omega})$ (a subgroup of $\Aut^*(\mathfrak{R}_{m,\omega})$ isomorphic to $\Sym(m)$ which permutes the colours)? We will show that $\Aut^*(\mathfrak{R}_{m,\omega})$ is the automorphism group
\index{group ! automorphism}%
 of the simple group $\Aut(\mathfrak{R}_{m,\omega})$ (so that the outer automorphism group of
this group is $\Sym(m)$).
\index{automorphism ! outer}%

\begin{theorem}
\label{splits}
The group $\Aut^*(\mathfrak{R}_{m,\omega})$ splits over $\Aut(\mathfrak{R}_{m,\omega})$ if and only if $m$ is odd.
\end{theorem}
\index{group ! split extension}%

\begin{theorem}
\label{autaut}
The group $\Aut^*(\mathfrak{R}_{m,\omega})$ is the automorphism group
\index{group ! automorphism}%
 of $\Aut(\mathfrak{R}_{m,\omega})$.
\end{theorem}

\subsection{Proof of Theorem~\ref{splits}}

We show first that the extension does not split
\index{group ! split extension}%
 if $m$ is even.  Let $\pi \in \Sym(m)$ be a fixed-point-free involution;
\index{group ! action ! fixed-point-free}%
 if $m$ is even there will be such a thing.  Suppose that the extension splits.  There is then an involution acting on $\mathfrak{R}_{m,\omega}$ inducing this colour permutation.  Since we are supposing that a complement exists, let $s$ be an element of this complement acting as $(1,2)(3,4)\cdots(m-1,m)$ on the colours. Then $s$ maps
the subgraph with colours $1,3,\ldots,m-1$ to its complement. But this is impossible, since the edge joining points in a $2$-cycle of $s$ has its colour fixed, for an involution simply swaps over the vertices of an edge.

Now suppose that $m$ is odd; we are going to construct a complement.

First, we show that there exists a function $f$ from pairs of distinct
elements of $\Sym(m)$ to $\{1,\ldots,m\}$ satisfying
\begin{itemize}
\item $f(x,y)=f(y,x)$ for all $x\ne y$;
\item $f(xg,yg)=f(x,y)^g$ for all $x\ne y$ and all $g$.
\end{itemize}
Because of the second property, it suffices to define $f(1,y)$ for $y\ne1$ arbitrarily subject to the
condition $f(1,x^{-1}) = f(1,x)^{x^{-1}}$. Note that this condition requires
$f(1,s)^s=f(1,s)$ whenever $s$ is an involution; this is possible, since
any involution has a fixed point (as $m$ is odd). Then we extend to all pairs
by defining $f(x,y)=f(1,yx^{-1})^x$.  We need to check that both properties are satisfied, and so no conflict arises.

Now we take a countable set of vertices, and let $\Sym(m)$ act semiregularly
\index{group ! permutation ! semiregular}%
on it, that is take the vertex set to be the union of countably many copies of $\Sym(m)$ with $\Sym(m)$ acting by right-multiplication on each copy.  Each orbit is naturally identified with $\Sym(m)$; we let $x_i$ denote
the element identified with $x \in \Sym(m)$ in the $i$th orbit, as $i\in\mathbb{N}$
(where orbits are indexed by natural numbers). Then we colour the edges
within each orbit (a copy of $\Sym(m)$) by giving $\{x_i,y_i\}$ the colour $f(x,y)$.  The second of the conditions above, the translation condition, ensures that we get the right colouring under $\Sym(m)$-action on pairs.  Now assume that colours within orbits $i$ and $j$ (each being a copy of $\Sym(m)$) are given.  For edges between orbits $i$ and $j$, with $i<j$, we colour $\{x_i,1_j\}$ randomly, and then give $\{y_i,z_j\}$ the image of the colour of $\{(yz^{-1})_i,1_j\}$ under $z$, that is act by translation on vertices and by permutation on colours.

Clearly the group $\Sym(m)$ permutes the colours of the edges
consistently, and in the same way as it permutes $\{1,\ldots,m\}$.

Next we show that a residual set of the coloured graphs we obtain are
isomorphic to $\mathfrak{R}_{m,\omega}$. We have to show that, given $m$ finite disjoint sets of
vertices, say $U_1,\ldots, U_m$, the set of graphs containing a vertex $v$
joined by edges of colour $i$ to all vertices in $U_i$  (for $i=1,\ldots,m$)
is open and dense. The openness is clear. To see that it is dense, note that
the $m$ finite sets are contained in the union of a finite number of orbits
(say those with index less than $N$); then, for any $i\ge N$, we are free to
choose the colours of the edges joining these vertices to $1_i$ arbitrarily.

Now by construction, the group $\Sym(m)$ we have constructed meets
$\Aut(\mathfrak{R}_{m,\omega})$ in the identity; so it is the required complement.

\bigskip

How close can we get when $m$ is even? The construction in the second part
shows that, if there is a group $G$ which acts as $\Sym(m)$ on the set
$\{1,\ldots,m\}$, in such a way that all involutions in $G$ have fixed points
on $\{1,\ldots,m\}$, then $G$ is a supplement for $\Aut(\mathfrak{R}_{m,\omega})$ in
$\Aut^*(\mathfrak{R}_{m,\omega})$ (that is, $G.\Aut(\mathfrak{R}_{m,\omega})=\Aut^*(\mathfrak{R}_{m,\omega})$), and $G\cap\Aut(\mathfrak{R}_{m,\omega})$ is the kernel of the action of $G$ on $\{1,\ldots,m\}$.

In the case that $m$ is even but not divisible by 8, we will see that $\exists H \leq \Aut^*(\mathfrak{R}_{m,\omega})$ which is a \emph{supplement} such that $|\Aut(\mathfrak{R}_{m,\omega}) \cap H| = 2$.  As the index is 2, $H$ is a double cover of $\Sym(m)$.  In this double cover of $\Sym(m)$, the fixed-point-free 
 \index{group ! action ! fixed-point-free}%
involutions lift to elements of order~$4$.  So we can just repeat the previous proof using a double cover in which fixed-point-free involutions lits to elements of order~$4$.

There are two double covers of $\Sym(n)$ for $n\ge4$,
described in~\cite[Chapter~2]{hh} and called there $\tilde S_m$ and $\hat S_m$.
In $\tilde S_m$, the product of $r$ disjoint transpositions lifts to an 
element of order~$4$ if and only if $r\equiv 1$ or~$2$ mod~$4$, while in 
$\hat S_m$, the condition is that $r\equiv 2$ or~$3$ mod~$4$.  (This is a fairly simple calculation.  The element
$(t_1t_3...t_{2r-1})^2 \in \tilde S_n$ can be calculated: each time we jump
one $t$ over another we get a factor $z$, and each $t^2$ gives a factor $z$, so
altogether we have $z^{r(r+1)/2}$; this is equal to $z$ if $r$ is congruent to
$1$ or $2$ mod~$4$. In $\hat S_n$, the calculation is similar with $s$ in place of
$t$, but $t^2=1$, so $(s_1s_3...s_{2r-1})^2$ is $z^{r(r-1)/2}$, which is equal
to $z$ if $r$ is congruent to $2$ or $3$ mod~$4$. So only the case $r=0$ mod~$4$ is
not covered. Now we have $m=2r$ (since we want a fixed-point-free
involution), so only the case $m=0$ mod~$8$ is not covered.)
 This shows that there is a supplement meeting $\Aut(\mathfrak{R}_{m,\omega})$ in a group of order~$2$
for $m$ even but not divisible by~$8$.

What happens in the remaining case, when $m$ is a multiple of $8$?  Is there a finite group $G$ with an epimorphism to $\Sym(m)$ such that the inverse image of a fixed-point-free involution is not an involution?  In other words, is there a finite supplement, and what is the smallest such?  This would be true only up to conjugacy as different colourings between different orbits $i$ and $j$ above give different conjugacy classes.

\subsection{Proof of Theorem~\ref{autaut}}

Since $\Aut(\mathfrak{R}_{m,\omega})$ is primitive 
\index{group ! permutation ! primitive}%
and not regular, its centralizer in the symmetric group is trivial; therefore $\Aut^*(\mathfrak{R}_{m,\omega})$ acts faithfully on $\Aut(\mathfrak{R}_{m,\omega})$ by conjugation. We have to show that there are no further automorphisms.

A permutation group $G$ of countable degree is said to have the \emph{small
index property}
\index{small index property (SIP)}%
 if any subgroup $H$ satisfying $|G:H|<2^{\aleph_0}$ contains
the pointwise stabilizer of a finite set; it has the \emph{strong small index property} 
\index{small index property, strong}%
 if any subgroup $H$ satisfying $|G:H|<2^{\aleph_0}$ lies between
the pointwise and setwise stabilizer of a finite set.

\head{Step 1} $\mathfrak{R}_{m,\omega}$ has the strong small index property.

This is proved by a simple modification of the  arguments for the case $m=2$.
The small index property was proved by Hodges \emph{et~al.}~\cite{hodges}, 
using a result of Hrushovski~\cite{hrush}; the strong version was a simple
extension due to Cameron~\cite{cam23}.
\index{Cameron, P. J.}%

Hrushovski 
\index{Hrushovski, E.}%
showed that any finite graph $X$ can be embedded into a finite
graph $Z$ such that all isomorphisms between subgraphs of $X$ extend to
automorphisms of $Z$. Moreover, the graph $Z$ is vertex-, edge- and 
nonedge-transitive. He uses this to construct a generic countable sequence
of automorphisms of $\mathfrak{R}$. To extend this  to $\mathfrak{R}_{m,\omega}$ is comparatively
straightforward. It is necessary to work with $(m-1)$-edge-coloured
graphs (regarding the $m$th colour as `transparent'). Now the arguments
of Hodges \emph{et~al.} and Cameron 
\index{Cameron, P. J.}%
go through essentially unchanged.

\head{Step 2} Since $\Aut(\mathfrak{R}_{m,\omega})$ acts primitively on the vertex set,
with permutation rank $m+1$, the vertex stabilizers are maximal subgroups of
countable index with $m+1$ double cosets (orbits on ordered pairs corresponding to equality and the $m$ edge colours). Moreover, any further subgroup of countable index has more than $m+1$ double cosets.

For let $H$ be a maximal subgroup of countable index. By the strong SIP, $H$
is the stabilizer of a $k$-set $X$. If $g$ maps $X$ to a disjoint $k$-set,
then $HgH$ determines the colours of the edges between $X$ and $X^g$, up to
permutations of these two sets. By universality, there are at least
$m^{k^2}/(k!)^2$ such double cosets. Now it is not hard to prove that
$m^{k^2}/(k!)^2>m$ for $k\ge2$. Hence we must have $k=1$.  (The number of double cosets cannot be changed by an automorphism, as the number $m+1$ is minimal).

\head{Step 3} It follows that any automorphism permutes the vertex
stabilizers among themselves, so is induced by a
permutation of the vertices which normalises $\Aut(\mathfrak{R}_{m,\omega})$. To finish
the proof, we show that the normaliser of $\Aut(\mathfrak{R}_{m,\omega})$ in the
symmetric group is $\Aut^*(\mathfrak{R}_{m,\omega})$.

This is straightforward. A vertex permutation which normalises
$\Aut(\mathfrak{R}_{m,\omega})$ must permute among themselves the $\Aut(\mathfrak{R}_{m,\omega})$-orbits on pairs of vertices, that is, the colour classes; so it belongs to $\Aut^*(\mathfrak{R}_{m,\omega})$.

\section{Introducing Some Groups}

In this section we will introduce some new groups related to automorphism groups
\index{group ! automorphism}%
 for random graphs.

First we borrow some concepts from~\cite{alcoma}, outlining their theory for a general countably infinite set $\Omega$, before we include any graph structure.  The set of \emph{near-bijections}~\label{NB}
\index{group ! near bijection}%
of $\Omega$ is given by
\[\NB(\Omega) = \{f: f\ \text{is a bijection}\ \Omega_1 \to \Omega_2\
\text{between cofinite subsets of}\ \Omega \}.\]
Regarded as sets we have that $\Sym(\Omega) \subseteq \NB(\Omega)$.
The set $\NS(\Omega)$~\label{NS} of \emph{near symmetries}
\index{group ! near symmetries}%
 is the set of
equivalence classes of the equivalence relation $\equiv$ on
$\NB(\Omega)$ where for bijections $f_1, f_2$,  $f_1 \equiv f_2$ if
there is a cofinite $\Omega' \subseteq \Omega$ with  $f_1|_{\Omega'}
= f_2|_{\Omega'}$.  This set forms a group, where the composition law
is induced by the composition of maps.  Composition $\circ$ of maps is well-defined, being independent of the choices made, and this induces a well-defined operation on equivalence classes by cofiniteness of operand sets.  For $(f_1 \circ f_2)(x)$ is defined if and only if both $x \in \dom(f_2)$ and $f_2(x) \in \dom(f_1)$.  Now the intersection of the two cofinite sets given by $\dom(f_1 \circ f_2) = f_2^{-1}(\dom(f_1) \cap \range(f_2))$ and $\range(f_1 \circ f_2)$ is cofinite.  Finally it is easy to see that $f_1 \equiv f'_1,  f_2 \equiv f'_2 \Leftrightarrow f_1 \circ f_2 \equiv f'_1 \circ f'_2$.

The sequence 
\[\FSym(\Omega) \to \Sym(\Omega) \to \NS(\Omega)\]
where the second map is a group homomorphism with kernel $\FSym(\Omega)$ is exact, so $\Sym(\Omega) / \FSym(\Omega) \le \NS(\Omega)$ and the elements of this quotient are the equivalence classes of permutations.  For $f \in
\NB(\Omega)$ with $\Omega_1, \Omega_2$ cofinite subsets of $\Omega$,
define the \emph{index} 
\index{index of near-bijection}%
of $f$ to be~\label{index}
\[ \ind f = |\Omega - \Omega_2| - |\Omega - \Omega_1|.\]
\begin{lemma}[Alperin, Covington and Macpherson]
\label{alpcovmac}
\index{Alperin, J. L.}%
\index{Covington, J.}%
\index{Macpherson, H. D.}%
The map $\ind: \NB(\Omega) \to \mathbb{Z}$ induces a mapping $\ind: \NS(\Omega) \to \mathbb{Z}$ which is a homomorphism onto $(\mathbb{Z}, +)$ with kernel $\Sym(\Omega) /
\FSym(\Omega)$.
\end{lemma}
From this lemma, which is proved in~\cite{alcoma}, it follows that   
\[1 \to \FSym(\Omega) \to \Sym(\Omega) \to \NS(\Omega) \to \mathbb{Z} \to 1.\]

Next we apply the above formalism to random graphs, firstly concentrating on the two-colour case.  Define a
\emph{near-automorphism}
\index{group ! near automorphism}%
 of $\mathfrak{R}$ to be an
element of $\NB(\mathfrak{R})$ which is an isomorphism from its domain
to its range; denote the set of such maps by $\NA(\mathfrak{R})$.~\label{NA}  It
is clear that a near-automorphism of $\mathfrak{R}$ has only finitely
many extensions to a maximal (with respect to inclusion)
near-automorphism of $\mathfrak{R}$.  Let $\NE(\mathfrak{R})$~\label{NE} be the
set of equivalence classes of elements of $\NA(\mathfrak{R})$ with the same equivalence relation as before.  The elements of $\NE(\mathfrak{R})$
\index{group ! equivalence classes of\\\quad near-automorphisms}%
 form a group with composition
being composition of maps.  Notice that $\NA(\mathfrak{R})$ is not a
permutation group on $\mathfrak{R}$, and that a maximal near-automorphism
is an automorphism if and only if it is a permutation, both necessary
and sufficient conditions following by definition.  

Let $\mathfrak{R}'$ be $\mathfrak{R}$ with edge $\{x, y\}$ deleted.  Then $\mathfrak{R}' \cong \mathfrak{R}$.  Let the isomorphism $f: \mathfrak{R}' \to \mathfrak{R}$ be such that $f: x \to u$, $f: y \to v$, $x \sim y, u \nsim v$, and otherwise $f$ is an automorphism.  So $f \in \NA(\mathfrak{R})$ is a permutation, and it will also be an element of
$\AAut(\mathfrak{R})$.  Therefore $f \in \NA(\mathfrak{R})$ is contained in both a permutation and in two maximal near-automorphisms (in one $\{u, v\}$ is an edge and in the other it is a non-edge), but $f \notin \Aut(\mathfrak{R})$.

Consider Figure~\ref{ininwgp} where the uppermost exact sequence for $\Omega$ instead of $\mathfrak{R}$ is taken from~\cite{alcoma}.

\begin{figure}[!h]
$$\xymatrix{
{1} \ar[r] & {\FSym(\mathfrak{R})} \ar[r] & {\Sym(\mathfrak{R})}
\ar@{->}[r]^{\theta} & {\NS(\mathfrak{R})} \ar[r]^{\quad \ind} & {\mathbb{Z}} \ar[r] & {1}\\
{1} \ar[r] & {1} \ar[r] & {\Aut(\mathfrak{R})} \ar[u]
\ar@{^{(}->}[r]  \ar@{^{(}->}[d] & {\NE(\mathfrak{R})}
\ar@{^{(}->}[u]^{\phi}\\
&& {\AAut(\mathfrak{R})} \ar@{^{(}->}[ur]^{\phi'}
}$$ 
\caption{The new groups $\NS(\mathfrak{R})$ and $\NE(\mathfrak{R})$ and their subgroup inclusions.}
\label{ininwgp}
\end{figure}

Choose $f \in \AAut(\mathfrak{R})$.  There is a finite set $S$ of edges and non-edges changed by the isomorphism
\[f : \mathfrak{R} \backslash \{\text{vertices in $S$}\} \to \mathfrak{R} \backslash \{\text{vertices in $f(S)$}\}.\]  
So $f \in \NA(\mathfrak{R})$.  Now $\phi': f \to [f]$.  Notice $\phi'(\AAut(\mathfrak{R})) \subset \NE(\mathfrak{R}) \cap
Im(\theta) = \NE(\mathfrak{R}) \cap \ker(\ind(\NE(\mathfrak{R})))$~\label{kernel}.  This inclusion is strict because not all isomorphisms $\mathfrak{R} \backslash \{x\} \to \mathfrak{R} \backslash \{y\}$,
where $x$ and $y$ are vertices of $\mathfrak{R}$ and both domain and
codomain are isomorphic to $\mathfrak{R}$, extend to near automorphisms of
$\mathfrak{R}$, but
only those which map the vertices attached to $x$ to the vertices
attached to $y$, apart from a finite number; only for these will the one-point extension
\index{one-point extension property}%
 property be valid.  For example if some of the adjacencies at $x$ are changed randomly then with probability 1, the resulting graph is isomorphic to $\mathfrak{R}$.  Similarly a switching at $x$ which interchanges all edges and non-edges at $x$ still gives a graph isomorphic to $\mathfrak{R}$.

Denote an equivalence class of near-automorphisms $h$ in
$\NE(\mathfrak{R})$ by $e = [h]_{NE}$.  Let $\phi:
\NE(\mathfrak{R}) \to \NS(\mathfrak{R}) : e \mapsto \phi(e) =
[h]_{NS}$ be a map which we will prove is a monomorphism.  That $\phi$
is well-defined, meaning that $[h]_{NE} = [h']_{NE} \Rightarrow
[h]_{NS} = [h']_{NS}$ follows from the definitions of
$\NE(\mathfrak{R})$ and $\NS(\mathfrak{R})$ and that the
equivalence class is the same in these groups, that
is the same cofinite operand for both groups.

Furthermore, if $e_1, e_2 \in \NE(\mathfrak{R})$ and $\Omega_1,
\Omega_2, \Omega_3$ are cofinite subsets of the vertex set of $\mathfrak{R}$ chosen as in Figure 3, then 

$\phi(e_1 e_2)
=\phi([h_1]^{1}_{NE}[h_2]^{2}_{NE}) =
\phi([h_1]^{3}_{NE}[h_2]^{3}_{NE}) = \phi([h_3]^{3}_{NE}) =$\\
$[h_3]^{3}_{NS} = [h_1]^{3}_{NS}[h_2]^{3}_{NS} = \phi(e_1) \phi(e_2)$,\\
where the superscript on $[h]$ indexes the equivalence class so that for example
$[h_1]^{3}_{NE} = {e_1}_{|\Omega_3}$ and $\Omega_3 = \Omega_1 \cap e_1^{-1}(\Omega_2 \cap
e_1(\Omega_1))$ (see Figure 3).  The homomorphism is clearly injective because
$\ker(\phi)=1$.  Whence

\begin{proposition}
The map $\phi : \NE(\mathfrak{R}) \hookrightarrow
\NS(\mathfrak{R})$ is a monomorphism.
\end{proposition}

\begin{figure}[htb]
\begin{center}
\input{test.latex}
\end{center}
\label{detom}
\caption{Determination of $\Omega_3$}
\end{figure}
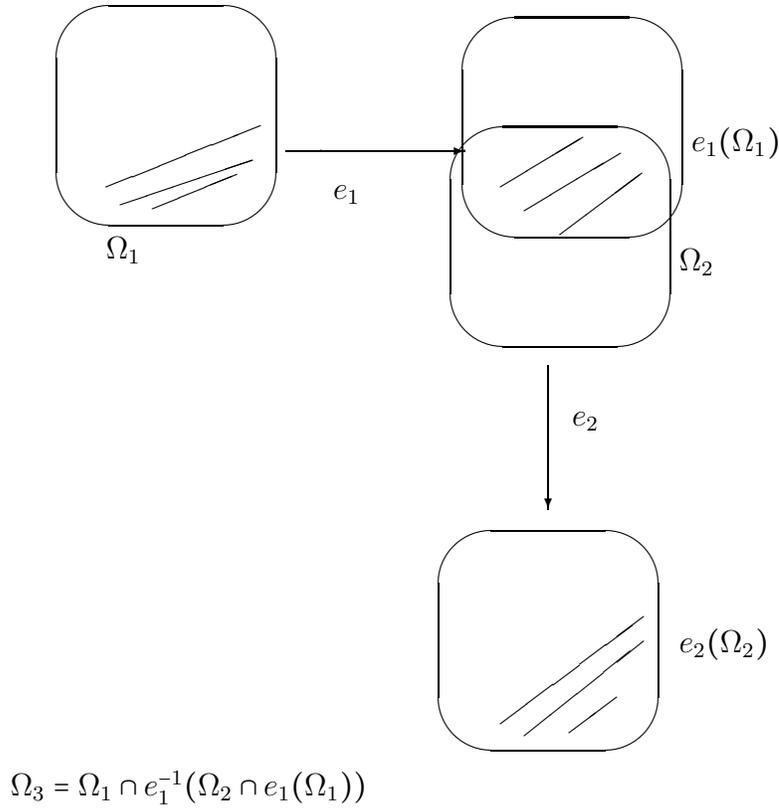

How does the map $\theta : \Sym(\mathfrak{R}) \to \NS(\mathfrak{R})$ behave when it is restricted
to $\Aut(\mathfrak{R})$?  The random graph
\index{graph ! random}%
 $\mathfrak{R}$ has no non-trivial finitary automorphisms, because a map from vertex $x$ to vertex $y$ moves an infinity of vertices attached to $x$ to those attached to $y$.  As $\Aut(\mathfrak{R}) \cap \FSym(\mathfrak{R}) = 1$,
\index{group ! finitary}%
 the map $\theta_{| \Aut(\mathfrak{R})}$ is a monomorphism into $\NE(\mathfrak{R})$.  However $\theta(\Aut(\mathfrak{R}))$ is not normal in
$\NE(\mathfrak{R})$ because as we will show below there is a
monomorphism from $\AAut(\mathfrak{R})$ to $\NE(\mathfrak{R})$, so that $\theta(\Aut(\mathfrak{R})) \triangleleft
\theta(\AAut(\mathfrak{R}))$, which because $\theta$ is a homomorphism
would give $\Aut(\mathfrak{R}) \triangleleft \AAut(\mathfrak{R})$, a contradiction because $\AAut(\mathfrak{R})$ is highly transitive (see Appendix~\ref{PreviousResults}), and so any non-trivial normal subgroup must also be highly transitive, which $\Aut(\mathfrak{R})$ is not.

Also $\phi'(\AAut(\mathfrak{R}))$ is not normal in
$\NE(\mathfrak{R})$ as can be seen from the following example.
Consider a vertex triple $\{x, y, z\} \in V(\mathfrak{R})$ with $x \sim y$, $x \sim z$, $y
\nsim z$.  Take $a \in \AAut(\mathfrak{R})$ mapping $x \sim y$ to $xa \nsim ya$ and
otherwise preserving all adjacencies, and also $e \in \NE(\mathfrak{R}) \backslash \phi'(\AAut(\mathfrak{R}))$ that moves vertex $y$.  We can choose $e$ to include the transposition $(y,z)$.  Then the effect of $e^{-1} a e$ acting on vertices $\{x, y, z\}$ is: $\{x \sim y\} \mapsto \{x \sim y\}$, $\{x \sim z\} \mapsto \{x \nsim z\}$, $\{y \nsim z\} \mapsto \{y \nsim z\}$.  Since $e \notin \AAut(\mathfrak{R})$, in moving $y$ it moves an infinite number of its neighbours.  So there are an infinite number of vertices $z$ whose adjacency with $x$ are changed by $e$.
Therefore $e^{-1} a e$ changes infinitely many edges on $x$.  So $e^{-1} a e \notin \AAut(\mathfrak{R})$.

\medskip

Define the \emph{vertex index}
\index{vertex index}%
$N(g)$~\label{vxindex} of an element $g$ of an acting group to be the number of vertices $v$ such that infinitely many adjacencies are changed at $v$ by $g$.  This index has a norm-like property $N(g_1 g_2) \le N(g_1) + N(g_2)$ and also $N(g^{-1}) = N(g)$.  We can then define two new groups
\[ \Aut_0 = \langle g : N(g) = 0 \rangle \]
\[ \Aut_{fin} = \langle g : N(g) < \infty \rangle. \]

The \emph{zero vertex index group} $\Aut_0$~\label{zerovx}
\index{group ! zero vertex index}%
is the group that changes infinitely many edges at no vertices, and the \emph{finite vertex index group} $\Aut_{fin}$
\index{group ! finite vertex index}%
is the group that changes infinitely many edges at finitely many vertices.

From the definitions, $\Aut_0 \subseteq \Aut_{fin}$.

It is clear that $\FSym(\mathfrak{R}) \lhd \Aut_{fin}$.
\index{group ! finitary}%
 Recall that $\NE(\mathfrak{R})$ is the set of equivalence classes of bijections on the cofinite subsets of the vertex set which are isomorphisms from their domain to their range; its elements form a group by composition of maps. 

Just as we were able to define an index map above, we can define another mapping $\ind': \NE(\mathfrak{R}) \to {\mathbb{Z}} : [w] \mapsto z$ or alternatively $\ind' f' = |\Omega - \range(f')| - |\Omega - \dom(f')|$ where $f'$ is a bijection between cofinite subsets of the vertex set $\Omega$.  As before, this mapping is a homomorphism onto $({\mathbb{Z}}, +)$ with kernel $\FSym(\mathfrak{R})$.
\index{group ! finitary}%
  The kernel $\ker(\ind')$ consists of all bijections such that $|\Omega - \range(f')| = |\Omega - \dom(f')|$.  The mapping $\tilde{\phi} : \Aut_{fin}(\mathfrak{R}) \to \NE(\mathfrak{R})$ is a group monomorphism because for $w \in \Aut_{fin}(\mathfrak{R})$, $\tilde{\phi}(w)$ is an automorphism on all but a finite vertex subset, and this subset can
be chosen to lie outside the domain and range of some $a \in \NA(\mathfrak{R})$ and so of $[a] \in \NE(\mathfrak{R})$.  From the equality $\Aut_{fin}(\mathfrak{R}) / \FSym(\mathfrak{R}) = \ker(\ind')$ it follows that

\centerline{$\Aut_{fin}(\mathfrak{R}) / \FSym(\mathfrak{R}) \lhd \NE(\mathfrak{R})$.}  

\smallskip

What else can we say about the quotient $\Aut_{fin}(\mathfrak{R}) / \FSym(\mathfrak{R})$?  It is clearly not maximal in $\NE(\mathfrak{R})$, because its quotient in this overgroup is $\mathbb{Z}$ which has one subgroup generated by each element.  By connectedness of the random graphs,
\index{graph ! random}%
 the action of $\FSym(\mathfrak{R})$ will always change the adjacencies on an infinite number of edges, so $\AAut(\mathfrak{R}) \cap \FSym(\mathfrak{R}) = \{1\}$.  However $\AAut(\mathfrak{R})$ and  $\FSym(\mathfrak{R})$ are not sufficient to generate  $\Aut_{fin}(\mathfrak{R})$, because their elements only ever act on finite vertex sets.  We can state this more formally as a proposition, for which we need the following claim.
\begin{claim}
For a vertex $v \in \mathfrak{R}$ let $N(v), \overline{N}(v)$~\label{neighbvx} be its neighbours and non-neighbours.  There exists an automorphism $g$ of $\mathfrak{R} \backslash \{v\}$ which interchanges the sets $N(v)$ and $\overline{N}(v)$.
\end{claim}

\begin{proof}
This is quite clear from elementary random graph properties.  Starting with two empty sets of vertices, build in parallel, two countably infinite graphs, which are isomorphic to $N(v) \cong \overline{N}(v)$.  By a further application of the $1$-point extension property, attach vertex $v$ to $N(v)$, then extend to construct the whole graph $\mathfrak{R}$.  (Alternatively, two infinite cliques  such as $N(v)$ and $\overline{N}(v)$ exist as subgraphs of $\mathfrak{R}$ by its universality).  Then the required $g$ exists by homogeneity of $\mathfrak{R}$.
\end{proof}

Then,

\begin{proposition}
$\langle \FSym(\mathfrak{R})), \AAut(\mathfrak{R})) \rangle \neq \Aut_{fin}(\mathfrak{R})$.
\end{proposition}

\begin{proof}
Suppose that $g \in \Aut(\mathfrak{R})$ as in the claim, and that $g_1 \in \FSym(\mathfrak{R}), g_2 \in \AAut(\mathfrak{R})$.  Let $S(g)$ denote the set of vertices $v$ such that infinitely many adjacencies at $v$ are changed by $g$.  Then $|S(g)| = 1$ and $S(g_2) = \emptyset$.  Also if $g_1 \neq 1$ then $|S(g)| \ge 2$, where $S(g_1)$ denotes the support of $g_1$.  Finally $S(g_1 g_2) = S(g_1)$.
\end{proof}

We pose an open question.  Is $\Aut_{fin}(\mathfrak{R}) / \FSym(\mathfrak{R})$ simple?

\bigskip

The above considerations generalize to $\mathfrak{R}_{m,\omega}$, as shown in Figure 4.  In this diagram, it is intended that each type of group appertaining to the $m$-coloured random graph
\index{graph ! random ! $m$-coloured}%
 is embedded in the corresponding $(m-1)$-coloured group.  So there are in fact $(m-1)$ layers of groups for $m$ colours.  We remark on some obvious embeddings:  

\begin{itemize}
\item[(i)]  The group $\Aut(\mathfrak{R}_{m, \omega})$ is a subgroup of $\Aut(\mathfrak{R}_{m-1, \omega})$.  The simplicity of all the automorphism groups
\index{group ! automorphism}%
 was proved in~\cite{truss1}.

\item[(ii)]  The groups in the inclusion $\AAut(\mathfrak{R}_{m-1, \omega}) \hookrightarrow
\AAut(\mathfrak{R}_{m, \omega})$ were studied in~\cite{truss2}.

\item[(iii)]  $\Aut(\mathfrak{R}_{m,\omega}) \hookrightarrow
\mathfrak{m}\Aut(\mathfrak{R}_{m,\omega}) \hookrightarrow
\Aut(\mathfrak{R}_{m,\omega})_{\mathfrak{m}} \cong \Aut(\mathfrak{R}_{m-1,\omega})$.  

Here $\mathfrak{m}\Aut(\mathfrak{R}_{m,\omega})$ denotes the automorphism
group acting on $\mathfrak{R}_{m,\omega}$ in which colours
$\mathfrak{m}$ and $\mathfrak{m-1}$ can be interchanged, and $\Aut(\mathfrak{R}_{m,\omega})_{\mathfrak{m}}$
means the automorphism group
\index{group ! automorphism}%
 acting on that $m$-coloured graph which
is colourblind in colours $\mathfrak{m}$ and $\mathfrak{m-1}$, that is does not
distinguish between these two colours.

By defining a $\mathbb{Z}_2$-valued parity map from $\mathfrak{m}\Aut(\mathfrak{R}_{m,\omega})$ to $\Aut(\mathfrak{R}_{m,\omega})$, according to the parity of the interchange of $\mathfrak{m}$ and $\mathfrak{m-1}$, we can see that
$\Aut(\mathfrak{R}_{m,\omega}) \triangleleft
\mathfrak{m}\Aut(\mathfrak{R}_{m,\omega})$.  However we claim that the image of the group embedding $\mathfrak{m}\Aut(\mathfrak{R}_{m,\omega})$ into
$\Aut(\mathfrak{R}_{m,\omega})_{\mathfrak{m}}$ is not normal in
$\Aut(\mathfrak{R}_{m,\omega})_{\mathfrak{m}}$.  For take a vertex triple $\{x, y, z\} \in V
(\mathfrak{R}_{m,\omega})$, and choose $a \in
\mathfrak{m}\Aut(\mathfrak{R}_{m,\omega})$ so as to stabilize
$\mathfrak{R}_{m,\omega}$ apart from a colour transposition
$\mathfrak{m} \mapsto \mathfrak{m-1}$ on edge $\{x, y\}$, and take $b
\in \Aut(\mathfrak{R}_{m,\omega})_{\mathfrak{m}}$ to include the $(y,
z)$ transposition.  Then $b^{-1} a b$ acting on the vertex triple stabilizes the colour on edge $\{x \sim y\}$ as $\mathfrak{m-1}$ and transposes the colour on edge $\{x \sim z\}$ from $\mathfrak{m}$ to $\mathfrak{m-1}$.  So $b^{-1} a b \notin \mathfrak{m}\Aut(\mathfrak{R}_{m,\omega})$, which
proves the claim.

\item[(iv)]  All of the groups $\NE(\mathfrak{R}_{m,\omega})$,
$\AAut(\mathfrak{R}_{m,\omega})$ and $\Aut_{fin}(\mathfrak{R}_{m,\omega})$
follow the same embedding pattern with respect
to $m$ as that of $\Aut(\mathfrak{R}_{m,\omega})$ in (i) above.

The following series of embeddings require further study:

(a)   $\NE(\mathfrak{R}_{m,\omega}) \hookrightarrow \NE(\mathfrak{R}_{m-1,\omega}) \hookrightarrow \ldots \NE(\mathfrak{R})$.  

(b)   $\Aut_0(\mathfrak{R}_{m,\omega}) \hookrightarrow \Aut_0(\mathfrak{R}_{m-1,\omega}) \hookrightarrow \ldots \Aut_0(\mathfrak{R})$.  

(c)   $\Aut_{fin}(\mathfrak{R}_{m,\omega}) \hookrightarrow \Aut_{fin}(\mathfrak{R}_{m-1,\omega}) \hookrightarrow \ldots \Aut_{fin}(\mathfrak{R})$.  

\item[(v)]  That $\Aut_{fin}(\mathfrak{R})$ is highly transitive
\index{group ! permutation ! highly transitive}%
 follows either from
the fact that it contains $\FSym(\mathfrak{R})$ or from it being an
overgroup of $\AAut(\mathfrak{R})$ which is itself highly transitive.

\item[(vi)]  For $m \ge 2$, consider the epimorphisms:
\[ \phi_m : \Aut_{fin}(\mathfrak{R}_{m,\omega}) \to \Aut_{fin}(\mathfrak{R}_{m,\omega}) / \FSym(\mathfrak{R}_{m,\omega}).\]
Given that $\Aut_{fin}(\mathfrak{R}_{m,\omega})$ acts as an automorphism-like group and because all of the finitary groups in the following series are equal, we have that
\[ \Aut_{fin}(\mathfrak{R}_{m,\omega}) / \FSym(\mathfrak{R}_{m,\omega}) < \Aut_{fin}(\mathfrak{R}_{{m-1},\omega}) / \FSym(\mathfrak{R}_{{m-1},\omega}) < \] 
\[ \ldots < \Aut_{fin}(\mathfrak{R}) / \FSym(\mathfrak{R}). \]
\end{itemize}

\clearpage

\begin{figure}[ht]
\vbox to \vsize{
 \vss
 \hbox to \hsize{
  \hss
  \rotatebox{90}{
   \vbox{
    $$
    \hss
    \xymatrix@R=30pt@C=5pt{
{} &&& {\Sym(\mathfrak{R})} \ar@{->}[rrrrrrr]^{\theta} &&&&&&& {\NS(\mathfrak{R})}  \ar@{->}[r]^{\quad \ind} & {\mathbb{Z}} \\
{} &&&&&&&&&& {\NE(\mathfrak{R})} \ar@{^{(}->}[u]^{\phi}  \ar@{->}[r]^{\quad \ind'} & {\mathbb{Z}} \\
{} &&& {\Aut_{fin}(\mathfrak{R})}  \ar@{^{(}->}[uu]   \ar@{^{(}->}[urrrrrrr]_{\tilde{\phi}} &&&&&&& {\quad \NE(\mathfrak{R}_{m,\omega})} \ar@{^{(}->}[u] \ar@{->}[r] & {\mathbb{Z}} \\
{} &&& {\Aut_{fin}(\mathfrak{R}_{m,\omega})}  \ar@{^{(}.>}[u]   \ar@{^{(}.>}[urrrrrrr]_{\tilde{\phi}}  && {\Aut_0} \ar@{^{(}->}[ull] \\  
{1} \ar@{->}[r] & {\FSym(\mathfrak{R})}   \ar@{^{(}->}[rrr] \ar@{^{(}->}[uurr]  \ar@{^{(}->}[uuuurr] &&& \AAut(\mathfrak{R}) \ar@{^{(}->}[ur]  \ar@{^{(}->}[uuurrrrrr]^{\phi'}  &  {\Aut_0(\mathfrak{R}_{m,\omega})} \ar@{^{(}.>}[ull] \ar@{^{(}.>}[u] \\
{1} \ar@{->}[r] & \FSym(\mathfrak{R}_{m,\omega})   \ar@{^{(}->}[rrr] \ar@{^{(}->}[uurr]  \ar@{^{(}->}[u] &&& \AAut(\mathfrak{R}_{m,\omega})  \ar@{^{(}.>}[uuurrrrrr]^{\phi'} \ar@{^{(}->}[u]   \ar@{^{(}.>}[ur] && \Aut(\mathfrak{R}) \ar@{^{(}->}[ull]  \ar@{^{(}->}[uuuurrrr] \\
{} &&&&&& \Aut(\mathfrak{R}_{m,\omega}) \ar@{^{(}->}[u] \ar@{^{(}->}[uuuurrrr]  \ar@{^{(}->}[ull]
}
    \hss
    $$
    \caption{The new groups $\NE(\mathfrak{R}_{m,\omega})$,\ $\NS(\mathfrak{R}_{m,\omega})$,\ $\Aut_{fin}$,\ $\Aut_0$ and their subgroup inclusions:\ variations on the theme of automorphism groups.}  
   }
  }
 \hss}
\vss}
\label{varfig}
\end{figure}
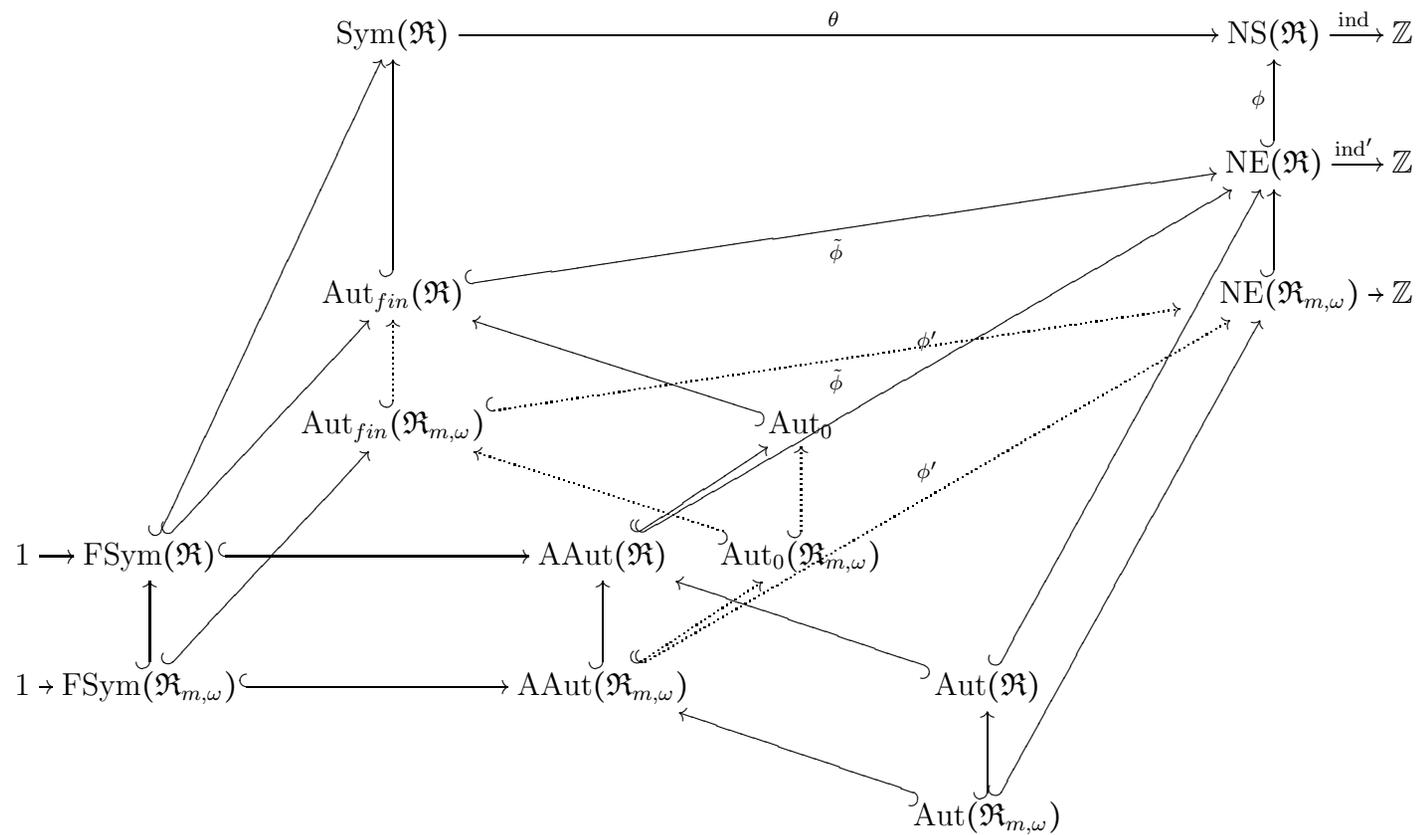

\chapter{Filters, Topologies and Groups from $\mathfrak{R}$}
\label{chap8}
\bigskip

Once the concept of infinity has been taken seriously, a human
dwelling can no longer be made of the universe.  The universe can
still be thought but it can no longer be imaged;  the man who thinks
it no longer really lives in it.
\begin{flushright}
Martin Buber
\end{flushright}

\bigskip

The role of relational structures
\index{relational structure}%
 (graphs, directed graphs, partial orders,
etc.) in the investigation of permutation groups is well-known. Of course, the
automorphism group
\index{group ! automorphism}%
 of a non-trivial relational structure cannot be highly
transitive; if a group is highly transitive, we would expect to have to
use ``infinitary'' structures such as filters
\index{filter}%
 and topologies in its study.
In fact, the gap is not so large. For example, Macpherson and 
Praeger~\cite{macprae}
\index{Macpherson, H. D.}%
\index{Praeger, C. E.}%
 showed that a permutation group of countable degree which
is not highly transitive
\index{group ! permutation ! highly transitive}%
 is contained in a maximal subgroup of the symmetric
group. As a reviewer of the paper said, ``Somewhat surprisingly, the proof 
is not entirely combinatorial, but also involves a little model theory,
\index{model theory}%
 in particular, an appeal to the Cherlin--Mills--Zil'ber theorem
\index{Cherlin--Mills--Zil'ber Theorem}%
 on $\aleph_0$-categorical strictly minimal sets.''
\index{strictly minimal set}%
  The other ingredients are filters and topologies, and indeed their main task is to show that such a
group preserves a non-trivial filter; they deduce this from the fact that
it preserves a non-trivial topology. (A more elementary proof of part of 
this theorem, avoiding the model theory,
\index{model theory}%
 was given in~\cite{cam10}.)

The purpose of this chapter is to carry further the investigation of topologies
and filters derived from relational structures,
\index{relational structure}%
 and their automorphism groups.
\index{group ! automorphism}%
  We give a flavour of some of the ideas that arise when filters are considered on random graph vertex sets rather than on structureless sets, for example what types of filters of random graph vertex subsets admit some of the groups that we have identified in the previous chapter.  

\section{The Random Graph and its Relations}
\index{graph ! random}%

We will study a filter and a topology defined naturally from $\mathfrak{R}$, generated by the vertex neighbourhoods in $\mathfrak{R}$, whose automorphism groups
\index{group ! automorphism}%
 contain the automorphism group of $\mathfrak{R}$. For such neighbourhood filters, $\mathfrak{R}$ has a `universal' property: any countable
graph whose neighbourhood filter is non-trivial contains $\mathfrak{R}$ as a spanning subgraph.
\index{graph ! spanning}%
  From the injection property of $\mathfrak{R}$ it follows that any finite vertex set of $\mathfrak{R}$ has a common neighbour, and this property characterizes the class of countable graphs containing $\mathfrak{R}$ as a spanning subgraph; see the appendix to the chapter for a proof of this property of $\mathfrak{R}$.  The neighbourhood topology
\index{topology ! neighbourhood}%
 is $T_1$ but not Hausdorff,
\index{separation condition ! $T_1$}%
\index{topology ! Hausdorff}%
\index{topology ! neighbourhood}%
 but appears to have some interesting properties, somewhat reminiscent of those of $\mathbb{Q}$.  We will also use the \emph{random bipartite graph} $\mathfrak{B}$,
 \index{graph ! random ! bipartite}%
  the unique graph with the property that, given a bipartition of a countable set
into two countable parts, if edges between sets of the bipartition
are chosen independently with probability~$1$, then the resulting
graph is isomorphic to $\mathfrak{B}$ with probability~$1$. 

Recall that S. Thomas \cite{thomas}
\index{Thomas, S.}%
 found all reducts of $\mathfrak{R}$: that is, all
closed subgroups of the symmetric group on the vertex set $V(\mathfrak{R})$ of $\mathfrak{R}$
which contain the automorphism group
\index{group ! automorphism}%
 of $\mathfrak{R}$. We will refer to his
result later. By contrast, the automorphism groups
\index{group ! automorphism}%
 of the filter and topology defined here are not reducts; they are highly transitive,
so their closures are the symmetric group $\Sym(V(\mathfrak{R}))$.

\section{Neighbourhood Filters}

A \emph{filter}
\index{filter}%
 on a set is a family $ \mathscr{F}$ of subsets of $V(\mathfrak{R})$ satisfying
\begin{itemize}
\item $X,Y\in\mathscr{F}$ implies $X\cap Y\in \mathscr{F}$;
\item $X\in \mathscr{F}$, $Y\supseteq X$ implies $Y\in\mathscr{F}$.
\end{itemize}

A filter $\mathscr{F}$ is \emph{trivial}
\index{filter ! trivial}%
 if it consists of all subsets of $V(\mathfrak{R})$; it is
\emph{principal}
\index{filter ! principal}%
 if it consists of all sets containing a fixed subset $A$ of
$\mathscr{F}$; and it is an \emph{ultrafilter}
\index{filter ! ultrafilter}%
\index{ultrafilter}%
 if, for any $X\subseteq V(\mathfrak{R})$, just one of $X$ and $V(\mathfrak{R})\setminus X$ belongs to $\mathscr{F}$. Ultrafilters are just maximal non-trivial filters; the axiom of choice
\index{axiom of choice}%
 implies that every non-trivial filter is contained in an ultrafilter.

Given a family $\mathcal{A}$ of subsets of $V(\mathfrak{R})$, the \emph{filter generated
by $\mathcal{A}$} is the set
\[\mathscr{F}=\{X\subseteq V: (\exists A_1, \ldots, A_n \in \mathcal{A})
(A_1\cap \cdots \cap A_n) \subseteq X\}.\]
Two families $\mathcal{A}_1$ and $\mathcal{A}_2$ generate the same filter if
and only if each member in $\mathcal{A}_2$ lies in the filter generated
by $\mathcal{A}_1$ (that is, contains a finite intersection of sets
of $\mathcal{A}_1$) and \emph{vice versa}.

\medskip

If $\mathscr{F}$ is a filter, then $\mathscr{F}\cup\{\emptyset\}$ is a topology.
\index{topology}%
  Since $V(\mathfrak{R})$ is countable, no non-discrete metric on $V(\mathfrak{R})$ can be complete.  If the topology is not discrete then the set of complements of finite unions of discrete sets is a filter. 

 In this chapter we will assume that the set $\Omega$ is countably infinite.  We claim that every nonprincipal ultrafilter $\mathscr{F}$
\index{ultrafilter ! nonprincipal}%
contains all cofinite sets in $\Omega$.  For a contradiction assume that $A$ is cofinite and $A \notin \mathscr{F}$ then $B = \Omega \backslash A \in \mathscr{F}$ and $B$ is finite.  Choose $B_0 \subseteq B$ minimal with respect to lying in $\mathscr{F}$ then $B_0 \subseteq X$, $\forall X \in \mathscr{F}$.  Suppose $|B_0| > 1$.  Then $B_0 \backslash \{x\} \notin \mathscr{F}$.  So $(\Omega \backslash \{B_0\}) \cup \{x\} \in \mathscr{F}$.  So $B_0 \cap (\Omega \backslash \{B_0\} \cup \{x\}) = \{x\} \in \mathscr{F}$.  So $\mathscr{F}$ is principal, and the claim is proved.

The filter of all cofinite subsets~\cite{cam10} is called the \emph{Fr\'echet filter}~\cite{jech}
\index{filter ! Fr\'echet}%
\index{Fr\'echet, M.}%
over an infinite set.  It is non-principal, and any ultrafilter that extends a Fr\'echet filter is also non-principal.

As an example of a filter supporting a specific group action,
\index{group ! action}%
 we show that if $\mathscr{F}$ is a nonprincipal ultrafilter then it admits $\FSym(\Omega)$.  Let $\mathscr{F}$ be such a filter.  If $A \in \mathscr{F}$ and $g \in \FSym(\Omega)$ then $|A \setminus A^g| < \omega$.  So $A \cap A^g$ is a cofinite subset of $A$, and we can write this as $A \cap B$ for some cofinite set $B$.  So $B \in \mathscr{F}$, so $A \cap B \in \mathscr{F}$ so $A^g \in \mathscr{F}$.

\bigskip

Let $\Gamma$ be a graph on a countable vertex set $V$. We define the
\emph{neighbourhood filter $\mathscr{F}_{\Gamma}$ of $\Gamma$}
\index{filter ! neighbourhood}%
 to be the filter generated by $\{\Gamma(v):v\in V\}$, where $\Gamma(v)$ denotes the neighbourhood of $v$ in $\Gamma$, the set of vertices adjacent to $v$.

\begin{proposition}
Suppose that $\Gamma$ has the property that each vertex has a non-neighbour.
Then the filter generated by the closed neighbourhoods
$\overline{\Gamma}(v)=\Gamma(v)\cup\{v\}$ is equal to $\mathscr{F}_{\Gamma}$.
\end{proposition}

\begin{proof} For we have $\Gamma(v) \subseteq \overline{\Gamma}(v)$, and, if $w$ is not adjacent to $v$, then $\overline{\Gamma}(v)\cap\overline{\Gamma}(w) \subseteq \Gamma(v)$.
\end{proof}

The condition on $\Gamma$ is necessary. If $\Gamma$ is the complete graph,
the closed neighbourhoods generate the filter on $\{V\}$,
\index{filter ! Fr\'echet}%
 while the open neighbourhoods generate the Fr\'echet filter of cofinite subsets of~$V$.
\index{filter ! Fr\'echet}%

Let $\mathfrak{R}$ denote the countable random graph.
\index{graph ! random}%

\begin{proposition}
\label{prop8.2}
The following three conditions on a graph $\Gamma$ are equivalent:
\begin{itemize}
\item[(a)] $\mathscr{F}_\Gamma$ is nontrivial;
\item[(b)] $\Gamma$ contains $\mathfrak{R}$ as a spanning subgraph;
\index{graph ! spanning}%
\item[(c)] $\mathscr{F}_\Gamma\subseteq\mathscr{F}_{\mathfrak{R}}$.
\end{itemize}
\end{proposition}

\begin{proof} A filter is trivial
\index{filter ! trivial}%
 if and only if it contains the empty
set. So $\mathscr{F}_\Gamma$ is non-trivial if and only if any finite number of
neighbourhoods have non-empty intersection. This is equivalent to the
statement that $\mathfrak{R}$ is a spanning subgraph of $\Gamma$; see the appendix to this
chapter.

So (a) and (b) are equivalent.

If $\Gamma$ contains $\mathfrak{R}$ as a spanning subgraph, then
$\mathfrak{R}(v)\subseteq\Gamma(v)$ for all $v$. So (b) implies (c). Conversely, $\mathscr{F}_{\mathfrak{R}}$
is non-trivial (by our proof that (b) implies (a)), so (c) implies (a).
\end{proof}

\head{Remarks}
This result shows that $\mathscr{F}_{\mathfrak{R}}$ is the unique maximal neighbourhood
filter. But this uniqueness is only up to isomorphism. So part (c) really
means that $\mathscr{F}_\Gamma$ is contained in a filter isomorphic to $\mathscr{F}_{\mathfrak{R}}$.

For example, it is possible to find two filters isomorphic to $\mathscr{F}_{\mathfrak{R}}$, one
contained in the other. For let $T$ be the random
$3$-colouring of the edges of the complete graph, with colours red, green and
blue. Let $\mathfrak{R}_1$ be the graph consisting of red edges, and $\mathfrak{R}_2$ the graph
consisting of red and green edges, in $T$. Then both $\mathfrak{R}_1$ and $\mathfrak{R}_2$ are
isomorphic to $\mathfrak{R}$. Since $\mathfrak{R}_1(v)\subseteq \mathfrak{R}_2(v)$, we have
$\mathscr{F}_{\mathfrak{R}_2}\subseteq\mathscr{F}_{\mathfrak{R}_1}$. We show that the inequality is strict.

The set $\mathfrak{R}_1(v)$ belongs to $\mathscr{F}_{\mathfrak{R}_1}$. Suppose that it belongs to
$\mathscr{F}_{\mathfrak{R}_2}$. Then there are vertices $w_1,\ldots,w_n$ such that
\[\mathfrak{R}_2(w_1)\cap \ldots\cap \mathfrak{R}_2(w_n)\subseteq \mathfrak{R}_1(v).\]
But, since the green graph is isomorphic to $\mathfrak{R}$, there is a vertex $x$
joined to all of $v,w_1,\ldots,  w_n$ by green edges; then $x$ belongs to the
left-hand expression of the displayed inclusion but not to $\mathfrak{R}_1(v)$, a
contradiction.

Similarly it can be shown that there are countable chains of filters
isomorphic to $\mathscr{F}_{\mathfrak{R}}$.

\section{Automorphism Groups of Neighbourhood Filters}
\index{topology ! neighbourhood}%

Clearly $\Aut(\mathfrak{R})$ is a subgroup of $\Aut(\mathscr{F}_{\mathfrak{R}})$. We will see in this section
that $\Aut(\mathscr{F}_{\mathfrak{R}})$ is much larger than $\Aut(\mathfrak{R})$.

For convenience let us recall the definition of a few groups. We say that a permutation $g$ changes the
adjacency of $v$ and $w$ if $(v\sim w)\Leftrightarrow(v^g\not\sim w^g)$. We
say that $g$ changes finitely many adjacencies at $v$ if there are only
finitely many points $w$ for which $g$ changes the adjacency of $v$ and $w$.

Let $C(g)$ be the set of pairs $\{v, w\}$ whose adjacency is changed by $g$.  Then $C(g^{-1}) = C(g)^{g^{-1}}$ and $C(gh) \subseteq C(g) \cup C(h)^{g^{-1}}$.

\begin{itemize}
\item $\DAut(\mathfrak{R})$ is the group of anti-automorphisms and automorphisms of $\mathfrak{R}$;
\index{graph ! anti-automorphism}%
\item $\SAut(\mathfrak{R})$ is the group of switching-automorphisms of $\mathfrak{R}$;
\item $\AAut(\mathfrak{R})$ is the group of permutations which change only finitely many adjacencies, (\emph{almost-automorphisms});
\item $\Aut_0(\mathfrak{R})$
\index{group ! zero vertex index}%
 is the group of permutations changing only finitely many adjacencies at any vertex (see Chapter~\ref{chap7});
\item $\Aut_{fin}(\mathfrak{R})$
\index{group ! finite vertex index}%
 is the group of permutations which change infinitely many adjacencies at only finitely many vertices (see Chapter~\ref{chap7});
\item $\FSym(V)$ and $\Sym(V)$ are the finitary symmetric group and the full
symmetric group on the set $V$.
\end{itemize}

By inspection, all these sets of permutations can truly be seen to form groups, as claimed.  (For $\AAut(\mathfrak{R})$,  $\Aut_0(\mathfrak{R})$ and $\Aut_{fin}(\mathfrak{R})$, use the above facts about $C(g)$).
The groups $\Aut_0(\mathfrak{R})$, $\Aut_{fin}(\mathfrak{R})$ and $\Aut(\mathscr{F}_{\mathfrak{R}})$ are not closed.  They are all highly transitive, so their closures are $\Sym(\mathfrak{R})$.

\begin{proposition}
\item{(a)} $\Aut_0(\mathfrak{R})\le\Aut(\mathscr{F}_{\mathfrak{R}})$;
\item{(b)} $\Aut_{fin}(\mathfrak{R})\not\le\Aut(\mathscr{F}_{\mathfrak{R}})$ and $\Aut(\mathscr{F}_{\mathfrak{R}})\not\le\Aut_{fin}(\mathfrak{R})$;
\item{(c)} $\FSym(V) \le \Aut_{fin}(\mathfrak{R})$, and
$\FSym(V)\cap\Aut_0(\mathfrak{R})=\{1\}$.
\item{(d)} $\FSym(V)\le\Aut(\mathscr{F}_{\mathfrak{R}})$.
\label{p1}
\end{proposition}

\begin{proof}
(a) Let $g \in \Aut_0(\mathfrak{R})$. It suffices to show that, for any vertex
$v$, we have $\mathfrak{R}(v)^g\in\mathscr{F}_{\mathfrak{R}}$. Now by assumption, $\mathfrak{R}(v)^g$ differs
only finitely from $\mathfrak{R}(v^g)$; let $\mathfrak{R}(v)^g\setminus \mathfrak{R}(v^g) = \{x_1,
\ldots, x_n\}$. Then we have
\[\mathfrak{R}(v^g)\cap \mathfrak{R}(x_1)\cap\cdots\cap \mathfrak{R}(x_n) \subseteq \mathfrak{R}(v)^g,\]
and we are done.

(b) Choose a vertex $v$, and  consider the graph $\mathfrak{R}'$ obtained by
changing all adjacencies at $v$. Then $\mathfrak{R}'\cong \mathfrak{R}$. Choose an
isomorphism $g$ from $\mathfrak{R}$ to $\mathfrak{R}'$; since $\mathfrak{R}'$ is vertex-transitive,
we can assume that $g$ fixes $v$. So $g$ maps $\mathfrak{R}(v)$ to
$\mathfrak{R}_1(v)=V(\mathfrak{R})\setminus(\{v\} \cup \mathfrak{R}(v))$. Clearly $g \in \Aut_{fin}(\mathfrak{R})$, since it changes only one adjacency at any point different from $v$.
But if $g\in\Aut(\mathscr{F}_{\mathfrak{R}})$, then we would have $\mathfrak{R}_1(v)\in \mathscr{F}_{\mathfrak{R}}$, a contradiction since $\mathfrak{R}(v)\cap \mathfrak{R}_1(v)=\emptyset$.

In the reverse direction, let $\mathfrak{R}''$ be the graph obtained by changing all adjacencies between non-neighbours of $v$.  Again $\mathfrak{R}'' \cong \mathfrak{R}$, and we can pick an isomorphism from $\mathfrak{R}$ to $\mathfrak{R}''$ which fixes $v$.  Now $g$ changes infinitely many adjacencies at all non-neighbours of $v$ (and none at $v$ or its neighbours).  Also, if $w$ is a non-neighbour of $v$, then $\mathfrak{R}(v) \cap \mathfrak{R}(w)^{g} = \mathfrak{R}(v) \cap \mathfrak{R}(w^{g})$, and so $g \in \Aut(\mathscr{F}_{\mathfrak{R}})$.

(c) Any non-identity finitary permutation belongs to
$\Aut_{fin}(\mathfrak{R})\setminus\Aut_0(\mathfrak{R})$. For if $g$ moves $v$, then $g$
changes infinitely many adjacencies at $v$ (namely, all $v$ and $w$,
where $w$ is adjacent to $v$ but not $v^g$ and is not in the support
of $g$). On the other hand, if $g$ fixes $v$, then $v$ changes the
adjacency of $v$ and $w$ only if $g$ moves $w$, and there are only
finitely many such $w$.

(d) If $g$ is a finitary permutation, then $\mathfrak{R}(v)^g$ contains a
cofinite subset of $\mathfrak{R}(v)$, which clearly contains an element of
$\mathscr{F}_{\mathfrak{R}}$.
\end{proof}

Note that $\Aut(\mathscr{F}_{\mathfrak{R}})$ is not equal to either $\Aut_0(\mathfrak{R})$ or
$\Aut_{fin}(\mathfrak{R})$.

The graph $\mathfrak{R}'$ in the proof of Proposition~\ref{p1}(b) is obtained
from $\mathfrak{R}$ by switching with respect to the set $\{v\}$; so the
permutation $g$ belongs to $S(\mathfrak{R})$. Thus $S(\mathfrak{R})\not\le\Aut(\mathscr{F}_{\mathfrak{R}})$. In fact, more is true:

\begin{proposition}
$\Aut(\mathscr{F}_{\mathfrak{R}})\cap \DAut(\mathfrak{R})=\Aut(\mathscr{F}_{\mathfrak{R}})\cap \SAut(\mathfrak{R})= \Aut(\mathfrak{R})$.
\end{proposition}
\index{graph ! anti-automorphism}%
\begin{proof} Any anti-automorphism $g$ of $\mathfrak{R}$ maps $\mathfrak{R}(v)$ to a set
disjoint from $\mathfrak{R}(v^g)$; so no anti-automorphism can belong to $\Aut(\mathscr{F}_{\mathfrak{R}})$.

Suppose that $g\in\Aut(\mathscr{F}_{\mathfrak{R}})$ is an isomorphism from $\mathfrak{R}$ to $\sigma_X(\mathfrak{R})$, where $\sigma_X$ denotes switching with respect to $X$. We may suppose that
$\sigma_X$ is not the identity, that is, $X \ne \emptyset$ and
$Y=V(\mathfrak{R}) \setminus X\ne\emptyset$. Choose $x$ and $y$ so that $x^g\in X$ and
$y^g\in Y$. Then
$\mathfrak{R}(x)^g\symd Y = \mathfrak{R}(x^g)$ and $\mathfrak{R}(y)^g\symd X  = \mathfrak{R}(y^g)$. Hence $\mathfrak{R}(x^g)\cap \mathfrak{R}(x)^g \subseteq X$ and $\mathfrak{R}(y^g) \cap \mathfrak{R}(y^g) \subseteq Y$. Hence
\[\mathfrak{R}(x^g)\cap \mathfrak{R}(x)^g \cap \mathfrak{R}(y^g)\cap \mathfrak{R}(y)^g = \emptyset,\]
a contradiction.
\end{proof}

\section{Fixed Points on Ultrafilters}
\index{filter ! ultrafilter}%
Let $\mathscr{F}$ be a filter on $\Omega$. Define $Aut(\mathscr{F})$ to be the group of automorphisms of $\mathscr{F}$ (permutations $g$ of $\Omega$ such that $g$ and $g^{-1}$ map $\mathscr{F}$ to itself), and $B(\mathscr{F})$ the set of permutations $g$ such that $\fix(g)\in\mathscr{F}$, where $\fix(g)$ is the set of fixed points of $g$.

\begin{proposition}
\begin{itemize}
\item[(a)] $B(\mathscr{F})$ is a normal subgroup of $Aut(\mathscr{F})$.
\item[(b)] If $\mathscr{F}$ is an ultrafilter then $B(\mathscr{F})=Aut(\mathscr{F})$.
\end{itemize}
\end{proposition}

\head{Proof} First we show that $B(\mathscr{F})$ is a group. Clearly it is closed under inversion. If $g_1,g_2\in B(\mathscr{F})$, then $\fix(g_1),\fix(g_2)\in \mathscr{F}$, and $\fix(g_1g_2)\supseteq \fix(g_1)\cap\fix(g_2)$.

Next we show that $B(\mathscr{F}) \subseteq Aut(\mathscr{F})$. Take $g\in B(\mathscr{F})$ and $Y\in\mathscr{F}$. Then $Yg\supseteq Y\cap\fix(g)$, so $Yg\in\mathscr{F}$.

Finally, suppose that $\mathscr{F}$ is an ultrafilter, and take $g\in Aut(\mathscr{F})$. Write $X=A\cup B\cup C\cup D$, where the four sets are chosen as follows. First, $A=\fix(g)$. Now, choose a point in each non-trivial cycle of $g$, and let $B$ consist of alternate points of the cycle beginning at the chosen point. (If the cycle has odd length $2n+1$, then we put just $n$ points into $B$.) Then $C=Bg$, and $D$ consists of the remaining points (one from each odd cycle of $g$). Now one of the sets $A,B,C,D$ must belong to $\mathscr{F}$, since it is an ultrafilter. But $Bg=C$, $Cg\subseteq B\cup D$, and $Dg\subseteq B$; so none of $B,C,D$ can belong to $\mathscr{F}$. Thus, $A\in\mathscr{F}$, showing that $g\in B(\mathscr{F})$, as required.

\begin{corollary}
A necessary condition for a  permutation group to fix an ultrafilter is that every element of the group has a fixed point.
\end{corollary}

This condition is not sufficient. If $K$ is an algebraically closed field, the group $\mathrm{PGL}(n,K)$,
\index{group ! projective linear}%
 acting on the $(n-1)$-dimensional projective space, has the property that every group element fixes a point (since every matrix has an eigenvector). But some elements fix only one point, and clearly the group preserves no principal filter. Note that, if $K=\mathbb{C}$, the projective space is compact. Note also that the automorphism group
\index{group ! automorphism}%
 of any non-principal ultrafilter contains the group of finitary permutations, and so is highly transitive. Also, such a group is a maximal subgroup of the full symmetric group.

\section{Sierpi\'nski's Theorem}

The group $H(\mathbb{Q})$ of autohomeomorphisms~\label{autohomeo}
\index{group ! autohomeomorphisms of $\mathbb{Q}$}%
of the rationals regarded as a topological group,
\index{topological group}%
 is both highly transitive and contains no non-trivial finitary permutations~\cite{neumann1}.
\index{group ! permutation ! highly transitive}%
Mekler
\index{Mekler, A.}%
\emph{et al.}~\cite{mek} proved that $\AAut(\mathfrak{R})$
cannot be embedded into $H(\mathbb{Q})$.
\index{group ! homeomorphisms of $\mathbb{Q}$}%

Truss
\index{Truss, J. K.}%
 proved~\cite{truss0} that the group $\Aut(\mathfrak{R})$ can be, and we outline his proof below.  He used the following idea: a \emph{permutation embedding} of a permutation group $G_{1}$ acting on a set $\Omega_{1}$ into a permutation group $G_{2}$ acting on a set $\Omega_{2}$ is a bijection from $\Omega_{1}$ onto $\Omega_{2}$ which induces a group monomorphism of $G_{1}$ into $G_{2}$.  
  
 Mekler gave~\cite{mekler} a necessary and sufficient condition for a countable permutation group to be embedded in $H(\mathbb{Q})$.  The 
reformulation of his criterion by $\Pi$eter Neumann~\cite{neumann1}
\index{Neumann, $\Pi$. M.}%
is: if $g_1, \ldots, g_n$ 
are finitely many members of a group $G$, then $\bigcap_{i = 1}^{n} \supp(g_i)$ is empty or infinite.  Certainly if $G \le 
H(\mathbb{Q})$ then the support of any homeomorphism is an open set
\index{open set}%
 in $\mathbb{Q}$ and any open set
\index{open set}%
  is empty or infinite.  Mekler's criterion,
\index{Mekler's criterion}%
which is similar to Truss's Lemma~\ref{infsuplem} in Appendix~\ref{PreviousResults}, implies that every countable subset of $\Aut(\mathfrak{R})$ embeds in $H(\mathbb{Q})$.  The group $\Aut(\mathfrak{R})$ is uncountable and Mekler's Theorem does not apply.

We will give $\Pi$eter Neumann's
\index{Neumann, $\Pi$. M.}%
 proof of Sierpi\'nski's characterisation of $\mathbb{Q}$~\cite{neumann1}, which we need in both the construction of the embedding which follows, and also in the next section.  A different proof is given in Sierpi\'nski's book~\cite[\S~59]{sierpinski1}.

We shall need some definitions.  A space is \emph{second countable}
\index{topological space ! second countable}%
if there is a countable base for the
topology; it is \emph{$0$-dimensional}
\index{topological space ! $0$-dimensional}%
if for any $x \in X$ and any open set
\index{open set}%
$U$ containing $x$ there is a clopen set
\index{clopen set}%
 $V$ such that $x \in V
\subseteq U$; $X$ is a \emph{$T_1$ space}
\index{topological space ! $T_1$}%
if singleton sets are closed, that
is if $\{x\} = \bigcap \{ U\ |\ x \in U \mbox{ and } U \mbox{ is closed in
} X \}$ for any $x \in X$; an \emph{isolated point}
\index{isolated point}%
$x$ is one for which the singleton set $\{x\}$ is open.  The topological space $X$ is metrizable.
\index{topological space ! metrizable}%

\begin{theorem}[Sierpi\'nski's Theorem]
\index{Sierpi\'nski's Theorem}%
\index{Sierpi\'nski, Waclaw}%
Let $X$ be a second countable, $0$-dimensional, $T_1$ topological space.
\index{topological space}%
  Then $X$ is homeomorphic with a subspace of $\mathbb{R}$.  If moreover $X$ is countable and has no isolated points then $X$ is homeomorphic with $\mathbb{Q}$.
\end{theorem}
\begin{proof}
The space $X$ has a countable family of subsets $C_1, C_2, \ldots$ such that

(1)  each $C_r$ is clopen;

(2)  if $x \in X$ then $\{x\} =   \bigcap \{ C_r\ |\ x \in C_r \}$;

(3)  the family is a base for the topology of $X$.

For a given $x \in X$ and family $\{C_r\}$, define
\begin{displaymath}
D_r(x) : = \left\{ \begin{array}{ll}
C_r & \text{if}\ x \in C_r\\
X - C_r & \text{if}\ x \notin C_r
 \end{array} \right.
\end{displaymath}

and clopen sets
\index{clopen set}%
 containing $x$, \begin{displaymath}
U_n(x) : =  \begin{array}{ll}
\bigcap_{r = 1}^{n} D_r(x).
 \end{array}
\end{displaymath}

An embedding of $X$ into $\mathbb{R}$ is sought.  Define the continuous map $f_r : X \to \{ 0, 2 \}$ by
\begin{displaymath}
f_r(x) : = \left\{ \begin{array}{ll}
2 & \text{if}\ x \in C_r\\
0 & \text{if}\ x \notin C_r.
 \end{array} \right.
\end{displaymath}

From (2), the function
\begin{displaymath}
f(x) : =  \begin{array}{ll}
\sum_{r = 1}^{\infty} f_r(x) \cdot 3^{-r}
 \end{array}
\end{displaymath}
is a Cantor set-valued $1$--$1$ map of $X$ into $\mathbb{R}$.  

To show that $f$ is continuous, let $m \in \mathbb{N}$ be such that $3^{-m} < \epsilon$ for all $\epsilon \in \mathbb{R}^{>0}$.  If $x \in U_m(x_0)$ then for any $x_0 \in X$, $f_r(x) = f_r(x_0)$ whenever $1 \leq r \leq m$, so that
\[ | f(x) - f(x_0) | \leq \sum_{m+1}^{\infty} 2 \cdot 3^{-r}  = 3^{-m} < \epsilon. \]

To show that $f^{-1} : f(X) \to X$ is continuous, for any $y_0 \in f(X)$ say $y_0 = f(x_0)$, let $U$ be any open subset of $X$ containing $x_0$.  From (3), $\exists n$ such that $x_0 \in C_n \subseteq U$, so that $U_n(x_0) \subseteq U$.  Define $\delta := 3^{-n}$.  If $y \in f(X)$, say $y = f(x)$ for $x \in X$, is such that $| y - y_0 | < \delta$,  and if $m$ is the least natural number such that $f_m(x) \neq f_m(x_0)$ (assuming that $x \neq x_0$), then
\[ | f(x) - f(x_0) | \geq 2 \cdot 3^{-m} - \sum_{m+1}^{\infty} 2 \cdot 3^{-r} = 3^{-m}, \]
implying that $m > n$ and so $x \in U_n(x_0)$.  So if $y \in f(X)$ and $| y - y_0 | < \delta$ then $f^{-1}(y) \in U$ and so $f^{-1}$ is continuous.  This proves that $f$ is a homeomorphism of $X$ to a subset of the Cantor set.  

Now further suppose that $X$ is countable and has no isolated points.  The next step is to verify that the family $\{C_r\}$ can further be chosen to satisfy

(4)  $\forall x \in X\ \forall m \in \mathbb{N}, \exists n \in \mathbb{N}$ such that $n \geq m, x \in C_{n+1}$ and $U_n(x) \nsubseteq C_{n+1}$;

(5)  $\forall x \in X\ \forall m \in \mathbb{N}, \exists n \in \mathbb{N}$ such that $n \geq m, x \notin C_{n+1}$ and $U_n(x) \cap C_{n+1} \neq \emptyset$.

Let $V_1, V_2, V_3, \ldots$ be the countable base for the topology
\index{topology}%
 on $X$, and enumerate the set of all pairs $(x, V_i)$ with $x \in V_i$ as $(x_1, X_1), (x_2, X_2),$ $(x_3, X_3), \ldots$.  Let $x_1 \in C_1 \subset X_1$ and $x_1 \notin C_2 \subset X_1$, where $C_2 \neq \emptyset$ and $C_1, C_2$ are clopen in $X$.  Next suppose that $C_1, C_2, \ldots, C_{2k-1}, C_{2k}$ have already been chosen.  Then $U_{2k+1}(x_{k+1})$ is open, contains $x_{k+1}$ and so contains $y \neq x_{k+1}$.  By (1) and (2) we can take $C_{2k+1}$ to be clopen in $X_{k+1}$, $x_{k+1} \in C_{2k+1}$, and $y \notin C_{2k+1}$.  So $U_{2k}(X_{k+1}) \subseteq C_{2k+1}$.  Now $\exists z \in U_{2k+1}(X_{k+1}) = U_{2k}(X_{k+1}) \cap C_{2k+1}$ where $z \neq x_{k+1}$.  Take $C_{2k+2}$ to be any clopen set
\index{clopen set}%
 containing $z$ but not $x_{k+1}$.

The constructed sequence $\{ C_i \}$ satisfies (1).  Furthermore,
\begin{displaymath}
\begin{array}{ll}
\quad \bigcap \{C_r |\ x \in C_r\}\\
\subseteq \bigcap \{C_r |\ x \in C_r \text{ and } r \text{ is odd}\}\\
\subseteq \bigcap \{X_s |\ x = x_s\}\\
\subseteq \bigcap \{V_i |\ x \in V_i\}\\
= \{x\}
\end{array}
\end{displaymath}
so (2) is satisfied.  If $U$ is any open set
\index{open set}%
 and $x \in U$ then there exists $i$ such that $x \in V_i$ and $V_i \subseteq U$; then $(x, V_i)$ is of the form $(x_n, V_n)$ for some $n$ and so $x_n \in C_{2n-1} \subseteq U$; so (3) is satisfied.

Next, let $x \in X$ and $m \in \mathbb{Z}^{+}$.  Choose $k \geq m$ so that $x = x_k$ (such an integer exists for otherwise $\{i |\ x \in V_i\}$ would be finite and $\{x\} = \bigcap \{V_i |\ x \in V_i\}$ would be open).  By construction, if $n := 2k$ then $x \in C_{n+1}$ and $U_n(x) \subseteq C_{n+1}$, whilst if $n := 2k+1$ then $x \neq C_{n+1}$ and $U_n(x) \cap C_{n+1} \neq \emptyset$.  Thus conditions (4) and (5) are satisfied.

Condition (4) (respectively (5)) shows that every member of $f(X)$ is a left (respectively right) limit in $f(X)$, and so every element of $f(X)$ is a two-sided limit in $f(X)$.  To see this for (4), let $x \in X$, $\epsilon > 0$, $m$ be such that $3^{-m} < \epsilon$, and $n \geq m$ such that $x \in C_{n+1}$ and $U_n(x) \nsubseteq C_{n+1}$, and let $y \in 
U_n(x) \backslash C_{n+1}$.  Since $y \in U_n(x)$ we have that if $1 \leq r \leq n$ then $f_r(y) = f_r(x)$.  But also $f_{n+1} (x) = 2$ and $f_{n+1} (x) = 0$.  So
\[  f(x) - f(y)  \geq 2 \cdot 3^{-(n+1)} - \sum_{n+2}^{\infty} 2 \cdot 3^{-r} = 3^{-(n+1)} > 0, \]
so that $f(x) > f(y)$, while at the same time
\[  f(x) - f(y)  \leq  \sum_{n+1}^{\infty} 2 \cdot 3^{-r} = 3^{-n} \leq 3^{-m} < \epsilon .\]
Thus $f(x) \in (f(x) - \epsilon, f(x))$, and since $\epsilon$ was arbitrary, $f(x)$ is a left limit in $f(X)$.

To complete the proof, note that if a set $S$ is a countable subset of $\mathbb{R}$ every point in which is a two-sided limit in $S$, then it is densely-ordered by $\leq$ and has no maximum nor minimum, and so by Cantor's theorem
\index{Cantor's theorem}%
 is order-isomorphic to $\mathbb{Q}$.  That every point is a two-sided limit in $S$ means that the order-topology
\index{topology ! order}%
  on $S$ is the same as that induced from $\mathbb{R}$.  So the order-isomorphism of $S$ to $\mathbb{Q}$ is a homomorphism.
\end{proof}

As a corollary Sierpi\'nski
\index{Sierpi\'nski, Waclaw}%
 showed that any countable metric space is homeomorphic to a subset of $\mathbb{Q}$.

There is yet another method for proving Sierpi\'nski's Theorem,
\index{Sierpi\'nski's Theorem}%
one that singles out the essential ingredient as being that the usual topology on $\mathbb{Q}$ has a countable base $\mathbb{B}$~\label{boolalg} of clopen subsets that under the operations of union, intersection and complementation in $\mathbb{Q}$ forms a Boolean algebra
\index{Boolean algebra}%
 all of whose non-empty members are infinite, and for which if $p \neq q$ with $p, q \in \mathbb{Q}$ then $\exists b \in \mathbb{B}$ with $x \in b$ and $y \notin b$.  As an example of such a generating set take the set of intervals of the form $(\pi + a, \pi + b) \cap \mathbb{Q}$ for $a, b \in \mathbb{Q}$.  It is the uniqueness up to isomorphism of this structure that leads to the uniqueness of the rational topology in the following way~\cite{truss0}.

If $\mathbb{B}_1$ is a Boolean algebra generated by a countable base of clopen sets
\index{clopen set}%
 for $X$.  Assuming that $X$ has no isolated points is equivalent to taking all nonempty members of $\mathbb{B}_1$ to be infinite.  If $\mathbb{B}_2$ is the corresponding Boolean algebra for the topology on $\mathbb{Q}$, derived from a countable base for clopen sets, then from~\cite[Theorem~2.1]{truss0} there is an isomorphism from $\mathbb{B}_1$ to $\mathbb{B}_2$ induced by a bijection from $X$ to $\mathbb{Q}$, which is the required homeomorphism.

\begin{theorem}[Truss]
\index{Truss, J. K.}%
The groups $\Aut(\mathfrak{R}_{m, \omega})$ for $2 \leq m \leq \aleph_{0}$ can be embedded in the group $H(\mathbb{Q})$.
\end{theorem}
\begin{proof}
By going colourblind in pairs of colours it is possible to embed $\Aut(\mathfrak{R}_{m_2, \omega})$ into $\Aut(\mathfrak{R}_{m_1, \omega})$ for $m_1 \leq m_2$, so we can restrict ourselves to $\mathfrak{R}$.  Let $\mathbb{B}$ be the Boolean algebra generated by neighbourhood sets~\label{neighb} $N(x) : = \{ z : z \sim x, \mbox{where } x, z \in \mathfrak{R}\}$.  The aim is to show that $\mathbb{B}$ is countable and that its elements are empty or infinite subsets of $\mathfrak{R}$, closed under the action of $\Aut(\mathfrak{R})$, such that for any distinct $x, y \in \Aut(\mathfrak{R})$, $\exists b \in \mathbb{B}$ with $x \in b$ and $y \notin b$.  For then it would follow that the topology with base $\mathbb{B}$ is homeomorphic to that on $\mathbb{Q}$, and therefore that $\Aut(\mathfrak{R})$ embeds into $H(\mathbb{Q})$.

Consider elements of $\mathbb{B}$ of the form 
\[ X = N(x_1) \cap \ldots N(x_k) \cap (\mathfrak{R} \backslash N(y_1)) \cap \ldots \cap (\mathfrak{R} \backslash N(y_l)).\]  
If $X \neq \emptyset$ then $x_i \neq y_j$ for each $i, j$.  If $A \subset \mathfrak{R}$ is any finite subset then by the ($*$)-property $\exists z\notin A$ such that $z \sim x_i$ and $z \nsim y_j$ for each $i, j$.  So $z \in X$ and so $X \backslash A \neq \emptyset$.  Since $A$ was an arbitrary finite set, $X$ must be infinite.  So  all non-empty members of $\mathbb{B}$ are infinite.

Furthermore, $\forall x, y \in  \mathfrak{R}$ and $g \in \Aut( \mathfrak{R}), y \in g N(x) \Leftrightarrow g^{-1} y \in N(x) \Leftrightarrow x \sim g^{-1}y \mbox{ in }  \mathfrak{R} \Leftrightarrow gx \sim y \mbox{ in } \mathfrak{R} \Leftrightarrow y \in N(g x)$, so that $g N(x) = N(g x) \in \mathbb{B}$.  So $\mathbb{B}$ is closed under the action of $\Aut(\mathfrak{R})$.

Finally, if $x \neq y$, it follows from the ($*$)-property there is $z$ joined to $x$ but not to $y$, so $x \in N(z)$ and $y \notin N(z)$.
\end{proof}

The group $\Aut^{*}(\mathfrak{R})$ of automorphisms and \emph{anti-automorphisms}
\index{graph ! anti-automorphism}%
(permutations which change all adjacencies), of $\mathfrak{R}$ is $2$-transitive,
\index{group ! permutation ! $2$-transitive}%
 contains $\Aut(\mathfrak{R})$ as a subgroup of index 2, but is not $3$-transitive because vertex triples containing 0 or 3 edges of $\mathfrak{R}$ are not equivalent to triples containing 1 or 2 edges.  To see that the embedding of $\Aut(\mathfrak{R})$ in $H(\mathbb{Q})$ does not extend to $\Aut^{*}(\mathfrak{R})$, consider a neighbourhood $N(U)$ of a finite set $U$ in which $\{z_1, z_2, z' \in N(U)\ |\ z_1 \nsim z_2, z_1 \nsim z', z_2 \sim z'\}$.  That such three points can be found with this configuration is made possible by repeated application of ($*$).  Then if $g \in \Aut(\mathfrak{R})$ such that $z_1 g = z_2, z' g = z_1, z_2 g = z', z g = z\ \forall z \in N(U) \backslash \{z_1, z_2, z'\}$ we have that the triple $\{u_1 \sim z_1, u_1 \sim z', z_1 \nsim z'\}$ is mapped to $\{u_1 \sim z_1, u_1 \sim z', z_1 \sim z'\}$.  So $g \in H(\mathbb{Q}) \backslash \Aut^{*}(\mathfrak{R})$.  To see that $\Aut^{*}(\mathfrak{R})$ cannot be embedded in $H(\mathbb{Q})$, note that the intersection of a closed set
\index{closed set}%
  (of the form a vertex together with its neighbours) and an open set
\index{open set}%
  (of the form a vertex together with its non-neighbours) is a single point, and as we already noted singleton sets can only support the discrete topology,
\index{topology ! discrete}%
 which is the largest topology containing all subsets (and in particular every point in the topological space) as open sets.

Finally, from~\cite{truss0} we have that $H(\mathbb{Q})$ can neither be embedded
in $\Aut(\mathfrak{R}_{m,\omega})$ nor in $\AAut(\mathfrak{R}_{m,\omega})$ for any $2 \leq m \leq \aleph_{0}$.  

\section{Topologies on Random Graphs}

As we noted at the beginning of the last section, Mekler~\cite{mekler} and Truss~\cite{truss0} gave a necessary and
sufficient condition (called the \emph{mimicking property} or \emph{strong Mekler criterion})
\index{mimicking property}%
\index{Strong Mekler Criterion}%
 for a countable group to be a subgroup of the
group of homeomorphisms of $\mathbb{Q}$: the intersection of the
supports of any finite number of elements should be empty or
infinite. We also noted that the group $\Aut(\mathfrak{R})$ is uncountable so the criterion does not apply, though the intersection property is still obeyed.  However, we will see that $\Aut(\mathfrak{R})$ is a subgroup of $\Aut(\mathbb{Q})$.  The argument involves constructing the topology of $\mathbb{Q}$ from
$\mathfrak{R}$.  We will also reveal another interesting topology.
\index{topology}%

If $\mathscr{F}$ is a filter,
\index{filter}%
 then $\mathscr{F}\cup\{\emptyset\}$ is a topology with the same automorphism group.
\index{group ! automorphism}%
  It is worth noting, however, that there are a couple of
topologies on $\mathfrak{R}$ related to $\mathscr{F}_{\mathfrak{R}}\cup\{\emptyset\}$ which are in some sense more natural.

In the first topology $\mathcal{T}$, a sub-basis for the open sets
\index{open set}%
 consists of the closed neighbourhoods of vertices.
\index{topology ! neighbourhood}%
 Thus the open sets are all unions of sets which are finite intersections of closed neighbourhoods.

The topology $\mathcal{T}$ is not Hausdorff:
\index{topology ! Hausdorff}%
 in fact, any two open sets have non-empty intersection. For it suffices to show this for basic open
sets; and the intersection of two finite intersections of neighbourhoods
is itself a finite intersection of neighbourhoods, and so is non-empty.

However, this topology does satisfy the $T_1$ separation condition:
\index{separation condition ! $T_1$}%
 for, given distinct points $x$ and $y$, there is a vertex $v$ joined to $x$
but not $y$, and so a neighbourhood containing $x$ but not $y$. Hence
all singletons (and so all finite sets) are closed.

We used closed neighbourhoods in the construction of $\mathcal{T}$.
In fact, open neighbourhoods would have given us the same topology,
\index{topology}%
 as we will now see.

Let $\mathfrak{B}$ denote the generic bipartite graph. 
\index{graph ! random ! bipartite}%
Consider the three topologies which have the following as the
points and sub-basic open sets:~\label{Ttop}
\index{sets ! basic open}%
\begin{itemize}
\item{$\mathcal{T}$:} points are vertices of $\mathfrak{R}$, sub-basic open sets are
open vertex neighbourhoods.
\item{$\mathcal{T}^*$:} points are vertices of $\mathfrak{R}$, sub-basic open sets are
closed vertex neighbourhoods.
\item{$\mathcal{T}^\dag$:} points are in one bipartite block in $\mathfrak{B}$, sub-basic open
sets are neighbourhoods of vertices in the other bipartite block.
\end{itemize}

\bigskip

\bigskip

\begin{proposition}
\label{threetops}
\begin{itemize}
\item[(a)] The three topologies defined above are all homeomorphic.
\item[(b)] The homeomorphism group of these topologies is highly transitive.
\end{itemize}
\end{proposition}
\index{topology}%

\begin{proof}
We construct two bipartite graphs $\mathfrak{B}_1$ and $\mathfrak{B}_2$
\index{graph ! bipartite}%
 as follows. The vertex set of each graph is $V(\mathfrak{R})\times\{0,1\}$; vertices $(v,0)$ and $(w,1)$ are adjacent if and only if
\begin{itemize}
\item $v\sim w$ in $\mathfrak{R}$ (for $\mathfrak{B}_1$);
\item $v=w$ or $v\sim w$ in $\mathfrak{R}$ (for $\mathfrak{B}_2$).
\end{itemize}

The characteristic property of $\mathfrak{R}$ shows that both bipartite graphs satisfy the characteristic property of the generic bipartite graph $\mathfrak{B}$; recall that $\mathfrak{B}$ is characterised as a countable bipartite graph by the following property:
\begin{quote}
Given two finite disjoint sets $X,Y$ of vertices in the same bipartite block,
there exists a vertex $z$ in the other block such that $z\sim x$ for all
$x\in X$, and $z\not\sim y$ for all $y\in Y$.
\end{quote}
It is clear from the property of $\mathfrak{R}$ that both $\mathfrak{B}_1$ and $\mathfrak{B}_2$ satisfy this condition, and so $\mathfrak{B}_1\cong \mathfrak{B}_2\cong \mathfrak{B}$.

It follows immediately that the three topologies are homeomorphic. Moreover,
the stabilizer of a bipartite block in $\mathfrak{B}$ acts on it as a group of
homeomorphisms of the topology $\mathcal{T}^\dag$, and this group is highly
transitive. (This follows from the homogeneity of $\mathfrak{B}$ as a graph with
bipartition: any two vertices in the same bipartite block have distance~$2$,
so any bijection between finite subsets of a bipartite block extends to an
automorphism of~$\mathfrak{B}$).
\end{proof}

\head{Remark} (a)  The topologies
\index{topology}%
 $\mathcal{T}$ and $\mathcal{T}^*$, though
homeomorphic, are not identical. Indeed, the identity map is a continuous
bijection from $\mathcal{T}^*$ to $\mathcal{T}$ but not a homeomorphism.

To see this, note first that, since the topology $\mathcal{T}^*$ is $T_1$,
\index{separation condition ! $T_1$}%
every singleton is closed, and so $\mathfrak{R}(v)=(\mathfrak{R}(v)\cup\{v\})\setminus\{v\}$ is
open in $\mathcal{T}^*$. It follows that any open set
\index{open set}%
 in $\mathcal{T}$ is also open in $\mathcal{T}^*$.

In the other direction, suppose that $\mathfrak{R}(v)\cup\{v\}$ is open in
$\mathcal{T}$. Then it is a union of basic open sets.
\index{sets ! basic open}%
 We can take one of these sets to be $\mathfrak{R}(v)$; let the other be $\bigcap_{x\in X}\mathfrak{R}(x)$ for some
finite set $X$. Then $v$ is joined to all vertices in $X$, but the remaining
common neighbours of these vertices are all in $\mathfrak{R}(v)$. So no point is joined
to all vertices in $X$ but not to $v$, a contradiction.

(b)  In terms of graph vertices and neighbourhoods, three of the separation axioms in the appendix, read as follows:-
\begin{itemize}
\item $T_0$ is $v \notin \Gamma(w) \wedge w \notin \Gamma(v)$;
\item $R_0$ is $v \notin \overline{\Gamma(w)} \wedge w \notin \overline{\Gamma(v)}$ but we can have $\Gamma(w) \subsetneq \Gamma(v)$ or $\Gamma(v) \subsetneq \Gamma(w)$;
\item $R_1$ is $v \notin \overline{\Gamma(w)} \wedge w \notin \overline{\Gamma(v)}$ $\&$ $\Gamma(v) \cap \Gamma(w) = \emptyset$;
\end{itemize}
So $\mathcal{T}$ is $T_0$ and $R_0$ but not $R_1$, whilst both $R_0$ and $R_1$ imply that $v \neq w$ and so $\mathfrak{B}_2$ is only $T_0$ and so whilst $\mathfrak{B}_1\cong \mathfrak{B}_2$, they are not identical.

\bigskip

The second topology $\mathcal{U}$~\label{U} is obtained by symmetrising this one
with respect to the graph $\mathfrak{R}$ and its complement $\mathfrak{R}^c$; in other words,
we also take closed neighbourhoods
\index{topology ! neighbourhood}%
 in $\mathfrak{R}^c$ to be open sets.
\index{open set}%
 So a basis for
the open sets consists of all sets of the form
\[Z(U,V)=\{z\in V(R): (\forall u\in U)(z\sim u) \wedge (\forall v\in V)(z \not \sim v)\}\]
for finite disjoint sets $U$ and $V$. Again it holds that all the non-empty
open sets are infinite. This time the topology is totally disconnected.
\index{topology ! totally disconnected}%
 For given $u\ne v$, there is a point $z\in Z(\{u\},\{v\})$; then the open neighbourhood of $z$ is open and closed in the topology and contains $u$ but not $v$.

By Sierpi\'nski's Theorem,
\index{Sierpi\'nski's Theorem}%
(see also \cite{neumann1}), this topology is homeomorphic to $\mathbb{Q}$.  So $\mathfrak{R}$ as a countable topological space is homeomorphic to $\mathbb{Q}$.

\begin{theorem}
Let $\mathcal{U}$ be a countable, second countable, totally disconnected
topological space
\index{topological group ! totally disconnected}%
with no isolated points. Then $\mathcal{U}$ is homeomorphic
to the usual topology on $\mathbb{Q}$.
\end{theorem}

There are several open questions about the topology $\mathcal{T}$.

\begin{itemize}
\item  Since $\mathcal{T}$ is a coarsening of $\mathcal{U}$, there must be an
identification of it with $\mathbb{Q}$ such that the open sets
\index{open set}%
 in $\mathcal{T}$ are
open in $\mathbb{Q}$. Can such an identification be found explicitly?
\item  Is there a characterisation of $\mathcal{T}$, along the lines of
Sierpi\'nski's Theorem?
\item  The homeomorphism group $\Aut(\mathcal{T}^\dag)$ has as a subgroup the group $\Aut'(\mathfrak{B})$ induced on a bipartite block of $\mathfrak{B}$ by its setwise stabilizer in the
automorphism group
\index{group ! automorphism}%
 of $\mathfrak{B}$. This group is highly transitive. Is it equal
to $\Aut(T^\dag)$ or not?
\item  How are the homeomorphisms of $\mathcal{T}$ and $\mathcal{T}^*$ related to the groups of Section 3?  Since $\mathscr{F}_{\mathfrak{R}}$ consists of all sets containing a non-empty $\mathcal{T}$-open set, we see that $\Aut(\mathcal{T}) \le \Aut(\mathscr{F}_{\mathfrak{R}})$; do any further relations hold?
\end{itemize}

We cannot answer these questions, but we present here a programme which might
lead to an affirmative answer for the third question (which would have implications for the second as well).

Call an open set $U$ \emph{full}
\index{open set ! full}%
 if, for all $x\notin  U$, the set $U\cup\{x\}$ is not open. By the argument used previously to show that the
topologies
\index{topology}%
 $\mathcal{T}$ and $\mathcal{T}^*$ are different, we see that any positive Boolean
combination of neighbourhoods in $\mathcal{T}^\dag$ (that is, any finite union of
basic open sets) is full.
\index{sets ! basic open}%
\begin{quote}
Is it true that these are the only full open sets in $\mathcal{T}^\dag$?
\end{quote}
If so, then we can recognise the basic open sets as the full open sets which
are not proper finite unions of full open sets; and then the neighbourhoods
are the basic open sets which are maximal under inclusion. So we can recover
the graph $\mathfrak{B}$ from the topology, and every homeomorphism is a graph
automorphism.

\section{The $\mathfrak{R}$-uniform Hypergraph $\mathfrak{RHyp}$}

Results of Claude Laflamme, Norbert Sauer and Maurice Pouzet
\index{Laflamme, C.}%
\index{Sauer, N.}%
\index{Pouzet, M.}%
 in~\cite{laflammea} communicated to us via Robert Woodrow
\index{Woodrow, R. E.}%
 concern the hypergraph $\mathfrak{RHyp}$ on the vertex set of $\mathfrak{R}$ whose edges are those sets of vertices which induce a copy of $\mathfrak{R}$.  Note that a cofinite subset of an edge is an edge.  
 
 In this section we will connect groups related to this hypergraph to the groups that have arisen in this and the previous chapter.
 
There are two interesting groups here; we shall call these the LPS groups.

\begin{itemize}
\item $\Aut(\mathfrak{RHyp})$;
\item $\Aut^{*}(\mathfrak{RHyp})$, the group of permutations $g$ with the property that, for every edge $E$, there is a finite subset $S \subset E$ such that $(E \backslash S) g$ and $(E \backslash S) g^{-1}$ are edges.  This is the group of almost automorphisms for $\mathfrak{RHyp}$.
\end{itemize}

In the definition of $\Aut^{*}(\mathfrak{RHyp})$ both conditions are necessary, for choose an infinite clique $C \subset \mathfrak{R}$, and partition $\mathfrak{R}$ into two edges $A$ and $B$.  Then it is shown in~\cite{laflammea} that there exists $g \in \Sym(\mathfrak{R})$ such that $(E)g$ is an edge for any any edge $E$, and such that $C g = A$.  But clearly $(A \backslash S)g^{-1}$ is not an edge for any (finite) $S$.

\begin{proposition}
\begin{itemize}
\item[(a)] $\Aut(\mathfrak{RHyp}) < \Aut^{*}(\mathfrak{RHyp})$.
\item[(b)] $\Aut_{0}(\mathfrak{R}) \leq \Aut(\mathfrak{RHyp})$ and $\Aut_{fin}(\mathfrak{R}) \leq \Aut^{*}(\mathfrak{RHyp})$.
\item[(c)] $\FSym(\mathfrak{R}) \leq \Aut^{*}(\mathfrak{RHyp})$ but $\FSym(\mathfrak{R}) \cap \Aut(\mathfrak{RHyp}) = 1$.
\end{itemize}
\end{proposition}

\begin{proof}
(a)  follows from (c).

(b)  If we alter a finite number of adjacencies at any point of $\mathfrak{R}$, the result is still isomorphic to $\mathfrak{R}$.  So induced copies of $\mathfrak{R}$ are preserved by $\Aut_{0}(\mathfrak{R})$.  Similarly, given an element of $\Aut^{*}(\mathfrak{RHyp})$, if we throw away the vertices where infinitely many adjacencies are changed, we have the case of  $\Aut(\mathfrak{RHyp})$.

(c)  If $g$ is finitary and $S$ is its support, then $g$ restricted to $E \backslash S$ is the identity for any edge $E$.

Choose any vertex $v$ and let $E$ be its set of neighbours in $\mathfrak{R}$ (this set is an edge of $\Aut(\mathfrak{RHyp})$).  Now, for any finitary permutation, there is a conjugate of it whose support contains $v$ and is contained in $\{ v \} \cup E$.  Then $Eg = E \cup \{v\} \backslash \{w\}$ for some $\{w\}$.  But the induced subgraph on this set is not isomorphic to $\mathfrak{R}$, since $v$ is joined to all other vertices.
\end{proof}

\medskip

It is obvious that if $G$ is any subgroup of $\Sym(\mathfrak{R})$ not containing $\FSym(\mathfrak{R})$, then $G . \FSym(\mathfrak{R})$ is a subgroup of $\Sym(\mathfrak{R})$ containing $G$, and is highly transitive.

\begin{lemma}  Any overgroup of $\Aut(\mathfrak{R})$ which is not contained in $\BAut(\mathfrak{R})$ is highly transitive.
\end{lemma}

\begin{proof}
Let $G \geq \Aut(\mathfrak{R})$, and suppose that $G$ is not $n$-transitive.  Let $O$  be an orbit of $G$ on $n$-tuples containing an $n$-clique.  Since the setwise stabilizer of an $n$-clique in $\Aut(\mathfrak{R})$ induces the symmetric group on it, we see that $O$ is symmetric, and so the set $O^{*}$ of $n$-sets supporting members of $O$ is not trivial.  But since $G \geq \Aut(\mathfrak{R})$ and $\mathfrak{R}$ is homogeneous, $O^{*}$ is specified by the set of isomorphism types of subgraphs induced; so $\Aut(O)$ is a reduct of $\Aut(\mathfrak{R})$.  Then $G \leq \Aut(O) \leq \BAut(\mathfrak{R})$ by Thomas' theorem.
\index{Thomas' Theorem}%
\end{proof}

The next result gives some information about how the LPS groups are connected to the reducts.

\begin{proposition}
\begin{itemize}
\item[(a)] $\DAut(\mathfrak{R}) < \Aut(\mathfrak{RHyp})$.
\item[(b)] $\SAut(\mathfrak{R}) \nleq \Aut(\mathfrak{RHyp})$.
\end{itemize}
\end{proposition}

\begin{proof}
(a)  Clear since $\mathfrak{R}$ is self-complementary.

(b)  We show that $\mathfrak{R}$ can be switched into a graph isomorphic to $\mathfrak{R}$ in such a way that some induced copy of $\mathfrak{R}$ has an isolated vertex after switching.  Then the isomorphism is a switching-automorphism but not an automorphism of $\mathfrak{RHyp}$.

Let $p, q$ be two vertices of $\mathfrak{R}$.  We shall work with the graph $\mathfrak{R}_1 = \mathfrak{R} \backslash \{ p \} \cong \mathfrak{R}$.  Let $A, B, C, D$ be the sets of vertices joined respectively to $p$ and $q$, $p$ but not $q$, $q$ but not $p$, and neither $p$ nor $q$.  Let $\sigma$ be the switching operation of $\mathfrak{R}_1$ with respect to $C$, and let $E = \{q\} \cup B \cup C$.  It is clear that, after the switching $\sigma$, the vertex $q$ is isolated in $E$.  So we must prove two claims:

\head{Claim 1:} $E$ induces a copy $\mathfrak{R}$.

For take $U, V$ to be finite disjoint subsets of $E$.  We may assume without loss of generality that $q \in U \cup V$.

\head{\quad Case 1:} $q \in U$.  Choose a witness for $(U, V \cup \{ p \})$ in $\mathfrak{R}$.  Then $z \notin p$ and $z \sim q$, so $z \in C$; thus $z$ is a witness for $(U, V)$ in $E$.

\head{\quad Case 2:} $q \in V$.  Now choose a witness for  $(U \cup \{ p \}, V)$ in $\mathfrak{R}$;  the argument is similar.

\head{Claim 2:} $\sigma(\mathfrak{R}_1) \cong \mathfrak{R}$.

Choose $U, V$ finite disjoint subsets of $\mathfrak{R} \backslash \{ p \}$.  Again, without loss, $q \in U \cup V$.  Set $U_1 = U \cap C, U_2 = U \backslash U_1$, and $V_1 = V \cap C, V_2 = V \backslash V_1$.

\head{\quad Case 1:} $q \in U$, so $q \in U_2$.  Take $z$ to be a witness for  $(U_2 \cup V_1 \cup \{ p \}, U_1 \cup V_2)$ in $\mathfrak{R}$.  Then $z \sim p, q$ so $z \in A$.  The switching $\sigma$ changes its adjacencies to $U_1$ and $V_1$, so in $\sigma(\mathfrak{R}_1)$ it is a witness for $(U_1 \cup U_2, V_1 \cup V_2)$.

\head{\quad Case 2:}  $q \in V$, so $q \in V_2$.  Now take $z$ to be a witness for $(U_1 \cup V_2, U_2 \cup V_1 \cup \{ p \})$ in $\mathfrak{R}$.  Then $z \sim q, z \nsim p$, so $z \in C$, and $\sigma$ changes its adjacencies to $U_2$ and $V_2$, making it a witness for $(U_1 \cup U_2, V_1 \cup V_2)$.
\end{proof}

\medskip

The following result shows that the LPS groups are incomparable to the automorphism group of the neighbourhood filter of $\mathfrak{R}$.

\begin{proposition}
$\Aut^{*}(\mathfrak{RHyp}) \nleq \Aut(\mathscr{F}_{\mathfrak{R}})$.
\end{proposition}

\begin{proof}
Choose $v \in \mathfrak{R}$, and suppose that $\mathfrak{R}^{\sigma(v)}$ is $\mathfrak{R}$ except for switching about $v$.  Then exists an isomorphism $g : \mathfrak{R}^{\sigma(v)} \to \mathfrak{R}$.  Since $\mathfrak{R}^{\sigma(v)}$ is vertex-transitive, we can assume that $v \in \fix(g)$.  So $g : \mathfrak{R}(v) \to \mathfrak{R}_1 : = V(\mathfrak{R}) \backslash (\{ v \} \cup \mathfrak{R}(v))$.  But $\mathfrak{R}(v)$ is an edge of $\mathfrak{RHyp}$.  So $g \in \Aut^{*}(\mathfrak{RHyp}) \backslash \Aut(\mathfrak{RHyp})$.  But then $g \notin \Aut(\mathscr{F}_{\mathfrak{R}})$, since otherwise $\mathfrak{R}_1(v) \in \mathscr{F}_{\mathfrak{R}}$, which is impossible since $\mathfrak{R} (v) \cap \mathfrak{R}_1(v) = \emptyset$.  So $\Aut^{*}(\mathfrak{RHyp}) \nleq \Aut(\mathscr{F}_{\mathfrak{R}})$.

\end{proof}

We note that $\FSym$ acts locally, having a local effect on $\mathfrak{R}$ or $\mathfrak{RHyp}$;  $\SAut(\mathfrak{R})$ acts locally having a global effect on $\mathfrak{R}$;  $\DAut(\mathfrak{R})$ acts globally having a global effect on $\mathfrak{R}$.

\bigskip

A recent collaboration~\cite{camlaf} has extended the work of this section.  In this paper, a further group is defined
\begin{itemize}
\item $\FAut(\mathfrak{RHyp})$, the group of permutations $g$ with the property that there is a finite subset $S \subset \mathfrak{R}$ such that for every edge $E$, both $(E \backslash S) g$ and $(E \backslash S) g^{-1}$ are edges.
\end{itemize}
and the following results are proved,

\begin{theorem}
\item{(i)}  $\Aut(\mathfrak{RHyp}) < \FAut(\mathfrak{RHyp}) < \Aut^{*}(\mathfrak{RHyp}).$
\item{(ii)}  $\Aut_{fin}(\mathfrak{R}) < \FAut(\mathfrak{RHyp}).$
\item{(iii)}  $ \FSym(\mathfrak{R}) < \FAut(\mathfrak{RHyp})$, but $ \FSym(\mathfrak{R}) \cap \Aut(\mathfrak{RHyp}) = 1$.
\item{(iv)}  $\Aut(\mathfrak{RHyp}) . \FSym(\mathfrak{R})  <  \FAut(\mathfrak{RHyp}).$
\item{(v)}  $\SAut(\mathfrak{R}) \leq \Aut^{*}(\mathfrak{RHyp})$.
\item{(vi)}  $\BAut(\mathfrak{R}) \leq \Aut^{*}(\mathfrak{RHyp})$.
\item{(vii)}  $\SAut(\mathfrak{R}) \nleq \FAut(\mathfrak{RHyp})$.
\item{(viii)}  $ \Aut(\mathscr{F}_{\mathfrak{R}}) \nleq \Aut^{*}(\mathfrak{RHyp})$.
\end{theorem}

\bigskip
\bigskip

We display the groups we have encountered in this and the previous chapter, together with some of their inclusions in the following diagram.  We make three observations:-

(i)  $\FSym(\mathfrak{R}) = B(\mathscr{F}_{Frechet})$.

(ii)  $\Aut(\mathcal{T}) = \Aut(\mathcal{T}^*) = \Aut(\mathcal{T}^\dag)$, by Proposition~\ref{threetops}.

(iii)   We have omitted the groups $\NE(\mathfrak{R})$ and $\NS(\mathfrak{R})$ that we found in the previous chapter, as they are not permutation groups.  It would be interesting though to investigate if and how they relate to the 14 minimal functions on the random graph found by Bodirsky and Pinsker in~\cite{bodirsky}.
\index{Bodirsky, M.}%
\index{Pinsker, M.}%

\bigskip
\head{Open Question} 
Is it true that $\Aut(\mathcal{T}) < \Aut(\mathscr{F}_{\mathfrak{R}})$?

\bigskip

Group inclusions involving $\FAut(\mathfrak{RHyp})$ can be found in~\cite{camlaf}.

\clearpage

\begin{figure}[ht]
\vbox to \vsize{
 \vss
 \hbox to \hsize{
  \hss
  \vbox{
    $$
    \hss
\xymatrix@R=30pt@C=5pt{
&&& {\Sym(\mathfrak{R}) =  \Aut(\mathscr{F}_{Frechet})} \ar@{-}[dd] \ar@{-}[drrrr] \\
&&&&&&& \Aut^{*}(\mathfrak{Hyp}) \ar@{-}[d] \\
&&&  \Aut(\mathscr{F}_{\mathfrak{R}}) \ar@{-}[dddd] &&&& \Aut(\mathfrak{Hyp}) \ar@{-}[d]  \\
&&&&&&& \Aut(\mathcal{U}) &   \BAut(\mathfrak{R}) \ar@{-}[luu]\\
& \Aut_{fin}(\mathfrak{R})   \ar@{-}[ddrr]  \ar@{-}[uuurrrrrr] &&&&&& {\quad\quad \DAut(\mathfrak{R})}  \ar@{-}[ddddllll] \ar@{-}[ur] \ar@{-}[u] & \SAut(\mathfrak{R})\ar@{-}[ddddlllll] \ar@{-}[u]\\
&&& {}  \\
& B(\mathscr{F}_{\mathfrak{R}}) \ar@{-}[d]  \ar@{-}[uuuurr] && \Aut_{0}(\mathfrak{R})  \ar@{-}[d]  \ar@{-}[uuuurrrr] &&&& \Aut(\mathcal{T})  \ar@{-}[d] \ar@{-}[ddllll]  \ar@{-}[uuuuuullll] \\
& \FSym(\mathfrak{R}) \ar@{-}[uuuuurr] \ar@{-}[ddrr]  && \AAut(\mathfrak{R})  \ar@{-}[d] &&&&  \Aut'(\mathfrak{B}) \ar@{-}[dllll] \\
&&& \Aut(\mathfrak{R}) \ar@{-}[d] \\
&&& (1)
}
    \hss
    $$
        \caption{Hasse Diagram for the groups considered}
        \label{HD1} 
   }
 \hss}
\vss}
\end{figure}
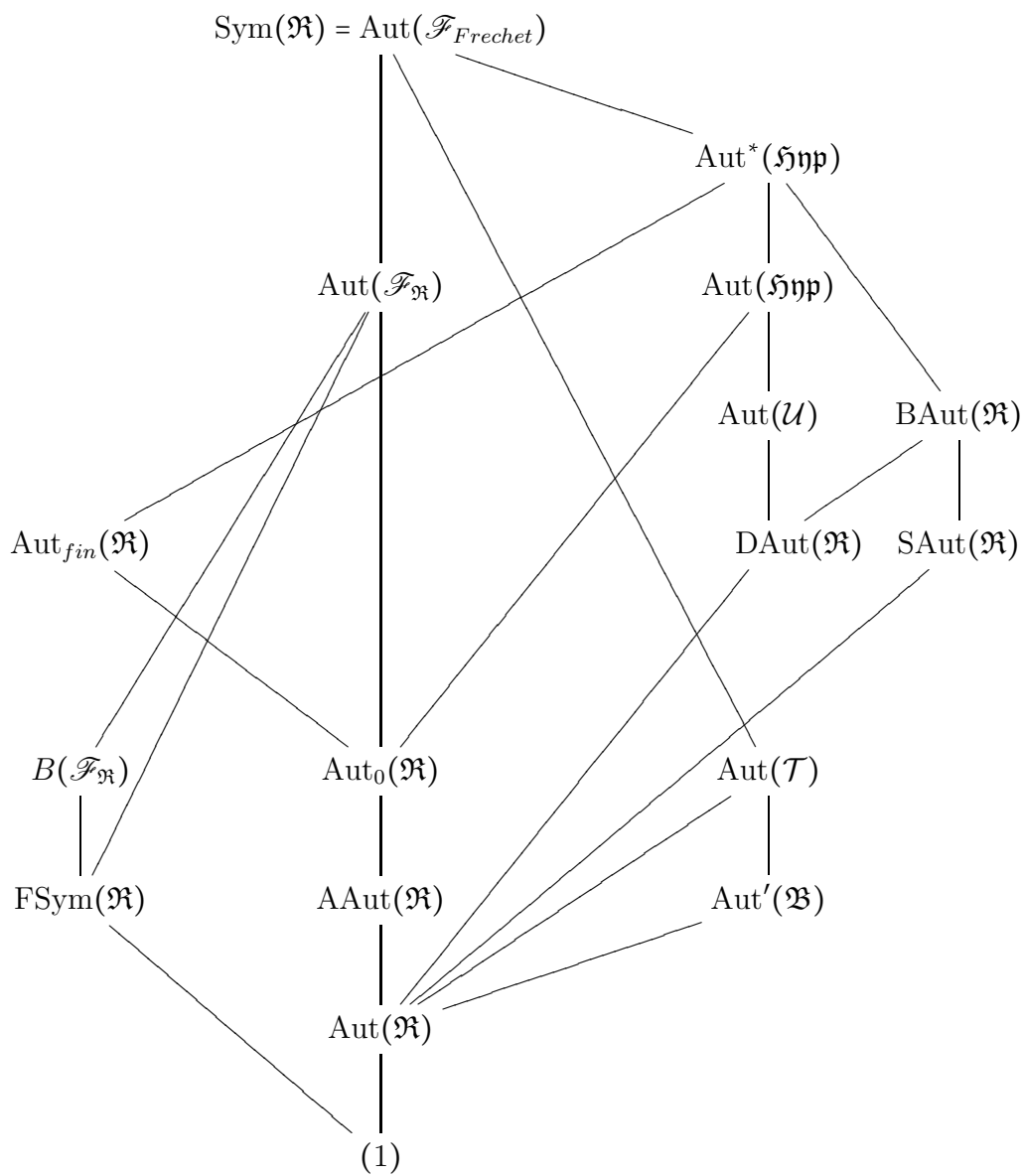
\index{Hasse diagram}%

\clearpage

\section{Stone--\v{C}ech Compactification}
\index{Stone--\v{C}ech compactification}%

In this section we see what the Stone--\v{C}ech compactification (maximal) compactification of a topological space\index{topological space}%
 looks like for neighbourhood filters
\index{filter}%
 on graphs.

Historically,  Grothendieck
\index{Grothendieck, A.}%
 viewed a space $X$ as being described by a topos of sheaves over $X$ (an object whose truth values correspond to the open sets of $X$), in which each point $x \in X$ has associated with it a filter of open neighbourhoods of $x$, these being the open subsets of $X$ containing $a$. This identifies a point with a filter.  This led Bourbaki to construct the completion of a uniform space and various compactifications including the Stone--\v{C}ech compactification~\cite{cartier}.

Recall the basic definitions.  Let $X$ be an infinite set. A \emph{filter} on $X$ is a family $\mathscr{F}$ of subsets of $X$ which is

\begin{itemize}
\item \emph{closed upwards}, i.e. $A\in\mathscr{F}$ and $B\supseteq A$ imply
$B\in\mathscr{F}$;
\item \emph{closed under intersection}, i.e.    $A,B\in\mathscr{F}$ implies $A\cap B\in\mathscr{F}$.
\end{itemize}

A filter $\mathscr{F}$ is \emph{non-trivial} if $\emptyset\notin\mathscr{F}$.  The word `filter' below will always mean `non-trivial filter'. (The only trivial filter is the power set of $X$.)
\index{filter ! trivial}%

An \emph{ultrafilter}
\index{filter ! ultrafilter}%
\index{ultrafilter}%
 is a filter $\mathscr{F}$ with the property that for any subset $A$, either $A\in\mathscr{F}$ or $X\setminus A\in\mathscr{F}$ (but not both, since $\mathscr{F}$ is non-trivial). Assuming the axiom of choice (as we do from now on), ultrafilters are the same as maximal filters.

Let $X^*$ be the set of ultrafilters on $X$. For any $x\in X$, the set $\mathscr{F}_x=\{A:x\in A\}$ is an ultrafilter (such ultrafilters are called \emph{principal}).
\index{ultrafilter ! principal}%
 The map $x\mapsto\mathscr{F}_x$ embeds $X$ into $X^*$.

Clearly any permutation group on $X$ induces a permutation group on $X^*$.

For any subset $A\subseteq X$, let $A^*$ be the set of all ultrafilters containing $A$. Clearly $X^*$ has the same meaning under this redefinition, and $\emptyset^*=\emptyset$. Moreover, $(A\cap B)^*=A^*\cap B^*$. (For an ultrafilter containing $A\cap B$ must contain $A$ and $B$, and conversely.)

Finally, we have $A^*\cap X=A$ for any $A\subseteq X$ (where $X$ is identified with the set of principal ultrafilters).

Thus the sets $A^*$, for all $A\subseteq X$, form the basis of a topology
\index{topology}%
 $\mathcal{T}$ on $X^*$ (that is, the open sets
\index{open set}%
 of $\mathcal{T}$ are all unions of sets of the form $A^*$), and the induced topology on $X$ is the discrete topology.
\index{topology ! discrete}%
 Permutations on $X$ induce homeomorphisms of this topology.

\begin{proposition}
The topological space $(X^*,\mathcal{T})$ is compact.
\end{proposition}
\index{topological space ! compact}%

\head{Proof} We have to show that any covering of $X^*$ by open sets has a finite sub-covering.

It is enough to consider coverings by basic open sets.
\index{sets ! basic open}%
  For given any covering $(U_i:i\in I)$, where $U_i=\bigcup_{j\in I_i}A_{i,j}^*$, we have a covering  by the basic open sets $A_{i,j}^*$ for all $i\in I$, $j\in I_i$. If such a covering has a finite sub-covering $A_{i_1,j_1}^*,\ldots,A_{i_n,j_n}^*$, then the original covering has a finite subcovering by $U_{i_1},\ldots,U_{i_n}$.

So suppose that $(A_i^*:i\in I)$ is a covering of $X^*$. Let us suppose, for a contradiction, that this covering has no finite sub-covering. The sets $(A_i:i\in I)$ form a covering of $X$. (For, if these sets don't contain a point $x$, then the principal ultrafilter $\mathscr{F}_x$ lies in no set $A_i^*$.) This covering has no finite sub-covering: for, if
\[A_1\cup\cdots\cup A_n=X,\]
then
\[A_1^*\cup\cdots\cup A_n^*=X^*,\]
contrary to assumption. Equivalently, the family $(A_i^c:i\in I)$ of complements of the sets $A_i$ has the property that any finite subfamily has non-empty intersection. Thus, the collection of all sets containing finite intersections of members of this family is a (non-trivial) filter. Let $\mathscr{F}$ be an ultrafilter containing it.

By assumption, $\mathscr{F}$ is covered by the sets $(A_i^*:i\in I)$; that is, $A_i\in\mathscr{F}$ for some $i$. But, by assumption, $\mathscr{F}$ contains $A_i^c$, and hence contains $\emptyset$, a contradiction.

So the result is proved.

\bigskip

We want to show that (a) $X^*$ is \emph{maximal}, (b) $X^*$ is \emph{unique}.

(a)  To show that $X^*$ is a maximal compactification of $X$ we can simply show that every compactification of $X$ is equivalent to a quotient space of $X^*$.  For this it suffices to find a continuous surjective closed map $f : X^* \to Y$, where $Y$ is an arbitrary compactification.

Certainly for any $A\subseteq X$, the map $f: A^* \to A^*\cap X=A$ is surjective.  If $A$ is a basic open set of the subcovering of $Y$, then $h_1: A^* \to A^*\cap Y=A$ is surjective.  This map is closed because the subcovering of $Y$ is finite.  Given the embedding $h: X \to X^*$ then it is clear that $h_{1}(f^{-1}(A)) = A$, and so $f$ is continuous.

(b)  Suppose that $X_i$ $(i = 1, 2)$ are two compactifications of $X$ and that $\kappa_i : X_i \to X^*$ are embeddings that extend the embedding $\kappa : X \to X^*$.  Then $\kappa_i$ maps $X$ onto $\kappa(X)$.  The $\kappa_i$ are continuous (because all functions between topological spaces
\index{topological space}%
 where the domain has the discrete topology
\index{topology ! discrete}%
 are continuous) and so map $X$ into the closure $\overline{\kappa(X)}$ of $\kappa(X)$.  But also the subspaces $\kappa_i(X_i)$ contains $X$ and are compact and so closed, and therefore they contain $\overline{X}$.  So $\kappa_i(X_i)=\overline{X}$, and $\kappa_2^{-1}\kappa_1$ is a homeomorphism stabilizing $X$.

\bigskip

\head{Remark} The space $X^*$ is the \emph{Stone--\v{C}ech compactification} of the discrete space $X$.  Our formulation of the Stone--\v{C}ech compactification is not the most general, but it is pertinent to a discussion in terms of filters.  The theory applies more generally to completely regular topological spaces and arbitrary metric spaces; see~\cite{johnstone}~\cite{walker}.

That any two Stone--\v{C}ech compactifications of the same topological space $X$ are homeomorphic, follows from (b) above by taking $X_i$ to be $X_i^*$ and noting that by uniqueness, $X_1^* \le X_2^* \le X_1^*$.

\begin{corollary}
$X$ embeds \emph{densely} in $X^*$, that is the closure of $X$ in $X^*$ is $X^*$.
\end{corollary}

\head{Proof}  We can restrict to working with basic open sets.
\index{sets ! basic open}%
  Take any such subset $A \subseteq X$.  Any open set
\index{open set}%
   in $X^*$ has a finite subcovering.  Each of the sets that form this subcovering contains a principal ultrafilter.  (That $X$ is countable and that any ultrafilter defining $A$ is closed under countable intersections implies that the ultrafilter is principal).  The principal ultrafilters $\mathscr{F}_A$ generated by $A$ all contain points of $A$.  So given the embedding $\iota: A \to \mathscr{F}_A$, every non-empty open set in $X^*$ contains $\iota(A)$ and so contains each $\mathscr{F}_A$ and so has a non-empty inverse image under $\iota$ and so contains $A$.  Therefore the image of $\iota$ is dense, and the result follows.  (In terms of closures, note that given $\kappa : X \to X^*$, if $\kappa(X)$ denotes a subspace of $X^*$ and $\overline{\kappa(X)}$ denotes its closure in $X^*$, then $\overline{\kappa(X)}$ is a compact Hausdorff space
\index{Hausdorff space}%
 and $\overline{X} = \overline{\kappa(X)}$; so $\overline{\kappa(X)}$ is a compactification of $\kappa(X)$).

\head{Caveat}

Here is another example of two different topologies on the same set giving rise to different properties.

Pestov~\cite{pestov}
\index{Pestov, V.}%
 proved the following theorem (see also~\cite{kechris}):

\begin{quote}
Let $G$ be a group of automorphisms of an infinite linearly ordered set $X$
acting transitively on $n$-subsets of $X$ for each $n$. Equip $G$ with the
topology of pointwise convergence.
\index{topology ! of pointwise convergence}%
 Then every continuous action of the topological group $G$ on a compact space has a fixed point. Equivalently,
the universal minimal $G$-flow is a singleton: $U(G)=\{{*}\}$.
\end{quote}

Take the linearly ordered set $X$ to be $\mathbb{Q}$; we may as well take
$G$ to be the group of all order-automorphisms of $\mathbb{Q}$.

Now take the discrete topology
\index{topology ! discrete}%
 on $\mathbb{Q}$, and let $U$ be its Stone--\v{C}ech compactification.
\index{Stone--\v{C}ech compactification}%
 Thus $U$ consists of all ultrafilters
\index{filter ! ultrafilter}%
  on the countable set $\mathbb{Q}$. A basis for the open sets
\index{open set}%
   of $U$ consists of the sets $\mathscr{F}(A)$ for $A\subseteq\mathbb{Q}$, where $\mathscr{F}(A)$ is the set of
ultrafilters containing $A$. (We have $\mathscr{F}(A\cap B)=\mathscr{F}(A)\cap \mathscr{F}(B)$, so the
set of all unions $\bigcup_{i\in I}\mathscr{F}(A_i)$ is a topology. It is well known
that it is compact.) The action of $G$ on $U$ is obviously continuous.

By Pestov's Theorem, $G$ fixes a point in $U$, that is, an ultrafilter $\mathscr{F}$
on $\mathbb{Q}$. 

But this is impossible. For let $A=\bigcup_{n\in\mathbb{Z}}[2n,2n+1)$.
Then exactly one of $A$ or its complement belongs to $F$. But the translation
$x\mapsto x+1$ in $G$ maps $A$ to its complement!

The answer to the apparent contradiction is that Pestov takes the usual permutation topology on the rationals whereas we have assumed the discrete topology.  This gives a different topology on the group of order-preserving permutations.  So the action of the group on the Stone--\v{C}ech compactification is not continuous for this topology.

Moreover, there is a natural compactification of $\mathbb{Q}$:  embed it in $\mathbb{R}$ and take the one-point compactification - this is a minimal compactification as opposed to Stone--\v{C}ech which is maximal.  Then the one added point is fixed.

\section{Other Graphs}

Let $\Gamma_i(x)$ denote the set of vertices at distance $i$ from a
vertex $x$ in the graph~$\Gamma$, and
$\overline{\Gamma}_k(x) = \bigcup_{i=0}^{k} \Gamma_i(x)$
denote the set of vertices within a distance $k$ from~$x$.

A graph $\Gamma$ with a path metric gives rise to a metric space $(V(\Gamma), d)$ where the distance $d(x, y)$ is the shortest path joining $x$ and $y$, and the largest value taken by $d$ is the \emph{graph diameter}.
\index{graph ! diameter}%

Which graphs have the property that the sets $\overline{\Gamma}_k(x)$ generate
a non-trivial filter, for some fixed $k$? Clearly we require that the diameter of such graphs to be
at most $2k$, else there are vertices $x$ and $y$ such that
$\overline{\Gamma}_k(x)\cap\overline{\Gamma}_k(y)=\emptyset$.  So defining the \emph{$k$-neighbourhood} of a vertex $v$ in $\Gamma$ to be the set of points distant at most $k$ from $v$, we require two $k$--neighbourhoods to not be disjoint.  Moreover, we can assume that the diameter is at least $k+1$, else every $k$-neighbourhood is the whole vertex set, $\overline{\Gamma}_k(x)=V(\Gamma)$ for all vertices $x$.

We make the following observation.  For every positive integer $d$, there is a countable homogeneous universal \emph{integral metric space}
\index{metric space ! integral}%
 (one with integer distances) of diameter $d$, unique up to isometry; see~\cite{cam8a}.  The metric is the path metric in a graph $M_d$ of diameter $d$.  Thus $M_2$ is the random graph $\mathfrak{R}$.
\index{graph ! random}%

It is easy to show that the filter generated by the $k$-neighbourhoods in $M_{2k}$ is isomorphic to $\mathscr{F}_{\mathfrak{R}}$.  

Is the following question true?  \emph{Let $\Gamma$ be a countable graph whose $k$-neighbourhoods generate a non-trivial filter.  Then $M_{2k}$ is a spanning subgraph of $\Gamma$}.
\index{graph ! spanning}%
\medskip

What can be said about other distance classes?

\section{Appendix: Graphs Spanned by $\mathfrak{R}$}

We first prove the following; (see also~\cite{camnes})

\begin{theorem}
A countable graph $\Gamma$ contains $\mathfrak{R}$ as a spanning subgraph if and only if any finite set of
vertices has a common neighbour in~$\Gamma$.
\end{theorem}

\begin{proof}
Let $V(\Gamma)=\{v_0,v_1,\ldots\}$ and $V(\mathfrak{R})=\{w_0,w_1,\ldots\}$. We
construct a bijection $\phi$ from $V(\mathfrak{R})$ and $V(\Gamma)$ by
back-and-forth,
\index{back-and-forth method}%
 which as so often with the random graph is the proof method of choice.

At even-numbered stages, choose the first unused vertex $w$ of $\mathfrak{R}$.
Let $U$ and $V$ be the sets of neighbours and non-neighbours of $w$
among vertices on which $\phi$ has been defined. Choose $v\in
V(\Gamma)$ joined to all vertices in $\phi(U)$, and extend $\phi$ to
map $w$ to $v$. In this extension, edges are mapped to edges.

At odd-numbered stages, choose the first unused vertex $v$ of
$\Gamma$. Let $U$ and $V$ be its neighbours and non-neighbours among
vertices in the image of $\phi$. Choose a vertex $w$ of $\mathfrak{R}$ joined
to no vertex in $\phi^{-1}(V)$, and extend $\phi$ to map $w$ to $v$.
Again, edges map to edges since the inverse image of a non-edge is a
non-edge. After countably many steps we have the required bijection
which takes edges of $\mathfrak{R}$ to edges of $\Gamma$.

Conversely, if $\Gamma$ contains $\mathfrak{R}$ as a spanning subgraph,
\index{graph ! spanning}%
 and $U$ is a finite set of vertices of $\Gamma$, then a non-neighbour to $U$
in $\mathfrak{R}$ (which certainly exists) is also a non-neighbour in $\Gamma$.
\end{proof}

Let $\Gamma$ be \emph{locally finite},
\index{graph ! locally finite}%
 that is, any vertex has only finitely many neighbours.  Embed $\Gamma$ as a spanning subgraph of $\mathfrak{R}$, which is possible by the preceding result. Then, on removing all the edges of $\Gamma$, we obtain a graph isomorphic to $\mathfrak{R}$.  So any locally finite graph is a spanning subgraph of $R$. 

This result is proved by showing that property $(*)$ can be verified without using any of the removed edges.

The statement in the first section of Chapter~\ref{homcochap} about Hamiltonian paths
\index{graph ! Hamiltonian path}%
 is a simple special case.

From this we can prove a remarkable decomposition theorem for $\mathfrak{R}$:

Let $\Gamma_1, \Gamma_2, \ldots$ be a sequence of locally finite graphs, each with at least one edge. Then we can decompose $\mathfrak{R}$ into spanning subgraphs isomorphic to $\Gamma_1, \Gamma_2, \ldots$.

To see this, we first enumerate all the edges of $\mathfrak{R}$. Since the automorphism group of $\mathfrak{R}$ is edge-transitive, we can find a spanning subgraph isomorphic to $\Gamma_1$ containing the first edge in the enumeration. Remove it; the resulting graph is still isomorphic to $\mathfrak{R}$, so we can find a spanning subgraph isomorphic to $\Gamma_2$ containing the first unused edge in the enumeration.  Iterating this procedure we find that every edge has been used.

So in particular, $\mathfrak{R}$ has a 1-factorisation, a Hamiltonian decomposition, and so on.

If we alter an arbitrary (as opposed to finite) set of edges, or switch with respect to an arbitrary set of vertices, the result is not always isomorphic to $\mathfrak{R}$; but it turns out that it is so `almost always', in either sense of the word.

\chapter{Monoids, Graphs and $\mathfrak{R}$}
\label{chapmon}
\bigskip

Rather, graphs provide a fundamental notational system for concepts and relationships that are not easily expressed in the standard mathematical languages of algebraic equations and probability calculus.
\begin{flushright}
Judea Pearl, \textit{Causality. Models, Reasoning, and Inference, Second Edition, p.138 -- Cambridge University Press (2009)}
\end{flushright}

\medskip

We have focussed thus far in this monograph on mappings that are isomorphisms between graphs.  In this chapter we turn our attention to aspects of homomorphisms between graphs in general, and the random graph in particular.  After an initial introductory section we shall derive results in separate directions, topological, model-theoretic and motivated by the theory of finite synchronizing groups we briefly look at infinite graph cores and hulls.
 
\section{Generalizations of Graph Homogeneity}

 \emph{Graph Homomorphism}.
\index{graph ! homomorphism}%
An \emph{arc}
\index{arc}%
 is a directed edge thus written as an ordered pair $(x,
y)$ instead of $\{x, y\}$.  In the study of graph homomorphisms it is
convenient to take \emph{digraphs} 
\index{digraph}%
 as the most basic objects.  Graphs with loops on vertices are then allowed, and a digraph is oriented
\index{graph ! oriented}%
if and only if it has no symmetric pairs of arcs.  A loop is the arc $(x, x)$.  A \emph{pseudograph}
\index{pseudograph}%
\index{graph ! pseudograph}%
 is a non-simple graph in which both graph loops and multiple edges are permitted, and a \emph{reflexive graph}
\index{graph ! reflexive}%
 is a pseudograph such that each vertex has an associated graph loop. 

A digraph homomorphism~\cite{hellnes}, $f: G \to H$, is a mapping $f(V): V(G) \to V(H)$ such that $(f(x), f(y)) \in E(H)$ whenever $(x, y) \in E(G)$.  Homomorphisms map \emph{paths} 
\index{graph ! path}%
 (a sequence of adjacent vertices and edges, with distinct vertices) to walks
\index{graph ! walk}%
 (a path in which vertices or edges may be repeated) in a graph; if $f$ is both vertex- and edge-bijective then it is an isomorphism.  Given that graph homomorphisms generalize isomorphisms, vertex colourings and arc directions, whilst preserving adjacency, it is not surprising that there are fruitful generalizations of the theory on which this treatise is based, which would be worthy of extension, in particular in light of~\cite{camnes}~\cite{camloc}, where a study is begun of relational structures
\index{relational structure}%
  (especially graphs and posets) which satisfy the analogue of homogeneity but for homomorphisms rather than isomorphisms.

We give the basic ingredients of this theory to motivate the reader.  Various notions of homogeneity arise, each being of the local-to-global variety, asserting that if two configurations have the same local structure in some sense, then they can be mapped to one another by a symmetry of the entire structure.  More precisely, a relational structure
\index{relational structure}%
 $M$ has property $XY$
\index{structure ! XY-property}%
 if every $x$-morphism between finite substructures of $M$ can be extended to a $y$-morphism from $M$ to $M$, where $(X, x)$ and $(Y, y)$ can be (I, iso), (M, mono), or (H, homo).  There are six properties of this kind that can be considered: HH, MH, IH, MM, IM, and II.  (Obviously a map cannot be extended to one satisfying a stronger condition.)  Note that II is equivalent to the standard notion of homogeneity.  These properties are related as in Figure~\ref{homtypes} with the strongest at the top:

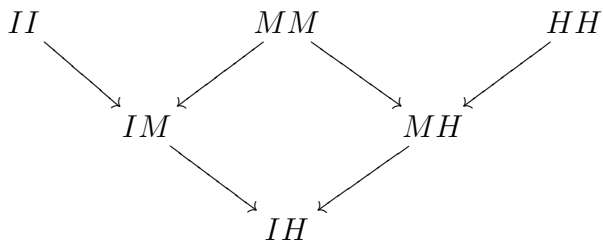
\begin{figure}[!h]
$$\xymatrix{ 
& {II} \ar@{->}[dr] && {MM} \ar@{->}[dr] \ar@{->}[dl] && {HH} \ar@{->}[dl]\\
&& {IM} \ar@{->}[dr] && {MH} \ar@{->}[dl]\\
&&& {IH}
}$$\caption{Relations between types of homogeneity}
\label{homtypes}
\end{figure} 

J. Ne\v{s}et\v{r}il and P. Cameron
\index{Ne\v{s}et\v{r}il, J.}%
\index{Cameron, P. J.}%
showed~\cite{camnes} that a countable MH graph either is an extension of
the random graph
\index{graph ! random}%
 $\mathfrak{R}$ (containing it as a spanning
\index{graph ! spanning}%
 subgraph), or has bounded claw size. Apart from disjoint unions of complete
graphs (containing no $K_1$ or $K_2$), no examples with bounded claw
size are known. Extensions of $\mathfrak{R}$ are MM and HH, and hence MH.  

All finite HH graphs are known. The known countable examples are disjoint unions of complete graphs of the same size, and graphs which contain the random graph
\index{graph ! random}%
 $\mathfrak{R}$ as an induced subgraph.  In~\cite{rusinov} a construction is given of graphs that do not contain $\mathfrak{R}$ as spanning subgraph, but nevertheless are HH.

The homogeneous (II) graphs were classified in~\cite{lach} by Lachlan and Woodrow,
\index{Lachlan, A. H.}%
\index{Woodrow, R. E.}%
but the classification of IH or IM graphs is open.

\bigskip

The graph $\mathfrak{R}$ has been proved \emph{not} to be characterized by the properties of universality and homomorphism homogeneity~\cite{camnes}.

\bigskip

\emph{Core-transitive Graphs}.
\index{graph ! core-transitive}%
 A homomorphism of graphs
\index{graph ! homomorphism}%
 is a mapping of vertices which carries edges to edges; its behaviour on non-edges is unspecified.  If there are homomorphisms in both directions the two graphs are said to be \emph{homomorphism-equivalent}.
\index{graph ! homomorphism-equivalent}%
 
Any homomorphism-equivalence class contains a unique smallest graph (up to isomorphism) called a \emph{core}.
\index{graph ! core}%
 For any graph $\Gamma$, the core of the equivalence class $[\Gamma]$ is unique up to isomorphism and is embeddable as an induced subgraph of $\Gamma$.  The core of $\Gamma$ has an embedding into $\Gamma$, and there is a \emph{retraction}
\index{graph ! retraction}%
 from a graph to its core, that is a homomorphism which fixes its image pointwise.  Alternatively, a \emph{retraction} of a graph is an endomorphism $f$ satisfying $f^2 = f$ (that is, $f$ is the identity on its image). The image of a retraction is a \emph{retract}.
\index{graph ! retract}%
  Given any map $f$ from a finite set to itself, some power of it is a retract.  For there is a number $n$ such that the image of any point under $f^n$ falls into the set of recurrent points of $f$; and $f^n$ induces a permutation of the set of recurrent points, so some power of it is the identity.

As we have already noted above, Sabidussi
\index{Sabisussi, G.}%
 showed~\cite{sabisussi1} that every vertex-transitive graph is a retract of a Cayley graph.  A graph is a core if it has no non-trivial retraction.  A \emph{graph endomorphism}
\index{graph ! endomorphism}%
is a homomorphism from a graph to itself.  The endomorphisms of a graph form a semigroup,
\index{semigroup}%
 and a graph is a core if and only if all its endomorphisms are automorphisms.  A retract is the image of an \emph{idempotent} endomorphism of the graph.  All cores of a finite graph are isomorphic, so we can speak of \emph{the} graph core.  Welzl
\index{Welzl, E.}%
 showed~\cite{welzl} that the core of a vertex-transitive graph is vertex-transitive, but stronger statements have been proved: for example, the core of a nonedge-transitive graph is either complete or the graph is its own core.   A graph $\Gamma$ is called \emph{core-transitive}
\index{graph ! core-transitive}%
 if any isomorphism between cores of $\Gamma$ can be extended to an automorphism of $\Gamma$.  Whilst any core is core-transitive, for graphs whose core is complete this is a strong condition.  We will shortly see a weakening of this condition, but first we have

\begin{proposition} If $\Gamma$ is vertex-transitive, then so is $Core(\Gamma)$. 
\label{vxtcore}
\end{proposition}
\begin{proof} Let $\iota : Core(\Gamma) \to \Gamma$ be an embedding as an induced subgraph, and $\rho : \Gamma \to Core(\Gamma)$ be a retraction.
\index{graph ! retraction}%
  Suppose that $\Gamma$ is vertex-transitive, and let $v, w$ be vertices of $Core(\Gamma)$.  If $g \in \Aut(\Gamma) :  v g = w$, then $g' = \iota g \rho$ is an endomorphism of $Core(\Gamma)$ taking $v$ to $w$.  Because $Core(\Gamma)$ is a core, then $g'$ is an automorphism of $Core(\Gamma)$.
\end{proof}

Similarly for other kinds of transitivity. 

The \emph{hull}
\index{graph ! hull}%
 of a graph $\Gamma$ is a graph containing $\Gamma$ as a spanning
\index{graph ! spanning}%
  subgraph, admitting all the endomorphisms of $\Gamma$, and having as core a complete graph of the same order as the core of $\Gamma$.  A graph $\Gamma$ is a core if and only if its hull is complete.

More on graph homomorphisms and cores can be found in~\cite{hahntardif} and in~\cite{hellnes}. 

A graph is a core if and only if all its endomorphisms are automorphisms.  Such graphs cannot be $MM$, $MH$ or $HH$, and the other three classes collapse into one: $II = IM = IH$.

Core-transitivity of a graph is weaker than homogeneity, and generalizing further, we can ask for graphs in which core isomorphisms extend to mono/homomorphisms of the entire graph.  

\smallskip

The $II$ graphs are rare and classified, whilst the class $IM$ requires more research.

\smallskip

If $\Gamma$ has $\mathfrak{R}$ as a spanning
\index{graph ! spanning}%
 subgraph, then it is $MH$ and so $IH$; in particular if $f$ is an embedding of any finite part $X'$ of $X \cup \overline{\mathfrak{R}} (x)$ into $\Gamma$, then $f$ can be extended into a monomorphism from $X' \cup w$ where the vertex $w \in \Gamma \backslash \{X \cup \overline{\mathfrak{R}} (x) \}$ into $Y' \cup z$ for some finite induced subgraph $Y'$ and $z \in \Gamma \backslash \{Y \cup \overline{\mathfrak{R}} (y) \}$.

\section{Graphs and Finite Transformation Monoids}
\index{monoid ! transformation}%
This section describes a pair of mappings between graphs and transformation monoids on the set 
$\{1, \ldots, n\}$, and some of their properties. 

Let $X$ be an $n$-element set, and $T(X)$ the \emph{monoid},
\index{monoid}%
 that is the (full transformation) semigroup
\index{semigroup}%
  with identity, of all maps $X \to X$.  Clearly $|T(X)| = n^n$ and $\Sym(X) \leq T(X)$.  For $f \in T(X)$, the \emph{rank} of $f$ is the cardinality of its image.  A monoid $M \leq T(X)$ is \emph{synchronizing}
\index{monoid ! synchronizing}%
 if there exists an $f \in M$ whose rank is 1.  Recent research interest has focussed on amongst other things, finding subsets $S \subseteq T(X)$ with $M := \langle S \rangle$. 

The set of endomorphisms of a graph $\Gamma$ is closed under composition and contains the identity; 
that is, it forms a \emph{monoid} $\End(\Gamma)$. An endomorphism is an automorphism if it is 
bijective.  If a graph is edge-transitive and not a core,
\index{graph ! core}%
\index{graph ! edge-transitive}%
 then its core is complete, and it has among its endomorphisms both automorphisms (of maximal rank) and proper colourings (of minimal rank).  
 
Two graphs $\Gamma_X$ and $\Gamma_Y$ are \emph{hom-equivalent}
\index{graph ! hom-equivalent}%
 if there are homomorphisms in both directions between them.  A graph is a \emph{core}
\index{graph ! core}%
 if it is not hom-equivalent to any graph with fewer vertices; and $\Gamma_Y$ is a core of $\Gamma_X$ if $\Gamma_Y$ is a core and is hom-equivalent to $\Gamma_X$. 
 
 We can recognise a core by the following property: 

\textbf{Fact 1} \emph{The graph $\Gamma$ is a core if and only if $\End(\Gamma) = \Aut(\Gamma)$.} 

For if $\Gamma_Y$ is a core and $f$ an endomorphism of $\Gamma_Y$, then the image of $f$ is hom-equivalent 
to $\Gamma_Y$ (the homomorphism in the other direction is just the embedding), so cannot be properly contained in $\Gamma_Y$. 

If $\Gamma_Y$ is a core of $\Gamma_X$, then a homomorphism from $\Gamma_Y$ to $\Gamma_X$ is injective, and a homomorphism from $\Gamma_X$ to $\Gamma_Y$ is surjective. (The first holds because, if $f : Y \to X$ is non-injective and $g : X \to Y$ is any homomorphism, then $f g$ is a non-injective endomorphism of $\Gamma_Y$. The second is similar.)  The homomorphisms in both directions between them are bijective and hence are isomorphisms, and so two hom-equivalent cores are isomorphic.  In other words, any graph has a unique core.

The \emph{clique number} 
\index{graph ! clique number}%
$\omega(\Gamma)$ is the cardinality of the largest complete induced subgraph of $\Gamma$; the \emph{chromatic number}
\index{graph ! chromatic number}%
\index{chromatic number}%
$\chi(\Gamma)$ is the smallest number of colours required to colour the vertices so that adjacent vertices get different colours.  (It is possible for a graph to have an infinite chromatic number, for example the Henson graphs $\mathfrak{H}_k$ $(2 < k < \omega)$~\cite[p.~81]{macneu}).
\index{graph ! Henson}%
  Clearly $\omega(\Gamma) \leq \chi(\Gamma)$, since vertices of a complete subgraph must get different colours.  

\textbf{Fact 2} \emph{The core of a graph $\Gamma$ is complete if and only if $\omega(\Gamma) = \chi(\Gamma)$.}

\begin{proof} $\omega(\Gamma) = m$ (that is an $m$-clique exists) if and only if there is a homomorphism $K_m \to \Gamma$, while $\chi(\Gamma) = m$ (that is an $m$-colouring exists) if and only if there is a homomorphism $\Gamma \to K_m$.  So $\omega(\Gamma) = \chi(\Gamma)  = m$ if and only if $\Gamma$ is homomorphically equivalent to $K_m$ (which is then necessarily the core of $\Gamma$). 
\end{proof}

In the other direction, let $M$ be a transformation monoid
\index{monoid ! transformation}%
 on the set $\{1, \ldots, n\}$, a submonoid of the full transformation monoid $T_n$. From $M$, we construct a graph as follows. Its vertex set is 
$\{1, \ldots, n\}$; two vertices $i$ and $j$ are joined by an edge if and only if there is no element $f \in M$ for which $i f = j f$. Denote this graph by $\Gr(M)$. 

A transformation monoid is synchronizing
\index{monoid ! synchronizing}%
 if it contains an element whose image has cardinality 1. 

\textbf{Fact 3} \emph{(a) $\Gr(M)$ is complete if and only if $M$ is a permutation group (that is, 
contained in the symmetric group).} 

(b) \emph{$\Gr(M)$ is null if and only if $M$ is synchronizing.} 

\begin{proof} (a) $\Gr(M)$ is complete if and only if no element of $M$ ever maps two points to the same place. 
(b) Let $f \in M$ be an element whose image is as small as possible. Then no two 
elements of the image of $f$ can be mapped to the same place; so they are pairwise 
adjacent. So, if $\Gr(M)$ is null, then the image of $f$ has cardinality 1. The converse 
is clear. 
\end{proof} 

\textbf{Fact 4} \emph{For any transformation monoid $M$, the graph $\Gr(M)$ has core a complete graph.} 

\begin{proof}  The argument in (b) above shows that the image of an element of $M$ of 
minimal rank is a complete subgraph of $\Gr(M)$. It is hom-equivalent to $\Gr(M)$ 
(the homomorphism in the other direction is just the embedding), and it is clearly a core. 
\end{proof} 

\textbf{Fact 5} \emph{For any transformation monoid $M$, 
(a) $M \leq \End(\Gr(M))$; 
(b) $\Gr(\End(\Gr(M))) = \Gr(M)$.}
 
\begin{proof} (a) Let $f$ be an endomorphism of $M$, and let $i$ and $j$ be adjacent in $\Gr(M)$. 
By definition, $i f \neq j f$. Could $i f$ and $j f$ be non-adjacent in $\Gr(M)$? If so, then 
there is an element $h \in \End(M)$ with $(i f )h = j f (h)$. But this contradicts the 
adjacency of $i$ and $j$, since $f h \in M$ by closure. 

(b) Suppose first that $i$ and $j$ are adjacent in $\Gr(M)$. Then no endomorphism 
of $\Gr(M)$ can collapse them, so they are adjacent in $\Gr(\End(\Gr(M)))$. 
Conversely, suppose that $i$ and $j$ are not adjacent in $\Gr(M)$. Then there is an 
element $f \in M$ satisfying $i f = j f$. By (a), $f \in \End(\Gr(M))$, and so $i$ and $j$ are 
non-adjacent in $\Gr(\End(\Gr(M)))$.
\end{proof} 

It is not true that $\End(\Gr(\End(\Gamma))) = \End(\Gamma)$ for all graphs $\Gamma$.  For let $\Gamma$ 
be the path of length 3, with just two automorphisms. It is easy to see that no endomorphism can identify the ends of the path, so that $\Gr(\End(\Gamma))$ is the 4-cycle, with eight automorphisms. 

\textbf{Fact 6} \emph{The maps $M \mapsto \End(\Gr(M))$ and $\Gamma \mapsto \Gr(\End(\Gamma))$ are idempotent.}
\begin{proof} This follows immediately from part (b) of the preceding Fact. 
\end{proof} 

Write $\Cl(M) = \End(\Gr(M))$, where $\Cl(M)$ denotes the \emph{closure} of the monoid.
\index{monoid ! closure}%
 Then $M \leq \Cl(M)$ and $\Cl(\Cl(M)) = \Cl(M)$, so $\Cl()$ is a closure operator on transformation monoids on $\{1, \ldots, n\}$.  We do not have a satisfactory description of the closed objects; but more on this below. 

In the other direction, let $Hull(\Gamma) = \Gr(\End(\Gamma))$.  Then we have that $Hull(Hull(\Gamma)) = Hull(\Gamma)$. The hull of a graph has the following properties: 

\textbf{Fact 7} \emph{(a) $\Gamma$ is a spanning
\index{graph ! spanning}%
 subgraph of $Hull(\Gamma)$ (that is, these graphs have the 
same vertex set, and every edge of $\Gamma$ is an edge of $Hull(\Gamma))$.}

\emph{(b) $\End(\Gamma) \leq \End(Hull(\Gamma))$ and $\Aut(\Gamma) \leq \Aut(Hull(\Gamma))$.} 

\emph{(c) $Core(Hull(\Gamma))$ is a complete graph on the vertex set of $Core(\Gamma)$.} 

\begin{proof} (a) If $i$ and $j$ are adjacent in $\Gamma$, then no endomorphism of $\Gamma$ can collapse $i$ 
and $j$, so they are adjacent in $\Gr(\End(\Gamma))$.
 
(b) Immediate from Fact 5(a). 

(c) The vertex set of $Core(\Gamma)$ cannot be collapsed by endomorphisms, so is a 
complete subgraph of $\Gr(\End(\Gamma)) = Hull(\Gamma)$. 
\end{proof} 

By (c), if $\Gamma$ is a hull, then $Core(\Gamma)$ is complete; but the converse is false. If $\Gamma$ is the path of length 3, then $Core(\Gamma)$ is a complete graph on two vertices, but $Hull(\Gamma)$ is the 4-cycle, by our previous argument. 

\textbf{Fact 8}  \emph{A transformation monoid $M$ is \emph{closed} (that is, satifies $M = \Cl(M))$ if and only if $M = \End(\Gamma)$ for some graph $\Gamma$ which is a hull (and in particular, whose core is complete).} 

\begin{proof} Suppose that $M$ is closed. Then $M = \End(\Gamma)$, where $\Gamma = \Gr(M)$; so $\Gamma = \Gr(\End(\Gamma)) = Hull(\Gamma)$. 
Conversely, if $\Gamma = Hull(\Gamma)$, then $\End(\Gamma) = \End(\Gr(\End(\Gamma))) = \Cl(\End(\Gamma))$. 
\end{proof} 

For any graph $\Gamma$,
\begin{itemize}
\item  $Hull(\Gamma)$ is complete if and only if $\Gamma$ is a core; 
\item  if $Hull(\Gamma) = \Gamma$ then $Core(\Gamma)$ is complete (but the converse is false).
\end{itemize}

\smallskip

\emph{Remark on closure.}
\index{closure axioms}%
 The axioms that a traditional closure operation $\cl$ should obey are~\cite{erne}:

(C1) $A \leq \cl(A)$\quad \emph{extensive};

(C2) $A \leq B \Rightarrow \cl(A) \leq \cl(B)$\quad \emph{isotone};

(C3) $A \leq \cl(\cl(A)) = A$ (or $\cl(\cl(A)) = \cl(A)$)\quad \emph{idempotent}.

We have seen that $M_1 \leq M_2 \Rightarrow \Gr(M_2) \leq \Gr(M_1)$.  Take $i, j \in \Gr(M_2)$ such that $i \sim j$.  Then $i \sim j$ also in $\Gr(M_1)$.  But if $\exists f \in \End(\Gr(M_2))$ such that $if =jf$ then $f \in \End(\Gr(M_1))$.  So $\Cl(M) \equiv \End(\Gr(M))$ satisfies (C1) and (C3), but is \emph{not} isotone, and cannot be a closure operator, as normally understood.

\smallskip

Finally, we reference an article~\cite{sauerst1} by Sauer and Stone
\index{Sauer, N.}%
\index{Stone, M.}%
 in which monoids and their local closures are studied.

\section{Hulls and Cores of Infinite Graphs}
\index{graph ! core}%

Any nonedge-transitive graph is either a core or a hull.  For if $f$ is an endomorphism of $\Gamma$ mapping two non-adjacent vertices in $\Gamma$ to non-adjacent vertices in $Hull(\Gamma)$ then $\Gamma = Hull(\Gamma)$.  Otherwise any pair of non-adjacent vertices can be mapped to a pair that are collapsed by $f$ and then application of $f$ maps them to a single vertex.  


We first make a few comments on infinite hulls. 
\index{graph ! hull}%

If $\Gamma$ is any countable graph which contains an infinite complete subgraph, then $\Gamma$ is a hull.  For let $K_{\aleph_{0}}$ be the infinite complete subgraph and $v$ and $w$ any two non-adjacent vertices in $\Gamma$.  Choose a bijection from $\Gamma - v$ to $K_{\aleph_{0}}$ (this is clearly a homomorphism). Now extend it to $\Gamma$ by mapping $v$ to the same place as $w$. The only pair not mapped to an edge is $v, w$ and by assumption it was not an edge in $\Gamma$.

This indicates that hulls, while easy to define in the infinite case (unlike cores), are not actually of any use.  So $\mathfrak{R}$ and all its overgraphs are hulls.  

\smallskip

Given a transformation monoid $M$ on $\Omega$, we can consider its closure in the topology of pointwise convergence: a sequence $(f_n)$ of elements of $M$ converges to the limit $f$ if, for all $x \in \Omega$, there exists $n_0$ such that for all $n \geq n_0$ we have $f_n(x) = f(x)$.
 
Now we say that a permutation group $G$ is \emph{strongly synchronizing}
\index{group ! permutation ! strongly synchronizing}%
 if, for any map $f$ which is not injective, the closure of $M = \langle G, f \rangle$ contains an element of rank 1. 

A characterization of synchronizing groups that works in the infinite case is: a permutation group $G$ is non-synchronizing if and only if there is a non-trivial hull $\Gamma$ on its vertex set such that $G \leq \Aut(\Gamma)$.  However in the case of infinite groups and monoids, quite unlike the finite case, there is a confluence of different concepts:

\begin{corollary}
The following are equivalent
\begin{itemize}
\item $G$ is synchronizing;
\item $G$ is strongly synchronizing;
\item $G$ is $2$-set transitive.
\end{itemize}
\end{corollary}





\bigskip

We next turn to cores, whose definition for infinite graphs is more problematic to define that for hulls.  A survey on graph homomorphisms and cores can be found in~\cite{hahntardif}, and more in~\cite{hellnes}.  

If $\Gamma$ is finite and $\exists f \in \End(\Gamma)$ which is non-injective then $\Gamma$ has a proper retract.  For $\exists n : f^n$ is a permutation of the image, and then $(f^n)^k$ is the identity on the image for some $k$, and so a retract.  We state the next result without proof.

\begin{proposition}
The following are equivalent for a graph $\Gamma$:

(a)  A graph $\Gamma$ is a core 

(b)  All its endomorphisms are injective.

(c)  It is a minimal retract.

\end{proposition}




For finite graphs both the existence of a minimal core and its uniqueness are clear, but for infinite graphs the question is more subtle.

As was pointed out in~\cite{hahntardif}, when an infinite graph has finite homomorphic images, compactness must be taken into account; for example an arbitrary graph admits a homomorphism to a finite graph if and only if its finite subgraphs do.  Furthermore, many of the properties that hold for finite graphs, fail in the case of infinite graphs.  For example, that a graph is a core if and only if every endomorphism is an automorphism was shown by Bauslaugh~\cite{bauslaugh}
\index{Bauslaugh, B.}%
not to hold for directed graphs.  It is also possible to construct graphs with no core in the sense that there is no minimal retract,
\index{graph ! retract}%
 as per the following example~\cite[p.~158]{hahntardif}.  From the complete graphs $K_k$ on $k$ vertices $\{v_1, \ldots, v_k\}$ $(k > 1)$, construct a graph $\gamma$ by identifying $k$ with $v_k$ for $k \in \mathbb{N}$.  The graph $\gamma$ has no core in the sense of minimal retract, because there are retraction maps $\Gamma_i \to \Gamma_j$ only if $i < j$ and so the retracts of $\gamma$ are the infinite components of $\gamma - k, k \in \mathbb{N}$.  Core-like properties of directed graphs and implications between these properties, are studied in~\cite{bauslaugh}, whilst other characteristics of cores of random graphs were examined in~\cite{bonatopralat}.

Two minimal requirements of the definition of a core are: 

(a) it should agree with the standard definition in the case of finite graphs; and 

(b) any complete graph is a core. 

We will test these putative definitions against this requirement. 

The first attempt is to take over the finite definition and ask that the core of a graph is the smallest graph hom-equivalent to it. This satisfies (a), but not (b). For infinite cardinals this is too crude: all countable graphs, 
for example, have the same number of vertices. 

The second attempt is to take the definition of a core as a graph all of whose endomorphisms are automorphisms. This also satisfies (a), but not (b); for an infinite complete graph has an endomorphism onto any infinite 
subgraph of the same cardinality. 

The third attempt builds on this. We say that a graph $\Gamma_X$ is a \emph{core} if every endomorphism is injective. This clearly satisfies (a) and (b). 

The fourth definition uses retracts. A graph $\Gamma_X$ is an \emph{R-core}
\index{graph ! R-core}%
 if its only retraction is the identity. This also satisfies (a) and (b). 

Clearly any core is an R-core. For if $\Gamma_X$ is a core, and $f$ is a retraction, then $(xf )f = xf$ for all $x \in X$; and so $xf = x$, since otherwise the points $x$ and $xf$ witness non-injectivity of $f$. So $f$ is the identity. 

\head{Open Question 1}  Is the converse true? That is, are the two notions of core equivalent? (It seems possible that X could be an R-core but have a non-injective endomorphism no power of which is a retraction.)

Note that the definition of a hull works without problems in the infinite case: a graph is a hull if and only if any two non-adjacent vertices are mapped to the same vertex by some endomorphism. Now the first definition of a core 
retains a property of finite graphs we saw in the preceding section: a graph is a core if and only if its hull is complete. 

Here is an example of a graph that is not a core.

Let $\Gamma_X$ be the disjoint union of complete graphs $K_n$ for all positive integers $n$.  What does a graph $\Gamma_Y$ hom-equivalent to $\Gamma_X$ look like?  Assume that $\Gamma_Y$ is a core.  Let $f : \Gamma_Y 
\to \Gamma_X$ and $g : \Gamma_X \to \Gamma_Y$ be homomorphisms; then $f g$ is injective, so in particular $f$ is injective.  Let $v$ and $w$ be vertices of $\Gamma_Y$ which belong to different components of $\Gamma_X$. There is an endomorphism $h$ of $\Gamma_X$ which fixes pointwise every component except that containing $vf$, and maps $vf$ to $wf$. Then $vf h = wf h$, so $vf hg = wf hg$; thus $\Gamma_Y$ has a non-injective endomorphism, contrary to assumption. 

\head{Open Question 2}  Is there an infinite graph which has no R-core?

In particular, is this true for the graph $\Gamma_X$ of the preceding paragraph? 

\head{Open Question 3}  Do there exist cores (in either sense) which are equivalent without being isomorphic? 

If $\Gamma_X$  and $\Gamma_Y$  are equivalent cores, then there exist injective homomorphisms $f : \Gamma_X \to \Gamma_Y$ and $g : \Gamma_Y \to \Gamma_X$. Now one might hope that a Cantor-Schr\"oder-Bernstein argument 
\index{Cantor-Schr\"oder-Bernstein Theorem}%
would produce an isomorphism betweem $\Gamma_X$ and $\Gamma_Y$ ; the problem is that the inverse of an injective homomorphism need not be a homomorphism. 

\begin{proposition}
A countably infinite graph containing an infinite clique is not a core unless it is complete.
\end{proposition}
\begin{proof}
Suppose that $v$ and $w$ are non-adjacent vertices.  There is a bijection $f$ from $\Gamma_X \backslash \{w\}$ to a clique $C$ of $\Gamma_X$.  Extend this map to $w$ by setting $wf = vf$.  The resulting map is an endomorphism of $\Gamma_X$ and is not injective.  

In particular, the countable random graph (the unique countable homogeneous universal graph) is not a core. 

In fact it is not an R-core either, though this is a little more diffcult. The countable random graph is characterised up to isomorphism by the property that, if $U$ and $V$ are disjoint finite sets of vertices, then there is a vertex $z$ 
joined to everything in $U$ and nothing in $V$. Now, if a countable graph is constructed by a random process such that, for any disjoint finite sets $U$ and $V$, there is an infinite set $Z$ such that the events that $z$ is joined to all of $U$ and none of $V$ are independent (for $z \in Z$ ) with fixed non-zero probability, 
then the resulting graph is almost surely isomorphic to $\mathfrak{R}$. 

Begin with countably many pairwise disjoint countable sets $\{A_i\}$, with a chosen point $a_i \in A_i$ for all $i$.  Now insist that the vertices $a_i$ are pairwise adjacent, and that there are no edges within any set $A_i$.  Then the above condition is satisfied.  (Given $U$ and $V$, let $Z$ consist of all points lying in sets $A_k$ containing none of $U \cup V$ except for the chosen points $a_k$ in these sets). So the resulting graph is isomorphic to the random graph. Now the map that collapses every set $A_i$ onto its representative point $a_i$ is obviously a retraction. 
\end{proof}

The argument actually shows that any countably infinite graph with an infinite clique is a hull. Clearly it works for any infinite graph having a clique of the same cardinality of the whole graph. 

\head{Open Question 4}  What happens for countable graphs which contain arbitrarily large finite cliques but no countable cliques? 

\head{Open Question 5} What happens for uncountable graphs whose largest cliques are countably infinite? 

For example, assuming the axiom of choice (AC),
\index{axiom of choice}%
  take a well-ordering of the real numbers, and define a graph in which $v$ and $w$ are joined if and only 
if the well-order and the usual order agree on $\{v, w\}$.  It is well-known that the largest clique (and the largest independent set) in this graph are both countably infinite.  Is this graph a core? 

\begin{proposition}
Let $\Gamma_X$ be a graph whose automorphism group is transitive on non-edges.  Then either $\Gamma_X$ is an core, or it is a hull. 
\end{proposition}
\begin{proof}
If there exists a non-injective endomorphism $f$, with $vf = wf$, then for any non-adjacent vertices $v'$ and $w'$ there is an endomorphism identifying them.  (Let $g$ be an automorphism mapping $\{v, w\}$ to $\{v', w' \}$, and let $f' = g^{-1} f g$.) 
\end{proof}

\begin{proposition}[Nick Gravin]
\index{Gravin, N.}%
 Let $\Gamma_X$ be an infinite graph with finite clique number $\omega(X)$, and suppose that $\Gamma_X$ is a hull.  Then $\chi(X) = \omega(X)$.
\end{proposition}
\begin{proof}
Given a finite subgraph $\Gamma_Y$ of $\Gamma_X$, if $\Gamma_Y$ is not complete, then there is an endomorphism $f_1$ of $\Gamma_X$ collapsing a non-edge of $\Gamma_Y$; if $\Gamma_Y f_1$ is not complete, 
there is an endomorphism $f_2$ collapsing a non-edge of $Y f_1$; and so on. We end with a homomorphism of $\Gamma_Y$ to a complete graph of size at most $\omega(X)$.  So $\chi(\Gamma_Y) \leq \omega(\Gamma_X)$ for any finite subgraph $\Gamma_Y$. A compactness argument shows that $\chi(\Gamma_X) \leq \omega(\Gamma_X)$, so equality holds. 

Hence, if $\Gamma_X$ is non-edge-transitive and has $\omega(\Gamma_X)$ finite and $\chi(\Gamma_X) \neq \omega(\Gamma_X)$, then $\Gamma_X$ is a core.  Examples include the Henson graph $\mathfrak{H}_k$ ($n \geq 3$).
\index{graph ! Henson}%
\end{proof}

\section{Fra\"{\i}ss\'e's Theorem for HH Structures}
\index{Fra\"{\i}ss\'e's Theorem}%
\index{relational structure}%

The two central results linking oligomorphic permutation groups,
\index{group ! permutation ! oligomorphic}%
 model theory
\index{model theory}%
  and combinatorics
\index{combinatorics}%
 are firstly Fra\"{\i}ss\'e's Theorem
\index{Fra\"{\i}ss\'e's Theorem}%
in which the amalgamation property to a large extent characterizes the classes of finite structures whose enumeration is equivalent to finding the number of orbits on $n$-sets, and secondly the Engeler--Ryll-Nardzewski--Svenonius Theorem.
\index{Engeler--Ryll-Nardzewski--Svenonius ! Theorem}%

In this section we prove a version of the former for HH relational structures
\index{relational structure}%
 in place of the usual universal homogeneous structures.  As we have already seen, both HH and MM relational structures are subclasses of MH structures, and the statements of this section have exact analogues for MM structures and the proofs of these are in~\cite{camnes}.
 
Concepts such age
\index{age}%
 and the joint embedding property (JEP; see Appendix~\ref{TheoryofRelationalStructures}) are retained from the usual theory.  A class $\mathcal{C}$ of finite relational structures satisfies the \emph{homo-amalgamation property} (HAP)
\index{homo-amalgamation property (HAP)}%
if:

for any $A, B_1, B_2 \in \mathcal{C}$, and any maps $f_i : A \to B_i$ (for $i = 1, 2$) such that $f_1$ is an embedding (an isomorphism to an induced substructure) and $f_2$ a homomorphism, there exists $C \in \mathcal{C}$ and homomorphisms $g_i : B_i \to C$ for $i = 1, 2$ such that $g_1 \circ f_1 = g_2 \circ f_2$ and $g_2$ is an embedding. 

The asymmetry between $B_1$ and $B_2$ is intentional.

The \emph{homo-extension property} of a structure $M$ with age  $Age(M)$ is defined thus:
\index{homo-extension property}%

if $B \in Age(M)$ and  $A$ is an induced substructure of $B$, then every homomorphism $A \to M$ extends to a homomorphism $B \to M$.

\begin{theorem}[Fra\"{\i}ss\'e's Theorem for HH Structures]
\label{fraiforhom}
\index{Fra\"{\i}ss\'e's Theorem}%

(a) A countable structure is HH if and only if it has the homo-extension property. 

(b) The age of any HH-structure has the homo-amalgamation property. 

(c) If a class $\mathcal{C}$ of finite relational structures
\index{relational structure}%
 is isomorphism-closed, closed under induced substructures, has only a countable number of isomorphism classes, and has the JEP and the HAP, then there is a countable HH structure $\mathcal{M}$ with age $\mathcal{C}$. 
\end{theorem}

\begin{proof}
(a) If $M$ has the homo-extension property, then any homomorphism from a finite substructure of $M$ can be extended point by point to a homomorphism of $M$.
 
Conversely, let $M$ is an HH structure. Let $B \in Age(M)$ and $A \subseteq B$; without loss of generality, $B \subseteq M$.  Then any homomorphism $f : A \to M$ extends to a homomorphism of $M$ whose restriction to $B$ is the required homomorphism $B \to M$.

(b) Suppose that $M$ is an HH structure and take $A, B_1, B_2, f_1, f_2$ to be as in the hypothesis of the HAP, with $A, B_1, B_2$ in $Age(M)$.  Assume without loss of generality that $B_1, B_2 \subseteq M$ and that $f_1$ restricted to $A$ is the identity.  Extend $f_1$ to a homomorphism $g$ of $M$; let $C = B_1 g$, $g_1 = g_{| B_1}$, and $g_2 = id_{| B_2}$.

(c)  Build the countable HH structure $M$ iteratively, supposing that $M_i \subset M$ has been constructed at step $i$.
 
At even $i$, use JEP to find structure $M_{i+1}$ such that $A, M_i \subseteq M_{i+1}$, for all $A \in \mathcal{C}$.

At odd $i$, choose $A, B \in \mathcal{C}$ with $A \subseteq B$.  Apply MAP to extend each homomorphism $A \to M$ to a homomorphism $B \to M'$ for some $M' \supseteq M$.  Successive application to each homomorphism $A \to M_i$ gives us a structure $M_{i+1}$ such that each homomorphism $A \to M_{i+1}$ extends to a homomorphism $B \to M_{i+1}$.  Arranging the steps so that every structure in $\mathcal{C}$ occurs at some even stage, and every pair $(A, B)$ at infinitely many odd stages, we arrive at a countable structure with age $\mathcal{C}$ and having the HH property, and hence is HH.
\end{proof}

Call a class $\mathcal{C}$ satisfying part (c) of the theorem, a \emph{homo-Fra\"{\i}ss\'e class},
\index{homo-Fra\"{\i}ss\'e class}%
and an HH structure with age $\mathcal{C}$ a \emph{homo-limit} of $\mathcal{C}$.  Unlike the usual form of Fra\"{\i}ss\'e's Theorem,
\index{Fra\"{\i}ss\'e's Theorem}%
 a class having the homo-amalgamation property does \emph{not} have a unique homo-limit, up to isomorphism; (the same is true of mono-amalgamation and mono-limit~\cite{camnes}).  For example, there are many examples of graphs containing $\mathfrak{R}$ as a spanning subgraph whose age is the class of all finite graphs.  However there is an equivalence relation which replaces isomorphism in our version of Fra\"{\i}ss\'e's Theorem.  We call two structures $M$ and $M'$ \emph{homo-equivalent}
\index{structure ! homo-equivalent}%
if
\begin{itemize}
\item $Age(M) = Age(M')$;
\item every embedding of a finite substructure $A$ of $M$ into $M'$ extends to a homomorphism from $M$ to $M'$, and the same with $M$ and $M'$ reversed.
\end{itemize}

\begin{proposition}
\label{eqheqahh}
(a) If $M$ and $M'$ are homo-equivalent structures and $M$ is an HH structure, then $M'$ is an HH structure.

(b)  Conversely, if $M$ and $M'$ are HH structures with $Age(M) = Age(M')$, then they are homo-equivalent.
\end{proposition}
\begin{proof}
Assume the hypotheses and take $A, B \in Age(M)$ with $A \subseteq B$, and let $f: A \to M'$ be a homomorphism.  Assume that $B \subseteq M$.  Let $A'$ be the image of $f$.  Since $Age(M) = Age(M')$, there is a copy $A''$ in $M$, that is exists a homomorphism $\phi : A \to A''$ and an isomorphism $g : A'' \to A'$ such that $g \circ \phi = f$.

Since $M$ has the HH property, $\phi$ extends to a homomorphism $\phi^* : M \to M$.  Let $B'' = \phi^*(B)$.  Also by assumption, the isomorphism $g$ extends to a homomorphism $g^* : M \to M'$.  Then $B_{| g^* \circ \phi^*}$ is a homomorphism $B \to M'$ extending $f$.  So $M'$ has the homo-extension property, and so is an HH structure.

(b)  Suppose that $M$ and $M'$ are HH structures with the same age.  If $A$ is a finite substructure of $M$ and $f : A \to M'$ is an embedding, for all $B \supseteq A$, by the homo-extension property in $M'$ we can extend $f$ to a homomorphism $B \to M'$.  So there is a homomorphism $M \to M'$ extending $f$.
\end{proof}

(This proposition is a generalization of the result that a countable first-order structure
\index{first-order structure}%
 $\mathcal{M}$ is $\aleph_0$-categorical
\index{aleph@$\aleph_0$-categorical}%
if and only if $\Aut(\mathcal{M})$ is oligomorphic~\cite[p.30]{cam6}.)

There is a partial order on the set of equivalence classes of an equivalence relation $M \leq_p M'$ between structures $M$ and $M'$ which holds if
\begin{itemize}
\item $Age(M') \subseteq Age(M)$; and
\item any embedding of a finite substructure $A$ into $M'$ extends to a homomorphism from $M \to M'$.
\end{itemize}

Note the reverse ordering of ages; for example in the case of graphs containing $\mathfrak{R}$ as a spanning subgraph the more extra edges are added, the smaller the age.  Part (a) of the last proposition shows that if $M$ is HH and $M \leq_p M'$ then $M'$ is also an HH structure.  

Finally, we mention that the following is proved in~\cite{camnes}

\begin{proposition} A countable graph $\Gamma$ satisfies $\mathfrak{R} \leq_p \Gamma$ if and only if $\mathfrak{R}$ is a spanning
\index{graph ! spanning}%
 subgraph of $\Gamma$.
\end{proposition}

If $M$ and $M'$ are HH structures and $Age(M) \supseteq Age(M')$, it does not follow that $M \leq_p M'$.  For example, if $M = \mathfrak{R}$ and $M'$ is the disjoint union of two infinite complete graphs; the map taking a non-edge in $M$ to a non-edge in $M'$ clearly cannot be extended. 

\section{Topological HH Monoids}

Oligomorphic permutation groups~\cite{cam6}
\index{group ! permutation ! oligomorphic}%
 were found to be a particularly interesting variety of infinite symmetric groups,
\index{group ! symmetric ! infinite}%
 due to the relative tameness of their properties as well as the connections with other fields, especially model theory.
\index{model theory}%
  One such feature is that the group can be studied in the form of its closure in the symmetric group in the topology of pointwise convergence.
\index{topology ! of pointwise convergence}%
 Thus we turn to characterizing closed (and in the next section oligomorphic) infinite monoids as a start to a more comprehensive study which could include finding monoids obeying primitivity and homogeneity properties and uncovering theorems on submonoids of closed monoids.  It may even transpire that certain naturally defined oligomorphic monoids have links with combinatorial enumeration problems via interesting integer sequences that count orbits of monoid functions, as has arisen in the study of inverse semigroups
\index{semigroup}%
  of partial bijections on a finite set, or in oligomorphic group theory.

A way of representing closed monoids was given in~\cite{camnes}, and we repeat it here.

For a countable set $X$, there is a natural topology
\index{topology ! natural}%
 on $X^X$, namely the product 
topology induced from the discrete topology
\index{topology ! discrete}%
 on $X$ . Thus the basic open sets
\index{sets ! basic open}%
 are of the form 
\[ \{ f \in X^X : f(x_i) = y_i\ for\ i = 1, \ldots , n \},\quad\quad\quad (\dagger)\] 
where $x_1, \ldots , x_n , y_1 , \ldots , y_n \in X$ and $x_1, \ldots, x_n$ are distinct.  Recall that in 
the induced topology on the symmetric group $\Sym(X)$, a permutation group $G$ is 
closed if and only if it is the automorphism group
\index{group ! automorphism}%
 of a homogeneous relational structure
\index{relational structure}%
  on X (see~\cite{cam6}). A similar observation holds for monoids:

\begin{proposition} 
\label{hhclosed}
(a) A submonoid $S$ of $X^X$ is closed in the product topology on 
$X^X$ if and only if $S$ is the monoid $End(\mathcal{M})$ of endomorphisms of an HH relational structure $\mathcal{M}$ on $X$. 

(b) A submonoid $S$ of the monoid of one-to-one maps $X \to X$ is closed in the 
product topology if and only if $S$ is the monoid of monomorphisms of an MM 
relational structure $\mathcal{M}$ on $X$. 
\end{proposition}
\index{relational structure}%

\begin{proof}
(a) For each $n$, and each $\overline{x} \in X^n$, we take an $n$-ary relation $R_{\overline{x}}$ defined by 
\[ R_{\overline{x}} (\overline{y}) \Leftrightarrow (\exists s \in S)(\overline{y} = s(\overline{x})). \] 
Let $\mathcal{M}$ be the relational structure with relations $R_{\overline{x}}$ for all $n$-tuples $\overline{x}$ (and all $n$).  We claim that $S$ acts as endomorphisms of $\mathcal{M}$, that $\mathcal{M}$ is HH, and $\End(\mathcal{M}) = S$. 

For the first point, take $s \in S$ and $\overline{y} \in X^n$ such that $R_{\overline{x}} (\overline{y})$ holds; we must show that $R_{\overline{x}} (s(\overline{y}))$ holds. But $\overline{y} = s' (\overline{x})$ for some $s' \in S$; then $s(\overline{y}) = ss' (\overline{x})$, so the assertion is true.
 
Next, let $f$ be a homomorphism between finite subsets of $X$, say $f(x_i) = y_i$ for $i = 1, \ldots , n$. Let $\overline{x} = (x_1, \ldots , x_n)$ and $\overline{y} = (y_1, \ldots, y_n)$. Now $S$ is a monoid and so contains the identity mapping. Thus, by definition, $R_{\overline{x}} (\overline{x})$ holds. Since $f$ is a homomorphism, $R_{\overline{x}} (\overline{y})$ holds. So by definition, there exists $s \in S$ such that $s(\overline{x}) = \overline{y}$.  Now $s$ is an endomorphism of $\mathcal{M}$ extending $f$. So $\mathcal{M}$ is HH.
 
Finally, to show that $\End(\mathcal{M}) = S$, we know already that $S \subseteq \End(\mathcal{M})$ and have to prove the reverse inclusion. We must take $h \in \End(\mathcal{M})$ and show that every basic neighbourhood of $h$ contains an element of $S$, so that $h$ is a limit point of $S$.  Since $S$ is assumed closed, we conclude that $h \in S$. Now each $n$-tuple $\overline{x}$ defines a basic neighbourhood of $h$, consisting of all functions $k$ such that $k(\overline{x}) = h(\overline{x})$. Now $R_{\overline{x}} (\overline{x})$ holds; since $h$ is a homomorphism, $R_{\overline{x}} (h(\overline{x}))$ also holds, and by definition of $R_{\overline{x}}$ this means that there exists $s \in S$ with $h(\overline{x}) = s(\overline{x})$, as required.
 
(b) The proof of this is entirely analogous, replacing homomorphisms by monomorphisms.

\end{proof}

However the relational structures
\index{relational structure}%
 constructed in the proof of this proposition have infinitely many relations of 
each arity.  We seek criteria that would enable us to recognize the monoids which are the endomorphism (or monomorphism) monoids of homogeneous structures with only finitely many relations of each arity (these would be the analogue of the closed oligomorphic permutation groups~\cite{cam6}),
\index{group ! permutation ! oligomorphic}%
 or even those with only finitely many relations altogether.  In the following section we will return to this topic.

\bigskip

Consider the following characterisation of topological groups~\cite{kayem}, repeated as Theorem~\ref{topaut} in Appendix~\ref{TopologyinPermutationGroups}. 

\begin{theorem}
\label{topaut}
\index{topological group ! Hausdorff}%
Let $G$ be a topological group.  Then $G = \Aut(M)$ for some countable
structure $M$ if and only if 
\begin{itemize}
\item[(a)] $G$ is Hausdorff;
\item[(b)] $G$ is complete;
\item[(c)] if $H < G$ is open then $|G : H| \le \aleph_0$;
\item[(d)] there is a countable family $\{H_i : i \in \mathbb{N}\}$ of
  subgroups of $G$ such that $\mathcal{B} = \{H_{0}^{g_0} \cap \ldots \cap
  H_{k}^{g_k} : k \in \mathbb{N}, g_0, \ldots, g_k \in G\}$ is a base
  of open subgroups of $G$, i.e. the set of cosets of elements of
  $\mathcal{B}$ forms a base for the topology on $G$.
\end{itemize}
\end{theorem}

That is, a topological group is the automorphism group
\index{group ! automorphism}%
 of a countable structure provided that the topology is complete and Hausdorff,
\index{topology ! Hausdorff}%
  and the family of open neighbourhoods of the identity has a countable base of open subgroups, each of which has index at most $\aleph_0$ in the automorphism group.  

\head{Open Question}  Much of group theory does not generalize in a straightforward way to monoids, but is there a analogue of this nice theorem characterizing topological monoids?
\index{topological ! monoid}%

\smallskip


We can rule out one potential condition immediately; we do not want any compactness condition on the topology because then the topology would be locally compact,
\index{topology ! locally compact}%
 which itself is sufficient to rule out oligomorphic structures by the Engeler-Ryll-Nardzewski-Svenonius Theorem.
\index{Engeler--Ryll-Nardzewski--Svenonius ! Theorem}%
  Since compact implies complete implies closed, it seems reasonable to consider completeness.  But in a complete metric space a subspace is complete if and only it it is closed.  Further, if $\mathcal{M}$ is a first-order structure then $\End(\mathcal{M})$ is closed; the converse is Proposition~\ref{hhclosed}, that is, a submonoid of $X^X$ is closed if and only if it is the endomorphism monoid of a first-order structure.

Take as basic open sets
\index{sets ! basic open}%
 of the topology, translates of submonoids (that is composition with elements of the monoid) generated by the functions $f$ given by $(\dagger)$ above.  Each finite $n$-tuple of elements of the countable set $X$ defines a basic neighbourhood of the function $f \in X^X$.  Open sets
\index{open set}%
  are countable unions of these basic open sets.

The natural metric derived from the permutation topology on the symmetric group, that of pointwise convergence,
\index{topology ! of pointwise convergence}%
 is simpler in the case of monoids, as the metric derived from the topology is complete even in the absence of inverses.  Suppose the operand is $\mathbb{N}$.  A countable topology is metrizable via a distance function, so a sequence ($m_n$) of monoid maps where $m_n = f_0\ f_1\ f_2\ldots f_{n-1}$ is a Cauchy sequence
\index{Cauchy sequence}%
 for the metric $d$ given by
\begin{displaymath}
d(m, m') = \left\{ \begin{array}{ll}
\text{0} & \text{if}\ m = m';\\
\text{$\frac{1}{2^i}$} & \text{if}\ m\ \&\ m'\ \text{agree on } \{0, \ldots, i-1\} \text{ and disagree on}\ i
 \end{array} \right.
\end{displaymath}
which begets the topology.  In the case of permutation groups, this sequence would still be a Cauchy sequence
\index{Cauchy sequence}%
 but one that fails to converge, as illustrated by noting that the pointwise limit of permutations $\pi_n = (0 1 \ldots n-1)$ is not a permutation, and so an extra relation for inverses is required~\cite{cam6}.  However in the case of monoids, the limit still lies in the monoid.  (Compare this with the fact~\cite{alcoma} that the outer automorphism group
\index{group ! automorphism}%
  of the quotient $\Sym(\Omega) / \FSym(\Omega)$ on a countably infinite set $\Omega$ is $\mathbb{Z}$).  Any metric space is Hausdorff.  
\index{Hausdorff space}%

In the case of topological groups, each open subgroup is also closed, being the complement of a union of cosets of basic open sets;
\index{sets ! basic open}%
 then choosing different basic open sets gives different topologies.  If we are to assume the same applies in the case of monoids then by Proposition~\ref{hhclosed} we are forced to choose endomorphism monoids of HH relational structures. 
\index{relational structure}%


For groups, the existence of an inverse ensures that cosets by a subgroup are disjoint and that there is a well-defined notion of index.  For monoids, we can define $|M : S|$ to mean that $M$ is the union of countably many right translates of $S$ but the cosets are not disjoint.

\bigskip

As a prelude to the next section we make the following comments.  

We know from the Upward L\"{o}wenheim--Skolem Theorem
\index{Upward L\"{o}wenheim--Skolem ! Theorem}%
\index{Skolem, T.}%
\index{L\"{o}wenheim, L.}%
 that if $\mathcal{M}$ is infinite, then its complete theory does not determine $\mathcal{M}$.  (However, the automorphism groups
\index{group ! automorphism}%
 of saturated structures can then be considered as invariants of the theory up to isomorphism; they are equal for bi-interpretable theories,
\index{theory ! bi-interpretable}%
 see Appendix~\ref{TopologyinPermutationGroups}).   The next best option is then to determine its model up to its cardinality.  Given this, a complete theory is $\aleph_0$-categorical 
\index{aleph@$\aleph_0$-categorical}%
if it has a unique countable model up to isomorphism.  $\aleph_0$-categoricity is proved using the back-and-forth method,
\index{back-and-forth method}%
 and it is not clear that when we extend the theory from the actions of groups to those of monoids that we will be able to use any more than forth.  

The back-and-forth argument, which can be traced back to Huntingdon~\cite{huntingdon}
\index{Huntingdon, E. V.}%
 and then the book of Hausdorff~\cite{hausdorff}
\index{Hausdorff, F.}%
 has been investigated, in particular with respect to the question of the sufficiency of going forth. 
\index{forth}%
  Generally, back-and-forth works when we begin with an equivalence relation (or `type')
\index{type}%
  on $n$-tuples ($\forall n \in \mathbb{N}$), possibly chosen from different structures in some class, satisfying
  
 $(\Diamond)$ If $\tp(\bar{a}) = \tp(\bar{b})$ and $x$ is any point of the structure $\mathcal{M}$ containing $\bar{a}$, then there is a point $y$ in the structure $\mathcal{N}$ containing $\bar{b}$ such that $\tp(\bar{a}, x) = \tp(\bar{b}, y)$.
 
 This guarantees both that all structures in the class are isomorphic and that tuples in a structure have the same type if and only if there is a type-preserving permutation of $\mathcal{M}$ which carries one to the other.  Then `types' are just orbits of a closed permutation group $G$ of $\mathcal{M}$,
\index{group ! permutation ! closed}%
  and we can dispense with structure and logic and define types to be orbits of a closed permutation group on tuples of a set $X$, that is on $\bigcup_{n \geq1} X^n$.  These are the isomorphism or first-order types, and if a structure needs to be taken into account, they lie in the canonical structure for this group, satisfying $(\Diamond)$.
  
 Say that \emph{forth suffices} for $G$ if the one-to-one type-preserving map from $X \to X$ is always onto, irrespective of the chosen enumerations, and that \emph{forth suffices} for $\mathcal{M}$ if it suffices for its automorphism group.  

Extend the usual definition of \emph{suborbit} of $G$
 \index{group ! permutation ! suborbit}%
from finite permutation group theory to be a pair $(\bar{a}, A)$, where $A$ is an orbit of the stabilizer of tuple $\bar{a}$.  If $G = \Aut(\mathcal{M})$ and $\mathcal{M}$ is $\aleph_0$-categorical,
\index{aleph@$\aleph_0$-categorical}%
 then $A$ is a minimal $\bar{a}$-definable set in $\mathcal{M}$.)  The suborbit $(\bar{a}, A)$ is \emph{trivial} if $A$ is a singleton which is a member of $\bar{a}$.  Say that $(\bar{b}, B)$ \emph{dominates} $(\bar{a}, A)$ if 
 
 (i)  each element of $\bar{a}$ occurs in $\bar{b}$;
 
 (ii)  $B \subseteq A$.
 
 The suborbit $(\bar{a}, A)$ is \emph{splittable} if there is a tuple $\bar{b}$ of points outside $A$ such that $G_{\bar{a}\bar{b}}$ is intransitive on $A$, and \emph{unsplittable} otherwise.  Roughly, forth suffices if there are many unsplittable suborbits and vice versa.  
 
 If every non-trivial suborbit dominating a suborbit $(\bar{a}, A)$ is splittable then forth does not suffice.   

Let $G$ be a permutation group acting on a set $X$.  The \emph{algebraic closure}
\index{algebraic closure}%
  $\acl(A)$ of a finite subset $A \subset X$ is the set of all those points of $X$ which lie in finite orbits of the pointwise stabilizer of $A$.  If for every finite set, $\acl(X) = X$ and $G$ has a primitive splittable suborbit then forth does not suffice.  The random graph
\index{graph ! random}%
 and Henson's graph $\mathfrak{H}_k$
\index{graph ! Henson}%
 are examples that fit this result.

Another example of an $\aleph_0$-categorical,
\index{aleph@$\aleph_0$-categorical}%
 homogeneous structure for which forth does not suffice is a countable dense ordered set without endpoints and with a distinguished dense subset whose complement is also dense.  This structure is naturally associated with the stabilizer of a point in the countable homogeneous local order.

For a closed and countable permutation group on a set $X$, forth suffices~\cite[p.129]{cam6}.  The gap between the necessary condition and the sufficient condition for forth to suffice given in~\cite{cam6}, has been narrowed by McLeish~\cite{mcleish}.

Extended discussions of this topic and its strong connection with Jordan groups can be found in~\cite{cam6} and~\cite{kayem}.


\section{Oligomorphic Monoids}
\index{monoid ! oligomorphic}%

The characterization of $\aleph_0$-categorical structures as given by the Engeler--Ryll-Nardzewski--Svenonius Theorem~\cite{engeler}~\cite{ryll} \cite{svenonius}
\index{Engeler--Ryll-Nardzewski--Svenonius ! Theorem}%
\index{Engeler, E.}%
\index{Ryll-Nardzewski, C.}%
\index{Svenonius, L.}%
 would seem to be a one-off result of the form ``axiomatizability = symmetry'', applying to countable structures and not to those of higher cardinality.  However there is a generalisation.  

We begin the section with the work of D. Ma\v{s}ulovi\'{c} and M. Pech
\index{Ma\v{s}ulovi\'{c}, D.}%
\index{Pech, M.}%
 on homomorphism-homogeneous structures and oligomorphic transformation monoids~\cite{masulovic}.
\index{structure ! homomorphism-homogeneous}%
 We end with the theorem that Bodirsky and Pinsker
\index{Bodirsky, M.}%
\index{Pinsker, M.}%
 have produced in~\cite{bodirsky} to characterize the endomorphism monoid of an $\aleph_0$-categorical structure.
 
There are several potential definitions of an oligomorphic monoid resulting from the fact that the ``orbit relation'' for a monoid is a preorder rather than an equivalence relation: it is reflexive and transitive but not necessarily symmetric.  Given a preorder $\to$ on a set $X$, the relation $\equiv$ defined by $x \equiv y$ if $x \to y$ and $y \to x$ both hold is an equivalence relation; its equivalence classes are partially ordered by $\to$. Thus we have the following three conditions of increasing strength for a transformation monoid $S$ on $X$: 

(a) the orbit preorder on $X^n$ has only finitely many connected components, for all $n$;
 
(b) the equivalence relation obtained from the orbit preorder on $X^n$ has only finitely many classes, for all $n$;
 
(c) the invertible elements (permutations) in $S$ form an oligomorphic permutation group on $X$.
 
There is an implication from (b) to (a) because each connected component of the orbit preorder contains one or more equivalence classes.  The equivalence relation is given by $\vec{x} \equiv \vec{y}$ if and only if there exist $f, h \in S$ with $\vec{x} f = \vec{y}$ and $\vec{y} h = \vec{x}$; if $\vec{x} f = \vec{y}$ and $f$ is invertible, then taking $h = f^{-1}$ this is satisfied, so (c) implies (b). 

\medskip

D. Ma\v{s}ulovi\'{c} and M. Pech~\cite{masulovic}
\index{Ma\v{s}ulovi\'{c}, D.}%
\index{Pech, M.}%
 found the most natural definition, whereby homogeneneous or $\aleph_0$-categorical
\index{aleph@$\aleph_0$-categorical}%
structures have as their analogues the class of countable homomorphism-homogeneous 
\index{structure ! homomorphism-homogeneous}%
 or countable \emph{weakly oligomorphic structures},
\index{structure ! weakly oligomorphic}%
  to be defined below.

A structure is \emph{homomorphism-homogeneous}
\index{structure ! homomorphism-homogeneous}%
 if every homomorphism between finitely generated substructures of the structure extends to an endomorphism of the structure~\cite{camnes}.

Let $\mathcal{R} = (R_i)_{i \in I}$ be a relational signature, and for an $\mathcal{R}$-formula $\phi(x_1, \ldots, x_n)$ and let $\mathcal{M} = \langle M, (R_i)_{i \in I}^{\mathcal{M}} \rangle$ be an  
$\mathcal{R}$-structure.  

A transformation monoid $M \subseteq X^X$ is \emph{oligomorphic}
\index{monoid ! oligomorphic}%
 if the action of $M$ on $X^n$ has only finitely many orbits for every $n \in \mathbb{N}$.

Just as the automorphism group of a homogeneous structure over a finite relational signature is oligomorphic, if $\mathcal{M}$ is a homomorphism-homogeneous structure over a finite relational language, then $\End(\mathcal{M})$ is an oligomorphic transformation monoid.

The finiteness of the signature of $\mathcal{M}$ ensures that it has only finitely many $n$-element substructures for all $n \in \mathbb{N}$.  For each pair of integers $m, k$ such that $m \geq k$ and each surjection $f : \{ 1, \ldots, m\} \to \{1, \ldots, k\}$ fix a right inverse $f^*$ of $f$, that is, a mapping $f^* : \{1, \ldots, k\} \to \{ 1, \ldots, m\}$ satisfying $f \circ f^* = \id$.  For example take $f^* (y) = \min(x : f(x) = y)$.

A signature $\mathcal{R}$ is \emph{residually finite} in an $\mathcal{R}$-structure if $\mathcal{R}_k$ is finite for every $k\in\mathbb{N}$, and $\forall n \in \mathbb{N}\ \exists l > n$ such that $\forall m \geq l$, every $k \leq n$, every relation symbol $R \in \mathcal{R}_m$ and every surjective mapping $f : \{1, \ldots, m\} \to \{1, \ldots, k\}$ there is a relation symbol $R^f \in \mathcal{R}_k$ such that

\[ \mathcal{M} \models \forall x_1 \ldots x_m \Big( ( \bigwedge_{(i, j) \in \ker(f)} x_i = x_j ) \Rightarrow  \]
\[   ( R^f(x_{f^*(1)}, \ldots, x_{f^*(k)}) \Leftrightarrow R(x_1 \ldots x_k) ) \Big) \]

Every finite relational signature $\mathcal{R}$ is residually finite in every $\mathcal{R}$-structure.  Also if $\mathcal{R}$ is residually finite in some $\mathcal{R}$-structure then $\mathcal{R}$ is countable.  If $\mathcal{R}$ is a countable relational structure and $\mathcal{M}$ an $\mathcal{R}$-structure such that the signature $\mathcal{R}$ is residually finite in $\mathcal{M}$ then for every $n \in \mathbb{N}$ there are, up to isomprphism, only finitely many $n$-element substructures of $\mathcal{M}$.

\begin{theorem}[Ma\v{s}ulovi\'{c} and M. Pech]
\index{Ma\v{s}ulovi\'{c}, D.}%
\index{Pech, M.}%
Let $X$ be an infinite set and let $M \subseteq X^X$ be a transformation monoid.  The following are equivalent:

(1)  $M$ is closed in $X^X$ and oligomorphic;

(2)  There is a homomorphism-homogeneous 
\index{structure ! homomorphism-homogeneous}%
 $\mathcal{R}$-structure $\mathcal{M}$ on $X$ for which $M = \End(\mathcal{M})$, where $\mathcal{R}$ is a countable relational signature which is residually finite in $\mathcal{M}$.
\end{theorem}

A relational structure $\mathcal{M}$ is \emph{weakly oligomorphic},
\index{structure ! weakly oligomorphic}%
 if $\End(\mathcal{M})$ is an oligomorphic transformation monoid.  Every relational structure with an oligomorphic automorphism group is weakly oligomorphic.

\smallskip

In order to state the oligomorphic monoid theorem equivalent to that of  Engeler--Ryll-Nardzewski--Svenonius
\index{Engeler--Ryll-Nardzewski--Svenonius ! Theorem}%
we need more definitions.

Let $T$ be a first order theory in relational signature $\mathcal{R}$.  Formulas $\phi(x_1, \ldots, x_n)$ and $\psi(y_1, \ldots, y_n)$ are \emph{equivalent in $T$}, written $\phi \equiv_{T} \psi$, if $T \models \forall \textbf{x} (\phi(\textbf{x}) \Leftrightarrow \psi(\textbf{x}))$.  Theory $T$ is said to have the \emph{positive Ryll-Nardzewski property} if for each $n \in \mathbb{N}$ there are only finitely many positive $\mathcal{R}$-formulas in variables $x_1, \ldots, x_n$ which are pairwise inequivalent in $T$.  An $\mathcal{R}$-structure $\mathcal{M}$ has this property if its first-order theory has it.

The \emph{complete positive $n$-type} of an $n$-tuple $\textbf{a} \in M^n$ in an $\mathcal{R}$-structure $\mathcal{M}$ is the set of all positive formulas that are satisfied by $\textbf{a} \in \mathcal{M}$.  An $\mathcal{R}$-structure $\mathcal{M}$ \emph{realizes} a set of positive formulas $\Phi(\textbf{x})$ if there exists an $n$-tuple $\textbf{x} \in M^n$ such that $\Phi$ is a subset of the complete positive $n$-type satisfied by $\textbf{a}$.  A \emph{complete positive type of a complete theory $T$} is a complete positive type of a tuple of some model of $T$.  Denote the set of all complete positive $n$-types of a theory $T$ by $S_{n}^{+}(T)$.  

A set of formulas $\Phi(\textbf{x})$ is \emph{principal} with respect to $T$ if there is a formula $\psi(\textbf{x}) \in \Phi(\textbf{x})$ such that $T \models \forall \textbf{x} ( \psi(\textbf{x}) \Rightarrow \bigwedge \Phi(\textbf{x}))$.  We say that $\psi$ \emph{generates} $\Phi$.  A first-order formula $\chi(\textbf{x})$ is a \emph{characteristic formula} for a set $\Phi(\textbf{x})$ of positive formulas with respect to $T$ if $T \models \forall \textbf{x} (\chi(\textbf{x}) \Leftrightarrow \bigwedge \Phi(\textbf{x}) \wedge \bigwedge \Phi^{\neg}(\textbf{x}))$, where

 $\Phi^{\neg}(\textbf{x}) = \{ \neg \psi(\textbf{x}) \text{is a positive formula and}\ \psi(\textbf{x}) \notin \Phi(\textbf{x}) \}$.

\begin{theorem}[Ma\v{s}ulovi\'{c} and M. Pech]
\index{Ma\v{s}ulovi\'{c}, D.}%
\index{Pech, M.}%
Let $T$ be a complete theory over a signature $\mathcal{R}$ and assume that $T$ has infinite models.  The following are equivalent:

(1)  Every countable model of $T$ is weakly oligomorphic.

(2)  There exists a countable model of $T$ which is weakly oligomorphic.

(3)  $T$ has the positive Ryll-Nardzewski property.

(4)  For every $n \geq 1$, every type in $S^{+}_n(T)$ has a characteristic formula.

(5)  For every $n \geq 1$, every type in $S^{+}_n(T)$ is finite.
\end{theorem}

Let $\mathcal{A}$ be an $\mathcal{R}$-structure and let $\textbf{a} \in A^n$.  The \emph{complete positive quantifier-free type} of $\textbf{a} \in \mathcal{A}$ is the set of positive quantifier-free formulas that are satisfied by $\textbf{a}$ in $\mathcal{A}$.

Ma\v{s}ulovi\'{c} and M. Pech
\index{Ma\v{s}ulovi\'{c}, D.}%
\index{Pech, M.}%
 also prove that an $\aleph_0$-categorical
\index{aleph@$\aleph_0$-categorical}%
 structure is both homogeneous and homomorphism-homogeneous if and only if it has quantifier elimination
\index{quantifier elimination (q.e.)}%
 where every positive formula reduces to a positive quantifier-free formula.

\medskip

We make some comments pertaining to random graphs.  Let $T$ be the first-order theory of graphs with the property that any finite set of vertices has a common neighbour. This requirement can be stated as
a countable sequence of $(\forall\exists)$-sentences.

Firstly, any two countable models of $T$ are mono-equivalent, for if $\mathcal{M}$ and $\mathcal{N}$ are
countable models, then ``forth'' gives a monomorphism from $\mathcal{M}$ into $\mathcal{N}$.  Secondly, not every graph mono-equivalent to a model of $T$ is a model of $T$; take, for example, the complete graph with an isolated vertex.

The countable models of $T$ are precisely the graphs containing the random graph
\index{graph ! random}%
 as a spanning
\index{graph ! spanning}%
 subgraph; so they are all HH and MM. However, these are not ``pure'' monoid properties since they also involve the graphs.

The theory of graphs in which every finite set of vertices has a common neighbour has the property that all its countable models are \emph{hom-equivalent} (that is, there exist homomorphisms in both directions between any two).  Let $T'$ be the first-order theory of graphs with the property that any finite clique has a common neighbour.  Both $T$ and $T'$ are examples of first-order theories which are not $\aleph_0$-categorical but for which all countable models are hom-equivalent.  Theory $T$ is an example of a first-order theory whose countable models are hom-equivalent but don't form a hom-equivalence class.


\smallskip

We turn to the general approach of Bodirsky and Pinsker.
\index{Bodirsky, M.}%
\index{Pinsker, M.}%

An operation $f : V^n \to V$ is a \emph{projection} if and only if there exists an $i \leq n$ such that $f(x_1, \ldots, x_n) = x_i$ for all $x_1, \ldots, x_n \in V$.  A \emph{clone}
\index{clone}%
 is a subset of the set of all finitary functions over a finite or infinite domain, which is closed under composition and containing all finitary projections - see also Chapter~\ref{chapFD} and Appendix~\ref{FurtherDetails}.  

The topologically closed permutation groups,
\index{group ! permutation ! closed}%
 the closed transformation monoids, and the closed clones,
\index{clone}%
 containing $\Aut(\mathfrak{R})$ all form complete lattices.  The atoms of these lattices are generated by functions on the random graph that are in a certain sense minimal.  An operation is \emph{minimal} if and only if it is non-trivial and all non-trivial functions $g$ it generates have at least the arity of $f$ and generate $f$.  The trivial functions are the automorphisms of $\mathfrak{R}$ possibly with additional dummy variables.  The group lattice has 5 atoms, corresponding to the 5 reducts of $\mathfrak{R}$ in Thomas' Theorem,
 \index{Thomas' Theorem}%
 whilst the other two lattices have infinite height.  
 
The lattice of reducts modulo the equivalence whereby two reducts are equivalent if and only if they first-order define one another, is anti-isomorphic to the lattice of closed permutation groups containing $\Aut(\mathfrak{R})$.  The finer lattice of reducts up to \emph{existential positive interdefinability} (see below), corresponds to the lattice of closed transformation monoids containing $\Aut(\mathfrak{R})$.  The lattice of reducts modulo the still finer equivalence of \emph{primitive positive interdefinability} (see below), corresponds to the lattice of closed clones containing $\Aut(\mathfrak{R})$. 
 
The set of all clones over a domain forms a complete lattice with respect to set-theoretical inclusion, and also forms a smaller lattice of clones which are closed in the natural topology
\index{topology ! natural}%
 of pointwise convergence 
\index{topology ! of pointwise convergence}%
 on the set $\mathcal{O}$ of all finitary operations on the domain closed under compositions and containing the projections.  The countable basis of sets take the form 
\[\mathcal{O}^{s}_{A} : = \{ f\ \epsilon\ \mathcal{O} : f_{| A} = s \}, \]
where for domain $D$, $A \subseteq D^n$ is finite and $s : A \to D$ is a finite function.  The universal algebraic name for clones which are closed subsets of $\mathcal{O}$ in this topology are called \emph{locally closed}
\index{clone ! locally closed (local)}%
 or \emph{local}.  This is the equivalent for clones of closure for permutation groups.  
\index{group ! permutation ! closed}%

A clone containing $\Aut(\mathfrak{R})$ is an atom in the lattice of local clones containing $\Aut(\mathfrak{R})$ if and only if there exists a minimal operation $f$ (defined in~\cite{bodirsky}) on $\mathfrak{R}$ such that the clone is the smallest local clone containing the set $\{ f \} \cup \Aut(\mathfrak{R})$, described as the local clone generated by $\{ f \}$ over $\Aut(\mathfrak{R})$.

Two structures $M_1$ and $M_2$ are \emph{first-order interdefinable} when each is first-order definable in the other.  Let $f : D^n \to D$ be an operation and let $R \subseteq D^m$ be a relation.  We say that $f$ \emph{preserves} $R$ if and only if $f(r_1, \ldots, r_n) \in R$ whenever $r_1, \ldots, r_n \in R$, where $f(r_1, \ldots, r_n)$ is calculated componentwise.   One statement of the Engeler--Ryll-Nardzewski--Svenonius Theorem,
\index{Engeler--Ryll-Nardzewski--Svenonius ! Theorem}%
 is that a relation $R$ is first-order definable in an $\aleph_0$-categorical structure
\index{aleph@$\aleph_0$-categorical}%
 $M$ if and only if $R$ is preserved by all automorphisms of $M$.  It follows that the reducts of an $\aleph_0$-categorical structure are in $1$--$1$ correspondence with the locally closed permutation groups containing $\Aut(M)$.

We need to define some concepts.  A \emph{polymorphism}
\index{polymorphism}%
 of $M$ is a homomorphism from a finite power $M^n$ to $M$, or simply a finitary operation preserving all relations of $M$.  A first-order formula is \emph{existential} if and only if it is of the form $\exists x_1, \ldots, x_k \psi$, where $\psi$ is quantifier-free.  A first-order formula is \emph{primitive positive} (respectively \emph{existential positive})
\index{formula ! primitive positive}%
\index{formula ! existential positive}%
 if and only if it contains no negations, disjunctions and universal quantifications (and it is existential).    As generalizations of the way that first-order definability characterizes automorphisms of an $\aleph_0$-categorical structure $M$,
\index{aleph@$\aleph_0$-categorical}%
 self-embeddings of an $\aleph_0$-categorical structure
 (respectively endomorphisms / polymorphisms) of $M$ can be used to characterize existential definability (respectively existential positive definability / primitive positive definability).  
 
Two structures  $M_1$ and $M_2$ are \emph{primitive positive interdefinable}
\index{structure ! primitive positive interdefinable}%
 if and only if every relation in $M_1$ has a definition by a primitive positive formula
\index{formula ! primitive positive}%
  in $M_2$ and vice versa; analogously for existential positive and existential interdefinability.  The following result is from~\cite{bodirskyn}
 
\begin{theorem}
A relation $R$ is primitive positive definable in an $\aleph_0$-categorical structure $M$ if and only if $R$ is preserved by the polymorphisms of $M$.
\end{theorem}

Finally we come to the theorem that describes existential and existential positive definability in an $\aleph_0$-categorical structure in terms of its endomorphism monoid.

\begin{theorem}[Bodirsky--Pinsker Theorem]
A relation $R$ has an existential positive (existential) definition in an $\aleph_0$-categorical structure if and only if $R$ is preserved by the endomorphisms (self-embeddings) of the structure.
\end{theorem}
\index{Bodirsky-Pinsker Theorem}%

From this it follows that the study of closed transformation monoids containing $\End(\mathfrak{R})$ is equivalent to the study of reducts of $\mathfrak{R}$ up to existential positive interdefinability.

\head{Open Question} Classify the locally closed transformation monoids that contain $\Aut(\mathfrak{R})$.

In~\cite{bodirskyj}, Bodirsky and Junker
\index{Bodirsky, M.}%
\index{Junker, M.}%
 extend the Ahlbrandt--Ziegler
\index{Ahlbrandt, G.}%
\index{Ziegler, M.}%
 analysis~\cite{ahlbrandtz} of interpretability in $\aleph_0$-categorical structures
\index{aleph@$\aleph_0$-categorical}%
 by showing that the way that interpretability is controlled by the automorphism group has analogies, namely in the way that existential interpretation is controlled by the monoid of self-embeddings and positive existential interpretation of structures without constant endomorphisms is controlled by the monoid of endomorphisms.

\head{Open Question} F. Point and M. Prest
\index{Point, F.}%
\index{Prest, M.}%
 proved in a study of the first-order theory of $\mathcal{R}$-modules over a ring $\mathcal{R}$~\cite[Lemma~1.1]{pointprest}, that there are forms of Morita equivalence
\index{Morita equivalence}%
  that are special cases of the concept of interpretations in model theory;
\index{model theory}%
 two generalizations are in~\cite{prest1} and~\cite{prest2}.  Is there a version of the theorem of Bodirsky and Pinsker
\index{Bodirsky, M.}%
\index{Pinsker, M.}%
 in this direction?

\section{Appendix: Previous Results on Monoids Supported by $\mathfrak{R}$}

We end the chapter by referencing some other publications that combine monoids and $\mathfrak{R}$.  

Part of the interest has been in designing random graph models for the (world wide) web graph.

The \emph{infinite locally random graph}
\index{graph ! random ! infinite locally}%
is the graph limit obtained from iterating the following construction: given a graph $\Gamma$, for each vertex $x$ and each subset $X$ of its closed neighborhood, add a new vertex $y$ whose neighbors are exactly $X$.  Charbit and Scott have classified the infinitely many isomorphism classes of limit graph~\cite{charbit}.

A semigroup
\index{semigroup}%
 $S$ is \emph{regular} if for all $f \in S$, there is a $g \in S$ such that $fgf = f$ and $gfg = g$.  In~\cite{bondel}, Bonato and Deli\'c, investigated $\End(\mathfrak{R})$, showing that it is neither regular nor satisfies the property that: for each \emph{singular} (that is, not onto) $f \in \End(\Gamma)$, $\exists n \geq 1$ such that $f$ is the product of $n$ idempotents, denoted $E(\End(\Gamma))$, of $\End(\Gamma)$.  They also show that $\mathbb{Q}$ and in particular, every countable linear order,
\index{linear order}%
 embeds in the set of idempotents in $\End(\mathfrak{R})$.  Another result of this paper~\cite[Proposition.~4.2]{bondel} is that for a graph $\Gamma$, there exists a retraction $f$ of $\mathfrak{R}$ such that $\im(e) \cong \Gamma$ if and only if $\Gamma$ is algebraically closed; they call a graph $\Gamma$ \emph{algebraically closed}
\index{graph ! algebraically closed}%
 if for each finite $S \subseteq V(\Gamma)$, there is a vertex $v \in V(\Gamma) \backslash S$ joined to all vertices from $S$.  Here, $\im(e)$ is the induced subgraph of $\mathfrak{R}$ on the vertex set of $e(V(\mathfrak{R}))$.

In~\cite{deldol}, Deli\'c and Dolinka prove that $\End(\mathfrak{R})$ is not simple, but rather has uncountably many ideals.

The relation $f \leq g$ defined by $fg = gf =f$ is an order on $E(\End(\Gamma))$.  Bonato proves~\cite{bona} that this order is universal, that is it embeds every countable order, and by a refinement, every countable preorder also.  This latter is constructed as follows: if $(A, \leq)$ is a preorder, then defining $x \sim y$ if $x \leq y$ and $y \leq x$, and writing the set of $\sim$-classes as $A / \sim$, gives an order $(A / \sim, \leq)$.

In~\cite{bonatopralat}, Bonato and Pralat define a core $H = Core(\Gamma(n, p(n)))$ of a random graph $\Gamma(n, p(n))$ \emph{great}
\index{graph ! core ! great}%
if for all $e \in E(H)$, there is a homomorphism from $H \backslash e$ to $H$ that is not onto.  For a large range of $p$ they prove that with probability tending to 1 as $n \to \infty$, $\Gamma \in \Gamma(n, p(n))$ is a core that is not great.

The proof that $\End(\mathfrak{R})$ is universal as a monoid~\cite{bondeldo} utilizes a chain of graphs whose union is $\mathfrak{R}$.  This extends to a chain of endomorphism monoids on those graphs that embeds $\End(\mathfrak{R})$, that is universally generates the full transformation monoid on a countable set.

Finally, Bonato
\index{Bonato, A.}
 proved~\cite{bona1} that for each finite core graph $\Gamma$, the class of all graphs admitting a homomorphism into $\Gamma$ is a \emph{pseudo-amalgamation class},
\index{amalgamation property ! pseudo-}%
(see Appendix~\ref{TheoryofRelationalStructures}).

\chapter{Random Graph Constructions} 
\label{rgconstr}

\bigskip

.... the sole end of science is the honor of the human mind, and that under this title a question about numbers is worth as much as a question about the system of the world.
\begin{flushright}
Carl Jacobi, \textit{Quoted in N Rose Mathematical Maxims and Minims (Raleigh N C 1988)}
\end{flushright} 

\medskip

The cult of Vishnu has some frail links with Vedic mythology.  He
there appears as a solar god who traverses the three worlds in three steps.
\begin{flushright}
\textit{New Larousse Encyclopedia of Mythology}
\end{flushright}

\medskip

The connections between different parts of number theory and (random)
graph theory in all likelihood lie in fertile ground (see for example Appendix~\ref{CategoryandMeasure}.  The first
two sections of this chapter illustrate the use of reciprocity
theorems and merely embellish the construction of $\mathfrak{R}$ given in~\cite{cameron}.  This is followed by observations on constructions from elsewhere in mathematics.

\section{A Number-Theoretic Construction of $\mathfrak{R}$}

A particularly concise presentation of $\mathfrak{R}$ is the following~\cite{bollobas1}
\[ \mathfrak{R} \cong \langle \omega, \{(i, j) : \text{the \emph{i}th prime divides \emph{j}, or vice versa}\} \rangle.\]

In~\cite{cameron} a construction of $\mathfrak{R}$ is given using the
original form of the quadratic reciprocity law
\index{quadratic reciprocity}%
due to Gauss.  Here we give an original generalization showing that the construction works for a
Gaussian algebraic number field $\mathbb{Q}[i]$.

Now to construct $\mathfrak{R}$.  Take as vertex set $\mathbb{P}$ the
set of all odd primes in $\mathbb{Z}$.  Join prime $p$ to prime $q$
with an edge if $\left( \frac{p}{q} \right) = 1$.  If $\mathbb{P}_1 =
\{p \equiv 1 (4)\}$ and $\mathbb{P}_{-1} = \{ p \equiv -1 (4)\}$, this
would give the structure $\mathbb{P}_1 \cup \mathbb{P}_{-1}$ where
$\mathbb{P}_1$ induces the random graph,
\index{graph ! random}%
 $\mathbb{P}_{-1}$ induces the random tournament
\index{tournament ! random}%
 and $(\mathbb{P}_1, \mathbb{P}_{-1})$ is the random
bipartite graph.
\index{graph ! random ! bipartite}%
In particular the partition into $\mathbb{P}_1$ and
$\mathbb{P}_{-1}$ is visible in the graph, with $\mathbb{P}_{-1} = \{$vertices $v
: \exists$ a directed edge containing $v\}$.  We get the same structures by taking the following vertex sets:
\[ \mathbb{P}_1\ =\ the\ set\ of\ primes\ which\ split\ in\ \mathbb{Q}[i].\]
\[ \mathbb{P}_{-1}\ =\ the\ set\ of\ primes\ which\ are\ inert\ in\ \mathbb{Q}[i].\]
\[ {2}\ =\ the\ prime\ which\ ramifies\ in\ \mathbb{Q}[i].\]

If we include $2$ we find that there is an edge from $p \to 2$ for all
$p$ and that there is an edge from $2 \to p$ iff $p \equiv \pm1 (8)$.
If we use a $2$-adic version $p \to 2$ iff $p \equiv 1 (4)$.  Finally we
get the following graph
$$
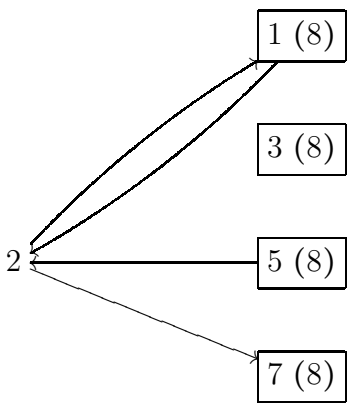
\xymatrix{ 
&&& *+[F]^{1\ (8)} \ar@/^/[ddlll]\\
&&& *+[F]^{3\ (8)}\\
{2} \ar@{->}[drrr] \ar@/^/[uurrr] &&& *+[F]^{5\ (8)} \ar@{->}[lll]\\
&&& *+[F]^{7\ (8)}
}$$

The construction of this section reinforces that the primes of a Gaussian number
field are randomly distributed.

The version of $\mathfrak{R}$ that is constructed using the Chinese Remainder Theorem and Dirichlet's Theorem~\cite{cameron} makes an unexpected appearance in Gareth Jones'
\index{Jones, G. A.}%
classification of the regular embeddings of complete bipartite graphs in orientable surfaces~\cite{jones1}.
\index{graph ! bipartite}%
 Here is how it arises.  Let $\to$ denote the binary relation on the set $\Pi$ of all prime numbers defined by $q \to p$ if and only if $q$ divides $p - 1$.  Regard $\Pi$ as a directed graph, with an arc from $q$ to $p$ whenever $q \to p$.  For each integer $n \ge 2$, let $\Pi_n$ denote the induced subgraph of $\Pi$ whose vertices are the prime factors $p_1, \ldots, p_l$ of $n$, formed by restricting the relation $\to$ to these primes.  In $\Pi$ there is an arc from the vertex $2$ to every other vertex.  Deleting this vertex and all incident arcs, and ignoring the direction of each remaining arc, we obtain an undirected graph $\Pi'$ whose vertices are the odd primes and edges joining $p$ and $q$ if and only if $q$ divides $p-1$ or vice versa.  Jones proves that $\Pi'  \cong \mathfrak{R}$.

\bigskip

\section{A Number-Theoretic Construction of $\mathfrak{R^{t}}$}

This follows the same procedure as in the previous section, and uses the theory of cubic reciprocity, elements of which are summarized in Appendix~\ref{NumberTheory}.

Let $U$ be a finite set of distinct cubic residues, and $V',\ W'$
be two finite sets of distinct cubic non-residues.  Define two sets by
$V = \{\omega . v_{j} : v_{j} \in V'\}$ and $W = \{\omega^{2} . w_{k} :
  w_{k} \in W'\}$, where $\omega$ denotes a cube root of unity.  

Now to construct $\mathfrak{R^{t}}$.  The vertex set is the set of
primary elements of $D$.  Join two such elements $\pi_1$ and $\pi_2$,
with a 
\begin{displaymath}
\left\{ \begin{array}{ll}
\text{non-edge} & \text{if}\ \chi_{\pi_{1}}(\pi_{2}) = 1\\
\text{single edge} & \text{if}\ \chi_{\pi_{1}}(\pi_{2}) = \omega\\
\text{double edge} & \text{if}\ \chi_{\pi_{1}}(\pi_{2}) = \omega^2,
 \end{array} \right.
\end{displaymath} where $\chi_{\pi_{1}}(\pi_{2})$ is the cubic residue character of $\pi_{2}$ modulo $\pi_{1}$.

The symmetry of the edge-joining condition is given by by the
statement of cubic reciprocity that can be found in Appendix~\ref{NumberTheory}.  Let $U, V, W$ be finite disjoint sets of primary ideals.  Choose $\alpha_{i}$ such that
$\chi_{u_{i}}(\alpha_{i}) = 1$ $\forall i$, $\beta_{j}$ such that
$\chi_{v_{j}}(\beta_{j}) = \omega$ $\forall j$, and $\gamma_{k}$ such that
$\chi_{w_{k}}(\gamma_{k}) = \omega^2$ $\forall k$.  That is, choose cubic residues
$\alpha_{i}\pmod{u_{i}}\ \forall u_{i} \in U$ (for example $\alpha_{i}
= 1$), cubic non-residues $\beta_{j}\pmod{v_{j}}\ \forall v_{j} \in V$,
and cubic non-residues $\gamma_{k}\pmod{w_{k}}\ \forall w_{k} \in W$.
The Chinese Remainder Theorem, of which there is a statement for general rings, implies that the congruences
 \index{congruence}%
\[ z \equiv \alpha_{i}\pmod{u_{i}}\ \forall u_{i} \in U \] 
\[ z \equiv \beta_{j}\pmod{v_{j}}\ \forall v_{j} \in V \] 
\[ z \equiv \gamma_{k}\pmod{w_{k}}\ \forall w_{k} \in W \] 
\[ z \equiv 1\pmod{3} \] 
have a unique solution $z\pmod{N}$, where $N = 3\ \prod u_{i}\  \prod
v_{j}\  \prod w_{k}$.

To capture the universality of $\mathfrak{R^{t}}$ we need to ensure an
infinite supply of vertices $z$ in the residue class.  This is given
to us by Dirichlet's Theorem which says that since $(z, N) = 1$, there
are an arithmetical progression of primes of the form $z + kN\ (k \in
\mathbb{N})$.  This theorem has far-reaching extensions, for example, to the Gaussian primes and to algebraic number fields; the Eisenstein integers form a commutative ring of algebraic integers in the algebraic number field $\mathbb{Q}\sqrt(-3)$.  This construction gives a graph isomorphic to
$\mathfrak{R^{t}}$ if we follow the ($*_t$) algorithm.

The construction of this section says that the primes in the ring
$\mathbb{Z}[\omega]$ of Eisenstein integers are randomly distributed.

\section{Universal Metric Spaces and $\mathfrak{R}$}

First we give an example of how to build $\mathfrak{R}$ from the points of a metric space.

\begin{theorem}
Let $M$ be any countable, universal, homogeneous, integral metric space.
\index{metric space ! integral}%
Putting $x \sim y$ if and only if $d(x, y)$ is odd gives $\mathfrak{R}$.
\end{theorem}

\begin{proof}
Consider the process of adding a new point $z$ with prescribed distances $d(z, a_i) = x_i$ from a given finite set of points $a_i, \ldots, a_n$.  These distances satisfy
\[ | x_i - x_j | \le d(a_i, a_j) \le x_i + x_j.\]
Assuming all $x_i$ are large enough, for example at least the diameter of $\{ a_i, \ldots, a_n \}$.  The upper bound is definitely satisfied, and the lower bound says
\[ x_i - d(a_i, a_j) \le x_j \le x_i + d(a_i, a_j). \]
The minimum value that $d$ takes is 1, so given $x_i$ there is always a range of at least 2 for the possible value of $x_j$, and we can choose its parity arbitrarily with $x_i \in \{m, m+1\}$ for some $m$. 
\end{proof}

If we take the countable universal homogeneous \emph{rational} metric space
\index{metric space ! integral}%
instead then we can also construct the random graph,
\index{graph ! random}%
 for example by partitioning the positive rationals into dense subsets $A, B$ and putting $x \sim y$ if and only if $d(x, y) \in A$~\cite{camver}.

The metric space equivalent of the random graph is uniquely defined by
a one-point extension property
\index{one-point extension property}%
 as it is for random graphs.  Some connections between metric spaces and random graphs are studied in~\cite[\S3.7]{cameron}, where metric space constructions are given for $\mathfrak{R}$, the almost-homogeneous random bipartite graph
\index{graph ! random ! bipartite}%
 $\mathfrak{B}$ and Henson's graph $\mathfrak{H}_k$. 
\index{graph ! Henson}%
 As pointed out in~\cite{cam22}, thinking of the edge-set of a simple graph on a given countable vertex set as a countable zero-one sequence,
\index{zero-one sequence}%
  leads to the fact that the class of graphs forms a complete metric space, of which the subclass of graphs satisfying the two-colour \emph{I-property} is a residual subset.
\index{I-property}%

A metric space $X$ is called a \emph{(generalized) Urysohn
space}~\cite{urysohn}
\index{Urysohn space}%
if whenever $A \subseteq X$ is a finite metric subspace of $X$
and $A' = A \cup \{a\}$ is an arbitrary one-point metric space
extension of $A$, the embedding $A \hookrightarrow X$ extends to an
isometric embedding $A' \hookrightarrow X$.  Up to isometry there is
only one universal complete separable Urysohn metric space, denote it
$\mathfrak{U}$, which contains an isometric copy of every separable
metric space.  Because there are an uncountable number of two-point isometric spaces let alone larger spaces, in order to construct a countable universal metric space as a Fra\"{\i}ss\'e limit
\index{Fra\"{\i}ss\'e class}%
 of the class of finite metric spaces, we need to impose some extra condition on distances otherwise there are too many spaces (there are uncountably many two-point spaces up to isometry).  There is a unique countable universal rational
metric space denoted $\mathbb{Q}\mathfrak{U}$ having all distances in
$\mathbb{Q}$, with $\mathfrak{U}$ being its completion.  The space $\mathfrak{U}$ is connected and locally connected and furthermore it is \emph{homogeneous}: every
isometry between two finite subspaces of $X$ extends to an
isometry of $\mathfrak{U}$ to itself; the same is true if we replace
finite by compact.  The group $\Aut(\mathfrak{U})$ is highly
transitive on isometry classes of finite subspaces, just as $\Aut(\mathfrak{R})$ is transitive on the set of
isomorphic finite subgraphs.  By comparison the Banach space
\index{Banach space}%
 $C[0,1]$ is universal but not homogeneous, whilst both the
infinite-dimensional Hilbert space
\index{Hilbert space}%
$\mathcal{H}$ and the unit sphere in $\mathcal{H}$ are homogeneous.
The infinite-dimensional Hilbert spaces play the role of Urysohn
metric spaces for the class of metric spaces embeddable into Hilbert
spaces.  

The way to show that $\mathfrak{U}$ is the `Random Polish Space'
\index{Polish space}%
is to choose $n$-tuples of random real numbers as distances from $n$-dimensional
Euclidean space
\index{Euclidean space}%
 $\mathbb{R}^n$, take a measure
\index{measure}%
for each dimension, and construct a
countable space one point at a time using a suitable measure on
extensions via a `random metric' between the $(n+1)^{\text{st}}$ point and
the first $n$ points, ensuring that the triangle inequality which determines a cone in $\mathbb{R}^n$ is satisfied each time.  Finally take the completion of this.  We can say more precisely what a random metric on $\mathbb{N}$ looks like.  Let $a_0, a_1$ be the first two points of the space.  Choose $d(a_0, a_1) = x_0^{(1)} \ge 0$, to lie in $\mathbb{R}^{+}$, where the superscript denotes dimension.  Next, choose $a_2$ so that  $d(a_0, a_2) = x_0^{(2)} = x \ge 0$ and $d(a_1, a_2) = x_1^{(2)} = y \ge 0$ satisfying
\[|x - y| \le x_0^{(1)} \le x + y.\]   Now choose $a_3$ so that  $d(a_0, a_3) = x_0^{(3)} = z_1 \ge 0$, $d(a_1, a_3) = x_1^{(3)} = z_2 \ge 0$ and $d(a_2, a_3) = x_2^{(3)} = z_3 \ge 0$ satisfying several inequalities of the form
\[|z_1 - z_2| \le x_0^{(2)} \le z_1 + z_2.\]  
Each cone at each new dimension depends on the previous one.  By choosing a wide range of reasonable measures on all cones, the completion of $\mathbb{N}$ equipped with a random metric is almost always isometric to $\mathfrak{U}$.  A. Vershik~\cite{vershik}
\index{Vershik, A. M.}%
attaches admissible vectors to distance matrices in order to mimic the one-point extension property.
\index{one-point extension property}%
  The indivisibility property of $\mathfrak{U}$ whereby if $B \subset \mathfrak{U}$ is an open ball
\index{open ball}%
then $\mathfrak{U} \backslash B \cong \mathfrak{U}$ is equivalent to removing a row and some columns in the distance matrices stabilizes $\mathfrak{U}$.

It is possible to construct $\mathfrak{R}$ from $\mathfrak{U}$.  An integer-valued construction of $\mathfrak{U}$ is given in~\cite[Theorem 3.11]{cameron}.  Any group acting on $\mathfrak{U}$ with a countable dense orbit preserves the structure of $\mathfrak{R}$ on the dense orbit~\cite[Proposition 7]{camver}.  Urysohn
space is generic with respect to Baire category, that is the class of
Urysohn spaces
\index{Urysohn space}%
 is residual in the class of all Polish spaces.
\index{Polish space}%

The homogeneity of $\mathfrak{U}$ implies that it is also $\aleph_0$-categorical,
\index{aleph@$\aleph_0$-categorical}%
because it is defined over a finite relational language~\cite{macth}.  Then by the
Engeler-Ryll-Nardzewski-Svenonius Theorem,
\index{Engeler--Ryll-Nardzewski--Svenonius ! Theorem}%
the group $\Aut(\mathfrak{U})$ is transitive on all isometric $n$-tuples because such a group can only have finitely many congruences.
 \index{congruence}%
   However within this constraint, every Polish metric space is isometric to the set $\fix(g)$ of fixed points of some $g \in \Aut(\mathfrak{U})$, so conjugacy in $\Aut(\mathfrak{U})$ is as complicated as it can be.  If $g \in \Aut(\mathfrak{U})$ is such that all its orbits are totally bounded then $\fix(g)$ is isometric to $\mathfrak{U}$.

If a metric space is separable, and therefore second-countable, then so is its automorphism group.
\index{group ! automorphism}%
  Both $\Aut(\mathbb{Q}\mathfrak{U})$ and $\BdAut(\mathbb{Q}\mathfrak{U})$ (the group of all \emph{bounded} isometries)
\index{bounded isometries}%
have been proven to be dense in $\Aut(\mathfrak{U})$~\cite{camver}.  Also
$\Aut(\mathbb{Q}\mathfrak{U}) < \Aut(\mathfrak{R})$; in fact it is shown in~\cite{camver} that there are $2^{\aleph_{0}}$ different reducts all of which are isomorphic to
$\mathfrak{R}$.  As $\mathfrak{U}$ is uncountable, any study
of reducts should be restricted to $\mathbb{Q}\mathfrak{U}$.  

Macpherson
\index{Macpherson, H. D.}%
has proved~\cite{macph} that any closed oligomorphic subgroup of $\Sym(\omega)$ contains a free subgroup
\index{group ! free}%
of infinite rank; this includes both $\Aut(\mathbb{Q}\mathfrak{U})$ and the group $\Iso(\mathfrak{U})$ of isometries of the Urysohn space, into which $\Aut(\mathbb{Q}\mathfrak{U})$ embeds densely.

Urysohn proved that there is no isometry of $\mathfrak{U}$ mapping \emph{infinite} subspaces to each other and Trofimov
\index{Trofimov, V. I.}%
 proved that there is no isomtery  of $\mathfrak{U}$ mapping two arbitrary \emph{countable} subsets to each other~\cite{trof1}.

In~\cite{vershik1} Vershik
\index{Vershik, A. M.}%
 proves the equivalence of the fact that $\Iso(\mathfrak{U})$ contains an everywhere dense locally finite subgroup with the the generalization of Hrushovski's theorem
\index{Hrushovski, E.}%
 (see also the work of Herwig and Lascar~\cite{herwig}~\cite{herwiglas})
\index{Herwig, B.}%
\index{Lascar, D.}%
 regarding the globalization of the partial isomorphisms of finite graphs~\cite{hrush} to metric spaces, namely that, for each finite metric space $X$ there exists another finite metric space $\overline{X}$ and isometric embedding $i$ of $X$ into $\overline{X}$ such that isometry $i$ induces the embedding of the group monomorphism of $\Iso(X)$ to $\Iso(\overline{X})$ and each partial isometry of $X$ (that is, isometry between two subsets of $X$) can be extended to a global isometry of $\Iso(X)$.  He also proves the following strengthening of homogeneity
 
\bigskip

\begin{theorem}
For any finite subset $X \subset \mathbb{Q}\mathfrak{U}$ there exists an isomorphic embedding $j : \Iso(X) \to \Iso(\mathbb{Q}\mathfrak{U})$ such that $(jg)(x) = g(x)$ for all $g \in \Iso(X), x \in X$.  Consequently, the space $\mathbb{Q}\mathfrak{U}$ is a countable union of finite orbits of the group $j\Iso(X)$.  
\end{theorem}

It follows that each finite metric space or Polish space can be embedded equivariantly into $\mathfrak{U}$.  Vershik goes on to remark that for graphs the words `random' and `universal' coincide despite the randomness in question not referring to one graph but indicative of the measure on the set of all graphs, but for continuous spaces the issue is more subtle and the role of a specific choice of measure is essential.

M. Doucha
 \index{Doucha, M.}%
has constructed an abelian separable invariant metric group $G$ that is universal and homogeneous; that is, every separable abelian topological group
\index{group ! abelian}%
\index{group ! topological}%
 equipped with an invariant compatible metric can be isometrically embedded into $G$, and any isometric homomorphism between two finitely generated subgroups of $G$ extends to an isometric automorphism
\index{group ! automorphism}%
 of $G$.  In particular, for every abelian Polish group
\index{group ! Polish}%
 there is a topologically isomorphic closed subgroup of $G$.  The construction is related to work of Cameron \index{Cameron, P. J.}%
and Vershik
\index{Vershik, A. M.}%
 imposing an abelian group structure on Urysohn space.
\index{Urysohn space}%
 Some more of Doucha's results on universal homogeneous Polish metric structures can be found in~\cite{doucha}.

Further properties of $\mathbb{Q}\mathfrak{U}$ are discussed
in~\cite{cam22} and~\cite{camver}; a construction of the rational Urysohn space is given in~\cite{hubickanes}.  See also the workshop proceedings on Urysohn space edited by A. Leiderman et al~\cite{leiderman}.
\index{Leiderman, A.}%
Separable and non-separable universal metric spaces are studied from a different viewpoint in the monograph~\cite{lipscomb}.

\section{Homogeneous Integral Metric Spaces}

Accounts of Cameron's
\index{Cameron, P. J.}%
 previous work on homogeneous integral metric spaces is given in~\cite{cameron} and ~\cite{cam8a}.

The unique countable homogeneous integral metric space $M_N$ of diameter $N$ is the same as the unique countable distance-homogeneous graph of diameter $N$ that has arisen in the work of Lawrence S. Moss,
\index{Moss, L. S.}%
so we give a synopsis of his relevant work on universal graphs of finite diameter, and then give a theorem that begins a classification of integral metric spaces.

In a graph $\Gamma$, denote the shortest length-$N$ path from $x$ to $y$ by $d_{\Gamma}(x, y) = N$, where an infinite length means that $x$ and $y$ belong to different connected components of $\Gamma$.  An isometric embedding
\index{graph ! isometric embedding}%
 between graphs is a distance-preserving map.  If $\overline{x}$ (resp. $\overline{y}$) is a length-$k$ path in $\Gamma_1$ (resp. $\Gamma_2$), write $\overline{x} =_{N} \overline{y}$ to mean that whenever  $1 \le i, j \le k$ and either $d_{\Gamma_1}(x_i, x_j) \le N$ or $d_{\Gamma_2}(y_i, y_j) \le N$, then $d_{\Gamma_1}(x_i, x_j) = d_{\Gamma_2}(y_i, y_j)$.  An isometric embedding of graphs is a stronger condition than an isomorphic embedding.  A graph $\Gamma$ is \emph{distance-homogeneous}~\cite{moss} 
\index{graph ! distance-homogeneous}%
if for every pair of tuples $\overline{x}$ and $\overline{y}$ of the same length such that $\overline{x} =_{\infty} \overline{y}$ there is a $g \in \Gamma : \overline{y} = \overline{x} g$.

In~\cite{moss1}~\cite{moss2} Moss showed that for each natural number $N$ there is a countable graph $\mathfrak{R}_N$~\label{mathfrak{RN}} of diameter $N$
\index{graph ! diameter N@diameter-$N$}%
and into which there is an isometric embedding of every countable graph of diameter $N$.  So, $\mathfrak{R}_0$ is a single vertex, $\mathfrak{R}_1$ is the complete graph on infinitely many vertices, and $\mathfrak{R}_2$ is the random graph $\mathfrak{R}$.
\index{graph ! random}%
  So the graphs $\mathfrak{R}_N$ are in a sense higher-dimensional analogues of $\mathfrak{R}$.  Each graph $\mathfrak{R}_N$ is the unique countable distance-homogeneous graph
\index{graph ! distance-homogeneous}%
of diameter $N$ which isometrically embeds every countable diameter-$N$ graph.  There is a universal countable connected graph $\mathfrak{R}_{\infty}$ first studied by Pach~\cite{pach}
\index{Pach, J.}%
and rediscovered in~\cite{moss} into which every countable connected graph is embedded.  Whilst the $\mathfrak{R}_N$ do not sit in each other in any obvious way and so there is no way of arbitrarily taking unions, in terms of the first-order theories
\index{first-order theory}%
 of these graphs,
\[ \Th(\mathfrak{R}_{\infty}) = \lim_{N \to \infty} \Th(\mathfrak{R}_N). \]

In other words, the first-order theory of the unique countable distance homogeneous countable graph which isometrically embeds every countable graph is the model completion of the theory of distanced graphs.  A first-order sentence is satisfied by $\mathfrak{R}_{\infty}$ if and only if it is satisfied by $\mathfrak{R}_N$ for almost all $N$.  There is one countable connected model but more than one uncountable model.

Corresponding to the ($*$)-condition for $\mathfrak{R}$, there are two equivalent characterizations of $\mathfrak{R}_N$, as follows.  The graph $\mathfrak{R}_N$ is the unique countable graph with the following property:

($*_N$)  Let $U_1, U_2, \ldots, U_{N-1}$ and $V$ be finite disjoint sets of vertices of $\mathfrak{R}_N$.  Suppose that whenever $1 \le i \le j \le N-1$, $x \in U_i$ and $y \in U_j$, we have that $j - i \le d_{\mathfrak{R}_N}(x, y) \le j + i$.  Also, whenever $1 \le i \le N-1$, $x \in U_i$ and $y \in V$ we have that $N - i \le d_{\mathfrak{R}_N}(x, y)$.  Then there exists in $z \in \mathfrak{R}_N$ such that (a) for all $i \le N$ and all $x \in U_i, d_{\mathfrak{R}_N}(x, z) = i$; and (b) for all $y \in V$, and all $d_{\mathfrak{R}_N}(y, z) \ge N$,

or equivalently

($*_N$)  Let $\overline{x}$ and $\overline{y}$ be $k$-tuples from $\mathfrak{R}_N$ and $\Gamma$ respectively, where $k \in \omega$ and $\Gamma$ are arbitrary, and suppose that $\overline{x} =_{2(N-1)} \overline{y}$.  Then $\forall y_1 \in \Gamma, \exists x_1 \in \mathfrak{R}_N$ such that $\overline{x}, x_1 =_{N-1} \overline{y}, y_1$. 

Every countable connected graph can be embedded into $\mathfrak{R}_{\infty}$ in $2^{\aleph_0}$ ways and $\mathfrak{R}_{\infty}$ has $2^{\aleph_0}$ many automorphisms.  A countable collection of copies of $\mathfrak{R}_{\infty}$ is the unique distance-homogeneous countably universal graph. 
\index{graph ! distance-homogeneous}%
The finite distance-homogeneous graphs have been classified~\cite{camab} and found to be precisely those graphs with the property that for every pair of $6$-tuples $\overline{x}$ and $\overline{y}$ such that $\overline{x} =_{\infty} \overline{y}$, there is a graph automorphism taking $\overline{x}$ to $\overline{y}$ pointwise.

Whilst the existence of universal countable graphs under isometric embeddings does not follow from general model theory,
\index{model theory}%
 there are characterizations related to the concept of existentially closed structures~\cite{moss} 
\index{graph ! n@$n$-existentially closed}%
(see Appendix~\ref{FurtherDetails}).

Finally, consider a graph $\Gamma$ with the following I-property ($I_{P}$)
\index{I-property}%

($I_{P}$)  If whenever $A$ and $B$ are connected graphs with $|A| \le P$ and $|B| = |A| + 1$, and $i : A \to \Gamma$ and $j : A \to B$ are isometric embeddings, then there is an isometric embedding $k : B \to \Gamma$ such that $k \circ j = i$.

The Johnson graph $J(n, m)$
\index{graph ! Johnson}%
 has vertices the $m$-subsets of $\{1, \ldots, n\}$, and two vertices connected if and only if their intersection has size $m-1$.  In~\cite{moss3}, Moss and Dabrowski
\index{Dabrowski, A.}%
 show that the Johnson graphs $J(n, 3)$ satisfy ($I_3$) whenever $n \ge 6$, and that $J(6, 3)$ is the smallest graph satisfying ($I_3$).

We say that a graph $\Gamma$ has the I-property ($LI_{P}$)
\index{I-property}%
if for all $v \in \Gamma$, the subgraph $\Gamma(v)$ induced by the neighbours of $\Gamma$ has property ($I_{P}$).  The Johnson graphs $J(n, m)$ are locally gridlike, and so do not satisfy ($LI_{P}$), but $\mathfrak{R}$ and $\mathfrak{R}_{\infty}$ satisfy  ($I_{3}$) and ($LI_{P}$) for all $P$.

\bigskip

Next we turn to reducts of integral metric spaces.  We can show that, in the universal
homogeneous metric space of diameter $n$, the ``distance at most $k$'' graph
is isomorphic to $\mathfrak{R}$ if and only if $k+1 \le n \le 2k$. More generally, we
can find necessary and sufficient conditions for a reduct of $M_n$ to be
isomorphic to $\mathfrak{R}$.  

Let $M_n$ be the unique countable homogeneous integral metric space of
diameter~$n$.
Thus, $M_2$ is the vertex set of the random graph
\index{graph ! random}%
 $\mathfrak{R}$ with the path metric.
 
\begin{theorem}
Let $\{1,\ldots,n\}=A\cup B$, where $A\cap B=\emptyset$. Form a graph $\Gamma$
on the vertex set $M_n$, by joining $x$ and $y$ if and only if $d(x,y)\in A$.
Then $\Gamma\cong \mathfrak{R}$ if and only if
$\{\lceil n/2\rceil,\ldots,n\}\not\subseteq A$ and
$\{\lceil n/2\rceil,\ldots,n\}\not\subseteq B$.
\end{theorem}

\begin{proof}
Since $\mathfrak{R}$ is self-complementary, we can interchange $A$ and $B$ where necessary.

Suppose first that $\{\lceil n/2\rceil,\ldots,n\}\subseteq B$. Then edges join
points at distance less than $n/2$; so $\Gamma$ has diameter greater than~$2$,
and cannot be isomorphic to $\mathfrak{R}$.

Now suppose that $\{\lceil n/2\rceil,\ldots,n\}$ is not contained in either
$A$ or $B$. We may suppose that $n\in B$. Let $a$ be the maximal element of
$A$; so $a>n/2$, and $a+1\in B$. Let $U$ and $V$ be finite disjoint sets of
points of $M_n$. We seek a point $z$ such that $d(z,u)=a$ for all $u\in U$,
and $d(z,v)=a+1$ for all $v\in V$. These requirements are consistent. For
the consistency conditions (instances of the triangle inequality) are
\begin{eqnarray*}
0\le d(u_1,u_2)\le 2a &\hbox{for}& u_1,u_2\in U,\\
1\le d(u,v)\le 2a+1 &\hbox{for}& u\in U,v\in V,\\
0\le d(v_1,v_2)\le 2a+2 &\hbox{for}& v_1,v_2\in V.
\end{eqnarray*}
These are satisfied since $2a\ge n$. So by the properties of $M_n$, the
point $z$ exists. Now by construction $z$ is joined to $U$ for all $u\in U$,
and not joined to $v$ for all $v\in V$. Thus, $\Gamma$ satisfies the
well-known characterisation of $\mathfrak{R}$.
\end{proof}

\begin{corollary}
The graph obtained from $M_n$ by joining points at distance at most $k$ is
isomorphic to $\mathfrak{R}$ if and only if $k+1\le n\le 2k$.
\end{corollary}

The model of  $\mathfrak{R}$ is a submodel of the model of $M_n$, so $\Aut(M_n) \le \Aut(\mathfrak{R})$, so 
the $\mathfrak{R}$-reducts lie between $\Aut(M_n)$ and $\Sym(M_n)$.

How many of the reducts of $M_k$ are isomorphic to $\mathfrak{R}$?  Note that these will also be reducts of $M_n$.  A space $M_k$ is a reduct of $M_n$ when it is possible to divide $n$-distances into $k$ classes such that you get a reduct isomorphic to $M_k$.

\smallskip

Which other reducts of $M_k$ for $k \geq 3$ are also reducts of $M_n$?  We need to search for higher-order relations on the points of $M_k$.  For  $M_3$, we require non-oriented ternary relations on triples of points.  

Let $A \subseteq M_n$ be a subset of points having only odd distances.  We need to consider triangles in which the perimeter is odd.  Given a subset $A$, what is the graph got by taking distances in $A$ to be proper reducts of the original metric space?  What is the two-graph?

There are two-graph-like conditions on triples directly derivable from the metric space that are not reducts of $\mathfrak{R}$.  Let $\{1,\ldots,n\}=A\cup A^c$, where $A^c$ denotes the complement of $A$ in the $n$-set.  Assume a ternary relation $(x, y, z)$ on $M_n$ if and only if there are an odd number of the distances in $A$.  The ternary relation defines a two-graph on triples of $M_n$ if and only if  all of the odd distances chosen belong to $A$, and 
$\{\lceil n/2\rceil,\ldots,n\}\not\subseteq A$ and
$\{\lceil n/2\rceil,\ldots,n\}\not\subseteq A^c$; this is the same as the two-graph on vertices of $\mathfrak{R}$.  However, with these criteria on $A$ and $A^c$, if one of the odd distances lies in $\{\lceil n/2\rceil,\ldots,n\}$ then the ternary structure on triples of $M_n$-points is no longer a two-graph; call it a \emph{$\mathsf{T}$-structure}.
\index{tstructure@$\mathsf{T}$-structure}%
(This is not to be confused with a different ternary relational construct called a T-structure that appears in~\cite[p.~68]{cam6}, which is associated with tournaments.)  We must prove that a $\mathsf{T}$-structure is a reduct of $M_n$.  Clearly, each automorphism of $M_n$ is also an automorphism of the $\mathsf{T}$-structure, so that $\Aut(M_n) \leq \Aut(\mathsf{T})$, so it remains to show that $\Aut(\mathsf{T})$ is closed, and we leave this as an open question for further work.

\head{Open Question}  Classify the reducts of integral metric spaces.




\index{model theory}%
 (see Appendix A.2).

\section{Miscellaneous Observations}

1.  There is a method of building $\aleph_0$-categorical structures
\index{aleph@$\aleph_0$-categorical}%
with a countably infinite domain by using finite substructures.  If $N$ is a finite substructure of a homogeneous structure $M$, then
all automorphisms of $N$ extend to automorphisms of $M$ fixing $N$; so
two tuples lying in the same orbit of $\Aut(N)$ also lie in the same
orbit of $\Aut(M)$.  If the converse holds we say that $N$ is a
\emph{finite homogeneous substructure}
\index{finite homogeneous substructure}%
 of $M$.  So a finite homogeneous substructure is itself a homogeneous structure.  Call $M$
\emph{smoothly approximable}
\index{smoothly approximable structure}%
 if it is the union of a chain of finite
homogeneous substructures of itself.  The canonical example of a smoothly approximable structure is a projective space of countable dimension over a finite field, which is the union of a chain of finite-dimensional subspaces.  Kantor, Liebeck and Macpherson classified~\cite{kantor}
\index{Kantor, W. M.}%
\index{Liebeck, M. W.}%
\index{Macpherson, H. D.}%
$\aleph_0$-categorical smoothly approximable structure with primitive and oligomorphic automorphism groups.
\index{group ! automorphism}%
  Hrushovski
\index{Hrushovski, E.}%
showed~\cite{hrushovski} that the assumption of primitivity is not necessary.  Smoothly approximable structures are further studied in~\cite{cherlinhru}.  The random bipartite graph
\index{graph ! random ! bipartite}%
 $\mathfrak{B}$ provides an example of a primitive $\aleph_0$-categorical permutation structure that is not smoothly approximable~\cite[p.~457]{kantor}.

\medskip

2.  We make one final observation about metric space theory
\index{metric space}%
with \cite[\S~3.7]{cameron} and~\cite[pp.232-234]{brco} in mind.  If we take a random graph
\index{graph ! random}%
 of diameter at most 2, and we assign distances between vertices $u, v$ as follows:
\begin{displaymath}
d(u, v) = \left\{ \begin{array}{ll}
\text{0} & \text{if}\ u = v\\
\text{1} & \text{if}\ u \sim v\\
\text{2} & \text{otherwise},
 \end{array} \right.
\end{displaymath}
then the distance-$1$ graph will be an integral metric space.
\index{metric space ! integral}%
Clearly, by inspection, this construction \emph{cannot distinguish} between
$\mathfrak{R}$, $\mathfrak{R}^{t}$ or any of the higher-adjacency
random graphs, because all distances are $\le 2$ regardless of the
adjacency.  This shortcoming can be overcome by introducing a theory
of coloured metric spaces, in which all distances have an upper bound, but two points are adjacent if and only if their
distance apart lies within the requisite upper bound and is
monochromatic.  This would give us as many different monochromatic random graphs
\index{graph ! random}%
 as there are colours.  Another possible way of overcoming the obstacle to modeling the higher-adjacency random graphs is to utilise the language of algebraic topology as is done in~\cite[p.60]{cameron}.  

\medskip

3.  The proof given in~\cite{cameron} that the injectivity property for
the random graph holds for the set-theoretic
construction of $\mathfrak{R}$ relies fundamentally on the axiom of
foundation (AF)
\index{axiom of foundation}%
not being violated.  A natural question arises as to
whether or not the result goes through for $\mathfrak{R^{t}}$ with an extension to the
anti-foundation Axiom (AFA).
\index{axiom of anti-foundation}%
This axiom is a prerequisite for the
existence of sets with $\in$ symmetrised, which would be required if in addition to the edge or non-edge possibilities, we were to take the third adjacency type to be a double edge, or set theoretically to demand both $x \in y$ and $y \in x$.

There is a loop version of the \emph{I-property}
\index{I-property}%
defining $\mathfrak{R^{t}}$ whereby an extra vertex with a loop, say $z^{l}$, must be found at each iteration.  To demonstrate that this property ($*_{t, l}$) does \emph{not} yield all possible
anti-founded sets, let $U, V, W$ be finite disjoint
sets of sets, and take
\[ z = \{ u_1, u_2, \ldots, u_p, w_1,\ldots, w_r, x_1,\ldots, x_q \} \]
where $x_j \notin \bigcup v_j, x_j \notin \{ v_1, \ldots, v_q \}$ for $(j = 1, \ldots, q)$, and
\[ z^{l} = \{ u_1, u_2, \ldots, u_p, w_1,\ldots, w_r, y_1,\ldots, y_q,
z^{l} \} \]
where $y_j \notin \bigcup \bigcup v_j, y_j \notin \bigcup
v_j, y_j \notin \{ v_1, \ldots, v_q \}$ for $(j = 1, \ldots, q)$.

Consider vertex $z$.  Certainly all members of $U$ and $W$ are in
$z$ if we consider all such $z$s.  Suppose $\exists v \in V$ joined to
$z$.  Then either $v \in z$; whence $v \in \{ u_1, \ldots, u_p \}$
or  $v \in \{ w_1, \ldots, w_r \}$ (both of which contradict our
disjointness assumption), whilst $ v \in x_j \notin \bigcup v_i $ is
also a contradiction.  Or, $ z \in v$ implying the contradiction $x_j
\in v$, for $x_j \notin \{ v_1, \ldots, v_q \} $, that is taking $v = v_j$
gives the contradiction $x_j \in \bigcup v_j$.

Now consider vertex $ z^{l} $.  Then $ z^{l} \in v $ implies $ y_j \in
z^{l} \in v $, contradicting $y_j \notin \bigcup v_j $.  Also $ y_j
\in z^{l} \in z^{l} \in v $ contradicts $ y_j \notin \bigcup \bigcup
v_j $.  Similarly $ v \in  z^{l} $ would imply $ v \in y_j \notin \{ v_1,
\ldots, v_q \} $.  However now we come to the crux of the matter.
Whilst $w_k \in z$ for all $k$, by assumption $w_1 \cap w_2 =
\emptyset$, so there is no $z$ for which $z \in w_1$ \emph{and} $z \in w_2$.

This example firstly shows that by weakening FA to AFA we are no longer able to satisfy ($*_{t, l}$), and in fact that AF is needed to verify ($*_{t, l}$).  Secondly the set-theoretic version of the defining condition for the triality graph
\index{graph ! triality}%
 is too strong to allow all possible graphical representations of all anti-founded sets, and we
cannot contend that $ZFA$ is the theory of orientations of
$\mathfrak{R^{t}}$ in the way that $ZF$ is the theory of
orientations of $\mathfrak{R}$~\cite{cameron}.  There are now several questions that can be asked.  Suppose we take a specific form of an anti-foundation axiom (as Barwise and Moss do in~\cite{barmos1}).
\index{Barwise, J.}%
\index{Moss, L. S.}%
Now given a countable model of ZFA with this axiom, we get a graph by symmetrizing $\in$ in that model.  Is the graph unique (up to isomorphism)?  If so, how can it be described?  If not, how many different graphs can arise thus?

Another way to proceed is to try to construct
a universal structure.  For example, take unary relations~\label{unaryr} $W_0, W_1,
W_2, \ldots$ where $W_i$ consists of $i$ vertices doubly connected to
a certain $z$, and insist that for all vertices, only one of the $W_i$
holds.  Does this give a universal structure?  Even so it may not be
homogeneous and amalgamation may fail.  So can we add further relations to
ensure that we get a Fra\"{\i}ss\'e limit?
\index{Fra\"{\i}ss\'e class}%

We make some further set-theoretic remarks.
\begin{enumerate}

\item  We can demonstrate simply the existence of random graphs
\index{graph ! random}%
 of arbitrarily large infinite cardinality, noting however that
identifying a canonical model of one, even with $\aleph_1$ vertices is
not necessarily straightforward.  Recall that every
ordinal number $\alpha > 0$ that is less than the smallest ordinal $\epsilon$ such that $\epsilon = \omega^{\epsilon}$ has a unique expansion in terms of a
decreasing finite sequence of ordinals $\alpha \geq 0 $,
\[ \alpha = \omega^{\alpha_{1}} + \omega^{\alpha_{2}} + \ldots + \omega^{\alpha_{k}}.\]
Now consider the base $2$ expansion of this so-called Cantor Normal Form of $\alpha$.
\index{Cantor Normal Form}%
\index{Cantor, G.}%
Take all ordinals having such an expression (hereditarily) less than $\omega_{\alpha}$ for some $\alpha$.  Label the vertices of Rado's construction of the random graph~\cite{cameron} with
ordinals $\alpha, \beta, \ldots $.  Join vertex $\beta$ to vertex
$\alpha$ (with $ \beta < \alpha$ as ordinals) if and only if $2^{\beta}$ appears in
the $\alpha$-expansion.  Now, for any finite disjoint sets of numbers
$U$ and $V$, we can add vertices to $U$ in order to achieve $max(U) >
max(V)$.  Then defining $z := \sum_{u \in |\alpha|}  2^{u}$ ensures that
$z$ is joined to every vertex in $U$ and to no vertex in $V$.  In view of this being an extension of Rado's
dyadic construction~\cite{rado}, there is a sense in which random graphs
constructed in this way could be considered to be the canonical
versions of higher infinite cardinality.

\item  The \emph{axiom of universality}
\index{axiom of universality}%
 (that every extensional graph
     has an injective decoration) is incompatible with the
     anti-foundation axiom (stating that every graph has a unique
     decoration), as the former yields many more non-well-founded sets;  for
     example, it implies the existence of collections of reflexive
     sets (of the form $x = \{ x \}$) of arbitrary cardinality.  Point
     (1) above indicates a constructive existence proof of such
     collections if the vertices $z$ are chosen to have loops.
     Notice that $\omega = 2^{\omega} (= n^{\omega})$, giving $\omega =
     \{ \omega \}$ for vertex $ \omega $.  Any vertex $z$ with a loop
     can be defined so that $ z = 2^{z}$.  \emph{Epsilon numbers}
\index{epsilon number}%
     are ordinals $ \epsilon $ such that  $\epsilon = {\omega}^{\epsilon}$ but may also be defined as ordinals $ \epsilon >
     \omega $ such that $ 2^{\epsilon} = \epsilon$.  The smallest ordinal satisfying the equation $\omega^{\alpha} = \alpha$ is $\alpha = \epsilon_0$, and so is the limit of the sequence $0, 1, \omega, \omega^{\omega}, \omega^{\omega^{\omega}}, \ldots$.  Since the class
     of epsilon numbers is of equicardinality with the class of
     ordinals, we can find as many vertices $ \epsilon_{\alpha} $ such
     that $ \epsilon_{\alpha} = \{ \epsilon_{\alpha} \} $ as
     there are ordinals $\alpha$.

\item  There is a generalization of the random graph
\index{graph ! random}%
 $\mathfrak{R}$ in another direction, to greater than the one dimension of
     $\mathfrak{R}$.  This is done by imagining an $n$-dimensional
     space, where each dimension has a coordinate with a dyadic
     expansion as per Rado's.  This transports us into
     the realm of \emph{simplicial complexes}.
     \index{simplicial complex}%
  These are natural generalizations of graphs; the first-order language
\index{first-order language}%
 $L_d$ of $d$-complexes resembles that of graphs, except that in place of the adjacency predicate $\sim$, it has $d$ predicates $\sim_1, \ldots, \sim_d$, where $\sim_i$ is an $(i+1)$-place predicate.  The \emph{dimension} of a simplex is one less than its cardinality, and a $d$-dimensional simplex reduces to a graph when $d = 1$.  It has been shown~\cite{blasshar} that every property of $d$-complexes expressible by an $L_d$-sentence holds of almost all or almost no $d$-complexes, just as it does for graphs.  In fact Rado's original paper~\cite{rado} has a result on universal simplicial complexes.

\end{enumerate}

\medskip

4.  Throughout this work we have assumed that the background set
theory is $ZF$.  Another possibility is to work within a set theory in
which the axiom of choice is false.  This allows us to define
so-called \emph{amorphous sets}~\cite{trussa}.
\index{amorphous set}%
A set is said to be \emph{amorphous} if it is infinite but is not the
disjoint union of two infinite sets.  This differs from the pigeonhole property
\index{pigeonhole property}%
obeyed by the countable homogeneous graphs that we have been studying, but amorphous sets appear to have some `pigeonhole-like' property: given any partition into two parts, one is finite and the other is amorphous.  In particular no countably
infinite set can be amorphous and so no countably infinite structure,
in particular $\mathfrak{R}$, can be built from such a set.  An equally immediate way to see that we
cannot construct the countably infinite random graph using an amorphous vertex set is to observe that one of the cycle types realized by
$\Aut(\mathfrak{R})$ has two infinite cycles and one fixed point (that is 
$\aleph_{0}^{2}.1^{1}$), which cannot be realized by any permutation
group acting on an amorphous set.  A more instructive argument runs as follows.  For a set $\Omega$, the canonical relational structure
\index{canonical structure}%
\index{relational structure}%
 $\mathcal{M}$ for a certain permutation group $G$ acting on $\Omega$ is such that
$\Aut(\mathcal{M})$ is the closure of $G$ in $\Sym(\Omega)$.  The $age$
of a relational structure $\mathcal{M}$ on $\Omega$ is the class of all finite
structures embeddable in $\mathcal{M}$ as induced substructures.  Regarding
Age$(\mathcal{M})$ as a tree with nodes in $\mathbb{N}$, let
$\mathcal{P(M)}$ be the set of infinite paths from the root.  For any
relational structure we can turn $\mathcal{P(M)}$ into a complete
metric space by defining a metric on $\mathcal{P(M)}$.  But if the
domain of the structure is an amorphous set, two-way infinite paths
are forbidden, and a metric cannot be defined for the whole space,
only a `partial' metric for finite paths.  Equally, we cannot apply Baire
category theory,
\index{Baire category theorem}%
 because we do not have residual sets, that is those
that contain an intersection of countably many dense open sets,
\index{dense open set}%
because the complement would be a union of countably many closed dense
sets.  Because the set of graphs on vertex set $\mathbb{N}$ which are
isomorphic to $\mathfrak{R}$ is residual in $\mathcal{P(\mathfrak{R})}$,
an amorphous set cannot be the domain of $\mathfrak{R}$, regarded as a
relational structure.  More generally, an amorphous set $\Omega$ does not obviously imply that either $\mathcal{P(M)}$ or Age$(\mathcal{M})$ is amorphous; the set of 2-element subsets of $\Omega$ is definitely not amorphous, since we can take $\alpha \in \Omega$ and split the pairs into those containing $\alpha$ and the rest.

In the following appendix we give some results obtained by working in a set theory in which the axiom of choice is false.

\section{Appendix: Bounded Amorphous Sets}

\subsection{The Theorem of Baer, Schreier and Ulam for Bounded Amorphous Sets}
\index{Baer, R.}%
\index{Schreier, J.}%
\index{Ulam, S.}%

An infinite set $U$ is amorphous if every subset of $U$ is finite or cofinite; an
amorphous set is \emph{bounded}
\index{amorphous set ! bounded}%
if there is a natural number $n$ such that
every non-trivial partition of $U$ has all but finitely many parts of 
size~$m$, for some $m\le n$.  In this appendix, as an exercise in working inside models of set theory in which the axiom of choice is false, we prove that within certain Fraenkel-Mostowski models of set theory,
\index{Fraenkel-Mostowski set theory}%
if $U$ is bounded amorphous, then the symmetric group on $U$, $\Sym(U)$, modulo the group of finitary permutations, $\FSym(U)$, is isomorphic to a finite group, and that every finite group can occur, and as a corollary, a variation of the Baer--Schreier--Ulam theorem.

\subsection{Permutations of a Bounded Amorphous Set}

A set is said to be \emph{amorphous} if it is infinite but is not the
disjoint union of two infinite sets. Amorphous sets are
\emph{Dedekind finite},
\index{Dedekind finite set}%
that is have no countably infinite
subset.  The presence of Dedekind finite infinite cardinals means that the
axiom of choice (AC) cannot hold. 
\index{axiom of choice}%

Truss
\index{Truss, J. K.}%
and co-workers~\cite{ctrussa1,ctrussa2,trussa} have
carried out an extensive study aimed at classification of amorphous sets
and their neighbours by examining the structure that an amorphous set can
carry, using finite permutation groups as invariants.  We will show that within certain Fraenkel-Mostowski models of set theory, that John Truss builds for his classification of bounded amorphous sets~\cite{trussa},  these groups have another role to play.

A permutation is \emph{finitary} if it fixes all but finitely many
points of $U$. The finitary permutations of $U$ form a group $\FSym(U)$,
the \emph{finitary symmetric  group}
\index{group ! finitary}%
on $U$, which is a normal subgroup
of $\Sym(U)$. A finitary permutation has a \emph{parity} (even or odd),
just as for a permutation of a finite set; the even permutations form a
subgroup of index~$2$ in $\FSym(U)$, the \emph{alternating group}
$\Alt(U)$.  A \emph{permutation group}~\cite{cam1} on $U$ is a subgroup of $\Sym(U)$. It is
\emph{transitive}
\index{group ! transitive}%
if the only fixed subsets of $U$ are $\emptyset$ and $U$,
and \emph{primitive}
\index{group ! primitive}%
if in addition the only fixed partitions of $U$ are
$\{U\}$ and the partition into singletons.

In ZFC, the factor group $\Sym(U)/\FSym(U)$ has cardinality~$2^{|U|}$, and its
normal subgroups form a chain. Indeed by the Theorem of Baer, Schreier and Ulam~\cite{baer}~\cite{schreierulam},
\index{Baer--Schreier--Ulam Theorem}%
the proper normal subgroups of
$\Sym(U)$ are (i) the trivial group; (ii) $\Alt(U)$; (iii) the bounded symmetric group $\BSym_\alpha(U)$~\label{BSym} consisting of all permutations moving fewer than $\alpha$ points, for each infinite cardinal $\alpha\le|U|$; (iv) $\Sym(U)$.

As we will see, the position is very different in the absence of the Axiom of
Choice!

A set is \emph{strictly amorphous}
\index{amorphous set ! strictly}%
if in addition to being amorphous, it is also the case that in any
partition of the set into infinitely many pieces, all except finitely
many are singletons.

\begin{theorem}[Jordan--Wielandt Theorem]
\index{Jordan--Wielandt Theorem}%
Let $U$ be an infinite set. Then a primitive permutation group on $U$
containing a non-identity finitary permutation must contain the alternating
group $\Alt(U)$.
\end{theorem}

We remark that the proof of this theorem (see~\cite[p.~166]{cam1}) makes
no use of the axiom of choice, and is valid in ZF.  With the help of this theorem, all transitive permutation groups on a strictly
amorphous set can be described. 

\begin{theorem}
Let $U$ be a strictly amorphous set. Then
\begin{enumerate}
\item $\Sym(U)=\FSym(U)$;
\item the only transitive subgroups of $\Sym(U)$ are $\Sym(U) \& \Alt(U)$.
\end{enumerate}
\end{theorem}

\begin{proof}
The proof of 1 is easy.

As for 2, if $G$ is transitive, then any invariant partition has all its parts of the
same size, and so (since $U$ is strictly amorphous) must be trivial. The
result now follows from the Jordan--Wielandt theorem.
\end{proof}

Let $U$ be an amorphous set, that is it is infinite but is not the
disjoint union of two infinite sets.  We consider
a partition of $U$ to be \emph{finitary} if all its parts are finite.
If $\pi$ is a finitary partition of $U$, there is a unique integer $n(\pi)$
such that all but finitely many parts of $\pi$ have size~$n(\pi)$. The number $n(\pi)$ is the \emph{gauge}
\index{gauge of partition}%
of $\pi$. We say that $U$ is \emph{bounded} or \emph{unbounded} according
as the gauges of its finitary partitions are. If $U$ is bounded, the
\emph{gauge} of $U$ is the largest gauge of a non-trivial partition of $U$.  If the gauge is $1$ then the set is strictly amorphous.

An example will illustrate what can happen when we try to extend this
last result to the permissible partitions of bounded sets.  Suppose first that the gauge is~$2$, so that there
is a partition of $U$ with all but finitely many points lying in parts of
size~$2$. Now there is a permutation of $U$ which interchanges the points in
each part of size~$2$ and fixes all other points. (No choice principle is
required for this). This element lies in $\Sym(U)$ but not in $\FSym(U)$,
so that the factor group has order at least~$2$ (and it is easy to see that
equality holds).

Next suppose that the bound is~$3$. There are two possible cycles on a set
of size~$3$, corresponding to its two cyclic orientations. If there is a
choice function for the set of orientations, we can use it to construct
a permutation in which all but finitely many points lie in $3$-cycles, so
that $|\Sym(U)/\FSym(U)|=3$. But if no such choice  function exists, then
$\Sym(U)=\FSym(U)$.

The final theorem in this section does not rely on any particular FM model.

\begin{theorem}
If a bounded amorphous set $U$ has no partition into $m$-sets and
finitely many other points, for $m >
n$, then 
\[ | \Sym(U) : \FSym(U) | \leq n. \]
\end{theorem}

\begin{proof}
Suppose all partitions of $U$ have gauge $\le n$, and suppose for
a contradiction that $g_1, g_2, \ldots, g_{n+1} \in \Sym(U)$ lie
in distinct right cosets of $\FSym(U)$.  Let $\pi_i$ be the
partition consisting of the orbits of $g_i$.  By Lemma~$2.1$ of~\cite{trussa} there
is a finitary partition $\pi$ such that each $\pi_i$ refines $\pi$.
By assumption, $\pi$ has gauge $\le n$.  Now for $i \neq j,\
g_ig_j^{-1}$ is not finitary; denote its orbit partition by
$\pi_{ij}$.  The partition $\pi$ can be chosen without loss of
generality to dominate all others including $\pi_{ij}$.  Hence $\exists
X \in \pi$ of size $\le n$ on which $g_ig_j^{-1}$ $(i \neq j)$ has no fixed points in $X$.  Pick $x \in X$.
Then $xg_i$ for $1 \le i \le n+1$ are distinct members of $X$, contradiction.
\end{proof}

\subsection{The Centralizer of the Gauge Group}

In this section we assume that we are working in the particular models of FM
mentioned in the Main Theorem (Theorem~\ref{mainth}) below.  Truss~\cite{trussa} showed that a bounded amorphous set of gauge~$n$ is
determined, up to a notion of equivalence which is somewhat difficult to
explain, by two parameters: the \emph{gauge group}
\index{gauge group}%
$G_0$ (a transitive subgroup of the symmetric group $\Sym(n)$), and the \emph{excess} (an
integer from the set $\{0,\ldots,n-1\}$). In this section, by examining
Truss's classification more closely, we will identify the finite group
$\Sym(U)/\FSym(U)$ with the centralizer of the gauge group.

There are two aspects of Truss's work which we depend on. First, he shows how the gauge group
and the excess are determined by the set~$U$.

Second,
he constructs a model of Fraenkel-Mostowski set theory
\index{Fraenkel-Mostowski set theory}%
 in which the set
of atoms is amorphous and has any prescribed gauge group and excess, and shows
how a corresponding model of ZF can be obtained using the method of
Jech and Sochor~\cite{js}.
\index{Jech, T.}%
\index{Sochor, A.}%

Note that, for any infinite set $U$, and any point $u\in U$, we have
\[\Sym(U)/\FSym(U) \cong \Sym(U\setminus\{u\})/\FSym(U\setminus\{u\}),\]
since $\FSym(U)$ is transitive and so meets every coset of
$\Sym(U\setminus\{u\})$ in $\Sym(U)$. So it is enough to consider sets
with excess equal to zero.

The gauge group is determined by $U$ (up to conjugacy), so we can
write it as $G(U)$.  Its construction is given
in~\cite[p.~199]{trussa}, where it is also shown that $G(U)$ acts
transitively on $\{0, 1,\ldots, n-1 \}$.  Truss further shows~\cite[Corollary~5.11]{trussa} that any bounded
amorphous set takes the form of the above construction.  In~\cite[Theorem 5.10]{trussa} Truss shows the type of transitive model of FM in which the sets of the main theorem lie.

The permutation group acting finitarily on elements of $U$ is $\mathcal{G} = (G\ Wr \Sym(\omega)) \cap \FSym(U) =
G\ wr \FSym(\omega)$ where $wr$ means the restricted wreath product.

\begin{theorem}[Main Theorem]
\label{mainth}
Let $\mathcal{N}$ be the transitive model of FM obtained using
the Truss construction~\cite{trussa} and let $U$ be a bounded amorphous set lying in
$\mathcal{N}$ and having gauge~$n=n(U)$ and gauge group
$G(U)$. Then
\[\Sym(U)/\FSym(U)\cong C_{\Sym(n)}(G(U)).\]
\end{theorem}

\begin{proof}
The proof of the theorem relies on the construction of
bounded amorphous sets by Truss.
\index{Truss, J. K.}%
In outline it begins with a model $\mathcal{M}$ of Fraenkel-Mostowski set theory with choice (FMC);
\index{Fraenkel-Mostowski set theory}%
a modification of ZFC allows a set $U$ of atoms. A \emph{permutation
model}
\index{permutation model}%
$\mathcal{N}$ is constructed using a permutation group $G$ on $U$
and a filter
\index{filter}%
$\mathcal{F}$ of subgroups of $G$. The Fraenkel-Mostowski
model $\mathcal{N}$
is defined recursively
by the rule
\[\mathcal{N}=\{x\in\mathcal{M}:x\subseteq\mathcal{N}\hbox{ and }G_{\{x\}}\in
\mathcal{F}\},\]
where $G_{\{x\}}$ denotes the setwise stabilizer of $x$ in $G$.  The
bounded amorphous sets of the theorem lie in a specific model
$\mathcal{N}$ of all `hereditarily symmetric' sets.  The paper~\cite{trussa} can be
referred to for a detailed exposition of the construction.  We need the construction to show that every finite group
can occur here.  Usually it is assumed that the filter $\mathcal{F}$
\index{filter}%
 contains the stabilizers of the elements of $U$. This guarantees that the elements
of $U$ (and hence $U$ itself) all belong to the model~$\mathcal{N}$.

It can also be shown that, under suitable hypotheses, most results
obtained by FM methods, for example in $\mathcal{N}$, can
be `transferred' to a model of ZF with similar properties. The first
general result that such a transfer was possible is the
theorem of Jech and Sochor~\cite{js}. In all the cases that we need, the
verification of the hypotheses has been done already (see Truss~\cite{trussa}),
so we have nothing to do in this connection.

The proof of our main theorem uses the following three propositions.  The proof of the first is in~\cite{tarzi} and that of the second is in~\cite{trussa} and uses a restricted wreath product.

\begin{proposition}
\label{firstprop}
Suppose that a permutation model $\mathcal{N}$ is defined by a set $U$ of
atoms, a permutation group $G$ on $U$, and a filter $\mathcal{F}$ of
subgroups of $G$ including the point stabilizers. Then $C_{\Sym(U)}(G)\in
\mathcal{N}$.
\end{proposition}

\begin{proposition}
\label{secondprop}
Let $G_0$ be a transitive subgroup of the symmetric group $\Sym(n)$ in
a model $\mathcal{M}$ of ZF, and
let $G=G_0\ wr\Sym(\omega)$, $U=n\times\omega$, and $\mathcal{F}$ the filter
generated by the point stabilizers in $G$. Then in the corresponding
permutation model $\mathcal{N}$, $U$ is a bounded amorphous set with gauge~$n$, gauge
group $G_0$, and excess~$0$.
\end{proposition}

\begin{proposition}
\label{isom}
With the notation of the preceding proposition,
\[ C_{\Sym(U)}(G) \cong C_{\Sym(n)}(G_0).\]
\end{proposition}

\begin{proof}
It is clear that, if $H_0=C_{\Sym(n)}(G_0)$, and we extend the action of
$H_0$ to $U=n\times\omega$ by letting it act on the first factor, then
$H_0\le C_{\Sym(U)}(G)$.

Conversely, take any element $h$ of $H=C_{\Sym(U)}(G)$. We claim first that
$h$ fixes every set $U_i=n\times\{i\}$ for $i\in\omega$. For there is a
factor $G_i$ of the base group of $G=G_0\ wr\Sym(\omega)$ which acts
transitively on $U_i$ and fixes all points of $U_j$ for $j\ne i$; if
$h$ maps $(x,i)$ to $(y,j)$ for some $j\ne i$, then the condition
$gh=hg$ for $g\in G_i$ would imply that every point of $U_i$ would be
mapped to $(y,j)$, a contradiction.

Now, using the fact that $h$ commutes with the top group $\Sym(\omega)$,
\index{group ! top}%
we find that $h$ acts in the same way on each set $U_i$; so $h\in H_0$,
and the proof is complete.
\end{proof}

By Proposition~\ref{secondprop}, the bounded amorphous set can be
realized as a countably infinite number of size $n$ pieces in the
model $\mathcal{M}$, and also $G_0$ is the gauge group $G(U)$. 

The theory of Truss requires that the models in which we are
working are transitive, that is that any member of a member of
$\mathcal{N}$ is a member of $\mathcal{N}$.  But this follows from the
fact that $\mathcal{N}$ contains all hereditarily symmetric sets.  

Let $g \in \Sym(U)$ in $\mathcal{M}$.  Now by definition $g \in \mathcal{N}$ if and
only if the stabilizer of $g$ in $G=G_0\ wr\Sym(\omega)$ is in the
filter
\index{filter}%
 $\mathcal{F}$.  A permutation is a set of ordered pairs $g =
\{(u, u^g) : u \in U\}$.  The stabilizer of $g$ is the set $\{h : (u,
u^g)^h = (v, v^g) \}$ for some $v$ in the operand.  But then $(u^h,
u^{gh}) = (v, v^g)$, and so $v = u^h$, $(u^{h})^{g}=v^g = u^{gh}$, that is
$h \in C_{G}(g)$.  It follows that $g \in \mathcal{N}$ if and only if
$C_{G}(g)$ contains $G_{(u_1, \ldots, u_n)}$ for some $u_1,
\ldots, u_n \in U$.

Now $\Sym(U) = \{g: C_{G}(g)
\ge G_{(u_1, \ldots, u_n)}\ \text{for\ some}\ u_1,
\ldots, u_n \in U\} = \{g: g \in C_{\Sym(U)}(G_{(u_1, \ldots, u_n)})\
\text{for\ some}\ u_1, \ldots, u_n \in U\}$, where in the first set
$\Sym(U)$ lies in $\mathcal{N}$ and in the second set
$\Sym(U)$ lies in $\mathcal{M}$.  Take $g \in \Sym(U)$ (in
$\mathcal{N}$).  Then $g$ fixes the set $\{u_1, \ldots, u_n\}$ and
fixes the set of fibres containing  $u_1, \ldots, u_n$.  On the
remaining fibres,
$g \in C_{\Sym(\tilde{U})}(\tilde{G})$, where $\tilde{G}$ is the group
acting on the remaining fibres as $\tilde{G} \cong G$, $\tilde{U}$ is
the set $U$ minus the fibres containing $\{u_1, \ldots, u_n\}$ and $\Sym(\tilde{U}) \in \mathcal{M}$.  By the proof
of Proposition~\ref{isom}, $g$ fixes all the fibres in $\tilde{U}$ and acts
on each as an element $\overline{g} \in C_{\Sym(n)}(G_0)$.  So we have a map (in $\mathcal{N}$) $\phi : \Sym(U) \to C_{\Sym(n)}(G_0):
g \mapsto \overline{g}$.  Here $\ker(\phi)$ consists of permutations
which, apart from finitely many fibres, are trivial, that is $\FSym(U)$.
We know that $\FSym(U) \in \mathcal{N}$ because if $g' \in
\FSym(U)$ then $g'$ moves a finite set of points $\{u_1, \ldots, u_n\}$
and $g'$ is centralized by $G_{(u_1, \ldots, u_n)} \in \mathcal{F}$.  The homomorphism $\phi$ is onto because by Proposition~\ref{isom}, for
every element of $C_{\Sym(n)}(G_0)$ there is an element of
$C_{\Sym(U)}(G)$ ($\Sym(U) \in \mathcal{M}$), and each of these lies
in $\Sym(U)$ ($\in \mathcal{N}$).  So
in $\mathcal{N}$, $\Sym(U)$ is an extension of $\FSym(U)$ by a group
isomorphic to $C_{\Sym(n)}(G_0) \cong C_{\Sym(U)}(G)$, where $\Sym(U)$
in the latter group lies in $\mathcal{M}$.  The product is
semidirect because $\FSym(U) \lhd \Sym(U)$ while
$C_{\Sym(U)}(G)$ is semiregular,
\index{group ! permutation ! semiregular}%
 so $\FSym(U) \cap C_{\Sym(U)}(G) =
1$.  Note that taking $n=0$ in
$\Sym(U)  = \{g: C_{\Sym(U)}(G_{(u_1, \ldots, u_n)})\}$ gives that
$C_{\Sym(U)}(G) \in \mathcal{N}$; this is also proved in Proposition~\ref{firstprop}.  Finally, both $C_{\Sym(n)}(G_0)$
and $C_{\Sym(U)}(G)$ are finite groups and so lie in $\mathcal{N}$ and
the isomorphism between them also lies in the model because from the
axioms of set theory, a bijection between finite sets lies in the
model.  Therefore $\Sym(U)/\FSym(U) \cong
C_{\Sym(n)}(G_0)$ and so Theorem~\ref{mainth} is proved.
\end{proof}

\begin{corollary}[Baer--Schreier--Ulam Theorem for bounded\\ amorphous sets]
\index{Baer--Schreier--Ulam Theorem}%

Every finite group arises as the quotient\\ $\Sym(U)/\FSym(U)$ for some bounded amorphous set~$U$.
\end{corollary}

\begin{proof}
The main theorem immediately gives that the quotient $\Sym(U) / \FSym(U)$ is finite
for the FM models of Theorem~\ref{mainth}.  Obviously, any finite group of order dividing~$n$ can be embedded as a semiregular
\index{group ! permutation ! semiregular}%
 subgroup in $\Sym(n)$ in some transitive model of FM.
\end{proof}

Usually the non-AC equivalent of an AC theorem is a more restrictive
result; there is a sense in which for this theorem the reverse is true.

More results can be found in~\cite{tarzi}, which contains a study of group actions
\index{group ! action}%
 on bounded and more general amorphous sets.

It is possible to build graphs from amorphous vertex sets, but they would have amorphous cardinality.  Theorems such as that of Lachlan-Woodrow
\index{Lachlan-Woodrow Theorem}%
\index{Lachlan, A. H.}%
\index{Woodrow, R. E.}%
that classifies countable homogeneous graphs (see Chapter~\ref{chapFD}) fail, not only because the notion of countable no longer exists, but also because in the absence of the axiom of choice it is no longer certain that finite partial automorphisms
\index{automorphism ! partial}%
 will extend to total automorphisms, and so the notion of homogeneity is no longer appropriate.  In fact it can be replace by ``local homogeneity''~\cite{trussa}.

The proof of~\cite[Theorem~2.1.2]{hod1} demonstrates that, without the axiom of choice, there can be amorphous sets.

Finally, we mention one problem where the above result has had a bearing.  The standard analysis showing that subgroups of bounded support are the only non-trivial normal subgroups of infinite symmetric groups,
\index{group ! symmetric ! infinite}%
 and that there are no nontrivial normal subgroups of small index, uses AC.  Bowler and Forster
\index{Bowler, N.}
\index{Forster, T. E.}
 proved~\cite{bowlerf} that if $|X| = \omega$ with $|X|=|X\times X|$ then $\Sym(X)$ has no nontrivial normal subgroups of small index.  The corollary that proves the Baer--Schreier--Ulam Theorem for bounded amorphous sets indicates the necessity of some extra condition; in this case $|X|=|X\times X|$ was adopted.  Their two-stage proof, first shows that every normal subgroup of $\Sym(X)$ (of small index) must contain every \emph{flexible} permutation,
\index{permutation ! flexible}%
that is, one that fixes at least $|X|$-many points.  Then they proved that the flexible permutations generate $\Sym(X)$.  Their interest was in the universes of typed set theories and of Quine's New Foundations (NF)
\index{Quine, W.}%
\index{NF set theory}%
 which satisfy the $|X|=|X\times X|$ condition, for if the symmetric group for such a universe has no non-trivial normal subgroups of small index then there is a simple theory of sets definable in such universes.  Note that AC fails in NF, though New Foundations with urelements (NFU) is a consistent system, and is consistent with Choice.  They applied their result to typed set theories, bounding the orbit sizes under the induced action of the symmetric group of the universe.  They derived a wellfoundedness result for small symmetric sets in these contexts, and a result limiting the sizes of wellfounded sets in NF.
 
For more on NF set theory see the book~\cite{forster} by Thomas Forster and the paper~\cite{holmes} by Randall Holmes.
\index{Holmes, R.}

\chapter{Further Directions}
\label{chapFD}
The will is infinite and the execution confined \ldots the desire is
boundless and the act a slave to limit.
\begin{flushright}
William Shakespeare, \textit{Troilus and Cressida, Act 3 scene 2}
\end{flushright}

\medskip

If the doors of perception were to be cleansed man would see everything as it truly is... Infinite. 

\begin{flushright}
William Blake, \textit{The Marriage of Heaven and Hell, (1790-3)\\ A
Memorable Fancy, plate 14} 
\end{flushright}

\bigskip

This chapter is an ensemble of notes together with suggested problems and speculations that are of varying degrees of difficulty and requiring differing levels of effort.  They are intended to motivate interest in the subject and encourage further work whilst appealing to differing tastes.  Some of the questions have a long
introduction both for guidance and as an opportunity to introduce or remind the reader of certain theories.  Not all will necessarily have positive solutions.  Some arise directly from some of the work in the main text whilst others intersect with rather different parts of mathematics.
\bigskip
\begin{enumerate}

\item  \emph{Classification of Countable Homogeneous Coloured Graphs}.  The
 Lachlan-Woodrow Theorem 
\index{graph ! random ! $m$-coloured}%
\index{Lachlan-Woodrow Theorem}%
\index{Lachlan, A. H.}%
\index{Woodrow, R. E.}%
for graphs~\cite{lach} states that:
\begin{theorem}
A countable homogeneous graph is isomorphic to one of the following:

(i) the disjoint union of $p$ complete $n$-vertex graphs, where at
least one of $p$ and $n$ is infinite, or its complement a $p$-partite
graph;
\index{graph ! ppartite@$p$-partite}%
(ii) a Henson graph $\mathfrak{H}_k$ ($n \geq 3$),
or the complement (see Chapter~\ref{chap2});
\index{graph ! Henson}%
(iii) the random graph $\mathfrak{R}$, (isomorphic to its complement).
\end{theorem}
\index{graph ! random}%

A graph $\Gamma$ is said to be \emph{distance-transitive}
\index{graph ! distance-transitive}%
 if , given any vertices $a, b, c, d$ such that $d(a, b)=d(c, d)$, there is an automorphism $\alpha$ of $\Gamma$ such that $\alpha(a)=c$ and $\alpha(b)=d$.  In terms of conditions, intermediate between the homogeneous and distance-transitive graphs are the \emph{connected-homogeneous graphs}, being those for which any isomorphism between \emph{connected} finite induced subgraphs extends to a graph automorphism.  The connected countably-homogeneous graphs were classified by Gray and Macpherson~\cite{graymac}.
\index{Macpherson, H. D.}%
\index{Gray, R.}%

Related to these are the \emph{distance-homogeneous graphs}~\cite{moss} \cite{camab} 
\index{graph ! distance-homogeneous}%
which are countable and homogeneous in a language that uses binary relations to encode distance within the graph, that is there is a binary predicate $P_k$ such that $P_k(v_1, v_2)$ holds if and only if $d(v_1, v_2) = k$.  Note that if $\Gamma$ is distance-homogeneous then $\Gamma(v)$ is homogeneous, since non-adjacent pairs in $\Gamma(v)$ are at distance 2.  So the Lachlan-Woodrow Theorem classifies the countable distance-homogeneous diameter 2 graphs.

Find a multicoloured analogues to these notions and classify countable homogeneous
$m$-coloured graphs.
\index{graph ! random ! $m$-coloured}%

\item \emph{Countable Homogeneous Metric Spaces}.
\index{metric space ! countable homogeneous}%
A graph of diameter $2$ is the same as a metric space in which the metric takes only the values 1 and 2. The graph $\mathfrak{R}$ is the unique countable homogeneous metric space with these properties. 
By the same methods we can construct countable universal 
homogeneous metric spaces with other sets of values of the 
metric: 
\begin{itemize}
\item $\{1, 2, . . . , d\}$ for any $d \geq 2$; 
\item the positive integers; 
\item the positive rationals. 
\end{itemize}
In the first two cases we can modify the construction to produce the analogue of Henson's graph
\index{graph ! Henson}%
(that is, one that has no equilateral triangles with side 1), or a bipartite graph (that is, in which all triangles have 
even perimeter). 
Problem: classify the countable homogeneous metric spaces?

\item   \emph{Hypermetric Space Theory}.
\index{hypermetric space}%
M. Deza
\index{Deza, M.}%
and co-workers have studied~\cite{dezagri} hypermetric spaces, which can be represented uniquely up to a multiple by a distance $d_{G,t}$ for a graph $G$ where $d_{G,t}(i,j) = 1$ or $t$ respectively according to whether $(i,j)$ is an edge or a non-edge in $G$.  Is there a universal hypermetric space?

\item \emph{Some Transitivity Properties of Reducts of $\mathfrak{R}_{m,\omega}$}.
\index{graph ! random ! $m$-coloured}%
  H. Wielandt
\index{Wielandt, H.}%
 introduced variations of $k$-transitivity with high transitivity being the strongest condition; it is also stronger that oligomorphicity.  In this item we outline two variations on the theme of multiple transitivity, one due to $\Pi$. M. Neumann
\index{Neumann, $\Pi$. M.}%
and the other due to Neumann and C. Praeger,
\index{Praeger, C. E.}%
and ask which random graph reducts admit either of them, or their relatives described below?  

A permutation group $G$ acting on a set $X$ is \emph{generously transitive}~\cite{neumann}
\index{group ! permutation ! generously transitive}%
if any two distinct points of $X$ are interchanged by some element of $G$.  More generally, a group $G$ is generously $k$-transitive
\index{group ! permutation ! generously $k$-transitive}%
 if for every $(k+1)$-element set $X$, $G^{X}_{\{X\}} = \Sym(k+1)$.  Generous transitivity is a strengthening of multiple transitivity so that 
 
\begin{proposition}
If $G$ is a generously $k$-transitive group on $|X| > k$ then it is $k$-transitive.
\end{proposition}

\begin{proof}
Let $G$ be a generously $k$-transitive group.  Then

(a)  $G$ is generously $(k-1)$-transitive, so $G^{X}_{\{X\}} = \Sym(k)$ if $|X|=k$.

(b)  If $X, Y$ are $k$-sets with $|X \cap Y| = k-1$ then $|X \cup Y| = k+1$ so there is an element $g \in G_{X \cup Y}$ mapping $X$ to $Y$.

(c)  Then by induction on $d$ where $X, Y$ are $k$-sets with $|X \cap Y| = k - d$, there exists $g \in G$ mapping $X$ to $Y$ (via $d-1$ intermediate sets which sequentially pairwise intersect at all but one of their points).

So $G$ is $k$-homogeneous (set-transitive).
\index{group ! permutation ! set-transitive}%
Finally (a) and (c) imply that $G$ is $k$-transitive.
\end{proof}
 
An \emph{orbital} of a permutation group
\index{group ! permutation ! orbital}%
$G$ on $\Omega$ is an orbit of $G$ on $\Omega \times \Omega$, and for a transitive group there is a bijection between its orbitals and the orbits of a point stabilizer; these latter are called \emph{suborbits}.
\index{group ! permutation ! suborbit}%
A $G$-orbit (subset of $\Omega \times \Omega$) is \emph{self-paired}
\index{group ! permutation ! self-paired orbit}%
if and only if for any (or all) vertex pairs $(v_1, v_2)$ there exists $g \in G$ such that $(v_1, v_2) g = (v_2, v_1)$.  For an edge-coloured graph an orbital is equivalent to an edge-colour.  A group $G$ is generously transitive if and only if it is transitive and all its suborbits are self-paired.  The
automorphism groups
\index{group ! automorphism}%
 of all homogeneous edge-coloured graphs, are generously
transitive because the edges are undirected so an interchange of a
vertex pair on an edge is an isomorphism that extends to an
automorphism of the whole graph.  So the groups
$\Aut(\mathfrak{R}_{m,\omega})$ are all generously transitive for all $m$.  This is not true for all homogeneous structures, and fails for example for $\mathbb{Q}$.  In fact the group $\Aut(\mathbb{Q}, \le)$~\label{autqle} of all order preserving automorphisms of the rationals is highly set-transitive (and primitive and oligomorphic) but is not even 2-transitive.

A permutation group is a \emph{three-star group}
\index{group ! permutation ! three-star}%
 if it induces a non-trivial group on each 3-element subset of points~\cite{neumpra}.  In brief then, a group is generously transitive if every pair of points of the operand can be transposed whilst it is three-star if it effects some non-trivial permutation of any triangle.  A primitive three-star group other than $\Alt(3)$ is generously transitive.  Three-star transitivity and 2-transitivity are weakenings of generous 2-transitivity in different directions and there is not always an implication from one to the other.  However generously transitive and 2-homogeneous implies 2-transitive.  Three-star is upwards closed for homogeneous structures in general but three-star does not necessarily imply four-star which is similarly defined, as is exemplified by the following.  Take a 4-ary relation $R$ such that $R(a, b, c, d)$ holds if and only if $a, b, c, d$ are distinct and $R(a, b, c, d) \Rightarrow \neg R(ga, gb, gc, gd)$ for any permutation $g$ of $a, b, c, d$, and for any 4 points it holds in one possible order.  If $M$ is the corresponding Fra\"{\i}ss\'e structure and $G = \Aut(M)$ which is 3-transitive then $G^{X}_{\{X\}} = 1$ for $|X| = 4$.  So we have constructed a 3-transitive group $G$ such that the setwise stabilizer of any 4-set is trivial. 

We begin applying this theory to random graph reducts.  The group $\Aut(\mathfrak{R})$ is vertex transitive but not 2-transitive; that it is generously transitive follows because it is transitive and all its suborbits
\index{group ! permutation ! suborbit}%
 are self-paired by definition of automorphism group.
\index{group ! automorphism}%
  The duality group $\DAut(\mathfrak{R})$ acts 2-transitively but not 3-transitively on $\mathfrak{R}$ because vertex triples containing 0 and 3 edges are not duality-equivalent to those containing 1 and 2 edges.  The switching group $\SAut(\mathfrak{R})$ acts 2-transitively but not 3-transitively on $\mathfrak{R}$ because vertex triples containing 0 and 2 edges  lie in different switching classes to those containing 1 and 3 edges.  

\begin{figure}[!h]$$\xymatrix{
{\bullet} \ar@{-}[d] \ar@{-}[drr] \ar@{-}[rr] && {\bullet} \ar@{-}[d] \ar@{-}[dll] &&& {\bullet} \ar@{-}[d] \ar@{-}[drr] && {\bullet} \ar@{-}[dll] \\
{\bullet} \ar@{-}[rr] && {\bullet} &&& {\bullet} \ar@{-}[rr] && {\bullet} 
}$$
\caption{Failure of 4-transitivity of $\BAut(\mathfrak{R})$}
\label{foursetdi}
\end{figure}
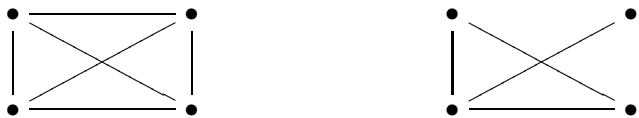

The biggest group $\BAut(\mathfrak{R}) = \langle \DAut(\mathfrak{R}), \SAut(\mathfrak{R}) \rangle$ acts 3-transitively because the four possible configurations on 3 vertices, shown in Figure~\ref{thstr}, are equivalent under $\BAut(\mathfrak{R})$,
but it does not act 4-transitively because for example the two 4-sets  in Figure~\ref{foursetdi} are not $\BAut(\mathfrak{R})$-equivalent:

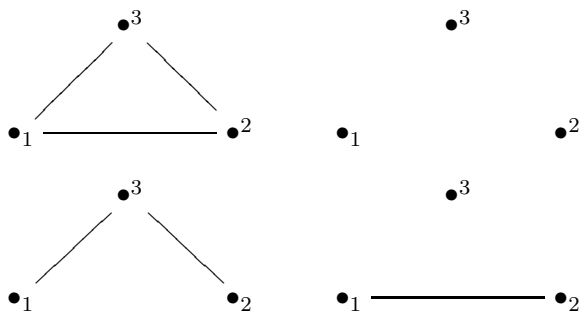
\begin{figure}[!h]$$\xymatrix{
& {\bullet}^3 \ar@{-}[dr] \ar@{-}[dl] &&& {\bullet}^3 \ar@{}[dr] \ar@{}[dl]\\
{\bullet}_1 \ar@{-}[rr] && {\bullet}^2 & {\bullet}_1 \ar@{}[rr] && {\bullet}^2
}$$
$$\xymatrix{
& {\bullet}^3 \ar@{-}[dr] \ar@{-}[dl] &&& {\bullet}^3 \ar@{}[dr] \ar@{}[dl]\\
{\bullet}_1 \ar@{}[rr] && {\bullet}_2 & {\bullet}_1 \ar@{-}[rr] && {\bullet}_2
}$$\caption{Three-star transitivity of $\Aut(\mathfrak{R})$}
\label{thstr}
\end{figure} 

\item  \emph{Pigeonhole Property}.
\index{pigeonhole property}%
  The pigeonhole property (PP) is a very strong indivisibility property of a relational structure, and states that whenever $X$ is the union of two disjoint induced substructures $Y$ and $Z$, at least one of $Y$ and $Z$ is isomorphic to $X$.  Induced means that the instances of relations in $Y$ are all instances in $X$ for which all the arguments lie in $X$.  The empty and singleton sets have this property, but no larger finite structure does, since it can always be divided into two strictly smaller parts.  Assuming the axiom of choice,
\index{axiom of choice}%
  any infinite structure has this property.
   
 The PP is often useful in theorem-proving.  For  example in~\cite{mozz} it is reported how the gap that initially had arisen in the proof~\cite{wiles}~\cite{taywil} of Fermat's Last Theorem was closed by  relating the infinite collection of Hecke rings
\index{Hecke ring}%
 by creating an infinite sequence of sets of pigeonholes and then showing that there must be objects arising in every set of pigeonholes.  
     
 To prove that $\mathfrak{R}$ has the PP using the $(*)$-condition, suppose that $\mathfrak{R}$ can be split into two pieces $Y$ and $Z$, neither of which is isomorphic to $\mathfrak{R}$.  Since $Y \ncong \mathfrak{R}$, there are finite sets $U_1$ and $V_1$ for which no ``correctly joined'' witness exists in $Y$.  Similarly there are sets $U_2$ and $V_2$ with no witness in $Z$.  Now put $U = U_1 \cup U_2$ and $V = V_1 \cup V_2$. Since the whole graph is isomorphic to $\mathfrak{R}$, there is a witness for $U$ and $V$; but, by assumption, this witness cannot lie in either $Y$ or $Z$, a contradiction.
 
 The PP for a non-trivial graph states that
     for every finite partition of its vertex set, the induced
     subgraph on at least one of the blocks is isomorphic to the
     graph.  Up to isomorphism, the three simple countably infinite
     graphs having PP are the complete graph $K_{\omega}$, its
     complement $\overline{K}_{\omega}$, and $\mathfrak{R}_{2, \omega}$.
\index{graph ! random}%
 To prove this, let $X$ be a countable graph with the PP. Partitioning $X$ into the set of isolated vertices and the rest, we conclude that either $X$ is null or it has no isolated vertices; we may assume the latter. Similarly we may assume that there is no vertex joined to all others.  Suppose that $X \ncong \mathfrak{R}$; let $U$ and $V$ be finite sets having no correctly joined witness, and suppose that $|U \cup V|$ is minimal subject to this.  We just concluded that this cardinality is at least 2.  Assume that $U \neq \emptyset$ and choose $u \in U$.  Then the set $W$ of vertices joined to everything in $U$ except $u$ and to nothing in $V$ is non-empty; let $Y = W \cup \{u\}$, and $Z$ the remaining points.  Then $Y$ is not isomorphic to $X$; for the non-existence of a witness for $U$ and $V$ shows that $u$ is joined to nothing in $Z$, so it is an isolated vertex in $Y$.  Also, $Z$ is not isomorphic to $X$, because the sets $U\backslash \{u\}$ and $V$ have no witness in $Z$ (all the witnesses lie in $Y$).

 If we allow our language to have $3$ binary relations we have seen that $\mathfrak{R}_{3, \omega} = \mathfrak{R^{t}}$ has the PP.

If an \emph{oriented graph},
\index{graph ! oriented}%
that is a directed graph not containing two vertices joined in both directions by edges, were to have the PP, then the undirected underlying graph is one of the three listed above.  So an oriented graph with the PP is null, a tournament,
\index{tournament}%
 or an orientation of $\mathfrak{R}$.

The countable random tournament,
\index{tournament ! random}%
arises with probability 1 if we choose the orientations independently at random.  It shares many properties with the random graph, including a very similar number-theoretic construction: take the primes congruent to $-1\pmod{4}$, and put an arc from $p$ to $q$ if $q$ is a quadratic residue mod $p$. A similar argument to the one with which we began shows that the random tournament has the PP.

Any countable totally ordered set whose order type is a limit ordinal (one with no greatest element) has the PP, for example, $\mathbb{N}$ with arcs directed from smaller to larger, or vice versa.  In~\cite{bonat1} it was shown that these are the only possibilities for tournaments with the PP.

In~\cite{bonat} the countable tournaments
\index{tournament}%
and orders with PP were classified.  It should be easy to show that $\mathfrak{R}_{m, \omega}$
  \index{graph ! random ! $m$-coloured}%
 has the PP for all $m \geq 2$.
     Further, all binary relational structures
\index{relational structure}%
   with $m$ different binary symmetric relations having PP have been classified and these correspond to multicoloured random graphs.
\index{graph ! random ! $m$-coloured}%
 We do not know if this has been done in the case of relations that are not symmetric.  Not every homogeneous graph has PP, so settling this question would not suffice to give the classification of the previous question.

 As for the orientations of $\mathfrak{R}$, there is a unique countable random oriented graph, which is an orientation of $\mathfrak{R}$, and which has the pigeonhole property.  Diestel et al.~\cite{diestel} constructed and characterised an acyclic orientation of $\mathfrak{R}$ which has the pigeonhole property, and showed that this exhausts the possibilities.  They classify the countably infinite oriented graphs which for every partition of their vertex sets into two parts, induce an isomorphic copy of themselves on at least one of the parts.  These are the edgeless graph,
\index{graph ! edgeless}%
 the random tournament,
\index{tournament ! random}%
 the transitive tournaments
\index{tournament ! transitive}%
 of order type $\omega^{\alpha}$, and two orientations of  $\mathfrak{R}$: the random oriented graph,
\index{graph ! random ! oriented}%
 and the newly-found random cyclic oriented graph.
\index{graph ! random ! cyclic oriented}%

A further generalization would ask for a structure with several (but only finitely many) binary relations.  Adding new relations if necessary, we can assume that (i) all the relations are non-empty, that (ii) every ordered pair of distinct points satisfies exactly one of the relations, and that (iii) each relation is either symmetric (a graph) or skew-symmetric (an oriented graph), and in the latter case its converse is also a relation in our set.  This accounts for the case of a general directed graph. There are two symmetric relations, ``not joined'' and ``joined in both directions'', and the remaining arcs and their reversals give a pair of skew-symmetric relations.

Suppose first that all the $m$ relations are symmetric.  Then, arguing as before, we can show that there is a unique structure, the random $m$-coloured complete graph.
\index{graph ! random ! $m$-coloured}%
 In this structure, each monochromatic subgraph is isomorphic to $\mathfrak{R}$.

Next, observe that in a set of relations as just defined, we can replace each converse pair of skew-symmetric relations by a symmetric relation without losing the PP.  When this is done, we obtain the random $m$-coloured complete graph.  Moreover, the original skew-symmetric relation is then an orientation of $\mathfrak{R}$ with the PP, so we know the possibilities for its structure.  One open question would be to put these pieces together. 

Here are some variations on the same theme.  Structures whose points are distributed among three pigeonholes, one of which contains a copy of the original have a form of the PP that is exactly equivalent to the original.  That is, if $X$ has the PP and is divided into three, then either the first pigeonhole has a copy of $X$, or the union of the other two does; in the latter case, one further application of the property gives the result. In the other direction, simply leave the third pigeonhole empty.

In~\cite{bonat1} a variant was studied which turns out to be quite different.  Suppose that, whenever the points are distributed among three pigeonholes, two of them together contain a copy of the original structure.  Further questions can be investigated in order to uncover generalized pigeonhole
 \index{pigeonhole property ! generalized}%
properties along the lines of~\cite{bonat1} where a relational structure $A$ is said to satisfy the $P(n, k)$ property if whenever the vertex set of $A$ is partitioned into $n$ nonempty parts, the substructure induced by the union of some $k$ of the parts is isomorphic to $A$.  

Another open direction for further study would be to investigate the PP of ternary relations.

A weaker version of the pigeonhole property asks only that, if the points of $X$ are distributed in two pigeonholes, then the structure in one of the pigeonholes contains a copy of $X$ as an induced substructure. Such an $X$ is called \emph{indivisible}; this property has been studied by El-Zahar and Sauer~\cite{elzahar},
\index{El-Zahar, M.}%
\index{Sauer, N.}%
 and others. There are too many of them for neat classifications.

Finally there may be another way of looking at the PP of random graphs worth studying, that is in terms of \emph{self-similarity}.
\index{self-similarity}%
  The countable complete, null and random graphs
\index{graph ! random}%
 are the only simple graphs (up to isomorphism) with the property that if the graph vertex set is partitioned into two parts, then the induced subgraph on one of the parts is isomorphic to the graph~\cite{cameron}.  Is there a nice connection between this property of (weak) self-similarity and the property of self-similarity that arises for linear fractals~\cite[Chapter 18]{mandelbrot}, Cantor sets
\index{Cantor set}%
\index{Cantor, G.}%
 and related objects?  Is there any random graph-like object whose evolution can be modeled as a self-similar stochastic process~\cite{embrechts} that picks up this weak self-similarity?

\item \emph{Paris-Harrington construction of $\mathfrak{R}$}.  The countable random graph is characterised by the property that, given any two finite sets $U$ and $V$ of vertices, there is a vertex $z$ joined to everything in $U$ and nothing in $V$.

Does there exist a construction of the random graph in which the witness $z$ is not bounded by any provably computable function?  Presumably this would be a graph whose isomorphism to $R$ would be unprovable in Peano arithmetic.  More specifically, is there a presentation of $\mathfrak{R}$ through which it can be constructed using a Paris-Harrington
\index{The Paris-Harrington Theorem}%
 numbering~\cite{parish} of vertices and edges, but that cannot be proved in Peano arithmetic
\index{Peano arithmetic}%
 to be $\mathfrak{R}$?

\item  \emph{Universal Sequences}.  Here is a construction of $\mathfrak{R}$: choose any universal set $S$ of positive integers, and build the graph whose vertex set is the set of all integers, two vertices being joined if their difference belongs to $S$.  We saw a construction of $\mathfrak{R}$ using prime numbers; here, perhaps, is another. Consider the zero-one sequence $s$ whose $n$th term is 0 if the $n$th odd prime is congruent to $1\pmod{4}$, or 1 if the $n$th odd prime is congruent to $3\pmod{4}$.  Is $s$ universal?

\item  \emph{Paradoxical Decompositions of $\mathfrak{R^{t}}$}.  There is an obvious decomposition of $\mathfrak{R^{t}}$ obtained by putting the vertices of the graph randomly into $l$ classes, where with probability 1 each class contains a copy of $\mathfrak{R^{t}}$.  This is easy because we do not require a group action.
\index{group ! action}%
  It is not paradoxical because we do not have a natural measure on $\mathfrak{R^{t}}$.

It has been said~\cite{frenchr} that the key step in the proof of the Banach-Tarski Theorem
\index{Banach-Tarski Theorem}%
\index{Tarski, A.}%
\index{Hausdorff, F.}%
is the Hausdorff paradox
\index{Hausdorff paradox}%
 which
     generates three disjoint subsets $I, J, K$ of the set of all
     transformations of the sphere using
\index{group ! modular}%
      $\PSL(2, \mathbb{Z}) = \langle
     \sigma, \rho : \sigma^{2} = \rho^{3} = 1 \rangle \cong \mathbb{Z}_{2} *
     \mathbb{Z}_{3}$, such that $\rho(I) = J,\ \rho^{2}(I) = K,\
     \sigma(I) = J \cup K$.  This paradox states that the sphere (minus
     a countable set) can be divided into three disjoint subsets of
     points $A, B, C$ such that $A, B, C$, and $B \cup C$ are pairwise
     congruent, mimicking the pairwise congruence
      \index{congruence}%
 of $I, J, K$ and $J
     \cup K$~\cite{wagon}.  Lindenbaum~\cite{linde} proved that no
     bounded set in the plane can have a paradoxical decomposition.
     As $\mathfrak{R^{t}}$ is not a planar graph, this theorem does not apply.  With Theorem~\ref{mdgptm} in mind, is there a paradoxical decomposition of this graph for the action of the modular group?

\item  \emph{Amenability of $\Aut(\mathfrak{R}_{m,\omega})$}.  There is a connection between Banach-Tarski paradoxes and the \emph{nonamenability} of isometry groups.
A discrete group $G$ is \emph{amenable}
\index{amenable}%
if there is a left-invariant measure
\index{measure ! invariant}%
 $\mu$ on $G$ which is finitely additive and has total measure 1. That is if there is a function \[\mu : \{\mbox{subsets of G} \} \to [0, 1]\] such that

(i)  $\mu(gA) = \mu(A)$ for all $g \in G$ and all subsets $A \subseteq G$,

(ii)  $\mu(G) = 1$, and

(iii)  $\mu(A \cup B) = \mu(A) + \mu(B)$ if $A$ and $B$ are disjoint subsets of $G$.

Finite and abelian groups are amenable but the free group of rank two is not.  Are any of the groups associated with the random graphs, in particular the automorphism groups,
\index{group ! automorphism}%
 amenable?

\item  \emph{A Conjugacy Class Representation of the Random Graph}.
\index{graph ! random}%
  Does there exist a countable group $G$ with the following properties:
\begin{enumerate}
\item  $G$ has just two non-identity conjugacy classes, each inverse-closed,
\item   the Cayley graph $\Cay(G,C)$ is isomorphic to the countable random graph, where $C$ is one of these classes.
\end{enumerate}

We give one clue as to how one approach to this problem, using HNN-extensions (see Appendix~\ref{PermutationGroups}), might proceed.

Recall from Appendix~\ref{LoopTheory} that the loop multiplication group
\index{loop ! multiplication group}%
 $\Mlt(L)$ of a loop $L$~\cite{niemenmaa} is generated by the set of all left and right translations, that is the permutations defined by $L_a (x) = ax$ and $R_a (x) = xa$ for every $x \in L$, and the inner mapping group
\index{loop ! inner mapping group}%
$I(L)$ of $L$ is the stabilizer of the neutral element $e$

For a loop $L$, the transversals $A = \{L_a : a \in L\}$ and $B = \{R_a : a \in L\}$ are said to be \emph{$I(L)$-connected} to $I(L)$ in $\Mlt(L) = \langle A, B \rangle$ if $[A, B] \leq I(L)$.

Suppose we have a loop multiplication group $\Mlt(L)$ whose pointwise stabilizer satisfies
\[ \Mlt(L)_{\alpha} \leq \Aut(\mathfrak{R}_{m, \omega})_{\alpha} \]
and has $m$ orbits on points $\beta \neq \alpha$ that is points joined to $\alpha$ by edges on $m$ colours.  If $L$ were a group and $\alpha = 1$, then $\Mlt(L)_{\alpha} = L$ acting by conjugation $l \mapsto a l b$ stabilizes 1 if and only if $ab = 1$ or $a = b^{-1}$.  So we require that $L$ has at least $m$ non-identity conjugacy classes, that is $\Mlt(L)_{\alpha}$ has at least $m$ orbits.  So is there a group $L$ which has $m$ non-identity conjugacy classes such that $\Mlt(L) \leq \Aut(\mathfrak{R}_{m, \omega})$?  Does the infinite Moufang Loop of Theorem~\ref{moulthm} satisfy analogous conditions?

Right loops can arise as algebraic structures on transversals, and this impacts on the question of multiplication group action
\index{group ! action}%
 of a quasigroup
\index{quasigroup}%
  or loop.  To within right loop isomorphism, every right loop can be obtained from the transversal to the stabilizer of the identity  in the right multiplication group of the loop.

Let $T$ be a \emph{right transversal}
\index{loop ! right transversal}%
 to a subgroup $H$ of loop multiplication group $G = \langle g \rangle$.  The natural $G$ action on $T$, denoted by $\star$ is given by:
\[\quad\quad\quad\quad  t \star g : = u\ \text{where}\ u\ \text{is the representative in}\ T\ \text{of the coset}\ Htg.\]
Restricted to $T$ itself it endows $T$ with a binary operation, so that $T$ is a \emph{magma},
\index{magma}%
denoted $\underline{T}$ which can be thought of as the \emph{right quotient} of $G$ by $H$.  

The following Lemma~\cite{niemenmaa}~\cite{phillips1} characterizes loop multiplication groups:

\begin{lemma} 
\label{loopmultgp}
A group $G$ is the multiplication group of some loop $L$ if and only if it has a subgroup $H$ and two right transversals $A$ and $B$ to $H$ in $G$, such that $[A, B] \leq H$, $\underline{A}$ coincides with $L$, the core of $H$ in $G$ is trivial, and $\langle A, B \rangle = G$, where $\langle A, B \rangle$ is the subgroup generated by $A$ and $B$.
\end{lemma}
We apply this lemma to a loop which we take to be a group.  Let $G_1 = \langle g \rangle$ be the multiplication group of a loop that is a group, and let $G_2 = \langle G_1, t : t^{-1} g t = g^{'} \rangle$ be an HNN-extension (see Appendix~\ref{PermutationGroups}), where $g, g^{'} \in G_1$ have the same order.  Keep $H = G_1$ constant in Lemma~\ref{loopmultgp}.  Then $\Core_{G_2} (G_1) = 1$.  

Elements in $G_2$ are words in $t, g \in G_2 \backslash H, H$.  Note that $t H$ is a left coset of $G_2$ relative to $H$ (or a left coset containing $t$, or a left coset generated by $t$); so words of type $t g h = g' t h$ are left cosets of $G_2$.

Furthermore if $A$ is a transversal for $H$ in $G_1$ and $U$ is a transversal for $G_1$ in $G_2$ then $AU$ is a transversal for $H$ in $G_2$.  To see this, first note that $G u \subseteq (Ha) u$ because every $g \in G$ must lie in a coset of $H$, so $g u = hau$ for some $a \in A$.  The opposite inclusion follows because $G = \bigcup H a$.  Firstly we must check that the HNN-extended $A^{*} = AU$ and $B^{*} = A_1 U_1$ are still transversals.  Secondly we must check that conjugacy classes remain at every HNN extension and that they are \emph{not} conjugated to each other.  Then, iterating the HNN-extension a countable infinity of times and taking the union of the groups gives $G_{\infty}$.  Outer automorphisms permute the conjugacy classes within the Cayley graph of $G_{\infty}$.

Another approach is to ask whether or not there is a $3$-colouring of the projective line over a field $\mathbb{Q}$,~\label{projln} $\mathbf{P}^{1}(\mathbb{Q})$, with all colour classes dense and such that the induced colouring of $\mathbb{Z}^2$ gives the triality graph?
\index{graph ! triality}%
 (We already know that $\mathfrak{R}$ can be derived from $\mathbb{Z}^2$ ~\cite{cam11}).

Note that if we are choosing our field to be $\mathbb{Q}$ then we do not need a density argument; density is a residual property, so is automatically present, and indeed any countable collection of residual properties can be attached. 

More generally, assume that there is a group $G$ for which $G \backslash \{1\} = A_1 \bigsqcup \ldots \bigsqcup A_m$, $A_i^{-1} = A_i$, $\Cay(G, A_1, \ldots, A_m) \cong \mathfrak{R}_{m,n}$.  Then proceed with either the further assumption that

(C)  If $a, b \in A_i$ then $a^n, b^n \in A_{i(n)}$ for some $i(n)$ for all $n \in \mathbb{N}$,

or the stronger assumption

(C')  For any $n \in \mathbb{N}$, $a, a^n$ have the same colour.

\item  \emph{Henson Graphs}.
Henson's graph $\mathfrak{H}_k$, for $k>2$,
\index{graph ! Henson}%
\index{Henson, C. W.}%
is the unique countable homogeneous graph whose finite subgraphs are the finite graphs containing no complete graph of size k.  The groups $\Aut(\mathfrak{H}_k)$ have been shown to be simple by Macpherson and Tent.
\index{Macpherson, H. D.}%
\index{Tent, K.}%
(i) Does $\Aut(\mathfrak{H}_k)$ have the small index property, that is, is it the case that a subgroup of countable index contains the pointwise stabilizer of a finite set?

(ii)  Are the groups $\Aut(\mathfrak{H}_k)$ and $\Aut(\mathfrak{H}_l)$ non-isomorphic for distinct $k$ and $l$?

(iii)  Is $\mathfrak{H}_k$ a Cayley graph for $k > 3$.

\item  \emph{Graph Limits and Random Graphs}.  In Appendix~\ref{FurtherDetails} we mention the relatively recent work on graph limits and the theory of exchangeable measures
\index{measure ! exchangeable}%
 which was begun by Bruno de Finetti
\index{de Finetti, B.}%
 in the 1930s.  Here we ask for links between this theory and the older theory of random graphs with forbidden subgraphs.  We motivate the question.

There are many models for such a graph.  For example, we could add vertices one at a time, and for each new vertex, choose its neighbour set randomly among all the possibilities which don't violate the restriction; in particular
we could pick a random independent set in the existing graph to be the neighbour set of the new vertex.  

Provided the restriction is such that this can always be done, then we have a probability measure on the set of all graphs on a given countable vertex set, with the property that the restriction holds with probability 1.

Alternatively, list the vertex pairs and consider one pair at a time.  If the restriction does not force the pair considered to be a non-edge or an edge then decide at random.  Again, provided we can continue indefinitely, this defines a probability measure as required.  Unfortunately, it may define many different measures; there is no obvious reason why different orderings of the pairs should give the same measure.  Two obvious test cases are lexicographic and reverse lexicographic ordering of the pairs.

The many variants all have the disadvantage that they depend on some ordering of the vertex set (or the set of pairs of vertices, or something).  To motivate finding a measure which did not require a choice of an ordering we ask:  what is the probability that two vertices are joined by an edge in a random triangle-free graph?
\index{graph ! triangle-free}%
  In the ordering models, this depends on the choice of ordering. In the model where the ordering of pairs begins $\{1,2\}, \{1,3\}, \{2,3\}, \ldots,$ the probability that $\{1,2\}$ is an edge is clearly $\frac{1}{2}$, and the same for $\{1,3\}$; but the probability that $\{2,3\}$ is an edge is only $\frac{3}{8}$ (since if we chose both $\{1,2\}$ and $\{1,3\}$, then $\{2,3\}$ would be forbidden).

Taking ``probability'' to mean ``limiting frequency'' suggests that the probability that any pair, say $\{1, 2\}$, is an edge should be the limiting edge-density in large triangle-free graphs.  Here it is convenient for the ``large triangle-free graphs'' to be labelled, that is, to have vertex set $\{1, \ldots,N\}$; then the required edge-density is equal to the proportion of these graphs in which $\{1,2\}$ is an edge.  More generally, we define a measure by saying that, if $\Gamma_H$ is a triangle-free graph (with vertex set $\{1, 2, \ldots, n\}$), then the ``frequency of $\Gamma_H$ in $N$-vertex graphs'' should be the proportion of triangle-free graphs on the vertex set $\{1, 2, \ldots, N\}$ which induce $\Gamma_H$ on $\{1, 2, \ldots, n\}$.  Now suppose that this frequency tends to a limit as $N \to \infty$.  Then, for any vertices $v_1, \ldots, v_n$ in the infinite set, we let the event that these vertices induce $\Gamma_H$ be this limiting frequency; this defines the required measure.

In this case, a theorem of Erd\H{o}s, Kleitman and Rothschild 
\index{Erd\H{o}s, P.}%
\index{Kleitman, D. J.}%
\index{Rothschild, B. L.}%
which says~\cite{erdos} that `almost all finite triangle-free graphs are bipartite', in the sense that the proportion of triangle-free graphs which are bipartite tends to 1 as the number of vertices tends to $\infty$; see $\S 6$ of Chapter~\ref{chap2}.  This implies that the probability that any odd number $n$ of vertices in the random graph
\index{graph ! random}%
 induce a cycle is zero; so the random infinite triangle-free graph is almost surely (that is, with probability 1) bipartite.
 
Consider next the age
\index{age}%
 of a graph $\Gamma$, that is the class of finite graphs embeddable in $\Gamma$ as induced subgraphs.  Whilst different infinite graphs may have the same age, there may be one of specific interest.  So for example, the class of all finite graphs is an age, and the most interesting infinite graph of which it is the age is the random graph.  For the finite triangle-free graphs, Henson's infinite graph
\index{graph ! Henson}%
 has this age.  So our earlier remarks show that Henson's graph is not ``the random triangle-free graph''.
 \index{graph ! triangle-free ! random}%
 (However, it is the typical triangle-free graph
\index{graph ! triangle-free}%
  in the sense of Baire category.)
\index{Baire category theorem}%

\head{Open Question}  Let $\mathcal{A}(\Gamma)$ be a fixed age of finite graphs, and let $\Gamma_H$ be a fixed graph on the vertex set $\{1, \ldots, n\}$.  For $N \geq n$, let $p_N (\Gamma_H)$ be the proportion of graphs in $\mathcal{A}(\Gamma)$ (on the vertex set $\{1, \ldots,N\}$) which induce $\Gamma_H$ on the first $n$ vertices.  Is it true that $p_N (\Gamma_H)$ tends to a limit $p(\Gamma_H)$ as $N \to \infty$?

If this is true, then we define a probability measure on graphs on a given countable vertex set by the rule that the event that a given $n$-tuple of vertices induces $\Gamma_H$ should be $p(\Gamma_H)$.  We can then ask a second question:

\head{Open Question}  What does the random graph in this model look like?

(Clearly, with probability 1, its age is contained in $\mathcal{A}(\Gamma)$; the containment may be proper, as in the case of the random triangle-free graph.)

As we mentioned the recent developments in infinitary limit objects for sequences of combinatorial objects are related to probabilistic results about exchangeable measures.  The simplest case is that of measures on $0$--$1$ sequences,
\index{zero-one sequence}%
 which we can regard as functions from the natural numbers to $\{0,1\}$.  Such a measure is called ``exchangeable''
\index{measure ! exchangeable}%
 if it is unaffected by arbitrary permutations of the natural numbers.

Exchangeability means that the probability measure is invariant under permutations, not any particular random object.  The simplest example of an exchangeable probability measure is given by tossing a biased coin countably often.  A slightly more complicated measure is defined as follows: take two biased coins with different probabilities of heads; choose one of the two according to some probability measure; then generate the random sequence using this coin.

This concept was introduced by de Finetti
\index{de Finetti, B.}%
 who worked on the foundations of probability in the 1930s.  He proved a representation theorem for exchangeable measures.  Any such measure is defined by a function $f$ from $[0,1] \times [0,1] \to [0,1]$ by the following rule.  Choose $u$ uniformly from $[0,1]$; for each natural number $n$, choose $u_n$ uniformly from $[0,1]$, with all choices independently; then let the $n$th term of the random binary sequence be chosen to be 1 with probability $f(u, u_n)$, independently for all $n$. It is an exercise for the reader to find the function $f$ representing each of the two examples above.

In the 1980s, this result was extended by Hoover and further by Aldous.
\index{Aldous, D.}%
\index{Hoover, D.}%
 We now define a probability measure on graphs on the vertex set $\mathbb{N}$ to be exchangeable if it is invariant under permutations of the natural numbers.  The measure giving the ``random graph with age contained in $\mathcal{A}(\Gamma)$'' described earlier is clearly exchangeable.
 
Aldous and Hoover gave a representation theorem for exchangeable random graph measures similar to de Finetti's
\index{de Finetti, B.}%
 except that the function replacing $f$ has three variables rather than two.

In an important special case, ergodicity can be used to show that the function is independent of the first parameter, so is really a symmetric function from $[0,1]^2$ to $[0,1]$.  As in de Finetti's case, the values of the function give the probabilities of seeing a given induced subgraph on a given set of vertices.  

In fact, what we have here is precisely a a limit of a sequence of finite graphs called a \emph{graphon}.
\index{graphon}%
 In brief, from the graph $\Gamma_n$ in a sequence $(\Gamma_n)$ of finite graphs, define a probability distribution $\mu_n$ on graphs with at most $k$ vertices, for any given $k$, by choosing $k$ random vertices of $\Gamma_n$ and taking the induced subgraph.  (The words ``at most'' allow for the event that the same vertex is picked more than once.)  The sequence $(\Gamma_n)$ converges if, for every fixed $k$, the probability measures $\mu_n$ converge to a limit $\mu$.  (This can also be expressed, using inclusion-exclusion, in terms of graph homomorphisms; the condition for convergence is that, for any fixed finite graph $\Gamma_H$, the probability $t(\Gamma_H,\Gamma_n)$ that a random map from $\Gamma_H$ to $\Gamma_n$ converges to a limit $t(\Gamma_H, \Gamma)$.)

A graphon is the ``limit'' of a convergent sequence of finite graphs.  But what is one?  It can be described by a function $W$ from $[0,1]^2$ to $[0,1]$, which is symmetric and measureable. The measure $\mu$ is defined from the graphon as follows: pick any real numbers $x_1, \ldots, x_k$ in the unit interval; restricted to these, $W$ defines a symmetric matrix of order $k$. Now join $i$ to $j$ with probability $W(x_i,x_j)$.

\head{Open Question}  What is the connection with the above question about random graphs with given age.

However, something is still missing.  The graphon corresponding to the random graph is the function $W$ which takes the constant value $\frac{1}{2}$.

\head{Open Question}  Can the remarkable properties of the random graph be extracted from the constant function  $\frac{1}{2}$? If so, how?

Also, changing $\frac{1}{2}$ to any other number between zero and one presumably gives a different graphon, but gives the same countable random graph!

\medskip

An automorphism of a graphon is a measure-preserving permutation of $[0,1]$ that preserves the value of the function almost everywhere.  After some standardization the automorphism group
\index{group ! automorphism}%
 of any graphon is compact.
\index{group ! compact}%
  The orbits of the automorphism group can be characterized by certain generalized degrees.
The tools required to prove these results include topologies on the set of points of a graphon.  Some properties of the graphon automorphism groups carry over to the limit objects of graph sequences; for example, the limit of graphs with vertex-transitive automorphism groups is a graphon with a point-transitive automorphism group.  In this joint work of Lov\'asz
\index{Lov\'asz, L.}%
 and Szegedy,
\index{Szegedy, B.}%
 the transition to the limit is not completely understood.

\bigskip

We end by noting an application of Petrov and Vershik's
\index{Petrov, F. V.}%
\index{Vershik, A. M.}%
 construction~\cite{petrovv} of probability measures on the set of graphs on a countable vertex set which is \emph{exchangable},
\index{graph ! random ! exchangeable}%
 that is, invariant under all permutations in the symmetric group, (see below) and concentrated on Henson's universal triangle-free graph.
\index{graph ! Henson}%
 Their approach of building topological or measure-theoretic triangle-free graphs on the real numbers means that triangles are not totally excluded, as long as the chance of getting one is zero; similarly, the universality axioms are not reqired to hold universally, only with high probability.)  Then they obtain a countable graph by taking countably many independent samples from a probability distribution on the real numbers.  The actual distribution used is not crucial; a Gaussian distribution suffices.  The resulting graph is isomorphic to Henson's graph with probability 1.

Using Aldous' Theorem on exchangeable measures Petrov and Vershik completely characterize all measures with the required property.  This probably discounts a finitary approach building a graph one vertex at a time.

Nate Ackerman, Cameron Freer and Rehana Patel 
\index{Ackerman, N.}%
\index{Patel, R.}%
\index{Freer, C.}%
 proved that there is an exchangeable measure concentrated on the structures isomorphic to a countably infinite structure $\mathcal{M}$ if and only if the stabiliser of a finite number of points in $\Aut(\mathcal{M})$ fixes no additional points, or more precisely,

\begin{theorem}[N. Ackerman, C. Freer, and R. Patel]
 Let $L$ be a countable relational language, $\mathcal{M}$ a countable $L$-structure with underlying set $\mathbb{N}$, and $S_L$ the measurable space of all $L$-structures with underlying set $\mathbb{N}$, equipped 
with the natural Borel $\sigma$-algebra (generated by sets of the form $\{\mathcal{M} \in S_L : \mathcal{M} \models R(n_1, \ldots, n_i) \}$, where $n_1, \ldots, n_i \in \mathbb{N}$ and $R \in L$.  The following are equivalent: 

(1) There is a probability measure on $S_L$, invariant under the action of $\Sym(\mathbb{N})$, 
that is concentrated on (the isomorphism class of) $\mathcal{M}$.

(2) ``Group-theoretic'' definable closure
\index{definable ! closure}%
 in $\mathcal{M}$ is trivial, that is for every finite 
tuple $\textbf{a} \in \mathcal{M}$, we have $dcl_\mathcal{M} (\textbf{a}) = \textbf{a}$, where the \emph{definable closure} $dcl_\mathcal{M} (\textbf{a})$ is the collection of $b \in \mathcal{M}$ that are fixed by all automorphisms of $\mathcal{M}$ fixing $\textbf{a}$ pointwise.
\end{theorem}

Part of the proof follows Petrov and Vershik in constructing a topological version of the structure in the real numbers, and then sampling countably many independent reals from a fixed distribution and taking the induced substructure.

For example, it holds for the random graph, or for the rational numbers as ordered set, but not for an infinite path or tree (where the stabiliser of two points fixes every point in the interval joining them).  For homogeneous structures in a finite relational language, having trivial $dcl$ is equivalent to having strong amalgamation.
\index{amalgamation property ! strong}%

\begin{corollary} 
Suppose further that $L$ is finite and $\mathcal{M}$ is homogeneous with age $\mathcal{A(\mathcal{M})}$.  Then the following are equivalent:\\
(1) There is an invariant measure on $S_L$ that is concentrated on $\mathcal{M}$.\\ 
(2) $\mathcal{A(\mathcal{M})}$ has the strong amalgamation property.
\end{corollary}
\item \emph{B-groups}.
\index{group ! B-group}%
A group $G$ is a B-group if every primitive permutation group containing the regular representation of $G$ is doubly transitive. There are many finite $B$-groups (including cyclic groups of composite order), but no countably infinite B-groups are known.

Does there exist a countably infinite B-group?

Let $G$ be the infinite dicyclic group,
\index{group ! dicyclic}%
 generated by two elements $a$ and $b$, where $b$ has order 4 and inverts $a$.  Is $G$ a B-group? (Note that the countable random graph is not a Cayley graph for $G$; so the standard method for showing a group is not a B-group fails).

\item  \emph{Group Actions on Random Hypergraphs}.  A \emph{tournament}
\index{group ! action}%
\index{tournament}%
 is a binary relation $b$ with trichotomy, that is exactly one of $b(x, y),\ x=y,\ b(y, x)$ holds for each pair $x, y$ of points.  There is no obvious way of defining a countable homogeneous
     random triality tournament, the tournament equivalent to
     $\mathfrak{R^{t}}$, because tournaments forbid symmetric joining
     relations.  However one way of defining a \emph{hypergraph}
     \index{hypergraph}%
 equivalent to $\mathfrak{R}_{m,\omega}$ is to consider the $m$-coloured $k$-complete
     $k$-uniform hypergraph, where each hyperedge comprises a
     $k$-subset of the vertex set in which every vertex is joined to
     every other vertex.  Investigate the properties of this object.  The I-property
\index{I-property}%
 for the countable universal homogeneous $k$-uniform hypergraph asserts that, given a finite set $A$ of vertices, and a set $B$ of $(k - 1)$-subsets of $A$, there is a vertex $z$ such that, for a $(k - 1)$-subset $K \subset A$, $K \cup \{ z \}$ is a hyperedge if and only if $K \in B$.
 
Produce a classification theorem that is the $m$-coloured equivalent of
     that of S. Thomas~\cite{thomas1},
     \index{Thomas, S.}%
 which classifies the reducts of
     the countable universal homogeneous random
     $k$-hypergraphs for all $k \ge 1$, which apart from the relevant automorphism and
     symmetric groups all turn out to be suitably-defined
     generalizations of switching groups.  Furthermore, the unique
     countable universal homogeneous $k$-hypergraph $\mathfrak{Hyp}$ has an
     infinite automorphism group
\index{group ! automorphism}%
      $G = \Aut(\mathfrak{Hyp})$ that is not
     $k$-transitive because it preserves $k$-hyperedges of
     $\mathfrak{Hyp}$~\cite{bollobas}.  However any two $k$-hypergraphs of size $k-1$
     are isomorphic because they have no hyperedges, so $G$ is
     $k-1$-transitive by the homogeneity of $\mathfrak{Hyp}$.  Find all homogeneous $k$-uniform hypergraphs for fixed $k$ and develop some theory of multicoloured complete $k$-uniform hypergraphs.

\item \emph{Switchings $\&$ Reducts of Random Homogeneous Tournaments}.
\index{tournament ! random}%
Classify the reducts of random homogeneous tournaments.  We give an introduction to this topic.

Two tournaments $T_1$ and $T_2$ on the same vertex set $X$ are said to be \emph{switching equivalent} if $X$ has a subset $Y$ such that $T_2$ arises from $T_1$ by switching all arcs between $Y$ and its complement $X \backslash Y$.  The abstract finite groups which are full automorphism groups
\index{group ! automorphism}%
 of switching classes of tournaments were characterized in~\cite{babaicam}.  In doing so, they describe two objects ``equivalent'' to switching classes of 
tournaments, one of which is the following.  A tournament can be regarded as an antisymmetric function $f$ from ordered pairs of distinct vertices to $\{\pm 1 \}$ (with $f (x, y ) = +1$ if and only if there is an arc from $x$ to $y$).   Switching with respect to $\{x\}$ corresponds to changing the sign of $f$ whenever $x$ is one of the arguments; and switching with respect to an arbitrary subset is performed by switching arcs between vertices in the subset and those in the complement.  Given a tournament $f$, define an \emph{oriented two-graph}
\index{two-graph ! oriented} 
to be a function $g$ on ordered triples of distinct elements satisfying
\[ g(x, y , z) = f (x, y )f (y , z)f (z, x).\] 
Then $g$ is alternating (for interchanging two arguments changes the sign) and satisfies the ``cocycle'' condition 
\[g(x, y , z)g(y , x, w)g(z, y , w)g(x, z, w) = +1.\]
Conversely, any oriented two-graph arises from a tournament in this way, and two tournaments yield the 
same oriented two-graph if and only if they are equivalent under switching.  Thus there is a natural bijection between switching classes of tournaments and oriented two-graphs; corresponding objects have the same automorphism group.
\index{group ! automorphism}%
  It transpires from Fra\"{\i}ss\'e's Theorem
\index{Fra\"{\i}ss\'e's Theorem}%
that both the class of switching classes of finite tournaments, and the class of ``local orders'' (that is, tournaments switching-equivalent to linear orders),
\index{linear order}%
 give rise to countably infinite structures with interesting automorphism groups.

Lachlan
\index{Lachlan, A. H.}%
determined~\cite{lachlan} that there are three countable homogeneous tournaments, (see also Cherlin~\cite{cherlina}):
\index{Cherlin, G. L.}%
the transitive tournament $\mathbb{Q}$, the homogeneous local order, and the homogeneous tournament $T$ containing all finite tournaments. The homogeneous oriented two-graph corresponds to the switching class of $T$.

Let us recall some homological algebraic theory~\cite{cartan1} applied to graph theory~\cite{cameron}~\cite{cam8a}.  Let $\Gamma$ be a graph and $\Lambda$ be an abelian group. The \emph{$0$-cochain} on $\Gamma$ (with coefficients in $\Lambda$) is a function from vertices of $\Gamma$ to $\Lambda$. A \emph{$1$-cochain} is a $\Lambda$-function from ordered edges, where the sign depends on edge-direction.  The \emph{$1$-coboundary} $\partial e$ of a $0$-chain $e$ is given by $\partial e (x, y) = e(y) - e(x)$, and the \emph{$2$-coboundary} of a $1$-chain $f$ is the function on oriented induced cycles whose value on a cycle $C$ is the sum of the values of $f$ on the ordered edges of $C$.  Two $1$-cochains have the same $2$-coboundary if and only if they differ by the $1$-coboundary of a $0$-chain. 

In the case that $\Gamma$ is complete and $\Lambda = C_2 = \{0, 1\}$, a $0$-cochain (respectively a $1$-cochain) is the characteristic function of a vertex (respectively edge) set of a graph.  If $e$ is a $1$-cochain, then $\partial e$ is a complete bipartite graph, and adding it to $f$ corresponds to switching the graph represented by $f$ with respect to the subset represented by $e$.  The $2$-coboundary of $f$ is the two-graph, the set of triples containing an odd number of graph edges.

Given a $1$-chain $f$ on $\Gamma$ with vertex set $V$, the graph $\Gamma_{f}$ has vertex set $V \times \Lambda$, where vertices $(x, \alpha)$ and $(y, \beta)$ are adjacent if $f(x, y) = \beta - \alpha$. The graph $\Gamma_f$ is a cover of $\Gamma$, with covering projection $p$ given by $p(x, \alpha) = x$.  If $f' = f + \partial e$ for some $0$-cochain $e$, then the map $(x, \alpha) \mapsto (x, \alpha + e(x))$ is an isomorphism from $\Gamma_f$ to $\Gamma_{f'}$.  So the covering graph depends only on the coboundary of $f$.  The following result is from~\cite{cameron}~\cite{cam8a}.

\begin{theorem}[P. J. Cameron]
\index{Cameron, P. J.}%
If $\Gamma$ is a countable graph whose age
\index{age}%
 has the strong amalgamation property,
\index{amalgamation property ! strong}%
 and $\Lambda$ a finite or countable abelian group, then there exists a universal homogeneous $2$-coboundary on $\Gamma$ (with coefficients in $\Lambda$).
\end{theorem}

It follows that the universal homogeneous two-graph
\index{two-graph ! universal homogeneous}%
 has an incarnation as the universal homogeneous 2-coboundary over $\mathbb{Z} / 2$ on the complete graph, meaning that it can be constructed as the coboundary of $\mathfrak{R}$, where a graph is regarded as the 1-cochain on the complete graph, edges and non-edges corresponding to values 1 and 0, respectively~\cite{cameron}~\cite{cam8a}.  Is it true that the universal co-chain on the complete graph over $\mathbb{Z} / 3$, which includes valuations over ordered triples, that is two edges joined at a vertex corresponding to the value 2, is a reduct of the random oriented graph
\index{graph ! oriented ! random}%
  or the universal $\mathbb{Z} / 3$ cochain on a complete graph?  

Assigning the graph edges as follows
\begin{figure}[!h]
$$\xymatrix{
& {\bullet}_{x} \ar@{->}[r]^{0} & {\bullet}_{y} && {\bullet}_{x} \ar@{->}[r]^{+1} & {\bullet}_{y} && {\bullet}_{x} \ar@{->}[r]^{-1} & {\bullet}_{y} 
}$$ 
\end{figure}
leads to the have the possibilities in Figure~\ref{3sw}.

There will be three different types of switching, between pairs of $0, +1$ and $-1$ configurations; for example,

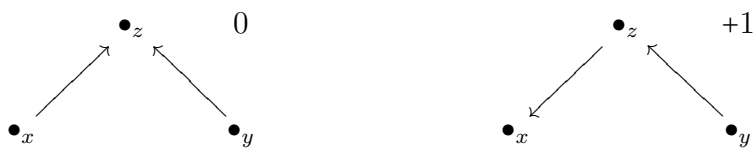
\begin{figure}[!h]
$$\xymatrix{
& {\bullet}_{z} & {0} &&&& {\bullet}_{z}  & {+1}\\
{\bullet}_{x} \ar@{->}[ur] && {\bullet}_{y} \ar@{->}[ul] &&& {\bullet}_{x} \ar@{<-}[ur] && {\bullet}_{y} \ar@{->}[ul] 
}$$\caption{Example of a $0$ - $+1$ Switching}
\end{figure} 

is a switch with respect $x$.

Is the automorphism group of this structure strictly larger than that of the random oriented graph?  This is the same as the universal homogeneous coboundary, so by definition it is a reduct, but is it the trivial reduct or are there non-trivial automorphisms?

\begin{figure}[!h]
$$\xymatrix{
& {\bullet}_{z} & {0} &&&& {\bullet}_{z}  & {+1}\\
{\bullet}_{x} && {\bullet}_{y} &&& {\bullet}_{x} \ar@{->}[rr] && {\bullet}_{y}\\
& {\bullet}_{z} & {0} &&&& {\bullet}_{z}  & {-1}\\
{\bullet}_{x} \ar@{->}[ur] && {\bullet}_{y} \ar@{->}[ul] &&& {\bullet}_{x} \ar@{->}[ur] && {\bullet}_{y} \ar@{<-}[ul] \\
& {\bullet}_{z} & {0} &&&& {\bullet}_{z}  & {0}\\
{\bullet}_{x} \ar@{<-}[ur] && {\bullet}_{y} \ar@{<-}[ul] &&& {\bullet}_{x} \ar@{->}[ur] && {\bullet}_{y} \ar@{<-}[ul] \ar@{<-}[ll] \\
& {\bullet}_{z} & {+1} &&&& {\bullet}_{z}  & {-1}\\
{\bullet}_{x} \ar@{->}[rr]  \ar@{->}[ur] && {\bullet}_{y} \ar@{<-}[ul] &&& {\bullet}_{x} \ar@{<-}[ur] && {\bullet}_{y} \ar@{->}[ul] \ar@{->}[ll]
}$$\caption{$\mathbb{Z} / 3$ Switching Configurations}
\label{3sw}
\end{figure}
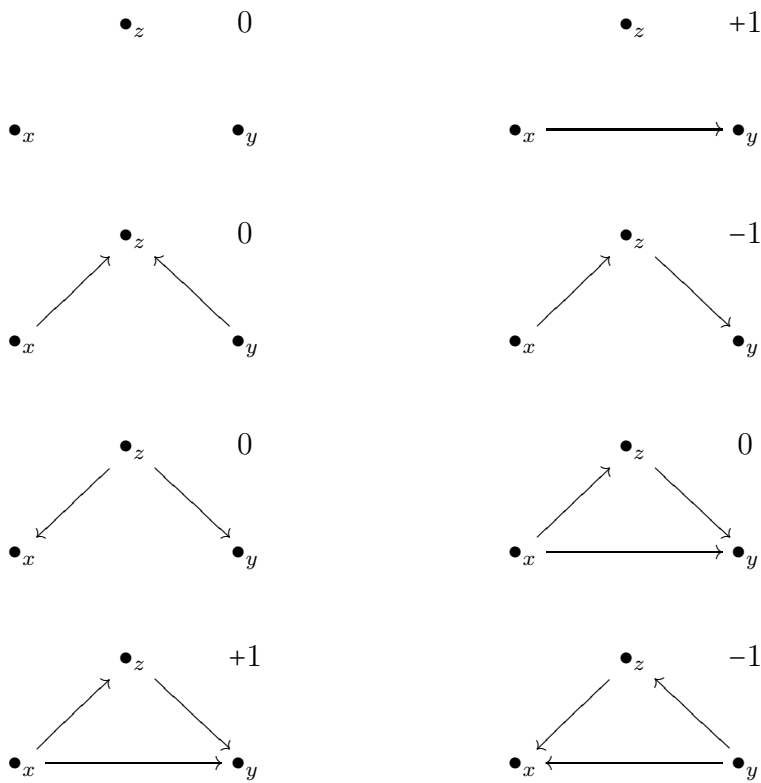 

\item  \emph{Infinite Association Schemes}.
\index{association scheme ! infinite}%
An \emph{association scheme}~\cite{bailey}
\index{association scheme}%
with $m$ classes on a finite set $\Omega$ is a partition of $\Omega \times \Omega$ into sets $\mathcal{C}_0, \mathcal{C}_1, \ldots, \mathcal{C}_m$ of binary relations (called \emph{associate classes})
\index{association scheme ! associate classes}%
such that 

$(i)$ $\mathcal{C}_0$ is the diagonal $\{(x, x) : x \in \Omega\}$ of $\Omega^2$

$(ii)$  $\mathcal{C}_i$ is symmetric for $i = 1, \ldots, m$;

$(iii)$  for all $i, j, k \in \{0, \ldots, m\}$ there is an integer $p^k_{ij}$ such that for all $(x, y) \in \mathcal{C}_k$, $p^k_{ij} = |\{z \in \Omega : (x, z) \in \mathcal{C}_i\ \text{and}\ (z, y) \in \mathcal{C}_j\}|$.  The number $p^k_{ij}$ is independent of $x$ and $y$.

If $(x, y) \in \mathcal{C}_i$ then $x$ and $y$ are called $i$-th associates.  If $G$ is transitive then the partition of $\Omega \times \Omega$ into $G$-orbits satisfies (i) and (iii); condition (ii) is satisfied if $G$ is generously transitive.
\index{group ! permutation ! generously transitive}%

There is a second definition of association scheme, given in terms of edge-coloured graphs.  An association scheme with $m$ associate classes on a finite set $\Omega$ is an edge-colouring of the complete undirected graph with vertex-set $\Omega$ by $m$ colours such that

$(iii)'$  for all $i, j, k \in \{0, \ldots, m\}$ there is an integer $p^k_{ij}$ such that whenever $\{x, y\}$ is an edge of colour
$k$ then $p^k_{ij} = |\{z \in \Omega :\ \{x, z\} \text{ has colour } i \text{ and } \{z, y\} \text{ has colour } j\}|$.

$(iv)'$  every colour is used at least once;

$(v)'$  there are integers $a_i$ for $i \in \{0, \ldots, m\}$ such that each vertex is contained in exactly $a_i$ edges of colour $i$.

The analogue of $(i)$ above is not required because every edge has two distinct vertices and the analogue of $(ii)$ is not required because we have specified an undirected graph.

\smallskip

Generously transitive groups are characterised
by the property that the partition of the edge set of the complete
graph on $X$ into $G$-orbits is an association scheme whose
\emph{Bose-Mesner algebra}
\index{Bose-Mesner algebra}%
 is the centralizer algebra
\index{centralizer algebra}%
 of $G$, the algebra of all matrices commuting with permutation matrices corresponding to elements of $G$.  If $rank(G) =
m+1$ then the association scheme has $m$ classes and $m+1$ symmetric binary relations on $X$, the first of which is the relation of equality.  

The matrix representation of a binary representation $\mathcal{C}$ on $\Omega$ has rows and columns indexed by $\Omega$ and entries $(x, y) =1$ if $(x, y) \in \mathcal{C}$ and $(x, y) = 0$ otherwise.  The representation of the association scheme
\index{association scheme ! representation of}%
 in terms of matrices $A_0, \ldots, A_m$ then translates to 

(BM1)  $A_0 + \ldots + A_m = J$, where $J$ is the all-1 matrix,

(BM2)  $A_0 = I$, where $I$ is the identity matrix,

(BM3)  $A_i^T = A_i$, for all $i$,

(BM4)  $\exists$ `intersection numbers' $p^k_{ij}$ such that \[A_i A_j = \sum_{k = 0}^m p^k_{ij} A_k.\]

Axiom (BM2) and (BM4) imply that $\{A_i\}$ spans an algebra (over $\mathbb{C}$).  This is the Bose-Mesner algebra.  Axiom (BM3) implies that the algebra is semisimple.  The matrices are symmetric and commuting so the vector space $\mathbb{R}^{\Omega}$ has an orthogonal decomposition into common eigenspaces of the matrices.  The identity $p^k_{ij} = p^k_{ji}$ implies that the $\{A_i\}$ commute pairwise.  That the $m$-dimensional span (over $\mathbb{R}$) of the matrices is closed under multiplication means that the algebra is commutative.

The theory of \emph{coherent configurations}
\index{coherent configuration (cellular ring)}%
 (c.c.) was develped by Higman
\index{Higman, D. G.}%
as operands of non-2-transitive finite permutation groups.  Boris Weisfeiler
\index{Weisfeiler, B.}%
 introduced the equivalent notion of cellular algebras in the USSR~\cite{Weisfeiler}, where they arise in the ``partition refinement'' approach to the graph isomorphism problem.  An association scheme is a symmetric coherent configuration.  For any permutation group $G$ acting on $\Omega$, the partition $\mathcal{K}(G)$ of $\Omega^2$ into orbits of $G$ is a c.c.  Note that $G$ is transitive if and only if $\mathcal{K}(G)$ is homogeneous (so that $A_0 = I$); also $G$ is 2-transitive if and only if $\mathcal{K}(G)$ attains its smallest dimension of 2 if and only if the c.c. is trivial.  A 2-transitive $G$
\index{group ! permutation ! $2$-transitive}%
 has two orbits on $\Omega^2$ with characteristc functions $I$ and $J - I$, so is \emph{AS-free},
 \index{group ! permutation ! AS-free}%
 that is the transitive group $G$ preserves no non-trivial association scheme on operand $\Omega$.
 
The set of partitions of $\Omega^2$ forms a lattice (with ``smaller'' = finer), called the \emph{partition lattice}.
\index{lattice ! partition}%
The set of coherent configurations is a meet-semilattice of the partition lattice.  Hence, given any family $F$ of subsets of  $\Omega^2$, there is a unique finest coherent configuration containing them, which we call the coherent configuration generated by $F$.  The partition into singletons forms the ``trivial'' coherent configuration.

More on links between association schemes and permutation groups can be found in~\cite{alejan} and~\cite{golf}.  

The first open question that we pose is find a definition of \emph{infinite} association schemes and develop a theory.  Clearly none of the  $p^k_{ij}$ can be infinite otherwise the theory falls down.
This may well require a ``local-finiteness'' condition, such as requiring that only a finite number of the $p^k_{ij}$ are non-zero.  If $G$ has only finitely many orbits on $\Omega^2$ then certainly only finitely many $p^k_{ij}$ are nonzero; this may be the only tractable case.  But there is thus far no theory for this situation.  

The condition that all $p^k_{ij}$ are finite is equivalent to requiring that the point stabilizer $G_{x}$ of the group $G$ acting on the underlying set $\Omega$ has all only finite orbits, which is the condition for the group to be locally compact.
\index{group ! locally compact}%
This equivalence should yield an interesting theory on further study.  There should also be a measure-theoretic
\index{measure}%
version of this question for which the measure satisfies $\mu\{z \in \Omega : (x, z) \in \mathcal{C}_i\ \text{and}\ (z, y) \in \mathcal{C}_j\} < \infty$ $\forall i, j, k$.  It may be possible to tie this up with the work of M\"{o}ller
\index{moller@M\"{o}ller, R.}%
 and co-workers on totally disconnected locally compact groups
\index{group ! locally compact}%
\index{group ! totally disconnected}%
 via graphs and permutations; see~\cite{malnic} and~\cite{kron} and references therein for details.

A comprehensive reference for the theory of association schemes is the book by R. A. Bailey~\cite{bailey},
\index{Bailey, R. A.}%
the notes by Chris Godsil
\index{Godsil, C. D.}%
or the paper~\cite{cam21}.

\item  \emph{Association Schemes and Galois Theory}.
\index{association scheme}%
\index{Galois ! theory}%
Building on earlier work, Chris Godsil
\index{Godsil, C. D.}%
used Galois theory to establish a correspondence between certain subfields of a splitting field
\index{splitting field}%
 and subschemes of association schemes, where the splitting field of the scheme is the extension of the rationals generated by the eigenvalues of the scheme.  Is there a further connection with switching groups regarded as Galois groups?  In particular, concretely realize the following correspondences:
$$
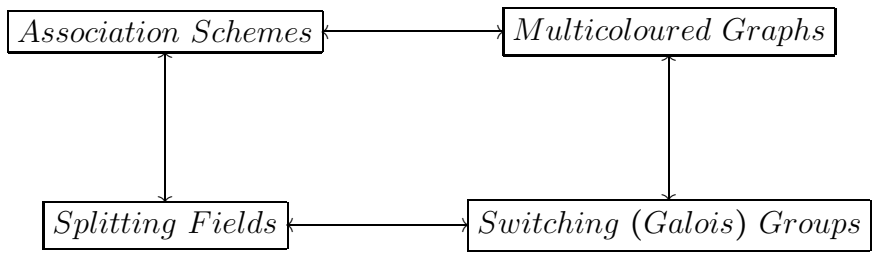
\xymatrix{ 
& *+[F]{Association\ Schemes}  \ar[rr] \ar[dd] 
&& *+[F]{Multicoloured\ Graphs} \ar[ll] \ar[dd]\\
&&\\
& *+[F]{Splitting\ Fields} \ar[rr] \ar[uu] &&
*+[F]{Switching\ (Galois)\ Groups} \ar[ll] \ar[uu]\\
}$$

\item  \emph{Invariant Relations of Permutation Groups}.
\index{group ! permutation ! invariant relations}%
  Investigate possible links between random graph automorphism groups,
\index{group ! automorphism}%
   finitary switching groups or other reducts and the work of Liebeck and Praeger~\cite{lipr}~\cite{prae}
\index{Praeger, C. E.}%
\index{Liebeck, M. W.}%
 on \emph{relation algebras} of
  binary relations.  They study the Galois connection
\index{Galois ! connection}%
 between the
  subalgebra lattice of the relation algebra
  \index{relation algebra}%
   of all binary relations
  on a finite set $\Omega$ and the subgroup lattice of subgroups of
  $\Sym(\Omega)$.  The Galois closed subgroups are the permutation
  groups that are \emph{2-closed} in the sense of Wielandt~\cite{wielandt1}, 
  \index{Wielandt, H.}%
  \index{group ! permutation ! $2$-closed}%
  that is the largest subgroup of $\Sym(\Omega)$ which preserves the orbits of $G$ in its coordinate-wise action on the set $\Omega \times \Omega$.  An example of a 2-closed permutation group
\index{group ! permutation ! $2$-closed}%
   is the full automorphism group
\index{group ! automorphism}%
   of a graph, because any vertex permutation which preserves the orbits on ordered pairs will preserve adjacency.
  (Groups of switching permutations on more than two colours are highly
  transitive have no group-invariant relational structure).
\index{relational structure}%
  
Alternatively investigate the reducts from the viewpoint of
  \emph{Krasner algebras}~\cite{farad},
\index{Krasner algebra}%
 these being the set of invariant relations on a set $\Omega$ which is a union of $k$-closed relation algebras defined on $\Omega^k$, $\forall k \ge 2$.  These can also be put in Galois correspondence
\index{Galois ! correspondence}%
 with the set of permutations on $\Omega$.

Given a set $S$, if the set $O$ of orbits of a permutation group $G$ acts diagonally and componentwise on $S^2$ to form an association scheme $(S, O)$, then the scheme $(S, O)$ is called a \emph{permutation group scheme}.
\index{association scheme ! permutation group}%
  This arises if and only if $(S, O)$ has a transitive automorphism group.
\index{group ! automorphism}%
    J. D. H. Smith
\index{Smith, J. D. H.}%
has given a higher-dimensional analogue of an association scheme called a \emph{superscheme}
\index{association scheme ! superscheme}%
which arises as the set of componentwise orbits on direct powers of $S$, of a transitive multiplicity-free permutation group $G$ acting on $S$.  The group and its action are recovered from the Bose-Mesner superalgebra
\index{Bose-Mesner algebra}%
of the superscheme, which is a `quantised' version of a Krasner (relation) algebra of the second kind.

Another direction of research which is relevant here is that concerning \emph{clones}.
A \emph{clone}
\index{clone}%
on a set $A$ is a collection of finitary operations on $A$ that
contains the projection operations and is closed under composition.  Clones are generalizations of monoids, considered as a set of selfmaps of the domain set which is closed under composition and contains the identity mapping.  For any set $F$ of operations on $A$, the clone $[F]$
generated by $F$ is the collection of all \emph{term operations} of
the algebra $\mathbb{A} = \langle A; F \rangle$, which is also called
the \emph{clone of} $\mathbb{A}$.  There is a Galois connection
\index{Galois ! connection}%
between the lattice of clones on a set $A$ and $\Sym(A)$ that is
determined by the relation that a permutation conjugates a clone onto
itself.  The Galois closed sets on the clone side are the lattices of
all clones that are closed under conjugation by all members of some
permutation group.  More on this can be found in the survey paper of
\'A. Szendrei~\cite{szendrei}.
\index{Szendrei, \'A.}%
Bodirsky and Pinsker
\index{Bodirsky, M.}%
\index{Pinsker, M.}%
 have initiated a study of clones~\cite{bodirsky} pertinent to random graph theory by studying the lattice of closed transformation monoids, and more generally clones closed under composition, that contain $\Aut(\mathfrak{R})$.


\item \emph{Links Between Switching and Named Groups}.  The switching group $C_{2} \Wr \Sym(n)$ occurs in many guises in such areas as matroid theory
\index{matroid theory}%
and Lie algebras 
\index{Lie algebra}%
as we now briefly recall.

We briefly indicate when the (slightly redefined from its appearance in Chapter 2) \emph{extended switching groups}
\index{group ! switching ! extended}%
defined by $S^{e}_{2,n}:=
S_{2,n+1} \sd \Sym(n)$ and acting on the set of $2$-coloured $n$-vertex random graphs may be isomorphic to hyperoctahedral groups which we denote $\Hyp(n)$.~\label{hyperoctahedral}
\index{group ! hyperoctahedral}%

The hyperoctahedral groups are special cases of groups that arise in
graph theory due to the fact that if $\Gamma_c$ is a connected graph,
\index{graph ! connected}%
then $\Aut(k \Gamma_c) = \Aut(\Gamma_c) \Wr \Sym(k)$, $k > 1$~\cite{hag}.  They also arise in studies of matroid automorphisms via signed graphs~\cite{fern}.
\index{graph ! signed}%
The $n$-dimensional hypercube can be regarded as a graph whose vertex set consists of $0$'s and $1$'s of length $n$.  The group $\Hyp(n)$ is the automorphism group
\index{group ! automorphism}%
 of this hypercube and can be represented by the group of signed permutations, that is vertex permutations of the hypercube which preserve adjacency.

The origin of $\Hyp(n)$ is as the Weyl group
\index{group ! Weyl}%
 $B_n := 2^n \sd \Sym(n)$ associated with one of four families of complex simple Lie
algebras also denoted $B_n$.  The hyperoctahedral group is the
automorphism group
\index{group ! automorphism}%
 of the root lattice of the Lie group $B_n$
\index{group ! Lie}%
with roots $\pm e_i$ and $\pm e_i \pm e_j$ for $1 \le i < j \le n$.  This
group has the presentation $\langle s_0, s_1, \ldots, s_{n-1} |
s_i^{2} = 1\ \forall\ i,\ s_i s_j = s_j s_i\ \text{if}\ | i - j | > 1,
s_i s_{i+1} s_i = s_{i+1} s_i s_{i+1}\ \text{if}\ i \ge 1,\
\text{and}\ s_1 s_0 s_1 s_0\\ = s_0 s_1 s_0 s_1 \rangle$.

Given the $B_n$ Coxeter-Dynkin diagram $(n \ge 2)$
$$\xymatrix{
\quad\quad\quad\ {\circ}^{e_1}_0 \ar@{=}[r] & {\circ}^{e_2 - e_1}_1 \ar@{-}[r] & {\circ}^{e_3 - e_2}_2 \ar@{-}[r] & {} {\ldots}  {\circ}^{e_{n-1} - e_{n-2}}_{n-2} \ar@{-}[r] & {\circ}^{e_{n} - e_{n-1}}_{n-1}
}$$
then $s_i$ is the reflection in the hyperplane perpendicular to the vector labelling the $i$th node.  Since the angle between $e_1$ and $e_2 - e_1$ is $135^{\circ}$, the angle between $e_i - e_{i-1}$ and $e_{i+1} - e_{i}$ is $120^{\circ}$, and all other angles are $90^{\circ}$, the orders of $s_0s_1$, $s_{i-1} s_i$ $(i \ge 2)$ and $s_i s_j$ $(|i - j| > 1)$ are 4, 3, and 2 respectively, giving the stated relations.

One way to view $B_n$~\cite{alamily} is as a subgroup of the group $\Sym(1, \ldots, n, \overline{1}, \ldots, \overline{n})$ such that for $h \in B_n$, $\overline{h(i)} = h(\overline{i})$ for $i = 1, \ldots, n$.  So $B_n = \langle s_0 = (1\ \bar{1}), s_i = (i\ i+1)(\overline{i}\ \overline{i+1} \rangle, (i = 1, \ldots, n-1) \rangle$.  The base group $2^{n}$ is generated by $s_0$ and its conjugates in $\Sym(n) = \langle s_1, \ldots, s_{n-1} \rangle$ given by $(2\ \overline{2}) = (1 \overline{1})^{s_1}, (3\ \overline{3}) = (2\ \overline{2})^{s_2}, \ldots$.

The Weyl group of interest to us is $D_n := 2^{n-1} \sd \Sym(n) =  \langle s_1 s_0 s_1 s_0, s_1, \ldots, s_{n-1} \rangle$.  

The elementary abelian group
\index{group ! elementary abelian}%
 $2^{n-1} = \langle t_i = (i\ \overline{i})(i+1\ \overline{i+1}) \rangle$ $(1 \leq i \leq n-1)$ and  $\Sym(n) = \langle s_1, \ldots, s_{n-1} \rangle$ acts on the $t_i$ generators by conjugation, for example $t_1^{s_2} = (1\ \overline{1})(2\ \overline{2})^{s_2} = (1\ \overline{1})(3\ \overline{3}) = t_1 t_2$.  Notice that regarding $B_n \le \Sym(2n)$ gives $D_n = B_n \cap \Alt(2n)$.

The $D_n$ action on $2^{n-1} = \{ (a_1, \ldots, a_n) |\ a_i \in \{0, 1\} \}$ is via conjugation, permuting the places.  Switching changes vertex \emph{pairs} so we are restricted to considering $n$ is odd so that $n-1$ is even.  

In conclusion $S^{e}_{2,n}$ and $D_n$ are dual $\Sym(n)$ modules that are isomorphic and have structure $2^{n-1} \sd \Sym(n)$ if and only if $n$ is odd.  Further we can only have as $(n-1)$-tuples $(0, \ldots, 0)$ and $(1, \ldots, 1)$ in the base module of $D_n$.

\bigskip

A \emph{pseudo-reflection}~\cite{rouquier} 
\index{pseudo-reflection}%
of a vector space over $\mathbb{C}$ is a
finite order automorphism (not necessarily order $2$) whose fixed
point set is a hyperplane.  The irreducible complex reflection groups
\index{group ! complex reflection}%
classified by Shephard
\index{Shephard, G. C.}%
and Todd~\cite{sheptodd}
\index{Todd, J. A.}%
are two infinite series $\Sym(n+1)$
and $G(p, q, n)$ $(p > 1, q \ge 1, n \ge 1$ and $q|p$) and $34$
exceptional groups.  The group $G(p, 1, n) = (\mathbb{Z}_{p})^n \sd
\Sym(n)$ is the group of $n \times n$ monomial matrices with $p$th
roots of unity as non-zero entries, where $\Sym(n)$ is the subgroup of
permutation matrices
\index{permutation matrix}%
and $(\mathbb{Z}_{p})^n$ is the subgroup of
diagonal matrices.  Letting $s_0 = \diag(\zeta, 1, \ldots, 1)$,~\label{diag} where
$\zeta$ is a primitive $p$th root of unity and $s_i = (i\ i+1)$ gives
that $G(p, 1, n)$ is generated by the pseudo-reflections $\{s_0, s_1,
\ldots, s_{n-1}\}$ satisfying:
 \begin{displaymath}
\textrm{braid\ relations} = \left\{ \begin{array}{ll}
s_0 s_1 s_0 s_1 = s_1 s_0 s_1 s_0\\
s_i s_j = s_j s_i & \textrm{if $|i - j| > 1$}\\
s_i s_{i+1} s_i = s_{i+1} s_i s_{i+1} & \textrm{for $i \ge 1$}
 \end{array} \right.
\end{displaymath}
\index{braid relations}%
 \begin{displaymath}
\textrm{finite\ order\ relations} = \left\{ \begin{array}{ll}
s_0^p = 1\\
s_i^2 = 1 & \textrm{for $i \ge 1$}
 \end{array} \right.
\end{displaymath}

So $G(2, 1, n) = B_n$ $(n \ge 2)$ and $G(2, 2, n) = D_n$ $(n \ge 3)$.
So for odd $n$, $S^{e}_{2,n}$ is isomorphic to a complex reflection group.

Investigate further the links between switching groups and other known groups.  One example of a welding of two types of groups is~\cite{fenn} which contains a study of the subgroup of the automorphism group
\index{group ! automorphism}%
 of the free group generated by the braid group
\index{group ! braid}%
 and the permutation group.

\item  \emph{Switchings of $\mathbb{L}_L$ Roots}.  A lattice
in a real $d$-dimensional vector space is \emph{Euclidean} or
\emph{Lorentzian}
\index{Lorentzian lattice}%
\index{Euclidean lattice}%
\index{Lattice ! Euclidean}%
 if the overlying vector space is Euclidean,
$\mathbb{R}^d$, or Minkowski, $\mathbb{R}^{d-1, 1}$.  Even Lorentzian
lattices, $II^{d-1,1}$, exist only in dimensions $d = 8n +2$, $n$ an
integer, and are defined as consisting of those $x$ for which (i)
either $x \in \mathbb{Z}^{d-1, 1}$ or $x-l \in \mathbb{Z}^{d-1, 1}$
and (ii) $x \cdot l \in \mathbb{Z}$ where $l = (\frac{1}{2},
\frac{1}{2}, \ldots, \frac{1}{2}; \frac{1}{2})$.  If $\mathbb{L}$ is
a Euclidean even self-dual lattice, then $\mathbb{L} \oplus II^{1,1} =
II^{8n + 1,1}$, where $\dim \mathbb{L} = 8n$.  It has been proven~\cite{conwaysloa}
that a set of simple roots of the even Lorentzian lattice $II^{25,1}$ is a set of
points that are isometric to the Leech lattice.
\index{lattice ! Leech}%
(A set of simple roots is a minimal set of vectors having the property that the
reflections in the hyperplanes perpendicular to them generate the Weyl group
\index{group ! Weyl}%
 of the lattice, that is the group of all reflections which
are automorphisms of the lattice).  In short, there is a Lorentzian form for $\mathbb{L}_L$  that lies in the even integral unimodular lattice $II_{25, 1}$~\cite[Chapters 26, 27, 30]{conwayslo}.  

 The lattice $II_{25, 1}$ contains a set of vectors that are isometric to the $\mathbb{L}_L$ under the metric $d(r, s)^2 = N(r-s)$, where
\[ x \cdot y = x_0y_o + \ldots + x_{n-1}y_{n-1} - x_n y_n,\quad N(x) = x \cdot x.\]
The reflection subgroup of $\Aut(II_{25, 1})$ is the Coxeter group
\index{Coxeter ! group}%
\index{group ! Coxeter}%
with a generator $R_r$ for each $\mathbb{L}_L$ vector $r$ presented by
\[ (R_r)^2 =1,\]
\[ (R_r R_s)^2 = 1\quad \text{if}\ N(r-s)=4,\]
\[ (R_r R_s)^3 = 1\quad \text{if}\ N(r-s)=6.\]
The group of all autochronous automorphisms (those not interchanging
the positive and negative time cones) of $II_{25, 1}$ is this Coxeter
group generated by reflections in the Leech roots
\index{lattice ! Leech ! roots}%
 extended by the group $\Co_{\infty}$ of its diagram automorphisms.  The direct product
of this extension and a central inversion $-1$ gives the group
$\Aut(II_{25, 1})$.  The group of graph automorphisms ($\cong
\Co_{\infty}$) transitively permutes the walls of the fundamental
region for the Coxeter group (which are in $1$--$1$ correspondence
with the Leech roots); it is the group of all automorphisms of
$\mathbb{L}_L$, including translations.  Does the above Coxeter subgroup of $\Aut(II_{25, 1})$ correspond to a specifiable subgroup of $S_{3, \omega}$?

\item  \emph{Switchings in Other Areas}.
\index{graph ! switching}%
This is a general question.  Where in mathematics, other than in graph theory, do switching operations
\index{switching ! operation}%
and switching groups arise?

\item  \emph{Kac-Moody algebras}. 
\index{Kac-Moody algebra}%
A graph is \emph{N-free}
\index{graph ! N-free@$N$-free}%
 if it does not contain the path of length $3$ as an induced subgraph. Let $c(n)$ be the number of connected N-free graphs on $n$ vertices.

Then the following identity holds:

\[ \prod_{n>0} (1-q^{n})^{-c(n)}  = 1 - q + 2 \sum_{n>0} c(n)q^{n} \] 

This somewhat resembles in form the denominator identity of a Kac-Moody Lie algebra.  Is there is a nice connection?  (For each N-free graph there is a complementary one, so there is a bijection between connected and disconnected N-free graphs).  

There may be an approach to solving this problem using the theory of free Lie algebras~\cite{bourbaki}.
\index{Lie algebra ! free}%

Another clue may be found in the formula that generates the number of unlabeled, rooted trees~\cite{wilf}.  A \emph{rooted tree}
\index{tree ! rooted}%
 is a tree whose vertices are unlabeled, except a distinguished one which is the `root'.  If $h(n)$ is the number of all rooted forests on $n$ vertices, then
\[ \sum_{n} h(n)x^{n}  =  \prod_{n>0} (1-x^{n})^{-t(n)}, \] 
where $t(n) = h(n, 1)$ is the number of rooted trees of $n$ vertices.  Then since the number of rooted forests on $n$ vertices equals the number of rooted trees on $n+1$ vertices, the formula becomes
\[ \sum_{n} t(n+1)x^{n}  =  \prod_{n>0} (1-x^{n})^{-t(n)}. \] 

Another occurrence of a similar identity, arising in the study of integer sequences~\cite{cam2e}.  Suppose that $x = (x_n)$ enumerates a class $C$, and the operator $S$ is defined on the set of sequences of non-negative integers by $Sx = (y_n)$ such that 
\[ \prod_{n>0} (1-t^{n})^{- x_n}  = 1 + \sum_{n>0} y_n t^{n}. \] 
Then $Sx$ enumerates the class of disjoint unions of members of $C$ where, for $Sx$, the order of the ``component'' members of $C$ is unimportant, for example, building a structure such as a graph from its connected components.  If $x$ enumerates the generators (homogeneous elements) of a graded algebra that is a polynomial ring, them $Sx$ gives the dimensions of the homogeneous components.  If $x$ is realized by a group $G$ then $Sx$ is realized by $G \Wr S$.  Finally from the same paper~\cite{cam2e}, if $x$ lists the dimensions of homogeneous components of a graded vector space $V$, then $S$ lists those of the symmetric algebra of $V$ and the operator $S^*$ defined by
\[ \prod_{n>0} (1+t^{n})^{x_n}  = 1 + \sum_{n>0} (S^*x)_n t^{n} \]
lists those of the exterior algebra of $V$.

Finally we mention a reference to the Hopf algebra
\index{Hopf algebra}%
 of rooted trees~\cite{brouder}, and one on pre-Lie algebras and rooted trees~\cite{chapoton} (the free pre-Lie algebra
\index{Lie algebra ! free}%
  on one generator has a basis given by rooted trees).  

\item \emph{Unified Triality}.  Consider the following trialities:

(a) 
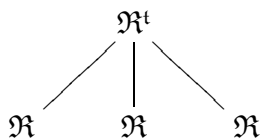
$$
\xymatrix{
& {\mathfrak{R^{t}}} \ar@{-}[d] \ar@{-}[dl] \ar@{-}[dr] \\
{\mathfrak{R}} & {\mathfrak{R}} & {\mathfrak{R}}
}$$

(b) 
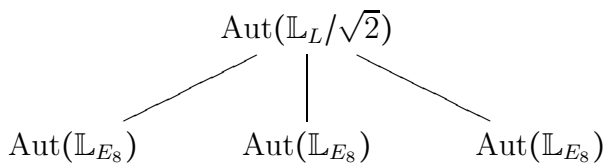
$$
\xymatrix{
& {\Aut(\mathbb{L}_{L}/ \sqrt{2})}  \ar@{-}[d] \ar@{-}[dl] \ar@{-}[dr]
  \\
{\Aut(\mathbb{L}_{E_8})} & {\Aut(\mathbb{L}_{E_8})} & {\Aut(\mathbb{L}_{E_8})}
}$$

(c) 
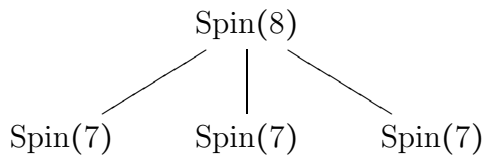
$$
\xymatrix{
& \Spin(8) \ar@{-}[d] \ar@{-}[dl] \ar@{-}[dr] \\
\Spin(7) & \Spin(7) & \Spin(7)
}$$
\index{group ! Spin(8)@$\Spin(8)$}%

(d) 
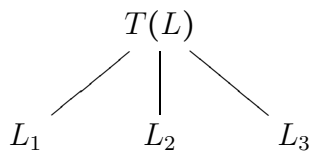
$$
\xymatrix{
& T(L) \ar@{-}[d] \ar@{-}[dl] \ar@{-}[dr] \\
L_1 & L_2 & L_3
}$$  

(a)  We can step down from  $\mathfrak{R^{t}}$ to one of the three  $\mathfrak{R}$s by going colourblind in any two of the three colours of $\mathfrak{R^{t}}$.

(b)    $\mathbb{L}_{L}$ is the Leech lattice and $\mathbb{L}_{E_8}$ is the $E_8$ lattice.  That three of the latter embed in one of the former, is clear from Lemma~\ref{lmlemma}.  Two lattices $\mathbb{L}, \mathbb{L'}$ are
\emph{isomorphic} if there is an orthogonal automorphism $u \in O(n, \mathbb{R})$ with $u(\mathbb{L})= \mathbb{L'}$.  If $\mathbb{L}=\mathbb{L'}$ then $u$ is an \emph{automorphism} of the lattice.  One further definition is worth stating as it mimics the homogeneity property of
\index{graph ! random}%
 random graphs: two overlattices $\mathbb{L}, \mathbb{L'}$ of a lattice $\Lambda$~\label{lamlattice} are \emph{isomorphic} if there is an automorphism of $\Lambda$ extending to an isomorphism of
$\mathbb{L}$ with $\mathbb{L'}$.

Modulo suitable normalization, $|\mathbb{L}_{L} / \mathbb{L}_{E_8}| = (2^{8})^3$~\cite{lepmeu}. 

In general, the root system of a lattice $\mathbb{L}$ spans a sublattice $\overline{\mathbb{L}}$ which has finite index in $\mathbb{L}$~\cite[p.~139]{milnorhus}.  In 24 dimensions the only exception is  $\mathbb{L}_L$ having  $\overline{\mathbb{L}_L} = 0$.  Furthermore, usually $\det(\overline{\mathbb{L}}) > 1$ so that $\overline{\mathbb{L}}$ is a proper sublattice of  $\mathbb{L}$.  This also has one exception, namely  $\mathbb{L} =  \overline{\mathbb{L}} = \mathbb{L}_{E_8} \oplus \mathbb{L}_{E_8} \oplus \mathbb{L}_{E_8}$, with root system $E_8 \oplus E_8 \oplus E_8$.  

(c)  The outer automorphisms of $\Spin(8)$ are the original well-known triality discovered by \'E. Cartan
\index{Cartan, \'E.}%
 in 1925, though with a geometric predecessor from 1903 due to E. Study.
\index{Study, E.}%



(d)  This refers to Grishkov's work~\cite{grishkov}
\index{Grishkov, A. N.}%
 on \emph{Lie algebras $\mathsf{L}$ with triality};~\label{mathsfL} see Appendix~\ref{LoopTheory}.
\index{Lie algebras with triality}%

Other than the 23 constructions of the Leech lattice based on a Niemeier lattice, there is also the Turyn construction using three copies of the $\mathbb{L}_{E_8}$, which follows the way that the binary Golay code
\index{Golay code}%
 can be constructed using three copies of the extended Hamming code, $H(8)$.
\index{Hamming code}%
   There are no doubt other examples, perhaps in fields such as design theory, Steiner triple systems and elsewhere; see Appendix~\ref{carsec}.

The substructures in our four above cases are not merely unions; a disjoint union of graphs is different to a direct sum of lattices, the latter being more like a transversal of graphs.  There are already approaches that connect some of the structures, such as~\cite{giudicia} where two infinite families of locally 3-arc transitive graphs
\index{graph ! locally $s$-arc transitive}%
 are constructed which have $\Aut(P\Omega^{+}(8,q))$ as their automorphism group.
 \index{group ! automorphism}%
   It is possible to argue that of the structures, the graph is the most basic.  Is there an ingredient or theory that describes the common form of the structures?

      
We already have some connections, which are summarized in Figure~\ref{mainfig}.  Here, $\mathbb{L}_{L}$ is the \emph{local} object corresponding to the global object $\mathfrak{R^{t}}$.  Note also that $\SL(2,2) \cong \PSL(2,2)$ and is its own Schur cover.  We list some associated open questions.

\emph{\textbf{Conjecture 1}}  \emph{Main Conjecture}.  We conjecture that it is possible to relate all these occurrences of triality, but the central question is how to prove this?  We already have some connections indicated in Figure~\ref{mainfig}, but is there an ingredient or theory that describes the common form of these structures?  

\emph{\textbf{Conjecture 2}}  \emph{Local-Global Connection}.   The lattice $\mathbb{L}_{L}$ is the \emph{local} object corresponding to the global object $\mathfrak{R^{t}}$.  Can we construct a Lie algebra with triality $T(L)$ to be the \emph{local} object corresponding to the global object comprising of the dense locally finite subgroup of $\Aut(\mathfrak{R^{t}})$?

A significant missing piece of the puzzle is how to link a purported Lie algebra $\mathsf{T(L)}$ to the triality graph $\mathfrak{R^{t}}$.

There are many approaches that connect some of the structures, such as Liebeck's~\cite{liebecka} which connects Lie algebras, two-graphs and permutation groups.  This may be a potential starting point.  In particular, is there a Liebeck Lie algebra on three colours having the form $\mathsf{T(L)} = \mathsf{L^{(2)}}_1 \oplus \mathsf{L^{(2)}}_2\oplus \mathsf{L^{(2)}}_3$, where $\mathsf{L^{(2)}}_i$ are isomorphic to a Liebeck-type Lie algebra on two colours?

 Take a finite graph $\Gamma = \{v_1, \ldots, v_n\}$, in which $N(v_i)$ denotes the neighbour set of vertex $v_i$, and suppose that $K$ is a field in which $|N(v_i)| - |N(v_j)| = 0$ for all $i, j$.  Define an algebra over $K$ with basis $\{v_1, \ldots, v_n\}$; Liebeck defines for $v_i, v_j \in V(\Gamma)$
\[\quad\quad v_i v_j : = \sum N(v_i) - \sum N(v_j), \]
that is the difference of the formal sums of neighbour sets of the vertices.  In characteristic zero, this Lie algebra derived from a regular graph is only nilpotent if $\Gamma$ is the null graph.

This definition, where the infinite sums are taken to have coefficients $1$ if the summed vertex is a neighbour and $0$ if it is a non-neighbour, does not readily extend to infinite graphs, so perhaps we can try a Lie algebra derived from field-valued functions on the vertices.  

We would first need to check the following properties are true of vertices in $\mathfrak{R}$ and $\mathfrak{R^{t}}$:-

(1)  Is a Lie bracket satisfied?

(2)  Does triality act in the right way, for example $\sigma$ effects $\mathfrak{b} \leftrightarrow \mathfrak{g}$, and $\rho$ effects $\mathfrak{b} \to \mathfrak{g} \to \mathfrak{r} \to \mathfrak{b}$ on edges in $\mathfrak{R^{t}}$?

(3)  Does the 3-colour case reduce to the 2-colour case?

(4)  Is the Lie algebra triality identity satisfied?

\emph{\textbf{Conjecture 3}}  \emph{Moufang Loops and Octonion Algebras}.  All octonion algebras satisfy the Moufang property so their loops of units are Moufang loops.  (In particular, the norm 1 split-octonions
\index{split-octonions}%
 over $K$ give rise to the triality group $P\Omega^{+}(8, K) \sd \Sym(3)$.)  Conversely, every known nonassociative simple Moufang loop,
\index{Moufang ! loop}%
 finite or infinite, arises as the central quotient of the norm 1 units from some octonion algebra~\cite{hallji2}.
\index{octonions}%
 Is this true of the Moufang loop constructed in Theorem~\ref{moulthm}, and how exactly does this connect $\mathfrak{R^{t}}$ with the octonions?

\emph{\textbf{Conjecture 4}}  \emph{Connection to $D_4(K)$}.
 \index{group ! D(4)@$D(4)$}%
There is a conjecture that a simple (infinite) group $G$ admits non-trivial triality $\langle \sigma, \rho \rangle$ if and only if $G \cong D(4, K)$, where $K$ is a field, and $\langle \sigma, \rho \rangle$ is the group of graph automorphisms of $G$.  Grishkov and A. V. Zavarnitsine make the following comment about this conjecture~\cite{grish1}:  ``This conjecture is both important and difficult. The only encouraging fact is that the corresponding problem for simple Lie algebras with triality was solved in the affirmative (see~\cite{grishkov})''. 

\emph{\textbf{Conjecture 5}}  \emph{Simple Centralizing Groups}.  
This conjecture is due to Nagy and P. Vojt\v{e}chovsk\'y~\cite{nagyvoj1}.  If $G$ is a simple group with triality $S = \langle \sigma, \rho \rangle$ such that $G = [G, S]$ and $\zeta_S(G) = 1$, then $C_G(\sigma) = \{ x \in G : x^{\sigma} = x \}$ and $C_G(\rho) = \{ x \in G : x^{\rho} = x \}$ are simple groups.

\emph{\textbf{Conjecture 6}}  \emph{Different Approaches to Triality}.  
Is there different approach that explicitly uses Clifford algebras?
\index{Clifford algebra}%


If so, can our lattice constructions of random graphs be applied to the formation by Wilson~\cite{wilsona}
\index{Wilson, R. A.}%
 of a 3-dimensional octonionic Leech lattice
\index{lattice ! Leech ! octonionic}%
 as the set of triples of octonions satisfying certain conditions, where the definition of the lattice is invariant under permutations of the three coordinates.  We note other attempts to build the Leech lattice from triples of integral octonions including~\cite{dixona}~\cite{dixonb}~\cite{elkies}; each of these can potentially be related to our work.

   A second plan of attack in proving the connection would be to
     adapt to our situation the theory of F. Zara~\cite{zara}
\index{Zara, F.}%
      who has shown how to
     construct Clifford algebras of Coxeter groups
     \index{Coxeter ! group}%
\index{group ! Coxeter}%
which are connected
     to Fischer systems
\index{Fischer system}%
\index{Fischer, B.}%
      and in turn to connected graphs.  In our case
     the group in question is the $3$-transposition group generated by
     the three types of switchings associated with
     $\mathfrak{R^{t}}$, and the connected graph is the affine diagram
      $\tilde{A_{2}}$, this being the smallest diagram \emph{not} associated with a
     finite Coxeter group, with the switchings acting upon or representing the vertices.
    
    A third possible approach would be to construct $\mathfrak{R^{t}}$
    using the definition of the icosian ring
\index{icosian ring}%
 discovered by P. L. H. Brooke
\index{Brooke, P. L. H.}%
     and given in~\cite[p.160]{wilsona}, which we now
    repeat.  Let $V$ be the 4-dimensional rational vector space spanned by vectors $G, H, I, J, K$ (permuted by $\Alt(5)$) subject to $G+H+I+J+K = 0$, with inner products $(X,X) = 2, (X, Y) = - \frac{1}{2}$ for $X \ne Y$.  The Clifford algebra $CA(V)$ is the quotient of the tensor algebra $1 \oplus V \oplus V^2 \oplus \ldots$ by the ideal $\langle v \otimes w + w \otimes v - 2(v,w) \rangle$ and has two 8-dimensional subalgebras over $\mathbb{Q}$.  It has an even subalgebra that is the image of $1 \oplus V^2 \oplus V^4$ and isomorphic to the icosian ring, and an odd algebra that is the image of $V \oplus V^3$.

\emph{\textbf{Conjecture 7}}  \emph{Ternary Quasigroups $\&$ the Modular Group}. In~\cite{smith}, J. D. H. Smith
\index{Smith, J. D. H.}%
\index{quasigroup}%
 unified the two definitions of $n$-quasigroups, the combinatorial one (amounting to Latin $n$-cubes), and the one given by $2n$ identities on $n + 1$ different $n$-ary operations, under the name of \emph{hyperquasigroup}.
\index{hyperquasigroup}%
   The combinatorial definition is a set $Q$ with an $n$-ary multiplication $\mu : Q^n \to Q: (x_n, \ldots, x_1) \mapsto x_n \ldots x_1 \mu$ such that, specification of any $n$ coordinates in an $(n+1)$-tuple satisfying $x_n \ldots x_1 \mu = x_0$ determines the remaining one uniquely.  For the free product
\index{group ! free product}%
  groups defined by $M_n := \langle \xi, \sigma | \xi^n = 1 = \sigma^2 \rangle,$
where $M_1 = C_2$, $M_2 = D_{\infty}$, $M_3 = \PSL(2, \mathbb{Z})$,
\index{group ! modular}%
 $\ldots$, there is an epimorphism $r : M_n \to \Sym(n+1): \xi \mapsto (1\ 2\ \ldots\ n), \sigma \mapsto (0\ 1)$.  Smith defines an \emph{$n$-space} $(G, \xi, \sigma)$ to be a set $G$ with a shift $\xi$ and an involution $\sigma$, such that $\xi : G \to G : g \mapsto \xi g$, $\sigma : G \to G : g \mapsto \sigma g$.  Then $\Sym(4) =  \langle \xi, \sigma | \xi^3 = \sigma^2 =  (\xi \sigma)^4 = 1 \rangle$ is a ternary space, where $\xi = (1\ 2\ 3)$ and $\sigma = (0\ 1)$ are left multiplications within $\Sym(4)$.  Smith shows that ternary quasigroups arise from actions of the modular group.  Can this be linked to our construction of $\mathfrak{R^{t}}$ as a homogenous Cayley graph for the index-3 subgroup $C_2 * C_2 * C_2$ of the modular group? Smith extended his work in~\cite{smith1} where his hyperquasigroups were characterized in terms of graphs, and the $\Sym(3)$-triality action was implemented in terms of the hyperquasigroups.

\clearpage

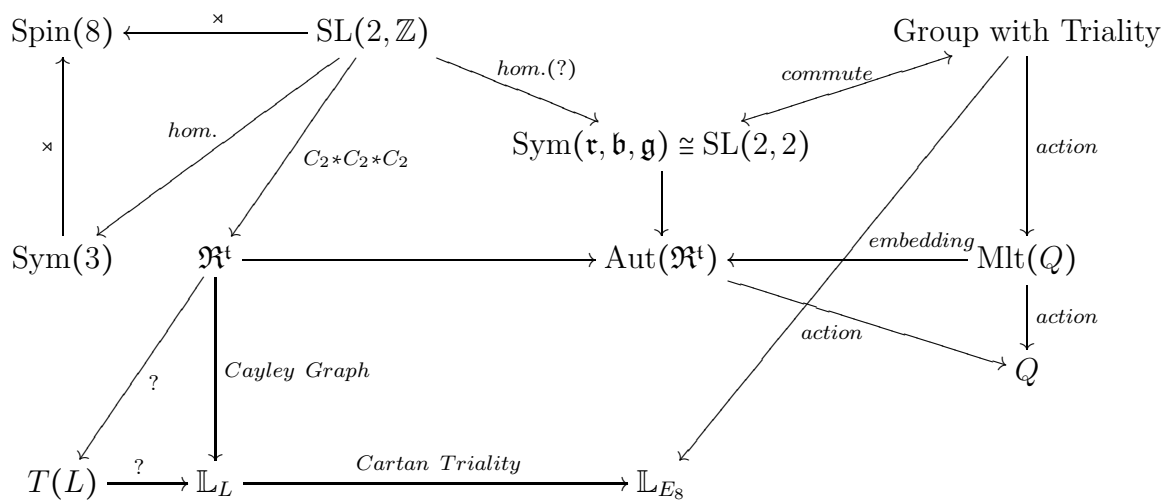
\begin{figure}[ht]
\vbox to \vsize{
 \vss
 \hbox to \hsize{
  \hss
  \rotatebox{90}{
   \vbox{
    \hss
\[
\xymatrix{
{\Spin(8)}  \ar@{<-}[rr]^{\sd} && {\SL(2, \mathbb{Z})}  \ar@{->}[ldd]^{C_2 \ast C_2 \ast C_2}  \ar@{->}[rd]^{hom. (?)} \ar@{->}[lldd]_{hom.}  &&  {\text{Group\ with\ Triality}} 
\ar@{->}[dd]^{action} \ar@{->}[ldddd]\\
&&& {\Sym(\mathfrak{r}, \mathfrak{b}, \mathfrak{g})  \cong \SL(2,2)} \ar@{<->}[ru]^{commute} \ar@{->}[d]\\
{\Sym(3)} \ar@{->}[uu]^{\sd} & {\mathfrak{R^{t}}}   \ar@{->}[rr] \ar@{->}[ldd]^{?} \ar@{->}[dd]^{Cayley\ Graph} &&  {\Aut(\mathfrak{R^{t}})}  
\ar@{->}[dr]_{action}    \ar@{<-}[r]^{\quad\quad\quad \quad\quad embedding}  &    {\Mlt(Q)} \ar@{->}[d]^{action}\\
&&&& Q\\
{T(L)}  \ar@{->}[r]^{?} &   {\mathbb{L}_{L}}   \ar@{->}[rr]^{Cartan\ Triality}    &&  {\mathbb{L}_{E_8}}
}
\]
    \hss
\caption{Diagram of Some Known Connections Related to Triality}
\label{mainfig}
        }
                }
  \hss          }
 \vss           }
\end{figure}

\clearpage

\item  \emph{Arnold's Trinities}.  Are there connections with Arnolds
\index{Arnold, V. I.}%
 programme of connecting \emph{trinities} that arise in different parts of mathematics, for example his attempt~\cite{arn} to colour braid groups
\index{group ! braid}%
 and complexify and quaternionize permutation groups?  Is there an analytical representation of permutations that can be analytically continued into the complex plane, perhaps via topological groups?

\item  \emph{Triality of $\mathfrak{R^{t}}$}.  As an extension of Theorem~\ref{outgpwttr}, investigate the possibility of a \emph{dense locally finite} subgroup $G$ of the group $\Aut(\mathfrak{R^{t}})$ such that $(G,
T(\mathfrak{R^{t}}))$ is a group with triality, with $T(\mathfrak{R^{t}})) \cong \Sym(\mathfrak{r}, \mathfrak{b}, \mathfrak{g})$.

\item \emph{Connection with modular group $(\gimel)$ action on $\hat{\mathbb{Q}}$}.  It is well-known that $\mathfrak{R}$ and $\mathbb{Q}$ have many links, for instance model-theoretic homogeneity,
$\aleph_0$-categoricity,
\index{aleph@$\aleph_0$-categorical}%
a topology
\index{topology}%
 on $\mathfrak{R}$ that is
homeomorphic to $\mathbb{Q}$ with the usual metric topology.  Also
$\mathbb{Q}$ is a direct limit
\index{direct limit}%
 of infinite cyclic groups
\index{group ! cyclic}%
 and $\mathfrak{R}$ is a Cayley graph
\index{graph ! Cayley}%
 for $C_{\infty}$.  The class of finite ordered graphs
\index{graph ! ordered}%
 has as Fra\"{\i}ss\'e limit
 \index{Fra\"{\i}ss\'e limit}%
 the random graph
\index{graph ! random}%
  with an appropriate linear ordering
\index{linear order}%
 isomorphic to $\mathbb{Q}$, which might be called the \emph{random ordered graph}~\cite{kechris}.
\index{graph ! random ! ordered}%
 Finally, Bodirsky and Pinsker
\index{Bodirsky, M.}%
\index{Pinsker, M.}%
  studied~\cite{bodirsky} the random ordered graph
\index{graph ! random ! ordered}%
 in which the vertex set is equipped with a total order to yield the unique countably infinite homogeneous graph containing all finite ordered graphs.  The order is then isomorphic to the order of the rationals.  Incidentally, Bodirsky, Pinsker and Pongr\'acz
\index{Bodirsky, M.}%
\index{Pinsker, M.}%
\index{Pongr\'acz, A.}%
  determined up to first-order interdefinability, the 42 reducts of the random ordered graph,
\index{graph ! random ! ordered}%
 this being the \emph{free superposition} (defined in the paper~\cite{bodirsky1}) of $(\mathbb{Q}, <)$ and the random graph, giving the unique countable homogeneous linearly ordered graph that embeds all finite linearly ordered graphs up to isomorphism.  This result is connected to the Kechris-Pestov-Todorcevic 
\index{Kechris, A. S.}%
\index{Pestov, V.}%
\index{Todorcevic, S.}%
theorem in topological dynamics~\cite{kechris}; see $\S~\ref{TheoryofRelationalStructures}$ of the Prerequisite Background Appendices.

Is there a connection between modular group action
\index{group ! action}%
 on $\mathfrak{R}$ and that on  $\hat{\mathbb{Q}} = \mathbb{Q} \cup \{\infty\}$?
 
For the reader's convenience we make some notes.  G. Jones
\index{Jones, G. A.}%
 and co-workers found~\cite{josiwi} that the action  of $\gimel$ on
 $\hat{\mathbb{Q}}$ where by convention
 $g \in \gimel : \infty \mapsto a/c$ and $a/c = \infty$ if
 $c = 0$ is
 transitive, and that the stabilizer of a point is $C_{\infty} = \Bigg\langle
\begin{pmatrix} 
1&1\\
0&1\\
\end{pmatrix} \Bigg\rangle$, which might be made to tally with the regular action of $C_{\infty}$ on $\mathfrak{R}$ if it stabilizes one colour in $\mathfrak{R^{t}}$.  Now if $\gimel_{\infty}$ denotes the stabilizer of $\infty$
then subgroups 
$H$ such that $\gimel_{\infty} < H \le \gimel$ produce
$\gimel$-invariant equivalence relations on $\hat{\mathbb{Q}}$.  
 
 Letting $\approx_n$ denote the $\gimel$-invariant equivalence relation
induced on $\hat{\mathbb{Q}}$ by the congruence subgroups
\index{group ! congruence subgroup}%
 $\gimel_{0}(n)$, it was found that for $p_1,
p_2 \in \hat{\mathbb{Q}}$,  $p_1 \approx_n p_2$ if and only if $p_1 - p_2 \equiv
0 \pmod{n}$, that is $p_1, p_2$ have the same reduction$\pmod{n}$.  For
example for $H = \gimel_{0}(n)$, the inclusion $H \leq \gimel$ is strict if $n > 1$, so $\gimel$ acts
imprimitively
\index{group ! permutation ! imprimitive}%
 on $\hat{\mathbb{Q}}$, the number of equivalence classes
being $|\gimel : \gimel_{0}(n)|$.  When $n$ is a prime $p$, there are
$p + 1$ blocks.  The action of $\gimel$ on the blocks is the same as
that of its quotient group $\PSL(2, p)$ on the projective line~\label{projln}
\index{projective line}%
 $GF(p) \cup \{\infty\}$.  For $p = 2$, $\gimel$ acts as $\PSL(2, 2) \cong
\Sym(3)$ on the blocks of imprimitivity $[0], [1]$ and $[\infty]$,
where $[a]$ denotes the block containing point $a$, and $[\infty] =
\{x / y \in \hat{\mathbb{Q}}\ |\ y \equiv 0 \pmod{n} \}$.  

Now take an infinite trivalent graph and label its vertices by those elements $p_i
\in \hat{\mathbb{Q}}$ such that for each pair $p_i, p_j \in
\hat{\mathbb{Q}}$ satisfying $p_i \approx_n p_j$ form the edge
$\{p_i, p_j\}$.  Firstly, prove that $\gimel$-invariance is equivalent to a
regular group action
\index{group ! action}%
 on graph vertices.  The $\gimel$-invariant equivalence relation  $\approx_2$ is then the stabilization of an edge rather than an edge colour.  Show that colouring the edges at random from  $\{\mathfrak{r} , \mathfrak{b} , \mathfrak{g}\}$ taking $p=2$ means that $\gimel$ acts as $\PSL(2, 2) \cong \Sym(\mathfrak{r} , \mathfrak{b} , \mathfrak{g})$ on the colours and as $\gimel_{0} (2)$ on the vertices.  By Theorem $4.2.5 (2)$ of~\cite{miyake} $\gimel_{0} (2)$ has index $3$ in $\gimel$ and so is the unique such group, so $\gimel_{0} (2) \cong C_2 * C_2 * C_2$, therefore $\mathfrak{R^{t}}$ can be a Cayley graph for $\gimel_{0} (2)$.
\index{graph ! Cayley}%

Another way to see the link is as follows~\cite{jones}.  Let $\gimel = \langle X, Y | X^2 = 1 = Y^3 \rangle$ where $X$ and $Y$ are the images in $\gimel$ of the matrices $\begin{pmatrix} 
0&-1\\
1&0\\
\end{pmatrix}$ and $\begin{pmatrix} 
0&1\\
-1&1\\
\end{pmatrix}$.  If $G$ is a subgroup of index $n$ in $\gimel$, then $\gimel$ acts transitively on the coset space $\gimel / G$.  Let the element $Z = (XY)^{-1} = \pm \begin{pmatrix}
1&1\\
0&1\\
\end{pmatrix}$ have $t$ cycles where $t$ is called the \emph{parabolic class number} of $G$.  Now $t$ is equal to the number of orbits of $G$ on $\hat{\mathbb{Q}}$, and also equal to the number of $G$--conjugacy classes of maximal parabolic subgroups in $G$.  Since $C_3$ is abelian,
\index{group ! abelian}%
 its order equals the number of its conjugacy classes.

\item  \emph{Connection with Conway's work on $\gimel_{0}(n)$}.
\index{Conway, J. H.}%
\index{group ! congruence subgroup}%
  If we consider the upper half complex plane then $\mathfrak{R^{t}}$ can be represented as a
random $3$-colouring of a tree whose leaves are the non-infinite edges of the
fundamental regions for the action of $\PSL(2, \mathbb{Z})$.
\index{group ! modular}%
  Can we use this observation, together with Conway's work~\cite{conway1} to help formulate some links between $\mathfrak{R^{t}}$ and the modular group?

We mention some elements of his work.  Two groups or lattices are \emph{commensurable}
\index{group ! commensurable}%
\index{lattice ! commensurable}%
 if their intersection has finite index in each of them.  

The \emph{commensurability subgroup} of $H < G$ is defined by $\Comm_G (H) : = \{g \in G | g^{-1}Hg\ \text{and}\ H\ \text{are commensurable}\}$.~\label{Comm_G (H)}

If $H_1, H_2 < G$ are
commensurable, then $\Comm_G (H_1)=\Comm_G (H_2)$.  Note the following
subgroup inclusions: $H \le N_G (H) \le \Comm_G (H) \le G$.  The modular group $\gimel$ is the automorphism group
\index{group ! automorphism}%
 of a $2$-dimensional fiducial lattice,
call it $\mathbb{L}_1 = \langle \underline{e}_1, \underline{e}_2 \rangle$ acting on it by
matrix multiplication, that is change of basis.  Conway studied the
groups commensurable with $\gimel$, for example the congruence subgroups
\index{group ! congruence subgroup}%
 $\gimel_{0}(n)$, in terms of their actions on the lattices commensurable with
$\mathbb{L}_1$.  

In a \emph{projective class} of lattices,
\index{lattice ! projective class}%
 $\langle \lambda \underline{e}_1, \lambda \underline{e}_2 \rangle$ and $\langle
\underline{e}_1, \underline{e}_2 \rangle$ are equivalent for all
$\lambda \in \mathbb{Q} \backslash \{0\}$.  For each integer $n$,
$\mathbb{L}_1$ and $\mathbb{L}_n = \langle n \underline{e}_1, \underline{e}_2 \rangle$
each contain a rescaled copy of the other to index $n$.  For each projective class
assign a node of a $(p+1)$-valent tree, with joins corresponding to
containments at prime index.  The stabilizer of any point is a
conjugate of $\gimel$, and the joint stabilizer of $\mathbb{L}_1$ and $\mathbb{L}_n$ is
$\gimel_0(n)$.  Any two commensurable lattices can by a basis
change become $\mathbb{L}_1$ and $\mathbb{L}_n$ for some $n$ and some conjugate of
$\gimel_0(n)$ stabilizes any two commensurable lattices.  A graph can be formed by joining each pair of lattices commensurable
with $\mathbb{L}_1 = \langle \underline{e}_1, \underline{e}_2 \rangle$ by a
line.  Finally, both $\gimel_0(n)$ and its normalizer in $\PSL(2, \mathbb{R})$ are
relevant to the moonshine properties of Thomson series,
\index{Thomson series}%
 and are required to describe monstrous moonshine.  
\index{monstrous moonshine}%

\item \emph{Almost Stabilizer of $\mathfrak{R}$}.
We will motivate and at the end state the open problem.

We mention at the outset some permutational similarities between $\Aut(\mathbb{Q}, <)$ and  $\Aut(\mathfrak{R})$:  they are both oligomorphic on their domains, primitive, transitive, and not even $2$-transitive.
\index{group ! permutation ! oligomorphic}%
\index{group ! permutation ! primitive}%
\index{group ! permutation ! $2$-transitive}%
 The differences include that $\Aut(\mathbb{Q}, <)$ is both torsion-free and highly homogeneous, whilst $\Aut(\mathfrak{R})$ is not.
\index{group ! torsion-free}%
\index{group ! permutation ! highly homogeneous}%

Macpherson
\index{Macpherson, H. D.}%
 has a construction of maximal subgroups of $\Aut(\mathbb{Q}, <)$ \cite[p.265]{macpherson}.  Maximal subgroups of $\Aut(\mathfrak{R})$ include stabilizers of finite sets of vertices, unordered edges and unordered non-edges.  (Observe that two obvious non-conjugate maximal subgroups of $\Aut(\mathbb{Q}, <)$ are the stabilizers of a rational and of an irrational; maximality of these groups follows from the primitivity of the corresponding permutation representations \cite[p.266]{macpherson}).

A collection $\mathscr{I}$ of subsets of a countably infinite set $\Omega$ is an \emph{ideal}
 \index{ideal}%
if 

(a)  $\emptyset \in \mathscr{I}, \Omega \notin \mathscr{I}$;

(b)  if $T \in \mathscr{I}$ and $S \subseteq T$ then $S \in \mathscr{I}$;

(c)  if $T, S \in \mathscr{I}$ then $T \cup S \in \mathscr{I}$.

An ideal $\mathscr{I}$ is \emph{maximal} if it is contained in no larger ideal; that is, for every $S \subset \Omega$, either $S \in \mathscr{I}$ or $\Omega \backslash S \in \mathscr{I}$.  Every ideal $\mathscr{I}$ has a dual \emph{filter},
\index{filter}%
 $\mathscr{F} := \{\Omega \backslash S : S \in \mathscr{I} \}$, and vice versa, and so the stabilizer of an ideal is the same as the stabilizer of a filter.  The dual of a maximal ideal is an \emph{ultrafilter}.
\index{ultrafilter}%
  If $\mathscr{H} \subseteq \mathcal{P}(\Omega)$ and $\Omega$ is not a finite union of members of $\mathscr{H}$, then the closure of $\mathscr{H}$ under subsets and finite unions is the ideal $\langle \mathscr{H} \rangle$ \emph{generated by} $\mathscr{H}$ and is the smallest ideal containing $\mathscr{H}$.

With this notation we have the following two definitions:-

\[ S_{\{\mathscr{H}\}} : = \{ g \in \Sym(\Omega) : (\forall T \subseteq \Omega)(T \in \mathscr{H} \leftrightarrow T^g \in \mathscr{H} )\}. \]
\[ S_{(\mathscr{I})} : = \{ g \in \Sym(\Omega) : \supp(g) \in \mathscr{I} \}. \]

Macpherson
\index{Macpherson, H. D.}%
proved the following result as an extension of a theorem of Richman~\cite{richman}, (see Appendix~\ref{PermutationGroups}):
\index{Richman, F.}%
\begin{theorem}
\label{macidealthm}
Let $\mathscr{I}$ be an ideal on $\Omega$ containing a moiety.  If $S_{\{\mathscr{I}\}}$ has three orbits on moieties of $\Omega$, then $S_{\{\mathscr{I}\}}$ is maximal in $\Sym(\Omega)$.
\end{theorem}

If $\Omega = \mathbb{Q}$ then one ideal that satisfies this theorem consists of all \emph{scattered} suborderings of $\mathbb{Q}$, that is, subsets which do not themselves embed a copy of $(\mathbb{Q}, \leq)$~\cite{macpherson}.

Covington, Macpherson and Mekler,
\index{Covington, J.}%
\index{Macpherson, H. D.}%
\index{Mekler, A.}%
 who were interested in ideals whose stabilizer is maximal, then proved~\cite{covingtonmac}

\begin{theorem}
Let $\mathcal{M}$ be a countably infinite indivisible homogeneous structure
\index{structure ! indivisible}%
 over a finite relational language, let $\Omega = M$ and $S = \Sym(M)$, and assume that $S \neq \Aut(\mathcal{M})$.  Define an ideal $\mathscr{I}$ to be the set of subsets $A \subseteq M$ such that the structure induced on $A$ does not embed a copy of $\mathcal{M}$.  Then
\begin{itemize}
\item[(a)]  $\mathscr{I}$ is an ideal on $\Omega$;
\item[(b)]  $\mathscr{I}$ contains a moiety of $\Omega$;
\item[(c)]  if the age
\index{age}%
 $\mathcal{A}(\mathcal{M})$ has the strong amalgamation property (SAP)
\index{amalgamation property ! strong}%
 then $S_{\{\mathscr{I}\}}$ has three orbits on moieties of $\Omega$ so is maximal in $\Sym(\Omega)$.
\end{itemize}
\end{theorem}

(A recent study has been made~\cite{laflamme1} in which a countable dimensional vector space over an $\mathbb{Q}$ provides an example of an age indivisible, but not weakly indivisible relational structure.)
\index{relational structure}%

In~\cite{brazil} and~\cite{macpherson}, non-maximal ideals are constructed, whose stabilizers $S_{\{\mathscr{I}\}}$ are maximal in $\Sym(\Omega)$.  In~\cite{covingtonmac} more such ideals were given, where the ideal was related to combinatorial structures with domain $\Omega$, such as $(\mathbb{Q}, <)$ and $\mathfrak{R}$.  The ideal is invariant under the automorphism group
\index{group ! automorphism}%
 of the structure, and the associated maximal subgroup of $\Sym(\Omega)$ can be viewed as an \emph{almost stabilizer}
\index{structure ! almost stabilizer}%
 of the structure, that is a group of \emph{almost automorphisms} of the structure.
\index{group ! almost automorphism}%
The last theorem applies to $\mathfrak{R}$ and to Henson graphs.
\index{graph ! Henson}%

By way of studying the extent to which non-isomorphic structures yield non-conjugate maximal subgroups, the following question has been asked~\cite{covingtonmac}: is it possible to put the structure of the random graph
\index{graph ! random}%
 on $\mathbb{Q}$ so that the maximal subgroup almost stabilizing $(\mathbb{Q},  \leq)$ is equal to the almost stabilizer of the random graph, in the present sense?

\item \emph{Random Cyclic Metric Space}.
\index{metric space ! random ! cyclic}%
A countable metric space admitting a cyclic automorphism permuting the 
points in a single cycle is specified by a function $f: \mathbb{N} \to \mathbb{R}_{+}$ which satisfies 
\[ \quad\quad\quad max_{x < n} ( | f(x) - f(n - x)| ) \leq f(n) \leq min_{y < n} (f(y) + f(n - y) ). \]
 
(Identify the vertex set with $\mathbb{Z}$ and set $d(x, y) = f(| x - y |)$.)  Assume that the values of $d(x, y)$ for $x \neq y$ lie in the set $\{1, 2, \ldots, d \}.$  If we choose $f(n)$ for $n = 1, 2, \ldots$, the displayed inequalities never conflict.  Suppose that we choose $f(n)$ randomly from the allowed values.  What can be said about 
the resulting 'random cyclic metric space'? (For $d = 2$ we obtain the path metric in the countable random graph $\mathfrak{R}$: see~\cite{rado}).

\item \emph{Limits of cubes}.
The celebrated Urysohn space
\index{Urysohn space}%
 is the completion of a countable universal homogeneous metric space which can itself be built as a direct limit
\index{direct limit}%
  of finite metric spaces.  In the paper~\cite{camtar1}, other examples are given of spaces constructed by forming the direct limit and completion, where the finite spaces are scaled hypercubes.
\index{hypercube}%
 The resulting countable spaces provide a context for a direct limit of finite symmetric groups with strictly diagonal embeddings (see~\cite{kroshko}), acting naturally on a module which additively is the ``Nim field''~\cite{conway2}
\index{Nim field}%
  (the quadratic closure of the field of order~$2$). Its completion is familiar in  another guise: it is the set of Lebesgue-measurable subsets of the unit interval modulo null sets with distance being the Lebesgue measure of the symmetric difference.  The isometry groups of these spaces are described as well as some interesting subgroups, and some generalizations are given.  Here we include some of the problems for further research that are in the mentioned paper.

(i)  The union of a chain of scaled hypercubes is the space $\mathcal{H}_{\omega}$~\label{hycube}, and this is represented on the set of non-negative integers by a recursively defined function.  Is it possible to give an explicit non-recursive formula for this function?

(ii)  \textbf{Cycles and Gray codes.}
\index{Gray code}%
The space $\mathbb{Q}\mathfrak{U}$ admits an isometry permuting all the points
in a single cycle. It follows that Urysohn space $\mathfrak{U}$ has an
isometry all of whose cycles are dense. See~\cite{camver}. Does anything
similar happen for the completion $\mathcal{H}$ of the infinite-dimensional cube $\mathcal{H}_{\omega}$~\label{comhycube}?

It is tempting to think that, even if such isometries don't exist, the
existence of Gray codes (that is, Hamiltonian cycles)
\index{Hamiltonian cycle}%
 in finite cubes should
imply the existence of something similar in the limit spaces (perhaps something
like a space-filling curve).

The space $\mathcal{H}$ cannot have a space-filling curve in the usual sense,
since it is not compact. We do not know how to proceed.

A related question would be the existence of something like a space-filling
curve in the ``middle level'' of $\mathcal{H}$, the set of points lying at
distance $\frac{1}{2}$ from $0$. Here, the existence of the analogous object
in finite cubes is a difficult combinatorial problem, as yet unsolved
(see~\cite{jrj}).

(iii)  \textbf{Other Hamming spaces.} 
Let $H(n,q)$ denote the Hamming space
\index{Hamming space}%
consisting of all $n$-tuples over an alphabet of size~$q$~\label{que}. As usual, the
distance between two $n$-tuples is the number of coordinates where they
differ. If we let $\mathcal{H}_n(q)$ denote the scaled Hamming space
$\frac{1}{n}H(n,q)$, then we have isometric embeddings
\[\mathcal{H}_1(q)\to\mathcal{H}_2(q)\to\mathcal{H}_4(q)\to\cdots\]
with union $\mathcal{H}_\omega(q)$ and completion $\mathcal{H}(q)$.
\begin{itemize}
\item Is there a convenient representation of $\mathcal{H}(q)$? What is the
structure of its isometry group?
\item Is it true that the set of points at distance $1$ from any point of
$\mathcal{H}(q)$ is isometric to $\mathcal{H}(q-1)$?
\item If we modify the embedding to take $\mathcal{H}_{n_i}(q)$ to
$\mathcal{H}_{n_{i+1}}(q)$, where $n_{i+1}/n_i=r_i$, with $(r_i)$ any given
sequence of integers greater than~$1$, is it the case that the completion
of the union is isometric to $\mathcal{H}(q)$, independent of the choice of
sequence $(r_i)$?
\end{itemize}

(iv)  \textbf{Philip Hall's locally finite group.} Philip Hall~\cite{hallp}
\index{Hall, P.}%
\index{group ! locally finite}%
\index{group ! Philip Hall's}%
constructed a universal homogeneous locally finite group as follows.
Embed $\Sym(n)$ into $\Sym(n!)$ by its regular representation, and
take the union of the sequence
\[\Sym(3)\to\Sym(6)\to\Sym(720)\to\cdots.\]
This group is countable and locally finite; it contains an isomorphic copy
of every finite group, and any isomorphism between finite subgroups is
induced by an inner automorphism of the group.

We can construct a limit of cubes to mirror this construction, so that
Hall's group acts on the union. Consider $\mathcal{H}_n$, with the
coordinates indexed $0,1,\ldots,n-1$ as usual. We will take the coordinates
of $\mathcal{H}_{n!}$ to be indexed by elements of $\Sym(n)$. Any
subset $S$ of $\{0,1,\ldots,n-1\}$ is mapped to the subset
\[\pi(S)=\{g\in\Sym(n):0^g\in S\}\]
of $\Sym(n!)$.

The embedding is an isometry, because $|\pi(S)|/n!=|S|/n$. Since
\[g\in\pi(S)\Leftrightarrow 0^g\in S \Leftrightarrow 0^{gh}\in S^h
\Leftrightarrow gh\in\pi(S^h),\]
we have $\pi(S^h)=\pi(S)h$, and so $\pi$ intertwines the natural action of
$\Sym(n)$ on $\{0,\ldots,n-1\}$ with its action on itself by right
multiplication. Hence Hall's group acts on the union of this chain of
cubes. We propose the name \emph{Hall cube}
\index{Hall cube}%
 for this space. Its completion is $\mathcal{H}$ (the construction above agrees with a more general construction given of $\mathcal{H}$).

What properties does the Hall cube have, and how does Hall's group act on it?

(v)  \textbf{Other embeddings of metric spaces.} We can play the same game
with other chains of finite metric spaces with lots of symmetry (for
example, scaled versions of distance-transitive graphs). 
\index{graph ! distance-transitive}%

One example involves the \emph{dual polar spaces}
\index{dual polar space}%
 $D_n(K)$ of type $D_n$ over a field~$K$ (see~\cite{camc}).
The points of such a space are the maximal totally singular subspaces of
a vector space of dimension~$2n$ over the field $K$ carrying the quadratic
form
\[Q(x_0,x_1,\ldots,x_{2n-1}) = x_0x_1+x_2x_3+\cdots+x_{2n-2}x_{2n-1}.\]
The distance between two subspaces is the codimension of their intersection.
There is a natural embedding of $D_n(K)$ into $D_{n+1}(K)$, as the
set of maximal totally singular subspaces containing a fixed $1$-dimensional
singular subspace. This embedding is the analogue of our embedding of
$H(n,2)$ into $H(n+1,2)$ without re-scaling.

A more interesting possibility would involve re-scaling, embedding
$D_n(K)$ into $D_{2n}(K)$. One possibility would be to let $\overline{K}$ be a quadratic
extension of $K$; then $D_n(K)$ is embedded in $D_n(\overline{K})$ (by tensoring the
underlying vector space with~$\overline{K}$), and $D_n(\overline{K})$ is
embedded in $\frac{1}{2}D_{2n}(K)$ by restriction of scalars. This is a
close analogue of the cubes considered in~\cite{camtar1}.

\item  \emph{Urysohn space}.
\index{Urysohn space}%

\emph{Metric Space Theory}.
\index{metric space}%
  Prove the conjecture that
\[  \Aut(\mathbb{Q}^m \mathfrak{U})  <  \Aut(\mathfrak{U}^{\mathbb{Z}^m}) < \Aut(\mathfrak{R}_{m,\omega}) \] 
where $\mathfrak{U}^{\mathbb{Z}^m}$ is the universal Urysohn space
\index{Urysohn space}%
 in which metrics take values in an $m$-coloured lattice $\mathbb{Z}^m$.

\emph{Isometries of $\mathfrak{U}$}.
\index{metric space}%
\index{Urysohn space}%
Any isometry of the universal rational metric space $\mathbb{Q} \mathfrak{U}$ can be extended to an isometry of its completion.  There is an isometry $\sigma$ of $\mathbb{Q} \mathfrak{U}$ permuting all its points in a single cycle (analogous to the cyclic automorphism
\index{automorphism ! cyclic}%
 of the random graph).  The isometry of $\mathfrak{U}$ induced by $\sigma$ has the property that all its orbits are dense.  The question is what other countable groups have this property?  All that is known so far is that the elementary abelian $2$-group
\index{group ! elementary abelian}%
 has this property but the elementary abelian $3$-group does not.

\emph{Abelian structure of $\mathfrak{U}$}.
\index{metric space}%
\index{Urysohn space}%
The closure of $\langle \sigma \rangle$ is an abelian group acting transitively on $\mathfrak{U}$ (so $\mathfrak{U}$ has an abelian group structure).  There are many such $\sigma$, and so the abelian group structure of $\mathfrak{U}$ is not canonical.  What isomorphism types of abelian groups can occur as the closure of $\langle \sigma \rangle$?  The closure of the countable elementary abelian $2$-group with dense orbits is an elementary abelian $2$-group
\index{group ! elementary abelian}%
 acting transitively on $\mathfrak{U}$.

\emph{Subgroups of Isometries of $\mathfrak{U}$}.
It is possible to ``randomly'' construct isometries of $\mathfrak{U}$ in which every orbit is dense~\cite{camver}.  The closure $\overline{\Iso(\mathfrak{U})}$ of the group generated by such an isometry is abelian and transitive, and so gives the Urysohn space an abelian group structure.  Can we make the choice such that $\overline{\Iso(\mathfrak{U})}$ has no non-trivial proper closed subgroup?

\emph{Ordered Metric Space Theory}.  Let $\Lambda$~\label{ordAbgp} be an ordered abelian group.
\index{group ! abelian}%
\index{metric space}%
\index{Urysohn space}%
  Is there a $\Lambda$-valued Urysohn space?

\item \emph{Two Galois Correspondences}.
\index{Galois ! correspondence}%
  A pair $({}^{\triangleright}, {}^{\triangleleft})$ of maps ${}^{\triangleright} : P \to Q$ and ${}^{\triangleleft} : Q \to P$ between two ordered sets $P$ and $Q$ is a \emph{Galois Connection}
\index{Galois ! connection}%
 between $P$ and $Q$ if, $\forall p \in P, q \in Q$, 
 \[ p^{\triangleright} \leq q \Leftrightarrow p \leq q^{\triangleleft}. \]

Whilst the two sets are usually defined to be posets, they can also be taken to be preorders,
\index{preorder}%
(that is reflexive and transitive binary relations) on a set~\cite[p.~173]{dav}, as follows.  Let $\leq_1$ and $\leq_2$ be two preorders on a set $P$ such that $x \leq_1 y$ and $x \leq_2 y$ imply $x = y$.  For $Y \subseteq P$ define,
 \[ Y^{\triangleright} := \{ x \in P : (\forall y) (y \leq_1 x \Rightarrow y \notin Y) \}, \]
 \[ Y^{\triangleleft} := \{ x \in P : (\forall y) (y \leq_2 x \Rightarrow y \notin Y) \}. \]
Let $Q \subseteq P$; $Q$ is a \emph{down-set} of $P$ if whenever $x \in Q$, $y \in P$ and $y \leq x$, we have $y \in Q$.  The pair $({}^{\triangleright}, {}^{\triangleleft})$ establishes a Galois connection between the $\leq_1$-down-sets of $P$ ordered by $\subseteq$ and the $\leq_2$-down-sets of $P$ ordered by $\supseteq$.

The reason that we need to consider preorders, is that the relation $\Gamma_1 \to \Gamma_2$, where $\to$ means `is homomorphic to', is a preorder on the class of graphs, but not a poset.  The relation $\to$ is neither symmetric nor anti-symmetric~\cite{hahntardif}.

Is there a possible application to the theory of monoids?

Consider finite transformation monoids on a fixed set of points with the graphs that can be defined on them, as we saw in a previous section.  

Use $\subseteq$ to denote inclusion as submonoids and $\leq$ to denote subgraphs in the sense of homomorphisms.  We know that $M_1 \subseteq M_2 \Rightarrow \Gr(M_2) \leq \Gr(M_1)$, because $v \sim w$ if $\nexists f \in M_2$ such that $f(v) = f(w)$; but  $\nexists f \in M_2  \Rightarrow \nexists f \in M_1$. 

Given that the symmetry relation $\Gamma_1 \leq \Gamma_2$ and $\Gamma_2 \leq \Gamma_1$ is only well-defined if the two graphs are homomorphic to each other, let $P$ be the set of all transformation monoids on a fixed set of points.  Allow $P$ to be endowed with preorders
\index{preorder}%
 $\subseteq$ and $\leq$, such that $M_1 \subseteq M_2 \Rightarrow \Gr(M_2) \leq \Gr(M_1)$, and $\Gr(M_2) \leq \Gr(M_1)$ if $\Gr(M_2)$ is a spanning
\index{graph ! spanning}%
  subgraph of $\Gr(M_1)$.  Then $M_1 \subseteq M_2 \Rightarrow \Gr(M_1) \subseteq \Gr(M_2)$.  So $\Gr(M_1) = \Gr(M_2)$.  However this does not imply $M_1 = M_2$; if for example $M_1, M_2$ are both permutation groups then $\Gr(M_1) = \Gr(M_2)$ are both null graphs.

Next let $\Gr(M_1) \subseteq \Gr(M_2)$ if $\Gr(M_1)$ is contained as a spanning subgraph of $\Gr(M_2)$ and $\Gr(M_1) \leq \Gr(M_2)$ if $\End(\Gr(M_1)) \subseteq \End(\Gr(M_2))$.  In the case that $\Gr(M_1)$ is a 3-path and $\Gr(M_2)$ is a 4-cycle, both $\Gr(M_1) \subseteq \Gr(M_2)$ and $\Gr(M_1) \leq \Gr(M_2)$ are true.

Can we generalise the construction mentioned at the beginning so as to ensure that the equality $\Gr(M_1) = \Gr(M_2)$ derived above is sufficient to get a Galois correspondence or the pair $({}^{\triangleright}, {}^{\triangleleft})$ sets up a Galois connection?

\smallskip

In the context of this question, we mention the Galois correspondence
\index{Galois ! correspondence}%
 found by Sauer and Stone~\cite{sauerst}
\index{Sauer, N.}%
\index{Stone, M.}%
 between semigroups
\index{semigroup}%
  of maps and (universal) algebras: a semigroup of endomorphisms is associated to each algebra over a set; to each semigroup of functions associate the algebra of finitary operations admitting the maps as homomorphisms. 
 
\item  \emph{Connections with Other Fields: K-Theory}.  Encouraged by the set, number and probability-theoretic constructions
of $\mathfrak{R}$ in~\cite{cameron}, we attempt to find other examples in
mathematics.  Whenever a universal property arises, there is a chance
to find a random edge-creating condition that gives $\mathfrak{R}$.  

Let $n \geqq 0$ be an integer, $S_{n}$ be the category of free abelian
groups
\index{group ! free abelian}%
$E_{n}$ of rank $n$ ($\cong \mathbb{Z}^{n}$), together with a
symmetric bilinear form  $E_{n}\ \times\ E_{n} \to \mathbb{Z} : (x,\ y)
\mapsto x \cdot y$, such that $E_{n}\ \cong\ Hom(E_{n},\
\mathbb{Z})$.  Write $E_{n}\ \cong\ E_{n}^{'}$ for two isomorphic
objects in $S_{n}$, and let $S = \bigcup_{n = 0} S_{n}$, $0 \le n
\le \infty$.  We say that $E, E^{'} \in S$ are \emph{stably\
  isomorphic} if there exists $F \in S$ such that $E \oplus F \cong
E^{'} \oplus F$.  This is an equivalence relation, defining an
additive composition law on $K_{+}(S) = S / \oplus$~\cite{serre}.

The Grothendieck group
\index{group ! Grothendieck}%
 $K(S)$ of $S$ with respect to $\oplus$ follows
from the semi-group $K_{+}(S)$, as  $\mathbb{Z}$ does from
$\mathbb{Z}^{+}$.  The universal property of $K(S)$ arises because
every product homomorphism $f_{s} \cdot f$, where $f_{s} : S \to K(S)$
and $f : K(S) \to A$ is an additive function on $S$, and $A$ is a
commutative group.

Define the Cayley graph
\index{graph ! Cayley}%
on $K(S)$ by $\Cay(K(S), \mathcal{S}^{-} \cup
\mathcal{S}^{+})$, where $\mathcal{S}^{\pm} \subset K_{\pm}(S)$.
One possible rule for edge-formation is to choose at random a finite subset $D = \{d_1, \ldots, d_k \} \subseteq \mathbb{N}$ which
are to be the dimensions of the elements in $K(S)$.  Then the $(*)$-condition to form an edge is for all $E, E^{'} \in K(S)$, join $E\
to\ E^{'}$ if and only if there exists an $F \in \mathcal{S}^{+}$
such that $\dim(E) - \dim(E^{'}) = \dim(F)$.

It is possible that this construction may be developed to forge links
between random graph theory and \emph{K-Theory}.
\index{K-Theory}%

\item  \emph{Algebraic Theory of Switching Groups}.  Inspection of the finite presentations of the symmetric and braid groups given below in Appendix D shows that there is a homomorphism $B_n \to \Sym(n) : b_i \mapsto s_i$.  A large part of the finite-dimensional representation theory of the symmetric group, using the ideas of partitions and Young diagrams, generalizes to representations of $B_n$.  For example, from the early work of Burau~\cite{burau} \index{Burau, W.}%
the reducible representation of $\Sym(n)$ over $\mathbb{Z}$ mapping $s_i \mapsto I_{i-1} \oplus   \begin{pmatrix} 
0&1\\
1&0
\end{pmatrix} \oplus I_{n-i-1}$ lifts to a reducible representation of $B_n$ over the ring of Laurent polynomials $\mathbb{Z}[t, t^{-1}]$ mapping $b_i \mapsto I_{i-1} \oplus   \begin{pmatrix} 
1-t&t\\
1&0
\end{pmatrix} \oplus I_{n-i-1}$.  When $t = 1$, the $B_n$ representation collapses to a $\Sym(n)$ one.

With this as motivation, we offer some speculative questions, which we divide up into four topics.  

 \emph{Algebra Homomorphisms}.  The group algebra of $\Sym(n)$ deforms to a Hecke algebra $H_n$,
\index{Hecke algebra}%
 which is a quotient of the group algebra of the braid group $B_n$.
 \index{group ! braid}%
   This is part of a sequence of algebra homomorphisms:
\begin{center}
$\mathbb{C} B_n \to W_n \to H_n \to A_n \to R_n \to \Sym(n)$
\end{center}
where $W_n$ is the Birman-Wenzl algebra, $A_n$ the Jones algebra and $R_n$ is generated by Burau matrices arising from the $(n-1)$-dimensional representation $\rho_n: B_n \to M_n(\mathbb{Z}[t, t^{-1}])$~\cite{birman}, and $M_n$ are the class of $n \times n$ matrices .  Each algebra supports a Markov trace (being a map $B_n \to M_n(\mathbb{Z}[t, t^{-1}])$ to a certain Laurent polynomial) and thereby yields a link-type invariant such as the Alexander, Jones or Kauffman polynomials.  A priori it seems that linking braidings, which require no two strands to pass through one another, with switchings, which have no such geometrical picture is unworkable.  However given the homomorphisms between direct products of braid groups and switching groups, is there a switching group algebra that would fit into the above sequence?  (We note that the homomorphism goes via a product of symmetric groups, making the braids, in all probability, irrelevant).

\emph{Switching Group Invariants}.  Most link invariants have been interpreted as Markov traces on the disjoint union $B_{\infty} = \bigsqcup_{k \ge 1} B_k$ of braid groups.
 \index{group ! braid}%
  A Markov trace is a ring-valued function $f: B_{\infty} \to R$, which is a class invariant in each $B_k$ and satisfies $f(\beta) = f(\mu^{\pm}(\beta))\ \forall \beta \in B_{\infty}$, where $\mu^{\pm}: B_k \to B_{k+1}: \beta \mapsto \beta b_k^{\pm}$.  For example the homomorphism $\theta : B_n \to \Sym(n)$ is defined by the sequence $1 \to P_n \to B_n \to \Sym(n) \to 1$ where in the \emph{coloured braid group} $P_n$, group multiplication preserves a well-defined assignment of colours to each of the $n$ strands.  This yields as a Markov trace the set $f(\beta) =$ the number of cycles in $\theta(\beta)$, and so the link invariant is the number of components.  Can such a Markov trace be produced for the switching group union $S_{\infty} = \bigsqcup_{n \ge 1} S_{m,n}$?

\emph{Braid Groups and Yang-Baxter Equation}.
 \index{group ! braid}%
 \index{Yang-Baxter equation}%
  There are a series of linear representations of braid groups $B_n \to \End(W^{\otimes n})$ associated to any finite-dimensional complex simple Lie algebra $\mathsf{L}$ and its finite-dimensional irreducible representations $\rho : \mathsf{L} \to \End(W)$.  Letting $\{I_{\mu}\}$ be an orthonormal basis of $\mathsf{L}$ with respect to the Cartan-Killing
\index{Cartan, \'E.}%
form enables the definition $\Omega_{\alpha \beta} = \Sigma_{\mu} \rho_{\alpha} (I_{\mu}) \otimes \rho_{\beta} (I_{\mu})$ where $\rho_{\alpha}$ is the non-trivial represention only on the $\alpha$-th factor.  Rational solutions of the classical Yang-Baxter equation lead to the \emph{infinitesimal pure braid relations}~\cite[pp.~730, 768]{kohno} \cite[p.~78]{kohno1}:
\begin{center}
$[\Omega_{\alpha \beta},  \Omega_{\alpha \gamma} + \Omega_{\beta \gamma}] =  [\Omega_{\alpha \beta} + \Omega_{\alpha \gamma}, \Omega_{\beta \gamma}] = 0\ \text{for}\ \alpha < \beta < \gamma$
\end{center}
\begin{center}
$[\Omega_{\alpha \beta}, \Omega_{\gamma \delta}] = 0\ \text{for distinct}\ \alpha, \beta, \gamma, \delta.$
\end{center}
As observed in Theorem~\ref{strg}, the commutator $[\sigma_{c,d,\{x\}}, \sigma_{d,e,\{y\}}]$ of two switchings permutes the ordered triple of colours $(c, d, e)$ to $(d, e, c)$.  With this interpretation the switchings satisfy the above relations.  Is there a genuine connection here?

\emph{Representations of Switching Groups}.
 \index{group ! switching}%
  Motivated by V. Jones' Hecke algebra representation theory for the braid group \cite{jonesvfr},
     \index{group ! braid}%
can we find a Hecke algebra
     \index{Hecke algebra}%
representation theory for any of the switching groups?  Given the  structural resemblance of the switching and the hyperoctahedral groups,
       \index{group ! hyperoctahedral}%
is there a geometrical interpretation for the switching groups that is somehow dual to that for the hyperoctahedral groups?

\end{enumerate}

\appendix

\chapter{Prerequisite Background}

In this appendix we establish some of the basic notation, and collect objects of study and theory on which the results of this monograph are built.  

\section{Permutation Groups}
\label{PermutationGroups}
A \emph{permutation group} $G$
\index{group ! permutation}%
on a set $\Omega$ is a subgroup of the symmetric group $\Sym(\Omega)$ of all permutations of $\Omega$.  More generally, a \emph{permutation representation}
\index{permutation representation}%
of $G$ on $\Omega$ is a homomorphism $G \to \Sym(\Omega)$.  Equivalently we can say that $G$ \emph{acts on} $\Omega$, or \emph{induces} a permutation group on $\Omega$, or that $\Omega$ is a \emph{$G$-space}.
\index{Gs@$G$-space}%
Most interesting properties of a group action are properties of its image; so the theories of group actions and permutation groups are almost identical.  Two $G$-actions on sets $\Omega_1$ and $\Omega_2$ are isomorphic if there is a bijection between the two domains such that the permutations induced by any element of $G$ correspond under the bijection.  The \emph{degree}
\index{group ! permutation ! degree}%
of a permutation group is the cardinality of the set $\Omega$.

Throughout this work we will assume the convention that groups act on the right. 

If there exists $g \in G$ such that $x g = y$ then $x \sim y$ is an equivalence relation on $\Omega$ whose equivalence classes are called \emph{orbits},
\index{group ! permutation ! orbit}%
and $G$ is \emph{transitive}
\index{group ! permutation ! transitive}%
if it has one orbit.  Every $G$-space can be uniquely expressed as a disjoint union of transitive $G$-spaces.  Any transitive action of $G$ is isomorphic to the action on the set of right cosets of a subgroup.  The actions on the sets of right cosets of subgroups $H$ and $K$ are isomorphic if and only if the subgroups $H$ and $K$ are conjugate.  Assume that $k \le |\Omega|$.  A group is \emph{$k$-transitive}
\index{group ! permutation ! $k$-transitive}%
if it is transitive on the set of ordered $k$-subsets of $\Omega$ and it is \emph{$k$-homogeneous}
\index{group ! permutation ! $k$-homogeneous}%
if it is transitive on the set of unordered $k$-subsets of $\Omega$.  A group is \emph{highly transitive}
\index{group ! permutation ! highly transitive}%
(respectively \emph{highly homogeneous})
\index{group ! permutation ! highly homogeneous}%
if it is $k$-transitive (respectively $k$-homogeneous) for every positive integer $k$.  If a group $G \le \Aut(\Gamma)$ acting on a graph $\Gamma$ is $2$-transitive
\index{group ! permutation ! $2$-transitive}%
 then the full automorphism group
 \index{group ! automorphism}%
  of the graph is the symmetric group on the vertex set, and the graph is either complete or null.  A theorem of Cameron~\cite{camaa}
\index{Cameron, P. J.}%
states that any highly homogeneous but not highly transitive permutation group is dense in the automorphism group
\index{group ! automorphism}%
 of one of the following: (i) a dense linear order;
\index{linear order}%
 the betweenness relation induced from a dense linear order; a dense circular order; the quaternary separation relation induced from a sense circular order.

The \emph{setwise stabilizer}
\index{group ! permutation ! setwise stabilizer}%
$G_{\{\Delta\}}$ of a subset $\Delta \in \Omega$ is the group whose elements are $\{g \in G : \Delta g = \Delta\}$.  More generally the \emph{pointwise stabilizer}
\index{group ! permutation ! pointwise stabilizer}%
$G_{(\Delta)}$ is the set of permutations that fix every point of $\Delta$.  In terms of subgroups, $G_{(\Delta)} \le G_{\{\Delta\}} \le G$.  If $\Delta = \{ x \}$, then its setwise and pointwise stabilizers coincide, and are denoted by $G_x$.

Let $H$ be a subgroup of the abstract group $G$.  The \emph{coset space}
\index{coset space}%
of $H$ in $G$ is the set  of right cosets of $H$ in $G$, and $G$ acts on it by $(H x) g = Hxg$.  If $G$ acts transitively on $\Omega$ then $\Omega$ is isomorphic to the coset space of $G_x$ in $G$ for a point $x \in \Omega$.  The permutation group $G$ is \emph{semiregular}
\index{group ! permutation ! semiregular}%
if only the identity has a fixed point; thus it is freely acting.  The centralizer of a transitive group is semiregular, and the centralizer of a semiregular group is transitive.  If $G$ is transitive and semiregular
\index{group ! permutation ! semiregular}%
 then $G$ is said to be \emph{regular}, or simply transitive or sharply $1$-transitive.
\index{group ! permutation ! regular}%
(A permutation is regular if all cycles in its canonical cycle decomposition have the same length).  In this case, $\Omega$ is isomorphic to $G$, $G$ acts on $\Omega$ by right multiplication, and this is called the \emph{right regular representation}.
\index{right regular representation}%
Also in this case, $G$ is isomorphic to its permutation group image and is said to act \emph{faithfully}
\index{group ! faithful action}%
on $\Omega$ because the kernel is the identity.  Cayley
\index{Cayley, A.}%
used this action to show that every group is isomorphic to a permutation group.  The centralizer in the symmetric group of the image of the right regular representation is the left regular representation in which $g \in G$ induces the permutation $x \mapsto g^{-1} x$.  The two representations commute with each other:\quad $g^{-1} (xh) = (g^{-1} x) h$.  The \emph{diagonal group}
\index{group ! diagonal}%
$G^* = G \times G$ acts transitively on $G$ as the product of the two regular actions, $g^{-1}xh$, with the left and right actions as normal subgroups, and the \emph{diagonal subgroup} $\{(g, g) : g \in G\}$, acting by conjugation on $G$ and stabilizing the identity.

The Orbit-Stabilizer Theorem
\index{Orbit-Stabilizer Theorem}%
says that if the finite group $G$ acts transitively on a set $\Omega$, and $x \in \Omega$, then $|G_x| = |G| / |\Omega|$.

Let $\fix(g)$ be the number of points of $\Omega$ fixed by $g \in G$.  We shall have cause to use the following theorem,

\begin{theorem}[Orbit-Counting Lemma]
\index{Orbit-Counting Lemma}%
Let $G$ be a permutation group on the finite set $\Omega$.  Then the number of orbits of $G$ on $\Omega$ is equal to the average number of fixed points of an element of $G$, that is,
\[ \frac{1}{|G|} \sum_{g \in G} \fix(g). \]
\end{theorem}

Recall some theory~\cite{cam1} of groups with regular normal subgroups.
\index{group ! regular normal subgroup}%
If $N$ is any group, and $H$ a subgroup of its automorphism group
\index{group ! automorphism}%
 $\Aut(N)$, then the semidirect product of $N$ by $H$ acts as a permutation group
on $N$, with $N$ as regular normal subgroup and $H$ as the stabilizer
of the identity: the element $hn$ of $N \sd H$ acts on $N$ as $x^h
n$ for $x, n \in N, h \in H$.  Conversely, if $G$ is any permutation
group on a set $\Omega$ with a regular normal subgroup $N$, then there
is a bijection from $\Omega$ to $N$, which takes the given action of
$N$ on $\Omega$ to its action on itself by right multiplication, and the given action of the subgroup $H = G_1$, the stabilizer of $1$ to its action by conjugation.  The
action of the point stabilizer on $\Omega$ is isomorphic
to its action on $N$ by conjugation.  The group $N
\sd \Aut(N)$ is called the \emph{holomorph}
\index{holomorph}%
of $N$, denoted $HOL(N)$, and $G \leq HOL(N)$.  The holomorph of $N$ is the normalizer of its right regular representation in $\Sym(N)$.

Suppose that $G$ is transitive on $\Omega$.  A \emph{congruence}
\index{congruence}%
is a $G$-invariant equivalence relation on $\Omega$, or thought of as a set of ordered pairs is a union of orbits of $G$ in its action on $\Omega \times \Omega$.  A group is \emph{primitive}
\index{group ! permutation ! primitive}%
if its only congruences are the trivial ones of equality and the universal relation with one equivalence class.  Otherwise the group is said to be \emph{imprimitive}.
\index{group ! permutation ! imprimitive}%
Finite primitive groups are small, and they are rare.  Examples of infinite primitive groups are the automorphism group
\index{group ! automorphism}%
 of the random graph
\index{graph ! random}%
  and of the Henson graphs
\index{graph ! Henson}%
   (see Chapters~\ref{chap2} and~\ref{chapFD}).  If $G$ is transitive but imprimitive, and $E$ is a non-trivial $G$-congruence on $\Omega$ then the set of $E$-classes is called a \emph{system of imprimitivity}
\index{system of imprimitivity}%
and its elements are called \emph{blocks of imprimitivity}.
\index{blocks of imprimitivity}%
\smallskip

A subgroup $H < G$ is \emph{maximal}
\index{group ! maximal}%
in $G$ if $\forall x \in G \backslash H$, $G := \langle H, x \rangle$. A transitive group $G$ is primitive if and only if $G_x$ is a maximal subgroup of $G$.  More generally, the congruences
\index{congruence}%
 form a lattice isomorphic to the lattice of subgroups lying between $G_x$ and $G$.  Finite transitive groups can be built from primitive groups,
\index{group ! permutation ! primitive}%
  which are classified by the O'Nan--Scott Theorem~\cite{liebprsa1}~\cite{scottll}.
\index{O'Nan--Scott Theorem}%
\index{O'Nan, M.}%
\index{Scott, L.}%
An analogous result is true for some but not all infinite permutation groups; it is true for example for oligomorphic groups because they only have finitely many congruences.  Examples of maximal subgroups of the automorphism group
\index{group ! automorphism}%
 of the random graph include stabilizers of finite sets of vertices, unordered edges and unordered non-edges.  Some properties of primitivity are given in

\begin{theorem}
\begin{itemize}
\item[(a)]  A $2$-transitive group is primitive.
\index{group ! permutation ! primitive}%
\index{group ! permutation ! $2$-transitive}%
\item[(b)]  A non-trivial normal subgroup of a primitive group is transitive.
\end{itemize}
\end{theorem}

The \emph{socle} $\soc(G)$
\index{group ! socle}%
of a finite group $G$ is the product of its minimal normal subgroups.  It is itself normal in $G$ and is a direct product of finite simple groups.
\index{group ! simple}%
If $G$ is primitive then either $\soc(G)$ is a product of isomorphic finite simple groups in a known permutation action, or $N = \soc(G)$ is simple and $G \le \Aut(N)$.

\medskip

A group $A$ is a \emph{B-group}
\index{group ! B-group}%
if any primitive permutation group
\index{group ! permutation ! primitive}%
 $G$ which contains $A$ acting regularly is $2$-transitive;
\index{group ! permutation ! $2$-transitive}%
 that is, if any overgroup of $A$ which kills all the $A$-invariant equivalence relations necessarily kills all the non-trivial $A$-invariant binary relations.  The letter B stands for Burnside,
\index{Burnside, W.}%
who showed that a cyclic group of prime-power but not prime order is a B-group.  The proof contained a gap which was subsequently fixed by Schur, who invented and developed Schur rings for this purpose.  The theory of Schur rings (or S-rings)
\index{Schur ring (S-ring)}%
is connected with many topics in representation theory, quasigroups,
\index{quasigroup}%
 association schemes, and other areas of mathematics; historically, it was an important source of ideas in these subjects.  The theory of S-rings and its connection with representation theory is described in Wielandt's
\index{Wielandt, H.}%
book~\cite{wielandt}.

Primitive groups
\index{group ! permutation ! primitive}%
 are comparatively rare.  For example, the set of numbers for which there exists a primitive group of degree $n$ other than $\Sym(n)$ and $\Alt(n)$ has density zero~\cite{camneutea}, and so the set of orders of non B-groups has density zero.  The precise result~\cite{camneutea} of Cameron, Neumann, Teague
\index{Cameron, P. J.}%
\index{Neumann, $\Pi$. M.}%
\index{Teague, D. N.}%
is that if $N(x)$ is the number of $n \leq x$ where $\exists H < G \neq \Alt(n), \Sym(n)$ with $H$ regular and $G$ primitive, then $\frac{N(x)}{x} \to 1$ as $x \to \infty$.

There are no known infinite B-groups.  A \emph{square-root set}
\index{square-root set}%
in a group $X$ is a set of the form
\[ \sqrt{a} = \{ x \in X : x^2 = a \}. \]
If $ a \neq 1$, it is called \emph{non-principal}.
\index{square-root set ! non-principal}%
One of the most powerful nonexistence theorems is the following result~\cite{camjohn},

\begin{theorem}
\label{bgpthm}

Let $A$ be a countable group with the following property:

\centerline{$A$ cannot be written as the union of finitely many translates of}
\centerline{non-principal square-root sets together with a finite set.}
Then $A$ is not a B-group.  More precisely, there exists a primitive
\index{group ! permutation ! primitive}%
 but not $2$-transitive group $G$
\index{group ! permutation ! $2$-transitive}%
 which contains a regular subgroup isomorphic to each countable group satisfying this condition.
\end{theorem}

The group $G$ in Theorem~\ref{bgpthm} includes the automorphism group of the ($2$-coloured) random graph.  The condition of this theorem is not very restrictive: any countable abelian group of infinite exponent satisfies the condition; and for any finite or countable group $A$, the direct product of $A$ with an infinite cyclic group satisfies it, so $G$ embeds every countable group as a semiregular subgroup.
\index{group ! permutation ! semiregular}%
  Automorphism groups
\index{group ! automorphism}%
   of multicoloured random graphs
\index{graph ! random ! $m$-coloured}%
    (see section below) are primitive but not $2$-transitive.
\index{group ! permutation ! $2$-transitive}%
\index{group ! permutation ! primitive}%
  So if a group embeds
  as a regular subgroup of such an automorphism group then it cannot
  be a B-group.  
  
\smallskip

The \emph{normal core}
\index{group ! normal core}%
 of a subgroup $H$ of a group $G$ is the largest normal subgroup of $G$ that is contained in $H$, or equivalently it is the intersection of the conjugates of H. More generally, the core of $H$ with respect to a subset $S \subseteq G$ is
 \[ \Core_S (H) := \bigcap_{s \in S} s^{-1} H s.\]
So the normal core is the core with respect to $S=G$. The normal core of a pointwise stabilizer acts as the identity on its entire orbit, and in a transitive action is precisely the kernel of the action.

A \emph{core-free}
\index{group ! core-free}%
subgroup is one whose normal core is the trivial subgroup, or equivalently, one that occurs as the stabilizer subgroup of a transitive, faithful group action.
\index{group ! action}%

\medskip

We next consider a concept that is weaker than primitivity, giving some of the interesting operand graphs; the groups and graphs are intended to be finite.  There are natural connections between \emph{quasiprimitive} permutation groups,
\index{group ! permutation ! quasiprimitive}%
all of whose non-trivial normal subgroups are transitive,~\cite{praeger}~\cite{praeger1} and arc-transitive, Cayley and bipartite graphs.
\index{graph ! arc-transitive}%
\index{graph ! bipartite}%
\index{graph ! Cayley}%
 A partition $\pi$ of a set $\Omega$ with $G \le \Sym(\Omega)$ being transitive, is called \emph{$G$-invariant}
\index{partition ! $G$-invariant}%
if for all $g \in G$ and $p \in \pi$~\label{partp}, we have that $p^g = \{x \in p\}$ again belongs to $\pi$.  To strengthen this, call a partition \emph{$G$-normal} if there exists $K \lhd G$ such that $\pi$ is the set of $K$-orbits in $\Omega$.  A primitive permutation group
\index{group ! permutation ! primitive}%
 has only trivial $G$-invariant partitions (singleton subsets and the universal partition $\{\Omega\}$),  but $G$ is quasiprimitive if and only if the only $G$-normal partitions are the trivial ones.  If $\pi$ has \emph{maximal blocks}, that is, if the only $G$-invariant partitions with $\pi$ as refinement are the universal partition $\{\Omega\}$ and $\pi$ itself, then $G$ induces a primitive permutation group
\index{group ! permutation ! primitive}%
 $G^{\pi} \cong G$ on $\pi$.  An overgroup of a quasiprimitive group
\index{group ! permutation ! quasiprimitive}%
  need not be quasiprimitive, unlike the equivalent property for primitive groups; if $G$ is imprimitive quasiprimitive and $\pi$ a non-trivial $G$-invariant partition then the stabilizer of $\pi$ in $\Sym(\Omega)$ has $\prod_{p \in \pi} \Sym(p)$  as a non-trivial intransitive normal subgroup.

An \emph{$s$-arc}
\index{arc}%
in a graph $\Gamma$ is an $(s+1)$-tuple $(v_0, v_1, \ldots, v_s)$ of vertices in $\Gamma$ such that $v_i \sim v_{i-1}$ $(1 \leq i \leq s)$ and $v_{j-1} \neq v_{j+1}$ $(1 \leq j \leq s-1)$.  Given $G \leq \Aut(\Gamma)$ we call $\Gamma$ \emph{locally $(G, s)$-arc transitive}
\index{graph ! locally $(G, s)$-arc transitive}%
if $\Gamma$ contains an $s$-arc and given any two $s$-arcs $\alpha, \beta$ starting at the same vertex $v$, there exists an element $g \in G_v : \alpha \to \beta$.  Call $\Gamma$ \emph{locally $s$-arc transitive}
\index{graph ! locally $s$-arc transitive}%
if it is locally $(G, s)$-arc transitive for some $G \leq \Aut(\Gamma)$.  A graph is \emph{$G$-locally primitive}
\index{graph ! gl@$G$-locally primitive}%
if the stabilizer $G_v$ of every vertex $v$ acts primitively on the set $\Gamma(v)$ of neighbours of $v$, and \emph{locally primitive}
\index{graph ! locally primitive}%
if it is $G$-locally primitive for some automorphism group
\index{group ! automorphism}%
 $G$.  If $\Gamma$ has valency
\index{graph ! valency}%
at least two then it is locally $(G, 2)$-arc transitive if and only if for every vertex $v$, $G_v$ acts $2$-transitively on the set $\Gamma(v)$.
\index{group ! permutation ! $2$-transitive}%
 So locally $(G, 2)$-arc transitive graphs are $G$-locally primitive.  Finite graphs of vertex valency at least two, admitting vertex-intransitive automorphism groups
\index{group ! automorphism}%
  $G$ are either locally $(G, s)$-arc transitive for $s \ge 2$ or $G$-locally primitive, are bipartite
\index{graph ! bipartite}%
with the two parts of the bipartition being the orbits of $G$, and $G$ being edge-transitive~\cite{giudici}.  Much of the research on locally $(G, s)$-arc transitive graphs has centered on the possible point stabilizers of two adjacent vertices.

If $\pi$ is a partition of $\Omega$, the \emph{quotient graph}
\index{graph ! quotient}%
of graph $\Gamma = (\Omega, E)$ relative to $\pi$ is defined to be the graph $\Gamma_{\pi} = (\pi, E_{\pi})$ where $\{\pi, \pi'\} \in E_{\pi}$ if and only if there exists $p \in \pi, p' \in \pi'$ such that $\{p, p'\} \in E$.  If $\pi$ is invariant under $G \le \Aut(\Gamma)$, then $G^{\pi} \le \Aut(\Gamma_{\pi})$.  If $N \lhd G = \Aut(\Gamma)$ which is intransitive on both parts of a bipartite graph, then taking quotients with respect to the orbits of $N$ preserves both local primitivity and local $s$-arc-transitivity and leads to the study of graphs where $G$ acts faithfully on both orbits and quasiprimitively
\index{group ! permutation ! quasiprimitive}%
 on at least one~\cite{giudici}.  If $N \lhd \Aut(\Gamma)$ acts intransitively on vertices, then the \emph{quotient graph}
\index{graph ! quotient}%
of graph $\Gamma_N$ has the $N$-orbits on $V(\Gamma)$ as the vertex set, and two $N$-orbits $B_1, B_2$ are adjacent in $\Gamma_N$ if and only if $\exists v \in B_1, w \in B_2$ such that $v \sim w$ in $\Gamma$.  A graph $\Gamma$ is a \emph{cover}
\index{graph ! cover}%
of $\Gamma_N$ if $|\Gamma(v) \cap B_2| = 1$ $\forall v \in B_1$ and edge $\{B_1, B_2\} \in \Gamma_N$.  If $\Gamma$ is a bipartite $(G, s)$-arc transitive graph, and $|G : G^{+}| = 2$ where $G^{+}$ fixes both parts of the bipartition then $\Gamma$ is locally $(G^{+}, s)$-arc transitive.  If $\Gamma$ is a nonbipartite $(G, s)$-arc transitive graph, then $\Gamma$ is a cover of $\Gamma_N$ and $\Gamma_N$ is a $(G/N, s)$-arc transitive graph.  So the problem becomes one of finding all examples where $G$ acts quasiprimitively on vertices and finding their covers.  

Finally we distill a construction, useful in the study of quasiprimitive groups, of an edge-transitive graph
\index{graph ! edge-transitive}%
from a given group $G$ and pair of subgroups $L, R$ such that $L \cap R$ is core-free in $G$; a subgroup $H$ of $G$ is \emph{core-free}
\index{group ! core-free}%
if $\cap_{g \in G} H^g = 1$.  Let $\Delta_1$ be the set $[G : L]$ of right cosets of $L$ in $G$, and $\Delta_2$ be the set $[G : R]$ of right cosets of $R$ in $G$.  The bipartite \emph{coset graph}
\index{graph ! bipartite}%
\index{graph ! coset}%
$\Cos(G, L, R)$ has  vertex set the disjoint union $\Delta_1 \sqcup \Delta_2$ such that two vertices $Lx$ and $Ry$ are adjacent if and only if $xy^{-1} \in LR$.  The group $G$ acts on the vertex set of graph by right multiplication, $G$ is edge-transitive, and $L$ and $R$ are the stabilizers of the adjacent vertices $L, R$ respectively.  
\begin{proposition}
Let $\Gamma = \Cos(G, L, R)$ for some group $G$ with subgroups $L, R$ such that $L \cap R$ is core-free, and $\Delta_1 = [G : L]$ and $\Delta_2 = [G : R]$. Then 
\begin{itemize}
\item[(1)]  $\Gamma$ is connected if and only if $\langle L, R \rangle = G$;
\item[(2)]  $G \le \Aut(\Gamma)$, $\Gamma$ is $G$-edge-transitive, and $\Delta_1$ and $\Delta_2$ are $G$-orbits on vertices;
\item[(3)]  $G$ acts faithfully on $\Delta_1$ and $\Delta_2$ if and only if both $L$ and $R$ are core-free,
\item[(4)]  $\Gamma$ is locally $(G, 2)$-arc transitive if and only if $L$ is $2$-transitive
\index{group ! permutation ! $2$-transitive}%
 on $[L : L \cap R]$ and $R$ is $2$-transitive on $[R : L \cap R]$. 
\end{itemize}
Conversely, if $\Gamma$ is a $G$-edge transitive but not $G$-vertex transitive graph, and $v \sim w$ then $\Gamma \cong \Cos(G, G_v, G_w)$.
\end{proposition}

In spite of the fact that the groups $\Aut(\mathfrak{R}_{m,\omega})$ $(m \ge 2)$ are simple, we have found overgroups which may potentially be quasiprimitive.
\index{group ! permutation ! quasiprimitive}%
  So our intention in giving the above concepts for the finite theory, has been to stimulate research on infinite quasiprimitive groups and their relations to infinite graphs.

We mention one other generalisationsof primitivity: a transitive permutation group $G$ of finite degree is \emph{semiprimitive}
\index{group ! permutation ! semiprimitive}%
 if it is not regular and if every normal subgroup of $G$ is transitive or semiregular~\cite{bermaroti}.
\index{group ! permutation ! semiregular}%
 
\smallskip

Let $S$ be a subset of a group $G$.  The \emph{Cayley graph}
\index{graph ! Cayley}%
$\Cay(G, S)$ is the directed graph with vertex set $G$, having directed edges $(g, sg)$ for all $g \in G$ and $s \in S$.  If $1 \notin S$, then the graph has no loops; if $s \in S \Rightarrow s^{-1} \in S$, then it is an undirected graph (that is, whenever $(g, h)$ is an edge, so is $(h, g)$, and we can regard edges as unordered pairs).  The graph $\Cay(G, S)$ is connected if and only if $S$ generates $G$.  The regular action of $G$ on the vertex set of $\Cay(G, S)$ by right multiplication is an automorphism group
\index{group ! automorphism}%
 of $\Cay(G, S)$.  Conversely when a graph $\Gamma$ admits an automorphism group $G$ acting regularly on the vertices, then $\Gamma$ is isomorphic to a Cayley graph for $G$.  (Choose a point $\alpha \in \Omega$, and take $S$ to be the set of elements $s$ for which $(\alpha, \alpha s)$ is an edge).

\begin{theorem}[Sabidussi's Theorem]~\cite{sabidussi}
\index{Sabidussi's Theorem}%
\index{Sabisussi, G.}%
A graph $\Gamma$ is a Cayley graph of a group $G$ if and only if there is a transitive and free action of $G$ on $\Gamma$.
\end{theorem}

Whilst not all vertex-transitive graphs are Cayley graphs, there is the following conjecture of McKay-Praeger 

\begin{conjecture} If
\[b_n := \frac{\sharp\ \text{ Cayley\ graphs\ on\ 
$\leq n$\ vertices }}{\sharp\ \text{Vertex-transitive\ graphs\ on\ 
$\leq n$\ vertices }}\]
then $b_n \to 1$ as $n \to \infty$. 
\end{conjecture}

A \emph{primitive Cayley graph}
\index{graph ! primitive Cayley}%
 is a graph $\Cay(G, S)$, whose automorphism group
 \index{group ! automorphism}%
  $\Aut(Cay(G, S))$ is vertex-primitive.  For example, if $S = G 
\backslash \{1\}$ then $\Cay(G, S)$ is the complete graph $K_n$, where $n = |G|$ and $\Aut (\Gamma) = \Sym (G) \cong Sym(n)$, and hence primitive.

Kr\"{o}n and M\"{o}ller
\index{kron@Kr\"{o}n, B.}%
\index{moller@M\"{o}ller, R.}%
have studied~\cite{kron} a generalization of Cayley graphs called \emph{rough Cayley graphs}
\index{graph ! rough Cayley}%
defined for compactly generated totally disconnected locally compact topological groups.
\index{group ! locally compact}%
\index{group ! totally disconnected}%
A connected graph $\Gamma$ is a rough Cayley graph of a topological group $G$ if $G$ acts transitively on $\Gamma$ and the vertex-stabilizers $U = G_x$ are compact open subgroups of $G$.  The vertex set is $G / U$ and two distinct left cosets $xU$ and $yU$ are adjacent if there are elements $g \in xU$ and $h \in yU$ such that $g$ and $h$ are adjacent in $G / U$. The $G$-action on $G / U$ induces an action of $G$ on $\Gamma$.  Conversely if $S = U \cup \{g_i(x) | i \in I\}$ where $U=G_x$ and $\{g_i(x) | i \in I\}$ are the neighbours of $x$ for $\{g_i(x)\} \in G$ generates an ordinary Cayley graph $\Gamma$, then the quotient graph (defined below) of $\Gamma$ with respect to $U$ gives a rough Cayley graph.  When a group is finitely generated, various of its properties are independent of the choice of finite generating set used to construct the Cayley graph; analogously a rough Cayley graph, generated by the union of a compact open subgroup
\index{group ! compact}%
\index{group ! open}%
and a finite set, is a quasi-isometric invariant of the group.
\index{quasi-isometry}%
Compact open subgroups are commensurable with each other and with their conjugates.

\smallskip

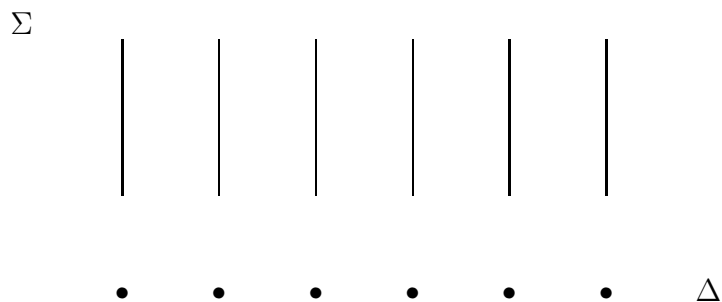
\begin{figure}[!h]
$$\xymatrix{ 
{\Sigma} & {}  & {} & {}  & {} & {}  & {}\\
{} & {}  & {} & {}  & {} & {}  & {}\\
& {} \ar@{-}[uu]  & {} \ar@{-}[uu] & {}
  \ar@{-}[uu] & {} \ar@{-}[uu] & {} \ar@{-}[uu] & {}\ar@{-}[uu]\\
{} & {\bullet}  & {\bullet} & {\bullet}  & {\bullet} & {\bullet}  & {\bullet} & {\Delta}
}$$\caption{Fibre diagram for defining wreath product}
\label{wrprod}
\end{figure}

Let $H$ and $K$ be permutation groups acting respectively on sets $\Sigma$ and $\Delta$.  Take $\Omega = \Sigma \times \Delta$ where fibres bijective with $\Sigma$ form a covering space of $\Delta$, as in Figure~\ref{wrprod}; this affords the so-called \emph{imprimitive action} of the wreath product.
\index{group ! permutation ! wreath product ! imprimitive action}%
The \emph{base group} $B$
\index{group ! base}%
is the cartesian product of $|\Delta|$ copies of $H$, that is one for each fibre of $\Omega$ (or each element of $\Delta$).  The \emph{top group} $K_1$
\index{group ! top}%
is the permutation group on $\Omega$ obtained by letting $K$ permute the fibres according to its given action on $\Delta$.  The \emph{wreath product}
\index{group ! permutation ! wreath product}%
$H \Wr K$ is the semi-direct product of $B$ and $K_1$.

Let $G$ be transitive but imprimitive on $\Omega$, $\Sigma$ be a congruence
 \index{congruence}%
 class, $H$ the permutation group induced on $\Sigma$ by its setwise stabilizer in $G$, $\Delta$ the set of congruence classes, and $K$ the group induced on $\Delta$ by $G$.  Then $G$ can be embedded in a natural way as a permutation group into $H \Wr K$.

The above \emph{full wreath product}
\index{wreath product ! full}%
is to be distinguished from the \emph{restricted wreath product}
\index{wreath product ! restricted}%
denoted $H \text{wr} K$, in which $B$ is a restricted cartesian power of $H$
by $\Delta$, that is the set of functions $\Delta \to H$ has value $1$
at all but a finite number of arguments.

If the congruences are partially ordered then a more general wreath product can be defined~\cite{cam6}.

In the \emph{product action}
\index{group ! permutation ! wreath product ! product action}%
of $H \Wr K$ on the set of functions $\phi : \Delta \to \Sigma$, that is on the set of global sections
\index{global section}%
(or transversals) of the fibre bundle, an element $f$ of the base group acts by $\phi \cdot f(i) = \phi(i) \cdot f(i)$, and the top group acts on the arguments of the functions by $\phi \cdot k(i) = \phi(i k^{-1})$.  If $|\Sigma|, |\Delta| > 1$ then $H \Wr K$ in its product action is primitive
\index{group ! permutation ! primitive}%
 if and only if $H$ is primitive but not regular on $\Sigma$, $\Delta$ is finite, and $K$ is transitive on $\Delta$.  (The only primitive regular finite groups are the cyclic groups $C_p$ of prime order $p$).
\index{group ! cyclic}%

\medskip

The \emph{finitary symmetric group}
\index{group ! symmetric ! finitary}%
$\FSym(\Omega)$ is the group of permutations of an infinite set $\Omega$ that move only finitely many points.  Because such permutations move only finitely many points, they can be assigned a parity, regardless of the number of fixed points.  This allows us to define the \emph{alternating group}
\index{group ! alternating}%
$\Alt(\Omega)$ of all even permutations.  

The Jordan--Wielandt Theorem
\index{Jordan--Wielandt Theorem}%
states that the only primitive
\index{group ! permutation ! primitive}%
 finitary groups on an infinite set $\Omega$ are $\FSym(\Omega)$ and $\Alt(\Omega)$.  We also mention Neumann's
\index{Neumann, $\Pi$. M.}%
extension to a structure theory for transitive finitary groups.

The Theorem of Baer, Schreier and Ulam~\cite{baer}~\cite{schreierulam},
\index{Baer--Schreier--Ulam Theorem}%
states that (using the axiom of choice, a statement of which is given in a later section)
\index{axiom of choice}%
for a set $\Omega$, not necessarily countable, the proper normal subgroups of $\Sym(\Omega)$ are (i) the trivial group; (ii) $\Alt(\Omega)$; (iii) the bounded symmetric group $\BSym_\alpha(\Omega)$ consisting of all permutations moving fewer than $\alpha$ points, for each infinite cardinal $\alpha\le|\Omega|$; (iv) $\Sym(\Omega)$.  In the case of a countable set $\Omega$ the factor groups are those in the normal series
\[ 1 \unlhd \Alt(\Omega) \unlhd \FSym(\Omega) \unlhd \Sym(\Omega) \]
 of simple groups.
 \index{group ! simple}%
A proof of the classification can be found in~\cite{scott}.

\medskip

An \emph{oligomorphic} permutation group
\index{group ! permutation ! oligomorphic}%
is one having finitely many orbits in its induced action on $n$-tuples for all $n$.  This concept has been central in uncovering connections between permutation groups and model theory.
\index{model theory}%
  A \emph{cofinitary} permutation is one that fixes only finitely many points.  A \emph{cofinitary permutation group}
\index{group ! permutation ! cofinitary}%
is one in which every non-identity element is cofinitary.  Primitive
\index{group ! permutation ! primitive}%
 cofinitary groups are too wild to enable a classification.

\medskip

A graph is \emph{locally finite} if all vertex valencies are finite, and \emph{vertex-symmetric}
\index{graph ! vertex-symmetric}%
if its automorphism group
\index{group ! automorphism}%
 is transitive on vertices.  Whilst the random graphs
\index{graph ! random}%
  of this monograph are far from being locally finite, there is a connection with other objects that are of interest.  The group of \emph{bounded automorphisms}
\index{group ! bounded automorphisms}%
(for which $d_{\Gamma} (v, gv) \le c$, for some $c \in \mathbb{N}$), is a normal subgroup of the overlying automorphism group.  For example, for $\mathbb{Z}^2$ the bounded automorphisms are the translations $(x, y) \mapsto (x + a_1, y + a_2)$.  A \emph{system of imprimitivity}
\index{system of imprimitivity}%
$\sigma$ for a group $G \leq \Aut(\Gamma)$ acting on $V(\Gamma)$ is a partition of $V(\Gamma)$ into blocks that are permuted by each automorphism; the automorphism group
\index{group ! automorphism}%
 of $\Gamma / \sigma$ induced by $G$ is denoted $G^{\sigma}$.  The work of Trofimov
\index{Trofimov, V. I.}%
on group actions
\index{group ! action}%
 on locally finite graphs, summarized in~\cite{trof}, includes the following results,

\begin{theorem}
Let $\Gamma$ be a connected locally finite vertex-symmetric graph.  An automorphism $g$ of $\Gamma$ is bounded if and only if there exists a system of imprimitivity $\sigma$ of $\Aut(\Gamma)$ on $V(\Gamma)$ with finite blocks such that $\langle g^{\Aut(\Gamma)} \rangle^{\sigma}$ is a finitely generated abelian group
\index{group ! abelian}%
with finite index $|\Aut(\Gamma / \sigma) : C_{\Aut(\Gamma / \sigma)} (\langle g^{\Aut(\Gamma)} \rangle^{\sigma})|$. 
\end{theorem}

\begin{theorem}
The automorphism group
\index{group ! automorphism}%
 of the previous theorem is bounded if and only if there is a system of imprimitivity of $\Aut(\Gamma)$ on $V(\Gamma)$ with finite blocks such that $\Gamma / \sigma$ is a Cayley graph of the abelian group
\index{group ! abelian}%
 $G^{\sigma}$ acting on it in the natural way if and only if $G^{\sigma}$ is a finitely generated free abelian group.
\index{group ! free abelian}%
\end{theorem}

A \emph{lattice}
\index{lattice}%
is precisely a connected locally finite graph with a vertex-transitive group of bounded automorphisms.

The class of all connected vertex-transitive graphs
\index{graph ! vertex-transitive}%
forms a locally compact, complete and totally disconnected metric space.
\index{metric space}%
Giudici, Li, Praeger, Seress and Trofimov
\index{Giudici, M.}%
\index{Li, C. H.}%
\index{Praeger, C. E.}%
\index{Seress, \'A.}%
\index{Trofimov, V. I.}%
studied~\cite{giudici1} the structure of graphs which are limit points of convergent sequences in this space of the subset of all finite graphs admitting a vertex-primitive automorphism group.
\index{group ! automorphism}%
  These \emph{limit graphs}
\index{graph ! limit}%
 differ from Fra\"{\i}ss\'e limits~\cite{frai}
\index{Fra\"{\i}ss\'e limit}%
  and from the theory of graphs as homogeneous structures~\cite{hrush}.  There is also a dependent concept of a limit of a sequence of automorphisms, called a \emph{limit automorphism}.
\index{automorphism ! limit}%
 A certain class of limiting graphs was shown to be Cayley graphs
\index{graph ! Cayley}%
of a finite-rank free abelian group.
\index{group ! free abelian}%
The Leech lattice
\index{lattice ! Leech}%
 can be realised as a limit graph, with the double cover of $Co_{1}$ as accompanying lattice reflection group, in the sense of~\cite{giudici1}, as can any root lattice of a crystallographic finite reflection group.

A different study of so-called \emph{graph limits} is mentioned in section 5 of this Appendix.

\medskip

A permutation group $(G, \Omega)$ is \emph{2-closed}
\index{group ! permutation ! $2$-closed}%
if any permutation of $\Omega$ which fixes every $G$-orbit on pairs belongs to $G$, or alternatively if and only if there is a complete edge-coloured digraph
\index{digraph}%
 $\Gamma$ with vertex set $\Omega$ for which $\Aut(\Gamma) = (G, \Omega)$.  The notion of $k$-closure is defined similarly: a group $G$ is $k$-closed if every permutation which fixes all $G$-orbits on $k$-sets belongs to $G$.  The group $G$ is \emph{closed}
\index{group ! closed}%
in the natural topology
\index{topology ! natural}%
 on $\Sym(\Omega)$ if and only if $G$ is $k$-closed for all $k \ge 2$.  

\medskip

There are many topics relating to infinite permutation groups that have an extensive literature, but are not utilised in this monograph, so we have omitted addressing them.  These include: Jordan groups~\cite{kayem}, and local and asymptotic behaviour of integer sequences enumerating orbits including studies of their growth rates and cycle index theory~\cite{cameron}.  There are others that have a confluence with topology, and we mention a select few of these in a later section.  

We end with an influential result proved in~\cite{macneu}:  if $\Omega$ is an infinite set then $\Sym(\Omega)$ is not the union of a chain of $|\Omega|$ or fewer proper subgroups.  An extension of the proof of this result gives the following one~\cite{bergman}: if $S$ is a generating set for $\Sym(\Omega)$ then there exists a positive integer $n$ such that every element of $S$ can be written as a group word of length at most $n$ in the elements of $S$.

\medskip

Thorough treatments of permutation groups can be found in~\cite{bhatmo} \cite{cam1} and~\cite{dixon}.  The book~\cite{meldrum} is devoted to a study of wreath products and its generalizations.  More on oligomorphic permutations groups can be found in~\cite{cam6}, and on cofinitary permutation groups in~\cite{cam5}.  Different aspects of the concept of closure of permutation groups are studied in~\cite{lipr},~\cite{prae} and~\cite{farad1}. 

\smallskip

Finally in this section, we mention the \emph{HNN-construction}~\cite{higmannn}
\index{HNN-construction}%
of an infinite group all of whose non-identity elements are conjugate, thus creating an infinite simple group.  The diagonal group $G^*$
\index{group ! diagonal}%
 stabilizes the identity, is transitive on the non-identity elements, so $G$ is $2$-transitive.
\index{group ! permutation ! $2$-transitive}%

If $A, B \leq G$ and $\phi : A \to B$ is an isomorphism then for $G' = \langle G, t : t^{-1}at = a\phi \mbox{ for all}\ a \in A \rangle$
\begin{itemize}
\item $G$ is embedded isomorphically in $G'$;
\item any finite-order element of $G'$ lies in a conjugate of $G$ (so that, if $G$ is torsion-free
 \index{group ! torsion-free}%
then so is $G'$);
\item $t^{-1}Gt \cap G = B$ and $tGt^{-1} \cap G = A$.
\end{itemize}

The group $G$ is called the \emph{base} of $G^*$, $t$ is called the \emph{stable letter}, $A, B$ the \emph{associated subgroups}, and $G^*$ the \emph{HNN-extension} of $G$ relative to $A, B$ and $\phi$.

If $1 \neq a, b \in G$ and $G$ is torsion-free, then $a$ is conjugate to $b$ in $G'$ and $G'$ is torsion-free.  Repeating this for each pair of non-identity elements of $G$, and putting $G_0 = G$ and $G_{n+1} = G_n^{\dagger}$ for all $n \in \mathbb{N}$, then $G^* = \bigcup_{n \in \mathbb{N}} G_n$ is countable, and any two non-identity elements of $G^*$ lie in $G_n$ for some $n$, so are conjugate in $G_{n+1}$.

Good discussions of this technique can be found in the book of Don Collins et al~\cite{collins} and in~\cite{lyndons}.

\section{Model Theory}
\label{ModelTheory}
\index{model theory}%

Model theory is the study of the relationship between the logical properties of sets of formulas and the mathematical properties of structures satisfying these formulas.  

Our focus is on \emph{first-order logic}.
\index{first-order logic}%
The concept of \emph{first-order property}
\index{first-order property}%
comes from logic.  Briefly, first-order \emph{languages}
\index{language}%
include logical connectives ($\neg, \Rightarrow, \vee, \wedge$), equality symbol ($=$), quantifiers ($\forall, \exists$), punctuation marks (parentheses, comma), countably many variables, symbols for functions, relations and constants.  

\emph{First-order}
\index{first-order structure}%
 means that the quantification is only over elements of the structure, whereas second order logic allows quantification over subsets and relations.  Tarski
\index{Tarski, A.}%
 proposed that models of all infinite cardinalities be allowed within first-order logic and this became the accepted view.

For example, to say that a graph is connected is not a first-order sentence,
\index{first-order sentence}%
 because a path connecting two vertices may be finite but arbitrarily long; providing a bound by specifying that the longest path has length $n$ makes it first-order.  A finite graph cannot satisfy the $(*)$-condition, but the statements of this condition are first-order.  Every vertex of the countable random graph
\index{graph ! random}%
 must have degree $\aleph_0$, because a degree-$n$ vertex cannot be a property of almost all graphs.

We refer to~\cite{bar}~\cite{hod1} for a detailed explanation of the term.

A \emph{term}
\index{term}%
is a variable, a constant symbol, or an $n$-ary function symbol of the form $F(t_1, \ldots, t_n)$ where the $t_i$ are terms.  An \emph{atomic formula}
\index{formula ! atomic}%
is an $n$-ary relation symbol of the form $R(t_1, \ldots, t_n)$, and a \emph{formula}
\index{formula}%
is either an atomic formula or of the form $(\phi \vee \psi)$, $(\phi \wedge \psi)$, $(\neg \phi)$, $(\phi \Rightarrow \psi)$, $(\phi \Leftrightarrow \psi)$, $(\forall x)\phi$, or $(\exists x) \phi$, where $\phi$ and $\psi$ are formulas and $x$ is a variable.

A variable $x$ is \emph{bound}
\index{variable ! bound}%
in a formula $\phi$ if it occurs in a subformula $(\forall x)\psi$ or $(\exists x)\psi$, and any other occurences are called \emph{free}.
\index{variable ! free}%
A \emph{sentence}
\index{sentence}%
is a formula with no free variables.

First-order logic
\index{first-order logic}%
 allows the combination of finitely many formulas with connectives, and the quantification over variables which range over the underlying set.  The first-order language
\index{first-order language}%
  of graph theory consists of formulas, which are finite strings of the above symbols together with the adjacency symbol ($\sim$), built inductively.  Sentences in first-order logic must be finite in length and disallow quantification over sets or relations.  So certain common properties, such as connectivity (which without a set quantifier fails the first of these conditions because it requires an infinite disjunction of sentences, one for each positive integer diameter), and $2$-colourability (which requires a colouring function), are not first-order expressible.

A graph is \emph{simple}
\index{graph ! simple}%
if it has only one type of edge or adjacency, and \emph{loopless}
\index{graph ! loopless}%
if no vertex is attached to itself by an edge.  The theory of simple loopless graphs is that of a structure with an irreflexive symmetric binary relation $\sim$ that satisfies the following sentences:
\[ \forall x (\neg x \sim x) \]
\[ \forall x \forall y ((x \sim y \wedge y \sim x) \vee (\neg x \sim y \wedge \neg y \sim x)). \]

A \emph{structure} $M$
\index{structure}%
over a first-order language
\index{first-order language}%
 $L$ is a set $X$ on which are defined constants (elements of $X$), $n$-ary functions (maps $X^n \to X$) and $n$-ary relations (subsets of $X^n$).  A \emph{valuation}
\index{valuation}%
assigns values to the variables in $X$ and truth values (`true' or `false') to the formulas.  The truth value assigned is independent of the values assigned to bound variables in the formula, so a sentence is either true or false in $M$ independent of the valuation.  If the sentence $\sigma$ is true in $M$ then $M$ is said to `satisfy' or `model' $\sigma$, which is written $M \models \sigma$.  If $M$ is a set of such sentences the we write $M \models \Sigma$.  A \emph{theory}
\index{theory}%
$\Th(M)$ of a structure $M$  is the set of all first-order sentences
\index{first-order sentence}%
 which are true in $M$.  

An \emph{expansion}
\index{expansion}%
or \emph{reduct}
\index{reduct}%
refer to the addition or removal of some new constants, functions or relations to a language or structure, but (in the case of structures) keeping the same domain.  In contrast, are notions of \emph{extension},
\index{extension}%
\emph{substructure}
\index{substructure}%
or \emph{submodel}
\index{submodel}%
where the language is retained but the domain is made smaller or larger.  The relations and functions of the submodel are the restrictions of those of the model.  In general mathematical objects can be represented as first-order structures
\index{first-order structure}%
 in different ways but each one will give rise to different substructures.

Thus far we have dealt with syntax and semantics in logic but a third aspect is deduction or inference.  The simplest deductive system specifies that particular formulas are \emph{axioms}
\index{axiom}%
whilst certain formulas can be deduced from others using \emph{rules of inference},
\index{rules of inference}%
such as \emph{modus ponens}
\index{modus ponens}%
where $\psi$ is inferred from $\phi$ and $\phi \Rightarrow \psi$.  In a deductive system a formula is a \emph{theorem}
\index{theorem}%
if it is either an axiom or deduced from other theorems using rules of inference.  Alternatively, a \emph{proof}
\index{proof}%
is a finite sequence of formulas each of which is either an axiom or is inferred from earlier formulas, and a \emph{theorem} is the last line of the proof.

A sentence $\sigma$ `is a consequence of' a set of sentences $T$, written $T \vdash \sigma$ if there is a finite list of sentences each deducible from the previous one, ending with $\sigma$.  A set of formulas is \emph{consistent}
\index{formula ! consistent}%
if no contradiction can be deduced from it.  If a theory has a model then it must be consistent.  A formula is \emph{logically valid}
\index{formula ! logically valid}%
if it is true in all valuations in all structures over the relevant language.  A deductive system is \emph{sound}
\index{system ! sound}%
if every theorem is logically valid and \emph{complete}
\index{system ! complete}%
if every logically valid formula is a theorem.  So a theory is complete if either $T \vdash \sigma$ or $T \vdash \neg \sigma$ for every sentence of $L$.  If all countable models of a consistent theory are isomorphic then the theory is complete.  An important theorem on completeness is the following,

\begin{theorem}[G\"{o}del--Henkin Completeness Theorem]
\index{godel@G\"{o}del--Henkin Completeness ! Theorem}%
\index{godel@G\"{o}del, K.}%
\index{Henkin, L.}%
There is deductive system for first-order logic which is sound and complete.

Moreover, this system has the following properties:
\begin{itemize}
\item[(a)] A sentence $\phi$ can be deduced from a set $\Sigma$ of sentences if and only if every structure satisfying $\Sigma$ also satisfies $\phi$.
\item[(b)] A set $\Sigma$ has a model if and only if it is consistent.
\end{itemize}
\end{theorem}

This theorem leads to the next two central results of model theory.
\index{model theory}%

\begin{theorem}[Compactness Theorem]
\index{Compactness Theorem}%
A set $\Sigma$ of sentences has a model if and only if every finite subset has a model.
\end{theorem}

\begin{theorem}[Downward L\"{o}wenheim--Skolem Theorem]
\index{Downward L\"{o}wenheim--Skolem ! Theorem}%
\index{Skolem, T.}%
\index{L\"{o}wenheim, L.}%
If a set of sentences in a countable language has a model, then it has a finite or countable model.
\end{theorem}

A language or structure is \emph{relational}
\index{relational language}%
\index{relational structure}%
if it contains no function or constant symbols.  Relational structures have the important property that every subset carries a substructure, for example an \emph{induced subgraph}
\index{graph ! induced subgraph}%
(one obtained by deleting vertices).  A structure is \emph{homogeneous}
\index{structure ! homogeneous}%
if every isomorphism between finite substructures extends to an automorphism of the entire structure.

The interplay of model theory
\index{model theory}%
 and permutation group theory has been very fruitful~\cite{kayem}.  The standard way of treating a permutation group $G$ acting on a set $\Omega$ as the automorphism group
\index{group ! automorphism}%
 of a relational structure with domain $\Omega$, is to utilize a relational language $L$ which has for every positive integer $k$ and each $G$-orbit of $\Omega^k$, a unique $k$-ary relation symbol $r$, and interpret $r(x_1, \ldots, x_k)$ to be true only if $(x_1, \ldots, x_k)$ lies in the orbit.  In this case $L$ is called the \emph{canonical language}
\index{canonical language}%
and the resulting structure $M$ as the \emph{canonical structure}
\index{canonical structure}%
for $G$.  If $\Omega$ is infinite then it is possible that $G < \Aut(M)$.  Moreover $\Aut(M)$ is the \emph{closure}
\index{group ! permutation ! closure}%
of $G$, that is they have the same orbits on the set of finite ordered subsets of $M$, (see the section on topology below).  For example assuming that $G$ is not $2$-homogeneous on $\Omega$, it can be realized as an automorphism group
\index{group ! automorphism}%
 of a graph on $\Omega$ which is neither null nor complete, where the orbit on unordered pairs can be taken to be the edge set; the group may not be the \emph{full} automorphism group of the graph.

A \emph{tuple}
\index{tuple}%
from the domain $M_s$ of a structure $M$ is an element of $M_{s}^k$ for some $k \in \mathbb{N}$, that is, a $k$-tuple for some $k$.  The \emph{type}
\index{type}%
$\tp_{M}(\bar{a})$ of a tuple $\bar{a} \in M_{s}^k$ is the set of all formulas $\phi(x_0, \ldots, x_{k-1})$ which are true in $M$ when $a_0, \ldots, a_{k-1}$ are substituted for $x_0, \ldots, x_{k-1}$; this is in fact a \emph{$k$-type},
\index{type@$k$-type}%
that is one of length $k$.  If $B \subset M_s$ then $\tp_{M}(\bar{a} / B)$ is the type $\bar{a}$ in the expanded structure $(M, b)_{b \in B}$.  A \emph{type over a theory} $T$
\index{type over a theory}%
is a set of formulas realizable as $\tp_{M}(\bar{a})$ for some $\bar{a} \in M \models T$, that is a maximal set of formulas consistent with $T$ when $x_0, \ldots, x_{k-1}$ are treated as new constant symbols.  This can also be called a \emph{type over $\emptyset$}
\index{type over $\emptyset$}%
because it has no additional parameters.  A \emph{type over a model} $M$
\index{type over a model}%
is a type over $\Th(M, a)_{a \in M}$.  The following definition is intended to capture the notion of a `large' model, as being one in which many types are realized.  A structure $M$ is \emph{$\lambda$-saturated}
\index{structure ! $\lambda$-saturated}%
where $\lambda$ is an infinite cardinal, if for every $A \subset M_{s}$ with $\card(A) < \lambda$,~\label{card} every type over $\Th(M, a)_{a \in M}$ is realized in $M_{s}$.  A complete theory has, up to isomorphism, at most one saturated model in a given cardinality.

A finite structure can be described by a single sentence.  The \emph{Upward L\"{o}wenheim--Skolem Theorem}
\index{Upward L\"{o}wenheim--Skolem ! Theorem}%
\index{Skolem, T.}%
\index{L\"{o}wenheim, L.}%
asserts that if a set of sentences has an infinite model, then it has models of arbitrarily large cardinality.  This follows easily from the Compactness Theorem.
\index{Compactness Theorem}%
  Cardinality of a structure is not a first-order property,
\index{first-order property}%
 so first-order axioms can at best specify a structure with the same cardinality as $M$, as being isomorphic to $M$.  When this holds and $|M| = \lambda$ we say that $M$ is \emph{$\lambda$-categorical}.
\index{theory ! $\lambda$-categorical}%
This concept applies also to a theory (a set of sentences).  A theory that is $\lambda$-categorical for all infinite $\lambda$ is called \emph{totally categorical}.
\index{theory ! totally categorical}%
Totally categorical theories are concerned with those infinite structures about which most can be said using first-order statements.  Modulo choice of language there are only countably many totally categorical theories and this motivates a programme to classify them~\cite{hodges1}.  It is known that, as long as a theory has infinite models, it cannot have a unique model; there 
will be models of arbitrarily large cardinality. The best we can do is ask 
that the theory is $\alpha$-categorical, where $\alpha$ is an infinite cardinal, meaning that 
there is a unique model of cardinality $\alpha$ up to isomorphism. By a theorem of 
Vaught, there are only two types of categoricity,
\index{categoricity}%
 countable and uncountable: 
the theorem of Vaught~\cite{vaught} states that if a theory in a countable language is $\lambda$-categorical for any uncountable cardinal $\lambda$, then it is $\mu$-categorical for any uncountable cardinal $\mu$.  So there is a dichotomy in types of categoricity: countable and uncountable.  Uncountable categoricity gives rise to a powerful structure theory, extending that for vector spaces over a fixed field or algebraically closed fields 
of fixed characteristic (where a single invariant, the rank, determines the model). According to Zil'ber's Trichotomy Theorem~\cite{zilber}
\index{Zil'ber, B. I.}%
\index{Zil'ber's Trichotomy Theorem}%
the class of all uncountably categorical structures divides, with respect to their properties, into three types: (1) field-like, (2) module-like, and (3) of disintegrated (discrete) type.  Countable categoricity is more akin to a symmetry condition.  One of the classical results that connect model theory
\index{model theory}%
 to permutation group theory, proved independently by Engeler~\cite{engeler}
\index{Engeler, E.}%
Ryll-Nardzewski~\cite{ryll}
\index{Ryll-Nardzewski, C.}%
and Svenonius~\cite{svenonius}
\index{Svenonius, L.}%
is

\begin{theorem}[Engeler--Ryll-Nardzewski--Svenonius Theorem]
\quad\\
The countable first-order structure
\index{first-order structure}%
 $M$ is countably categorical if and only if $\Aut(M)$ is oligomorphic.  Moreover, if these conditions hold then two tuples lie in the same orbit of $\Aut(M)$ if and only if they have the same type (that is, they satisfy the same $n$-variable first-order formulas).
\end{theorem}

As we mentioned in the preface, this theorem implies that axiomatisability of a structure is equivalent to the existence of a very large automorphism group,
\index{group ! automorphism}%
 a situation reminiscent of Felix Kelin's
\index{Klein, F.}%
 Erlanger Programm~\cite{klein}. 

A more substantial statement of this theorem, including a list of equivalences is given in~\cite{hod1}.  The theorem indicates that the symmetry groups of $\aleph_{0}$-categorical
\index{aleph@$\aleph_0$-categorical}%
structures are rich in structure.  Each isomorphism type of $n$-element substructures contributes at most $n!$ orbits on $M^{(n)}$, leading to the following corollary,

\begin{corollary}
Let $M$ be a countable homogeneous relational structure
\index{relational structure}%
 having only a finite number of isomorphism types of $n$-element substructures for all $n$.  Then $M$ is $\aleph_0$-categorical.  In particular, this holds if there are only finitely many relation symbols in the language of $M$.
\end{corollary}
The multicoloured random graphs
\index{graph ! random ! $m$-coloured}%
 that are the subject of this work exemplify this result.

In general, primitivity
\index{group ! permutation ! primitive}%
of a permutation group is not first-order, but it is if there is a bound on the number of $2$-types in $\Omega$ (that is orbits on pairs)~\cite[p.21]{kayem}.  From Higman's work~\cite{higman}
\index{Higman, D. G.}%
it follows that $G$ is primitive
\index{group ! permutation ! primitive}%
 if and only if each $G$-orbit on unordered pairs, thought of as the edge set of an undirected graph, is path-connected.  If $G$ is thought of as a first-order structure
\index{first-order structure}%
 and has finitely many orbits on $\Omega^2$, for example if we assume a graph of bounded diameter $\le k$ (which is a first-order property),
\index{first-order property}%
  then the claim follows.  Of course primitivity is not first-order if no such bound exists.

\medskip

We mention a general question.  Which permutation groups of countable degree are automorphism groups
\index{group ! automorphism}%
 of relational structures
\index{relational structure}%
  over finite relational languages?  Frucht~\cite{frucht} showed that every abstract group is the automorphism group of some simple undirected graph.

\medskip

We end with two notions that will be used in the section on topology in permutation groups below, and which arise in the classification of countable totally categorical structures;  let $M$ and $N$ be two such structures.  

Totally categorical theories are determined purely by their countable models.  A \emph{definable $k$-ary relation}
\index{definable ! $k$-ary relation}%
in a structure $M$ is a set $\phi(M^k)$ where $\phi(x_0, \ldots, x_{k-1})$ is a formula on $L$ meaning $\{ \bar{a} : M \models \phi(\bar{a}) \}$; if there are no parameters from $M$, we say that the structure is \emph{$\emptyset$-definable}.
\index{definable ! $\emptyset$-definable}%
The structures $M$ and $N$ are \emph{bi-definable}
\index{bi-definable structure}%
if there is a bijection from $\dom(M)$ to $\dom(N)$ taking $\emptyset$-definable relations to $\emptyset$-definable relations; this induces an isomorphism of automorphism groups,
\index{group ! automorphism}%
 $\Aut(M) \cong \Aut(N)$, as permutation groups.

A weaker notion than bi-definability is \emph{bi-interpretability}.
\index{bi-interpretable structure}%
\index{structure ! bi-interpretable}%
We refer to~\cite{kayem} for proper definition, but note the following:

\begin{theorem}
Let $M$ and $N$ be countable structures.  Then
\begin{itemize}
\item[(A1)]  $\Aut(M) \cong \Aut(N)$ as topological groups if and only if $M$ and $N$ are (infinitarily) bi-interpretable.
\item[(A2)]  If $M$ and $N$ are both $\aleph_0$-categorical then $\Aut(M) \cong \Aut(N)$ as topological groups if and only if they are finitarily\\ bi-interpretable.
\end{itemize}
\end{theorem}

\bigskip

The theorems of Fra\"{\i}ss\'e's (see later appendix) and Engeler--Ryll-Nardzewski--Svenonius
\index{Engeler--Ryll-Nardzewski--Svenonius ! Theorem}%
 are connected through the concept of quantifier elimination,
\index{quantifier elimination (q.e.)}%
  so we say a little about this.
  
A substructure $M$ of a structure $N$ is \emph{elementary} if the identity map on $M$ is elementary; $N$ is an \emph{elementary extension} of $M$.  A complete first-order theory
\index{first-order theory}%
 $T$ has \emph{quantifier elimination (q.e.)}
\index{quantifier elimination (q.e.)}%
 if, given any formula $\phi(\bar{x})$, there is a quantifier-free formula
\[ (\forall \bar{x}) (\phi(\bar{x}) \leftrightarrow \psi(\bar{x}))\]
provable from $T$.  By the Completeness Theorem,
\index{godel@G\"{o}del--Henkin Completeness ! Theorem}%
 this is equivalent to saying that the displayed sentence holds in every model of $T$.  A structure $M$ has q.e. if its theory does.  If a theory has q.e. then tuples have the same type if and only if they are isomorphic.  For a countable and $\aleph_0$-categorical 
\index{aleph@$\aleph_0$-categorical}%
structure q.e. is equivalent to homogeneity~\cite[p.44]{cam6}.  The concept of \emph{existential closure} 
\index{existentially closed (e.c.)}%
 (see Appendix~\ref{FurtherDetails}) is a generalization of quantifier elimination -- a substructure $M$ of a structure $N$ is \emph{existentially closed (e.c.)} in $N$ if, for every tuple $\bar{a} \in M$ and existential formula $\phi(\bar{x})$, $\phi(\bar{a})$ holds in $M$ if and only if it holds in $N$.  A first-order theory $T$ has q.e. if and only if every model of $T$ is e.c. in a homogeneous model of $T$.
 
Let $\mathbf{K}$ a class of $L$-structures for a first-order language $L$.
\index{first-order language}%
 A set $\Phi$ of formulas of $L$ is an \emph{elimination set} for $\mathbf{K}$ if the formulas are equivalent (Boolean combinations of each other) in every structure in $\mathbf{K}$.

One of the uses of \emph{quantifier elimination} is that structures in $\mathbf{K}$ can be classified up to elementary equivalence by inspection of equivalent sentences in the elimination set.  As a special case, if $\mathbf{K}$ is the class of all models of a first-order theory $T$ and all sentences in $\Phi$ are either deducible from $T$ or inconsistent with $T$, then all models of $T$ are elementarily equivalent, and so $T$ is a complete theory.  Another use is in the description of elementary embeddings; if $\Phi$ is an elimination set for $\mathbf{K}$, then the elementary maps
\index{elementary map}%
 between structures in $\mathbf{K}$ are precisely those \emph{homomorphisms} which preserve $\psi$ and $\neg \psi$ for every formula $\psi \in \Phi$.  

A theory with q.e. is \emph{model-complete},
\index{theory ! model-complete}%
 that is where every embedding of one model of $T$ into another is an elementary map, that is preserves all first-order sentences.  (By a result of Abraham Robinson,
\index{Robinson, A.}%
 it suffices that every embedding preserves existential formulas.)  The converse is false, but $T$ has q.e. if and only if it is model-complete and the class of substructures of all models of $T$ has the amalgamation property.  

The theory of a countable $\aleph_0$-categorical structure is model-complete if and only if it is equivalent to a class of $(\forall\exists)$-sentences.   Further, this is true if the structure is also homogeneous.

The random graph $\mathfrak{R}$ has quantifier elimination,
\index{quantifier elimination (q.e.)}%
 but it was shown in~\cite{bodirsky} that its reducts do not, though its reducts are model-complete.  However, they also show that
 
\begin{theorem}
Every reduct of $\mathfrak{R}$ has quantifier elimination if it is expanded by all relations with an existential definition in the reduct.
\end{theorem}
 
\bigskip

\emph{The Random Graph and Stability}

Central to model theory
\index{model theory}%
 is the classification of mathematical structures using logical formulas, that is, which logical sentences are true in the structures, as well as the ways of constructing models of given sentences.  

For example, the first-order theory
\index{first-order theory}%
 of the random graph
\index{graph ! random}%
  is equivalent to the theory

\[ T = \forall x\ \neg (x \sim x) \cup \forall x y (x \sim y \leftrightarrow y \sim x) \cup \{ \phi_{a, b} : a, b < \omega \}, \]
where $\phi_{a, b} = \forall x_1 \ldots x_a y_1 \ldots y_b ( \bigwedge_{i < j} x_i \neq x_j \wedge \bigwedge_{i < j} y_i \neq y_j  \wedge \bigwedge_{i, j} x_i \neq y_j \longrightarrow \exists z ( \bigwedge_{i} (z \neq x_i \wedge z \sim x_i)\wedge \bigwedge_{j} (z \neq y_j \wedge \neg (z \sim y_j))))$.

Shelah introduced the idea of \emph{stability}
\index{stability}%
 as an aid to counting the number of non-isomorphic models of a theory in a given cardinality, going on to classify all complete first-order theories~\cite{shelah}.
\index{first-order theory}%
   A complete theory $T$ in a first-order language
\index{first-order language}%
  $L$ is \emph{unstable}
\index{theory ! unstable}%
if there is a formula $\phi(\bar{x}, \bar{y})$ of $L$ and a model $M$ of $T$ containing tuples of elements $\bar{a}_i$ $(i < \omega)$ such that
\[ \text{for all } i, j < \omega, M \models \phi(\bar{a}_i, \bar{a}_j) \Leftrightarrow i < j.\]
$T$ is \emph{stable} if it is not unstable.  The structure $M$ is stable or unstable according to whether or not $T$ is.

The definable sets of a structure are a function of the language used.  Building on earlier ideas of Morley~\cite{morl}, Shelah used the complexity of the definable sets to count the number of models.  If a theory is not stable then heuristically, its models are too numerous to classify, while a stable theory has a chance of classification.  For example, if the definable sets are nested in a complicated fashion, it is less likely that there will be isomorphisms between models.

A formula $\phi(\overline{x}, \overline{y})$ has the \emph{strict order property} (SOP)
\index{strict order property}%
 (for a complete theory $T$) if in every model $M$ of $T$, $\phi$ defines a partial order on the set of $n$-tuples with an infinite (or arbitrarily long finite) chain(s).  The theory $T$ has the SOP if some formula has the SOP for $T$.

A formula $\phi(\overline{x}, \overline{y})$ has the \emph{independence property} (IP)
\index{independence property}%
 (for a complete theory $T$) if in every model $M$ of $T$ there is, for each $n < \omega$, a family of tuples $\bar{b}_{0}, \ldots, \bar{b}_{n-1} \in M$ such that for every subset $X$ of $n$ there is a tuple $\bar{a} \in M$ for which $M \models \phi(\bar{a}, \bar{b}_{i})$ if and only if $i \in X$.  The theory $T$ has the IP if some formula has the IP for $T$.

Shelah defined a theory $T$ to be \emph{unstable} if there is a formula interpreting either the IP or the SOP.  The prototypical theory having SOP is the theory of dense linear orders
\index{linear order}%
 without endpoints, $(\mathbb{Q}, <)$.  The theory of the (uniform) random graph
\index{graph ! random}%
 is unstable, having the IP but not the SOP.

\bigskip

More on model theory
\index{model theory}%
 can be found in standard texts such as~\cite{hod1}, or its shorter version~\cite{hod2}.  A survey of the first-order theory
\index{first-order theory}%
 of graphs appears in~\cite{cam7b}.  More on the model theory
\index{model theory}%
 of groups and automorphism groups
\index{group ! automorphism}%
  can be found in~\cite{evans1}.

\vspace{1cm}

\section{Category and Measure}
\label{CategoryandMeasure}
Let $(X, d)$ be a complete metric space.  A subset of $X$ is \emph{dense}
\index{dense set}%
if its closure is $X$, that is if it meets every open set.
\index{open set}%
  A subset is
\emph{residual}
\index{residual set}%
or \emph{comeagre}
\index{sets ! comeagre}%
or \emph{of second category}
\index{second category}%
if it contains the intersection of countably many dense open sets.
\index{dense open set}%
The complement of such a set is called \emph{meagre}
\index{meagre}%
or \emph{of the first category}.
\index{first category}%

The \emph{Baire category theorem}~\cite{baire}
\index{Baire category theorem}%
\index{Baire, R.}%
states that a residual set in a complete metric space is non-empty.  (See also Osgood~\cite{osgood}).
\index{Osgood, W.}%

A residual set is dense; in other words it has a non-empty intersection with every dense open set.
\index{dense open set}%
Moreover, any countable collection of residual sets has non-empty (even residual) intersection.  For this reason we regard residual sets as ``large'', containing ``almost all'' of the space.  (For example Cantor took uncountable sets (such as the reals) to be large and to contain finite or countable subsets.  But there are two families of uncountable small sets:

$\bullet$\  the \emph{null sets} with respect to Lebesgue measure (those sets which can be covered by a countable union of intervals with total length at most $\epsilon$, for any $\epsilon > 0$);

$\bullet$\  the \emph{meagre sets} in the topology
\index{topology}%
 (see below), (those contained in a countable union of closed sets
\index{closed set}%
 with empty interior).

Both classes include the countable set.  Sierpi\'nski
\index{Sierpi\'nski, Waclaw}%
 showed that there is a permutation of the real numbers which maps null sets to meagre sets; Erd\H{o}s
\index{Erd\H{o}s, P.}%
 showed that we can take this permutation to be an involution interchanging the two families. So, as families of sets, they are isomorphic~\cite{oxt}.)
 
It follows that if we can show that the set of elements having property $\mathcal{P}$ is residual, then some element will have property $\mathcal{P}$. 

Metric spaces relevant to our study typically arise as follows.  A point in the space is determined by a countable sequence of choices.  The larger the agreement between the initial segment of the choice sequences determining two points, the closer are the two points.  If points $x$ and $y$ differ first in the $n^{\text{th}}$ term say, any decreasing function can determine their distance apart, for example $d(x, y) = \frac{1}{2^n}$.

We can see how this applies to the set of paths in a tree.  The nodes of the tree occur on `levels' indexed by the first $n$ natural numbers $0, \ldots, n-1$.  Define a metric on this set by the rule that two paths which agree up to level $n$ and diverge thereafter should have distance $\frac{1}{2^n}$.  This choice of distance function gives a space whose diameter is 1, but as long as a decreasing function of $n$ is used its precise form is not crucial.  With this metric the set of paths forms a complete metric space.

Before defining measures we mention some spaces on which a they may be defined.  A \emph{field} $\mathcal{F}(X)$~\label{mathcalF(X)} 
\index{field}%
on a space $X$ is a nonempty collection of subsets of $X$ such that
\begin{itemize}
\item[(i)]  if $O_1, O_2 \in \mathcal{F}(X)$ then $O_1 \cup O_2 \in \mathcal{F}(X)$;
\item[(ii)]  if $O \in \mathcal{F}(X)$ then $X \backslash O \in \mathcal{F}(X)$.
\end{itemize}
If in addition,
\begin{itemize}
\item[(iii)]  whenever $O_i \in \mathcal{F}(X)$ $\forall i = 1, 2, \ldots$ we have $\bigcup_{i=1}^{\infty} O_i \in \mathcal{F}(X)$,
\end{itemize}
then the collection $\mathcal{F}(X)$ is called a \emph{$\sigma$-algebra}.
\index{sigma@$\sigma$-algebra}%
The space of all subsets is the largest $\sigma$-algebra.  Another is the set $\mathcal{F}'(X)$ of all sets that can be described by finite-length expressions in the space $X$ together with any collection $G$ of subsets of $X$, using unions and complements; we say that $\mathcal{F}'(X)$ is the \emph{field generated} by $G$.  The intersection $\mathcal{F}_{G}(X)$ of all $\sigma$-algebras  on $X$ containing $G$, that is \emph{generated by} $G$, is the smallest $\sigma$-algebra on $X$ containing $\mathcal{F}'(X)$.  If $G$ has a finite set of generators then $\mathcal{F}'(X) = \mathcal{F}_{G}(X)$.

A \emph{topological space}
\index{topological space}%
$(X, T)$~\label{(X, T)} is a space $X$ together with a topology $T = T(X)$. A \emph{topology}
\index{topology}%
$T$ for $X$ is a set of subsets of $X$ such that
\begin{itemize}
\item[(i)]  $\emptyset, X \in T$,
\item[(ii)]  if $\{O_i : i \in \mathcal{I} \} \subset T$ is any collection of members of $T$ then $\bigcup_{i \in \mathcal{I}} O_i \in T$,
\item[(iii)]  if $\{O_n : n = 1, \ldots, N\}$ is any finite set of members of $T$ then $\bigcap_{n = 1}^{N}  O_n \in T$.
\end{itemize}
The sets of $T$ are called \emph{open} sets.
\index{open set}%

Let $(X, T)$ be a topological space.  Then the set $\mathcal{B}(X)$~\label{mathcal{B}(X)} of \emph{Borel subsets}
\index{Borel set}%
of $X$ is the $\sigma$-algebra $\mathcal{F}_{T(X)}(X)$ generated by the open subsets $T(X)$ of $X$.

A (finite, positive) \emph{measure}
\index{measure}%
$\mu$ on a space $X$ is a function $\mu: \mathcal{F}(X) \to [0, \infty)$, where $\mathcal{F}(X)$ is a field, and 
\[ \sum_{n=1}^{\infty} \mu(O_n) = \mu \left( \bigcup_{n=1}^{\infty} O_n \right)\]
whenever $\{O_n \in \mathcal{F}(X) : n = 1, 2, \ldots\}$ is a sequence such that for all $n, m \in \mathbb{N}$ with $n \neq m$, $\bigcap_{n = 1}^{\infty}  O_n \in \mathcal{F}(X)$ and $O_n \cap O_m = \emptyset$.

One interpretation of $\mu$ is as a probability law in which $\mu(E)$ measures the probability that a randomly chosen point $x \in X$ belongs to $E$.  A \emph{probability} (or \emph{normalized}) measure
\index{measure ! probability (normalized)}%
on $X$ is one for which $\mu(X) = 1$.  The measure for a Borel $\sigma$-algebra is called a \emph{Borel measure}.
\index{measure ! Borel}%
Denote the set of Borel measures on $X$ by $\mathcal{M}(X)$.  The \emph{support}
\index{measure ! support}%
of a Borel measure $\mu$ on a metric space $(X, d)$ is the set of points $x \in X$ such that $\mu(O_x) > 0$ whenever $O_x$ is an open set
\index{open set}%
 that contains $x$.  Sometimes a measure can be described by a \emph{density function},
\index{density function}%
and it may even be permissible that this function be piecewise continuous with discontinuities appearing on certain sets.

The next theorem states that a measure on a field uniquely extends to one on the $\sigma$-algebra generated by the field, and that the measure on the $\sigma$-algebra can be evaluated using only its values on the field.

\begin{theorem}
Let $X$ be a space, $\mathcal{F}'(X)$ be a field on $X$ and $\mathcal{F}(X)$ be the smallest $\sigma$-algebra on $X$ that contains $\mathcal{F}'(X)$.  Let $\mu' : \mathcal{F}'(X) \to [0, \infty)$ be a measure.  Then there exists a unique measure $\mu : \mathcal{F}(X) \to [0, \infty)$ such that $\mu(B) = \mu'(B)$ for all $B \in \mathcal{F}'(X)$.
Moreover for all $A \in \mathcal{F}(X)$,
\[ \mu(A) = \inf \left\{ \sum_{n=1}^{\infty} \mu(B_n) : A \subset \bigcup_{n = 1}^{\infty} B_n, B_n \in \mathcal{F}'(X)\ \forall n = 1, 2, \ldots \right\} . \]
\end{theorem}

A \emph{measureable set}
\index{measureable set}%
is one to which a measure can be assigned, that is, a member of a field or $\sigma$-algebra.  The next theorem states that a continuous transformation of a Borel measure is also Borel, because Borel sets
\index{Borel set}%
 are generated by open sets
\index{open set}%
  and the inverse transformations of open sets, under continuous transformations, are open sets.

\begin{theorem}
Let $\nu \in \mathcal{M}(X)$ be a Borel measure and let $f : X \to X$ be a continuous function.  Then there exists on $X$ a unique Borel measure $\mu \in \mathcal{M}(X)$ such that
\[ \mu(B) = \nu(f^{-1}(B))\ \forall B \in \mathcal{B}(X). \]
\end{theorem} 
The measure $\mu$ is called the \emph{transformation of the measure} $\nu$
\index{measure ! transformation of}%
by the function $f$, and denoted by $f(\nu)$ or $f \circ \nu$.

Let $(X, d)$ be a metric space and let $f : X \to X$ be a continuous function.  A measure $\mu \in \mathcal{M}(X)$ is said to be invariant under $f$ if and only if
\[ \mu(B) = \mu(f^{-1}(B))\ \forall B \in \mathcal{B}(X). \]
Such a measure is called an \emph{invariant measure}
\index{measure ! invariant}%
of the transformation $f$, and the condition can be written $f(\mu) = \mu$.

\begin{theorem}
If $f, g : X \to X$ are continuous functions and the measure $\mu \in \mathcal{M}(X)$ is invariant under $f$ then the measure $g(\mu)$ is invariant under $g \circ f \circ g^{-1}$.
\end{theorem} 

\smallskip

The two notions of ``largeness'' that appear in topology
\index{topology}%
 and Baire category theory
\index{Baire category theorem}%
 are different.  For example consider zero-one sequences
 \index{zero-one sequence}%
  and define the \emph{upper density}
\index{upper density ! of zero-one sequences}%
of a zero-one sequence to be  $ \lim \sup_{n \to \infty} \frac{s(n)} {n}$
where $s(n)$ is the number of ones among the first $n$ terms of $s$, and the \emph{lower density} to be defined similarly as $\lim \inf_{n \to \infty} \frac{s(n)} {n}$.
\index{lower density ! of zero-one sequences}%
According to the Law of Large Numbers,
\index{Law of Large Numbers}%
the set of sequences having density $\frac{1}{2}$ has measure 1.  However in the metric defined above, the set of sequences with upper density 1 and lower density 0 is residual.  The ``large'' sets, analogous to residual sets, are sets of measure $1$.  So, for example, any countable intersection of sets of measure $1$ has measure $1$, and any set of measure $1$ has cardinality $2^{\aleph_0}$.

In probability theory, a $\sigma$-algebra is often referred to as a \emph{$\sigma$-field}.
\index{sigmaf@$\sigma$-field}%
 A \emph{probability space}
\index{probability space}%
is a triple $(X, \mathcal{F}(X), \mu)$, where $X$ is a set, $\mathcal{F}(X)$ is a $\sigma$-field of subsets of $X$, $\mu$ a non-negative measure on $X$, and $\mu(X) = 1$.  Then $\mu$ is determined by a function $X \to [0, 1], y \to \mu(\{y\})$, $\mu(Y) = \sum_{y \in Y} \mu(\{y\}), Y \subset X$.  A real-valued \emph{random variable}
\index{random variable}%
is a measurable real-valued function on a probability space.

In a probability space (a measure space of total measure 1), sets of measure $1$
\index{measure}%
are considered ``large''.  They share the property
with residual sets that the intersection of countably many such sets
is again of measure $1$, and hence non-empty.  The complementary
statement is that the union of countably many null sets is a null
set.  Any non-empty open set
\index{open set}%
 has positive measure; so a set of measure
$1$ is dense.

The notions of small sets that we mentioned at the beginning of the section can be widely generalized.  

In any measure space, in particular in any probability space, the null sets form a family of small sets in our earlier sense, and the property that the union of countably many null sets is null, is a simple consequence of Kolmogorov's
\index{Kolmogorov, A. N.}%
 axioms for probability.

On the other hand, in any metric space, the \emph{interior} of a set $A$ is the set of points having a neighbourhood contained in $A$; a set is \emph{open} if it is equal to its interior; a set is \emph{closed} if its complement is open; and a set is \emph{meagre} if it is contained in a countable union of closed sets
\index{closed set}%
 with empty interior. It is clear that a countable union of meagre sets is meagre and a re-statement of the Baire category theorem
\index{Baire category theorem}%
 would be: a complete metric space is not meagre.

So the meagre sets in any complete metric space satisfy the requirements for a family of small sets.  The ``large'' sets will then be the comeager ones.  A set is \emph{comeagre}
\index{sets ! comeagre}%
 if and only if it contains a countable intersection of open dense sets (where a set is \emph{dense} if every neighbourhood of every point in the space meets it).

There is a specific example which can be used to give a characterization of $\mathfrak{R}$ and so is instructive to elucidate it.   Let $X$ be the set of all infinite sequences of zeros and ones. We define the appropriate structures on $X$ as follows:

$\bullet$\  as a probability space,
\index{polynomial algebra}%
 we regard elements of $X$ as recording infinitely many tosses of a fair coin;

$\bullet$\  as a metric space, two sequences are close if they agree on a long initial segment (for example, if the first disagreement is in position $n$, we take the distance to be $\frac{1}{n}$).

The metric space is not only complete, but in fact, it satisfies a strengthening of the triangle inequality known as the hypermetric inequality:
\[ d(x,z) \leq \max( d(x,y),d(y,z) ). \]

A subset $A \subset X$ is \emph{open} if and only if it is \emph{finitely determined}, that is, every sequence in $A$ has a finite prefix all of whose continuations belong to $A$; and a set $A$ is \emph{dense} if and only if it is \emph{always reachable}, that is, every finite sequence has a continuation belonging to $A$. Thus, any countable intersection of sets which are finitely determined and always reachable is non-empty.

These two conditions are purely combinatorial, and the topology
\index{topology}%
 has been removed.  The sequence space obeys the Baire category theorem.
\index{Baire category theorem}%

For measure, we have in fact nothing new.  Regarding an infinite binary sequence as the base 2 ``decimal'' expansion of a real number, we have a bijection between the sequence space $X$ and the unit interval, apart from a null set (the rationals with 2-power denominators have two represerntations); the bijection transforms the coin-tossing measure to Lebesgue measure on the interval. However, the topology of the two metric spaces is quite different. We can't expect to make infinitely many cuts in the unit interval without changing its topology drastically.

A property which holds for a large set is regarded as being a ``typical'' or ``generic'' property of elements of $X$, that is, almost all elements of $X$ have this property. Here is a case where measure and category give different views of what the typical element looks like:

$\bullet$\  The Law of Large Numbers
\index{Law of Large Numbers}%
 asserts that amost all binary sequences (in the sense of measure, that is, the complement of a null set) have limiting density 1/2.

$\bullet$\  On the other hand, almost all sequences (in the sense of category) have the property that they do not have a density; indeed, the $\lim$ $\inf$ of the density of the prefixes is 0 and the $\lim$ $\sup$ is 1. (To prove this, show that, for any $n$, the set of sequences having a prefix of length at least $n$ with density less than $\frac{1}{n}$ is obviously finitely determined, and is always reachable by simply appending enough zeros to the given finite sequence.)

But there are many cases where they agree. For example, call a sequence \emph{universal} if it contains every finite sequence of zeros and ones as a consecutive subsequence. To see that universal sequences exist simply concatenate all finite sequences.  The category argument to seeing that universality is a ``large'' property, so that almost every sequence (in either sense) is universal, runs as follows.   There are only countably many finite sequences, so it suffices to show that the set of infinite sequences containing a given finite sequence s is finitely determined and always reachable.  For the first, if $s$ occurs in a given sequence, then just choose a prefix containing $s$.   For the second, given any prefix, append $s$ to it, and every continuation will contain $s$.

Here is another example. Take a countably infinite set, say the natural numbers, and write down all the 2-element subsets in a sequence. Then there is a natural bijective correspondence between graphs on this vertex set, and our sequence space $X$: the $n$th pair of vertices is joined by an edge if the $n$th term of the sequence is 1, and not joined it it is 0. In the probabilistic model, we are choosing a random graph,
\index{graph ! random}%
 by tossing a fair coin once for each pair of vertices to decide whether to put an edge or not. What do almost all graphs look like?  Well its $\mathfrak{R}$, so that almost all graphs are isomorphic to $\mathfrak{R}$.

Before turning to $\mathfrak{R}$,  we give a generalisation of the binary sequence space.

Let $T$ be a countable rooted tree. Its nodes come on levels indexed by the natural numbers (including zero). The root is the unique node on level 0; each node other than the root has a unique predecessor on the preceding level, and at least one but at most countably many successors on the following level. A branch of $T$ is an infinite path, starting at the root and containing one node on each level. (So if a branch contains a node, it must contain its predecessor). Let $X(T)$ be the set of all branches of $T$.

If each node has two successors, we have the infinite binary tree; labelling them by 0 and 1, we have a bijection between the set of branches and the set of infinite binary sequences.

If the number of successors of any node is finite, then we may define a probability measure on $X$ in a simple way. For any node $v$, the measure of the set of branches containing $v$ is the reciprocal of the product of the numbers of successors of nodes on the path from the root to $v$. We then extend to all measurable sets in the standard way. More intuitively, we take a random walk where, arriving at a node, we then proceed to one of its successors, each being equally likely.  There are other measures which could be used.

The category model is simpler, and does not require that each node has only finitely many successors. We define the distance between two branches to be $\frac{1}{n}$ if they first diverge at level $n$; this makes $X$ a complete metric space, so the Baire category theorem applies. Just as in the sequence case, it is purely combinatorial.

A set $A$ of branches is \emph{open} if and only if it is \emph{finitely determined}, that is, for any branch $x$ in $A$, there is a node $v$ on $x$ such that all branches containing $v$ lie in $A$. And $A$ is \emph{dense} if and only if it is \emph{finitely determined}, that is, any node lies on a branch in $A$. So our large sets are again those which contain countable intersections of finitely determined, always reachable sets.

We shall return to the concepts of finitely determined and always reachable in the next section.

\smallskip

A measure on the random graph
\index{graph ! random}%
 is a function from the countable set of definable vertex subsets to $[0, 1]$.  The finitely additive probability measures on the definable subsets of the random graph which are invariant under the graph automorphism group
\index{group ! automorphism}%
 have been shown to be integrals of Bernoulli measures
\index{Bernoulli measure}%
arising from the coin-flipping model of random graph construction~\cite{albert}.  (Furthermore a random system can be formed from, the Bernoulli (left) shift map $T$,
\index{Bernoulli shift}%
 on the space $X = \{0, 1\}^{\mathbb{Z}} \equiv 2^{\mathbb{Z}}$ of infinite 0-1 sequences (or equivalently, the space of all sets of integers), together with uniform product probability measure, where $T(x_n)_{n \in \mathbb{Z}} := (x_{n+1})_{n \in \mathbb{Z}}$ is invertible and both it and its inverse are measurable.  The Bernoulli shift is an example of a \emph{strongly mixing system}
\index{strongly mixing system}%
 which for measurable sets $E$ and $F$ for which $\mu(T^n E \cap F) \to \mu(E) \mu(F)$ as $n \to \infty$, that is the shifted sets become asymptotically independent of the unshifted sets~\cite{tao}.)  This was extended to structures which do not have the independence property in~\cite{ensley}.  An alternative random graph measure is studied in~\cite[p.~112]{cam6}, where it is used to indicate the probability that a sequence in the random graph is isomorphic to a given finite subgraph.  The two types of measure are contrasted in~\cite{albert} together with a study of other graph measures, such as those on triangle-free graphs.
\index{graph ! triangle-free}%

Petrov and Vershik~\cite{petrovv}
\index{Petrov, F. V.}%
\index{Vershik, A. M.}%
study invariant measures
\index{measure ! invariant}%
 on the set of universal countable graphs.  The infinite symmetric group of a fixed countable vertex set $V$ acts naturally on the set $\mathcal{G}$ of simple graphs.  They equip the set $\mathcal{G}_V$ of all graphs with the \emph{weak topology} whose base is formed by the collections of sets of graphs that have a given induced graph structure on a given finite vertex set.  This permits consideration of Borel sets, $\sigma$-fields and probability measures on $\mathcal{G}_V$.  They prove the existence of $\Sym(V)$-invariant probabilistic Borel ergodic (invariant) measures
\index{Borel measure}%
  on the set of universal graphs and universal Henson graphs 
\index{graph ! Henson}%
with the set of vertices $\mathbb{Q}$ and with shift invariant graph structure, and then extend the graph structure on the whole line $\mathbb{R}$.  Taking $V = \mathbb{N}$, the space $\mathcal{G}_V$ supports $\Sym(\mathbb{N})$, with orbits a class of isomorphic graphs and the stabilizer of a given graph, the automorphism group of the graph.  This action extends to one on invariant Borel probability measures on the spaces of graphs and matrices, and is transitive on universal graphs.  A `random graph'
\index{graph ! random}%
 in a given category is then a $\Sym(\mathbb{N})$-invariant Borel probability measure on the set of graphs that is concentrated on the set of universal graphs of this category.  The Erd\H{o}s-R\'enyi random graphs are the Bernoulli measures on the space of adjacency matrices
\index{adjacency matrix}%
  of Erd\H{o}s-R\'enyi 
\index{Erd\H{o}s, P.}%
\index{R\'enyi, A.}%
 random graphs,  with the distribution $(p, 1? p), 0 < p < 1$, for each entry.  Our case of $p = \frac{1}{2}$ this Bernoulli measure is the weak limit of the uniform measures on the sets of finite graphs with $n$ vertices as $n \to \infty$.   There are uncountably many invariant ergodic measures on the set of $K_n$-free graphs for $n > 2$.
 
They define two new concepts, that of \emph{measurable universal graphs} and \emph{topological universal graphs}.
\index{graph ! topological universal}%
 Vertices of the latter lie in a Polish space.
\index{Polish space}%
  They build a graph on the additive group $V = \mathbb{R}$, $E = \{(x, y) : |x - y| \in Z \subset (0, + \infty) \}$ and prove the following
\begin{theorem}
There is a universal topological graph (resp., universal topological triangle-free graph) with vertex set the additive group $\mathbb{R}$ and with graph structure which is invariant under the additive action of the group $\mathbb{R}$ on itself.
\end{theorem}
In summary, the route they take is universal Borel graph with measures $\to$ topologically universal graph ($\to$ homogeneous topologically universal graph for $\mathfrak{R}$ and $K_3$-free cases) $\to$ randomization in vertices $\to$ invariant measures on the set of countable universal graphs $\to$ randomization in edges $\to$ the list of all invariant measures on the set of countable universal (or $K_n$-free universal) graphs.

\smallskip

With an eye to future research we record some of Terence Tao's
\index{Tao, Terence}%
remarks on structure and randomness in ergodic theory and graph theory~\cite{tao}.  Tao notes the dichotomy between these two properties in 

\begin{itemize}
\item Combinatorial number theory, (finding patterns in unstructured dense sets (or colourings) of integers);
\item Ergodic theory and in particular multiple recurrence theory, (where a discrete dynamical system acts on probability spaces
\index{polynomial algebra}%
 yielding patterns in positive-measure sets (or more specifically, measure-preserving actions of $\mathbb{Z}$) or probability preserving systems;
\item Graph theory, (specifically finding patterns in large unstructured dense graphs); and
\item Ergodic graph theory, (finding patterns within a measure-preserving action of the infinite permutation group,
\index{group ! permutation ! infinite}%
 as one way in which to study infinite limits of graphs.
\end{itemize}

The first two on this list are discrete and finitary and connected by a Cayley graph construction, whilst the second two are continuous and infinitary.  Tao writes

``On the other hand, we have some very rigorous connections between combinatorial number theory and ergodic theory, and also (more recently) between graph theory and ergodic graph theory, basically by the procedure of viewing the infinitary continuous setting as a limit of the finitary discrete setting. These two connections go by the names of the \emph{Furstenberg correspondence principle} and the \emph{graph correspondence principle} respectively. These principles allow one to tap the power of the infinitary world (for instance, the ability to take limits and perform completions or closures of objects) in order to establish results in the finitary world, or at least to take the intuition gained in the infinitary world and transfer it to a finitary setting. Conversely, the finitary world provides an excellent model setting to refine one's understanding of infinitary objects, for instance by establishing quantitative analogues of `soft' results obtained in an infinitary manner.''

These comments, further discussed in~\cite{tao}, echo part of the philosophy behind our project. 

\smallskip

As we often use the concept of Cayley graph
\index{graph ! Cayley}%
in the main text, we mention two ways of putting metric spaces
\index{metric space}%
on them.  For a finitely-generated group $G$ with identity element $e$ define $\forall g \in G$ the \emph{length} $l_S(g) := \min(n)$ such that $g = s_1 s_2 \ldots s_n$ where the $s_i$ belong to a finite set $S$ of generators.  The Cayley graph
\index{graph ! Cayley}%
$\Cay(G, S)$ of $G$ can be turned into a metric space by defining the distance function $d_S: G \times G \to \mathbb{R}^{+}$ where $d_S(g_1, g_2) = l_S(g_1^{-1} g_2)$.  

Each edge of the Cayley graph can be turned into a metric space isomorphic to the unit interval $[0, 1]$ such that the left action of $G$ gives isometries between edges.  Thus $\Cay(G, S)$ can be turned into an arc-connected metric space, and the inclusion $G \subset \Cay(G, S)$ is an isometry~\cite[p.300]{ghysha}.

\smallskip

We end by noting that whilst there is a natural topology on the symmetric group, there is no natural measure.

\smallskip

A standard introduction to measure and category is~\cite{oxt}.


\section{Fra\"{\i}ss\'e's Theory of Relational Structures}
\label{TheoryofRelationalStructures}
Recall that a structure is \emph{homogeneous}
\index{structure ! homogeneous}%
if every isomorphism between finite substructures extends to an automorphism of the entire structure.  Also recall that the \emph{theory}
\index{theory}%
of a structure $M$, $\Th(M)$, is the set of all first-order sentences
\index{first-order sentence}%
 which are true in $M$.  An $\aleph_0$-categorical relational structure $M$
\index{aleph@$\aleph_0$-categorical}%
is homogeneous if and only if, modulo $\Th(M)$, every formula with parameters in $M$ is equivalent to a quantifier-free formula~\cite{hod1}.  Fra\"{\i}ss\'e
\index{Fra\"{\i}ss\'e, R.}%
produced a theory~\cite{frai} that allows construction of a great variety of unique homogeneous structures.  The multicoloured random graphs
\index{graph ! random ! $m$-coloured}%
 introduced in the next section are examples.  The basic concept in Fra\"{\i}ss\'e's work 
\index{Fra\"{\i}ss\'e, R.}%
is the \emph{age}
\index{age}%
Age$(M)$ of a relational structure $M$, which is the class of all finite structures (over the same logical language) which are embeddable in $M$ as induced substructures.  

The main result applicable to random graphs is
\begin{theorem}[Fra\"{\i}ss\'e's Theorem]
\index{Fra\"{\i}ss\'e's Theorem}%
Let $\mathit{L}$ be a purely relational first-order language.
\index{first-order language}%
  A class Age$(M)$ of finite $\mathit{L}$-structures is the age
\index{age}%
 of a countable homogeneous relational structure
\index{relational structure}%
  $M$ if and only if it satisfies the following properties:
\begin{itemize}
\item[(A1)]  Age$(M)$ is closed under isomorphism,
\item[(A2)]  Age$(M)$ contains at most a countable number of non isomorphic structures,
\item[(A3)]  Age$(M)$ has the \emph{Hereditary Property},
\index{hereditary property}%
that is it is closed under taking induced substructures,
\item[(A4)]  Age$(M)$ has the \emph{Joint Embedding Property},
\index{joint embedding property}%
that is, given $A, B \in$ Age$(M)$ then there are a structure $C \in$ Age$(M)$ and embeddings $f_1: A \rightarrow C$ and embeddings $g_1: B \rightarrow C$,
\item[(A5)]  Age$(M)$ has the \emph{Amalgamation Property},
\index{amalgamation property}%
that is, given $A, B_1,$ $B_2 \in$ Age$(M)$ and embeddings $f_i: A \rightarrow B_i$ for $i = 1,2$, there exists $C \in$ Age$(M)$ and embeddings $g_i: B_i \rightarrow C$ for $i = 1,2$ such that $f_1 g_1 = f_2 g_2$,
\end{itemize}
then the structure $M$ is unique up to isomorphism.
\end{theorem}

We refer to the amalgamation class of structures that satisfy this theorem as the \emph{Fra\"{\i}ss\'e Class}
\index{Fra\"{\i}ss\'e class}%
(which we denote $\mathcal{A}$), and the homogeneous structure $M$ guaranteed by this theorem as the \emph{Fra\"{\i}ss\'e limit}.
\index{Fra\"{\i}ss\'e limit}%

If the language $\mathit{L}$ is finite then there are only countably many isomorphism classes of finite $\mathit{L}$-structures.  If Age$(M)$ is defined by a set of universal $\mathit{L}$-sentences then Age$(M)$ has the Hereditary Property.  Much of the effort to prove a structure is homogeneous with the required age is consumed with verification of (A5).

Properties of relational structures
\index{relational structure}%
 such as their uniqueness and universality can be demonstrated using the back-and-forth method;
\index{back-and-forth method}%
 in fact this is how Fra\"{\i}ss\'e's Theorem
\index{Fra\"{\i}ss\'e's Theorem}%
  is proved.  A related object, the homogeneous universal directed graph, is the Fra\"{\i}ss\'e limit
\index{Fra\"{\i}ss\'e limit}%
   for the class of all
directed graphs; its existence can be proved by using the extension
property, which implies both its universality and its
homogeneity~\cite{hubicka}.  

An important strengthening of the amalgamation property is the following,

\begin{itemize}
\item[(A5')]  Age$(M)$ has the \emph{Strong Amalgamation Property},
\index{amalgamation property ! strong}%
that is, given $A, B_1, B_2 \in \mathcal{A}$ and embeddings $f_i: A \rightarrow B_i$ for $i = 1,2$, there exists $C \in$ Age$(M)$ and embeddings $g_i: B_i \rightarrow C$ for $i = 1,2$ such that $f_1 g_1 = f_2 g_2$ and, if $b_1 \in B_1, b_2 \in B_2$ and $b_1 g_1 = b_2 g_2$, then there exists $a \in A$ such that $b_i = a f_i$ for $i = 1,2$.
\end{itemize}

The difference between amalgamation and strong amalgamation is that the former allows the possibility that when $B_1$ and $B_2$ are glued over a common substructure $A$, some additional points also become identified, whereas the stronger version guarantees that amalgamation can be performed without such extra identifications.  Therefore, the region labelled $A''$ in Figure~\ref{TAP} can be chosen to be empty.  Equivalently, in the Fra\"{\i}ss\'e Limit $M$,
\index{Fra\"{\i}ss\'e limit}%
 the automorphism group
\index{group ! automorphism}%
  fixing any finite set of points has no further fixed points.  There are structures that are homogeneous but for which strong amalgamation fails, such as ``treelike objects''.   The strong amalgamation property is useful for producing constructions; one of its known properties~\cite{cam2d}~\cite{cameron} is the following theorem:

\begin{theorem}
Let $M$ be a countable homogeneous structure.  Then the following conditions are equivalent:
\begin{itemize}
\item[(a)]  The age Age$(M)$ has the strong amalgamation property;
\item[(b)]  the stabilizer in $\Aut(M)$ of any finite number of points has no additional fixed points;
\item[(c)]  the stabilizer in $\Aut(M)$ of any finite number of points has no additional finite orbits;
\item[(d)]  $M \backslash \{x\} \cong M$ for any point $x$ of $M$.
\end{itemize}
\end{theorem}

This proposition implies, for example, that the orbits of the stabilizer of a tuple are the infinite open intervals in $\mathbb{Q}$ lying in between the points of the tuple of stabilized points.  (For every structure $M$ obeying this proposition there is a homogeneous structure $M'$ with age
\index{age}%
 $Age(M')$ consists of labelled members of $Age(M)$, and $\Aut(M') \leq \Aut(\mathbb{Q})$.)

\bigskip

A homogeneous structure has the strong amalgamation property if and only if deleting any finite number of points of its points gives a structure isomorphic to the original.  This is the analogue of the  coherence condition in Alexander Gnedin's
\index{Gnedin, A.}%
 coherence theory of regeneration in random combinatorial structures~\cite{gnedin},
\index{regeneration in random structures}%
  in which removing something from a combinatorial object leaves a copy of what is there before.  A single structure cannot have this property, nor can an infinite chain of structures in any interesting way, so Gnedin considered probability distributions on partitions of a positive integer.
 
A collection of probability distributions on the set of all partitions $\Pi_n$ of the integer $n$ ($\forall n$) is (i) \emph{coherent} if subtracting one from a random part of an element of $\Pi_n$ gives a random element of $\Pi_{n-1}$; (ii) \emph{regenerative} if removing an entire part of an element of $\Pi_n$ gives a random element of the set of partitions of what remains.

There are some distributions to be specified here.  For coherence, we could choose the part uniformly, or with probability proportional to its size.

Distributions are described uniformly in two steps, which can be illustrated with an example.  Firstly divide the unit interval randomly then colour each interval with a different colour (\emph{Kingman's paintbox}).  Secondly, choose $n$ independent random points from the unit interval according to some possibly new and unrelated distribution.  Then we have a partition of $n$ given by the number of points of each colour in the selection.

This gives rise to the distributions which are the cycle lengths of a random permutation.  That is, the probability of an element of $\Pi_n$ is the probability that a random element of the symmetric group $\Sym(n)$ has those cycle lengths.

Gnedin
\index{Gnedin, A.}%
 has found explicit formulae for the distribution on partitions in many cases.

Now, the random graph is produced by a very simple probability distribution, and has very strong regeneration properties: remove a finite set of vertices or edges, or a finite set of vertices and their common neighbours, or the edges of a locally finite subgraph, and what we obtain is still isomorphic to the random graph.

We have already mentioned the work of Petrov and Vershik
\index{Petrov, F. V.}%
\index{Vershik, A. M.}%
 who constructed~\cite{petrovv} an exchangeable measure on countable graphs which is concentrated on Henson's homogeneous universal triangle-free graph
\index{graph ! triangle-free}%
\index{graph ! Henson}%
  (that is, the random graph is isomorphic to Henson's graph with probability 1), and a generalisation of this due to Ackerman, Patel and Freer,
\index{Ackerman, N.}%
\index{Patel, R.}%
\index{Freer, C.}%
   who gave a necessary and sufficient condition on a structure $M$ for there to be an exchangeable measure concentrated on $M$.  Their condition is precisely that the stabiliser of a finite set in $\Aut(M)$ has only infinite orbits on the remaining points.

Now there is an even closer connection with what Gnedin was talking about. The method used by these five researchers was to define a structure on the real numbers which is almost of the kind desired (thus, in the Petrov-Vershik case, the measure of the set of triangles is zero), and then obtain the countable structure by choosing countably many real numbers independently and taking the induced structure on them.

We ask an open question.  What is the deeper connection here with Gnedin's work on probability distributions?

\smallskip

A countable homogeneous structure $M$ has the \emph{strong embedding property}
\index{structure ! strong embedding property}%
 if for all $A \in$ Age$(M)$ and $x \in A$ and embeddings $e : A \backslash {x} \to M$ there are infinitely many different extensions of $e$ to embeddings of $A$ into $M$.  The following result is from~\cite{elzahar}:

\begin{theorem}
Let $M$ be a countable homogeneous structure.  Then the following conditions are equivalent:
\begin{itemize}
\item[(a)]  $M$ has the strong amalgamation property;
\item[(b)]  $M$ has the strong embedding property;
\item[(c)]  $M$ is strongly inexhaustible.
\end{itemize}
\end{theorem}

\clearpage

\begin{figure}[!h]
$$\xymatrix{
& {\bullet} \ar@{-}[dddrrr]^{B_1} \ar@{-}[dl]\\
{\bullet} \ar@{-}[dddrrr]\\
& {}\\
&& {}  \ar@{-}[rr]_{A} &&  {\bullet} {}  \ar@{-}[dl]\\
&&& {\bullet}
}$$
$$\xymatrix{ 
&&& {\bullet}  \ar@{-}[dr]\\
&&&& {\bullet}  \ar@{-}[ddll]\\
& {\bullet} \ar@{-}[rr]_{A}  \ar@{-}[dr] \ar@{-}[uurr]^{B_2} && {}\\
&& {\bullet}
}$$\end{figure}

\begin{figure}[!h]
$$\xymatrix{
&& {\bullet} \ar@{-}[dl]  \ar@{-}[ddrr]^{C}\\
& {\bullet} \ar@{-}[dddrrr] \ar@{-}[dl]\\
{\bullet} \ar@{-}[dddrrr] &&&& {\bullet} \ar@{-}[dr]\\
& {} &&&& {\bullet}\\
&& {\bullet} \ar@{-}[uurr] \ar@{-}[rr]_{A}^{A''} &&  {\bullet} {}  \ar@{-}[dl] \ar@{-}[ur]\\
&&& {\bullet}
}$$\caption{The Amalgamation Property}
\label{TAP}
\end{figure}
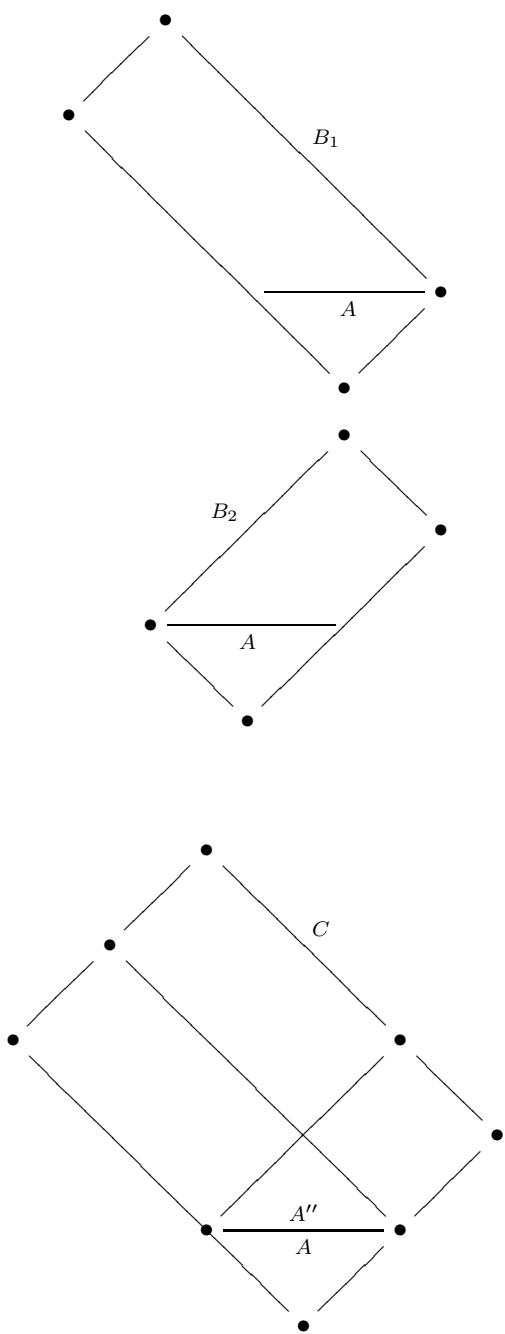

\clearpage

The Strong Amalgamation Property has been further strengthened by Macpherson and Tent~\cite{mactent}, 
\index{Macpherson, H. D.}%
\index{Tent, K.}%
as follows.  A Fra\"{\i}ss\'e class
\index{Fra\"{\i}ss\'e class}%
 $\mathcal{A}$ has the \emph{Free Amalgamation Property}
\index{amalgamation property ! free}%
if, whenever $B_1$ and $B_2$ are structures in $\mathcal{A}$ with a common substructure $A$, there is an amalgam $C$ of $B_1$ and $B_2$ such that
\begin{itemize}
\item[(i)]  the intersection of $B_1$ and $B_2$ in $C$ is precisely $A$ (that is, strong amalgamation);
\item[(ii)] every instance of a relation in $C$ is contained in either $B_1$ or $B_2$.
\end{itemize}
This holds for graphs: we can choose to make the amalgam so that there are 
no edges between $B_1 \backslash A$ and $B_2 \backslash A$.  They prove the following: 

\begin{theorem}
Let $\mathcal{A}$ be a nontrivial Fra\"{\i}ss\'e class
\index{Fra\"{\i}ss\'e class}%
 (that is, there are some non-trivial relations) with the free amalgamation property, and $G$ the automorphism group
\index{group ! automorphism}%
  of its Fra\"{\i}ss\'e limit.
\index{Fra\"{\i}ss\'e limit}%
 Then $G$ is simple.
\end{theorem}

The trivial case must be excluded, since then $G$ is the symmetric group, which is not simple. 

There is a generalization of homogeneity in relational structures,
\index{relational structure}%
 called \emph{pseudo-homogeneity}
\index{structure ! pseudo-homogeneous}%
whose study was initiated by Calais~\cite{calais}
\index{Calais, J. P.}%
and continued by R. Fra\"{\i}ss\'e
\index{Fra\"{\i}ss\'e, R.}%
who showed~\cite{frai} that if $\mathcal{C}$ is a class of countable structures closed under both isomorphisms and taking induced substructures, then a structure $M \in \mathcal{C}$ is \emph{universal pseudo-homogeneous} if there is a subclass of finite substructures $\mathcal{C'} \subset \mathcal{C}$ satisfying (A4), (A5), and $\mathcal{C'}$ and $\mathcal{C}$ satisfy the following \emph{cofinality condition}: for each finite $A \in \mathcal{C}$, there is a finite $B \in \mathcal{C'}$ so that $A \leq B$.  

The subclass $\mathcal{C}$ is called a \emph{pseudo-amalgamation class}~\cite{bona1}.  There is an infinite universal pseudo-homogeneous graph sharing some of the properties of $\mathfrak{R}$.  Taking $\mathcal{C}$ to be the class of all countable graphs and $\mathcal{C'}$ to be the class of all finite graphs, makes $\mathfrak{R}$ universal pseudo-homogeneous.

We mention just one more variation on the basic theory that has been studied.  Counting problems for orbits on sets and tuples of oligomorphic permutation groups are identical with those for unlabelled and labelled structures in so-called \emph{oligomorphic Fra\"{\i}ss\'e classes}
\index{Fra\"{\i}ss\'e class, oligomorphic}%
 of relational structures,
\index{relational structure}%
 which satisfy a stronger version of condition (A2):
\begin{itemize}
\item[(A2')]  Age$(M)$ contains only finitely many $n$-element structures up to isomorphism for all $n$,
\end{itemize}

These include large numbers of combinatorially important classes of structures, such as graphs, directed graphs, tournaments, partially ordered sets, $k$-edge-coloured graphs, and graphs with a fixed bipartition.  The condition certainly holds if the relational language has finitely many relations, for example, the random graph
\index{graph ! random}%
 or the Fra\"{\i}ss\'e limit
\index{Fra\"{\i}ss\'e limit}%
 $(\mathbb{Q}, <)$ of the class of finite totally ordered sets.
 
\bigskip

We mention that there is an unpublished version of Fra\"{\i}ss\'e's Theorem
\index{Fra\"{\i}ss\'e's Theorem}%
 in the language of Category Theory due to J. Covington.
\index{Covington, J.}%

\vspace{1cm}

The age
\index{age}%
 Age$(M)$ of a relational structure
\index{relational structure}%
  $M$ can be regarded as an example of two different theories:

either

(i)  a \emph{species}
\index{species}%
in the sense of Joyal~\cite{joyal}~\cite{bergeron}~\cite{joyal},
\index{Joyal, \'A.}%
which is a formal power series where the coefficients are combinatorial objects instead of numerical coefficients, though the former can be reduced to the latter,

or 

(ii)  as a tree which leads to the study of \emph{ubiquity}~\cite{cameron}:
\index{ubiquity}%
a countable homogeneous relational structure
\index{relational structure}%
 $M$ is \emph{ubiquitous in category}
\index{ubiquitous in category}%
if a residual set of structures younger than $M$ are isomorphic to $M$.

We develop the second of these theories, as it is used in the main text.  The nodes of Age$(M)$ thought of as a tree occur on `levels' indexed by non-negative integers.  The empty set lies at the root of tree.  All structures for which Age$(M)$ has point set $\{1, \ldots, n\}$ occur on tree level $n$.  Each node on level $n + 1$ is adjacent to the unique level-$n$ node obtained by deleting the point $n + 1$, that is the induced substructure on $\{1, \ldots, n\}$.  An infinite path in the tree corresponds to a structure $S$ on $\mathbb{N}$ whose induced substructure on $\{1, \ldots, n\}$ belongs to Age$(M)$ for all $n$.  By property (A3) above, every finite substructure of $S$ lies in Age$(M)$, so $S$ is \emph{younger than} $M$, that is Age$(S) \subseteq$ Age$(M)$.

Conversely, if the points of a countable structure $S$ younger than $M$ are indexed by the natural numbers, then the induced substructure of $S$ on $\{1, \ldots, n\}$ is a node of the tree for each $n$, and nodes on consecutive levels are adjacent.  So $S$ is represented by a path in the tree.  Different enumerations of the points of $S$ give different paths.

So  the infinite paths from the root represent the (labelled) structures younger than $M$.

As explained in the previous section, a metric can be defined on the set $\mathcal{Y}(M)$ of paths by the rule that two paths which agree up to level $n$ and diverge thereafter should have distance $\frac{1}{2^n}$.  Any decreasing function can be chosen here, and this one gives the space a diameter of 1.  With this metric $\mathcal{Y}(M)$ is a complete metric space.

A set $\mathcal{S}$ of paths (that is structures younger than $M$) is \emph{open} if and only if it is \emph{finitely determined}, meaning that for any $S \in \mathcal{S}$, there exists $n$ such that any structure in $\mathcal{Y}(M)$ which induces the same structure on $\{1, \ldots, n\}$ as $S$, belongs to $\mathcal{S}$.  For the ball of radius $\frac{1}{2^n}$ containing $S$ consists of all structures agreeing with $S$ on $\{1, \ldots, n\}$.

A set $\mathcal{S}$ is \emph{dense} if and only if it is \emph{always reachable}, meaning that any finite structure is induced by some member of $\mathcal{S}$.  In other words, $\mathcal{S}$ meets all open balls,
\index{open ball}%
whose definition in given in a later section.

We will continue with the discussion of ages and ubiquity after defining random graphs in the next section.

The homogeneous structures among various classes that have been determined include posets by Schmerl~\cite{schmerl},
\index{Schmerl, J. H.}%
 tournaments by Lachlan~\cite{lachlan}
\index{Lachlan, A. H.}%
 permutations by Cameron~\cite{cam12}
\index{Cameron, P. J.}%
and digraphs by Cherlin~\cite{cherlinb}.
\index{Cherlin, G. L.}%

There are two noteworthy recent developments in the theory of homogeneous structures and Fra\"{\i}ss\'e classes.  \index{Fra\"{\i}ss\'e class}%
The first one makes connections with Ramsey theory;
\index{Ramsey theory}%
in particular, if the age
\index{age}%
 of a structure with an explicit or implicit order is a Ramsey class then it is a Fra\"{\i}ss\'e class~\cite{hubicka}.  The second one makes connections with topological dynamics,
\index{topological dynamics}%
 for example in~\cite{kechris}~\cite{kechrisros} it is shown that the automorphism groups
\index{group ! automorphism}%
  of homogeneous structures have very strong amenability properties.  We refer to~\cite[\S3.4]{cam23a} for a brief account of these.

\smallskip

Fra\"{\i}ss\'e's
\index{Fra\"{\i}ss\'e, R.}%
 original theory~\cite{frai} has been greatly extended in a number of directions, a sample of which are~\cite{camnes}~\cite{camloc}~\cite{covington}~\cite{drostegir}~\cite{elzahar} \cite{higmansc} \cite{hod1} \cite{hrushov} \cite{hubicka}~\cite{jonsson} \cite{jonsson1}~\cite{kechris}~\cite{pouzet}.  Hrushovski's
\index{Hrushovski, E.}%
theory in essence requires only some instances of amalgamation to be possible, enabling the construction of many new types of structures and permutation groups; a detailed discussion of it is given in Wagner's
\index{Wagner, F. O.}%
expository article in~\cite{kayem}.  Homogeneous structures arise in many areas of mathematics; for example the connection with Ramsey theory
\index{Ramsey theory}%
 is studied by Hubi\v{c}ka and Ne\v{s}et\v{r}il in~\cite{hubicka}~\cite{nesetril1},
\index{Ne\v{s}et\v{r}il, J.}%
\index{Hubi\v{c}ka, J.}%
  and in~\cite{hubickanes1} they prove the existence of countably universal structures for every class of structures determined by finitely many forbidden homomorphisms.  Arbitrary relational structures on an infinite groundset satisfying a finiteness hypothesis stating that there are only finitely many isomorphism types of substructures of size a certain infinite cardinality condition is studied in~\cite{gibson}.

Several parts of model theory
\index{model theory}%
 are increasingly becoming standard tools in various areas of mathematics.  For example, for constructions of special Banach spaces in Banach space theory
\index{Banach space}%
 and in Metric space theory
\index{metric space}%
 and oscillation of such spaces.  There are papers of Uspenskij
\index{Uspenskij, V.}%
 and others working in the theory of large group actions.  
\index{group ! action}%

Especially noteworthy is the fundamental paper by Kechris, Pestov and Todorcevic~\cite{kechris}
\index{Kechris, A. S.}%
\index{Pestov, V.}%
\index{Todorcevic, S.}%
 because it opened up a new research direction, and we give a brief introduction to one of the main results.

A class $\mathbf{C}$ has the \emph{Ramsey property} if for any two structures $A, B \in \mathbf{C}$ and any positive integer $m$ there exists a structure $C \in \mathbf{C}$ such that if the substructures of $C$ isomorphic to $A$ are coloured with $m$ colours, then there is a substructure isomorphic to $B$ such that all its A-substructures have the same colour. If $\mathbf{C}$ consists of sets with no structure then this property reduces to the classical form of Ramsey's theorem.
\index{Ramsey theorem}%

Recall that the automorphism group of a countable structure $\mathcal{M}$ is a topological group
\index{group ! topological}%
 with the topology of pointwise convergence. A \emph{flow}
\index{flow}%
 for a topological group
\index{group ! topological}%
 $G$ is a continuous action of $G$ on a compact Hausdorff space.
\index{Hausdorff space}%
 The group $G$ is \emph{extremely amenable}
\index{group ! extremely amenable}%
 if every flow for $G$ has a fixed point. More generally, any group has a \emph{universal minimal flow}: it is extremely amenable if the universal minimal flow is a space with a single point.

The theorem of Kechris, Pestov and Todorcevic
\index{Kechris, A. S.}%
\index{Pestov, V.}%
\index{Todorcevic, S.}%
 mentioned above states:

The age of a countable homogeneous structure with a total order is a Ramsey class if and only if its automorphism group is extremely amenable.

This was an unexpected confluence of homogeneous structures, Ramsey theory and topological dynamics. 
\index{topological dynamics}%
The countable dense linear order
\index{linear order}%
 without endpoints (the rational numbers as ordered set) are an example of a homogeneous structure with this property.

Lionel Nguyen Van Th\'e, together with Yonatan Gutman and Todor Tsankov
\index{Nguyen Van Th\'e, L.}%
\index{Tsankov, T.}%
\index{Gutman, Y.}%
 showed~\cite{nguyen} that a Fra\"{\i}ss\'e class
\index{Fra\"{\i}ss\'e class}%
 of structures has an expansion (obtained by adding extra relations) with the Ramsey property if and only if the automorphism group
\index{group ! automorphism}%
 of its Fra\"{\i}ss\'e limit
\index{Fra\"{\i}ss\'e limit}%
 has an extremely amenable subgroup with precompact quotient.

The Glasner-Weiss theorem
\index{Glasner-Weiss theorem}%
\index{Glasner, E.}%
\index{Weiss, B.}%
 states that the universal minimal flow for the symmetric group of countable degree is the action on the class of linear orders
\index{linear order}%
 of the countable set.  Tsankov
\index{Tsankov, T.}%
 classified the minimal flows for this group.  The arguments required for this classification were produced by Claude Frasnay
\index{Frasnay, C.}%
 in the 1960s, in a totally different context.

\bigskip

More topics were discussed at the 2nd Workshop on Homogeneous Structures held in Prague in July 2012.

\clearpage

\section{The Random Graph $\mathfrak{R}$}
\label{TheRandomGraph}
\index{graph ! random}%
In 1963, Erd\H{o}s and R\'enyi~\cite{er} proved the following remarkable theorem.
\index{Erd\H{o}s, P.}%
\index{R\'enyi, A.}%
\begin{theorem}
There is a countable graph $\mathfrak{R}$ with the following property:
Suppose that a graph on a given countable vertex set is chosen at random by
selecting edges independently with probability $\frac{1}{2}$ from the set of
unordered pairs of vertices. Then, with probability~$1$, the resulting random
graph is isomorphic to $\mathfrak{R}$.
\end{theorem}

This \emph{countable random graph} $\mathfrak{R}$ is the prototype and
motivation for what we have described in this book. Accordingly, we give a
brief introduction to this graph and its properties in this Appendix, referring
to~\cite{cameron} for further details.

The proof of the Erd\H{o}s--R\'enyi theorem
is surprisingly easy; we now sketch it. It depends on the following
property $(*)$ which may or may not be satisfied by a given graph:
\begin{quote}
For any finite disjoint sets $U$, $V$ of vertices, there is a vertex $z$ 
which is joined to every vertex in $U$ and to none in $V$.
\end{quote}
Now the Erd\H{o}s--R\'enyi theorem follows immediately from the following two
facts, whose proofs we sketch:

Fact 1: With probability~$1$, the random countable graph satisfies
$(*)$.

Fact 2: Any two countable graphs satisfying $(*)$ are
isomorphic.

Proof of Fact 1: Since there are only a countable number of 
disjoint pairs $(U,V)$ of finite sets, and a countable union of null sets
is null, it suffices to show that the failure of $(*)$ for a particular
choice of $U$ and $V$ has probability~$0$. But if $z_1,z_2,\ldots,z_N$ are
vertices outside $S=U\cup V$, the probability that none of $z_1,\ldots,z_N$
is ``correctly'' joined to $S$ is $(1-1/2^|S|)^N$, which tends to~$0$ as
$N\to\infty$, as required.

Proof of Fact 2: Let $X_1$ and $X_2$ be countable graphs satisfying
condition $(*)$; assume that we have an enumeration of the vertices of each
of these graphs.  We construct an isomorphism between them by the method of
back-and-forth,
\index{back-and-forth method}%
 as the union of a sequence of finite partial isomorphisms.
Start with the empty partial isomorphism $f_0$. At an odd-numbered stage, let
$x$ be the first vertex of $X_1$ not in the domain of the current partial
isomorphism $f_n$; use property $(*)$ of the graph $X_2$ to find a vertex
$y$ such that the extension $f_{n+1}$ of $f_n$ which maps $x$ to $y$ is a
partial isomorphism. At an even-numbered stage, let $y$ be the first vertex
of $X_2$ not in the range of $f_n$, and use property $(*)$ for $X_1$ to find
a suitable pre-image $x$.

The proof gives further information about the random graph.

\begin{theorem}
\begin{enumerate}
\item $\mathfrak{R}$ is \emph{universal}: that is, every finite or countable
graph can be embedded as an induced subgraph of $\mathfrak{R}$.
\item $\mathfrak{R}$ is \emph{homogeneous}: that is, every isomorphism between
finite induced subgraphs of $\mathfrak{R}$ can be extended to an automorphism
of $\mathfrak{R}$.
\end{enumerate}
\end{theorem}

\begin{proof} For both parts, we use a modification of the 
``back-and-forth'' method. For the first, note that going in the forward
direction only requires that $(*)$ holds in $X_2$, and if we take all our
steps in this direction we construct an embedding (rather than an isomorphism)
from an arbitrary graph to $X_2$. For the second, given a finite partial
isomorphism $g$ of $\mathcal{R}$, take $X_1=X_2=\mathcal{R}$, and modify
the back-and-forth by starting with $g$ rather than the empty map.
\end{proof}

\medskip

The Erd\H{o}s--R\'enyi theorem gives a non-constructive existence proof for the
graph $\mathfrak{R}$: an event with probability~$1$ is certainly non-empty!
In order to give an explicit construction, all we need to do is to verify that
the constructed graph has property $(*)$. We give three such constructions.
The first was given  by Rado~\cite{rado} in 1964.

\head{First construction} The vertex set is the set $\mathbb{N}$ of
natural numbers (including zero). Given two distinct vertices $x$ and $y$,
suppose (without loss of generality) that $x<y$; then join $x$ to $y$ if and
only if the $x$th digit in the base~$2$ representation of $y$ is equal to~$1$.

Let $U$ and $V$ be finite disjoint sets of vertices. By enlarging $U$ if
necessary, we can assume that $\max(U)>\max(V)$. Now it is straightforward
to show that $z=\sum_{u\in U}2^u$ satisfies the requirements of $(*)$.

\head{Second construction} The vertex set is the set of primes congruent
to $1$ mod~$4$. Join vertices $p$ and $q$ if and only if $p$ is a quadratic
residue mod~$q$. By the Law of Quadratic Reciprocity, this is a symmetric
adjacency relation.

Let $U$ and $V$ be finite disjoint sets of vertices. Choose a quadratic
non-residue $a_v$ for each $v\in V$. By the Chinese Remainder Theorem, there
is a solution to the simultaneous congruences
\begin{center}
\begin{tabular}{ll}
$x \equiv 1~\hbox{mod}~u$ & for all $u\in U$;\\
$x \equiv a_v~\hbox{mod}~v$ & for all $v\in V$;\\
$x \equiv 1~\hbox{mod}~4$;&\\
\end{tabular}
\end{center}
the solution is unique modulo $N=4\prod_{u\in U}u\prod_{v\in V}v$. By
Dirichlet's Theorem, this arithmetic theorem contains a prime, which gives
the required value of $z$.

\head{Third construction} By the downward L\"owenheim--Skolem theorem of
first-order logic,
\index{first-order logic}%
 if the Zermelo-Fraenkel axioms for set theory are
consistent, then they have a countable model; that is, there is a countable
set $X$, with a binary relation $\in$ on it, for which the axioms are satisfied.
We take the vertex set to be $X$, and join $x$ to $y$ if and only if either
$x\in y$ or $y\in x$ (in other words, adjacency is symmetrised membership).

Let $U$ and $V$ be finite disjoint subsets of $X$. Using the Pairing, Union
and Foundation axioms, it is easy to construct an element $z$ which contains
all the elements of $U$ and none of the elements of $V$.

Note that only the three axioms mentioned are used in this proof. This means
that we obtain the random graph,
\index{graph ! random}%
 whether or not the other axioms (such as Choice or Infinity) are satisfied. Indeed, Rado's construction can be regarded as obtained by this method using a hereditarily finite model of set theory
(in which all sets are finite).

\head{Remark} Since $\mathfrak{R}$ is homogeneous, its automorphism
group is very large: it acts transitively on vertices, edges, non-edges, and
indeed on any isomorphism type of finite subgraphs. But each of our three
constructions of $\mathfrak{R}$ has the property that no non-trivial
automorphisms are ``obvious'' from the construction.

\section{$\aleph_0$-categoricity}
\label{Aleph0categoricity}

Let $L$ be a countable first-order language.
\index{first-order language}%
 A theory in $L$ is said to be
\emph{$\aleph_0$-categorical} if it has a unique countable model; and a
countable $L$-structure is \emph{$\aleph_0$-categorical} if its theory is,
that is, if it is the unique countable model of its theory (up to isomorphism).

Two examples of $\aleph_0$-categorical theories are:
\begin{itemize} 
\item the rational numbers $\mathbb{Q}$, as ordered set (this follows
from Cantor's theorem
\index{Cantor's theorem}%
 characterising $\mathbb{Q}$ as the unique countable
dense linear order
\index{linear order}%
 without endpoints);
\item the random graph $\mathfrak{R}$
\index{graph ! random}%
 (this follows because condition $(*)$
in the preceding section is first-order, being the conjunction of infinitely
many $\forall\exists$ sentences, one for each pair of values of $|U|$ and
$|V|$).
\end{itemize}

On the other hand, Peano's axioms
\index{Peano's axioms}%
 for the natural numbers, or the
Zermelo--Fraenkel
\index{Zermelo-Fraenkel set theory}%
 axioms for set theory, are not $\aleph_0$-categorical
(assuming that they are consistent), since
in each case there are many countable models. (For Peano's axioms, these
are the ``non-standard models of arithmetic''; for ZF there are, for example,
models with and without the axiom of choice, and by the downward 
L\"owenheim--Skolem theorem there are countable models in each case.)

The remarkable theorem of Engeler, Ryll--Nardzewski and Svenonius
\index{Engeler--Ryll-Nardzewski--Svenonius ! Theorem}%
 characterises the $\aleph_0$-categorical countable structures in terms of their symmetry.

\begin{theorem}
A structure $M$ over the first-order language
\index{first-order language}%
 $M$ is $\aleph_0$-categorical
if and only if $\Aut(M)$ has only finitely many orbits on $M^n$ for each
natural number $n$.
\end{theorem}

A permutation group is said to be \emph{oligomorphic} if it has only finitely
any orbits on the $n$th Cartesian power of the permutation domain for each
natural number $n$. (So a structure $M$ is $\aleph_0$-categorical if and only
if its automorphism group
\index{group ! automorphism}%
 is oligomorphic.) It is easy to see that we may
equivalently ask for finitely many orbits on the set of $n$-tuples of distinct
elements, or on the set of $n$-element subsets for each $n$.

To illustrate, we look again at our two examples of $\aleph_0$-categorical
structures.

\begin{itemize}
\item $\Aut(\mathbb{Q},<)$ has just one orbit on $n$-element subsets of
$\mathbb{Q}$, for each $n$: for given two such subsets, it is easy to find
a piecewise-linear order-preserving
\index{linear order}%
 permutation carrying one to the other.
\item Since $\mathfrak{R}$ is homogeneous, the number of orbits on $n$-element
subsets is equal to the number of graphs with $n$ vertices, which is clearly
finite. So $\Aut(\mathfrak{R})$ is oligomorphic.
\end{itemize}

Extending the argument in the second example, we see that if $M$ is any
countable homogeneous structure, then the number of $\Aut(M)$-orbits on
$n$-tuples is equal to the number of $n$-element substructures of $M$ up to
isomorphism. If the language $L$ is a finite relational language, then this
number is finite. We conclude that
\begin{quote}
A countable homogeneous structure over a finite relational language is
$\aleph_0$-categorical.
\end{quote}

We see, for example, that a countable model of ZF is not $\aleph_0$-categorical,
but that if we symmetrise the membership relation, the resulting graph is
$\aleph_0$-categorical (and indeed we get the same graph, no  matter what
model of ZF we chose).

As for the first example,  the $\aleph_0$-categoricity of the rationals with the usual ordering $(\mathbb{Q},<)$ which was proved by Cantor,
\index{Cantor, G.}%
was generalized by by Skolem~\cite{skolem}
\index{Skolem, T.}%
who showed that it is possible to partition $\mathbb{Q}$ by adding finitely many dense and cofinal subsets.  The Covington graph $\mathfrak{N}$
\index{graph ! Covington}%
 can be described in terms of a dense codense colouring of the rationals~\cite{cam11}.  We conjecture that it is possible to realize it as a Cayley graph that can be constructed from a conjugacy-class generating set.  (By contrast a countable $N$-free graph admits at most one conjugacy class of cyclic automorphisms)~\cite[p.96]{cam6}.  It may also be possible to build other structure-free graphs, such as square-free $W$-free graphs, using conditions similar to those parameterizing $N$-free graphs~\cite[p.97]{cam6}.

We end with some pedagogical remarks regarding the categoricity
\index{categoricity}%
 of random graphs.
\index{graph ! random ! categoricity}%
Whilst the two-coloured random graph is known to be $\aleph_{0}$-categorical
\index{aleph@$\aleph_0$-categorical}%
it is not $\aleph_{1}$-categorical.  The higher-adjacency random graph $\mathfrak{R}_{\omega,\omega}$ is not $\aleph_0$-categorical if we regard colours as ``named'' relations, but it is if we simply take the structure as a partition of edges into colour classes.  It is unlikely to be $\aleph_1$-categorical.  Because $\mathfrak{R}_{m,\omega}$ is not $\aleph_{1}$-categorical there is no unique uncountable random graph, so it is not clear what an object like $\mathfrak{R}_{\kappa_1, \kappa_2}$ where the $\kappa_i$~\label{kappa} are chosen uncountable cardinals, actually is but this may be an interesting question for set-theorists.  

Another discussion of $\aleph_0$-categoricity can be found in~\cite{kayem}.  We refer also to~\cite{bodirskyj} for extensions of interpretability in $\aleph_{0}$-categorical structures.
\index{aleph@$\aleph_0$-categorical}%

\section{Mathematics and Logic}
\label{MathematicalLogic}%
Logic plays a dual role in mathematics; it provides the foundations on which the subject is built, and it is a branch of mathematics in its own right developing by using the common culture of mathematics, and makes its own very important contributions to this culture.  Let's expand on this observation.

The \emph{Deduction Theorem}
\index{Deduction Theorem}%
 asserts that, if a proposition $q$ can be deduced from a set $S$ of propositions together with an extra proposition $p$, then the proposition $p \models q$ can be deduced from $S$ alone. The proof, of course, is by induction on the number of steps in the assumed proof of $q$.

Now a formula $p$ is logically valid if it is true in every interpretation, and is a logical consequence of a set $S$ if it is true in all interpretations which make every formula in $S$ true.

Finally we can say what soundness and completeness are. A system is sound if all its theorems are logically valid, and is complete if all its logically valid formulas are theorems. A sound system is free of contradiction, since a contradiction cannot be logically valid. Clearly soundness is an essential requirement, and completeness a desirable requirement, for any logical system which is to serve as a basis of mathematics.  

Two very important and familiar logical systems are sound and complete. The first of these is propositional (or Boolean) logic, where atomic propositions have no internal structure, but are simply statements which are either true or false and are combined with connectives such as ``and'', ``or'', ``not'', ``implies''.  We can decide whether a compound proposition is logically valid by constructing a truth table for it; there are standard deduction systems for propositional logic which prove precisely the logically valid formulas.

Secondly there is first-order logic,
\index{first-order logic}%
 which is the formalism in which most of mathematics is naturally cast. Its structure has constants, functions, and variables which are combined into `atomic' propositions by means of relations (including the relation of equality); these are then joined by connectives as in the Boolean case, or bound by quantifiers over variables such as ``for all x, $\ldots$'' or ``there exists $x$ such that $\ldots$'' G\"odel proved that first-order logic is complete.  \emph{G\"odel's Completeness Theorem}
\index{godel@G\"odel's Completeness Theorem}%
says that in first-order logic, truth and provability are equivalent.  Higher-order logics have sound deductive systems, but no such system can be complete; (though second-order logic is complete for Henkin semantics).

\emph{G\"odel's Incompleteness Theorem}
\index{godel@G\"odel's Completeness Theorem}%
 asserts something more specific. Any logical system provided with axioms which attempts to describe the natural numbers with the operations of addition and multiplication and the relation of order, if it is consistent, must be ``incomplete'', in the sense that there are true statements about the natural numbers which are not provable. Of course, by the completeness theorem, any statement which is true in all models of the axioms is provable. So this is actually a constructive theorem: it says that any axiom system which attempts to describe the natural numbers, if consistent, will always have `non-standard models', structures which satisfy the axioms but fail to have some property which is true of the natural numbers.

This is quite subtle. It is known that the theory of the natural numbers, with the order relation and the operation of addition but without multiplication is complete: this is \emph{Pressberger arithmetic}~\cite{presburger}.
\index{Presburger arithmetic}%
\index{Presburger, Moj\'zesz}%
  Somehow throwing in multiplication as well makes the difference, despite the truism that ``multiplication can be defined in terms of addition''.  Not to mention \emph{Skolem Arithmetic},
\index{Skolem arithmetic}%
\index{Skolem, T.}%
 where the operation of multiplication is kept but addition is thrown away, which is also complete.  Moj\'zesz Presburger proved Presburger arithmetic to be:
consistent, complete and decidable, that is, there exists an algorithm which decides whether any given statement in Presburger arithmetic is true or false. Robinson arithmetic (Q),
\index{Robinson arithmetic}%
 is a finitely axiomatised fragment of Peano arithmetic (PA),
\index{Peano arithmetic}%
 first set out in~\cite{robinson} by Raphael Robinson.
\index{Robinson, R. M.}%
 This is essentially PA without the axiom schema of induction. Since it is weaker than PA, it is incomplete in the sense of G\"{o}del,
\index{godel@G\"{o}del, K.}%
 but crucially this weak finitely axiomatised theory is already incompletable and essentially undecidable.  More on these matters can be found in~\cite{boolos} and~\cite{marker}.

The desirable properties of a logical system, soundness and completeness, mean that a formula $f$ is provable if and only if it is true in all interpretations.  More generally, a formula $f$ is provable from a set $S$ of formulas if it is true in all interpretations satisfying all the formulas in $S$.

This last property has an important consequence: a set $S$ of formulas is consistent (no contradiction can be proved from it) if and only if it is satisfiable (true in some interpretation).  For a contradiction is true in no interpretation, so is deducible from $S$ if and only if $S$ is true in no interpretation.

Consistency was a matter of great signficance to Hilbert and his followers.  Part of G\"odel's Incompleteness Theorem states that a logical system which is at least strong enough to include the theory of the natural numbers (with order, addition, and multiplication) cannot prove its own consistency.  For example, Peano's axioms,
\index{Peano's axioms}%
 the commonest such system for the natural numbers, cannot prove their own consistency.

The vast majority of mathematicians accept the existence of the natural numbers.  It is easy to see that they do indeed satisfy Peano's axioms.
\index{Peano's axioms}%
 So, as stated above, this shows that the axioms are consistent.

Another take on this is to consider the Zermelo-Fraenkel (ZF) axioms for set theory.
\index{Zermelo-Fraenkel set theory}%
 Set theory is generally accepted as the standard foundation, in which everything that mathematicians study can be built.  The Zermelo-Fraenkel axioms are strong enough to prove the consistency of the Peano axioms because the natural numbers can be constructed in set theory.

Of course this only pushes the difficulty elsewhere, since the ZF axioms cannot prove their own consistency. But ZF with a ``large cardinal axiom'' proves the consistency of ZF. The higher we go, the smaller proportion of mathematicians accept the axioms, of course.  While few doubt the existence of the natural numbers, ``inaccessible cardinals'' are more problematic.

We mention in passing a popular fad, which is to regard category theory
\index{category theory}%
 rather than set theory as the natural foundation for mathematics. But the foundational difficulties with category theory (the large cardinal axioms that have to be assumed) are much stronger than for set theory!

Briefly, the construction of numbers in set theory. Numbers are for counting; like the standard metre, the standard number 2 should be a 2-element set against which other sets can be compared to see if they have 2 elements. The standard set 0 must be the empty set (it follows from the Axiom of Extensionality that there is only one empty set, so we have no choice). Then the standard 1-element set should be the only one we have at this point, namely the set {0}; the standard 2-element set {0,1} (in other words, {0,{0}}); and so on. Each number is the set of its predecessors.

To return to the incompleteness theorem, the existence of true but unprovable statements in a theory at least as strong as Peano's.
\index{Peano's axioms}%
 We have to look at what appear to be ways to circumvent this. If $f$ is true but unprovable, no contradiction would be introduced by simply adding $f$ as a new axiom.  But G\"odel's argument still applies, and gives a true but unprovable statement in the new theory; and so on ad infinitum.

Adding all true statements as axioms also fails, because one of our requirements for a logical system is that there should be a mechanical procedure for recognising the axioms; and there is no mechanical procedure for recognising true statements about the natural numbers. (This was Alan Turing's
\index{Turing, A.}%
 great contribution; but in a sense, it must be so, else G\"odel's theorem would be false. Turing gave a precise definition of a ``mechanical procedure'' by inventing the Turing machine, the theoretical computer which is the prototype for our real computers.)

What do these true but unprovable statements look like? They have the property that G\"odel required of them; but, if written out as logical statements in terms of order, addition and multiplication on the natural numbers, they are entirely unnatural (and dependent on the particular axiom system chosen).

There was considerable excitement in the 1970s when Jeff Paris
\index{Paris, J.}%
 and Leo Harrington
\index{Harrington, L.}%
 found the first example of a ``natural'' true but unprovable statement. (There was also Gerhard Gentzen's
\index{Gentzen, G.}%
 1943 direct proof of the unprovability of $\epsilon_0$-induction in Peano arithmetic,
\index{Peano arithmetic}%
 and Goodstein's Theorem).
\index{Goodstein's Theorem}%
  Their statement was a variant of Ramsey's theorem,
\index{Ramsey theory}%
 to which we shall return shortly.  (We mention Harvey Friedman's
\index{Friedman, H.}%
  work on \emph{reverse mathematics},
\index{reverse mathematics}%
 whose aim is to derive the axioms of mathematics from the theorems considered to be necessary. One aspect of this is his \emph{Boolean relation theory}, which attempts to justify large cardinal axioms by demonstrating their necessity for deriving certain important propositions.)
 
It is very nearly true that no statement which mathematicians had previously invested time and effort in trying to prove has ever turned out to be undecidable. (The exception to this is the Continuum Hypothesis,
\index{Continuum Hypothesis}%
 the statement that there is no set intermediate in size between the natural numbers and the real numbers.  In 1960, Paul Cohen
\index{Cohen, P.}%
showed that it is independent of the ZF axioms. This means that, rather like geometry where we can study either Euclidean
\index{Euclidean geometry}%
 or non-Euclidean geometry (or both), we can do mathematics in which the Continuum Hypothesis is true, or in which it is false.)  So the syntactic aspect of logic is not much of a concern. 

However results such as G\"odel's Completeness Theorem for first-order logic
\index{first-order logic}%
 have important consequences for ordinary mathematics. Here are two of them, which can be regarded as the portals of model theory.
\index{model theory}%
 (A model for a set $S$ of formulas is a structure in which the formulas of $S$ are true. This is entirely a semantic notion.)  Let $S$ be a set of sentences (formulas with no free variables) in a first-order language
\index{first-order language}%
  which has at most countably many symbols in its alphabet.

\emph{The Compactness Theorem}: If every finite subset of $S$ has a model, then $S$ has a model.
\index{Compactness Theorem}%
\emph{The Downward L\"owenheim--Skolem Theorem}: If $S$ has a model, then it has a model which is at most countably infinite.
\index{Downward L\"{o}wenheim--Skolem ! Theorem}%

We comment briefly on the proofs. According to soundness and completeness, the statement ``S has a model'' is equivalent to ``S is consistent''.  Now if $S$ is inconsistent, then a contradiction can be proved from it; since proofs are finite, the contradiction must follow from a finite subset of $S$. This proves compactness. The Downward L\"owenheim--Skolem Theorem relies on the fact that, in the proof of completeness, it is necessary to construct a model for a consistent set of sentences; observing the proof we see that the model constructed is at most countably infinite.

Ramsey's Theorem
\index{Ramsey theory}%
 has a rather complicated statement. The easiest finite form of it is the party theorem: given any number $k$, there is a number $n=R(k)$ such that, among any given $n$ people, there are either $k$ mutual friends or $k$ mutual strangers.  (We assume that any two people are either friends or strangers. The proof that $R(3)=6$ is well-known.) This can be formulated (and proved) in Peano arithmetic.
\index{Peano arithmetic}%
 But there is an infinite form: given an infinite number of people, there are either infinitely many mutual friends or infinitely many mutual strangers. This is a theorem of set theory which can be proved from the ZF axioms. Now the Compactness Theorem
\index{Compactness Theorem}%
  shows that the finite party theorem is a consequence of the infinite one.

The \emph{Paris-Harrington Theorem}~\cite{parish},
\index{The Paris-Harrington Theorem}%
 alluded to before, is a generalised version of Ramsey's Theorem.
\index{Ramsey theory}%
   The `party' version says that, given $k$, there is a number $n=PH(k)$ such that, given $n$ people $p_0, \ldots, p_{n-1}$, there is a set of at least $k$ people, mutual friends or mutual strangers, so that if the least-numbered person in the set is pm, then there are more than $m$ people in the set.

The theorem can be deduced from the infinite Ramsey theorem by exactly the same compactness argument that gives the finite Ramsey theorem. But, although it can be formulated in the language of Peano arithmetic, it cannot be proved from Peano's axioms for the interesting reason that the appropriate ``Paris-Harrington function'' grows so rapidly that it cannot be expressed in terms of the functions of arithmetic!

We complete this section with another connection, one that makes transparent the relevance of the above to our study in this manuscript.

By the Upward L\"{o}wenheim--Skolem Theorem,
\index{Upward L\"{o}wenheim--Skolem ! Theorem}%
 a theory which has infinite models will have arbitrarily large infinite models; so no set of first-order axioms can uniquely specify an infinite structure. The next best thing is to specify a structure once the size is given. The case of most relevance to us is the following.

A set $S$ of first-order sentences
\index{first-order sentence}%
 is \emph{countably categorical}
\index{countably categorical}%
 if it has only one countable model up to isomorphism.  A countable structure with countably categorical theory is thus one which can be completely specified by a set of first-order axioms together with the extra assertion that it is countable.  The prototype Cantor's theorem
\index{Cantor's theorem}%
 characterising $\mathbb{Q}$ as the unique countable dense total order without endpoints. (The conditions of being totally ordered and without endpoints are first-order statements.)

Now we repeat the theorem proved independently in 1959 by three people, Engeler, Ryll-Nardzewski, and Svenonius:
\index{Engeler, E.}%
\index{Ryll-Nardzewski, C.}%
\index{Svenonius, L.}%
\index{Engeler--Ryll-Nardzewski--Svenonius ! Theorem}%

A countable structure $M$ has countably categorical theory if and only if its automorphism group
\index{group ! automorphism}%
 has only finitely many orbits on the set of $n$-tuples of elements of $M$ for all natural numbers $n$.

\emph{In other words, axiomatisability (by countability and first-order axioms) is equivalent to symmetry (for the condition of the theorem asserts that such a structure has a huge and rich group of automorphisms).  This equivalence of axiomatisability and symmetry is surely one of the most remarkable facts in mathematics.}

Oligomorphic permutation groups which have only finitely many orbits on $n$-tuples for all $n$, different orbits representing different ``shapes'' of $n$-tuples, and only a few (i.e. finitely many) such shapes are required to satisfy the condition.  These groups are essentially just the automorphism groups
\index{group ! automorphism}%
 of structures with countably categorical theories.  Given an oligomorphic permutation group, and consider the number of its orbits on $n$-tuples. Call a sequence of natural numbers realisable if it arises in this way.

Which sequences are realisable? We are very far from a complete answer to this question; but here are two general remarks.

By the Downward L\"owenheim--Skolem Theorem, if a sequence is realisable, then it is realisable by a countable group acting on a countable set.  So this theory is not plagued with difficulties about ``large cardinals''.  By the Compactness Theorem,
\index{Compactness Theorem}%
 a sequence is realisable if and only if every (finite) initial subsequence of it is an initial subsequence of some realisable sequence. So realisability is a ``finitely determined'' property.

\section{Previous Results in Multicoloured Random Graph Theory}
\label{PreviousResults}
\index{graph ! random ! $m$-coloured}%

Ours is not the first study of multi-adjacency random graphs.  They
have been analysed before by J. K. Truss
\index{Truss, J. K.}%
as countable universal homogeneous
$C$-coloured graphs $\Aut(\Gamma_{C})$~\label{Gamma_C} for $|C| \ge 2$~\cite{truss1}~\cite{truss2},
using detailed permutation group theory.  Whilst our
results are complementary to these works, it is appropriate to list
their main findings. 

Among the results of~\cite{truss1} are:-

(i)  $\Aut(\Gamma_{C})$ is simple for each $C$.
\index{group ! simple}%
(Clearly for $|C| = 1$,
$\Aut(\Gamma_{C}) \cong \Sym(\aleph_0)$ which has $\FSym(\aleph_0)$ as
a normal subgroup).

(ii) Non-identity members of $\Aut(\Gamma_{C})$ have infinite support
$\forall\ C$.  (Clearly not true for $|C| = 1$).

(iii) Different cycle types arise, depending on $C$.

(iv) For each $C$ with $2 \le |C| \le \aleph_{0}$, $\exists g \in
\Aut(\Gamma_{C})$ of cycle type $(k_{\infty}, k_{1}, k_{2}, k_{3},
\ldots) = (n, 1, 0, 0, \ldots)\ \forall n \ge |C|$, but for $m < |C|$
none of type $(m, 1, 0, 0, \ldots)$.  Hence $\Aut(\Gamma_{C}) \cong
\Aut(\Gamma_{C'})$ as \emph{permutation} \emph{groups}
\emph{iff} $|C|=|C'|$.
\index{group ! permutation}%
(M. Rubin~\cite{rubin}
\index{Rubin, M.}%
has shown that the $\Aut(\Gamma_{C})$ are non-isomorphic for different values of $|C|$).

(v)  For each $C$ with $2 \le |C| < \aleph_{0}$, $\Gamma_{C}$
is $\aleph_0$-categorical.
\index{aleph@$\aleph_0$-categorical}%

(vi) The group $\Aut(\Gamma_{\aleph_0})$ has a subgroup
isomorphic to $\Sym(\aleph_0)$.

(vii)  For $|C| \ge 1$, $\Aut(\Gamma_{C})$ contains two subgroups
isomorphic to $\Aut(\Gamma_{C+1})$ whose intersection is isomorphic to
$\Aut(\Gamma_{C+2})$.

(viii)  There is a subgroup chain $\Aut(\Gamma_{1}) = H_{1}> H_{2}>
\ldots >  H_{\infty}$ such that for all $n \le \infty$, $H_{n} \cong
\Aut(\Gamma_{n})$ and $H_{\infty} = \cap_{n<\infty} H_{n}$.  (Clearly
any extra structure imposed on $\Gamma_{n}$ restricts the class of
automorphisms, so defines a subgroup).

(ix)  For any $C, C'$ with $1 \le |C|, |C'| \le \aleph_{0}$,
$\Aut(\Gamma_{C})$ is isomorphic to a subgroup of $\Aut(\Gamma_{C'})$.  So for
example, each of $\Aut(\Gamma_{1})$ and $\Aut(\Gamma_{2})$ can be embedded as a subgroup of the other; note however that $\Aut(\Gamma_{2})$ is simple whilst
$\Aut(\Gamma_{1})$ is not.

\smallskip

In~\cite{truss1} Truss showed that five conjugates suffice to prove that 
$\Aut(\Gamma_C)$ is simple, that is any non-identity element can be 
expressed as the product of five conjugates of any other non-identity 
element.  In~\cite{truss5} this was improved to three conjugates.  We give 
an outline of the form of the original proof and refer the reader to the articles for a detailed discussion.

Let $C$ be the set of edge-colours and $F : E(\Gamma_{C}) \to C$ be the colouring function on the edges.  Let $\Sigma$ be the set of all $g \in \Aut(\Gamma_C)$ which are an infinite product of disjoint infinite cycles (with no finite cycles) for which: if $\alpha$ maps a finite subset of $\Gamma_C$ into $C$ then $\exists x \in \Gamma_C \backslash \dom(\alpha)$ such that $x g^n \notin \dom(\alpha)$ $(\forall n \in \mathbb{Z})$ and $F\{x, y\} = \alpha(y)$ $(\forall y \in \dom(\alpha))$.  Call any $x$ that satisfies the formula $(\exists x)(\forall y \in \dom(\alpha)) F\{x, y\} = \alpha(y)$ a \emph{witness} for the formula.  Witnesses for the formula lie outside any cycles of elements of $\Sigma$ containing members of $\dom(\alpha)$.  Firstly two lemmas are proved:
\begin{lemma} If $1 \neq g \in \Aut(\Gamma_C)$ then $\{x \in \Gamma_C : g x \neq x\}$ is infinite.
~\label{infsuplem}
\end{lemma}
So non-identity members of $\Aut(\Gamma_C)$ have infinite support.  Truss comments that this appears to be the only point in the proof of the simplicity of $\Aut(\Gamma_C)$ where $|C| \ge 2$ is used.

\begin{lemma} If $1 \neq g \in \Aut(\Gamma_C)$, $n \in \mathbb{N}$, $\alpha$ maps a finite subset of $\Gamma_C$ into $C$, and $A$ is a finite subset of $\Gamma_C$, then $\exists x \in \Gamma_C$ such that
\[ (\forall i) ( -n \le i \le n) \Rightarrow x g^i \notin A \cup \dom(\alpha)),\ x g \neq x \]
and
\[ (\forall y \in \dom(\alpha))\ F\{x, y\} = \alpha(y). \] 
\end{lemma}
So infinitely many of the witnesses are moved by $g$.

Then the following two results are proved:
\begin{theorem} If $1 \neq g_1, g_2 \in \Aut(\Gamma_C)$ then there is a conjugate $h$ of $g_1$ such that $g_2 h \in \Sigma$.
\end{theorem}

So there are conjugates $g_1, g_2$ of $g^{-1}$ such that $g^{-1} g_1, h g_2 \in \Sigma$.

\begin{theorem} If $g_1, g_2, g_3 \in \Sigma$ then there are conjugates $h_1, h_2, h_3$ of $g_1, g_2, g_3$ such that $h_3 h_2 h_1 = 1$. 
\end{theorem}

So there are conjugates $h_1, h_2$ of $g^{-1} g_1$ and $h_3$ of $h g_2$ such that $h_3 h_2 h_1 = 1$.

Putting these together, $h_3 = h_1^{-1} h_2^{-1} = 1$, so $h_3$ is the product of four conjugates of $g$.  Hence so is $h g_2$, and so if $g$ and $h$ are non-identity elements of $ \Aut(\Gamma_C)$, then $h$ is the product of five conjugates of $g$.  So $\Aut(\Gamma_C)$ is simple.

\smallskip

In~\cite{truss2}, $\Aut(\Gamma_{C})$ is extended to a highly
transitive permutation group $\AAut(\Gamma_{C})$ on the same set by
considering ``almost automorphisms''
\index{group ! almost automorphism}%
 of $\Gamma$, these being permutations $g$
of $\Gamma_{C}$ for which the set of two-element subsets $\{x, y\}$ of
$\Gamma_{C}$ such that $\{x, y\}$ is a different colour to $\{gx, gy\}$ is
finite; that is permutations of $\Gamma_{C}$ which preserve the colour
of all but a finite set of edges.

Among the results of~\cite{truss2} are:-
\begin{itemize}
\item[(i)]  The lattice of non-trivial normal subgroups $N$ of
$\AAut(\Gamma_{C})$ are isomorphic to the lattice of subgroups of the
free abelian group of rank $n$, where $n = |C|-1$.
\index{group ! free abelian}%
In particular $\AAut(\Gamma_{C})$ has a unique minimal non-trivial normal subgroup.
\item[(ii)]  Any non-trivial normal subgroup $N$ of $\AAut(\Gamma_{C})$
contains $\Aut(\Gamma_{C})$.
\item[(iii)]  Cycle types occurring in $\Aut(\Gamma_{C})$ and
$\AAut(\Gamma_{C})$ are the same apart from those which are the product
of finitely many cycles.
\end{itemize}

A classification of all the cycle types occurring in the automorphism groups
\index{group ! automorphism}%
 of countable homogeneous graphs can be found in~\cite{lovelltruss}.

\smallskip

A word which is an element of a free group
\index{group ! free}%
 on a finite set for which every equation of the form $w = g$ for an element $g$ of a given group $G$ is soluble, is said to be \emph{universal} for $G$.  In~\cite{droste1} Droste
\index{Droste, M.}%
and Truss prove special cases of the conjecture that a word is universal for $G = \Aut(\Gamma_{C})$ if and only if it cannot be written as a proper power, a result already established for infinite symmetric groups.
\index{group ! symmetric ! infinite}%

\medskip

Finally, M. Rubin
\index{Rubin, M.}%
has proved the following results:-
\begin{itemize}
\item[(i)]  If the automorphism group
\index{group ! automorphism}%
 of a countable homogeneous graph
(the varieties are classified by the Lachlan-Woodrow Theorem)
\index{Lachlan-Woodrow Theorem}%
\index{Lachlan, A. H.}%
\index{Woodrow, R. E.}%
is primitive then it is simple.
\item[(ii)]  The stabilizer $\Aut(\Gamma_{C})_{(X)}$ is a simple group
\index{group ! simple}%
for finite $X \subseteq \Gamma_{C}$.
\end{itemize}

\bigskip

Combinatorial aspects of random graph colourings are explored in~\cite{nelsw}, but these refer to finite graphs in contrast to our emphasis on the infinite. 

\bigskip

\section{Further Details on Random Graphs}
\label{FurtherDetails}

Let $V(\Gamma)$ denote the vertex set of graph $\Gamma$ and $E(\Gamma)$ its edge set.  The \emph{trivial graph}
\index{graph ! trivial}%
consists of just one vertex.
 A bijection $\alpha: V(\Gamma) \to V(\Gamma')$ between two graphs
$\Gamma$ and $\Gamma'$ is an \emph{isomorphism}
\index{graph ! isomorphism}%
if $\{u, v\} \in E(\Gamma) \Leftrightarrow \{\alpha(u), \alpha(v)\} \in E(\Gamma')$.
An \emph{automorphism}
\index{graph ! automorphism}%
is a self-isomorphism.  The set of automorphisms of $\Gamma$ under the operation of composition form a group, denoted $\Aut(\Gamma)$, and its elements permute the vertices
of $\Gamma$ so as to leave invariant its edges.  For more on graph
automorphisms see~\cite{cam12a}~\cite{laurisc}.  A \emph{graph invariant}
\index{graph ! invariant}%
is a map which assigns equal values to isomorphic graphs.  One of the simplest graph invariant is the \emph{diameter}~\label{diam}
\index{graph ! diameter}%
of a connected graph, which is the maximal distance between pairs of its vertices.

Presently we shall explain why the \emph{random graph} $\mathfrak{R}$
\index{graph ! random}%
is the unique countably infinite graph~\cite[p.~241]{bollobas} whose defining relation is either given by the ($*$)-condition:

($*$)  If $U$ and $V$ are finite
disjoint sets of vertices of $\mathfrak{R}$, then there exists in
$\mathfrak{R}$ a vertex $z$ joined to every vertex in $U$ and to no vertex in $V$;

or equivalently, by the one-point extension property:
\index{one-point extension property}%

$(\dagger)$  If $A \subset B$ are finite graphs with $A$ the induced subgraph of $B$ obtained by deleting one vertex, then any embedding of $A$ into graph $\mathfrak{R}$ can be extended to an embedding of $B$ into $\mathfrak{R}$.

Condition $(\dagger)$ is called the \emph{I-property}
\index{I-property}%
because it is a form of injectivity.  Clearly,  $(\dagger)$ is equivalent to the property ($*$), which is just its interpretation for graphs.  The back-and-forth method
\index{back-and-forth method}%
shows the uniqueness of a countable graph satisfying them.  An example of the back-and-forth argument proving uniqueness, and its extension proving homogeneity for a slightly different random graph can be found in Chapter~\ref{fstchap}.

Property ($*$) is a conjunction of one first-order sentence
\index{first-order sentence}%
 per application (or per value of $|U \cup V|$), so $\mathfrak{R}$ is $\aleph_{0}$-categorical and $\Aut(\mathfrak{R})$ is oligomorphic.  Property $(\dagger)$ implies that $\mathfrak{R}$ is \emph{universal}, that is that any finite or countable subgraph of $\mathfrak{R}$ is embeddable in $\mathfrak{R}$ as an induced subgraph.  One further property of $\mathfrak{R}$ is its \emph{homogeneity}, meaning that any isomorphism between finite subgraphs of $\mathfrak{R}$ can be extended to an automorphism of $\mathfrak{R}$.  Homogeneity is a measure of the symmetry of a mathematical object.  It is quite a general symmetry property that for graphs includes as special cases the properties of vertex-, edge- and path-transitivity.  The graph $\mathfrak{R}$ is the unique countable homogeneous universal graph. 

An $\aleph_0$-categorical structure
\index{aleph@$\aleph_0$-categorical}%
is universal, as can be proved~\cite{cam3} using K\"onig's Infinity Lemma~\cite{konig}.  

A homogeneous structure is universal, as proved in~\cite{cam3}.

P. Komj\'ath
\index{Komj\'ath, P.}%
 and J. Pach
\index{Pach, J.}%
survey~\cite{komjathpach} results in the theory of universal graphs with different properties, including examples of classes of graphs for which there is no universal element, for example, the class of $K_{\omega}$-free countably infinite graphs.

A graph is \emph{$n$-existentially closed} or \emph{$n$-e.c.}
\index{graph ! n@$n$-existentially closed}%
if for every pair of subsets $U$, $W$ of the vertex set $V(\Gamma)$ of the graph such that $U\cap W=\emptyset$ and $|U|+|W|=n$, there is a vertex $v\in V(\Gamma) \backslash (U\cup W)$ such that all edges between $v$ and $U$ are present and no edges between $v$ and $W$ are present.  So another name for property ($*$) is the \emph{existentially closed}
\index{existentially closed (e.c.)}%
or \emph{e.c.} adjacency property.  This leads to stability properties such as for each $v \in \mathfrak{R}$, $\mathfrak{R} \cong \mathfrak{R} \backslash \{v\}$.  The $n$-e.c. property first arose in the discussion of random graphs
\index{graph ! random}%
 and the zero-one law
\index{zero-one law}%
for first-order sentences~\cite{glebskii} (see below),
\index{first-order sentence}%
 but it is more general.  Let $M$ be a homogeneous relational structure
\index{relational structure}%
  with age Age$(M)$.  Then $(\dagger)$ is the e.c. property that holds in $M$ for any finite substructures $A, B \subseteq$ Age$(M)$ of $M$ with $|B| = |A| + 1$.  We say that a finite $S$ is \emph{$n$-e.c.} with respect to Age$(M)$ if \begin{itemize}
\item[(i)]  $S \in$ Age$(M)$;
\item[(ii)]  $(\dagger)$ holds for all $A, B \in$ Age$(M)$ with $|A| \le n$;
\end{itemize}
The following question arises: when can we find finite $n$-e.c. structures for a Fra\"{\i}ss\'e class
\index{Fra\"{\i}ss\'e class}%
 Age$(M)$, and how large do they have to be in terms of $n$?
For graphs, a lower bound $2^{n} + n$ is known, and an upper bound has been derived from Paley graphs.
\index{graph ! Paley}%
Almost all random graphs are $n$-e.c.~\cite{camerstark}.

The book by Higman and Scott~\cite{higmansc}
\index{Higman, D. G.}%
\index{Scott, E.}%
is an account of the theory of existentially closed groups.

The structure of the graph $\mathfrak{R}$ can be generalized by increasing the number of different adjacency types from two to any number.  In this manuscript we shall mostly be interested in the three-adjacency case.  These three can more conveniently be thought of as colours, red, blue and green ($\mathfrak{r}, \mathfrak{b}, \mathfrak{g}$) thereby giving the \emph{triality
graph} $\mathfrak{R^{t}}$,
\index{graph ! triality}%
which is the 3-colour generalization of $\mathfrak{R}$, with a modified injectivity property which is an interpretation of $(\dagger)$ for 3-coloured complete graphs:

($*_{t}$)  If $U$, $V$ and $W$ are finite
disjoint sets of vertices of $\mathfrak{R^{t}}$, then there exist in
$\mathfrak{R^{t}}$, a vertex $z$,
joined to every vertex in $U$ with a red edge, to every vertex in
$V$ with a blue edge, and to every vertex in $W$ with a green edge.  

The $m$-coloured random graph for $m \ge 3$
\index{graph ! random ! $m$-coloured}%
edge colours is similarly defined and is denoted $\mathfrak{R}_{m,\omega}$, where $\omega$ denotes a countable infinity of vertices; so $\mathfrak{R^{t}} = \mathfrak{R}_{3,\omega}$.

The graphs $\mathfrak{R}_{m,\omega}$ are also universal and homogeneous, and the proofs of this are similar to those given in the main text for $\mathfrak{R^{t}}$.  Another measure of symmetry is the size of an object's automorphism group.
\index{group ! automorphism}%
  Most graphs, certainly finite ones, are \emph{asymmetric},
\index{graph ! asymmetric}%
that is they possess a trivial automorphism group
\index{group ! automorphism}%
 consisting solely of the identity element.  By contrast the automorphism group of $\mathfrak{R}$ and its multicoloured
\index{graph ! random ! $m$-coloured}%
 versions are gigantic, being uncountably infinite in size, thereby indicating a colossal degree of symmetry possessed by these graphs.

In short then, the multicoloured graph $\mathfrak{R}_{m,\omega}$ has
similar structural properties to the two-coloured graph $\mathfrak{R}$.  However there are differences in the structure of the groups supported by the graphs, as will unfold in the main text.  Whenever we say `the random graph',
\index{graph ! random}%
 we mean the two-coloured random graph.  A random countable graph on two edge colours, that is with two adjacency relations, is isomorphic to $\mathfrak{R}$ with probability $1$.

It is useful to imagine the `evolution' of a random structure as a set of $n$ points develops a structure, which in the case of random graphs is the addition of edges at random; this is how P. Erd\H{o}s and A. R\'enyi 
\index{Erd\H{o}s, P.}%
\index{R\'enyi, A.}%
 originally viewed random processes on large discrete structures.  The main aim of random graph theory is to determine at what stage of this evolutionary process is the graph likely to develop a certain property.  A way of establishing a formalism to study this precisely is to consider the property of a `typical' graph in a probability space
\index{probability space}%
of graphs of a particular type.   The simplest probability space consists of all graphs on $n$ vertices, each of which is assigned the same probability, though there are many such \emph{models}
\index{graph ! random ! models}%
of random graphs~\cite[Chapter 2]{bollobas1}.  

To highlight the evolutionary process, consider a \emph{random graph process}~\cite[p.~42]{bollobas1}
\index{graph ! random ! process}%
on $V(\Gamma) = \{0, \ldots, n\}$, or simply a \emph{graph process}, which is a Markov chain
\index{Markov chain}%
$\tilde{G} = (\Gamma_t)_{0}^{\infty}$, whose states are graphs on $V(\Gamma)$.  Beginning with the \emph{empty graph}
\index{graph ! empty}%
$\Gamma_{0} = (V(\Gamma)= \emptyset, E(\Gamma) = \emptyset)$, form the graph $\Gamma_t$ from $\Gamma_{t-1}$ for $1 \le t \le {n \choose 2}$ by the addition of an edge, all new edges being equiprobable.  Then $\Gamma_t$ has exactly $t$ edges, and if $t = {n \choose 2}$ then $\Gamma_t$ is the complete graph on $n$ vertices, and if $t > {n \choose 2}$ then this is also adopted.

An alternative way of establishing the process of random edge-acquisition is to say that a graph process is a sequence $(\Gamma_t)_{t}^{N} = 0$ such that
\begin{itemize}
\item[(i)]  each $\Gamma_t$ is a graph on $V(\Gamma)$,
\item[(ii)]  $\Gamma_t$ has $t$ edges for $t = 0, 1, \ldots, N$ and
\item[(iii)]  $\Gamma_0 \subset \Gamma_1 \subset \ldots$.
\end{itemize}
If we give all members of the set $\mathcal{\tilde{G}}$ of all $N!$ graph processes the same probability, then we turn it into a probability space.
\index{polynomial algebra}%
  The graph $\Gamma_t$ is called the \emph{state} of the process $\tilde{G} = (\Gamma_t)_{0}^{N}$ at time $t$.  

There is a measure-preserving map
\index{measure-preserving map}%
between $\mathcal{\tilde{G}}$ and the set $\mathcal{G}_{M,n}$ of simple complete $n$-vertex $M$-edge graphs, in which the graphs have equal probability.  This is defined by $\tilde{G} = (\Gamma_t)_{0}^{N} \to \Gamma_{M,n}$, so that we can identify the set of graphs obtained at time $t = M$ with $\mathcal{G}_{M,n}$, that is a random graph
\index{graph ! random}%
 from $\mathcal{G}_{M,n}$ is the same as the state of a graph process $\tilde{G}$ at time $M$.

A \emph{random bipartite graph process}
\index{graph ! random ! bipartite process}%
has also been defined and studied~\cite[p.~171]{bollobas1}. 

It should be mentioned that the graphs in the classical theory often have a high degree of symmetry and succumb to probabilistic, algebraic and analytic methods for dealing with the questions that arise.  More recently, massive and complex real-world graphs have become popular and these have the opposite properties~\cite{chung}, including, 
\begin{itemize}
\item[(i)]  sparsity, (the number of edges is at most a constant multiple of the number of vertices);
\item[(ii)]  small world phenomena, (any two vertices are connected by a short path, and two vertices having a common neighbour are more likely to be neighbours);
\item[(iii)]  power law degree distribution, (the number of vertices with degree $n$ is proportional to $n^{- \beta}$ for a constant $\beta$.
\end{itemize}

A \emph{space of random structures}
\index{space ! of random structures}%
is a sequence of probability spaces.
\index{polynomial algebra}%
  A probability distribution can be attached to the class of $n$-element structures for each $n \in \mathbb{N}$.  The point set of the structure can either be labeled or unlabeled and in the latter case the class is an isomorphism class.  
 
If $r_n$ denotes the proportion of graphs on $n$ vertices having a graph-theoretic property $\mathcal{P}$ and if $\lim_{n \to \infty} r_n = 1$, then we say that \emph{almost every graph has property $\mathcal{P}$}.  The intuition of the probabilistic existence proof of the random graph
\index{graph ! random}%
 given by Erd\H{o}s
\index{Erd\H{o}s, P.}%
and R\'enyi~\cite{er}
\index{R\'enyi, A.}%
is that for increasingly large structures  taken from a space of random structures, the properties of interest will either almost always hold or almost never hold. 
This idea of a \emph{zero-one law}
\index{zero-one law}%
goes back to Kolmogorov,
\index{Kolmogorov, A. N.}%
and in a graph-theoretic context it says that either almost every finite graph with a certain property $\mathcal{P}$ satisfies the law or almost none does.  

In order to make this idea more precise, the property $\mathcal{P}$ is taken to be monotone increasing, that is one for which a graph satisfies $\mathcal{P}$ whenever one of its subgraphs does.  The probabilistic approach~\cite{alonspencer}
\index{Alon, N.}%
\index{Spencer, J.}%
is one of the most effective methods in combinatorics.
\index{combinatorics}%

Let $M_n$ be a random structure of size $n$ over a first-order language
\index{first-order language}%
 $L$ and define for any sentence $\sigma$ of $L$,
\[ Pr_n(\sigma) = Pr(M_n \models \sigma) \]
and if it exists
\[ Pr(\sigma) = \lim_{n \to \infty} Pr_n(\sigma). \]

\bigskip

The zero-one law
\index{zero-one law}%
 for the first-order theory
\index{first-order theory}%
  of random graphs
\index{graph ! random}%
 is due to Glebskii \emph{et al.}~\cite{glebskii}
\index{Glebskii, Y.}%
and Fagin~\cite{fagin}
\index{Fagin, R.}%
and can be stated as a corollary to the following theorem:

\bigskip

\begin{theorem}
Almost every finite graph satisfies a first-order property if and only if $\mathfrak{R}$ does.
\end{theorem}

\begin{proof}
Let $\sigma$ be a first-order sentence
\index{first-order sentence}%
 in the language of graphs.  Denote by $\Psi_k$ the axiom resembling the $k$-e.c. property, with $U$ and $V$ both being sets of at most $k$ vertices.  A probabilistic argument similar to that given in Chapter~\ref{fstchap} for three colours, shows that $\Pr(\Psi_k) = 1$ for each $k$.  A graph $\Gamma$ is e.c., that is satisfies $(\dagger)$, if $\Gamma \models \Psi_k$ for all $k$.  We know that such a countable graph is unique up to isomorphism and is $\mathfrak{R}$.  Therefore the $\Psi_k$ are axioms for $\mathfrak{R}$, because $\mathfrak{R}$ is characterized by ($*$).  Take any sentence $\sigma$ such that $\mathfrak{R} \models \sigma$.  Then $\Psi_k \vdash \sigma$.  So $\{ \Psi_1, \ldots, \Psi_{k_0} \} \vdash \sigma$, where $k_0$ is finite.  Then $\Pr(\sigma) = 1$.

Conversely, if $\mathfrak{R} \nvDash \sigma$, then $\mathfrak{R} \models \neg \sigma$, so $\Pr(\neg \sigma) = 1$, so $\Pr(\sigma) = 0$.

\end{proof}

\begin{corollary}
If $\sigma$ is a first-order sentence
\index{first-order sentence}%
 in the language of graphs, and random graphs are chosen uniformly, then $Pr(\sigma)$ exists and equals either $0$ or $1$.
\end{corollary}
\index{graph ! random}%

\begin{proof}
From the theorem, $Pr(\sigma) = 1$ if and only if $\mathfrak{R} \models \sigma$.

If it is not the case that $\mathfrak{R} \models \sigma$ then $\mathfrak{R} \models \neg \sigma$ and so $Pr(\sigma) = 0$.
\end{proof}

A unique random countable structure in Fagin's
\index{Fagin, R.}%
proof does not always exists, but the proof applies to the uniform space of relational structures.
\index{relational structure}%
  Fagin's proof fails if the language has a function symbol; for example if $c$ is a constant and $W_i$ is a unary relation, then the sentence $c \in W_i$ has probability $\frac{1}{2}$.  There are instances when the zero-one law
\index{zero-one law}%
fails but $\Pr(\sigma)$ exists for each first-order sentence $\sigma$.
\index{first-order sentence}%

More basic results on first-order theory of random graphs are due to Blass and Harary~\cite{blasshar}.
\index{Blass, A.}%
\index{Harary, F.}%

 Cameron and Martins~\cite{cammar} have an application that relates oligomorphic automorphisms to first order sentences:
\begin{theorem}[Cameron--Martins]
Given a finite collection of finite graphs, and the subsets of vertices of a random graph
\index{graph ! random}%
 $\Gamma$ (but not $\Gamma$ itself) that induce those graphs, it is almost always possible to uniquely reconstruct a class of graphs equivalent to $\Gamma$.
\end{theorem}
The uniqueness in the theorem is up to equivalence, referring to one of five natural equivalence relations/automorphisms on graphs given by Thomas' Theorem
 \index{Thomas' Theorem}%
on the classification of the reducts of  $\mathfrak{R}$.

The result is obtained by combining the above theorem of Fagin and the Engeler--Ryll-Nardzewski--Svenonius Theorem~\cite{engeler}~\cite{ryll}~\cite{svenonius}.
\index{Engeler--Ryll-Nardzewski--Svenonius ! Theorem}%

The theory of graphs together with the ($*$)-condition for every natural number has no finite models and all countable models are isomorphic, so the theory is complete~\cite{vaughta}~\cite{gaifman}.

Whilst in our work the vertex-joining probabilities are always constant, an obvious generalization, which has generated a vast amount of mathematics, is assuming a variable edge probability.
 An important result~\cite{spencer}
\index{Spencer, J.}%
is that the probability $p_n(A)$ of a random
graph having a \emph{first-order property}
\index{first-order property}%
 $A$, almost never or
almost surely holds, that is $\lim_{n \to \infty} p_n(A) = 0$ or
$1$.  If $p(n)$ denotes the vertex-dependent connection probability then a
\emph{threshold function}
\index{threshold function}%
 $p_0(n)$ for $A$ is such that if $p(n) \ll
p_0(n)$ then $A$ almost never holds and whenever $p(n) \gg p_0(n)$
then $A$ holds almost surely.  Of particular significance is the phenomenon of phase transitions, an
example of which occurs in the emergence of giant components
\index{giant component}%
 in a random graph
\index{graph ! random}%
 as it grows in size when vertex pairs are joined independently with probability $p$.  For example, denoting by $\Gamma(n, p)$ the
random graph on $n$ vertices with connection probability $p$, which is
also a probability space,
\index{polynomial algebra}%
 Erd\H{o}s and R\'enyi
\index{Erd\H{o}s, P.}%
\index{R\'enyi, A.}%
discovered a global change in the nature of $\Gamma(n, p)$ near $p = 1/n$.  When $p = (1 -
\epsilon) /n$ the components are all small (the largest having size
$O(\ln n)$) and are all trees or
unicyclic, but when  $p = (1 + \epsilon) /n$ a giant component with
far more edges than vertices emerges.  Benny Sudakov
\index{Sudakov, B.}%
  and Michael Krivelevich
\index{Krivelevich, M.}%
 have given~\cite{sudakov} a new simple, proof of the theorem about the threshold for a ``giant component'' of a random graph, using the \emph{depth-first}
\index{depth-first}%
  technique used in computer science.

The book~\cite{spencer} discusses the work of J. Spencer
\index{Spencer, J.}%
and S. Shelah
\index{Shelah, S.}%
and its main result is that when $p = n^{- \alpha}$, if $\alpha$
is irrational then there is a zero-one law,
\index{zero-one law}%
but when $\alpha$ is any rational in $(0,1)$ then there is no such
law.  In fact for some first-order properties,
\index{first-order property}%
 there is an infinite number of threshold functions.  

When the random graph
\index{graph ! random}%
 has a uniform distribution, as has been assumed in this book, then the random graphs have very different properties to the case where random graphs have joining probability is of the form $p = n^{- \alpha}$, for irrational $\alpha$.  The differences are listed by Baldwin~\cite{baldwin}
\index{Baldwin, J.}%
In particular, the random graph with a uniform distribution is unstable (in fact, a prototypical theory with the independence property)
\index{independence property}%
 and is $\aleph_0$-categorical,
\index{aleph@$\aleph_0$-categorical}%
whereas the random graph with $p = n^{- \alpha}$, for irrational $\alpha$ is stable and is not $\aleph_0$-categorical.

J. Cohen
\index{Cohen, J.}%
~\cite{cohenj}
has written an account of the importance of threshold phenomena in
random structures for explaining the physical theory of phase
transitions, where sudden changes of phase result from gradual changes
of a parameter like temperature.  He speculates on possible
applications of random combinatorial structures to telecommunications,
neurobiology and the origin of life.  Thresholds are neither peculiar
to a particular definition of randomness nor to graphs; rather they arise in a
variety of random combinatorial structures in the limit of large
size.  See~\cite[p.~40]{bollobas1} for more on threshold functions.

We must now note that the graphs we study in this monograph have uniform probability, that is constant and independent of $n$, and for these there are no threshold phenomena.  Although there is no known general criterion that characterizes which spaces of random structures satisfy zero-one laws,
\index{zero-one law}%
 Compton~\cite{compton}  
\index{Compton, K.}%
has found that \emph{uniform} spaces of random structures obey zero-one laws as long as the number of such structures does not grow too quickly.  

Furthermore some ``almost all'' theorems cannot be easily obtained.  For example, neither rigidity nor hamiltonicity can be deduced from any first-order property
\index{first-order property}%
 of almost all graphs~\cite{blasshar}, even though it is known that almost all graphs are rigid and hamiltonian.
\index{graph ! Hamiltonian}%
\index{graph ! rigid}%

Let $a(n) = \sum_{n=0}^{\infty} \frac{a_n}{n!} x^n$ be an exponential generating function where $a_n$ is the number of labeled structures on an $n$-element set (see Appendix~\ref{EnumerationandReconstruction}).  Let $b(n) = \sum_{n=0}^{\infty} b_n x^n$ be a generating function where $b_n$ is the number of isomorphism classes of structures of size $n$.  Compton calls a uniform space of random structures \emph{slow-growing}
\index{structure ! slow-growing}%
if $a(x)$ (or $b(x)$) has a non-zero radius of convergence.  Two of his results are

\begin{theorem}
A slow-growing uniform labeled space of random structures closed under disjoint union and taking components has a zero-one law
\index{zero-one law}%
for first-order logic
\index{first-order logic}%
 if and only if, for any $m$,
\[ \lim_{n \to \infty} \frac{a_{n-m} / (n - m)!} {a_n / n!} = \infty.\]
\end{theorem}

\begin{theorem}
A slow-growing uniform unlabeled space of random structures closed under disjoint union and taking components has a zero-one law
\index{zero-one law}%
 for first-order logic if and only if, for any $m$,
\[ \lim_{n \to \infty} \frac{b_{n-m}} {b_n} = 1.\]
\end{theorem}

Uniform graphs are fast-growing, but still have a zero-one law.

\bigskip

There have been many variations on the theme of random graphs.
\index{graph ! random}%
  We survey an example of one such, studied in~\cite{bonatojanssen}.  The \emph{infinite random $n$-ordered graphs}, $R^{(n)}$,
\index{graph ! random ! $n$-ordered}%
is defined as the limit of a certain chain of finite graphs.  Let $R_0 \sim K_n$.  For some $t \ge 0$, define a finite overgraph of $R_0$, by adding a new vertex $x_S$ to each $n$-vertex subgraph $S$ of $V(R_t)$, thus giving $R_{t+1}$.  Then $R^{(n)} = (\cup_{t \in \mathbb{N}} V(R_t), \cup_{t \in \mathbb{N}} E(R_t)) = \lim_{t \to \infty} R_t$.  For the random graph $\mathfrak{R}$, simply add $x_S$ to \emph{every} subset $S$ of $V(R_t)$, and not just those of cardinality at most $n$.

\bigskip

A graph is \emph{strongly $n$-e.c.}
\index{graph ! n@$n$-existentially closed ! strongly}%
if for every $n$-element subset $U$, and finite $W$ of the vertex set $V(\Gamma)$ of the graph, there is a vertex $z\in V(\Gamma) \backslash (U\cup W)$ such that all edges between $v$ and $U$ are present and no edges between $v$ and $W$ are present.  A strongly $n$-e.c. graph is infinite and strongly $m$-e.c. for all $m < n$.

\bigskip
\bigskip

\begin{theorem}  Fix a positive integer $n$.
\begin{itemize}
\item[(1)]  The graph $R^{(n)}$is strongly $n$-e.c, but not $(n+1)$-e.c;
\item[(2)]  If $\Gamma$ is a strongly $n$-e.c. graph, then $R^{(n)}$is an induced subgraph of $\Gamma$;
\item[(3)]  $R^{(n)} \leq R^{(n+1)}$;
\item[(4)]  $\lim_{n \to \infty} R^{(n)} = \mathfrak{R}$;
\item[(5)]  $\Aut(R^{(n)})$ embeds all countable groups, including  $\Sym(\omega)$. 
\end{itemize}
\end{theorem}

In a variation on this random graph process, new vertices have a higher probability of being added to older than younger vertices. 

\bigskip

The notion of \emph{graph limits}
\index{graph ! limit}%
 which was introduced by Borgs, Chayes, Lov\'asz, S\'os and Vesztergombi (see~\cite{borgs})
\index{Lov\'asz, L.}%
\index{Borgs, C.}%
\index{Chayes, J.}%
\index{S\'os, V. T.}%
\index{Vesztergombi, K.}%
 and further studied in~\cite{lovasz} and~\cite{diaconisjan} by B. Szegedy, P. Diaconis and S. Janson.
\index{Szegedy, B.}%
\index{Janson, S.}%
\index{Diaconis, P.}%
 The limits of sequences of (dense) graphs
 \index{graph ! dense}%
 are not graphs, but rather symmetric measurable functions $W: [0, 1]^2 \to [0, 1]$~\label{W} (the function $W$ is \emph{symmetric}
 \index{symmetric function}%
if $W(x, y) = W(y, x)$).  This limit object determines all the limits of subgraph densities and conversely each such function arises as a limit object.  The notion is roughly that a growing sequence of finite graphs $\Gamma_n$ converges if, for any fixed graph $\Gamma$, the proportion of copies of $\Gamma$ in $\Gamma_n$ converges.  If $\gamma$ is a fixed simple graph and $\homo(\gamma, \Gamma)$~\label{homo} denotes the number of homomorphisms of $\gamma$ to $\Gamma$ then the \emph{homomorphism density}
\index{homomorphism density}%
defined by~\label{tgraph}
\[ t(\gamma, \Gamma) = \frac{\homo(\gamma, \Gamma)}{|V(\Gamma)|^{|V(\gamma)|}} \]
is the probability that a random mapping $V(\gamma) \to V(\Gamma)$ is a homomorphism.  Both $t(\cdot, \Gamma)$ and $\homo(\cdot, \Gamma)$ are multiplicative; a graph parameter $f$ is \emph{multiplicative}
\index{graph ! multiplicative parameter}%
if $f(\Gamma_1 \sqcup \Gamma_2) = f(\Gamma_1) f(\Gamma_2)$
where $\Gamma_1 \sqcup \Gamma_2$ denotes the disjoint union of two graphs $\Gamma_1$ and $\Gamma_2$.  Whilst the $n$-vertex random graphs
\index{graph ! random}%
 converge with probability 1 uniquely up to automorphism, the limit object for random graphs of density $p$ on $n$ nodes is the constant function $p$, and the homomorphism densities have different limits depending on density.  Every such $W$ gives rise to a general model of random graphs, called \emph{$W$-random},
\index{graph ! wrandom@$W$-random}%
 and every random graph model satisfying some natural criteria can be obtained in this way for an appropriate $W$.

By defining for a simple graph $\gamma$
\[ t(\gamma, W) = \int_{{[0, 1]}^{V(\gamma)}} \prod_{i, j \in E(\gamma)} W(x_i, x_j)\,dx,\] 
the homomorphism density can be thought of as a ``moment'' of $W$.  Borgs et al prove~\cite{borgs1} that every bounded symmetric function $W$ is determined by its moments up to a measure preserving transformation of the variables.

It was also discovered~\cite{diaconisjan} that the characterization of Lov\'asz et al of the graph limits is essentially equivalent to the characterization of \emph{exchangeable random infinite graphs}
\index{graph ! random ! exchangeable}%
 by Aldous and Hoover.
\index{Aldous, D.}%
\index{Hoover, D.}%
  A random infinite graph in the set of all labeled countable graphs, is \emph{exchangeable} if its distribution is invariant under every permutation of the vertices; it is sufficient here to consider only finite vertex permutations.  Exchangeable random graphs
\index{graph ! random ! exchangeable}%
   converge to a limiting object which may be thought of as a probability measure on infinite random graphs.  There is a $1$--$1$ correspondence between infinite exchangeable random graphs and distributions on the space of proper graph limits.

Our work has focussed on the classical random graph and its derivatives, introduced by Erd\H{o}s and A. R\'enyi~\cite{er}.
\index{Erd\H{o}s, P.}%
\index{R\'enyi, A.}%
 These are homogeneous in the sense that all their vertices come in one type, and the vertex degrees
tend to be concentrated around a typical value.  Recently `inhomogeneous' random graph models
\index{graph ! random ! inhomogeneous}%
 have been studied with a view to modeling real-world phenomena, for example in~\cite{soderberg}.  Others have power-law distributions.  Amongst these are the inhomogeneous models of sparse graphs introduced~\cite{bollobas2} by Bollob\'as, Janson and Riordan
\index{Bollob\'as, B.}%
\index{Janson, S.}%
\index{Riordan, O.}%
 which are closely related to the convergent sequences of dense graphs.
 \index{graph ! dense}%

In the context of the above work on graph limits, Freedman, Lov\'asz, and Schrijver,
\index{Lov\'asz, L.}%
\index{Freedman, M.}%
\index{Schrijver, A.}%
introduced the concept of a \emph{graph algebra}~\cite{freedman}.
\index{graph ! algebra}%
 
Another structure that has arisen is the notion of \emph{graphon}
\index{graphon}%
 introduced by Lov\'asz and Szegedy in~\cite{lovasz}.  Razborov~\cite{razborov}~\cite{razborov1}
\index{Razborov, A. A.}%
has an alternative formulation in terms of the related concept of \emph{flag algebras}~\cite{razborov}.
\index{flag algeras}%
A third equivalent formulation is Tao's notion of a permutation-invariant measure
\index{measure ! invariant}%
 space~\cite{tao1}~\cite[p.~171]{tao}, which takes a sequence of increasingly large but still dense graphs whose limiting object is a probability space 
\index{probability space}%
 $(X, \mathcal{B}, \mu)$ together with a $\Sym(\infty)$-action on $\mathbb{Z}$.
\index{group ! symmetric ! infinite}%

We summarize the work of David Aldous~\cite{aldous}
\index{Aldous, D.}%
 and co-workers on exchangeability and continuum limits of discrete random structures.

For any kind of object, there is a \emph{random} such object.  An $S$-valued random variable $X$ has a distribution $\dist(X)$, the induced probability measure on $S$.  A probability measure on $S^{\infty}$ is \emph{exchangeable} if it is invariant under the action of each map $\tilde{\pi} : S^{\infty} \to S^{\infty} : (s_i) \mapsto (s_{\pi(i)})$, where $\pi$ is a finite permutation of $\mathbb{N} = \{1, 2, \ldots\}$.  A sequence of random variables obeys this strong symmetry condition of exchangeability if order is not important.  Any countable mixture of independent indentically-defined (IID) sequences is exchangeable; de Finetti's theorem
\index{de Finetti's theorem}%
\index{de Finetti, B.}%
 is the converse.

One way of examining a structure with a probability measure is to sample $k$ IID random points (say graph vertices) and consider an induced $k$-element substructure relating the random points.  Within the limit of random structures, the $k$-elements are exchangeable and the distributions are consistent as $k$ increases, and so can be used to define an infinite structure.  Thus, exchangeable representations of random structures provide a way to derive the $n \to \infty$ continuum limit of discrete random finite $n$-element structures.

Aldous goes on to show how the theory can be applied to continuum scale-invariant random spatial networks, and in particular how beginning with a square lattice it is possible to refine it so that it can be used to define (by continuity) routes between points in $\mathbb{R}^2$.  External randomization gives the process further invariance properties.

\bigskip

Now that we have defined random graphs, we can discuss them in the context of ages and ubiquity.
 
Let $\mathcal{Y}(\mathfrak{R})$ be the set of all graphs on the vertex set $\mathbb{N}$.  The Baire category version of the Erd\H{o}s--R\'enyi
\index{Erd\H{o}s, P.}%
\index{R\'enyi, A.}%
statement that almost all countable graphs are isomorphic to $\mathfrak{R}$, is:
\begin{theorem}
The set of graphs on $\mathbb{N}$ which are isomorphic to $\mathfrak{R}$ is residual in $\mathcal{Y}(\mathfrak{R})$.
\end{theorem}
\begin{proof}
For given finite disjoint sets $U$ and $V$, the set $\mathcal{S}(U, V)$ of graphs for which there exists $z$ joined to every vertex in $U$ and to none in $V$, is open and dense
\index{dense open set}%
in $\mathcal{Y}(\mathfrak{R})$.  The result will then follow, because a graph is isomorphic to $\mathfrak{R}$ if and only if it lies in the intersection of all these sets.

$\mathcal{S}(U, V)$ \emph{is open}:  If $S$ is a graph in this set, $z$ is the witnessing vertex, and $m = \max(U \cup V \cup \{z\})$, then any graph agreeing with $S$ on $\{1, \ldots, m\}$ is in $\mathcal{S}(U, V)$.

$\mathcal{S}(U, V)$ \emph{is dense}:  If we have a graph on $\{1, \ldots, n\}$, where without loss of generality $\max(U \cup V) \le n$, then we can join the next vertex $n + 1$ correctly to $U$ and not $V$, and guarantee membership in $\mathcal{S}(U, V)$.
\end{proof}

The principle behind this proof can be generalized.  A sentence of first-order logic
\index{first-order logic}%
 is called \emph{inductive}
\index{sentence ! inductive}%
if for a quantifier-free $\sigma$, it takes the form $(\forall x_1, \ldots, x_m)(\exists y_1, \ldots, y_n)\sigma$.

\begin{theorem}
Let $M$ be a countable relational structure,
\index{relational structure}%
 and $\Sigma$ a countable set of inductive sentences in the language of $M$.  Then the set of models of $\Sigma$ is residual in $\mathcal{Y}(M)$.
\end{theorem}
\begin{proof}
A countable intersection of residual sets is residual, so it suffices to prove this for a single sentence.  Now follow the argument used above for the inductive sentence $\sigma_{m, n}$ in $\mathfrak{R}$.
\end{proof} 

\begin{corollary}
Let $M$ be a countable relational structure.  Then the set of structures $S$ satisfying Age$(S)$ = Age$(M)$ is residual in $\mathcal{Y}(M)$.
\end{corollary}

This corollary states that almost all structures younger than $M$ actually have the same age.  It is proved by taking, for each finite structure $X$ of $M$, the existential sentence saying that $X$ occurs as a substructure.

\begin{corollary}
If $M$ is a countable homogeneous relational structure, then $M$ is ubiquitous in category.
\end{corollary}

Since almost all structures in $\mathcal{Y}(M)$ have the same age
\index{age}%
 as $M$, we merely observe that $M$ is characterized by the I-property, which can be expressed as a countable number of inductive sentences.

A relational structure
\index{relational structure}%
 $M$ is \emph{universal}
\index{structure ! universal}%
if every structure which is younger than $M$ is embeddable in $M$.  If $S$ is embeddable in $M$ then $S$ is younger than $M$, but the converse is false.  For example, if $M$ is a two-way infinite path then Age$(M)$ consists of all disjoint unions of finite paths and $\mathcal{Y}(M)$ consists of all disjoint unions of finite or countable paths.  The union of two two-way infinite paths is not embeddable in $M$.  However the following holds~\cite{cam3}:
\begin{theorem}
If a countable relational structure is either homogeneous or $\aleph_0$-categorical then it is universal.
\end{theorem}
\index{aleph@$\aleph_0$-categorical}%
\bigskip

One graph invariant that has drawn attention is the distinguishing number, introduced in~\cite{albertson}.

The \emph{distinguishing number}
\index{graph ! distinguishing number}%
 of a graph $\Gamma$ is the smallest positive integer $r$ such that $\Gamma$ has a labeling of its vertices with $r$ labels for which there is no non-trivial automorphism of $\Gamma$ preserving these labels.  Imrich, Klav\v{z}ar and Trofimov~\cite{imrich} computed the distinguishing number of some infinite graphs, showing in particular that $\mathfrak{R}$ has distinguishing number 2.   This result has been generalised by Laflamme et al.~\cite{laflamme}, who show that the automorphism groups
\index{group ! automorphism}%
  of many countable homogeneous relational structures
\index{relational structure}%
   have distinguishing number either 2 or $\infty$, including countable homogenous graphs (both directed and undirected), the universal poset, and the countable dense
\index{linear order}%
 linear order $(\mathbb{Q}, \leq)$. Moreover, they show that any countable homogeneous structure satisfying the \emph{Free Amalgamation Property}
\index{amalgamation property ! free}%
 (see Appendix~ref{TheoryofRelationalStructures}) has distinguishing number 2.  Bonato and Deli\'c~\cite{bondelic}
\index{Bonato, A.}
\index{Deli\'c, D.}
 further showed that any countable relational structure satisfying a particular adjacency property (which they call the \emph{weak-e.c.property})
\index{existentially closed (e.c.) ! weak}%
 has distinguishing number 2.  More on this topic can be found in the paper by Robert Bailey and Peter Cameron~\cite{bailcam}.
\index{Bailey, Robert F.}
\index{Cameron, P. J.}

The subject of random graphs has been an exciting research area for over half a century, and it would require an enormous treatise to cover the variations on the theme together with applications. 

\bigskip

We end the section with some comments on the automorphism group
\index{group ! automorphism}%
 of the random graph, reserving most of the discussion for the main text.  This group embeds all finite or countable groups.  All its cycle structures of were classified by Truss~\cite{truss3}.
\index{Truss, J. K.}%
The following two theorems~\cite{bhatmaca}~\cite{bhatmac} of Bhattacharjee
\index{Bhattacharjee, M.}%
and Macpherson
\index{Macpherson, H. D.}%
demonstrate remarkable properties of subgroups of $\Aut(\mathfrak{R})$.
\begin{theorem}
There exist $g_1, g_2 \in \Aut(\mathfrak{R})$ such that
\begin{itemize}
\item[(a)]  $g_1$ has a single cycle on $\mathfrak{R}$ which is infinite;
\item[(b)]  $g_2$ fixes a vertex $v$ and has two cycles on the remaining vertices (namely, the neighbours and non-neighbours of $v$);
\item[(c)]  the group $\langle g_1, g_2 \rangle$ is free and is transitive on vertices, edges, and non-edges of $\mathfrak{R}$, and each of its non-identity elements has only finitely many cycles on $\mathfrak{R}$.
\end{itemize}
\end{theorem}

This theorem is proved by building the permutations $g_1$ and $g_2$ as limits of partial maps constructed in stages.  Parts (a) and (b) follow from Truss' classification
\index{Truss, J. K.}%
of cycle types, about which more will be said in the main text.  To prove (c) requires a lot of work.

\begin{theorem}
There is a locally finite group $G$ of $\Aut(\mathfrak{R})$ which acts homogeneously (that is, any isomorphism betwen finite subgraphs can be extended to an element of $G$).
\end{theorem}

This theorem uses a result of Hrushovski~\cite{hrush}
\index{Hrushovski, E.}%
on extending partial automorphisms of graphs.
\index{automorphism ! partial}%
\bigskip

We have had a fair bit to say about the concept of switching of graphs, and we end this section with an application of the classical theory of Seidel switching.
\index{Seidel switching}%
\index{switching! Seidel}%

There is a problem in extremal graph theory linked to a problem in classical geometry through the work of Gromov.
\index{Gromov, M.}%
  I. B\'ar\'any
\index{B\'ar\'any, I.}%
 showed~\cite{barany} that, for any dimension $d$, there is a constant $c_d$ such that, given any finite set of points in general position in $d$-dimensional Euclidean space,
\index{Euclidean space}%
 there is a point (not necessarily in the set) which is contained in a proportion at least $c_d$ of all the $d$-simplices spanned by the given points.  The problem is to calculate the value of $c_d$, which is known to be $2 / 9$ for $d=2$, but little more is known.

Gromov introduced a topological method, which depended on a function $f_d$.  The definition for $d=2$ which we give presently, easily generalizes and depends on Seidel switching of graphs.  A graph is called \emph{Seidel-minimal}
\index{graph ! Seidel-minimal}%
 if it has the smallest number of edges of any graph in its Seidel switching class (equivalently, given any cut, there are at least as many non-edges as edges crossing it).  Any graph has a two-graph,
\index{graph ! two-graph}%
\index{two-graph}%
 or coboundary, which is an invariant of the switching class.  Now $f_2(\alpha)$ is defined to be the limit inferior of the density of the coboundary of a Seidel-minimal graph of density at least $\alpha$, where here density means the proportion of edges or triples out of the total possible.  Gromov gave a lower bound for $c_d$ depending on the functions $f_e$ for all $e \le d$.  So better bounds on $f_2$ could improve lower bounds for $c_d$ for all $d>2$.

This is exactly what B\'ar\'any and L. Mach have done, using flag algebra
\index{flag algeras}%
 computations.  They have improved the lower bound for $c_3$ from 0.06332 to 0.07509. (The best upper bound is 0.09375.)  See~\cite{kral} for references including the work of Gromov.

\bigskip

Bodirsky and Pinsker
\index{Bodirsky, M.}%
\index{Pinsker, M.}%
 have produced a detailed study~\cite{bodirsky} of what they call \emph{minimal functions} on the random graph, proving using a Ramsey-type theorem,
 \index{Ramsey theory}%
  that there is a system of 14 non-trivial finitary functions on the random graph such that any nontrivial function on the random graph generates one of the functions of this system by composition with automorphisms and by topological closure, and that the system is minimal in the sense that no subset of the system has the same property.  
 
They derive Thomas' Theorem on reduct
\index{reduct}%
  classification
\index{Thomas' Theorem}%
 and prove some refinements of this theorem.  They also classify the minimal reducts closed under \emph{primitive positive}
\index{formula ! primitive positive}%
  definitions, and prove that all reducts of the random graph and of the linear order
\index{linear order}%
 of the rationals $(\mathbb{Q}; <)$ are model-complete.  
\index{theory ! model-complete}%
  A structure is model-complete
\index{structure ! model-complete}%
   if and only if its first-order theory is model-complete.

\bigskip

The substantial treatise by Bollob\'as~\cite{bollobas1}
\index{Bollob\'as, B.}%
deals with the combinatorial theory of random graphs, and the book by Spencer~\cite{spencer}
\index{Spencer, J.}%
is an account of the interplay of logic, probability theory and random graphs.  The article by Winkler~\cite{winkler}
\index{Winkler, P.}%
has more on the relationship between random structures and zero-one laws
\index{zero-one law}%
including those in higher logics than first-order.  The article by Baldwin~\cite{baldwin}
\index{Baldwin, J.}%
summarizes the use of logic for studying limit laws on finite probability spaces,
\index{polynomial algebra}%
 and shows how Urysohn space
\index{Urysohn space}%
arises as the limit of finite spaces with probability measures as compared with Vershik's
\index{Vershik, A. M.}%
construction in which measures lie on infinite sets.


\section{Topology in Permutation Groups}
\label{TopologyinPermutationGroups}

To a metric space $(X, d)$ is associated a \emph{natural topology}
\index{topology ! natural}%
$T_d(X)$ in which a set $O \subset X$ is called \emph{open} if and only if for every $x \in O$ there is a real number $r > 0$ such that
\[ B(x, r) := \{ y \in X : d(y, x) < r\} \subset O,\]
where $B(x, r)$ is the \emph{open ball}
\index{open ball}%
of radius $r$ centred at $x$.  So the \emph{topological space} $(X, T_d(X))$
\index{topological space}%
is associated with the metric space $(X, d)$.
 
The natural topology on permutation groups is the \emph{topology of pointwise convergence}
\index{topology ! of pointwise convergence}%
in which for a sequence $g_n$ of permutations on $X$, $\lim_{n \to
  \infty} (g_n) = g$ if and only if $\forall x_i \in X,\ \exists n_0
\in \omega$ such that $\forall n > n_0,\ x_i g_n = x_i g$.  At the
heart of the topology of pointwise convergence
\index{topology ! of pointwise convergence}%
 is the idea that permutations agree on points that are close to each other. The space
of zero-one sequences
\index{zero-one sequence}%
 is a complete metric space and its topology is
that of pointwise convergence.

Endowing a group $G$ with this topology turns it into a \emph{topological group},
\index{group ! topological}%
that is multiplication and inversion are continuous, so that if $g_n \to g$ and $h_n \to h$ then $g_n h_n \to gh$ and $g_n^{-1} \to g^{-1}$.

A basis of open sets
\index{open set}%
 in this topology consists of the cosets of
pointwise stabilizers of finite tuples; the subgroups involved are
called \emph{basic open subgroups}.
\index{group ! basic open}%
The topology can be derived from a metric.  The distance function defined by:
 
\begin{displaymath}
d(g, h) = \left\{ \begin{array}{ll}
0 & \text{if $g = h;$} \\
\frac{1}{2^i} & \text{if $jg = jh$ and $jg^{-1} = jh^{-1}$ for all $j < i$,} \\
  & \text{but $ig \neq i h$ and $ig^{-1} \neq ih^{-1}$,}
 \end{array} \right.
\end{displaymath}
turns $\Sym(X)$ into a complete metric space.

\medskip

The \emph{closure}
\index{group ! permutation ! closure}%
of a permutation group is the set of all permutations which are limits of sequences of group elements.  The closure of $G < \Sym(X)$ acting on $\Omega$ is the largest group having the same orbits on $\Omega^n$ as $G$ for all $n$.  An \emph{automorphism} of a first-order structure
\index{first-order structure}%
 $M$ is a permutation $g$ that preserves all relations, fixes constants and commutes with functions by $f(x_1g, \ldots, x_ng) = f(x_1, \ldots, x_n)g$.  Denote the group of such automorphisms by $\Aut(M)$.
\index{group ! automorphism}%
The topologies defined on $\Sym(X)$ are the same as that defined on
automorphism groups.
\index{group ! automorphism}%

A space admitting a countable dense set
\index{dense set}%
is called \emph{separable}.
\index{topological space ! separable}%
If $M$ is a countably infinite first-order structure and $G = \Aut(M)$ then either 

$|G| \le \aleph_0$,\\
or

$|G| = 2^{\aleph_0}$, the first alternative holding if and
only if the stabilizer of some tuple is the identity.  Because $M$ is
countable this topology is metrizable (via the above distance function), and $G$ is a complete separable
metric space.

In the latter case where $|G| = 2^{\aleph_0}$, the identity and hence every point is a limit point.  In the first case, $|G|$ is equal to the number of tuples in its orbit, so is at most countable and the induced topology on $G$ is \emph{discrete},
\index{topology ! discrete}%
that is every subset is open; it follows that every subset is also closed.  Taking the discrete topology
\index{topology ! discrete}%
on a set $X$ means that a sequence of elements of $X$ converges if and only if its limit is constant.

Closed subgroups of the symmetric group have tamer properties than arbitrary subgroups.  Evans
\index{Evans, D. M.}%
has proved~\cite{evans} that if $G$ and $H$ are closed subgroups of $\Sym(X)$, and $H \le G$, then either $|G : H| \le \aleph_0$ or $|G : H| = 2^{\aleph_0}$, and the former holds if and only if $H$ contains the stabilizer in $G$ of some tuple.   The following result is proved in~\cite{cameron}:
\begin{theorem}
Let $G$ be a permutation group acting on a countable set $\Omega$.  The following are equivalent:
\begin{itemize}
\item[(a)] $G$ is a closed subgroup of $\Sym(\Omega)$;
\item[(b)] $G = \Aut(M)$ for some first-order structure $M$ on $\Omega$;
\item[(c)] $G = \Aut(M)$ for some relational structure $M$ on $\Omega$.
\end{itemize}
\end{theorem}
\index{first-order structure}%
\index{relational structure}%

\medskip

The following~\cite{cam4} lists some of the topological properties of subgroups.  
\begin{theorem}
\label{toplist}
Let $G < \Sym(\mathbb{N})$.  Then,
\begin{itemize}
\item[(a)] $G$ is open if and only if it contains the pointwise stabilizer of a finite set;
\index{group ! open}%
\item[(b)] $G$ is closed if and only if it is the automorphism group
\index{group ! automorphism}%
 of some first-order structure on $\mathbb{N}$;
\index{group ! closed}%
\item[(c)] $G$ is discrete if and only if the pointwise stabilizer of some finite set is trivial;
\index{group ! discrete}%
\item[(d)] $G$ is compact if and only if it is closed and all its orbits are finite;
\index{group ! compact}%
\item[(e)] $G$ is locally compact if and only if it is closed and all orbits of the stabilizer of some finite set are finite.
\index{group ! locally compact}%
\end{itemize}
\end{theorem}

Another characterisation of a topological group is~\cite{kayem}

\begin{theorem}
\label{topaut}
\index{topological group ! Hausdorff}%
Let $G$ be a topological group.  Then $G = \Aut(M)$ for some countable
structure $M$ if and only if 
\begin{itemize}
\item[(a)] $G$ is Hausdorff;
\item[(b)] $G$ is complete;
\item[(c)] if $H < G$ is open then $|G : H| \le \aleph_0$;
\item[(d)] there is a countable family $\{H_i : i \in \mathbb{N}\}$ of
  subgroups of $G$ such that $\mathcal{B} = \{H_{0}^{g_0} \cap \ldots \cap
  H_{k}^{g_k} : k \in \mathbb{N}, g_0, \ldots, g_k \in G\}$ is a base
  of open subgroups of $G$, i.e. the set of cosets of elements of
  $\mathcal{B}$ forms a base for the topology on $G$.
\end{itemize}
\end{theorem}

The topology defined at the beginning of this section varies according
to the precise basic open subgroups chosen.  It is always Hausdorff
\index{topology ! Hausdorff}%
(so any Cauchy sequence
\index{Cauchy sequence}%
converges to at most one point) but not in
general compact.  If the automorphism group
\index{group ! automorphism}%
 of a structure is oligomorphic then the topology is not locally compact.

\medskip

Bergman
\index{Bergman, G.}%
and Shelah~\cite{bergshel}
\index{Shelah, S.}%
have classified the closed subgroups of the countably infinite symmetric group $\Sym(\Omega)$,
\index{group ! symmetric ! infinite}%
 as follows:
\begin{theorem}
For $G_1, G_2\leq \Sym(\Omega)$ write $G_1\approx G_2$ if there exists a finite set $U\subseteq \Sym(\Omega)$ such that $\langle G_1\cup U \rangle = \langle G_2\cup U \rangle$.  The subgroups closed in the function topology on $\Sym(\Omega)$ lie in one of four equivalence classes under this relation, depending on which of the following statements about pointwise stabilizer subgroups $G_{(S)}$ of finite subsets $S \subseteq\Omega$ holds: 
\begin{itemize}
\item[(i)]  For every finite set $S$, the subgroup $G_{(S)}$ has at least one infinite orbit in $\Omega$;
\item[(ii)]  There exist finite sets $S$ such that all orbits of $G_{(S)}$ are finite, but none such that the cardinalities of these orbits have a common finite bound. 
\item[(iii)]  There exist finite sets $S$ such that the cardinalities of the orbits of $G_{(S)}$ have a common finite bound, but none such that $G_{(S)}=\{1\}$. 
\item[(iv)]  There exist finite sets $S$ such that $G_{(S)}=\{1\}$.
\end{itemize}
\end{theorem}
Note that $\Aut(\mathfrak{R})$ satisfies case (i); in fact, for any finite set $S$ all the orbits of $G_{(S)}$ outside $S$ are infinite, and have the property that the induced subgraph is isomorphic to $\mathfrak{R}$ and $G_{(S)}$ acts homogeneously on it.  (If $|W| = n$, then the pointwise stabilizer of $W$ has $2^n$ orbits outside $W$: for each subset $U \subset W$, there is an orbit consisting of the points witnessing condition $(*)$ for $(U, W \backslash U)$).  So~\cite[Theorem 1.1]{macneu} there is a finite set $B$ such that $\langle \Aut(\mathfrak{R}), B \rangle = \Sym(\mathfrak{R})$.  From~\cite{galvin}, $B$ can be a single element.

\medskip

We now turn to reconstruction problems in model theory:
\index{model theory}%
 what can we know about a saturated structure $M$ or its theory $\Th(M)$ from its automorphism group
\index{group ! automorphism}%
 $\Aut(M)$?  As the next theorem~\cite{hodgesa} indicates, we can pin down the theory up to \emph{bi-interpretability},
\index{theory ! bi-interpretable}%
a concept which we now define.

 A \emph{permutation structure}
\index{permutation structure}%
is a pair $\langle W; \Aut(W) \rangle$
where $W$ is a non-empty set~\cite{evans1}.  Two such structures are
\emph{bi-interpretable}
\index{bi-interpretable structure}%
\index{structure ! bi-interpretable}%
 if their automorphism groups
\index{group ! automorphism}%
  are isomorphic as topological groups.
\index{topological group}%
 Two countable totally categorical structures $M$ and $N$ are called \emph{bi-interpretable},
\index{structure ! bi-interpretable}%
that is they have $\Aut(M) \cong \Aut(N)$ as topological groups if and only if there are interpretations $\iota_1$~\label{iotaint} of $M$ in $N$ and $\iota_2$ of $N$ in $M$, both without parameters, such that the natural isomorphisms $M \cong \iota_1\iota_2(M)$ and $N \cong \iota_2 \iota_1(M)$ are definable without parameters in $M$ and $N$ respectively.

\begin{theorem}
If $T$ and $T'$ are $\aleph_0$-categorical theories and $M$ and $N$ are their countable models then the following are equivalent:
\begin{itemize}
\item[(a)]  there is a bicontinuous isomorphism from $\Aut(M)$ to $\Aut(N)$;
\item[(b)]  $T$ and $T'$ are bi-interpretable.
\end{itemize}
\end{theorem}
Therefore $\Aut(M)$ as a topological group tells us everything about $T$.
Let $M$ and $N$ be relational structures.
\index{relational structure}%
  If $M$ and $N$ are saturated models of the same theory $T$, then it is possible to reconstruct $\Aut(M)$ and $\Aut(N)$ as topological groups from each other.  It is possible to reconstruct~\cite{lascar} the category of all models of a countably saturated model $M$ from the topological group $\Aut(M)$, with elementary maps
\index{elementary map}%
 as morphisms.  (A map between two $L$-structures $M$ and $N$ is \emph{elementary} if, for any tuple $\overline{a}$ in $M$ and formula $\phi(\overline{x})$ of $L$, if $\phi(\overline{a})$ holds in $M$, then $\phi(\overline{a}f)$ holds in $N$).

 A topological group is the automorphism group
\index{group ! automorphism}%
\index{topological group}%
  of a first-order structure
\index{first-order structure}%
  if it satisfies the criteria of Theorem~\ref{topaut}.  It
 is possible to reconstruct an $\aleph_0$-categorical structure from
 its automorphism group
\index{group ! automorphism}%
  as a permutation group.  From the work of
 D. Evans
\index{Evans, D. M.}%
and P. Hewitt~\cite{evhe}
\index{Hewitt, P. R.}%
it is not possible to reconstruct an $\aleph_0$-categorical theory from its pure automorphism group;
\index{group ! automorphism}%
 they exhibited two $\aleph_0$-categorical structures whose automorphism groups are isomorphic as abstract groups but not as topological groups.
\index{aleph@$\aleph_0$-categorical}%
\index{topological group}%

Classifying totally categorical structures up to their abstract automorphism group
\index{group ! automorphism}%
 is the same as classifying them up to their topological automorphism group,
because the topology may be obtained completely from an abstract
group.  This follows from the fact that such structures have the
\emph{small index property},
\index{small index property (SIP)}%
because the topology on $\Aut(N)$ can be
recovered from $\Aut(N)$ as an abstract group~\cite{hodges1}.  Before we give a brief introduction to this property, we mention another approach to reconstruction of various structures, which is that of M. Rubin~\cite{rubin}.
\index{Rubin, M.}%
He found a condition related to the definability, in $\Aut(N)$, of point stabilizers, which implies that $\Aut(N)$ determines $N$ up to bi-interpretability,
\index{bi-interpretable structure}%
 or in some instances, up to bi-definability; such a structure is said to have a \emph{weak $\forall\exists$ interpretation}.
\index{weak $\forall\exists$ interpretation}%
Slightly weaker versions of these interpretations for a range of relational structures are known.
\index{relational structure}%

\medskip

A topological group
\index{topological group}%
 $G$, or structure $M$ where $G = \Aut(M)$, is said to have the \emph{small index property} (SIP)
\index{small index property (SIP)}%
if for all subgroups $H < G$, $H$ is open if and only if $|G : H| \le \aleph_0$.  

For an infinite group $G$ and an $n$-tuple $\bar{\mu} := (1, \ldots, n)$, write
$\bar{\mu} g$ for $(1g, \ldots, ng)$ and $G_{\bar{\mu}} := \{g \in G : \bar{\mu} g = \bar{\mu} \}$.
The group $G$ is a topological group
\index{topological group}%
 with its basic open sets
\index{sets ! basic open}%
 being the cosets of the $G_{\bar{\mu}}$, for all $n$ and $\bar{\mu}$.  Left
and right cosets are equivalent here, for $G_{\bar{\mu}g} =
g^{-1}G_{\bar{\mu}}g$.   

A subset of a topological space is \emph{dense}
\index{topological space ! dense}%
if it meets every nonempty open set.
\index{open set}%
  A \emph{Polish space}
\index{Polish space}%
is a topological space that is separable and completely metrizable,
\index{topological space ! completely metrizable}%
meaning that its topology is induced by a complete metric.  A topological group is \emph{Polish}
\index{group ! topological ! Polish}%
if the underlying topology of the group is that of a Polish space, with the open subgroups forming a base of open neighbourhoods of $1$.  For each $\bar{\mu}$, there is a bijection between the
right cosets of $G_{\bar{\mu}}$ in $G$ and the set $\{\bar{\mu} g : g \in
G\}$.  So $G_{\bar{\mu}}$ has small index in $G$, and so therefore does
any open subgroup of $G$.  For a topological group each of whose subgroups have SIP, the open subgroups are identifiable by their index alone, and then the open sets
\index{open set}%
 are unions of cosets of open subgroups.  So the topology
\index{topology}%
  is determined by the group structure.

Suppose $G_1$ and $G_2$ are Polish groups, $G_1$ has the small index property and $\phi : G_1 \to G_2$.  If $\phi$ is a homomorphism of abstract groups then it is continuous; if $\phi$ is an isomorphism of abstract groups then it is a homeomorphism and $G_2$ also has the small index property.

Cherlin
\index{Cherlin, G. L.}%
and Hrushovski
\index{Hrushovski, E.}%
discovered that there are structures that fail the SIP~\cite{kayem}~\cite{lascar}.  Examples of such groups are any meagre subgroup (that is one with comeagre complement) of a Polish group~\cite{hodges},
\index{group ! Polish}%
 or any of its cosets, all of which have index $2^{\omega}$ in the overgroup.  In~\cite[p.~108]{cameron} another structure is given, consisting of an interchangeable
complementary pair of $k$-uniform hypergraphs
\index{hypergraph}%
 (each hyperedge being a $k$-vertex subset) which has an automorphism group
\index{group ! automorphism}%
  with quotient $\mathbb{Z}_{2}^{\omega}$.

Countable totally categorical structures necessarily have the
SIP (see~\cite[p.~115]{kayem}.)  A consequence of failure of SIP is that there
will exist two models of the relevant structure whose
automorphism groups are isomorphic as abstract groups but not as
topological groups,
\index{topological group}%
 so the two models fail to be bi-interpretable.  
\index{bi-interpretable structure}%

W. Hodges
\index{Hodges, W.}%
\emph{et al.}~\cite{hodges} proved that the countable random graph
\index{graph ! random}%
 has the SIP, and Cameron
\index{Cameron, P. J.}%
proved~\cite{cam23} that it also has the \emph{strong small index property}
\index{small index property, strong}%
meaning that any subgroup of $\Aut(\mathfrak{R})$ with index less than $2^{\aleph_0}$ lies between the pointwise and setwise stabilizers of some finite tuple.

\smallskip

Gartside and Knight
\index{Gartside, P. M.}%
and 
\index{Knight, R. W.}%
 define~\cite{gartside} a Polish group $G$ to be \emph{almost free}
\index{group ! almost free}%
if a residual subset of the $n$-tuples of elements of G freely generate a free group.  They give a number of equivalent characterisations of such groups, and show:

\begin{theorem}
Closed oligomorphic groups are almost free. 
\end{theorem}

\smallskip

Another concept used in the study of automorphism groups
\index{group ! automorphism}%
 is relevant, and was researched by Truss.
\index{Truss, J. K.}%
 If $M$ is a countable saturated structure then $g \in \Aut(M)$ is \emph{generic}
\index{automorphism ! generic}%
if the conjugacy class of $g$ is comeager in $\Aut(M)$, in the pointwise convergence topology.
\index{topology ! of pointwise convergence}%
  The random graph
\index{graph ! random}%
 as well as many other structures admit such automorphisms~\cite{truss3}~\cite{hrush}, and this has been useful in proving the small index property
\index{small index property (SIP)}%
for them, see for example~\cite{hodges}.  Generic sequences of automorphisms are finite sequences $(g_1, \ldots, g_n)$ for which the set $\{(\alpha g_1 \alpha^{-1}, \ldots, \alpha g_n \alpha^{-1}) : \alpha \in \Aut(M)\}$ is comeager on a nonempty open subset of the product of $n$ copies of $\Aut(M)$.  For such sequences to exist for all $n$, firstly automorphisms must amalgamate, and secondly if $f_1, \ldots, f_n$ are finite partial elementary maps
\index{elementary map}%
 from finite subsets of $M$ onto finite subsets of $M$, then there exists a finite $A \subset M$ containing all the domains and ranges of the $f_i$'s and elementary permutations $g_1, \ldots, g_n$ of $A$ which extend the $f_i$'s.  

The generic elements of $\Aut(\mathfrak{R}_{m,\omega})$ have cycle type $\prod_{1 \le k < \infty} k^{\infty}$, reflecting that no finite amount of information about $g$ can prevent $x$ from lying in a finite-length cycle, and so $g$ has no infinite cycles.  For a structure $M$, the necessary condition for generic elements to exist in $G = \Aut(M)$ is that an appropriate cofinal subset of the natural family of finite approximations to $g \in G$ satisfies the relevant amalgamation and joint embedding properties.
\index{joint embedding property}%
For $\mathfrak{R}_{m,\omega}$ this cofinal subset may be taken to be the set of $p$ for which $\dom(p) = \range(p)$.  Macpherson
\index{Macpherson, H. D.}%
 gave the following example.  Let $g \in \Aut(\mathfrak{R}_{m,\omega})$, $\alpha \in \Sym(m)$ and function $F$ be such that $F \{xg, yg \} = (F\{x, y\})\alpha$ for all $x \neq y$ in $\mathfrak{R}_{m,\omega}$.  If $m$ is finite then each conjugacy class of the full symmetry group $\Sym(m)$ has a corresponding generic automorphism, whilst if $m$ is infinite there will be a generic automorphism.  Generic automorphisms of the random (countable, universal, homegeneous) tournament
\index{tournament ! random}%
are also known.  The group $\Aut(\mathbb{Q}, <)$ however is torsion-free
\index{group ! torsion -free}%
 and therefore no element of this group can have finite cycles of length greater than 1, and so in this case  the generic elements look different.  A condition for genericity
\index{generic structure}%
\index{structure ! generic}%
  of automorphisms of the $m$-coloured version of the rationals in $\Aut(\mathbb{Q}, <)$ is known, and refinements of the concept of genericity were examined by Truss.  An example of a non-closed permutation group
\index{group ! permutation ! non-closed}%
   having generic members is the point-stabilizer $\Aut(\mathbb{Q})_{x}$.
  
A notion of genericity 
\index{generic structure}%
\index{structure ! generic}%
 for endomorphisms of homogeneous structures has been defined and studied~\cite{lockett}.
 
The notion of Baire category works better than measure here, for while there are a countable set of finite approximations, it is unclear how measures should be applied.

In~\cite{kechrisros}, Kechris
\index{Kechris, A. S.}%
and Rosendal
\index{Rosendal, C.}%
extend the concept of generics that Hodges et al.~\cite{hodges} employ for automorphism groups
\index{group ! automorphism}%
 to more general Polish groups.  A Polish group has \emph{ample generic}
\index{ample generic element}%
elements if for all finite $n$ there exists a comeagre orbit for the (diagonal) conjugacy action of $G$ on $G^n : g \cdot (g_1, \ldots, g_n) = (g g_1 g^{-1}, \ldots, g g_n g^{-1})$.  This is clearly a stronger condition than just having a comeagre conjugacy class; for example $\Aut(\mathbb{Q}, <)$ has the latter but not the former.  The \emph{cofinality}
\index{group ! cofinality}%
of a group $G$ is the least cardinality of a well-ordered chain of proper subgroups whose union is $G$.  Again generalizing results of~\cite{hodges} they show
\begin{theorem}
A Polish group with ample generic elements is not the union of countably many non-open subgroups (or even cosets of subgroups).
\end{theorem}

If $G$, a closed subgroup of $\Sym(\omega)$, is oligomorphic, then by a result of Cameron, any open subgroup of $G$ is contained in only finitely many subgroups of $G$, so if $G$ has ample generics then it has uncountable cofinality.

That $\Sym(\omega)$ has ample generics has continuity consequences.  For example any unitary representation of $\Sym(\omega)$ on separable Hilbert space acts by homeomorphisms on some locally compact Polish space or by isometries on some Polish metric space, then it does so continuously.

\begin{theorem}~\cite{kechrisros}
$\Sym(\omega)$ has precisely two separable group topologies, namely the trivial one and the usual Polish topology.
\end{theorem}

The following theorem~\cite{macth} connects the existence of a comeagre conjugacy class in a Polish group and actions of the group on trees.
\begin{theorem}[Macpherson--Thomas]
\index{Thomas, S.}%
\index{Macpherson, H. D.}%
Let $G$ be a Polish group with a comeagre conjugacy class.  Then $G$ cannot be written as a free product
\index{group ! free product}%
 with amalgamation.
\end{theorem}
It is also true that~\cite{kechrisros} every element of such a group $G$ is a commutator and that it does not have $\mathbb{Z}$ as a homomorphic image.

We can also ask what group is generated by a `typical' $n$-tuple of elements. 
Dixon
\index{Dixon, J. D.}%
 proved that almost all pairs of elements of the finite symmetric group $Sym(n)$ generate $Sym(n)$ or the alternating group $Alt(n)$, later~\cite{dix} proving an analogue for the symmetric group of countable degree: almost all pairs of elements (in the sense of Baire category, that is, a residual set) generate a highly transitive free subgroup.  (The existence of highly transitive free groups was first shown by McDonough).
\index{McDonough, T. P.}%
These results have been extended to wider classes of groups. 

\medskip

Our final topic is that of large subgroups of infinite symmetric groups~\cite{macpherson}~\cite{cam10},
\index{group ! symmetric ! infinite}%
 of which subgroups of small index is a part, and in particular proper \emph{maximal subgroups}.
\index{maximal subgroup}%
A subgroup of a group $G$ is maximal in $G$ if and only if the action of $G$ by right multiplication on the right cosets of the subgroup is primitive.  
\index{group ! permutation ! primitive}%
  So maximal subgroups and primitive group actions are essentially the same thing.  The maximal subgroups of the symmetric group were classified in~\cite{liebprsa}.  We can only indicate a few of the highlights in this subject.  The operand $|\Omega|$ is not restricted to being countable.  

If $H$ is a maximal subgroup of $\Sym(\Omega)$ and $\aleph_{0} \le \alpha \le |\Omega|$ then either $\BSym_{\alpha}(\Omega) \le H$, in which case $H$ is highly transitive; or  $\BSym_{\alpha}(\Omega) H = \Sym(\Omega)$, in which case for some $\Delta \subset \Omega$ with $|\Delta| < \alpha$, $\Sym_{\{\Delta\}}$ induces $\Sym(\Omega \backslash \Delta)$~\cite{macneu}.

A \emph{filter} $\mathcal{F}$
\index{filter}%
on a set $\Omega$ is a subset of $\mathcal{P}(\Omega)$ containing $\Omega \backslash \emptyset$, which is closed upwards and under finite intersections. By replacing every set in the filter with its complement, we get an \emph{ideal} 
\index{ideal}%
on $\Omega$, which is a subset of $\mathcal{P}(\Omega)$ containing $\emptyset \backslash \Omega$, and is closed downwards and under finite unions.  An \emph{ultrafilter}
\index{ultrafilter}%
is a filter which is maximal with respect to inclusion, and a \emph{principal ultrafilter}
\index{ultrafilter ! principal}%
is one that has a least element.  Two ultrafilters are \emph{equivalent}
\index{ultrafilter ! equivalent}%
if a permutation of $\Omega$ maps one to the other.  Two ultrafilters are \emph{uniform}
\index{ultrafilter ! uniform}%
if all its elements have the same cardinality.  If $S \subseteq \mathcal{P}(\Omega)$ and $G \le \Sym(\Omega)$ then we can define $G_{S} := \{ g \in G : \forall s \in S \Leftrightarrow s g \in S \}$.  Also $\Sym(\Omega)_{(\mathcal{F})} : = \{ g \in \Sym(\Omega) : \fix(g) \in \mathcal{F} \}$, from which it is immediate that $\Sym(\Omega)_{(\mathcal{F})} \unlhd \Sym(\Omega)_{\{\mathcal{F}\}}$.  The next theorem is in~\cite{dixnt},
\begin{theorem}[Dixon, Neumann and Thomas]
\index{Dixon, J. D.}%
\index{Neumann, $\Pi$. M.}%
\index{Thomas, S.}%
Suppose $G \le \Sym(\Omega)$ and $|\Sym(\Omega) : G| < 2^{\kappa}$.  Let
\[ \mathcal{F} := \{\Gamma \subseteq \Omega : \exists \Delta \subseteq \Gamma\ \text{with}\ |\Omega \backslash \Delta| = \kappa\ \text{and}\ \Sym_{(\Delta)} \le G \}. \]
Then
\begin{itemize}
\item[(1)]  $\mathcal{F}$ is a filter on $\Omega$ and contains a moiety of $\Omega$;
\item[(2)]  $\Sym_{(\mathcal{F})} \le G \le \Sym_{\{\mathcal{F}\}}$.
\end{itemize}
\end{theorem}
A structural similarity between $\Sym(\Omega)_{\{\mathcal{F}\}}$ and $\Sym(\Omega)$ is suggested by the next result~\cite{macpherson},
\begin{theorem}
\label{bthm}
Let $\mathcal{F}$ be a uniform ultrafilter $\mathcal{F}$ on $\Omega$,and $|\Omega| = \kappa$.  Then $\BSym_{\kappa}(\Omega) \unlhd \Sym_{\{\mathcal{F}\}}$, and $\Sym_{\{\mathcal{F}\}} / \BSym_{\kappa}(\Omega)$ is simple.
\end{theorem}

Richman and Ball provided the first results on maximal subgroups of infinite symmetric groups.
\index{group ! symmetric ! infinite}%
 Among the results proved by Richman~\cite{richman}
\index{Richman, F.}%
are
\begin{itemize}
\item[(a)]  For an ultrafilter $\mathcal{F}$ on $\Omega$, $\Sym_{(\mathcal{F})} = \Sym_{\{\mathcal{F}\}}$ and is maximal in $\Sym(\Omega)$;
\item[(b)]  $\Sym(\Omega)$ is highly transitive on the right cosets of $\Sym_{\{\mathcal{F}\}}$;
\item[(c)]  the setwise stabilizer of a finite set of equivalent ultrafilters is maximal in $\Sym(\Omega)$ (including simply the setwise stabilizer of a finite set);
\item[(d)]  the almost stabilizer of a partition into a finite number of equi-cardinal sets is maximal in $\Sym(\Omega)$ (a partition $\pi$ is \emph{almost stabilized}
\index{almost stabilizer}%
by $g$ if, for all $A \in \pi$, there exists $B \in \pi$ such that the symmetric difference of $Ag$ and $B$ is finite).
\end{itemize}

Ball showed~\cite{ball} that for an infinite set $\Omega$,  the setwise stabilizer of a finite set is a maximal subgroup in $\Sym(\Omega)$, and that any transitive maximal subgroup contains the finitary symmetric group.

If $H < \Sym(\Omega)$ is not highly transitive, then it is contained in a maximal subgroup~\cite{macprae}.  If $H < \Sym(\Omega)$ is such that $\langle H, B \rangle = \Sym(\Omega)$ for some $B$ with $|B| \le |\Omega|$, then $\exists g$ such that $\langle H, g \rangle = \Sym(\Omega)$~\cite{galvin}.

Another topic studied within this group is \emph{maximal closed groups}
\index{group ! maximal closed}%
that is groups $G < \Sym(\omega)$ such that $G$ is closed as a permutation group and $\Sym(\omega)$ has no proper closed subgroup properly containing $G$.  Known examples of these groups suggest that this concept is the right analogue of maximal subgroups of finite symmetric groups.  In order to classify these groups, it is natural to reduce them to primitive
\index{group ! permutation ! primitive}%
 versions, but whereas the O'Nan--Scott Theorem requires the existence of a socle, these $G$ may not have them.  What is adopted instead is the idea of a \emph{minimal closed normal subgroup};
\index{group ! minimal closed normal subgroup}%
 any automorphism group
\index{group ! automorphism}%
  of an $\aleph_0$-categorical structure $M$ possesses one, and they act on $M$ either regularly or oligomorphically or neither.

In a different direction, we mention a result dealing with closed groups on $\mathfrak{R}$.  In~\cite{macphwood}, Macpherson and Woodrow
\index{Macpherson, H. D.}%
\index{Woodrow, R. E.}%
consider closed permutation groups
\index{group ! permutation ! closed}%
 $G$ on a countably infinite set $X$ and study the permutation groups (setwise stabilizers) $(G\sb{\{A\}};A)$ induced on moieties $A \subset X$; here the assumption that $G$ is closed on $X$ can be characterized by the property that $G$ is the full automorphism group
\index{group ! automorphism}%
 of some first order structure on $X$.  The group $\Aut(\mathfrak{R})$ is shown to have interesting universality properties with respect to such subsets. 
 For instance, they showed that if $(H, B)$ is any closed permutation group of countable degree, then there is a moiety $A$ of the vertex set of $\mathfrak{R}$ such that $\Aut(\mathfrak{R})_{(A)} = 1$ and the permutation groups $(H, B)$ and $(\Aut(\mathfrak{R})_{\{A\}}, A)$.

\bigskip

Original sources of results on maximal subgroups of infinite symmetric groups include~\cite{brazil}~\cite{covingtonmac}~\cite{macneu}~\cite{macprae}~\cite{richman}
\index{group ! symmetric ! infinite}%

We refer the reader to~\cite{cam6}~\cite{cam4}~\cite{cam10}~\cite{kayem}~\cite{macpherson} and references therein for more on the general topic of topology in permutation groups and proofs of some of the assertions.

\section{Polynomials}
\label{Polynomials}

\head{Graded Algebras.}
\index{graded algebra}%

Let $G$ be an infinite permutation group acting on an infinite set
$\Omega$.  Cameron
\index{Cameron, P. J.}%
discovered~\cite{camaa}~\cite{cam2a} an algebra $\mathcal{A}$
\index{graph ! random ! polynomial algebra}%
\index{polynomial algebra}%
called the \emph{age algebra}
\index{age algebra}%
 which encodes information about the action of $G$ on finite subsets of
$\Omega$.   Let $K$ be a field of characteristic $0$, for example
$\mathbb{Q}$ or $\mathbb{C}$.  Define a $K$-vector space of functions
from $n$-element subsets ${\Omega \choose n}$ of $\Omega$ to $K$ which
are invariant under the natural action of $G$ on $n$-subsets.  Let $V_{n}$ be the $n$th homogeneous component, which is the space of functions from ${\Omega \choose n}$ to $K$, with pointwise addition.

Define $\mathcal{A} := \oplus_{n \ge 0}
V_{n}$, where $V_{n}$ is the space of $G$-fixed points, and multiplication for $f \in V_{m}, g \in V_{n}$ by 
\[ (f \cdot g ) (X) = \sum_{Y \subseteq X}  f(Y) \cdot g(X
  \setminus Y) \]
for $|Y| = m, |X| = n+m$ and this defines $fg \in V_{m+n}$.  If $\mathcal{A}^{G} := \oplus_{n \ge 0} V_{n}^{G}$ is the algebra of $\mathcal{A}$-fixed points then $\dim(V_{n}^{G}) = f_n$ if this number is finite.  

By extending this multiplication linearly to the whole domain we get an $\mathbb{N}$-graded, commutative, associative algebra $\mathcal{A}$, where $f$ and $g$ can be interpreted as the coefficients in the Hilbert series
\index{Hilbert series}%
 of the graded algebra.  It has as
identity the function in $V_{0}$ with value $1$, and as
non-zero-divisor
\index{zero-divisor}%
the function $e : V_{1} \to 1$.  Furthermore, all functions
with finite supports are nilpotent and so zero-divisors, because $f^{n}
= 0$ if the support of $f$ fails to have $n$ disjoint sets.  The age algebra
\index{age algebra}%
 of a relational structure
\index{relational structure}%
  is graded and connected, and so is finitely generated if and only if it is a Noetherian ring~\cite{atiyahmac}.
\index{Noetherian ring}%
 
The structure of $\mathcal{A}^{G}$ is known in certain cases:
\begin{itemize}
\item[(1)]  If $G$ is highly homogeneous
\index{group ! highly homogeneous}%
(e.g. $G = \Sym(\Omega)$), then $\mathcal{A}^{G}$ is the polynomial algebra
\index{polynomial algebra}%
 in one generator $e$, the constant function with value $1$ in $V_1$.
\item[(2)]  If $G = G_1 \times G_2$ in its intransitive action,
\index{group ! intransitive}%
then $\mathcal{A}^{G} = \mathcal{A}^{G_1} \otimes_{\mathbb{C}} \mathcal{A}^{G_2}$.  So for $G = \Sym(\Omega)^n$ in the intransitive action, $\mathcal{A}^{G}$ is a polynomial ring
\index{graph ! random ! polynomial ring}%
in $n$ generators of degree $1$.
\item[(3)]  If $G = \Aut(M)$, where $M$ is the Fra\"{\i}ss\'e limit
\index{Fra\"{\i}ss\'e limit}%
of the class $\mathcal{C}$ of all finite simple graphs, that is $M = \mathfrak{R}$, then the generators of the polynomial algebra
\index{polynomial algebra}%
 $\mathcal{A}^{G}$ correspond to connected graphs.
\index{graph ! connected}%
This example can be generalized~\cite{cam7} to other structures.
\end{itemize}

If $\mathcal{A}^{G}$ is a polynomial algebra
\index{polynomial algebra}%
 with $a_n$ homogeneous generators of degree $n$ for each $n$, then the sequences $(f_n)$ and $(a_n)$ are related by
\[ \sum_{n\ge0} f_n  x^{n} = \prod_{k \ge 1} (1 - x^k)^{- a_k}. \]
These sequences have combinatorial interpretations independently of whether or not it is known that $\mathcal{A}^{G}$ is a polynomial algebra.
\index{polynomial algebra}%
  For example, if $G$ is the group of switching automorphisms of the random graph $\mathfrak{R}$ (see main text), then $f_n$ is the number of $n$-vertex \emph{even graphs}
\index{graph ! even}%
(graphs with all vertex degrees even)
so $a_n$ is the number of $n$-vertex \emph{Eulerian graphs}
\index{graph ! Eulerian}%
(connected even graphs)~\cite{cama}~\cite{ms}.  Given that $\mathcal{A}^{G} \subseteq \mathcal{A}^{\Aut(\mathfrak{R})}$, and that $\mathcal{A}^{\Aut(\mathfrak{R})}$ is a polynomial ring, then $\mathcal{A}^{G}$ is an integral domain, but it is an open question as to whether or not it is a polynomial ring.

A more detailed account of the above theory can be found in~\cite{cam8}.  Julian Gilbey
\index{Gilbey, J. D.}%
has explored what happens when $G$ has no finite orbits~\cite{gilbey}.  He shows that if the permutation group in the Cameron's algebra has no finite orbits, then no homogeneous element of degree one is a zero-divisor of the algebra.  He conjectures that the algebra is an integral domain if, in addition, the group is oligomorphic, proving special cases in the form of wreath products and showing that the algebras corresponding to oligomorphic groups 
\index{group ! permutation ! oligomorphic}%
are polynomial algebras.
\index{polynomial algebra}%

It was conjectured that if $G$ has no finite orbits then $\mathcal{A}^{G}$ is an integral domain,
\index{graph ! random ! integral domain}%
this has recently been proved by M. Pouzet~\cite{pouzet1}.
\index{Pouzet, M.}%
using the following approach.  The \emph{kernel}
\index{relational structure ! kernel}%
 of a relational structure $\mathcal{M}$ is the subset $\ker(\mathcal{M})$ of $x \in \Omega$ $(|\Omega| = \infty)$ such that $\mathcal{A} (\mathcal{M}_{| \Omega \backslash \{x\}}) \neq \mathcal{A}(\mathcal{M})$.  If  $\ker(\mathcal{M}) \neq \emptyset$, pick $x \in \ker(\mathcal{M})$ and $S \in \Omega^{< \omega}$ such that $\mathcal{M}_{| S} \in \mathcal{A}(\mathcal{M}) \backslash \mathcal{A}(\mathcal{M}_{| \Omega \backslash \{x\}})$.  Let $T \in \Omega^{< \omega}$.  Set $f(T) := 1$ if $\mathcal{M}_{| T}$ is isomorphic to $\mathcal{M}_{| S}$, and $0$ otherwise.  Then $f^2 := f \cdot f = 0$. 

\begin{theorem}[Pouzet]
\label{pouzettheorem}
Let $\mathcal{M}$ be a relational structure
\index{relational structure}%
 with possibly infinitely many non-isomorphic types of $n$-element substructures. The age algebra
\index{age algebra}%
 $\mathbb{C} \mathcal{A} (\mathcal{M})$ is an integral domain if and only if $\ker(\mathcal{M}) = \emptyset$.
\end{theorem}

Thus in the $G$-action on $\Omega$ encoded by $\mathcal{M}$, $\ker(\mathcal{M})$ is the union of the finite $G$-orbits of the one-element sets.  So if $G$ has no finite orbit, $\ker(\mathcal{M}) = \emptyset$.  Therefore $\mathcal{A}^{G}$ is an integral domain.  If $\mathcal{A}^{G}$ is an integral domain then this property is preserved by passing to overgroups or transitive extensions.
\index{transitive extension}%

The kernel of a relational structure
\index{relational structure}%
 $\mathcal{M}$ is empty if and only if for every finite subset $S \subset \Omega$ there is a disjoint subset $T$ such that the restrictions $\mathcal{M}_{| S}$ and $\mathcal{M}_{| T}$ are isomorphic. The ages
\index{age}%
 of such relational structures have the \emph{disjoint embedding property}
\index{disjoint embedding property}%
meaning that two arbitrary members of the age can be embedded into a third in such a way that their domain are disjoint.  Ages with the disjoint embedding property are called \emph{inexhaustible}
\index{age ! inexhaustible}%
 and relational structures whose age is inexhaustible are called age-inexhaustible.
\index{relational structure ! age-inexhaustible}%
 Finally relational structures with finite kernel are \emph{almost age-inexhaustible}.
\index{relational structure ! almost age-inexhaustible}%
 If $G$ is oligomorphic then the kernel of the relational structure is finite.

\begin{lemma}~\cite{pouzeta}
Let $\mathcal{A}$ be an infinite inexhaustible age. Then: 

(i) For every age
\index{age}%
 $\mathcal{A'}$ included into $\mathcal{A}$ there is an infinite strictly increasing sequence of ages included into $\mathcal{A}$ such that
$\mathcal{A'} = \mathcal{A}_0 \subset \ldots \subset \mathcal{A}_n \subset \ldots$

(ii) The height of $\mathcal{A}$, is a limit ordinal, provided that it is defined. 
\end{lemma}

A relation is \emph{quasi-inexhaustible} if its kernel is finite.  Two relations with the same age have isomorphic kernels, so we say that the age $\mathcal{A}$ is \emph{quasi-inexhaustible}
 \index{age ! quasi-inexhaustible}%
if it is the age of a quasi-inexhaustible relation.  If the kernel is empty then quasi-inexhaustible reduces to inexhaustible.  R. Fra\"{\i}ss\'e
\index{Fra\"{\i}ss\'e, R.}%
 showed that ages are nonempty ideals of the poset of finite relational structures, considered up to isomorphism and ordered up by embeddability.  In~\cite{pouzetsobrani}, Pouzet and Sobrani
\index{Pouzet, M.}%
\index{Sobrani, M.}%
 study the poset of ages that lie in between two given ages, proving among other things that if this poset is infinite then it contains an infinite chain.

\bigskip

The concept of fixed points of an algebra leads on to the topic of polynomial invariants.

\bigskip

\head{Polynomial Invariants.}

We recall the basic idea of polynomial invariants
\index{polynomial invariant}%
of finite groups \cite{smithl}~\cite{benson}~\cite{stanl1}.  

A polynomial $f(x_1, \ldots, x_n)$ in the
polynomial algebra
\index{polynomial algebra}%
 $K[x_1, \ldots, x_n]$ in $n$ variables over field $K$, is called
\emph{symmetric}
\index{symmetric polynomial}%
if for $\pi \in Sym(n)$, $f(x_{\pi(1)}, \ldots, x_{\pi(n)}) = f(x_1,
\ldots, x_n)$.  Let $K[V]$ be a finite
dimensional $K$-vector space on which a finite group $G$ acts by linear
substitutions via the representation homomorphism $\rho : G \to
\Aut(V) = GL(V)$.  Equivalently $V$ is a finitely generated
$K[G]$-module, where $K[G]$ is a group ring.  So if $x_1, \ldots, x_n$
is a basis for the dual space $V^{*} = \Hom_{K}(V, K)$,~\label{Hom} then $K[V^{*}] =
K[x_1, \ldots, x_n] = K \oplus V^{*} \oplus S^{2}(V^{*}) \oplus
S^{3}(V^{*}) \oplus \ldots$, where $S^{m}(V^{*})$, the symmetric \emph{m}th power of
$V^{*}$, comprises degree \emph{m} homogeneous
polynomials in $x_1, \ldots, x_n$.  For example $S^{2}(V^{*})$ has a
basis consisting of the monomials $x_i x_j$ for the ${n+1 \choose 2}$
choices of $i$ and $j$.  The vector space dimension of $S^{m}(V^{*})$
is ${n+m-1 \choose m}$.  Regarding each $x_i$ as having degree one makes $K[V]$ a
graded ring.  If the $G$-action on $K[V]$ is defined by
$(gf)(v):=f(g^{-1}v)$, then the set of all fixed points of this action
is the \emph{ring of invariants}
\index{ring of invariants}%
 $K[V]^{G}$.  If $G =
\Sym(n)$, then $K[V]^{\Sym(n)} = K[e_1, \ldots, e_n]$, where the $e_i$
are the elementary symmetric polynomials in the $x_i$ $(1 \le i \le
n)$.   Also if $H$ is a finite permutation group then $\mathcal{A}^{\Sym(\mathbb{N}) \Wr H}$ is isomorphic to the ring of invariants of the linear group $G$, where the wreath product is considered in its imprimitive action.
\index{group ! permutation ! imprimitive}%
  In particular, if $G = \Sym(\mathbb{N}) \Wr \Sym(n)$, then $\mathcal{A}^{G}$ is a polynomial ring in generators of degree $1, 2, \ldots, n$.

The \emph{profile}~\cite{pouzetthiery}~\cite{pouzetthieryII}
\index{structure ! profile}%
of a relational structure $M$ 
\index{relational structure}%
is the function $\phi_{M}$ which counts for every integer $n$ the (possibly infinite) number of substructures of $M$ induced on non-isomorphic $n$-element subsets.  The profile depends only on the age
\index{age}%
 $\mathcal{A}$ of $M$, and we can associate many graded algebras with $M$ such that $\phi_{M}$ is the Hilbert function.
\index{Hilbert function}
One such graded algebra is Cameron's
\index{Cameron, P. J.}%
above-mentioned age algebra
\index{age algebra}%
 $K \mathcal{A}(M)$ over a field $K$.  The profile in this case is the dimension of the homogeneous component of degree $n$, and so $H_{\phi_{M}} := \sum^{\infty}_{n=0} \phi_{M}(n) Z^n$ is the Hilbert series
\index{Hilbert series}%
 of $K \mathcal{A}(M)$.  For a universal graph, the age is all finite graphs and the profile is $\phi_n \sim 2^{n(n-1)/2} /n!$.  For the random graph $\mathfrak{R}$,
\index{graph ! random}%
 the age is the free commutative algebra generated by the connected graphs and the Hilbert series is $H_{\mathfrak{R}} (Z) = \prod_{d \geq 1} \frac{1}{(1 - Z^d)^{c_n}}$, where $c_n$ is the number of $n$-vertex connected graphs.
 
For structures obeying certain conditions the \emph{invariant ring}
\index{invariant ring}%
of a finite permutation group was found by Pouzet and Thi\'ery~\cite{pouzetthiery}
\index{Pouzet, M.}%
\index{Thi\'ery, N. M.}%
to be isomorphic to some $K \mathcal{A}$.  If a group $G$ acts on $\{1, \ldots, k\}$, the Hilbert function of the subalgebra $K[X]^G$ of the polynomials in the algebra $K[X] = K[X_1, \ldots, X_k]$ which are invariant under $G$-action is called the \emph{orbital profile}.
\index{profile ! orbital}%
It counts the number of orbits of $n$-element subsets; oligomorphic groups
\index{group ! oligomorphic}%
have orbital profiles taking only finite values.   

Let $\mathcal{M} := (\Omega, (R_i)_{i \in I})$ be a relational structure,
\index{relational structure}%
 where the $R_i$ are $n_i$-ary relations on the domain $\Omega$.  To each size-$d$ subset $S \subset \Omega$, associate the monomial $x^{d(S)} := \prod_{i \in X} x_i^{d_i(S)}$, where $d_i(S) = | S \cap \Omega_i |$ for all $i \in X$.  To each orbit of sets, associate the unique maximal \emph{leading monomial} $\lm(S)$,
\index{leading monomial}%
 where $S$ ranges through the orbit.  In order to prove that the Hilbert generating series of a certain profile is a rational function of a certain form, Pouzet and Thi\'ery~\cite{pouzetthiery} endow the set of leading monomials with and \emph{ideal} structure in an appropriate polynomial ring.  In particular, if  $\lm_C$ is the set of leading monomials with chain $C = (\empty \subsetneq S_1 \subsetneq \ldots \subsetneq S_l \subsetneq X)$ of subsets as support, they realized $\lm_C$ as the linear basis of some \emph{ideal} of a polynomial ring, so that the generating series of $\lm_C$ is realized as a Hilbert series.
\index{Hilbert series}%
 There is a natural polynomial ring embedding $K[S_1, \ldots, S_l] \to K[X] : S_j \mapsto \prod_{i \in S_j} x_i$.  Infinite relational structures with a constant profile, equal to 1, were called \emph{monomorphic}.
\index{relational structure ! monomorphic}%
  Relations admitting a finite monomorphic decomposition have age algebras
 \index{age algebra}%
  that are graded subalgebras of finitely generated polynomial algebras.
\index{polynomial algebra}%
    When all monomorphic components are infinite, the monomial ideal spanned by monomials $m = S_1^{r_1} \ldots S_l^{r_l}$ such that $d_i(m) > |\Omega_i|$ for some $i$, is the trivial ideal $\{0\}$.

Much is known about the growth of a profile~\cite{cam6} but less is known about `local' conditions relating 
individual values of $\phi_n$, however it can be proved either using Ramsey-type theorem or using finite combinatorics
\index{combinatorics}%
 and linear algebra that $\phi_n \le \phi_{n+1}$.  So the profile of an infinite relational structure $\mathcal{M}$ is nondecreasing.  Furthermore, if the relations constituting $\mathcal{M}$ have bounded arity or the kernel of $\mathcal{M}$ is finite, its growth rate is either polynomial or faster than every polynomial~\cite{pouzeta}.  If $\mathcal{M}$ is an infinite relational structure
\index{relational structure}%
 with a finite monomorphic decomposition, then its generating Hilbert series is a rational fraction and its orbital profile is bounded by some polynomial, that is $\phi_{\mathcal{M}} \sim a n^{k-1}$~\cite{pouzetthiery}.

A conjecture of Cameron that remains open is that the orbital profile of an $\aleph_0$-categorical
\index{aleph@$\aleph_0$-categorical}%
structure is polynomial provided it is bounded by a polynomial.

Pouzet's Theorem mentioned above has the following consequence: if $G$ has no finite orbits then $\phi_{m+n} \ge \phi_{m} + \phi_{n-1}$.  This is because multiplication induces a map from the Segre variety 
\index{Segre variety}%
(the rank 1 tensors modulo scalars) in $V_m \otimes V_n$ into $V_{m+n}$ modulo scalars; so the dimension of $V_{m+n}$ is at least as great as that of the Segre variety~\cite{cam8}.

If $G$ has no finite orbits then $e$ is prime in $\mathcal{A}^{G}$, that is $\mathcal{A}^{G} / \langle e \rangle$ is an integral domain.  This implies that the profile of the relational structure satisfies $\psi_{m+n} \ge \psi_{m} + \psi_{n-1}$, where $\psi_{n} = \phi_{n+1} - \phi_{n}$.  (This follows by an application of the Segre variety argument to $\mathcal{A} / \langle e \rangle$, whose $n$th homogeneous component is $V_{n+1} / e V_n$, with dimension $\phi_{n+1} - \phi_{n}$).

Our final result on profiles is the following result which appears in~\cite{pouzet1}, to which we refer for details.

\begin{proposition}
\label{pouprop}
Let $\mathcal{A}$ be the age
\index{age}%
 of a relational structure
\index{relational structure}%
  $\mathcal{M}$ such that the profile
of $\mathcal{M}$ takes only finite values and $K \mathcal{A}$ be its age algebra.
\index{age algebra}%
 If $\mathcal{A'}$ is an initial segment of $\mathcal{A}$ then: 

(i) The vector subspace  $\mathscr{J} := K(\mathcal{A} \backslash \mathcal{A'})$ spanned by $\mathcal{A} \backslash \mathcal{A'}$ is an ideal of $K \mathcal{A}$.
 
Moreover, the quotient of $K \mathcal{A}$ by $\mathscr{J}$ is a ring isomorphic to the ring $K \mathcal{A'}$. 

(ii) If this ideal is irreducible then $\mathcal{A'}$ is a subage of $\mathcal{A}$. 

(iii) This is a prime ideal if and only if $\mathcal{A'}$ is an inexhaustible age.
\end{proposition}

\smallskip

We end with an application~\cite[p.95]{neusel} of invariant theory is towards counting \emph{weighted graphs},
\index{graph ! weighted}%
that is graphs which have a weight $m_{ij} \in \mathbb{C}$ attached to every edge.  A \emph{weighted isomorphism} between two weighted graphs is a bijection $\phi$ mapping edges onto edges of the same weight.  Counting invariants of the same representation type gives the number of isomorphism classes of $n$-vertex weighted graphs.  Briefly, $\Sym(n)$ acts on the polynomial ring $\mathbb{C}[x_{ij} | 1 \le i < j \le n]$ in $n \choose 2$ variables by acting on the indices,
\begin{displaymath}
s(x_{ij}) = \left\{ \begin{array}{ll}
x_{s(i)s(j)} & \text{if $s(i) < s(j),$} \\
x_{s(j)s(i)} & \text{otherwise.}
 \end{array} \right.
\end{displaymath}
If the basis elements $x_{ij}$ are arranged in a symmetric matrix
\begin{displaymath}
X_{ij} = \left\{ \begin{array}{ll}
x_{ij} & \text{for $i < j,$} \\
x_{ji} & \text{for $i > j,$} \\
0 & \text{otherwise}
 \end{array} \right.
\end{displaymath}
with zeroes on the diagonal, then $\Sym(n)$ acts by conjugation, $sX_{ij}s^{-1}$, where we have identified $s \in \Sym(n)$ with its image under the defining representation.  If $\underline{m} = (m_{ij})$ then two weighted graphs $\Gamma, \Gamma'$ are isomorphic if and only if $p(\underline{m}) = p(\underline{m'})$ for all invariant polynomials $p \in \mathbb{C}[x_{ij} | 1 \le i < j \le n]^{\Sym(n)}$, where $s \in \Sym(n)$ such that $\underline{m} = s(\underline{m'})$ gives the graph isomorphism.

\bigskip

\section{Enumeration, Growth Rates and Reconstruction}
\label{EnumerationandReconstruction}

We will give the briefest of introductions to these three topics.

\emph{Enumeration}

Let $M$ be a countable homogeneous structure with age
\index{age}%
 $\mathcal{A}$ and let $G = \Aut(M)$.  By homogeneity, two finite subsets of $M$ lie in the same $G$-orbit if and only if the induced substructures are isomorphic.  So the sequence enumerating \emph{unlabelled} $n$-element members of $\mathcal{A}$ (that is, up to isomorphism) is identical with the sequence enumerating the $G$-orbits on unordered $n$-element subsets of $M$.  Similarly, the number of  \emph{labelled} $n$-element members of $\mathcal{A}$ (that is, members of $\mathcal{A}$ on the set $\{1, 2, \ldots, n\}$) is equal to the number of $G$-orbits on ordered $n$-tuples of distinct elements $M$. 

Let $G$ be an oligomorphic permutation group acting on each of the sets $\Omega^n$ (all $n$-tuples of elements of $\Omega$), $\Omega^{(n)}$ (all $n$-tuples of distinct elements of $\Omega$), and $\Omega^{\{n\}}$ (all $n$-element subsets of $\Omega$).  If $F_{n}^{*}, F_{n}$ and $f_n$ denote these numbers of orbits then $f_n \leq F_n \leq n! f_n$, since each orbit on $n$-sets corresponds to at least one and most $n!$ orbits on $n$-tuples.  A permutation group is \emph{$n$-transitive}
\index{group ! permutation ! $k$-transitive}%
 if $F_n = 1$,  and \emph{$n$-homogeneous}
\index{group ! permutation ! $k$-homogeneous}%
if $f_n = 1$.

Also,
\[ F_{n}^{*} = \sum^{n}_{k=1} S(n, k) F_{k}, \]
where $S(n, k)$ is the Stirling number of the second kind,
\index{Stirling number ! of second kind}%
 the number of partitions of an $n$-set into $k$ parts~\cite{camtay}.  For an orbit $(\alpha_1, \ldots, \alpha_n)^G$ on $n$-tuples determines, and is determined by, a partition of $\{1, \ldots, n\}$ into $k$ parts (where $i$ and $j$ lie in the same part if $\alpha_i = \alpha_j$) and an orbit on $k$-tuples of distinct elements.  The \emph{exponential generating function} 
\index{exponential generating function}%
is given by $F(t) = \sum \frac{F_n t^n}{n!}$.  The series $F(t)$ for a direct product (in the intransitive action) or a wreath product (in the imprimitive action)
\index{group ! permutation ! imprimitive}%
 can be calculated from those of the factors:
\[ F_{G \times H} (t) =  F_{G} (t) \times F_{H} (t), \]
\[ F_{G \Wr H} (t) =  F_{H} (t) (F_{G} (t) - 1). \]

If $S = \Sym(\omega)$ and the sequence $x$ is realized by a group $G$, then $S x$ is realized by $G \Wr S$.

Let $G$ act transitively on $\Omega$, and $G_{\alpha}$ be the stabilizer of the point $\alpha \in \Omega$.  It can be shown that 
\[ F_{G_{\alpha}} (t) = {\operatorname{d}\over\operatorname{d}t} F_G (t), \]
or equivalently differentiating an exponential generating function corresponds to shifting the sequence terms one place to the left,
\[ F_{n} (G_{\alpha}) = F_{n+1} (G). \]
So there is an equivalence between $G_{\alpha}$-orbits on $n$-tuples and of $G$ orbits on $(n+1)$-tuples.  If $G$ is intransitive, then the derivative of $F_G (t)$ gives the sum of $F_{G_{\alpha}} (t)$, over a set of representatives $\alpha$ of the orbits of $G$.

For any oligomorphic group $G$, \emph{generalized Stirling numbers}
\index{Stirling number ! generalised}%
denoted by $S[G](n, k)$, can be defined~\cite{cam6a} that obey
\[ \sum^{n}_{k=1} S[G](n, k)  =  F_{n}(G \Wr S),\]
and have the composition property
\[ \sum^{n}_{l=k} S[G](n, l) S[H](l, k)  =  S[G  \Wr H](n, k).\]
The last result can be expressed in terms of infinite lower triangular matrices of generalized Stirling numbers.  There is also a statement in terms of Joyal's theory of species~\cite{bergeron}~\cite{joyal}.
\index{species}%
\index{Joyal, \'A.}%

There is a linear analogue of the combinatorics
\index{combinatorics}%
 of sets and functions that applies to vector spaces over finite fields and linear transformations in which the Gaussian (or $q$-binomial coefficient) replaces the binomial coefficient~\cite{camtay}.  These coefficients enumerate the number of $k$-dimensional subspaces of $n$-dimensional $GF(q)$-vector space~\cite{kacc}.

\bigskip

\emph{Growth Rates}

There are general theorems on growth rates of counting sequences, many of them requiring a primitive permutation group
\index{group ! permutation ! primitive}%
 $G = \Aut(M)$, for example if $M = \mathfrak{R}$.  Though primitivity is not a first-order concept, it can be completely described if $G$ is oligomorphic.

If is known that $f_{n+1} \geq f_n$ for all $n$, that is, the sequence $(f_n)$ is non-decreasing, as proved in~\cite{cam6}.  Macpherson
\index{Macpherson, H. D.}%
 showed there to be a gap in growth rates between constant and exponential:
 
\begin{theorem}
If $G$ is an oligomorphic permutation group which is primitive but not highly set-transitive then there exists a constant $c$ such that $f_n > c^n$ for all sufficiently large $n$.
\end{theorem}

Macpherson proves the theorem holds with $c = 2^{\frac{1}{5}} - \epsilon$.  There is a variant of this theorem that does not require primitivity due to Pouzet~\cite{pouzeta1} and Macpherson~\cite{macpha}

\begin{theorem}
Let $G$ be an oligomorphic permutation group.

(a)  $c_1 n^d \leq f_n \leq c_2 n^d$ for all $n$, where $d$ is a non-negative integer and $c_1, c_2$ positive constants; or

(b) $f_n \geq exp(n^{\frac{1}{2} - \epsilon})$ for all sufficiently large $n$.
\end{theorem}

There are known examples in the growth rate spectrum lying between $exp(n^{\frac{p}{p+1} - \epsilon})$ and $exp(n^{\frac{p}{p+1} + \epsilon})$, as well as faster than any fractional exponential $exp(n^{c})$ for any $c < 1$ but slower than straight exponential $c^n$ for $c > 1$.  Macpherson has studied~\cite{macph1a} permutation groups of rapid growth, that is faster than exponential, and established connections with model-theoretic properties related to stability.

Macpherson's bound was improved by Merola~\cite{merola},
\index{Merola, F.}%
 who shows that if $G$ is an infinite primitive permutation group
\index{group ! permutation ! primitive}%
  which is not highly homogeneous (so for some $n$ has more than one orbit on $n$-sets) then the number of orbits on $n$-tuples of distinct elements is, for large enough $n$, bounded below by $\frac{c^{n} n!}{p_G(n)}$, where $p_G(n)$ is a $G$-dependent polynomial, and the result is proved with the constant $c \approxeq 1.172$; she proves that 2 is an upper bound for $c$.  The proof is an induction on the degree of transitivity, considering the cases when $G$ acts on a graph (when $G$ is not $2$homogeneous), or a tournament (when $G$ is $2$-homogeneous but not 2-transitive), or a Steiner system (one of the possibilities in the $2$-transitive not $2$-primitive case).  Merola finds lower bounds for the numbers of labelled size $n$ substructures of these structures.

Therefore counting orbits of an oligomorphic permutation group on ordered or unordered $n$-sets is the same as enumerating labelled or unlabelled structures in a Fra\"{\i}ss\'e class~\cite{cam6a}.
\index{Fra\"{\i}ss\'e class}%

More recently, the following has been proved~\cite{pouzetthiery}
\index{Pouzet, M.}%
\index{Thi\'ery, N. M.}%

\begin{theorem}[Pouzet, Thi\'ery]
Let $\mathcal{M} : = (X, (\rho_i)_{i \in I})$ be a relational structure
\index{relational structure}%
 on an infinite set $X$.  Let $\phi_{\mathcal{M}}$ be a non-decreasing profile.  If either the signature $\mu := (n_i)_{i \in I}$ is bounded or the kernel $\ker(\mathcal{M})$ of $\mathcal{M}$ is finite, then the growth of $\phi_{\mathcal{M}}$ is either polynomial or as fast as every polynomial
\end{theorem}

Conditions equivalent to the profile of $\mathcal{M}$ being bounded are given in~\cite{pouzetthiery}.

\bigskip

\emph{Reconstruction}
\index{reconstruction}%
\index{growth rates}%

The Vertex Reconstruction Conjecture of Kelly~\cite{kelly} 
\index{Kelly, P. J.}%
is one of the most famous unsolved problems in graph theory, and is a statement about graphs being uniquely determined by their subgraphs.  A \emph{vertex-deleted subgraph}
\index{graph ! vertex-deleted}%
is an induced subgraph formed by deleting one vertex from a graph $\Gamma$.  If $D(\Gamma)$ is the \emph{deck}, 
\index{graph ! deck}%
 that is the set of all vertex-deleted subgraphs, then call two graphs having the same deck \emph{hypomorphic}.
\index{graph ! hypomorphic}%
The conjecture states that 

\begin{conjecture}
Any two hypomorphic graphs on at least three vertices are isomorphic.
\end{conjecture}

For an algebraic formulation of the conjecture, take $G = \Aut(\mathfrak{R})$, so that a basis for the $G$-invariant functions $V_{n}^{G}$ (where $G$ acts on $V_{n}$ naturally: $f^g (X) = f(X g^{-1})$), consists of the characteristic functions $\chi_{\Gamma}$ of all unlabelled $n$-vertex graphs $\Gamma$.  Since $e \in V_{1}^{G}$ is not a zero-divisor, multiplication by $e$ is a monomorphism from $V_{n}^{G} \to V_{n+1}^{G}$ for all $n$.  Consider the map $\theta = \theta_n : V_{n}^{G} \to V_{n-1}^{G}$ which is dual to multiplication by $e$:
\[ \theta(\chi_{\Gamma}) = \sum_{x \in \Gamma} \chi_{\Gamma \backslash \{x\}}. \]
The map $\theta$ takes the graph $\Gamma$ to the sum of the graphs in its deck of vertex-deleted subgraphs.  Since $e$ is one-to-one, $\theta$ is onto.  More generally, if $G$ is any oligomorphic group then $G$-orbits on $n$-sets replace the characteristic functions.

Variants of this conjecture for graphs and other classes of structures are discussed in~\cite{cam6a}~\cite{cam6b} and~\cite{laurisc}.  Also noteworthy is the reformulation of the edge-reconstruction conjecture by V. Mnukhin~\cite{mnukhin1}~\cite{mnukhin2}, in terms of orbit reconstructions of an arbitrary finite permutation group, based on embedding the set of orbits in an \emph{orbit algebra}.
\index{Mnukhin, V.}%
 Other approaches are listed in~\cite[p.~67]{cameron}, amongst which we mention that of Joyal's 
\index{Joyal, \'A.}%
theory of combinatorial species~\cite{bergeron}~\cite{joyal}.

\section{Cartan Triality and Eight-Dimensional Exotica}
\label{carsec}
\index{Cartan, \'E.}%
The triality property of $\mathfrak{R^{t}}$ has been sufficiently significant in directing one of our lines of research as to warrant a section outlining the basic motivating finite theory whose infinite version we used to demonstrate the connection.

This section is a brief introduction to the exceptional $\Sym(3)$ outer automorphisms
\index{automorphism ! outer}%
that arises for the groups $\PSO(8)$
\index{group ! orthogonal ! projective special}%
 and $\Spin(8)$
 \index{group ! spin}%
 \index{group ! Spin(8)@$\Spin(8)$}%
(but not for $\SO(8)$,
\index{group ! SO(8)@$\SO(8)$}%
which is the light cone restriction of the Lorentz group
$\SO(9,1)$).  Though discovered by \'E. Cartan~\cite{cartan}, in fact the oldest phenomenon associated with
triality predates Cartan and is connected to complex quadrics~\cite{adam}~\cite{cam20}; see~\cite{study}.  We will follow the usual textbook account where it will transpire that the Clifford algebra
\index{Clifford algebra}%
for a given real or complex vector space has an associated quadratic form which contains in its
multiplicative group a subgroup which is a double cover of the orthogonal group of automorphisms of the space preserving the form.  We summarize the theory as expounded in~\cite{chevalley}.

Within the main text, we shall prove that the action of the colour group
$\Sym(\mathfrak{r}, \mathfrak{b}, \mathfrak{g})$ on the automorphism group
\index{group ! automorphism}%
 $\Aut(\mathfrak{R^{t}})$ of the triality graph,
\index{graph ! triality}%
 in a sense originates in a generalization of the \emph{Cartan triality} group action.
\index{Cartan triality}%
\index{group ! action}%

Vectors are preserved under a rotation of $2\pi$.  Spinors are objects that require a $4\pi$ rotation to return to their original value; this explains their occurrence in describing half-integer spin elementary particles, fermions.  That the angular momenta of the emitter has half-integral units (as multiples of Planck's constant) of the quantized spin can be discerned by the emitted photons.  This appears in Dirac's spinor solutions of his equations for the electron. It also indicates why we may expect spaces of both such objects to arise in representations of the rotation group.  In 8 dimensions, we have an extra symmetry that is unique to spaces of this dimensionality.

Let $V$ be an $8$-dimensional vector space over $\mathbb{R}$ for which we can take
without loss of generality the quadratic form to be $q= x_1 x_2 + x_3 x_4 +
x_5 x_6 + x_7 x_8$ with two maximal totally singular subspaces being
$\langle e_1, e_2, e_3, e_4 \rangle$ and $\langle e_5, e_6, e_7, e_8
\rangle$ in terms of standard basis vectors.  This leads to a \emph{Clifford
algebra}
\index{Clifford algebra}%
denoted by $CA = CA(r_1, \ldots, r_n)$, $r_i \in \mathbb{R}$ (field of char $\neq
2$), which is the free $\mathbb{R}$-algebra $\mathbb{R}\{x_1, \ldots, x_n\}$ over indeterminates
${x_1, \ldots, x_n}$ satisfying relations $x^2 = q(x) . 1$.

For example over the field $\mathbb{R}$, $CA(-1) = \mathbb{C}$.  The two-sided ideal of the
relevant tensor algebra (to be defined below)
$I(q) \subset T(V)$ is generated by even degree elements.  This
induces a $\mathbb{Z} / 2 \mathbb{Z}$ grading on the Clifford algebra
$CA = CA^+ \oplus CA^-$ with $CA^+ \cdot CA^+ \subset CA^+$, $CA^- \cdot CA^-
\subset CA^+$, $CA^+ \cdot CA^- \subset CA^-$, $CA^- \cdot CA^+ \subset
CA^-$.  The dimension $2^{n-1}$ subalgebra $CA^+$ is spanned by
products of an even number of elements in $V$.

It transpires~\cite{fultha} that the Lie algebra
\index{Lie algebra}%
$\mathfrak{so}(q) = \mathfrak{so}(n, \mathbb{R})$ embeds inside the Lie
algebra of the even part of the Clifford algebra and that the Clifford
algebra can be identified with one or two copies of matrix algebras.  Now let $S$ be the space of spinors or spin representations
\index{spinor (spin representation)}%
 of $CA$ with dim $S = 16$ and $S = S^{+} \oplus S^{-}$.
Then on each of the half-spin representations of
$\mathfrak{so}(n, \mathbb{R})$, written $S^{+}$, $S^{-}$ there is a quadratic form isomorphic to the original form on
$V$.  Algebraic triality results from the order $3$ mapping $J: V \to S^{+} \to S^{-} \to V$.

For any $2n$-dimensional vector space $V$, the Clifford algebra
has dimension $2^{2n} \cong M_{2^n \times 2^n} (\mathbb{R})$, a space of matrices.  Then $S = S^{+}
\oplus S^{-}$ gives dim $S^{+}\ =$ dim $S^{-} = 2^{n-1}$.  Triality
occurs only when $2n' = 2^{n'-1}$, that is only when $n' = 4$ that is $n = 8$.

The following is a flow-chart that gives the key ingredients of the derivation of this algebraic triality together with an explanation of the terms:
\begin{figure}[!h]$$\xymatrix{
{V\ :\ \mbox{8 dimensional vector space}} \ar[d]\\
{T\ :\ \mbox{tensor algebra}} \ar[d]\\
{CA\ (= T/I)\ :\ \mbox{Clifford algebra of } q \mbox{ with } |CA|=2^8} \ar[d]^{\rho_{s}\
\mbox{(spinor representation)}}\\
{S\ :\ \mbox{spinor space of the } CA \mbox{ representation with } |S|= 2^4} \ar[d]\\
{CG\ :\ \mbox{Clifford group of } q} \ar[d]\\
{CG^{+}_{0}\ = CG_{0} \cap CG^{+}\ :\ \mbox{Reduced Clifford group }} \ar[d]^{\chi}\\
{\chi(CG^{+}_{0}) = \Omega^{+}(8, \mathbb{R})\ :\ \mbox{Reduced orthogonal group}}
}$$\caption{Pictorial synopsis of 8--dimensional Cartan triality}  
\end{figure}
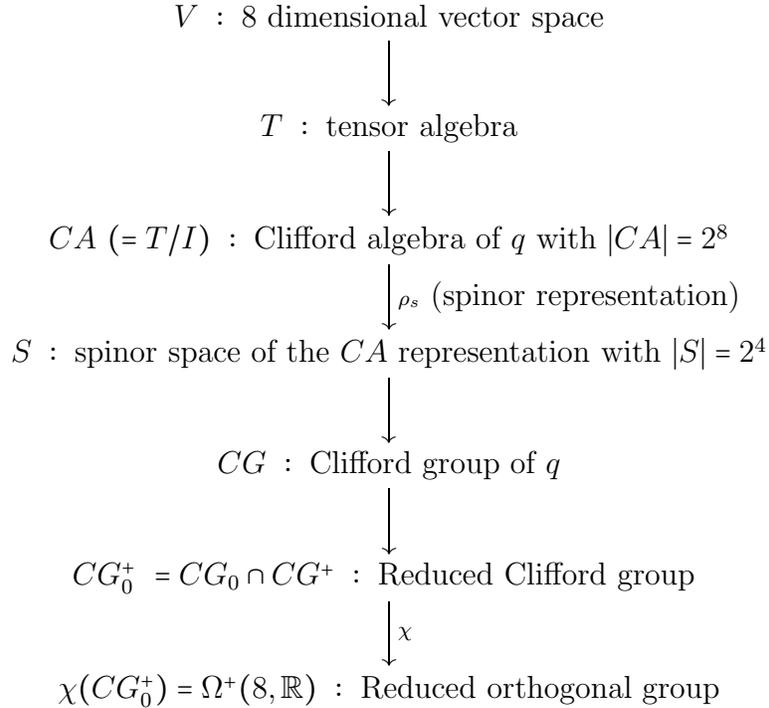

$\bullet$\ $T = T^+ \oplus T^- = \sum_{k\ even} \otimes^k V \oplus
\sum_{k\ odd} \otimes^k V = F \oplus V \oplus (V \otimes V) \oplus \ldots$

$\bullet$\ The ideal $I$ of $T$ is generated by $x \otimes x - q(x) \cdot 1$.

$\bullet$\ $CG = \langle \text{invertible\ elements}\ s \in CA: sxs^{-1} \in V,\
\forall x \in V \rangle$.

$\bullet$\ $CG^+ = CG \cap CA^+$ is the special Clifford group
of $q$, where $CA^+ = T^+ / (I \cap T^+)$.

$\bullet$\ The `main anti-automorphism' of $CA$ is a linear map on $T$
given by $\alpha(x_1 \otimes \ldots \otimes x_k) = x_k \otimes \ldots \otimes x_1$.

$\bullet$\ $CG_{0} = \{s \in CG : \alpha(s) s = 1 \}$, $\alpha$
is the main anti-automorphism of $CA$.

$\bullet$\ $CG^{+}_{0}$ has three inequivalent $8$-dimensional representations
via $\rho$ on $V$ and the half-spinor spaces $S^+, S^-$.

$\bullet$\ The linear map $\chi: x \mapsto s x s^{-1}$ is a vector
representation of $CG$.

$\bullet$\ $|O^{+}(8, \mathbb{R}) : \chi(CG^+)| = 2$, and $\Omega^{+}(8, \mathbb{R})$ is the
commutator subgroup of $O^{+}(8, \mathbb{R})$.

$\bullet$\ If $char(F) \neq 2$ then $|\chi(CG^+)| = |SO^{+}(8, \mathbb{R})|$.  If $char(F) = 2$ then every element of $O^{+}(8, \mathbb{R})$ has unit determinant; the
determinant of a reflection is $-1$. 

$\bullet$\ $\rho_{s}(CG^{+}_{0})$ preserves a quadratic form $j$ on $S$.

$\bullet$\ There is an order $3$ mapping $J: V \to S^+ \to S^- \to V$ such
that:
\[ j(J(x)) = q(x) \]
\[ j(J(u)) = j(u) \]
\[ q(J(u')) = j(u'), \]
for $x \in V, u \in S^+, u' \in S^-$, and for some quadratic form $j$.

$\bullet$\ the map $J$ yields a \emph{triality} map $\tau_{CA}$ that gives outer
automorphisms of $P\Omega^{+}(8, \mathbb{R})$.  As mentioned at the beginning of the section, there is a geometrical version of this algebraic triality due to E. Study,
\index{Study, E.}%
 mapping $1$-dimensional singular subspaces to $4$-dimensional totally singular subspaces; an account of this is given in~\cite{cam20}.

\bigskip

Triality is only non-trivial for orthogonal groups in 8 dimensions.  Put slightly differently, it is only in 8 dimensions that there is a Moufang loop
\index{Moufang ! loop}%
 available to create a non-trivial triality; see Appendix~\ref{LoopTheory}.

Note from~\cite{conwaysm} that $\Spin(7) \cap \Spin(7) \subset \Spin(8)$
\index{group ! Spin(8)@$\Spin(8)$}%
 for any two of the
three pairs of $\Spin(7)$.   If $\alpha, \beta, \gamma \in \SO(8)$
\index{group ! SO(8)@$\SO(8)$}%
 then the surjection $\Spin(8) \to
\SO(8): \{(\alpha, \beta, \gamma), (- \alpha, - \beta, \gamma) \mapsto
\gamma\}$ is $2$--$1$.  In fact one way of knowing that the Cartan triality
\index{Cartan triality}%
automorphisms are outer is because they act nontrivially on the centre of $\Spin(8)$ which is 
\[\{(1,1,1),\ (1,-1,-1),\ (-1,1,-1),\ (-1,-1,1)\}.\]

The group $\Spin(7)$ is `almost' as big as $SO(8)$ in the sense that $\Spin(7)$ also acts transitively
on spheres in $\mathbb{R}^{8}$~\cite[p.283]{harvey}.  The elements of
$\Spin(7)$ are those elements of $\Spin(8)$ that stabilize a vector in
one of the three eight-dimensional vector spaces $V^{8}, S^{+}, S^{-}$ representing the three irreducible
representations
\index{representation ! irreducible}%
 of $\PSO(8)$.  The three copies of $\Spin(7)$ in $\Spin(8)$
have as elements the triples $(\alpha, \beta, \gamma)$ for which
$\gamma$ (respectively $\alpha, \beta$) fix
$1$~\cite[p.~95]{conwaysm}.  There are $3 \choose 2$ ways of doing this and the triality automorphism $\tau_C$ (where $\tau_C^3 = 1$) maps each of the 3 non-isomorphic copies of $\Spin(7)$ into $\Spin(8)$.  

\bigskip

Let us discuss some of the occurences of the Cartan outer automorphisms.  An outer automorphism is a coset of inner automorphisms
\index{automorphism ! inner}%
\index{automorphism ! outer}%
and it permutes the conjugacy classes of representations of a group.
The Lie group $D_4$
\index{group ! D(4)@$D(4)$}%
\index{group ! Lie}%
is the compact, connected and simply-connected group $\Spin(8)$, and has the largest outer automorphisms of all the Lie groups, these being isomorphic to $\Sym(3)$.  The exceptional Jordan algebra
\index{Jordan ! exceptional algebra}%
 $\mathfrak{h}_3(\mathbb{O})$ is a $27$-dimensional real vector space consisting of all $3 \times 3$
Hermitian symmetric matrices over $\mathbb{O}$ of the form
$ \begin{pmatrix} 
a&z&\bar{y}\\
\bar{z}&b&x\\
y&\bar{x}&c
\end{pmatrix}$ with $a, b, c \in \mathbb{R}$ and $x, y, z \in
\mathbb{O}$.  If $\Spin(8)$ is embedded in $F_4 = \Aut(\mathfrak{h}_3(\mathbb{O}))$ then the $\Sym(3)$ action can be effected by conjugation by elements of $F_4$~\cite{adams}.  The finite versions of Cartan triality
\index{Cartan triality}%
arise by taking the simple groups
\index{group ! simple}%
$G=D_4(\mathbb{F}_q)=P\Omega^{+}(8, \mathbb{F}_q)$
over finite fields.  If $\mathbb{L}_{E_8}$ is the root lattice of the $E_8$ Lie algebra and $W_{\mathbb{L}_{E_8}}$ is its Weyl group, then the determinant
one subgroup of $\Aut(\mathbb{L}_{E_8})/\{\pm 1\}$ is isomorphic to $SO^{+} (8,
\mathbb{Z} / 2 \mathbb{Z})$, or alternatively $W_{\mathbb{L}_{E_8}}/\zeta(W_{\mathbb{L}_{E_8}})
\cong W_{E_8}/\{\pm 1\} \cong O^{+} (8, 2)$. The group $P\Omega^{+}(8, \mathbb{F}_q)$ which has index $2$ in $PSO^{+}(8, \mathbb{F}_q) = SO^{+}(8, \mathbb{F}_q) / (\{\pm\} \cap SO^{+}(8, \mathbb{F}_q))$, has its own Cartan
triality outer automorphism group.
\index{automorphism ! outer}%

The group $\Spin(n)$ is the simply connected covering group of $\SO(n)$
    and the group homomorphism $\pi : \Spin(n) \to \SO(n)$
    has kernel of order $2$ if $char(F) \neq 2$ and kernel of order
    $1$ if $char(F) = 2$.  The reduction mod
    $2$ homomorphism $\SL(2, \mathbb{Z}) \to \SL(2, 2) \cong
    \Sym(3)$,
\index{group ! SL()@$\SL(2, \mathbb{Z})$}%
    has as kernel matrices congruent to the identity matrix mod $2$.  There is an action of $\SL(2, \mathbb{Z})$ on $\Spin(8)$ giving the
    semidirect product $Spin(8) \sd \SL(2, \mathbb{Z})$ that
    characterises \emph{S-duality},
    \index{S-duality}%
conjectured to be an exact symmetry of full string theory~\cite{gaunt}~\cite{seibwit}, whereby the $\Spin(8)$ vector $(v)$, spinor $(s)$
\index{spinor (spin representation)}%
 and conjugate spinor $(c)$ representations of the number of magnetic and electric charges, $(n_m, n_e)$, these being coprime integers, are
    determined by the following hypermultiplet of states:
\[ (0, 1) \to \mathbf{8_v};\   (1, 0) \to \mathbf{8_s};\   (1, 1) \to
    \mathbf{8_c}. \]
That the irreducible representations, $\mathbf{8_v}, \mathbf{8_s}, \mathbf{8_c}$, of the Lie algebra $\mathfrak{so}_8$ are related to one another by the triality automorphism was known to E. Study~\cite{study}.  It may be worth pursuing a new graph-theoretic approach to the physics
using the idea that (0, 1), (1, 0) and (1, 1) represent the $3$ types of
adjacency of $\mathfrak{R^{t}}$.

The group $\Spin(8) \sd \Sym(3)$ is the automorphism group
\index{group ! automorphism}%
 of many different algebras, and yields invariant subgroups of exceptional type $G_2$~\cite{knus}.
\index{group ! exceptional}%

We list some of the theories where the appearance of Cartan triality
\index{Cartan triality}%
 has been essential to a description of the phenomena.  

In~\cite{shankar} an application of triality is given to solvable models, one a quantum Ising model of spins in statistical mechanics and the other the $O(8)$ Gross-Neveu model in quantum field theory; it is shown that this triality is the origin of the fermionic solutions.

In~\cite{godolsc} the basic representation of the $E_8$ Kac-Moody algebra
\index{Kac-Moody algebra}%
 is constructed out of $16$ real Fermi fields, themselves constructed out
of $8$ real Bose fields, and the triality properties of the
eight-dimensional orthogonal group were used
in a crucial way.  In an analysis of the algebraic structure of
fermions interrelated by the vertex operator construction~\cite{godnaolrsc}, triality
is realized via orbit triples, just as it will be in our Leech lattice ($\mathbb{L}_L$)
\index{lattice ! Leech}%
construction of $\mathfrak{R^{t}}$, (see the main text).

More recently the triality automorphism has had other manifestations.  For
example, as an automorphism of both the Griess algebra~\cite{gri}
\index{Griess Jnr., R. L.}%
 whose automorphism group
\index{group ! automorphism}%
  is the Monster group $\mathbb{M}$
\index{group ! Monster}%
 and of the Monster module
$V^{\natural}$~\cite{frenkellm},
\index{Monster ! module}%
as well as inside the lattice vertex operator
algebra (VOA) constructed from a marked binary doubly-even self-dual
code~\cite{griessho}.  (A marked, length $2n$ binary code is a partition of the $2n$
coordinates into $n$ sets of size $2$).  The results
of~\cite{frenkellm} giving the natural representation of $\mathbb{M}$
as a conformal field theory were extended~\cite{dolangm}~\cite{dolangm1} to produce untwisted
and $\mathbb{Z}_2$-twisted theories constructed from even self-dual
lattices or binary codes which also possess a triality.  In each case
the twisting or intertwining effected by the triality is essential to
the construction.  A $\mathbb{Z}_3$ analogue of the
$\mathbb{Z}_2$-twisted construction of~\cite{frenkellm} is given
in~\cite{xu} where a ternary moonshine space is constructed which is
closed under intertwining operators and in some instances is an
operand under automorphic braid group action.  
\index{group ! braid}%
\index{group ! action}%
\head{Codicil}

The objects and phenomena of most depth,
beauty and importance in mathematics are often the exceptional ones; see the review article~\cite{baez}.  

Consider the following:

\centerline{Cayley-Graves Octonions $\mathbb{O}$}
\index{octonions}%
\centerline{Exceptional Lie Groups (e.g. $G_2$, $F_4$)}
\index{group ! exceptional}%
\centerline{Moufang Projective Plane $\mathbb{O}\mathbb{P}^2$}
\centerline{Exceptional Jordan Algebra $\mathfrak{h}_3(\mathbb{O})$}
\index{Jordan ! exceptional algebra}%
\centerline{Cartan Triality Outer Automorphism}
\index{Cartan triality}%
\centerline{Exceptional Hopf Map}
\index{Hopf, H.}%
\centerline{Bott Periodicity}
\index{Bott, R.}%
\centerline{Leech Lattice $\mathbb{L}_L$}
\index{lattice ! Leech}%
\centerline{Fischer-Griess Monster Group $\mathbb{M}$}
\index{Griess Jnr., R. L.}%
\index{Fischer, B.}%

On this list are a number system, a small set of continuous groups, a projective
plane, an algebra, a coset of a group action,
\index{group ! action}%
 a fibration, a topological property, a lattice
and a finite group.  In one sense or another, these are all considered to be exceptions in their class and it would be interesting to find more among other branches of mathematics.  For example consider singularity theory
which amongst other things classifies singular objects in terms of
their complexity.  Some of the types are \emph{stable} singularities
which are immutable to small
perturbations, whilst \emph{simple} ones arise as a finite number of
non-equivalent cases under small perturbations, and Kleinian
singularities
\index{group ! Kleinian singularity}%
are associated with finite subgroups of $\SL(2, \mathbb{C})$.  Arnold's
\index{Arnold, V. I.}%
 strange duality in this field has been related to $\mathbb{L}_L$~\cite{ebel}, so
perhaps this points to a candidate special singularity.

As an example of the almost inexhaustible number of applications of $D_4$
\index{group ! D(4)@$D(4)$}%
 is the construction by Kantor~\cite{kant}
\index{Kantor, W. M.}%
of a geometry which is almost a building whose diagram is the extended $D_4$, in which a middle node has four neighbours.

J. Milnor
\index{Milnor, J.}%
 proved in $1956$ that there are $28$ different smooth
structures on the manifold $S^7$, with $\mathbb{Z}_{28}$ acting on
them.  The octonionic Hopf map is the projection $p_{\mathbb{O}} :
E_{\mathbb{O}} \to \mathbb{O}\mathbb{P}^1$, with total space
$E_{\mathbb{O}}$ consisting of all unit vectors in $\mathbb{O}^2$, so
is a $2.8-1=15$-dimensional sphere $S^{15}$, and with base space $S^8$.  The projection is a Hopf bundle of $7$-spheres over
$\mathbb{O}\mathbb{P}^1$ with a unit $7$-sphere in each fibre.  Bott Periodicity is the property of the homotopy groups of the
topological group
\index{topological group}%
 $O(\infty)$ (a direct limit
\index{direct limit}%
 of Orthogonal groups
$O(n)$ as $n \to \infty$) given by $\pi_{i+8}(O(\infty)) \cong
\pi_{i}(O(\infty))$.  The above objects are inter-related and are all exceptional in some real sense~\cite{slob} \cite{stil}.  This is a subjective list and perhaps other entries can be added.  The
octonions
\index{octonions}%
 are thought to be the linchpin of
exceptionality~\cite{baez} \cite{stil} and it
appears that Cartan triality
\index{Cartan triality}%
 is a signature.  The concrete
realization of $\Spin(8)$ using the octonions has several
applications~\cite[p.279]{harvey}.  The two exceptional tesselations of $\mathbb{R}^2 = \mathbb{C}$ are both built
on the lattice of Eisenstein integers~\cite{stil},
\index{Eisenstein integers}%
in terms of which we will construct
$\mathfrak{R^{t}}$ in one of the chapters.  Analogously in the lattice of `integer octonions' in
$\mathbb{R}^8$, the neighbours of each lattice point form the
Gosset polytope which is not regular but has $E_8$ as its symmetry
group.  Finally, the idea that octonions can be made to be commutative
and associative within the context of a symmetric monoidal
category~\cite{alberm} by transferring the nonassociativity to a
suitable product operation, may have an important part to play in future developments.

Whilst dualities are in plentiful supply in
mathematics, trialities require a subtleness that makes them rarer.
Equally, outer automorphisms
\index{automorphism ! outer}%
 are much less frequently occurring than the
inner ones got by conjugation.  Any division algebra yields a
triality, so trialities potentially occur in dimensions $1, 2, 4, 8$.  B. Kostant
\index{Kostant, B.}%
found that the outer automorphisms of $\Spin(8)$ are the \emph{reason}
for the existence of the five exceptional Lie groups~\cite{kost}~\cite{rota}.
\index{group ! exceptional}%
 It would be nice to amass sufficient evidence to add $\mathfrak{R^{t}}$ as an asterisked footnote to the above list.

Another piece of evidence for the theory that is often asserted that the exceptional Lie groups all exist because of the octonions
\index{octonions}%
   is the Freudenthal-Tits
\index{Tits, J.}%
  magic square,
\index{Freudenthal-Tits magic square}%
   and its variations; see for example~\cite{barton}~\cite{freudenthal}~\cite{freudenthal1}~\cite{titsa}~\cite{vinberg}.  

\section{Loop Theory and Groups with Triality}
\label{LoopTheory}

A \emph{magma}
\index{magma}%
is a set with a single binary operation.  A \emph{semigroup}
\index{semigroup}%
is a magma whose binary operation is associative.  Let $L$ be a \emph{quasigroup}
\index{quasigroup}%
that is a \emph{magma} with unique division, and $\Mlt(L)$ the \emph{multiplication group}
\index{loop ! multiplication group}%
of $L$, which is generated by the set of all left and right \emph{translation maps}
\index{loop ! translation map}%
 that is the permutations defined by $L_a (x) = ax$ and $R_a (x) = xa$ for every $x \in L$.  Clearly, $M(L)$ acts transitively on $L$ and the stabilizers of the elements of $L$ are conjugated in $L$.   Another definition of a \emph{quasigroup} is a magma $(G, \cdot)$ in which the maps $L(a) : G \to G$ and $R(a) : G \to G$ are bijections for all $a \in G$.  (We should say that \emph{unique division} means that $\forall l_1, l_2$ in a quasigroup $(L, \circ)$ there exist unique elements $x, y \in L$ such that: $l_1 \circ x = l_2$ and $y \circ l_1 = l_2$.  The unique solutions to these equations are written $x = l_1 \backslash l_2$ and $y = l_2 \slash l_1$. The operations $\backslash$ and $\slash$ are called, respectively, left and right division.)  
  
For a group $G$, the right multiplication group is isomorphic to $G$ (thus proving Cayley's Theorem.)  For an arbitrary quasigroup, the right multiplication group will be larger than $G$ (it will act transitively but not regularly on $G$).  In fact, for a random finite quasigroup, the right multiplication group is the symmetic group (all permutations of $G$) with high probability.  Some of the best structural results on quasigroups come from studying their multiplication groups.

A \emph{loop}
\index{loop}%
of order $n$ is a set $L$~\label{Lloop} of $n$ elements with an identity element $e$ satisfying $ex = xe = x\ \forall x \in L$ and such that for $a, b \in L$, the equations $ax = b$ and $y a = b$ each has a unique solution in $L$.  The multiplication table is a Latin square,
\index{Latin square}%
which is a square array of $n^2$ cells containing the numbers $1, \ldots,n$ in such a way that no number appears twice in either the same row or in the same column.  The Latin square may be regarded as a `global' object as compared to a `local' loop.  The binary operation defined by a Latin square need not be associative.  If the Latin square is symmetric, then the corresponding loop is commutative.  A quasigroup is precisely a structure whose operation table is a Latin square;
\index{Latin square}%
  it is a loop precisely when, assuming that the first element is the identity, the first row and column are the same as the row and column labels.

If $L$ is a \emph{loop},
\index{loop}%
that is a quasigroup
\index{quasigroup}%
 with a neutral/identity element $e$ and inverses, then the \emph{inner mapping group}
\index{loop ! inner mapping group}%
$I(L)$ of $L$ is the stabilizer of $e$; this is the analogue for loops of the inner automorphism group
\index{group ! automorphism}%
 of a group.  The loop $L$ is an abelian group
\index{group ! abelian}%
if and only if $I(L) = 1$.  If $L$ is a group, then $I(L)$ consists of the inner automorphisms of $L$. 

Note one difference between the theories of semigroups and quasigroups.  In a semigroup, the associative law holds, so we can write products unambiguously.  In a quasigroup or loop, we would have to specify which of the possible bracketings of a product is intended.

For a loop $L$, the transversals $A = \{L_a : a \in L\}$ and $B = \{R_a : a \in L\}$ are \emph{$I(L)$-connected} to $I(L)$ in $\Mlt(L) = \langle A, B \rangle$ if $[A, B] \leq I(L)$.
 
Many of the concepts that arise in group theory have their equivalents in loop theory.  For example, a non-empty subset $H$ of a set $G$ is a \emph{subloop}
\index{loop ! subloop}%
 of a quasigroup $(G, \cdot)$ if $(H, \cdot)$ is a loop.  If $H$ is a subloop
\index{loop ! subloop}%
of a loop $L$, and $L$ has a left and right decomposition modulo $H$ then $H$ is a \emph{normal subloop}
\index{loop ! normal subloop}%
if $l H = H l$, $(l_1 H) l_2 = l_1 (H l_2)$ and $l_1 (l_2 H) = (l_1 l_2) H$.  A \emph{simple loop}
\index{loop ! simple}%
 is one with no non-trivial proper homomorphic images, that is it has no non-trivial proper normal subloops; alternatively phrased, if every surjective loop homomorphism is either bijective or has image the identity.  Let $H$ be a subloop and $K$ a subset of a loop $L$.  Then $K$ is a \emph{left (right) coset modulo} $H$ if $K = aH$ $(K = Ha)$ for some $a \in L$.  The \emph{exponent} $exp(L)$ of a loop $L$ is the smallest positive integer $n$ for which $l^n = 1$ for all $l \in L$.

For a quasigroup
\index{quasigroup}%
 $L$ it is well-known that $\Mlt(L)$ is primitive on $L$ if and only if $L$ is simple.  In~\cite{phillips3}
\index{Phillips, J. D.}%
\index{Smith, J. D. H.}%
 prove that $L$ is simple if and only if $\Mlt(L)$ is quasiprimitive
\index{group ! permutation ! quasiprimitive}%
 on $L$; (a quasiprimitive group has only transitive nontrivial normal subgroups, see Appendix~\ref{PermutationGroups}).  So transitive imprimitive quasiprimitive permutation groups cannot be multiplication groups of quasigroups~\cite{praeger1}.
\index{quasigroup}%
  Contrast this with the observation that for almost all quasigroups, $\Mlt(L) = \Sym(L)$.
\index{group ! symmetric}%
This fact follows either from the result of {\L}uczak and Pyber~\cite{luczak}
\index{Luczak@{\L}uczak, T.}%
\index{Pyber, L.}%
 that for almost all $g \in \Sym(L)$ the only transitive subgroups of $\Sym(L)$ containing $g$ are $\Sym(L)$ and perhaps $\Alt(L)$, or the result of H\"{a}ggkvist and Janssen~\cite{haagkvist} that $\Mlt(L)$ is almost never contained in $\Alt(L)$.   Other than books on loops, a good reference for multiplication groups of loops is~\cite{niemenmaa}.

A \emph{quasigroup homotopy}
\index{quasigroup ! homotopy}%
 from quasigroup $G_1$ to $G_2$ is a triple $(\alpha, \beta, \gamma)$ of maps from $G_1$ to $G_2$ such that
\[ \alpha(x) \beta(y) = \gamma(xy) \]
for all $x, y \in G_1$. A quasigroup homomorphism is just a homotopy for which the three maps are equal.
An \emph{isotopy}
\index{quasigroup ! isotopy}%
 is a homotopy for which each of the three maps $(\alpha, \beta, \gamma)$ is a bijection. Two quasigroups are \emph{isotopic} if there is an isotopy between them. In terms of Latin squares, an isotopy $(\alpha, \beta, \gamma)$ is given by a permutation of rows $\alpha$, a permutation of columns $\beta$, and a permutation on the underlying element set $\gamma$.  An \emph{autotopy}
\index{quasigroup ! autotopy}%
 is an isotopy from a quasigroup to itself. The set of all autotopies of a quasigroup form a group with the automorphism group
\index{group ! automorphism}%
 as a subgroup.  Each quasigroup is isotopic to a loop.

A \emph{Moufang loop,}~\label{QMoufangloop}
\index{Moufang ! loop}%
 named after Ruth Moufang~\cite{moufang}
 \index{Moufang, R.}%
 and denoted $Q$, is a set of elements with a
binary operation, a unique product $xy$ for every two elements $x, y$, an identity element $1$, a unique inverse $x^{-1}$ of $x$ such that $x^{-1} \cdot x = 1 = x \cdot x^{-1}$, that satisfy a weak form of associativity given by any of the following equivalent conditions:
\[ [x(zx)]y = x [z(xy)]\quad\quad\quad\quad\quad  [(xz)x]y = x[z(xy)] \]
\[ (xy) (zx) = x [(yz) x]\quad\quad\quad\quad\quad  [(yx)z]x = y [x(zx)]. \]

The last two conditions can be equivalently written as the following multiplication laws~\cite[p.~74]{conwaysm}
\[ (xzx) \cdot x^{-1}y = x(zy) \]
\[ (yz)x = yx^{-1} \cdot (xzx). \]

Nonassociativity does not arise in compositions of pairs of Moufang elements.  (Note that for permutation groups composition is always associative, and it only becomes necessary to invoke the associativity axiom for abstract groups.)

 For loops in general, $I(L)$ is not an automorphism group
\index{group ! automorphism}%
  but it is for commutative Moufang loops,
\index{Moufang ! commutative loop (CML)}%
 that is loops satisfying $xx \cdot yz = xy \cdot xz$.

If $M$ is a simple nonassociative Moufang loop, and $L(x):= m \to xm, R(x):= m \to mx$ $(m \in M)$, ~\label{L(M)} then $G(M) : = \langle L(x), R(x): x \in M \rangle$ is a group with triality.  Every Moufang loop with trivial nucleus has an associated group with triality; for more such results see~\cite{phillips}~\cite{phillips1}.

Commutative Moufang loops (CML)
\index{Moufang ! commutative loop (CML)}%
 can be defined by the relation
\[ x^2 (yz) = (xy) (xz) \]
and have the property that all cubes are central, that is the mapping
 $l \mapsto l^3$ is a centralizing endomorphism.  Every loop-isotope
 of a CML also has centrally endomorphic cubing whilst commutativity
 itself is not an isotopic invariant.  R. H. Bruck
\index{Bruck, R.H.}%
 proved that any Moufang loop with centrally endomorphic cubing is
 either an isotope of a CML or is an index-3 normal subloop of such an isotope.  It is possible for elements of a
 CML to have infinite order.  If $L$ is commutative then the
 square-root set of the identity is a subloop $L'$ consisting of all
 elements of order $1$ or $2$, and any other square-root set is empty
 or a coset of $L'$.  (This follows because $x^{-1} y^{-1} x^{-1} y =
 1 \Rightarrow y = x y x \Rightarrow y z = (((x) y) x) z \Rightarrow y
 z = (x y) (x z) \Rightarrow x^2 (y z) = (x y) (x z)$, satisfying the
 condition for a CML), and the last equation implies that the elements
 of $L'$ have order 1 or 2.  Furthermore $L'$ is a subloop because of
 closure of pre- and post-multiplication of any element of $L'$ by
 another).  So if $|L : L'| = \infty$, for example if $L'$ has unbounded exponent, then $L$ is $\mathfrak{R}_{m}$-genic.
\index{graph ! Rmgenic@$\mathfrak{R}_{m}$-genic}%

The \emph{nucleus}
\index{Moufang ! nucleus}%
 $\Nuc(Q)$~\label{Nuc} of a Moufang loop $Q$ is defined by $\Nuc(Q) = \{x \in Q :
\forall y, z \in Q, (xy)z = x (yz)\ \&\ cyc.\}$ If $C(Q) : = \{x \in Q
: \forall y \in Q, xy = yx \}$ is the \emph{Moufang centre}
\index{Moufang ! Moufang centre}%
 then the \emph{centre}
\index{Moufang ! centre}%
 is $\zeta(Q) := C(Q) \cap \Nuc(Q)$.  Finally
if $(L_1, \cdot)$ and $(L_2, \circ)$ are two loops and $\phi: L_1 \to
L_2$ such that $(x \cdot y) \phi = x\phi \circ y \phi$, then $\phi$ is
a \emph{homomorphism}
\index{Moufang ! homomorphism}%
 of $(L_1, \cdot)$ into $(L_2, \circ)$.  A homomorphic image of a loop
 is not always a loop~\cite[p.~28]{pflug}, but we shall work only with those which are. 

There is an associated incidence structure called a \emph{3-net}~\cite{bruck},
\index{ta@3-net}%
where the disjoint line classes are broken down as $\mathcal{L} =
\mathcal{L}_1 \sqcup \mathcal{L}_2 \sqcup
\mathcal{L}_3$.  Lines in each class are disjoint; lines
from different classes meet in one point.  The line classes are
subsets of the point set $\mathcal{P} = Q \times Q$, and are given by
$\mathcal{L}_1 = \{\{(x, c)\}:\ x \in Q\}:\ c \in Q\}$, $\mathcal{L}_2
= \{\{(c, y)\}:\ y \in Q\}:\ c \in Q\}$, $\mathcal{L}_3 = \{\{(x,
y)\}:\ x, y \in Q\}:\ xy=c\}:\ c \in Q\}$; the element $c$ is a
constant.  A \emph{Bol reflection} 
\index{Bol, reflection}%
about a fixed axis, say $x = \text{constant}$, is a particular permutation of the point set that stabilizes the
line forming the axis; this is described in~\cite{hallnagy}.  Bol proved that a 3-net is coordinatised by a Moufang loop if and only if any Bol reflection is a collineation.  A \emph{Latin square design}
\index{Latin square design}%
 $D$ is a pair of points $P$ and lines $A$ (subsets of $P$) such that (i) $P$ is the
disjoint union of 3 parts $R, C, E$, (ii) every line $l \in A$ contains
exactly three points, meeting each of $R, C, E$ once, (iii) every pair of
points from different parts belong to exactly one line.  Equivalently,
it is a Latin square where the entry in row $x$ and column $y$ is $z$
$\Longleftrightarrow \langle x, y, z \rangle \in A$.  The partial
automorphism $\tau_x \in \Aut(D)$ acting as $\tau_x(y)=z,
\tau_x(z)=y$ $\Longleftrightarrow \langle x, y, z \rangle \in A$
extends to a permutation on all $P$ such that $\{\tau_p |\ p \in P,
\tau_p \in \Aut(D)\}$ is a set of involutions.  If $p$ and $q$ belong
to different parts of $P$ and $\tau_p, \tau_q \in \Aut(D)$, then $(
\tau_p \tau_q)^3=1$.  So $\tau_x$ extends to a permutation on all
$P$; it is a \emph{central automorphism}
\index{loop ! central automorphism}%
\index{automorphism ! central}%
of $L$ with centre $x$.  More on Latin square designs can be found in~\cite{hallji1}~\cite{hallji1}. 

Bol proved the following~\cite{bol}:

\begin{theorem}
Let $(L, \cdot)$ be a loop.  Then for every point $x$, $\tau_x \in \Aut(D)$ if and only if
\[ (xa)(bx) = (x (ab))x,\quad\quad x, a, b \in L. \]
\end{theorem}

In other words Moufang loops satisfy this identity.  Thus to every Moufang loop there is an associated automorphism group generated by a conjugacy class of involutions obeying

\begin{proposition}
In $\Aut(D)$ there is at most one central automorphism $\tau_p$ with center $p$ for each $p \in \mathcal{P}$. If $\tau_p \in \Aut(D)$, then it has order 2. If $\tau_p, \tau_q \in \Aut(D)$ with $p, q$ in different fibres, then $\tau_p \tau_q$ has order 3 and $\langle \tau_p, \tau_q \rangle \cong \Sym(3)$.  If this is the case, then 
there is a unique conjugacy class $T$ of central automorphisms in $\Aut(D)$.
\index{automorphism ! central}%
\end{proposition}

The converse relates Moufang loops to \emph{groups with triality}
\index{group ! with triality}%
 $G = \langle T \rangle$ as studied by Glauberman~\cite{glauberman}
\index{Glauberman, G.}%
 and Doro~\cite{doro}.  It states~\cite{hallji2},
\index{Doro, S.}%

\begin{theorem}
Let $T$ be the conjugacy class of involutions in $G = \langle T \rangle$, and let $\pi : G \to \Sym(3)$ be an epimomorphism.  Assume that $G$ satisfies the hypothesis that if $\pi(t) \neq \pi(r)$ then $|\pi(t)\pi(r)| = 3$, for all $t, r \in T$.

Then there is a Moufang loop $(q, \cdot)$ with
\[ G / \zeta(G) \cong \Aut(D(Q, \cdot))^0, \]
where the class $T$ maps bijectively to the class of central automorphisms
\index{automorphism ! central}%
 of $\Aut(D(Q, \cdot))^0$, the subgroup of $\Aut(D(Q, \cdot))$ generated by all central automorphisms.
\end{theorem}

There is a category equivalence between loops and Latin square designs, a special case of which is a category equivalence between Moufang loops and those Latin square designs that admit all possible central automorphisms.
\index{automorphism ! central}%
 These are in turn equivalent to an appropriate category of groups with triality~\cite{hallji2}.

The reason for the difference between Moufang loops and groups with triality is that different groups with triality correspond to the same central Latin square design (one admitting a central automorphism at every point).  The problem has to do with centre of groups with triality.  To every such group $(G, T, \pi)$ there is a central quotient $(G^A, T^A, \pi^{A})$ where the centre is factored out, and a \emph{universal group with triality}
\index{group ! with triality ! universal}%
 $(G^U, T^U, \pi^{U})$ which is a maximal (so universal) central extension.  The distinction between these two classes gives the equivalence between the category of Moufang loops and the category of \emph{universal} groups with triality.

In Chapter~\ref{homcochap} we gave a pr\'ecis of the theory of groups with triality as expounded in~\cite{hallnagy}, via its relation to Moufang loop theory, elements of which was used in the main text, but we need to repeat the definition.

Doro~\cite{doro}
\index{Doro, S.}%
defined a pair $(G, S)$~\label{(G,S)} to be a \emph{group with triality}
\index{group ! with triality}%
if $G$ is an abstract group, $S \le \Aut(G), S = \langle
\sigma, \rho\ |\ \sigma^2 = \rho^3 = (\sigma\rho)^2 = 1 \rangle \cong
\Sym(3)$, and the \emph{triality identity}
\index{triality identity}%
given by $[g, \sigma]
[g, \sigma]^{\rho} [g, \sigma]^{\rho^2} = 1$ or equivalently $(\sigma^g
(\sigma \rho^2))^3 = 1$ holds for all $g \in G$.  The group with triality of a group $H$ is (essentially) the wreath product $H \wr \Sym(3)$.

Two loops $L_1, L_2$ are \emph{isotopic}
\index{loop ! isotopic}%
if there exist bijections $\alpha, \beta, \gamma: L_1 \to L_2$ with $x^{\alpha} y^{\beta} = (x y)^{\gamma}$.  Isotopic groups are isomorphic, but this is not true in general for loops or Moufang loops.

\begin{theorem}[J. I. Hall and G. P. Nagy]
\index{Hall, J. I.}%
\index{Nagy, G. P.}%
The following are equivalent
\begin{itemize}
\item[(i)]  Groups $(G, S)$ with triality and centre $\zeta(GS)={1}$;
\item[(ii)]  Latin square designs
\index{Latin square design}%
in which $\tau_p$ extends to a central automorphism
\index{automorphism ! central}%
 for every point $p$;
\item[(iii)]  Moufang $3$-nets.
\index{Moufang ! 3-net}%
\item[(iv)]  Isotopy classes
\index{Moufang ! isotopy class}%
 of Moufang loops.
\end{itemize}
\end{theorem}

In view of Theorems~\ref{outgpwttr} and~\ref{autsymaut} we state for comparison the following result~\cite{nagyvoj} which parallels Theorem~\ref{centthree},

\begin{theorem}[Nagy and Vojt\v{e}chovsk\'y]
\index{Nagy, G. P.}%
\index{Vojt\v{e}chovsk\'y, P.}%
\label{nagyvothm}
Let $Q$ be a Moufang loop
\index{Moufang ! loop}%
and $\mathcal{N}$ be its associated 3-net.  Let $G_0$ be the group of collineations generated by the Bol reflections of $\mathcal{N}$, $G$ the direction preserving part of $G_0 = G \sd S$, and $S \cong \Sym(3)$ the group generated by the Bol reflections whose axis contains the origin of $\mathcal{N}$.  Then $\Aut(Q) \cong C_{\Aut(G)} (S)$.
\end{theorem}

Using the correspondence between simple Moufang loops and simple groups with triality that Doro~\cite{doro} and Glauberman~\cite{glauberman} had discovered, of which by Theorem~\ref{outgpwttr} the group $\Aut(\mathfrak{R^{t}})$ is an example, Liebeck
\index{Liebeck, M. W.}%
 proved~\cite{liebeck} that a finite simple Moufang loop is either associative (and so is a finite simple group),
\index{group ! simple}%
or is isomorphic to a Paige loop~\cite{paige}.
\index{Paige loop}%
In addition to the normed division algebras, the other composition algebras over $\mathbb{R}$, are the split-complex numbers, the split-quaternions and the split-octonions.
\index{split-octonions}%
  Let $Oct(K)$ be the eight-dimensional algebra of split-octonions over the field $K$.  Its units form a Moufang loop, and $SOct(K)$ is the subloop of norm-$1$ octonions.  Factoring out the normal subloop $\{\pm 1\}$ of $SOct(K)$ gives a \emph{Paige loop} $PSOct(K)$.  A Paige loop is a simple Moufang loop that is isomorphic to all its loop isotopes.  

\begin{theorem}[Grishkov and Zavarnitsine~\cite{grishkovza}]
\index{Grishkov, A. N.}%
\index{Zavarnitsine, A. V.}%
Every finite Moufang loop $Q$ contains a uniquely determined normal series
\[ 1 \leq \Gr(Q) < Q_0 \leq Q \]
such that $Q / Q_0$ is an elementary abelian $2$-group,
\index{group ! elementary abelian}%
 $Q_0 / \Gr(Q)$ is the direct product of simple Paige loops $Q(q)$ (where $q$ may vary), the composition factors of $\Gr(Q)$ are groups, and $\Gr(Q / \Gr(Q)) = 1$.
\end{theorem}

For $\mathbb{F}_q$ odd, the Paige loop $PSOct(\mathbb{F}_q)$ has a two-fold extension isomorphic to the loop $PGL(\mathbb{O}(\mathbb{F}_q))$, where $\mathbb{O}(\mathbb{F}_q)$ is the Cayley algebra over $\mathbb{F}_q$.  Denoting this extension by $Q(\mathbb{F}_q) . 2$, define

\begin{displaymath}
\widehat{Q(\mathbb{F}_q)} = \left\{ \begin{array}{ll}
Q(\mathbb{F}_q) . 2 & \textrm{if\ $q$\ is\ odd,}\\
Q(\mathbb{F}_q) & \textrm{if\ $q$\ is\ even.}
 \end{array} \right.
\end{displaymath}
The group $\InnDiag(P\Omega^{+}(8, \mathbb{F}_q))$ of inner-diagonal automorphisms of $P\Omega^{+}(8, \mathbb{F}_q)$ is a group with triality $S$ corresponding to $\widehat{Q(q)}$, where $S$ is the group of graph automorphisms of $P\Omega^{+}(8, \mathbb{F}_q)$.  Furthermore the factor group

\begin{displaymath}
\InnDiag(P\Omega^{+}(8, \mathbb{F}_q)) / P\Omega^{+}(8, \mathbb{F}_q) = \left\{ \begin{array}{ll}
1 & \textrm{if\ $q$\ is\ even,}\\
\cong \mathbb{Z}_2 \times \mathbb{Z}_2 & \textrm{if\ $q$\ is\ odd,}
 \end{array} \right.
\end{displaymath}
and in the latter case $S$ acts nontrivially on $\mathbb{Z}_2 \times \mathbb{Z}_2$.

Whilst Doro proved
\index{Doro, S.}%
 that every Moufang loop can be obtained from a suitable group with triality, in general this group is not uniquely determined by the Moufang loop.  In~\cite{grish1}, Grishkov and Zavarnitsine describe all possible groups associated with a given Moufang loop.   They also introduce several universal groups with triality and discuss their properties.

A group $G$ with $S \leq \Aut(G)$ is \emph{S-simple} if the identity and $G$ are the only $S$-invariant normal subgroups of $G$.  The group $G$ is \emph{triality-simple}
\index{group ! triality-simple}%
 if it is $S$-simple for $S \cong \Sym(3)$ and additionally the group $G . S$ is a group with triality with respect to the conjugacy class containing the transpositions of $S$.

A \emph{locally finite field}
\index{locally finite field}%
 is one that is isomorphic to a subfield of the algebraic closure of a finite field $\mathbb{F}_p$ for some prime $p$.  Liebeck's results were extended to locally finite simple Moufang loops by J. I. Hall~\cite{hallji},
\index{Hall, J. I.}%
 as follows,
\begin{theorem}
If $G$ is a nonabelian locally finite triality-simple group, then $G . S$ is one of:
\index{group ! locally finite}%

(a)  $N \wr \Sym(3)$ for a nonabelian locally finite simple group $N$,

(b)  $P\Omega^{+}(8, K) \sd \Sym(3)$ for a locally finite field $K$. 
\end{theorem}

\begin{theorem}
A locally finite simple Moufang loop
\index{Moufang loop}%
 is either associative (and so is a locally finite simple group) or is isomorphic to a Paige loop over a locally finite field.
\end{theorem}

This result applies to both finite and infinite groups $G$.

Hall has a different but equivalent formulation of the theory of groups with triality to that of Doro, in terms of a universal triality group; in~\cite{hallji} he proves that a Moufang loop $(Q , \cdot)$ is locally finite if and only if the associated group with triality $\Aut(D(Q , \cdot))^0$ is locally finite.  All locally finite fields are countable, and a finite-dimensional (matrix) algebra over a countable field is countable.  This has the corollary~\cite{hallji2} that an uncountable locally finite simple Moufang loop is associative and so is a locally finite simple group.
\index{group ! locally finite}%
\index{group ! simple}%
 Also that the octonions and the associated Paige loop
\index{Paige loop}%
  over a locally finite field are countable.   

There are many Moufang loops related to the 8-dimensional real octonion algebra,
\index{octonions}%
 for example the multiplicative loop of non-zero elements in $\mathbb{O}$.  A nice reference is~\cite{nagyvoj1}.  We can say more.  Let $M$ be a Moufang loop with an \emph{isotopy} (with a slightly different but equivalent definition to the one above),
\index{quasigroup ! isotopy}%
 that is, a triple of maps $\alpha, \beta, \gamma : M \to M$ such that $x^{\alpha} y^{\beta} z^{\gamma} = 1$ whenever $xyz = 1$.  For example, if $u$ is an octonion of norm 1, so that $u^{-1} = \overline{u}$, and with left, right and bi-multiplications given respectively by $L_u (x) = ux$, $R_u (x) = xu$ and bi-multiplication $B_u (x) = uxu$ for every $x \in L$, then it follows from the Moufang law that $(L_u, R_u, B_{\overline{u}})$ is an isotopy.  Now if $u$ is purely imaginary of norm 1, then $B_1 B_u$ acts as minus the reflection in $u$ on purely imaginary octonions.  These generate $\SO(7)$.  Also bi-multiplications do not generally stabilize 1, so they generate at least $SO(8)$.
\index{group ! SO(8)@$\SO(8)$}%
 So the group of isotopies is approximately the spin group
\index{group ! spin}%
 (double cover of the orthogonal group) with a triality automorphism given by the map $(\alpha, \beta, \gamma) \mapsto (\beta, \gamma, \alpha)$.

Moufang proved that the loop of units in any alternative algebra satisfies (one of the equivalent forms of) the Moufang identity.  Moufang loops, like composition algebras, obey the \emph{alternative laws}
\index{alternative laws}%
\[ (xy)x = x(yx),\quad\quad x(xy) = x^2y,\quad\quad (xy)y = xy^2.\quad\quad \]
Any octonion algebra is alternating, so the units of norm 1 in the split-octonians
\index{split-octonions}%
 (over any field) form a Moufang loop. 
 
The connection between Moufang loops and Cartan triality
\index{Cartan triality}%
\index{Cartan, \'E.}%
is further illustrated by the following
results~\cite{barlottistra}:  (a) if $\mathbb{O}^{*}$ denotes the
Moufang loop of all Cayley numbers
\index{Cayley numbers}%
of norm one and $Z = \{\pm 1\}$ is
its centre, then the group generated by all right and left
translations of the factor loop $\mathbb{O}^{*} / Z$ is isomorphic to
$\PSO(8, \mathbb{R})$ and is closed in the whole homeomorphism group of
$\mathbb{O}^{*}$; (b) the group $\Delta$ of all continuous
direction-preserving collineations of $\mathbb{O}^{*} / Z$ is
isomorphic to $\PSO(8, \mathbb{R})$; the stabilizer $\Delta_L$ of
$\Delta$ on any line $L$ in $\mathcal{N}(Q)$ is isomorphic to $\PSO(7,
\mathbb{R})$.  The stabilizer in $\Delta$ of every point is isomorphic
to the compact exceptional Lie group
\index{group ! Lie}%
$G(2)$;
\index{group ! G(2)@$G(2)$}%
(c) the finite simple orthogonal group $P\Omega^{+}(8, q)$
\index{group ! orthogonal}%
 is the group with triality
\index{group ! with triality}%
associated with the finite simple nonassociative Moufang loops~\cite{grish}, or $P\Omega^{+}(8, q) \sd \Sym(3)$ according to~\cite{hallga}.

\head{Lie Algebras with Triality}

This refers to the work of Grishkov~\cite{grishkov}
\index{Grishkov, A. N.}%
 who defined a \emph{Lie algebra $\mathsf{L}$ with triality}
\index{Lie algebras with triality}%
 to be one on which $\Sym(3) = \langle \sigma, \rho : \sigma^2= 1 = \rho^3, \sigma \rho \sigma=\sigma^2 \rangle$ acts as an automorphism group
\index{group ! automorphism}%
  such that \[(x^\sigma-x) +(x^\sigma- x)^\rho+ (x^\sigma-x)^{\rho^2} =0\quad \forall x \in \mathsf{L}.\]  This is the additive version of the multiplicative triality identity for groups.
\index{triality identity}%
 Examples include $D_4$
\index{group ! D(4)@$D(4)$}%
  and any Lie algebra of an algebraic or Lie group
\index{group ! Lie}%
 with triality. For any Lie algebra $\mathsf{L}$, we can construct the Lie algebra with triality $\mathsf{T(L)} = \mathsf{L}_1 \oplus \mathsf{L}_2\oplus \mathsf{L}_3$, where $\mathsf{L}_1$, $\mathsf{L}_2$ and $\mathsf{L}_3$ are isomorphic to $\mathsf{L}$ has triality where $\sigma$ transposes $\mathsf{L}_1$ and $\mathsf{L}_2$ and stabilizes $\mathsf{L}_3$, and $\rho \in C_3.$  A Lie algebra with triality is \emph{standard} if it is isomorphic to some $\mathsf{T(L)}$ or its invariant subalgebra.  An algebra $A$ is perfect if $A^2 = A$.  Grishkov proves that a perfect finite-dimensional Lie algebra over an algebraically closed field of characteristic zero which has triality is an extension of a Lie algebra with triality of type $D_4$ by a standard Lie algebra with triality.  A Lie algebra $\mathsf{L}$ has triality if and only if $\mathsf{L} = \mathsf{L}_0\oplus \mathsf{L}_2$, where $\mathsf{L}_0 = \{x \in \mathsf{L} : x^{\lambda}, \forall \lambda \in \Sym(3)$ is $\Sym(3)$-invariant, and $\mathsf{L}_2$ is a sum of irreducible 2-dimensional $\Sym(3)$-modules.  A finite or infinite-dimensional simple Lie algebra with triality is of type $D_4$.
\index{group ! D(4)@$D(4)$}%

\section{Set Theory}
\label{SetTheory}

We omit introductory discussions of Zermelo-Fraenkel set theory, details of which can be found in standard texts such as~\cite{devlin} or~\cite{drakesi}.  We need an extension of the following:

\head{Axiom of Foundation}.  A set contains no infinitely descending membership sequence.  (Alternatively, a set contains a membership minimal element).
\index{axiom of foundation}%

There are two axioms that are featured in our work:

\head{(a) Axiom of Anti-Foundation}.
There are many versions of anti-foundation, and rather than give one of the statements or an extensive discussion, it suffices for our purposes to note that it represents a weakening of foundation by allowing the existence of set inclusions of the form $x \in \{x , y\}$ and infinite membership sequences of the form $\ldots \in x_4 \in x_3 \in x_2 \in x_1$.  Accounts of anti-founded set theory can be found in~\cite{acz} and~\cite{barmos1}. 

The set theory of the binary membership relation $\epsilon$ satisfies the Zermelo-Fraenkel axioms.  In~\cite{cameron} it was observed that first-order set theory can
be viewed as the theory of certain special orientations of the random graph.
\index{graph ! random}%
 It might be hoped that an analogous proposition would demonstrate that
first-order $ZFA$ set theory, this being the Zermelo-Fraenkel axioms
with anti-foundation
\index{axiom of anti-foundation}%
replacing foundation, is the theory of special orientations of $\mathfrak{R^{t}}$.  However as we will show in one of the chapters, there is no straightforward extrapolation and that the triality graph
\index{graph ! triality}%
 \emph{fails} to be a model of ZFA set theory.

A graph can be formed from a (vertex) set $\{x, y, \ldots \}$ by
attaching $y \to x$ with an edge if and only if $x \in y$.  

A \emph{directed}
\index{graph ! random ! directed}%
random graph must allow for four possible directed types of edges:

$$
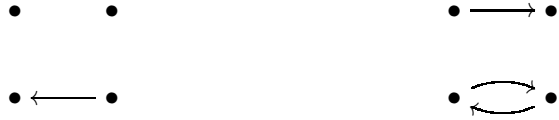
\xymatrix{{\bullet} & {\bullet} &&&& {\bullet} \ar[r] & {\bullet}}$$

$$\xymatrix{{\bullet} & {\bullet} \ar[l] &&&& {\bullet} \ar@/^/[r] & {\bullet} \ar@/^/[l]}$$
as well as vertices with loops, regarded as self-adjacencies.  The fourth type is disallowed in an \emph{oriented graph}
\index{graph ! oriented}%
which is a directed graph having no symmetric pair of directed edges.  In terms of the set-theoretic construction of $\mathfrak{R}$, symmetrising the binary set-membership relation $\in$ reduces this to \emph{two} possible vertex pairs, edges and non-edges, because again the fourth type is disallowed by the axiom of foundation.
\index{axiom of foundation}%

However,
if we allow anti-foundation then $\in$ is no longer asymmetric, and we can no
longer regard set membership as an oriented graph.  The symmetry of
$\in$ gives double-edges and its reflexivity permits vertices with
loops.  

In terms of sentences of first-order logic,
\index{first-order logic}%
we can define the triality graph
\index{graph ! triality}%
 with loops and directed edges, where the 3 adjacency types are non-edge, edge and double-edge, as the conjunction of an infinite
number of first-order sentences of the form:
$$  \forall u_1,\ldots,u_p,v_1,\ldots,v_q,w_1,\ldots,w_r\ \exists
z,z^{l} $$ 
$$ \bigcup_{\substack{1\le i\le p\\ 1\le j\le q\\
1\le k\le r}}\left( u_i \in z \wedge \lnot (z \in u_i) \right)
\vee \left(z \in  u_i \wedge \lnot (u_i \in z) \right) \wedge \left(
\lnot (v_j \in z \vee z \in v_j) \right)$$ $$\wedge$$
$$\left( u_i \in z^{l} \wedge \lnot (z^{l} \in u_i) \right)
\vee \left(z^{l} \in  u_i \wedge \lnot (u_i \in z^{l}) \right) \wedge \left(
\lnot (v_j \in z^{l} \vee z^{l} \in v_j) \right)$$
$$\wedge \left( w_k \in z \wedge z \in w_k \right) \wedge \left(
w_k \in z^{l} \wedge z^{l} \in w_k \right) \wedge \lnot (z \in z) \wedge (z^{l} \in z^{l}).$$

The superscript $l$ denotes loops.  Because in the I-property we have equal adjacency status for
vertices both with and without a loop, edge-complementation in a graph is independent of loops.  

\bigskip

\head{ (b) axiom of choice}.
\index{axiom of choice}%
If $C$ be a collection of nonempty sets, then we can choose a member from each set in that collection.  In other words, there exists a function $f$ defined on $C$ with the property that, for each set $S$ in the collection, $f(S)$ is a member of $S$.
 
In an appendix to one of the chapters, we build a model of set theory in which the axiom of choice is taken to be false, in particular leading to a modified statement of the Baer--Schreier--Ulam Theorem, as it applies in this model. 

An account of permutations and the axiom of choice is given in the article by Truss
\index{Truss, J. K.}%
in~\cite{kayem}.  The many weaker and stronger versions of the axiom of choice are discussed in the book~\cite{howardrub} by Howard and Rubin.
\index{Howard, P.}%
\index{Rubin, J. E.}%

\bigskip

\section{Number Theory}
\label{NumberTheory}

\head{Quadratic Number fields}
\index{quadratic number field}%
We give a prologue on quadratic number fields to introduce the terms that we
will need in the text.

A reformulation of the law of
quadratic reciprocity for quadratic number fields has been influential in
the development of $L$-function theory; it is given by the functional equation
\[ \zeta(s)\ L(\chi, s) = \zeta_{K}(s). \] 

The terms are as follows:

$\zeta(s) = \sum_{n \geq 1} n^{-s}$, (integers $s \geq 1$) is the
Riemann zeta function;

$L(\chi, s) = \sum_{n \geq 1} \chi(n) n^{-s}$, is the Dirichlet
$L$-function;

$\chi: \mathbb{Z} \to \mathbb{C}$ is a Dirichlet character
\index{character ! Dirichlet}%
 modulo $q
\geq 1$, i.e. a group homomorphism $(\widetilde\chi):
(\mathbb{Z}/q\mathbb{Z})^{*} \to (\mathbb{C})^{*}$, such that $\chi(x)
= (\widetilde\chi)\ (x \bmod{q})$ for $(x, q) = 1$ and $\chi(x) = 0$ if
$(x, q) \neq 1$;

$\zeta_{K}(s) = \sum_{\mathfrak{a}} (N\mathfrak{a})^{-s}$, is the
Dirichlet zeta function, (which is the only analogue of $\zeta(s)$
which satisfies the required Euler product);

$K = \mathbb{Q}[(\chi(-1) q)^\frac{1}{2}]$ is the quadratic field, so $K$ is
imaginary if $\chi(-1) = -1$ and real if $\chi(-1) = 1$;

$\mathfrak{a} \subset \mathfrak{O}_{K}$, denotes non-zero integral
ideals in the ring of integers in $K$.  The norm of the ideal is
$N\mathfrak{a}$.

In the formula, $\chi \bmod{q}$ is a non-trivial \emph{primitive}
quadratic character
\index{character ! quadratic}%
 of \emph{conductor} $q$, i.e. there does not exist
$\widetilde q\ |\ q$, $\widetilde q < q$ and a character of
$(\mathbb{Z}/ \widetilde q\ \mathbb{Z})^{*}$ such that $\chi(n) =
\widetilde\chi (n \bmod{\widetilde q})$ for $(n, \widetilde q) = 1$.  If $q$ is
prime, then any non-trivial character is primitive modulo $q$.  The
conductor identifies the `bad' primes and measures how bad they are.

The three types of primes are encapsulated in the Euler product
expression
\[ \zeta_{K}(s) = \prod_{p\ split} (1 - p^{-s})^{-2}  \prod_{p\ inert} (1
- p^{-2s})^{-1} \prod_{p\ ramified} (1 - p^{-s})^{-1} \] 
which is equivalent to the above functional equation with the characterisation
\begin{displaymath}
\left\{ \begin{array}{ll}
p\ is\ ramified & \mathit{iff\ \chi(p)=\ 0\  iff\ (p|q)} \\
p\ is\ split & \mathit{iff\ \chi(p) =\ 1\  iff\ (\frac{D}{p}) = 1} \\
p\ is\ inert & \mathit{iff\ \chi(p) = -1\  iff\ (\frac{D}{p}) = -1} \\
 \end{array} \right.
\end{displaymath}
where D is the discriminant of $K$, and the Legendre symbol $(a/p)$
has the value $1$ if $a$ is a quadratic residue $\bmod{p}$ and $-1$ if
$a$ is a quadratic nonresidue $\bmod{p}$, and zero if $p|a$.  The terms
need explanation.  If $\mathfrak{O}$ is the ring of integers of an algebraic number field
$K$, and $\mathfrak{p}$ is a non-zero prime ideal, $\mathfrak{O} /
\mathfrak{p}$ is a finite field, so it contains $p^{f}$ for some
rational prime $p$ and some $f > f_{\mathfrak{p}} > 0$.  Since
$\mathfrak{p} \supset (p)$, let $\mathfrak{p}^{e}$ be the exact power
of $\mathfrak{p}$ dividing $(p)$.  A prime with $e_{\mathfrak{p}} > 1$ is
called \emph{ramified}.

If number fields $K \supset k$ have respective prime ideals in integer
rings $\mathfrak{P} \in \mathfrak{O}$ and $\mathfrak{p} \in
\mathfrak{o}$ with $\mathfrak{P}$ a prime factor of $\mathfrak{p}\
(=\mathfrak{P} \cap k)$, the \emph{splitting group} $Z$ is the
subgroup of $\Gal(K/k)$ fixing $\mathfrak{P}$.  The kernel of the
epimorphism $Z \to
\Gal((\mathfrak{O}/\mathfrak{P})/(\mathfrak{o}/\mathfrak{p}))$ is the 
\emph{inertia group} $T \vartriangleleft Z$.

We require the following form of the Chinese Remainder Theorem~\cite[p.12]{swin}:
\index{Chinese Remainder Theorem}%
let $\mathfrak{a}_{1}, \ldots, \mathfrak{a}_{m}$ be
non-zero pairwise coprime integral ideals and let $\alpha_{1}, \ldots,
\alpha_{m} \in \mathfrak{o}$ (a ring of integers).  Then $\exists
\alpha \in \mathfrak{o}$ such that $\alpha \equiv \alpha_{\mu} \bmod{\mathfrak{a}_{\mu}} (\mu = 1, \ldots, m)$.  There is also a generalization of Dirichlet's theorem on primes in arithmetic
progression~\cite[p.94]{swin}.  If $q$ is squarefree and odd,
then $\chi(q) = \prod_{p|n} (\frac{n}{p})$ is a unique primitive character
\index{character ! primitive}%
modulo $q$.  This is clearly quadratic$\mod{q}$, and since
$p|q>2$, there are quadratic residues and non-residues$\mod{p}$, so
$\chi$ cannot be induced from $q|p$ for any $p$.  Then $\chi(-1)
\equiv q \pmod{4}$ so $K = \mathbb{Q}(\sqrt q)$ if $q \equiv 1 \pmod{4}$ and $K =
\mathbb{Q}(\sqrt -q )$ if $q \equiv 3 \pmod{4}$, with discriminant $D =
\chi(q) \cdot q$.  Taking $q$ to be prime $>2$, gives $\chi(-1) =
(-1)^{(q-1)/2}$.  Then quadratic reciprocity says that $\left(
\frac{-1}{p} \right) = 1$ iff $p \equiv 1 (4)$, that $\left(
\frac{2}{p} \right) = 1$ iff $p \equiv \pm1 (8)$ and that $\left( \frac{p}{q} \right) = \left( \frac{D}{p}\right) = \left(
\frac{(-1)^{(q-1)/2}}{p} \right) \left( \frac{q}{p} \right) =
(-1)^{(p-1)(q-1)/4} \left( \frac{q}{p} \right)$.

A similar analysis goes through if $q$ is even.

\bigskip
\bigskip

\head{Cubic Reciprocity}
\index{cubic reciprocity}%

This is a preamble outlining the theory of cubic reciprocity~\cite{ire}.
 Consider the ring $D=\mathbb{Z}[\omega]$, where
$\omega = (-1 + \sqrt{-3}) / 2$, which is a unique factorisation
domain and a Euclidean domain.
\index{Euclidean domain}%
  Define the norm of $\alpha = a + b \omega \in
\mathbb{Z}[\omega]$ by $N\alpha = \alpha \overline{\alpha} =
a^{2}-ab+b^{2}$, where $\overline{\alpha}$ means the complex conjugate of
$\alpha$.  Call the primes in $\mathbb{Z}$ rational primes and those
in $D$ simply primes (so $7 = (3 + \omega) (2 - \omega)$ is not
prime).  Say that an ideal $I$ is principal if $I = (a)$ for $a \in
I$.  Then $a$ and $b$ are associate if and only if $(a) = (b)$.
 Recall the following propositions:
\begin{proposition}
For all primes $\pi \in D$, there is a rational prime $p$ such that
\begin{displaymath}
N\pi = \left\{ \begin{array}{ll}
p & \text{if $\pi$ is not associate to a rational prime} \\
p^{2} & \text{if $\pi$ is associate to p.}
 \end{array} \right.
\end{displaymath}
\end{proposition}
\begin{proposition}
If $\pi \in D$ and $N\pi=p$, a rational prime, then $\pi$ is a prime
in $D$.
\end{proposition}
\begin{proposition}
Suppose $p$ and $q$ are rational primes.  If $q\ \equiv\ 2\ (mod\ 3)$, then
$q$ is prime in $D$.  If $p\ \equiv\ 1\ (mod\ 3)$, then $p = \pi
\overline{\pi}$, where $\pi$ is a prime in $D$.
\end{proposition}

So if $q$ is prime, then $N(q) = q^2 \equiv 1\ (3)$.  If $p = \pi
\bar{\pi}$, then $N(p) = p^2 = N(\pi) N(\overline{\pi})$ so $N(\pi) = p
\equiv 1\ (3).$

The \emph{units} $\alpha \in D$ (for which $\exists \beta \in D$ such that
$\alpha \beta = 1$, or equivalently $N\alpha = 1$) are $\pm 1,\ \pm \omega,\ \pm \omega^{2}$.  The congruence
 \index{congruence}%
 classes modulo a non-zero, nonunit $\gamma \in D$ is called the residue class
ring modulo $\gamma$, denoted $D / \gamma D$.  For a prime $\pi$, $D /
\pi D$ is a finite field with $N\pi$ elements, so the
multiplicative group $(D / \pi D)^{*}$ has order $N\pi-1$, and
since $\{1,\ \omega,\ \omega^{2}\} \cong C_{3}$, it follows that $3 |
  N\pi-1$.  Most importantly for us, if $\pi \ne 3$, then the residue
  classes of $1,\ \omega$, and $\omega^{2}$ are distinct in $D / \pi D$.

\head{Definition} If $N\pi \ne 3$, the \emph{cubic residue character}
of $\alpha$ modulo $\pi$ is given by

(a)  $\chi_{\pi}(\alpha) := (\alpha / \pi)_{3} = 0$ if $\pi|\alpha$

(b)  $\alpha^{(N\pi-1)/3} \equiv (\alpha/\pi)_{3}\ (mod\ \pi) = 1$ or
$\omega$ or $\omega^{2}\ (mod\ \pi)$.

Now $\chi_{\pi}(\alpha) = 1$ or $\omega$ or $\omega^{2}\ (mod\ \pi)$ by
definition and equals $1$ $\emph{iff}$ $\alpha$ is a cubic residue.

\head{Definition} For a prime $\pi \in D$, $\pi$ is $\emph{primary}$
if $\pi \equiv 2\ (mod\ 3)$.

Every non-zero element of $D$ has six associates.  If $N\pi = p = 1\ (mod\ 3)$, only one
of the associates of $\pi$ is primary, and we may identify the two
fields $D / \pi D$ and $\mathbb{Z}/ p \mathbb{Z}$ by mapping the coset
of $n \in \mathbb{Z}_{p}$ to the coset of $n \in D / \pi D$. 

The Law of Cubic Reciprocity says that if $\pi_{1}$ and $\pi_{2}$ are
primary, and $N \pi_{1}, N \pi_{2} \neq 3$, and $N \pi_{1} \neq N
\pi_{2}$, then 
\[  \chi_{\pi_{1}}(\pi_{2}) = \chi_{\pi_{2}}(\pi_{1}).  \]  

This says that $\pi_{1}$ is a cubic residue$\pmod{\pi_{2}}$ \emph{iff}
$\pi_{2}$ is a cubic residue$\pmod{\pi_{1}}$.


\chapter{Diagrams of Algebraic Structures}

\vspace*{1cm}

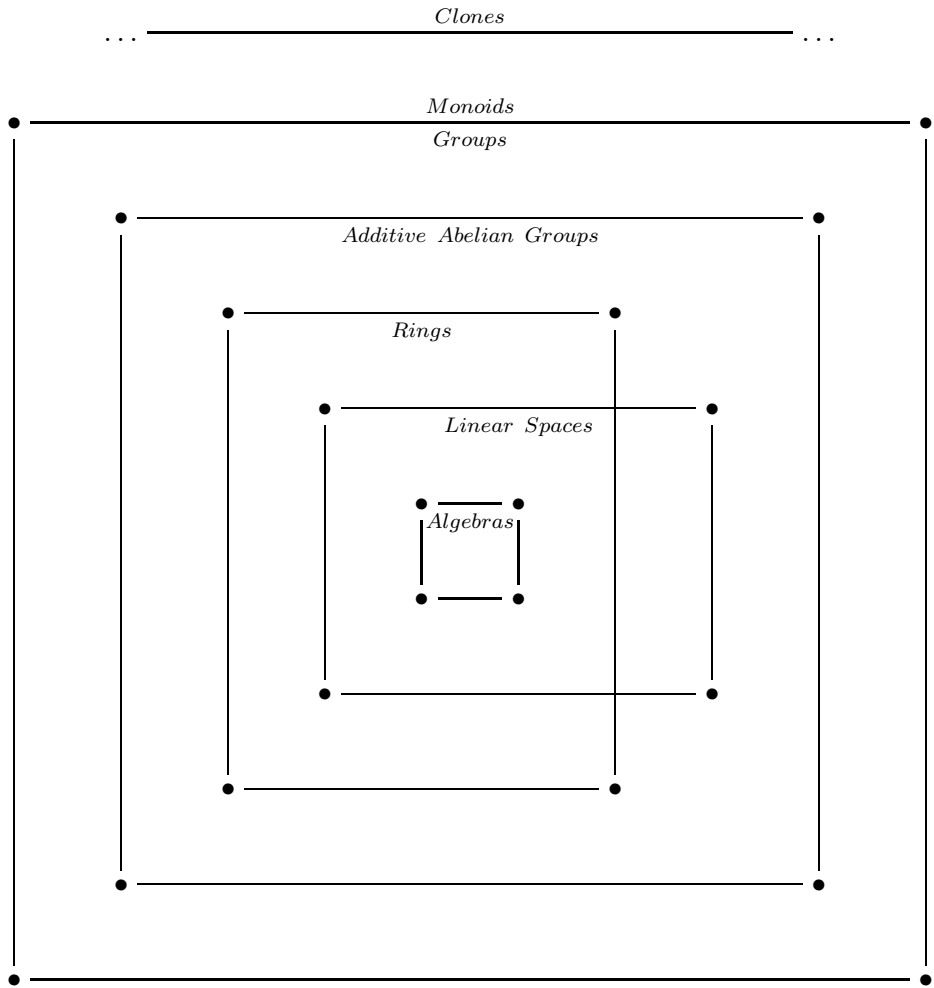
\begin{figure}[!h]
$$\xymatrix{
& {\ldots}  \ar@{-}[rrrrrrr]^{Clones} &&&&&&& {\ldots}\\
{\bullet} \ar@{-}[ddddddddd] \ar@{-}[rrrrrrrrr]_{Groups}^{Monoids} &&&&&&&&& {\bullet} \ar@{-}[ddddddddd]\\
&{\bullet} \ar@{-}[ddddddd] \ar@{-}[rrrrrrr]_{Additive\ Abelian\ Groups} &&&&&&& {\bullet} \ar@{-}[ddddddd]\\
&&{\bullet} \ar@{-}[ddddd] \ar@{-}[rrrr]_{Rings} &&&& {\bullet} \ar@{-}[ddddd]\\
&&& {\bullet} \ar@{-}[ddd] \ar@{-}[rrrr]_{Linear\ Spaces} &&&& {\bullet} \ar@{-}[ddd]\\
&&&& {\bullet} \ar@{-}[d] \ar@{-}[r]_{Algebras} & {\bullet} \ar@{-}[d]\\
&&&& {\bullet} \ar@{-}[r] & {\bullet}\\
&&& {\bullet} \ar@{-}[rrrr] &&&& {\bullet}\\
&& {\bullet} \ar@{-}[rrrr] &&&& {\bullet}\\
&{\bullet} \ar@{-}[rrrrrrr] &&&&&&& {\bullet}\\
{\bullet} \ar@{-}[rrrrrrrrr] &&&&&&&&& {\bullet}
}$$\caption{Interrelationship of Algebraic Structures}
\end{figure}

\clearpage

\vspace*{2.5cm}

\begin{figure}[!h]
$$\xymatrix{
& {Magmas}\\
{Monoids\ \&\ Semigroups}  \ar@{->}[ur] && {Loops\ \&\ Quasigroups}  \ar@{->}[ul] \\
& {Groups} \ar@{->}[ur]  \ar@{->}[ul] 
}$$\caption{Further Algebraic Interrelationships}
\end{figure}
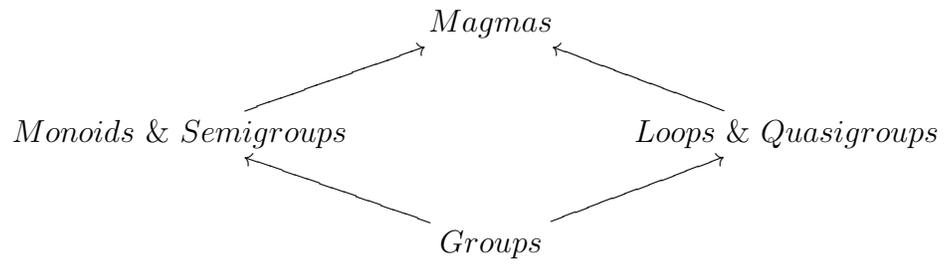

\chapter{Diagrams of Graph Types}

Let $\mathscr{X}$ be some class of graphs.  

Define a graph $\Gamma$ to be \emph{$\mathscr{X}$-homogeneous}
\index{graph ! homogeneous}%
 if any isomorphism between induced $\mathscr{X}$-subgraphs of $\Gamma$ extends to an automorphism of $\Gamma$.
 
Define a graph $\Gamma$ to be \emph{$\mathscr{X}$-transitive}
\index{graph ! transitive}%
 if whenever $A$ and $B$ are isomorphic $\mathscr{X}$-subgraphs of $\Gamma$, there is an automorphism of $\Gamma$ carrying $A$ to $B$.  Hence transitive is the weaker version of a given homogeneous property.

If $\Gamma$ is finite, define $\Gamma$ to be \emph{combinatorially $\mathscr{X}$-homogeneous}
\index{graph ! combinatorially homogeneous}%
 if whenever $A$ and $B$ are isomorphic $\mathscr{X}$-subgraphs of $\Gamma$ then $|N(A)| = |N(B)|$, where $N(x)$ is the set of neighbours of $x \in V(\Gamma)$, and we define $N(X) := \bigcap_{x \in X} N(x)$.

At the end of this section we give a more complete list of definitions.  Note that\\
$\text{$\mathscr{X}$-homogeneous} \Rightarrow \text{$\mathscr{X}$-transitive} \Rightarrow \text{combinatorially $\mathscr{X}$-homogeneous}$.

If $\mathscr{X}$ is the set of all graphs on at most $5$ vertices then for finite graphs $\text{combinatorially $\mathscr{X}$-homogeneous} \Rightarrow \text{$\mathscr{X}$-homogeneous}$.

We use \emph{$k$-homogeneous}, \emph{$k$-transitive}, and so on if $\mathscr{X}$ is the set of all graphs on at most $k$ vertices.  Usually when we use the $k$-prefix we intend a family of properties rather than a single one.

\head{Example}

If $\mathscr{X}$ is the set of connected graphs 
\index{graph ! connected}%
 and $\mathscr{X} \subseteq \mathscr{Y} \Rightarrow  \text{$\mathscr{Y}$-property} \Rightarrow  \text{$\mathscr{X}$-property}$, where `property' means homogeneous, transitive or combinatorially homogeneous.

We also have for finite graphs that
$\text{$(k+1)$-CS-homogeneous} \Rightarrow \text{$(k+1)$-CS-transitive} \Rightarrow \text{$k$-distance-transitive}$, or the stronger version\\
$\text{CS-homogeneous} \Rightarrow \text{CS-transitive} \Rightarrow \text{distance-transitive}$.
Further finite graph connections are given in the next two Figures.
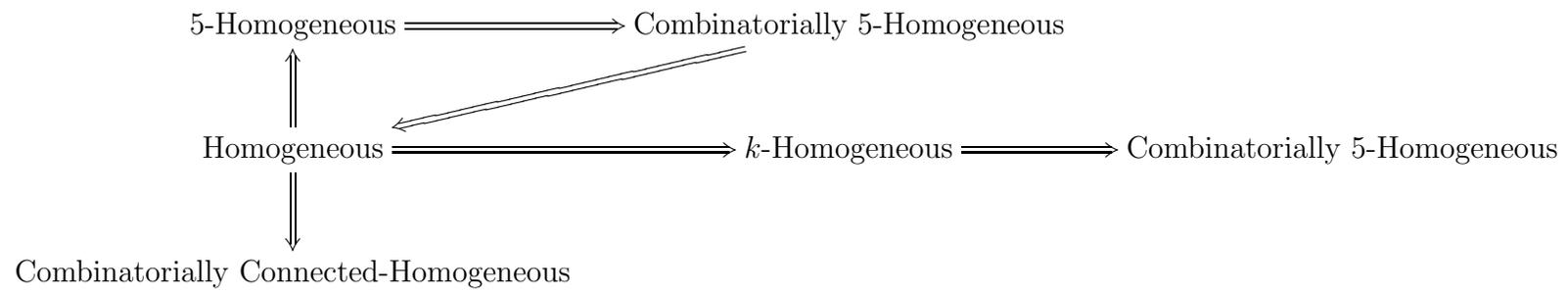
\begin{figure}[ht]
\vbox to \vsize{
 \vss
 \hbox to \hsize{
  \hss
  \rotatebox{90}{
  \vbox{
    $$
    \hss
\xymatrix@R=30pt@C=5pt{
\text{$5$-Homogeneous} \ar@{=>}[rr]  && \text{Combinatorially $5$-Homogeneous}  \ar@{=>}[dll]\\
\text{Homogeneous} \ar@{=>}[u] \ar@{=>}[d]  \ar@{=>}[rr]  && \text{$k$-Homogeneous} \ar@{=>}[rr]  && \text{Combinatorially $5$-Homogeneous}\\
\text{Combinatorially Connected-Homogeneous} 
}
    \hss
    $$
   }
  }
 \hss}
\vss}
\end{figure}

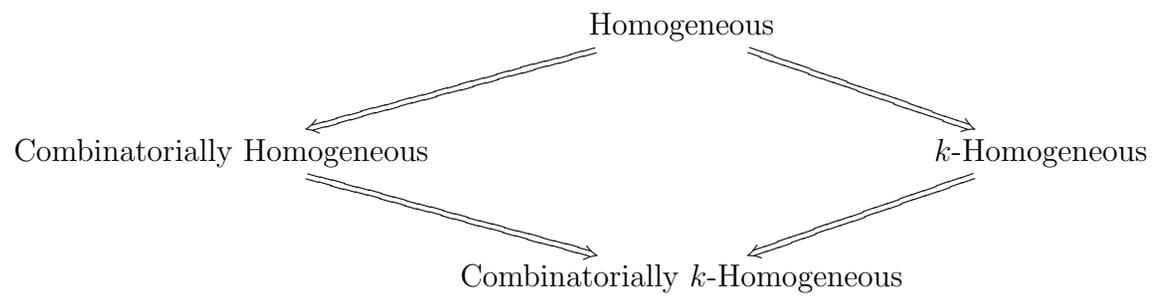
\begin{figure}[ht]
\vbox to \vsize{
 \vss
 \hbox to \hsize{
  \hss
  \rotatebox{90}{
  \vbox{
    $$
    \hss
\xymatrix@R=30pt@C=5pt{
& \text{Homogeneous}  \ar@{=>}[dl] \ar@{=>}[dr] \\
\text{Combinatorially Homogeneous}  \ar@{=>}[dr]  && \text{$k$-Homogeneous} \ar@{=>}[dl] \\
& \text{Combinatorially $k$-Homogeneous} 
}
    \hss
    $$
   }
  }
 \hss}
\vss}
\end{figure}

\clearpage

For infinite graphs we have the following:\\

$\text{homogeneous} \Rightarrow {\ldots} \Rightarrow \text{$k$-homogeneous} \Rightarrow \text{$(k-1)$-homogeneous} \Rightarrow {\ldots} \Rightarrow \text{$2$-homogeneous} \Rightarrow \text{$1$-homogeneous (= vertex-transitive)}$.  

\bigskip
\bigskip

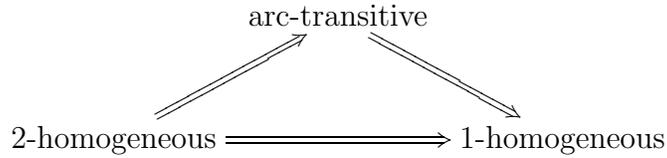
\begin{figure}[ht]    
$$
\xymatrix@R=30pt@C=5pt{
& \text{arc-transitive}  \ar@{=>}[dr]\\
\text{$2$-homogeneous}  \ar@{=>}[ur]  \ar@{=>}[rr]  && \text{$1$-homogeneous}
}
$$
\end{figure} 

\bigskip
\bigskip

A distance-regular graph has neither a symmetry nor a distance invariance, but rather a counting condition.  Hence,

\bigskip

\begin{figure}[!ht]
$$
\xymatrix@R=30pt@C=5pt{
&&& \text{distance regular} \ar@{=>}[dll] \ar@{=>}[dr] \\
& \text{distance degree regular} &&& \text{$F$-geodesic}
}  
$$
        \caption{Inclusion relations for some countable graph classes with a counting condition} 
        \label{IR2} 
\end{figure}
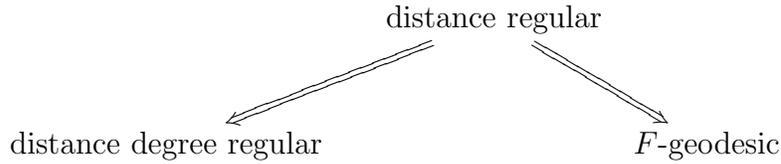  

\bigskip
\bigskip

Some of the relations whose definition we have given are more important than others.  Among the more obscure are:-

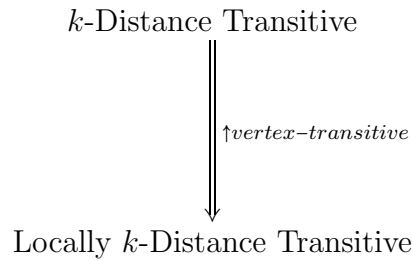
\begin{figure}[ht]
$$
\xymatrix@R=30pt@C=5pt{
&&& \text{$k$-Distance Transitive}  \ar@{=>}[dd]^{\uparrow vertex-transitive}\\
&\\
&&& \text{Locally $k$-Distance Transitive} 
}
$$
\end{figure} 

\clearpage


\begin{figure}[ht]
\vbox to \vsize{
 \vss
 \hbox to \hsize{
  \hss
  \vbox{
    $$
    \hss
\xymatrix@R=30pt@C=5pt{
&&& \text{Homogeneous (II)} \ar@{=>}[d] \ar@{=>}[ddr]  \ar@{=>}[ddll] \ar@{=>}[r]  & \text{$XY$-Transitive}  \\
&   \text{Edge-Transitive} && \text{Arc-Transitive} \ar@{=>}[d]  \ar@{=>}[ll] \ar@{=>}[r] &  \text{Vertex-Transitive} \ar@{=>}[d] \\
&  \text{Metrically Homogeneous} \ar@{=>}[d] && \text{Connected-Homogeneous} \ar@{=>}[d] & \text{$k$-arc Transitive} \ar@{=>}[d]_{k = 2} \\
& \text{Metrically $k$-Homogeneous} \ar@{=>}[d] \ar@{=>}[rr]^{\quad k = 2} && \text{Distance-Transitive}  \ar@{<=}[dr]   & \text{Triangle-free or Complete} \\
& \text{Locally Homogeneous}  \ar@{=>}[d] && \text{} \ar@{=>}[u] & \text{$3$-CS-Homogeneous} \\
&  \text{Locally P-Homogeneous} && \text{}\\
&&& \text{$k$-CS Transitive}  \ar@{<=}[d]  \ar@{=>}[uuu] \\
&&& \text{$k$-CS-Homogeneous}\ \text{$(k\geq2)$}  \ar@{=>}[uuur]_{k = 3}   \ar@{<=}[d] \\
&&& \text{$k$-Homogeneous}
}
    \hss
    $$
        \caption{Inclusion relations for some countable graph classes with symmetry or distance invariance showing the centrality of Homogeneity} 
        \label{IR1} 
   }
 \hss}
\vss}
\end{figure}
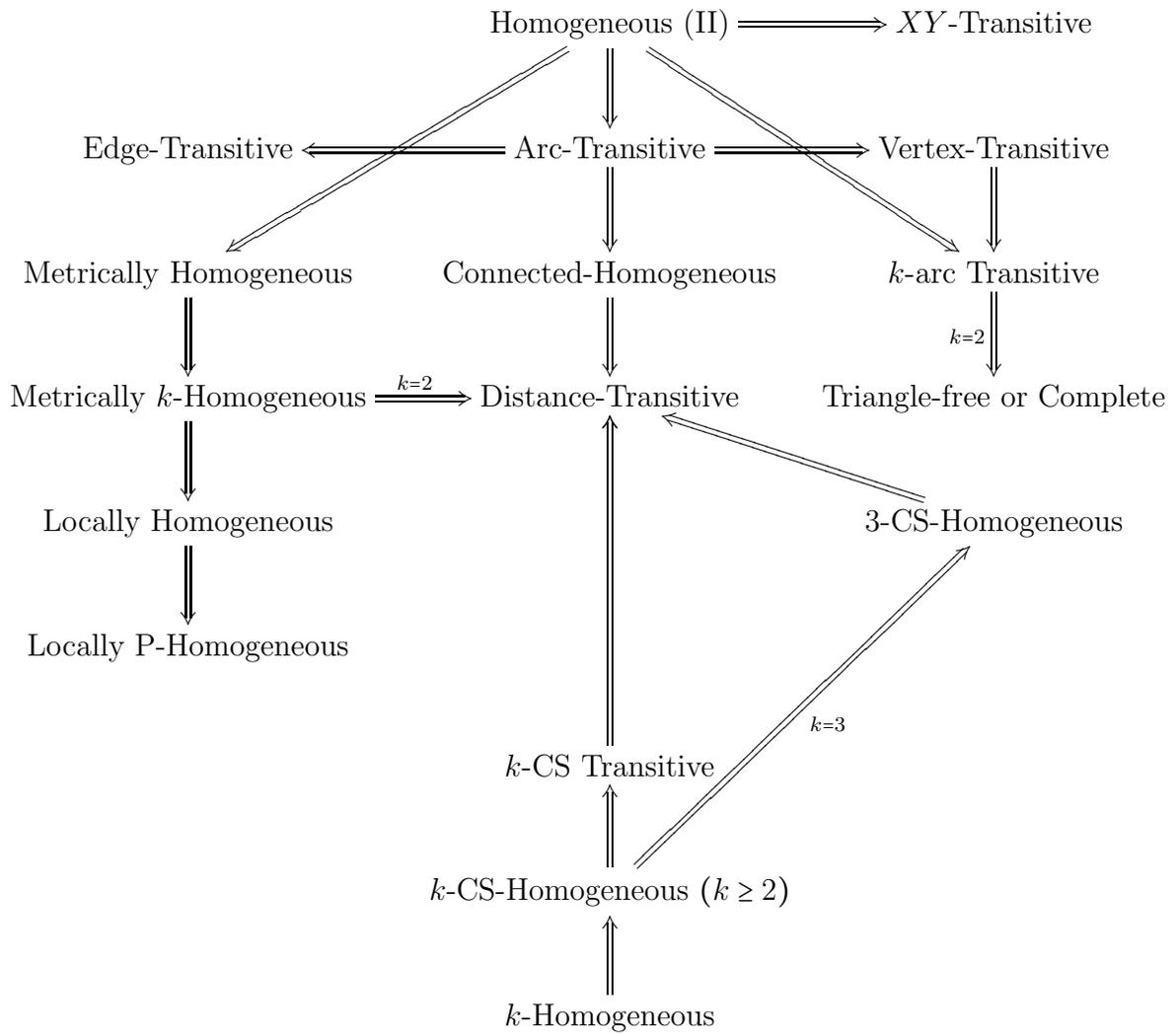

\clearpage

\textbf{Definitions}

In Figure~\ref{IR1}, transitive means transitive on tuples.  Also homogeneous means set-transitive; sometimes this is referred to as weakly homogeneous.
\begin{itemize}
 \item[(a)]  A graph is \emph{homogeneous}
\index{graph ! homogeneous}%
(respectively \emph{$k$-homogeneous}) if every isomorphism between finite substructures (respectively on at most $k$ vertices) extends to an automorphism of the entire structure; see Appendix~\ref{ModelTheory}.  

There is a different concept that was proved to be equivalent to this in~\cite{enomoto} for finite, simple, undirected graphs; it is defined as follows.  For a subset $X \subset V(\Gamma)$, and set of neighbours $N(x)$ of $x \in V(\Gamma)$, define $N(X) = \bigcap_{x \in X} N(x)$.   A finite graph is said to be \emph{combinatorially homogeneous}
\index{graph ! combinatorially homogeneous}%
 if $|N(X)| = |N(X')|$ for any isomorphic vertex-subgraphs on $X(\Gamma)$ and $X'(\Gamma)$.
 \item[(b)]  A graph is \emph{combinatorially connected-homogeneous}
\index{graph ! combinatorially connected-homogeneous}%
if whenever two subgraphs $X(\Gamma)$ and $X'(\Gamma)$ are connected and isomorphic, then for any subset $Y \backslash X$ there is an isomorphism $\alpha$ of $X(\Gamma)$ onto $X'(\Gamma)$ such that $|N(Y)| = |N(Y^{\alpha})|$.  
\item[(c)]  A graph is \emph{metrically homogeneous}
\index{graph ! metrically homogeneous}%
 if any isometry between finite induced subgraphs extends to an automorphism of the graph;
\item[(d)]  A graph is \emph{metrically $k$-homogeneous}
\index{graph ! metrically $k$-homogeneous}%
 if any isometry between $k$-element induced subgraphs extends to an automorphism of the graph;
\item[(e)]  A graph is \emph{locally homogeneous}
\index{graph ! locally homogeneous}%
 if any automorphism of induced subgraphs of the graph extends to an automorphism of the graph.  A variation is that a graph is  \emph{locally P-homogeneous}
\index{graph ! locally P-homogeneous}%
 if any automorphism of induced subgraphs of the graph that lie in P extends to an automorphism of the graph; \item[(f)]  A graph is \emph{$XY$-transitive}
\index{graph ! XY-property}%
 if every $x$-morphism between finite subgraphs $X$ and $Y$ can be extended to a $y$-morphism from the graph to itself, where $(X, x)$ and $(Y, y)$ can be (I, iso), (M, mono), or (H, homo).  There are six properties of this kind that can be considered: HH, MH, IH, MM, IM, and II.  Property II is equivalent to the standard notion of homogeneity. 
\item[(g)]  A graph is called \emph{$k$-arc-transitive} 
\index{graph ! karctransitive@$k$-arc-transitive}%
for a positive integer $k$, if it has a walk of length $k$ with specified initial point in which no line succeeds itself, and if there is always a graph automorphism of  sending each $k$-arc onto any other $k$-arc.  Note that $1$-arc-transitive is equivalent to \emph{symmetric};
\index{graph ! symmetric}%
 a graph is symmetric if, for all vertices $u, v, x, y$ of $\Gamma$, such that $u \sim v$ and $x \sim y$ are adjacent, there exists $g \in \Aut(\Gamma)$ satisfying $g(u) = x$ and $g(v) = y$.  Note also that
 
$\text{distance-transitive} \Rightarrow \text{symmetric} \Rightarrow \text{vertex-transitive.}$

Furthermore $2$-arc-transitive implies triangle-free.
\index{graph ! triangle-free}%

\item[(h)]  A graph is \emph{connected-homogeneous}
\index{graph ! connected-homogeneous}%
 if any isomorphism between connected finite induced subgraphs extends to an automorphism of the graph;
\item[(i)]  A graph is \emph{distance transitive} 
\index{graph ! distance-transitive}%
 if, given any vertices $u, v, u', v'$ such that $d(u, v) = d(u', v')$, there is an $\alpha \in \Aut(\Gamma)$ such that $\alpha (u) = u'$ and $\alpha (v) = v'$.  Note that distance transitive is equivalent to metrically $2$-homogeneous and a homogeneous graph is a distance-transitive graph with distance at most 2;
\item[(j)]  A graph is \emph{locally $k$-distance transitive} 
\index{graph ! distance-transitive ! locally}%
 if the stabilizer of a vertex $v$ acts transitively on the set of vertices at distance $1$ from $v$ and on the set of vertices at distance $k$ from $v$~\cite{devillers}.  If the graph is also vertex-transitive then it is said to be $k$-distance transitive;
\item[(k)] A graph is \emph{$k$-CS-transitive},
\index{graph ! kcstransitive@$k$-CS-transitive}%
for a positive integer $k$, if for any two connected isomorphic induced subgraphs $A, B$ of $\Gamma$, each of size $k$, there is an automorphism of $\Gamma$ taking $A$ to $B$;
\item[(l)]  A graph is called \emph{$k$-CS-homogeneous} 
\index{graph ! kcshomogeneous@$k$-CS-homogeneous}%
if any isomorphism between two connected induced subgraphs of size $k$ extends to an automorphism;
\item[(m)]  A regular connected graph is \emph{distance-regular} 
\index{graph ! distance-regular}%
if, (a) for each vertex $u$ the distance partition $\partial_u$ of $\Gamma$ relative to $u$ is equitable, and (b) the quotient $\Gamma \backslash \partial_u$ is independent of $u$;
\item[(n)]  A graph with diameter $d$ is \emph{distance degree regular} 
\index{graph ! degree-regular}%
if, the number of vertices at a distance $i$ from a given vertex depends only on $i$;
\item[(o)]  A connected graph with diameter $d$ and function $F : \{1, \ldots, d \}$ $\mapsto \mathbb{N}$ is \emph{$F$-geodesic} 
\index{graph ! Fgeodesic@$F$-geodesic}%
if, given any pair of vertices at a distance $i$ from each other have exactly $F(i)$ shortest paths joining them.  If $F(i) = 1$ for all $i$, then $\Gamma$ is called \emph{geodetic}.
\index{graph ! geodetic}%
\end{itemize}

\chapter{Some Separation Axioms for Topological Spaces}
\index{topological space ! $R_i$}%
\index{topological space ! $T_i$}%

Two subsets (or points) in a topological space $X$ are \emph{topologically distinguishable} 
\index{sets ! topologically distinguishable}%
 if at least one of them has a neighbourhood that is not a neighbourhood of the other. They are \emph{separated}
\index{sets ! separated}%
 if each of them has a neighbourhood that is not a neighbourhood of the other; that is, neither belongs to the closure of the other.

We list a few of the restrictions that can be imposed on a topological space in order of increasing strength, though there are varying definitions in the literature.

$\bullet$\ A topological space is a \emph{$T_0$ space} if one point out of any two distinct points has a neighbourhood which does not contain the other point, that is they are topologically distinguishable.

$\bullet$\ A topological space is an \emph{$R_0$ space} if any two topologically distinguishable points are separated.

$\bullet$\ A topological space is a \emph{$T_1$ space} if any two distinct points are separated.  A space is $T_1$ if and only if it is both $T_0$ and $R_0$. 

$\bullet$\ A topological space is an \emph{$R_1$ space} if any two topologically distinguishable points are separated by neighbourhoods, that is have disjoint neighbourhoods. An $R_1$ space must also be $R_0$.

$\bullet$\ A topological space is a \emph{Hausdorff space}
\index{Hausdorff space}%
 (or a \emph{$T_2$ space}) if any two distinct points have disjoint neighbourhoods which do not contain the other point.  A space is $T_2$ if and only if it is both $T_0$ and $R_1$. 

$\bullet$\ A space is called a \emph{regular space} if whenever $S$ is a closed set
\index{closed set}%
 and $x$ is a point which is not in $S$ then there are disjoint open sets
\index{open set}%
 $U$ and $V$ such that $x \in U$ and $S \subset V$, that is they are separated by neighbourhoods.   (There are also closed neighbourhoods that separate.) A regular $T_0$ space is called a \emph{$T_3$ space}.

$\bullet$\ A space is called a \emph{normal space} if any two disjoint closed subsets of X are separated by neighbourhoods.  A normal $T_1$-space is called a \emph{$T_4$ space}.

For a topological space we have the implications
\[ T_4 \Rightarrow T_3 \Rightarrow T_2 \Rightarrow T_1 \Rightarrow T_0.\]
The $T_0$ axiom can not only be added to a property, but also subtracted from a property.  An example of a regular space that is neither $T_1$ nor $T_2$ is the indiscrete topology $(X, \{\emptyset, X\})$.  A regular space where all singletons are closed is $T_2$.

\chapter{Tables}



\begin{figure}[ht]
\begin{tabular}{|l||c|c|c|c|}
\hline
\quad\quad\quad\quad\quad\quad Groups&Notation&Chapter\\ 
&&No.\\
\hline
\hline
1. Random Bipartite Graph&&\\
Groups with vertex colours&$G(\mathfrak{B}^{v})$&\ref{chap2}\\
\hline
2. Multicoloured Switching Groups&$S_{m,n}$&\ref{chap2}\\
and Duality Groups (and Variations)&$D_{m,n}$&\\
\hline
3. Finitary Switching Group&$S^{f}_{m,\omega}$&\ref{chap3}\\
\hline
4. Enhanced Switching Group&$ES_{m, \omega}$&\ref{chap3}\\
\hline
5. Closure of Switching Group&$S_{m, \omega}^{cl}$&\ref{chap3}\\
\hline
6. Coloured Tsaranov Group&$\Ts(\mathfrak{R}_{m,n})$&\ref{chap4}\\
\hline
7. Coloured Coxeter Group&$\Cox(\mathfrak{R}_{m,n})$&\ref{chap4}\\
\hline
8. Group of Equivalence Classes&$\NS(\mathfrak{R}_{m,n})$&\ref{chap7}\\
of Equivalence Relation on $\NB(\mathfrak{R}_{m,n})$&&\\
\hline
9. Group of Equivalence Classes&$\NE(\mathfrak{R}_{m,n})$&\ref{chap7}\\ 
of Near-Automorphisms $\NA(\mathfrak{R}_{m,n})$&&\\
\hline
10. Zero Vertex Index Group&$\Aut_0(\mathfrak{R}_{m,n})$&\ref{chap7}\\
\hline
11. Finite Vertex Index Group&$\Aut_{fin}(\mathfrak{R}_{m,n})$&\ref{chap7}\\
\hline
12. Automorphism Group of&$\Aut(\mathcal{F}_{\mathfrak{R}})$&\ref{chap8}\\
Filters on $\mathfrak{R}$&&\\
\hline
13. Group of Auto-homeomorphisms&$\Aut(\mathcal{T})$&\ref{chap8}\\
of $\mathcal{T}$ Topology&&\\
\hline
14. Group of Auto-homeomorphisms&$\Aut(\mathcal{T}^*)$&\ref{chap8}\\
of $\mathcal{T}^*$ Topology&&\\
\hline
15. Group of Auto-homeomorphisms&$\Aut(\mathcal{T}^{\dagger})$&\ref{chap8}\\
of $\mathcal{T}^{\dagger}$ Bipartite Graph Topology&&\\
\hline
16. Group of Automorphisms&$\Aut(\mathfrak{RHyp})$&\ref{chap8}\\
of $\mathfrak{R}$-uniform hypergraph $\mathfrak{RHyp}$ &&\\
\hline
17. Group of Almost Automorphisms&$\Aut^{*}(\mathfrak{RHyp})$&\ref{chap8}\\
of $\mathfrak{R}$-uniform hypergraph $\mathfrak{RHyp}$ &&\\
\hline
\end{tabular}
   \caption{Table of Some Groups Arising}  
\end{figure}

\begin{figure}[ht]
\vbox to \vsize{
 \vss
 \hbox to \hsize{
  \hss
  \rotatebox{90}{
   \vbox{
    \hss
\[
\begin{array}{|c|c|}
\hline
S_{2,n} \cong (C_2)^{n-1}&S_{m,n} \cong (\Alt(m))^{n(n-1)/2}\sd (C_2)^{n-1}\\
\hline
\SAut(\mathfrak{R})\ is\ a\ reduct&\SAut(\mathfrak{R}_{m,\omega})\ are\ not\ reducts\\
\hline
S_{2,n}\ has\ many\ switching\ classes&S_{m,n}\ is\
transitive\ on\ \mathcal{G}_{m,n}\\
\hline
\SAut(\mathfrak{R})\ reducts\ are\
$2$-transitive&\SAut(\mathfrak{R}_{m,n})\ reducts\ are\ highly\ transitive\\
\hline
G^W_{2,n}\ is\ not\ even\ transitive\ on\ \mathcal{G}_{2,n}&G^W_{m,n}\
acts\ primitively\ on\ \mathcal{G}_{m,n}\\
\hline
S^*_{2,n}\  is\ not\ even\ transitive\ on\ \mathcal{G}_{2,n}&S^*_{m,n}\
acts\ primitively\ on\ \mathcal{G}_{m,n}\\
\hline
S_{2,n}\ has\ no\ transitive\ extension&S_{m,n}\ weakly\ transitively\ extends\\
\hline
\Aut(\mathfrak{R})\ is\ three-star\
transitive&\Aut(\mathfrak{R}_{m,\omega})\ is\ not\ three-star\ transitive\\
\hline
\SAut(\mathfrak{R})\ is\ generously\
$2$-transitive&\SAut(\mathfrak{R}_{m,\omega})\ is\ not\ generously\ $2$-transitive\\
\hline
\SAut(\mathfrak{R})\ is\
$2$-transitive&\SAut(\mathfrak{R}_{m,\omega})\ are\ highly\ transitive\\
\hline
There\ is\ no\ polynomial&\\algebra\ Galois\ Correspondence&H \le \Sym(m-1) \leftrightarrow \mathcal{A}^{D^H_{m,\omega}}\\
\hline
\mathcal{A}(\mathfrak{R}_{2, \omega})^{H_{r}}\ is\ polynomial\
algebra&\mathcal{A}(\mathfrak{R}_{m, \omega})^{H_{r}}\ not\ polynomial\
algebra\\
\hline
The\ proportion\ of\ switching&The\ proportion\ of\ switching\\
classes\ of\ 2-colour\ graphs&classes\ of\ m-colour\ graphs\\
with\ n\ vertices\ which\ have\ non-trivial&with\ n\ vertices\ which\
have\ non-trivial\\
automorphisms\ \to\ 0\ as\ n\ \to \infty&automorphisms\ \to \infty\
as\ n\ to\ \infty\\
\hline
\end{array}
\]
    \hss    
\caption{Table of Differences between $\mathfrak{R}$ and $\mathfrak{R}_{m,\omega}\ (m \ge 3)$}
        }
                }
  \hss          }
 \vss           }
\end{figure}

\clearpage

\begin{figure}[ht]
\vbox to \vsize{
  \vss
  \hbox to \hsize{
   \hss
    \vbox{
     \hss
\[
\begin{tabular}{|c|c|}
\hline
$\Sym(n)$&$B_n$\\
\hline
$s_i^2 = (i\ i+1)^2 = 1 (i=1, \ldots, n-1)$&$b_i\ (i=1, \ldots, 
n-2)$\\
\hline
$s_i s_j = s_j s_i$\ if $|i-j| \ge 2$&$b_i b_j = b_j b_i$\ if 
$|i-j| \ge 
2$\\
\hline
$s_i s_{i+1} s_i = s_{i+1} s_i s_{i+1}$&$b_i b_{i+1} b_i = b_{i+1} 
b_i 
b_{i+1}$\\
\hline
\hline
$\Hyp(n)$&$S_{m,n}$\\
\hline
$h_i^{2} = 1\ (i=0, \ldots, n-1) $&$\sigma_{i,j}^2=1 (1 \le i \ne j \le m)$\\
\hline
$h_i h_j = h_j h_i\ \text{if}\ | i - j | > 1$&$\sigma_{i,j}\sigma_{i,k}=\sigma_{j,k}=\sigma_{i,k}\sigma_{i,j}$\\
\hline
$h_i h_{i+1} h_i = h_{i+1} h_i h_{i+1}$\ $(i \ge 1)$&$\sigma_{i,j}\sigma_{k,l}\sigma_{i,j}=\sigma_{k,l}\sigma_{i,j}\sigma_{k,l}$\\ 
$(i \ge 1)$&$|i-j|,|k-l|  \ge 2$\\
\hline
$h_1 h_0 h_1 h_0 = h_0 h_1 h_0 h_1$&\\
\hline
\end{tabular}
\]
     \hss
\caption{Presentations of $\Sym(n),\ B_n,\ \Hyp(n),\ S_{m,n}$}
         }
   \hss          }
  \vss           }
\end{figure}


\clearpage

\chapter{Quotations}

In this section we give some quotations that underlie the philosophy behind the Further Directions chapter.

\medskip

If man were restricted to collecting facts the sciences were only a sterile nomenclature and 
he would never have known the great
       laws of nature. It is in comparing the phenomena with each
       other, in seeking to grasp their relationships, that he is led
       to discover these laws... 
\begin{flushright}
Pierre-Simon Laplace
\end{flushright}

\medskip

He who wants to unlock secrets should not lock himself away in one
area of science, but should maintain connections with its other areas
as well.
\begin{flushright}
J. Hadamard
\end{flushright}

\medskip

the mathematical facts worthy of being studied are those which, by
their analogy with other facts, are capable of leading us to the
knowledge of a mathematical law, just as experimental facts lead us to
the knowledge of a physical law.  They are those which reveal to us
unsuspected kinship between other facts, long known, but wrongly
believed to be strangers to one another.
\begin{flushright}
H. Poincar\'e
\end{flushright}

In proportion as science develops, its total comprehension becomes
more difficult; then we seek to cut it in pieces and to be satisfied
with one of these pieces: in a word, to specialize.  If we went on in
this way, it would be a grievous obstacle to the progress of science.
As we have said, it is by unexpected union between its diverse parts
that it progresses.
\begin{flushright}
H. Poincar\'e
\end{flushright}


Mathematics, physics, chemistry, astronomy, march in one front. Whichever lags behind is drawn after. 
\begin{flushright}
K. Schwarzschild
\end{flushright}


As every mathematician knows, nothing is more fruitful than these
obscure analogies, these indistinct reflections of one theory into
another, these furtive caresses, these inexplicable disagreements; also nothing gives the researcher greater pleasure. 
\begin{flushright}
A. Weil
\end{flushright}
\smallskip

I think chance is a more fundamental conception than causality
\begin{flushright}
Max Born, \emph{Natural Philosophy of Cause and Chance}, (1948).
\end{flushright}

\smallskip

I don't want to belong to any club that will accept me as a member.
\begin{flushright}
Groucho Marx, \emph{Groucho and me. An autobiography.}  (1959), 320.
\end{flushright}

\smallskip
Technical skill is mastery of complexity while creativity is mastery of simplicity.
Catastrophe Theory, 1977. 
\begin{flushright}
E. C. Zeeman, \emph{Catastrophe Theory}, (1977).
\end{flushright}

\smallskip

-- that's the beauty of pure mathematics -- one is always discovering new
connections between things that looked unconnected.
\begin{flushright}
F. J. Dyson
\end{flushright}

\smallskip

The same pathological structures that the mathematicians invented to break loose from $19$th--century naturalism turn out to be inherent in familiar objects all around us in Nature.
\begin{flushright}
F. J. Dyson, \emph{Characterizing Irregularity}, Science 200 (1978), 677-8.
\end{flushright}

\smallskip

Putting together widely differing facts drawn either from the
experimental sciences or from within mathematics itself is one of the
essential ingredients of mathematics.  We must have people who try to
connect up different parts of mathematics as well as those who
restrict themselves to one area and try to get as far as possible in
that direction.
\begin{flushright}
M. Atiyah
\end{flushright}

\smallskip

The experience of past centuries shows that the development of
mathematics was due not to technical progress (consuming most of the
efforts of mathematicians at any given moment), but rather to
discoveries of unexpected interrelations between different domains
(which were made possible by these efforts).  $\ldots$ Growing specialisation and bureaucratic subdivision of mathematics
into small domains becomes an obstacle to its development $\ldots$
This lack of understanding of the interrelations between different
domains of mathematics originates from the disasterous divorce of
mathematics from physics in the middle of the 20th century, and from
the resulting degeometrisation of mathematical education.
\begin{flushright}
V. I. Arnold
\end{flushright}

\medskip

But, to be sure, this tendency is a manifestation of the unity of mathematics and, equally evidently, its emergence was a direct or indirect effect of theoretical physics.  In fact, in the second half of the XXth Century it began to swing sideways \emph{away from analysis and the linear equations of mathematial physics towards algebra (including algebraic geometry) and towards probability}, whereupon the special role of probability theoretical constructions in modern theoretical physics, and, above all, in statistical physics and quantum-field theory became undeniable.  $\ldots$ We have in mind the fact that when one studies, for example, symmetric groups of large degree or Lie groups of large rank, one encounters the same phenomenon as the Law of Large Numbers
\index{Law of Large Numbers}%
in probability theory, although, a priori, the formulation of the problem did not involve any probability!  $\ldots$ Here, the similarity to the Law of Large Numbers is not an approximate analogy, but a completely precise fact[.]  $\ldots$ the study of random algebraic systems are of particular interest.   $\ldots$ General trends in the development of mathematics are determined not only, and even not so often, by major problems remaining from the past; it has become commonplace to speak about the role of Hilbert's problems in the development of XXth Century mathematics.  Internal trends  play a much more substantial role, together with the internal logic of the development of the science, which is determined not only by unsolved problems, however attractive they might appear, but by the interweaving of the interests, methods and applications of different areas.  However, to attempt to predict these trends at a remote date is almost hopeless.
\begin{flushright}
A. M. Vershik
\end{flushright}

\medskip

It is also vital to always keep moving. The risk otherwise is to confine oneself in a relatively small area of extreme technical specialization, thus shrinking one's perception of the mathematical world and its bewildering diversity. $\ldots$ In other words there is just "one" mathematical world, whose exploration is the task of all mathematicians, and they are all in the same boat somehow. 
\begin{flushright}
A. Connes 
\end{flushright}

\medskip

The thing that \emph{doesn't} fit, is the thing that is most interesting.  The part that \emph{doesn't} go according to what you would expect.
\begin{flushright}
R. Feynman
\end{flushright}


Be prepared to consider crazy ideas, radical ideas $\ldots$ it doesn't follow that every crazy idea is a good one, but some crazy ideas are good ones and you want to think carefully before discarding them entirely.
\begin{flushright}
M. Atiyah
\end{flushright}

\clearpage

\backmatter

\chapter*{Glossary of Notation}

The page or chapter number refers to the first (or only) occurrence of the notation.  We often exclude notation that is used only within a short proof, example or paragraph.  Greek letter entries are at the end.

\begin{displaymath}
\begin{array}{lllllll}
\text{$a(n)$} &&& \text{exponential generating function}\hfill 
\mbox{ Appendix~\ref{EnumerationandReconstruction}}\\
a_n &&& \text{number of labeled structures}\\
\text{} &&& \text{on an $n$-element set}\hfill \mbox{ Appendix~\ref{EnumerationandReconstruction}}\\
\mathcal{A(M)} &&& \text{age of a relational structure}\ \mathcal{M}\hfill \mbox{Chapter~\ref{ringchap}}\\
\text{Age(M)} &&& \text{age of a relational structure M}\hfill \mbox{ Appendix~\ref{TheoryofRelationalStructures}}\\
\mathcal{A} &&& \text{polynomial algebra}\hfill \mbox{ Appendix A, Chapter~\ref{ringchap}}\\
\mathcal{A}^{G} &&& \text{algebra of $\mathcal{A}$-fixed points}\hfill \mbox{ Appendix A, Chapter~\ref{ringchap}}\\
\text{$\Alt(\Omega)$} &&& \text{alternating group of}\\
\text{} &&& \text{permutations on the set $\Omega$}\hfill  \mbox{ throughout}\\
\text{$\Aut$} &&& \text{automorphism group} \hfill \mbox{Page~\pageref{Aut}}\\
\text{$\AAut$} &&& \text{group of almost automorphisms} \hfill \mbox{ throughout}\\
\text{$\Aut^{H}(\mathfrak{R}_{m,\omega})$} &&& \text{group of vertex-permutations}\\
\text{} &&& \text{inducing colour-permutations $H$}\hfill \mbox{Page~\pageref{auth}}\\
\text{$\Aut(\mathfrak{R^{t}_{\mathfrak{i}}})$} &&& \text{automorphism group of $\mathfrak{R^{t}}$ but}\\ \text{} &&& \text{colourblind in $\mathfrak{i}$ and another colour}\hfill  \mbox{Page~\pageref{autri}}\\
\text{$\Aut(\mathbb{Q}, \le)$} &&& \text{group of all order preserving}\\ 
\text{} &&& \text{automorphisms of $\mathbb{Q}$} \hfill \mbox{Page~\pageref{autqle}}\\
\text{$\Aut^{*}(\mathfrak{R}_{m,\omega})$} &&& \text{group of all automorphisms}\\ 
\text{} &&& \text{of $\mathfrak{R}_{m,\omega}$} \hfill \mbox{Page~\pageref{allaut}}\\
\text{$\Aut_0$} &&& \text{zero vertex index group}\hfill   \mbox{Page~\pageref{zerovx}}\\
\text{$\Aut_{fin}$} &&& \text{finite vertex index group}\hfill   \mbox{Page~\pageref{zerovx}}\\
AGL(V) &&& \text{affine general linear group} \hfill \mbox{Page~\pageref{AGL(V)}}\\
b(n) &&& \text{generating function}\hfill \mbox{ Appendix~\ref{EnumerationandReconstruction}}\\
b_n &&& \text{number of isomorphism classes}\\
\text{} &&& \text{of structures of size $n$}\hfill \mbox{ Appendix~\ref{FurtherDetails}}\\
\mathcal{B} &&& \text{base of open subgroups of group}\hfill \mbox{ Appendix~\ref{TopologyinPermutationGroups}}\\
\mathcal{B}(X) &&& \text{Borel subsets of $X$}\hfill \mbox{ Appendix~\ref{CategoryandMeasure}}\\
B(x, r)  &&& \text{open ball of radius $r$ centred at $x$}\hfill  \mbox{ Appendix~\ref{TopologyinPermutationGroups}}\\
\mathbb{B} &&& \text{Boolean algebra}\hfill \mbox{Page~\pageref{boolalg}}\\
\end{array} 
\end{displaymath}

\clearpage

\begin{displaymath}
\begin{array}{lllllll}
\text{$\BA(\Gamma)$} &&& \text{group of biggest automorphisms of $\Gamma$} \hfill \mbox{Page~\pageref{dsbgps}}\\
\mathfrak{B} &&& \text{random bipartite graph} \hfill  \mbox{Page \pageref{randbip}}\\
\mathfrak{B}^{*} &&& \text{bipartite complement of $\mathfrak{B}^{*}$} \hfill  \mbox{Page \pageref{bipcomp}}\\
\text{$B_m$} &&& \text{Braid group}\hfill \mbox{Page~\pageref{braidgroup}}\\
\text{$\mathcal{B}_{c}$} &&& \text{set of $n$-vertex complete}\\
\text{} &&& \text{bipartite graphs}\hfill \mbox{Page~\pageref{mathcal{B}_{c}}}\\
\text{$\BSym_{\alpha}(\Omega)$} &&& \text{group of all permutations moving}\\
\text{} &&& \text{fewer than $\alpha$ points}\\
\text{} &&& \text{for each infinite cardinal $\alpha\le|\Omega|$}\hfill \mbox{Page~\pageref{BSym}}\\
\mathbb{C} &&& \text{field of complex numbers}\hfill  \mbox{ throughout}\\
c(x, y) &&& \text{colour of the edge between}\\
\text{} &&& \text{vertices $x$ and $y$}\hfill \mbox{Page~\pageref{edgecol}}\\
\text{$\Cr$} &&& \text{cartesian product of groups}\hfill  \mbox{Page~\pageref{Crprod}}\\
\text{$\Co$} &&& \text{Conway group}\hfill  \mbox{Page~\pageref{Conway}}\\
\text{$\Cox$} &&& \text{Coxeter group}\hfill  \mbox{Page~\pageref{cox}}\\
\text{$\Cox^{+}$} &&& \text{even part of Coxeter group}\hfill  \mbox{Page~\pageref{cox+}}\\
CA &&& \text{Clifford algebra}\hfill  \mbox{ Chapter~\ref{chapFD}, Appendix~\ref{carsec}}\\
CG &&& \text{Clifford group}\hfill  \mbox{ Appendix~\ref{carsec}}\\
CG_{0}^{+} &&& \text{Reduced Clifford group}\hfill  \mbox{ Appendix~\ref{carsec}}\\
\text{$Cay(G, S)$} &&& \text{Cayley graph on $G$}\hfill \mbox{ Appendix~\ref{PermutationGroups}}\\
\text{$\mathcal{C}_2$} &&& \text{$2$-dimensional simplicial complex}\hfill \mbox{Page~\pageref{2dimcomp}}\\
\text{$\mathcal{C}$} &&& \text{set of induced cycles in $(\Gamma, f)$}\hfill \mbox{Page~\pageref{indcyc}}\\
\text{$\mathcal{C}$} &&& \text{set of edge colours}\hfill \mbox{Page~\pageref{colrset}}\\
\text{$\mathcal{C}_m$} &&& \text{set of $m$ edge colours}\hfill \mbox{ throughout}\\
char(K)  &&& \text{characteristic of field $K$}\hfill  \mbox{ throughout}\\
\text{$\card(A)$}  &&& \text{cardinality of a set $A$}\hfill  \mbox{Page~\pageref{card}}\\
\text{$\Comm_G (H)$} &&& \text{commensurability subgroup}\\
\text{} &&& \text{of $H < G$}\hfill \mbox{Page~\pageref{Comm_G (H)}}\\
D &&& \text{Latin square design}\hfill  \mbox{ Appendix~\ref{LoopTheory}}\\
d(x, y)  &&& \text{distance between $x$ and $y$}\hfill  \mbox{ throughout}\\
\dim(\cdot)  &&& \text{dimension}\hfill  \mbox{ throughout}\\
\text{$\dom(f)$} &&& \text{domain of a map $f$}\hfill \mbox{ throughout}\\
\text{$\diag$} &&& \text{diagonal matrix}\hfill \mbox{Page~\pageref{diag}}\\
\text{$\diam(\Gamma)$} &&& \text{denotes the diameter of graph $\Gamma$}\hfill \mbox{Page~\pageref{diam}}\\
\text{$\DA(\Gamma)$} &&& \text{duality automorphism group of $\Gamma$} \hfill \mbox{Page~\pageref{dsbgps}}\\
\text{$\Dr$} &&& \text{direct product of groups}\hfill  \mbox{Page~\pageref{Drprod}}\\
E(\Gamma) &&& \text{edge set of graph $\Gamma$}\hfill  \mbox{ throughout}\\
E = E(\Gamma, f) &&& \text{adjacency matrix of signed graph}\hfill  \mbox{Page~\pageref{admx}}\\
e_n &&& \text{number of Euler graphs}\\
\text{} &&& \text{on $n$-vertex graphs} \hfill  \mbox{Page~\pageref{tnswcl}}\\
\end{array} 
\end{displaymath}

\clearpage

\begin{displaymath}
\begin{array}{lllllll}
e &&& \text{edge set of graph}\hfill  \mbox{Page~\pageref{edgeofgraph}}\\
\mathbb{F}, GF(q) &&& \text{finite fields}\hfill  \mbox{ throughout}\\
f, \phi &&& \text{functions}\hfill  \mbox{ throughout}\\
\text{$f_m$~\label{f_m}} &&& \text{colouring function on edges of $\mathfrak{R}_{m,\omega}$}\hfill \mbox{Page~\pageref{fm}}\\
\text{$f_X$} &&& \text{signed switching class}\hfill \mbox{Page~\pageref{ssclass}}\\
F(t_1, \ldots, t_n) &&& \text{$n$-ary function symbol}\hfill  \mbox{ Appendix~\ref{ModelTheory}}\\
f(x_1, \ldots, x_n) &&& \text{polynomial}\hfill \mbox{ Appendix~\ref{Polynomials}}\\
\mathscr{F} &&& \text{filter}\hfill \mbox{ throughout}\\
\text{$\mathscr{F}_c$} &&& \text{cofinite filter}\hfill \mbox{ throughout}\\
\text{$\mathscr{F}_{\Gamma}$} &&& \text{neighbourhood filter on graph $\Gamma$}\hfill \mbox{Chapter~\ref{chap8}}\\
\mathcal{F}(X) &&& \text{field or $\sigma$-algebra on a space $X$}\hfill \mbox{ Appendix~\ref{CategoryandMeasure}}\\
\text{$\FSym(\Omega)$} &&& \text{finitary symmetric group of}\\
\text{} &&& \text{permutations on the set $\Omega$}\hfill  \mbox{ throughout}\\
\text{fix} &&& \text{fixed points}\\
\text{} &&& \text{of group action}\hfill \mbox{Chapter~\ref{homcochap}}\\
G &&& \text{group}\hfill  \mbox{ throughout}\\
G_{\{\Delta\}} &&& \text{setwise stabilizer of permutation group}\hfill  \mbox{ throughout}\\
G_{(\Delta)} &&& \text{pointwise stabilizer of permutation}\\
\text{} &&& \text{group}\hfill  \mbox{ throughout}\\
G_{\bar{\mu}} &&& \text{stabilizer of tuple} \hfill  \mbox{ throughout}\\
\text{$(G, \Omega)$} &&& \text{permutation group $G$}\\
\text{} &&& \text{on a set $\Omega$}\hfill  \mbox{ Appendix~\ref{PermutationGroups}}\\
g_n &&& \text{sequence of group elements}\hfill  \mbox{ throughout}\\
\mathcal{G}_{m,n} &&& \text{set of simple complete $n$-vertex}\\
\text{} &&& \text{$m$-edge graphs} \hfill \mbox{Page~\pageref{mathcal{G}_{m,n}}}\\
\text{$\tilde{G}$} &&& \text{Markov chain whose states are graphs } \hfill \mbox{ Appendix~\ref{FurtherDetails}}\\
\text{$\mathcal{\tilde{G}}$} &&& \text{set of all graph processes} \hfill \mbox{ Appendix~\ref{FurtherDetails}}\\ 
G^{cl} &&& \text{closure of group $G$} \hfill \mbox{Page~\pageref{G^{cl}}}\\
G^W_{m,n} &&& \text{extended switching groups}\hfill \mbox{Page~\pageref{extsw}}\\
G_0 &&& \text{amorphous set gauge group}\hfill   \mbox{ Chapter~\ref{rgconstr}}\\
\text{$\mathfrak{i}G(\mathfrak{R^{t}})$} &&& \text{group acting on $\mathfrak{R^{t}}$}\\
\text{} &&& \text{stabilizing colour $\mathfrak{i}$}\hfill  \mbox{Page~\pageref{stabi}}\\
\text{$3G(\mathfrak{R^{t}})$} &&& \text{group acting on $\mathfrak{R^{t}}$}\\
\text{} &&& \text{preceded by action of $C_3$}\hfill  \mbox{Page~\pageref{stabi}}\\
\text{$[G,G]$} &&& \text{commutator subgroup of $G$}\\
\text{} &&& \text{derived group of $G$}\hfill \mbox{Page~\pageref{dergp}}\\
\text{$\Gal$} &&& \text{Galois group} \hfill \mbox{ Appendix~\ref{NumberTheory}}\\
(G, S) &&& \text{group with triality} \hfill \mbox{Page~\pageref{(G,S)}}\\
G(M) &&& \text{group generated by L(x) \& R(x)} \hfill  \mbox{Page~\pageref{(G,S)}}\\
\text{$\gimel$} &&& \text{modular group}\hfill \mbox{Page~\pageref{gimel}}\\
\end{array} 
\end{displaymath}

\clearpage

\begin{displaymath}
\begin{array}{lllllll}
HOL(G) &&& \text{holomorph of a group G}\hfill  \mbox{ Appendix~\ref{PermutationGroups}}\\
Hx &&& \text{coset space; set of right}\\
\text{} &&& \text{cosets of $H$ in $G$}\hfill  \mbox{ Chapter~\ref{profilsec}}\\
H(\mathbb{Q}) &&& \text{autohomeomorphism group of $\mathbb{Q}$}\hfill  \mbox{ Page~\pageref{autohomeo}}\\
\mathfrak{H}_k &&& \text{homogeneous universal}\\
\text{} &&& \text{$K_k$-free graph}\hfill \mbox{Page~\pageref{K_k}}\\
\mathfrak{h}_3(\mathbb{O}) &&& \text{exceptional Jordan algebra}\hfill  \mbox{ Appendix~\ref{carsec}}\\
\text{$H^1(\mathcal{C}_2,\mathbb{Z}_{2})$} &&& \text{cohomology group}\hfill \mbox{Page~\pageref{cohgp}}\\
\text{$\Hyp(n)$} &&& \text{Hyperoctahedral group}\hfill \mbox{Page~\pageref{hyperoctahedral}}\\
\text{$\mathcal{H}_{\omega}$} &&& \text{direct limit of hypercubes}\hfill \mbox{Page~\pageref{hycube}}\\
\text{$\mathcal{H}$} &&& \text{completion of $\mathcal{H}_{\omega}$}\hfill \mbox{Page~\pageref{comhycube}}\\
\text{$\Hom$} &&& \text{vector space homomorphism}\hfill \mbox{Page~\pageref{Hom}}\\
\text{$\homo$} &&& \text{number of graph homomorphisms}\hfill \mbox{Page~\pageref{homo}}\\ 
\mathcal{I} &&& \text{index set}\hfill \mbox{ Appendix~\ref{CategoryandMeasure}}\\
I(Q) &&& \text{two-sided ideal}\hfill  \mbox{ Appendix~\ref{carsec}}\\
\text{$\im(f)$} &&& \text{image of a map $f$}\hfill  \mbox{Page~\pageref{image}}\\
\text{$\ind(f)$} &&& \text{index of a map $f$}\hfill  \mbox{Page~\pageref{index}}\\
\text{$\Inn(\cdot)$} &&& \text{group of inner automorphisms} \hfill \mbox{Page~\pageref{Inn}}\\
J &&& \text{rim of signed graph}\hfill \mbox{Page~\pageref{rim}}\\
j &&& \text{quadratic form}\hfill \mbox{ Appendix~\ref{carsec}}\\
J &&& \text{order 3 mapping}\hfill \mbox{ Appendix~\ref{carsec}}\\
K &&& \text{field}\hfill  \mbox{ throughout}\\
K[x_1, \ldots, x_n] &&& \text{polynomial algebra in $n$ variables}\\
\text{} &&& \text{over field $K$}\hfill \mbox{ throughout}\\
K[V] &&& \text{finite-dimensional $K$-vector space}\hfill  \mbox{ throughout}\\
K[V]^{G} &&& \text{ring of invariants of $K$-vector}\\
\text{} &&& \text{space under $G$-action}\hfill  \mbox{ throughout}\\
K[G] &&& \text{group ring}\hfill  \mbox{ throughout}\\
\text{$K_k$} &&& \text{complete graph on $k$ vertices}\hfill \mbox{Page~\pageref{K_k}}\\
\text{$\ker(\cdot)$} &&& \text{kernel of a map or permutation}\hfill  \mbox{ throughout}\\
L &&& \text{first-order language}\hfill  \mbox{ Appendix~\ref{ModelTheory}}\\
L &&& \text{loop}\hfill  \mbox{ Page~\pageref{Lloop}}\\
L(x) &&& \text{left-acting bijection}\\
\text{} &&& \text{on a loop element $x$}\hfill  \mbox{Page~\pageref{L(M)}}\\
\text{$\LMlt(Q)$} &&& \text{Left bijections}\\
\text{} &&& \text{on a Moufang loop $Q$}\hfill \mbox{Page~\pageref{LMlt}}\\
\mathbb{L} &&& \text{lattice}\hfill \mbox{Page~\pageref{mathbb{L}}}\\
\mathbb{L}_{E_8} &&& \text{root lattice of the $E_8$ Lie algebra}\hfill  \mbox{ throughout}\\
\mathbb{L}_{L} (\mathbb{L}_L^{\mathbb{C}}) &&& \text{(complex) Leech lattice}\hfill  \mbox{ throughout}\\
\mathcal{L} &&& \text{3-net}\hfill  \mbox{ Appendix~\ref{LoopTheory}}\\
\end{array} 
\end{displaymath}

\clearpage

\begin{displaymath}
\begin{array}{lllllll}
\mathsf{L} &&& \text{Lie algebra}\hfill \mbox{Page~\pageref{mathsfL}}\\
M, \mathcal{M} &&& \text{a structure}\hfill  \mbox{ throughout}\\
\mathcal{M}(X) &&& \text{set of Borel measures on $X$}\hfill \mbox{ Appendix~\ref{CategoryandMeasure}}\\
M_n &&& \text{random structure of size $n$}\hfill  \mbox{ Appendix~\ref{ModelTheory}}\\
M_s &&& \text{domain of a structure}\hfill  \mbox{ Appendix~\ref{ModelTheory}}\\
(M, b)_{b \in B}  &&& \text{expanded structure}\hfill  \mbox{ Appendix~\ref{ModelTheory}}\\
M_{2^n \times 2^n} (\mathbb{R}) &&& \text{a space of matrices}\hfill \mbox{ Appendix~\ref{carsec}}\\ 
\mathbb{M} &&& \text{Fischer-Griess Monster group}\hfill \mbox{ throughout}\\
\text{$\Mlt(Q)$} &&& \text{multiplication group}\\
\text{} &&& \text{of a Moufang loop $Q$}\hfill \mbox{Page~\pageref{Mlt}}\\
\mathbb{N} &&& \text{commutative monoid of}\\ 
\text{} &&& \text{natural numbers}\hfill  \mbox{ throughout}\\
N(g) &&& \text{vertex index of element $g$}\hfill \mbox{Page~\pageref{vxindex}}\\
N(v) &&& \text{neighbours of vertex $v$}\hfill \mbox{Page~\pageref{neighbvx}}\\
N(U) &&& \text{neighbourhood set for a}\hfill\\
\text{} &&& \text{finite vertex set $U$} \hfill \mbox{Page~\pageref{neighb}}\\
\mathfrak{N} &&& \text{Covington graph} \hfill \mbox{Page~\pageref{cogr}}\\
\mathbf{N} &&& \text{class of N-free graphs} \hfill \mbox{Page~\pageref{clnfr}}\\
\mathbf{N}_{f} &&& \text{class of finite N-free graphs} \hfill \mbox{Page~\pageref{clnfr}}\\
\mathcal{N} &&& \text{permutation model}\hfill  \mbox{ Chapter~\label{rgconstr}}\\
\text{$\Nuc(Q)$} &&& \text{nucleus of a Moufang loop $Q$}\hfill \mbox{Page~\pageref{Nuc}}\\
\NB(\Omega) &&& \text{group of near bijections of a set $\Omega$} \hfill \mbox{Page~\pageref{NB}}\\
\NS(\Omega) &&& \text{group of near symmetries of a set $\Omega$} \hfill \mbox{Page~\pageref{NS}}\\
\text{$NA(\Gamma)$} &&& \text{group of near automorphisms of a graph $\Gamma$} \hfill \mbox{Page~\pageref{NA}}\\
\text{$NE(\Gamma)$} &&& \text{group of equivalence classes of}\\
\text{} &&& \text{elements of $\NA(\Gamma)$} \hfill \mbox{Page~\pageref{NE}}\\
\mathbb{O} &&& \text{nonassociative normed division}\hfill\\
\text{} &&& \text{algebra of octonions}\hfill \mbox{ throughout}\\
O, O_i &&& \text{open sets}\hfill \mbox{ Appendices~\ref{CategoryandMeasure}}\\
\mathbb{OP}^2 &&& \text{Moufang projective plane}\hfill \mbox{ Appendix~\ref{carsec}}\\
O(n, K) &&& \text{$n$-dimensional orthogonal}\\
\text{} &&& \text{group over field $K$} \hfill  \mbox{ throughout}\\
O(\infty) &&& \text{direct limit of orthogonal groups $O(n)$}\hfill \mbox{ Appendix~\ref{carsec}}\\
\mathfrak{O} &&& \text{point at infinity}\hfill \mbox{Page~\pageref{ptatinf}}\\
\text{$\Out(\Gamma)$} &&& \text{group of outer automorphisms of $\Gamma$} \hfill \mbox{Page~\pageref{outaut}}\\
Pr(\sigma) &&&  \text{$\lim_{n \to \infty} Pr_n(\sigma)$}\hfill   \mbox{ Appendix~\ref{FurtherDetails}}\\
\text{$\mathbf{P}^{1}(K)$} &&& \text{projective line over a field $K$} \hfill \mbox{Page~\pageref{projln}}\\
Pr_n(\sigma) &&& \text{probability sentence $\sigma$}\\
\text{} &&& \text{is true in $M_n$} \hfill  \mbox{ Appendix~\ref{FurtherDetails}}\\
\end{array} 
\end{displaymath}

\clearpage

\begin{displaymath}
\begin{array}{lllllll}
\text{P} &&& \text{parity on a complete}\\
\text{} &&& \text{$3$-coloured graph} \hfill  \mbox{Page~\pageref{P-equivalent}}\\
\text{$P_m$} &&& \text{parity on an $m$-coloured graph}\\
\text{} &&& \text{generalizing P parity} \hfill  \mbox{Page~\pageref{P_m-equivalent}}\\
\text{$\mathsf{P}_m$} &&& \text{parity on class of complete}\\
\text{} &&& \text{$m$-coloured graphs} \hfill  \mbox{Page~\pageref{mathsf{P}_m}}\\
\mathcal{P}  &&& \text{property}\hfill  \mbox{ throughout}\\
\text{$p$} &&& \text{finite partial bijection}\\
\text{} &&& \text{from $\mathfrak{R}_{m,\omega}$ to $\mathfrak{R}_{m,\omega}$} \hfill \mbox{Page~\pageref{fpb}}\\
\text{$p$} &&& \text{part of a partition}\hfill  \mbox{Page~\pageref{partp}}\\
\mathcal{P}(\Omega) &&& \text{power set of $\Omega$}\hfill  \mbox{ throughout}\\
\text{$PSL$} &&& \text{projective special linear group} \hfill  \mbox{ throughout}\\
\text{$PSO$} &&& \text{projective special}\\
\text{} &&& \text{orthogonal group} \hfill  \mbox{ throughout}\\
\mathbb{Q} &&& \text{field of rational numbers}\hfill  \mbox{ throughout}\\
\mathbb{Q}[i] &&& \text{Gaussian algebraic number field}\hfill \mbox{ throughout}\\
q  &&& \text{quadratic form}\hfill  \mbox{ Appendix~\ref{carsec}}\\
Q &&& \text{Moufang loop}\hfill  \mbox{Page~\pageref{QMoufangloop}}\\
q &&& \text{size of Hamming space alphabet}\hfill \mbox{Page~\pageref{que}}\\
\mathbb{R} &&& \text{field of real numbers}\hfill  \mbox{ throughout}\\
\mathbb{R}\{x_1, \ldots, x_n\} &&& \text{free $\mathbb{R}$-algebra}\hfill  \mbox{ Appendix~\ref{carsec}}\\
r(x_1, \ldots, x_n) &&& \text{$n$-ary relation symbol}\hfill  \mbox{ Appendix~\ref{ModelTheory}}\\
R(t_1, \ldots, t_n) &&& \text{$n$-ary relation symbol}\hfill  \mbox{ Appendix~\ref{ModelTheory}}\\
\mathfrak{R} &&& \text{2 edge-colour countably infinite}\\
\text{} &&& \text{random graph} \hfill  \mbox{Page \pageref{mathfrak{R}}}\\
\mathfrak{R^{t}} &&& \text{triality graph, 3 edge-colour countably}\\ 
\text{} &&& \text{infinite random graph} \hfill \mbox{Page \pageref{mathfrak{Rt}}}\\
\mathfrak{r}, \mathfrak{b}, \mathfrak{g} &&& \text{red, blue and green}\hfill\\
\text{} &&& \text{triality graph edge colours} \hfill \mbox{ throughout}\\
\mathfrak{R}_{m,\omega} &&& \text{$m$-colour countably infinite}\\
\text{} &&& \text{random graph} \hfill  \mbox{Page \pageref{mathfrak{R_mn}}}\\
\mathfrak{R}^{v} &&& \text{random graph with coloured vertices} \hfill \mbox {Page~\pageref{rangphcolv}}\\
r_n &&& \text{proportion of $n$-vertex}\\
\text{} &&& \text{graphs with specified property} \hfill \mbox{ Appendix~\ref{FurtherDetails}}\\
R(x) &&& \text{right-acting bijection on a loop}\hfill  \mbox{ Page~\pageref{L(M)}}\\
R_k &&& \text{relations in $\Gamma$}\hfill \mbox{Page~\pageref{Relations}}\\
\text{$\range(f)$} &&& \text{range of a map $f$}\hfill \mbox{ throughout}\\
\mathfrak{R}_N &&& \text{universal graph of fixed}\\
\text{} &&& \text{finite diameter $N$} \hfill  \mbox{Page \pageref{mathfrak{RN}}}\\
\end{array} 
\end{displaymath}

\clearpage

\begin{displaymath}
\begin{array}{lllllll}
\text{$\RMlt(Q)$} &&& \text{Right bijections}\\
\text{} &&& \text{on a Moufang loop $Q$}\hfill \mbox{Page~\pageref{RMlt}}\\
S &&& \text{a generating set for $\Sym(\Omega)$}\hfill  \mbox{ Appendix~\ref{PermutationGroups}}\\
S &&& \text{a structure}\hfill  \mbox{ Appendix~\ref{ModelTheory}}\\
S &&& \text{spinor space of $CA$}\hfill  \mbox{ Appendix~\ref{carsec}}\\
\text{S} &&& \text{switching parity on a complete}\\
\text{} &&& \text{$3$-coloured graph} \hfill  \mbox{Page~\pageref{S-equivalent}}\\
S &&& \text{subset of $\mathbb{N}$}\hfill  \mbox{Page~\pageref{SsubsetN}}\\
S^{\pm} &&& \text{half-spin representations of $\mathfrak{so}(n, \mathbb{R})$}\hfill  \mbox{Appendix~\ref{carsec}}\\
\mathcal{S}  &&& \text{a set of graph paths}\hfill  \mbox{ Appendix~\ref{TheoryofRelationalStructures}}\\
\mathcal{S}_{\omega}  &&& \text{group inducing finitary edge-colour}\\
\text{} &&& \text{permutations}\hfill  \mbox{ Chapter~\ref{chap3}}\\
\mathcal{S}(U, V) &&& \text{set of graphs satisfying $(*)$}\hfill  \mbox{ Appendix~\ref{FurtherDetails}}\\
\text{$S_{m,n}$} &&& \text{switching group} \hfill \mbox{Page~\pageref{switchinggroup}}\\
\text{$S_{m, \omega}^{cl}$} &&& \text{closure of $S_{m, \omega}$} \hfill \mbox{Page~\pageref{clswitchinggroup}}\\
\soc(G) &&& \text{socle of}\\  
\text{} &&& \text{a group G}\hfill  \mbox{ Appendix~\ref{PermutationGroups}, Chapter~\ref{chapFD}}\\
\text{$Spin$} &&& \text{spin group} \hfill  \mbox{ throughout}\\
\text{$\SA(\Gamma)$} &&& \text{group of switching automorphisms of $\Gamma$} \hfill \mbox{Page~\pageref{groupofswaut}}\\
\sqrt{a} &&& \text{a square-root set in a}\\
\text{} &&& \text{group}\hfill  \mbox{ Appendix~\ref{PermutationGroups}, Chapter~\ref{homcochap}}\\
\text{$\Sym(\Omega)$} &&& \text{symmetric group of}\\
\text{} &&& \text{permutations on the set $\Omega$}\hfill  \mbox{ throughout}\\
\text{$\mathsf{S}_{n}$} &&& \text{group of switchings on $\mathfrak{R}_{m, n}$}\\
\text{} &&& \text{as products of elementary switchings}\hfill \mbox{Page~\pageref{mathsf{S}_{n}}}\\
\text{$\Suz$} &&& \text{Suzuki group}\hfill  \mbox{Page~\pageref{Suz}}\\
t_i &&& \text{a term}\hfill  \mbox{ Appendix~\ref{ModelTheory}}\\
T &&& \text{theory}\hfill  \mbox{ Appendix~\ref{ModelTheory}}\\
T &&& \text{topology}\hfill  \mbox{ Appendix~\ref{TopologyinPermutationGroups}}\\
T(V) &&& \text{tensor algebra}\hfill  \mbox{ Appendix~\ref{carsec}}\\
\Th(M) &&& \text{theory of a structure M}\hfill  \mbox{ Appendix~\ref{ModelTheory}}\\
\Th(M, a)_{a \in M}  &&& \text{type over a model $M$}\hfill   \mbox{ throughout}\\
\tp_{M}(\bar{a}) &&& \text{type of $\bar{a} \in M_{s}^k$}\hfill  \mbox{ Appendix~\ref{ModelTheory}}\\
\tp_{M}(\bar{a} / B) &&& \text{type of $\bar{a}$ in $(M, b)_{b \in B}$}\hfill  \mbox{ Appendix~\ref{ModelTheory}}\\
\mathcal{T}(\Gamma) &&& \text{two-graph of $\Gamma$}\hfill  \mbox{ Page~\pageref{two-graph}}\\
t_n &&& \text{number of switching classes}\\
\text{} &&& \text{for $n$-vertex graphs} \hfill  \mbox{Page~\pageref{tnswcl}}\\
\text{$\TAut(\Gamma)$} &&& \text{group of all automorphisms of $\Gamma$} \hfill \mbox{Page~\pageref{colrset}}\\
T_d(X) &&& \text{natural topology}\hfill \mbox{ Appendix~\ref{TopologyinPermutationGroups}}\\
\text{$\Ts$} &&& \text{Tsaranov group}\hfill  \mbox{Page~\pageref{Ts}}\\
\end{array} 
\end{displaymath}

\clearpage

\begin{displaymath}
\begin{array}{lllllll}
\text{$(T, r)$} &&& \text{rooted tree}\hfill \mbox{Page~\pageref{roottree}}\\
\mathcal{T}, \mathcal{T}^*, \mathcal{T}^\dag &&& \text{topologies on $\mathfrak{R}$}\hfill  \mbox{Page~\pageref{Ttop}}\\
\text{$t(\gamma, \Gamma)$} &&& \text{graph homomorphism density}\hfill  \mbox{Page~\pageref{tgraph}}\\
U, V &&& \text{finite disjoint sets of vertices}\\
\text{} &&& \text{in $\mathfrak{R}$}\hfill   \mbox{ throughout}\\
U, V, W &&& \text{finite disjoint sets of vertices}\\
\text{} &&& \text{in $\mathfrak{R^{t}}$}\hfill   \mbox{ throughout}\\
U &&& \text{amorphous set}\hfill   \mbox{ Chapter 8}\\
\mathcal{U} &&&  \text{topology on $\mathfrak{R}$}\hfill  \mbox{Page~\pageref{U}}\\
V_n &&& \text{$n$th homogeneous component}\\
\text{} &&& \text{of vector space}\hfill  \mbox{ Appendix~\ref{Polynomials}}\\
V^{*} &&& \text{dual vector space}\hfill  \mbox{ throughout}\\
V(\Gamma) &&& \text{vertex set of graph $\Gamma$}\hfill  \mbox{ throughout}\\
v &&& \text{vertex in $V(\Gamma)$}\hfill  \mbox{ throughout}\\
v^{\mathfrak{i}} &&& \text{vertex with colour $\mathfrak{i}$}\hfill  \mbox{ Chapters~\ref{fstchap},~\ref{chap2}}\\
V^{\natural} &&& \text{Monster module}\hfill  \mbox{ throughout}\\
W_i &&& \text{unary relation}\hfill \mbox{ Appendix~\ref{FurtherDetails}, Page~\pageref{unaryr}}\\
\text{$\W(\Gamma, f)$} &&& \text{Weyl group of root lattice $\mathbb{L}(\Gamma, f)$}\hfill \mbox{Page~\pageref{weylgp}}\\
W_{\mathbb{L}_{E_8}} &&& \text{Weyl group of the $E_8$ Lie algebra}\hfill  \mbox{ throughout}\\
X &&& \text{set}\hfill  \mbox{ throughout}\\
x_i &&& \text{elements of a set}\hfill  \mbox{ throughout}\\
X &&& \text{a space}\hfill \mbox{ Appendix~\ref{CategoryandMeasure}}\\
(X, d) &&& \text{metric space}\hfill \mbox{ throughout}\\
(X, T) &&& \text{topological space}\hfill \mbox{ Appendix~\ref{CategoryandMeasure}}\\
(X, T_d(X)) &&& \text{topological space}\hfill \mbox{ Appendix~\ref{CategoryandMeasure}}\\
(X, \mathcal{F}(X), \mu) &&& \text{probability space}\hfill \mbox{ Appendix~\ref{CategoryandMeasure}}\\ 
\mathcal{Y}(M)  &&& \text{a set of paths in a}\\
\text{} &&& \text{structure $M$}\hfill  \mbox{ Appendix~\ref{TheoryofRelationalStructures}}\\
\mathbb{Z} &&& \text{ring of integers}\hfill  \mbox{ throughout}\\
\mathbb{Z}_p &&& \text{ring of integers modulo $p$}\hfill  \mbox{ throughout}\\
\mathbb{Z}[\omega] &&& \text{ring of Eisenstein integers}\hfill \mbox{ Appendix~\ref{carsec}, Chapter~\ref{homcochap}}\\
\end{array} 
\end{displaymath}

\clearpage

\begin{displaymath}
\begin{array}{lllllll}
\alpha &&& \text{main anti-automorphism of $CA$} \hfill \mbox{ Appendix~\ref{carsec}}\\
\beta &&& \text{automorphism or}\\
\text{} &&& \text{almost automorphism} \hfill \mbox{ Chapter~\ref{homcochap}}\\
\Gamma_{C} &&& \text{countable universal homogeneous}\\
\text{} &&& \text{$C$-coloured graph} \hfill \mbox{Page~\pageref{Gamma_C}}\\
(\Gamma, f) &&& \text{signed graph}\hfill  \mbox{Page~\pageref{Gammaf}}\\
\Gamma(n, p) &&& \text{random graph on $n$ vertices}\\
\text{} &&& \text{with connection probability $p$} \hfill \mbox{ Appendix~\ref{FurtherDetails}}\\
\text{$\Gamma(S), \Gamma(G)$} &&&
\text{Cayley graphs} \hfill \mbox{Chapter~\ref{homcochap}}\\
\text{$[\Gamma(G)]^2$} &&& \text{2-complex of Cayley graph} \hfill \mbox{Page~\pageref{Gamma(G)]^2}}\\
\gamma &&& \text{first invariant of $\mathcal{T}(\Gamma)$} \hfill \mbox{Page~\pageref{fstinv}}\\
\text{$\symd$} &&& \text{symmetric difference} \hfill \mbox{ throughout}\\ 
\text{$\Delta_1, \Delta_2$} &&& \text{bipartite blocks of $\mathfrak{B}$} \hfill
\mbox{Page~\pageref{deltabip}}\\
\zeta(\cdot)  &&& \text{centre of a group or loop} \hfill \mbox{ throughout}\\
\theta &&& \text{homomorphism (variously iso-}\\
\text{} &&& \text{mono- \& epimorphism)} \hfill \mbox{ throughout}\\
\Theta &&& \text{theta function} \hfill \mbox{Page~\pageref{Thetafn}}\\
\iota &&& \text{inclusion map} \hfill \mbox{Page~\pageref{iota}}\\
\iota_i &&& \text{interpretations of}\\
\text{} &&& \text{structure of $M$ in $N$}\hfill \mbox{Page~\pageref{iotaint}}\\
\kappa &&& \text{uncountable cardinal}\hfill \mbox{Page~\pageref{kappa}}\\
\Lambda &&& \text{lattice or}\\
\text{} &&& \text{projective class of lattices}\hfill   \mbox{Page~\pageref{lamlattice}}\\
\Lambda &&& \text{ordered abelian group}\hfill \mbox{Page~\pageref{ordAbgp}}\\
\mu, \nu &&& \text{measures of a set}\hfill \mbox{ Appendix~\ref{CategoryandMeasure}}\\
\xi &&& \text{switching operation}\hfill \mbox{Page~\pageref{switop}}\\
\pi &&& \text{permutation}\hfill   \mbox{Page~\pageref{piperm}}\\
\pi &&& \text{set partition}\hfill   \mbox{Page~\pageref{partition}, Appendix~\ref{PermutationGroups}}\\
\text{$\rho$} &&& \text{representation homomorphism}\hfill \mbox{ Appendix~\ref{carsec}}\\
\text{$\rho_s$} &&& \text{spinor representation}\hfill \mbox{ Appendix~\ref{carsec}}\\
\text{$\rho$} &&& \text{3-cycle of Group with Triality}\hfill \mbox{ Chapter~\ref{homcochap} \S4-5}\\
\text{$\sigma, \sigma_i, \sigma'$} &&& \text{cyclic automorphisms}\hfill \mbox{ Chapter~\ref{homcochap} \S1}\\
\text{$\sigma$} &&& \text{2-cycle of Group with Triality}\hfill \mbox{ Chapter~\ref{homcochap} \S4-5}\\
\sigma &&& \text{a sentence}\hfill  \mbox{ Appendix~\ref{ModelTheory}}\\
\Sigma &&& \text{set of sentences}\hfill  \mbox{ Appendix~\ref{ModelTheory}}\\
\sigma_{c,d,Y} &&& \text{switching operation on graph $Y$}\\
\text{} &&&  \text{of edge-colours $c$ \& $d$}\hfill \mbox {Page~\pageref{sigmacdy}}\\
\end{array} 
\end{displaymath}

\clearpage

\begin{displaymath}
\begin{array}{lllllll}
\text{$\sigma$} &&& \text{isometry}\hfill \mbox{ Chapter 10}\\
\tau &&& \text{3-cycle (automorphism)}\hfill  \mbox{Page~\pageref{tau}}\\ 
\tau_C &&& \text{Cartan triality 3-cycle}\hfill  \mbox{Page~\pageref{tau}}\\
\tau_x &&& \text{partial automorphism}\hfill  \mbox{ Appendix~\ref{LoopTheory}}\\
\tau_{CA} &&& \text{triality outer}\\
\text{} &&& \text{automorphism of $P\Omega^{+}(8, \mathbb{R})$}\hfill \mbox{ Appendix~\ref{carsec}}\\
\phi_{M} &&& \text{profile of relational structure $M$}\hfill \mbox{ Appendix~\ref{Polynomials}}\\
\Phi(g)  &&& \text{triality identity for $g \in G$}\hfill \mbox{Page~\pageref{trid}}\\
\Upsilon &&& \text{SML graph}\hfill  \mbox{Page~\pageref{smlgraph}}\\ 
\chi_{\pi_{1}}(\pi_{2}) &&& \text{cubic residue character}\\ 
\text{} &&& \text{of $\pi_{2}$ modulo $\pi_{1}$}\hfill \mbox{ Appendix~\ref{NumberTheory}}\\
\Psi_k &&& \text{axiom resembling}\\
\text{} &&& \text{the $k$-e.c. property}\hfill  \mbox{ Appendix~\ref{FurtherDetails}}\\
\omega &&& \text{cube root of unity}\hfill \mbox{Chapter~\ref{homcochap}}\\
\Omega &&& \text{a set}\hfill  \mbox{ throughout}\\
\Omega^k &&& \text{$k$-tuple of elements}\\
\text{} &&& \text{from a set $\Omega$}\hfill \mbox{ throughout}\\
\Omega^{+}(8, \mathbb{R}) &&& \text{Reduced orthogonal group, }\\
\text{} &&&\text{commutator subgroup of $O^{+}(8, \mathbb{R})$}\hfill \mbox{ Appendix~\ref{carsec}}\\
\end{array} 
\end{displaymath}

\printindex

\end{document}

%% file: test.latex
\setlength{\unitlength}{3947sp}%
\begingroup\makeatletter\ifx\SetFigFont\undefined%
\gdef\SetFigFont#1#2#3#4#5{%
  \reset@font\fontsize{#1}{#2pt}%
  \fontfamily{#3}\fontseries{#4}\fontshape{#5}%
  \selectfont}%
\fi\endgroup%
\begin{picture}(4275,5032)(676,-4544)
\thinlines
\put(1651,-211){\oval(1382,1382)}
\put(4201,-286){\oval(1382,1382)}
\put(4127,-973){\oval(1382,1382)}
\put(4051,-3511){\oval(1382,1382)}
\put(3526,-436){\makebox(1.6667,11.6667){\SetFigFont{5}{6}{\rmdefault}{\mddefault}{\updefault}.}}
\put(2626,-436){\makebox(1.6667,11.6667){\SetFigFont{5}{6}{\rmdefault}{\mddefault}{\updefault}.}}
\put(2401,-436){\vector( 1, 0){1125}}
\put(4051,-1786){\vector( 0,-1){900}}
\put(3751,-661){\line( 5, 3){518.382}}
\put(3901,-811){\line( 5, 3){606.618}}
\put(4126,-961){\line( 4, 3){516}}
\put(1568,-793){\line( 5, 2){530.172}}
\put(1364,-772){\line( 3, 1){832.500}}
\put(1276,-661){\line( 5, 2){969.828}}
\put(4183,-4087){\line( 4, 3){300}}
\put(3901,-4111){\line( 5, 4){750}}
\put(3751,-4036){\line( 4, 3){900}}
\put(2701,-736){\makebox(0,0)[lb]{\smash{\SetFigFont{12}{14.4}{\rmdefault}{\mddefault}{\updefault}$e_1$}}}
\put(1276,-1111){\makebox(0,0)[lb]{\smash{\SetFigFont{12}{14.4}{\rmdefault}{\mddefault}{\updefault}$\Omega_1$}}}
\put(4951,-436){\makebox(0,0)[lb]{\smash{\SetFigFont{12}{14.4}{\rmdefault}{\mddefault}{\updefault}$e_1(\Omega_1)$}}}
\put(4876,-1186){\makebox(0,0)[lb]{\smash{\SetFigFont{12}{14.4}{\rmdefault}{\mddefault}{\updefault}$\Omega_2$}}}
\put(676,-4486){\makebox(0,0)[lb]{\smash{\SetFigFont{12}{14.4}{\rmdefault}{\mddefault}{\updefault}$\Omega_3 = \Omega_1 \cap e_1^{-1}(\Omega_2 \cap e_1(\Omega_1))$}}}
\put(4201,-2161){\makebox(0,0)[lb]{\smash{\SetFigFont{12}{14.4}{\rmdefault}{\mddefault}{\updefault}$e_2$}}}
\put(4876,-3586){\makebox(0,0)[lb]{\smash{\SetFigFont{12}{14.4}{\rmdefault}{\mddefault}{\updefault}$e_2(\Omega_2)$}}}
\end{picture}